\pdfoutput=1
\documentclass[a4paper]{willowtreebook}
\Title{Introduction to Cartan Geometries}
\Author{\texorpdfstring{Benjamin \scotsMc{}Kay}{Benjamin McKay}}
\BibliographyFile{introduction-to-Cartan-geometries}
\makeatletter
\DeclareRobustCommand{\scotsMc}{\scotsMcx{c}}
\DeclareRobustCommand{\scotsMC}{\scotsMcx{\textsc{c}}}
\DeclareRobustCommand{\scotsMcx}[1]{%
  M%
  \raisebox{\dimexpr\fontcharht\font`M-\height}{%
    \check@mathfonts\fontsize{\sf@size}{0}\selectfont
    \kern.3ex\underline{\kern-.3ex #1\kern-.3ex}\kern.3ex
  }%
}
\expandafter\def\expandafter\@uclclist\expandafter{%
  \@uclclist\scotsMc\scotsMC
}
\makeatother

\DeclareUnicodeCharacter{2227}{\ensuremath{\wedge}}
\DeclareUnicodeCharacter{2207}{\ensuremath{\nabla}} %∇
\DeclareUnicodeCharacter{00A3}{\ensuremath{\LieDer}} %£
\DeclareUnicodeCharacter{02E9}{\ensuremath{\hook}} %˩
\DeclareUnicodeCharacter{2308}{\ensuremath{\lceil}}
\DeclareUnicodeCharacter{2309}{\ensuremath{\rceil}}
\DeclareUnicodeCharacter{230A}{\ensuremath{\lfloor}}
\DeclareUnicodeCharacter{230B}{\ensuremath{\rfloor}}
\DeclareUnicodeCharacter{03F5}{\ensuremath{\varepsilon}}
\DeclareUnicodeCharacter{2282}{\ensuremath{\subset}}%⊂
\DeclareUnicodeCharacter{03C6}{\ensuremath{\phi}}%ϕ
\DeclareUnicodeCharacter{03D5}{\ensuremath{\varphi}}%ϕ
\DeclareUnicodeCharacter{03B9}{\ensuremath{\iota}} %ι
\DeclareUnicodeCharacter{2297}{\otimes}
\DeclareUnicodeCharacter{1D524}{\ensuremath{\LieG}} %𝔤
\DeclareUnicodeCharacter{1D525}{\ensuremath{\LieH}} %𝔥
\DeclareUnicodeCharacter{1D529}{\ensuremath{\LieL}} %𝔩
\DeclareUnicodeCharacter{1D52A}{\ensuremath{\LieM}} %𝔪
\DeclareUnicodeCharacter{03C9}{\ensuremath{\omega}}%ω
\DeclareUnicodeCharacter{222B}{\int}%
\DeclareUnicodeCharacter{221E}{\infty}%
\DeclareUnicodeCharacter{2211}{\sum}%
\DeclareUnicodeCharacter{03B1}{\alpha}%
\DeclareUnicodeCharacter{03B2}{\beta}%
\DeclareUnicodeCharacter{02B9}{\ensuremath{'}}%
\newunicodechar{ℝ}{\mathbb{R}}
\newunicodechar{ϕ}{\varphi}
\newunicodechar{⟨}{\left<}
\newunicodechar{⟩}{\right>}
\usepackage{hyperref}
\usepackage{mathrsfs}
\usepackage{standalone}
\usepackage{graphicx}
\usepackage{epigraph-keys}
\usepackage{xstring}
\usepackage{bm}
\usepackage{tikz-cd}
\usepackage{pgfplots}
\pgfplotsset{compat=1.17}
\usetikzlibrary{%
decorations.markings,
decorations.pathmorphing,
backgrounds,
shapes,
arrows,
calc,
hobby}
\usetikzlibrary {arrows.meta,bending}
\tikzset
  {midarrow/.style={decoration={markings,mark=at position 0.5 with
     {\arrow[xshift=2pt]{Latex[length=4pt,#1]}}},postaction={decorate}}
  }
\usepackage{nicematrix}
\NewDocumentCommand\bark{}{\ensuremath{\bar{k}}}
\NewDocumentCommand\C{e{^}}%
{%
	\ensuremath%
	{%
		\IfValueTF{#1}%
		{%
			\mathbb{C}^{#1}%
		}%
		{%
			\mathbb{C}%
		}%
	}%
}%
\NewDocumentCommand\R{e{^}}%
{%
	\ensuremath%
	{%
		\IfValueTF{#1}%
		{%
			\mathbb{R}^{#1}%
		}%
		{%
			\mathbb{R}%
		}%
	}%
}%
\NewDocumentCommand\Z{e{^}}%
{%
	\ensuremath%
	{%
		\IfValueTF{#1}%
		{%
			\mathbb{Z}^{#1}%
		}%
		{%
			\mathbb{Z}%
		}%
	}%
}%
\NewDocumentCommand\Quat{e{^}}%
{%
	\ensuremath%
	{%
		\IfValueTF{#1}%
		{%
			\mathbb{H}^{#1}%
		}%
		{%
			\mathbb{H}%
		}%
	}%
}%
\NewDocumentCommand\Romanbar{m}{%
\relax\ifmmode%
\IfStrEqCase{#1}{
{1}{I}%
{2}{I\kern-2.5pt I}%
{3}{I\kern-2.5pt I\kern-2.5pt I}%
{4}{I\kern-1pt V}%
{5}{V}%
{6}{V\kern-2.6pt I}%
{7}{V\kern-2.6pt I\kern-2.5pt I}%
{8}{V\kern-2.6pt I\kern-2.5pt I\kern-2.5pt I}%
{9}{I\kern-2.5pt X}%
{10}{X}%
{11}{X\kern-2.8pt I}%
}%
\else%
\IfStrEqCase{#1}{
{1}{I}%
{2}{I\kern-0.57pt I}%
{3}{I\kern-0.57pt I\kern-0.57pt I}%
{4}{I\kern-1pt V}%
{5}{V}%
{6}{V\kern-.49pt I}%
{7}{V\kern-.49pt I\kern-0.57pt I}%
{8}{V\kern-.49pt I\kern-0.57pt I\kern-0.57pt I}%
{9}{I\kern-0.66pt X}%
{10}{X}%
{11}{X\kern-1pt I}%
}%
\fi%
}
\NewDocumentCommand\LT{m}{\ell_{#1}}
\NewDocumentCommand\RT{m}{r_{#1}}
\newcommand{\shapeOp}{\ensuremath{\Romanbar{2}}}
\NewDocumentCommand\Proj{m}{\ensuremath{\mathbb{P}^{#1}}}
\NewDocumentCommand\RP{m}{\ensuremath{\mathbb{RP}^{#1}}}
\NewDocumentCommand\CP{m}{\ensuremath{\mathbb{CP}^{#1}}}
\NewDocumentCommand\GL{m}{\ensuremath{\operatorname{GL}_{#1}}}
\NewDocumentCommand\Lie{mo}{%
\ensuremath{%%
\IfValueTF{#2}%
{%%%
{\mathfrak{\MakeLowercase{#1}}}_{#2}
}%%%
{%%%
\mathfrak{\MakeLowercase{#1}}
}%%%
}%%
}%
\NewDocumentCommand\CompleteCone{o}{\IfValueTF{#1}{\mathscr{K}_{#1}}{\mathscr{K}}}
\NewDocumentCommand\barCompleteCone{o}{\IfValueTF{#1}{\bar{\mathscr{K}}_{#1}}{\bar{\mathscr{K}}}}
\NewDocumentCommand\LieGL{m}{\Lie{GL}[#1]}
\NewDocumentCommand\LieSO{m}{\Lie{SO}[#1]}
\NewDocumentCommand\SL{m}{\ensuremath{\operatorname{SL}_{#1}}}
\NewDocumentCommand\LieSL{m}{\Lie{SL}[#1]}
\NewDocumentCommand\PSL{m}{\ensuremath{\mathbb{P}\!\operatorname{SL}_{#1}}}
\NewDocumentCommand\PGL{m}{\ensuremath{\mathbb{P}\!\operatorname{GL}_{#1}}}
\NewDocumentCommand\PO{m}{\ensuremath{\mathbb{P}\!\operatorname{O}_{#1}}}
\NewDocumentCommand\SU{m}{\ensuremath{\operatorname{SU}_{#1}}}
\NewDocumentCommand\SUL{}{\ensuremath{\operatorname{SU}_2^{\ell}}}
\NewDocumentCommand\SUR{}{\ensuremath{\operatorname{SU}_2^r}}

\NewDocumentCommand\LieSUR{}{\ensuremath{\mathfrak{su}_2^r}}
\NewDocumentCommand\Aff{m}{\operatorname{Aff}_{#1}}
\NewDocumentCommand\Euc{m}{\operatorname{Euc}_{#1}}
\NewDocumentCommand\Lm{smm}{\ensuremath{\Lambda^{\IfStrEq{#2}{*}{\complex}{#2}}\IfBooleanTF{#1}{\!\left(#3\right)}{#3}}}
\NewDocumentCommand\transpose{m}{{#1}^t}
\DeclareMathOperator{\Ad}{Ad}
\DeclareMathOperator{\ad}{ad}
\NewDocumentCommand\LieDer{}{\ensuremath{\mathcal L}}
\NewDocumentCommand\hook{}{\ensuremath{\mathbin{ \hbox{\vrule height1.4pt
        width4pt depth-1pt \vrule height4pt width0.4pt depth-1pt}}}}
\NewDocumentCommand\SO{m}{\ensuremath{\operatorname{SO}_{#1}}}
\NewDocumentCommand\Un{m}{\ensuremath{\operatorname{U}_{#1}}}
\NewDocumentCommand\PU{m}{\ensuremath{\mathbb{P}\!\operatorname{U}_{#1}}}
\NewDocumentCommand\UnR{}{\ensuremath{\operatorname{U}^r_1}}
\NewDocumentCommand\LieUn{m}{\ensuremath{\mathfrak{u}_{#1}}}
\NewDocumentCommand\LieUnR{}{\ensuremath{\mathfrak{u}^r_1}}
\NewDocumentCommand\Orth{m}{\ensuremath{\operatorname{O}_{#1}}}
\NewDocumentCommand\CO{m}{\ensuremath{\operatorname{CO}_{#1}}}
\NewDocumentCommand\Spin{m}{\ensuremath{\operatorname{Spin}_{#1}}}
\NewDocumentCommand\ip{mm}{\left<#1,#2\right>}
\NewDocumentCommand\Hom{mm}{\ensuremath{\operatorname{Hom}({#1},{#2})}}
\NewDocumentCommand\homotopyGroup{mm}{\ensuremath{\pi_{#1}({#2})}}
\NewDocumentCommand\fundamentalGroup{m}{\ensuremath{\homotopyGroup{1}{#1}}}
\NewDocumentCommand\complex{}{\ensuremath{\scriptscriptstyle{\bullet}}}
\NewDocumentCommand\MakeLie{m}{\expandafter\def\csname Lie#1\endcsname{\Lie{#1}}}
\NewDocumentCommand\lb{smm}{\ensuremath{\left[{#2}\IfBooleanT{#1}{,}{#3}\right]}}
\NewDocumentCommand\Aut{o}%
{%%
\ensuremath{%%
\IfValueTF{#1}%
{%%%
\operatorname{Aut}_{#1}
}%%%
{%%%
\operatorname{Aut}%
}%%%
}%%
}%%
\NewDocumentCommand\VertAut{o}%
{%%
\ensuremath{%%
\IfValueTF{#1}%
{%%%
\operatorname{Aut}^|_{#1}
}%%%
{%%%
\operatorname{Aut}%
}%%%
}%%
}%%
\NewDocumentCommand\amal{mmm}{\ensuremath{{#1}\mathbin{\times^{#2}} \! #3}}
\NewDocumentCommand\vbTM{}{\vb{{\mathfrak{g}/{\mathfrak{h}}}}}
\NewDocumentCommand\vbg{}{\bm{\mathfrak{g}}}
\NewDocumentCommand\vbk{}{\bm{\mathfrak{k}}}
\NewDocumentCommand\vbh{}{\bm{\mathfrak{h}}}
\NewDocumentCommand\framebundle{e_m}{\mathscr{F}_{#2\IfValueT{#1}{,{#1}}}}
\NewDocumentCommand\vb{m}{\ensuremath{\bm{#1}}}
\NewDocumentCommand\vc{sm}{\ensuremath{\IfBooleanTF{#1}{\overrightarrow{#2}}{\vec{#2}}}}
\NewDocumentCommand\Sym{se{^}m}%
{%
\IfBooleanTF{#1}%
{%%
\IfValueTF{#2}%
{%%%
\operatorname{Sym}^{#2}%
}%%%
{%%%
\operatorname{Sym}%
}%%%
\left(#3\right)%
}%%
{%%
\IfValueTF{#2}%
{%%%
\operatorname{Sym}^{#2}%
}%%%
{%%%
\operatorname{Sym}%
}%%%
#3%
}%%
}%
\NewDocumentCommand\SymZ{se{^}m}%
{%
\IfBooleanTF{#1}%
{%%
\IfValueTF{#2}%
{%%%
\operatorname{Sym}_0^{#2}%
}%%%
{%%%
\operatorname{Sym}_0%
}%%%
\left(#3\right)%
}%%
{%%
\IfValueTF{#2}%
{%%%
\operatorname{Sym}_0^{#2}%
}%%%
{%%%
\operatorname{Sym}_0%
}%%%
#3%
}%%
}%

\NewDocumentCommand\rightMC{o}%
{%
\IfValueTF{#1}{\mu_{#1}}{\mu}%
}%

\makeatletter
\DeclareRobustCommand{\cev}[1]{%
  \mathpalette\do@cev{#1}%
}
\newcommand{\do@cev}[2]{%
  \fix@cev{#1}{+}%
  \reflectbox{$\m@th#1\vec{\reflectbox{$\fix@cev{#1}{-}\m@th#1#2\fix@cev{#1}{+}$}}$}%
  \fix@cev{#1}{-}%
}
\newcommand{\fix@cev}[2]{%
  \ifx#1\displaystyle
    \mkern#23mu
  \else
    \ifx#1\textstyle
      \mkern#23mu
    \else
      \ifx#1\scriptstyle
        \mkern#22mu
      \else
        \mkern#22mu
      \fi
    \fi
  \fi
}

\makeatother

\NewDocumentCommand\rvc{sm}{\ensuremath{\IfBooleanTF{#1}{\overleftarrow{#2}}{\cev{#2}}}}
\NewDocumentCommand\fl{m}{\ensuremath{\mathbf{e}^{#1}}}
\def\lst{A,F,G,H,K,L,M,N,Z}
\makeatletter
\@for\i:=\lst\do{\expandafter\MakeLie \i}
\makeatother
\NewDocumentCommand\G{}{\mathscr{G}}
\NewDocumentCommand\g{}{\texttt{\textup{g}}}
\NewDocumentCommand\Greg{}{\mathscr{G}^r}
\NewDocumentCommand\Mreg{}{M^r}
\NewDocumentCommand\bbb{m}{#1^{\flat}}
\NewDocumentCommand\bboundary{}{\(\flat\)-boundary}
\NewDocumentCommand\Bun{}{E}
\NewDocumentCommand\redComplement{}{\mathfrak{m}}
\NewDocumentCommand\map{omm}{\ensuremath{\IfValueTF{#1}{{#2}\xrightarrow{#1}{#3}}{{#2}\to{#3}}}}
\NewDocumentCommand\mapto{omm}{\ensuremath{\IfValueTF{#1}{{#2}\xmapsto{#1}{#3}}{{#2}\mapsto{#3}}}}
\NewDocumentCommand\ExtAut{m}{\operatorname{Aut}'_{#1}}
\NewDocumentCommand\Conn{m}{\mathscr{A}_{#1}}
\NewDocumentCommand\cohomology{mm}{H^{#1}\left(#2\right)}
\NewDocumentCommand\nForms{mm}{\Omega^{\IfStrEq{#1}{*}{\complex}{#1}}_{#2}}
\NewDocumentCommand\Gr{smm}%
{%
\ensuremath{\operatorname{Gr}_{#2}\!\IfBooleanTF{#1}{\left({#3}\right)}{#3}}%
}%
\NewDocumentCommand\Gau{}{T}
\newcommand{\bunderline}[2][4]{\underline{#2\mkern-#1mu}\mkern#1mu }
\NewDocumentCommand\otM{}{\bunderline{M}}
\NewDocumentCommand\otx{}{\bunderline{x}}
\NewDocumentCommand\ota{}{\bunderline{a}}
\NewDocumentCommand\ote{}{\bunderline{e}}
\NewDocumentCommand\otv{}{\bunderline{v}}
\NewDocumentCommand\otp{}{\bunderline{p}}
\NewDocumentCommand\otsigma{}{\bunderline{\sigma}}
\NewDocumentCommand\otgamma{}{\bunderline{\gamma}}
\NewDocumentCommand\otG{}{\bunderline{\mathscr{G}}}
\NewDocumentCommand\eep{m}{e^{\scalebox{.5}{+}}_{#1}}
\NewDocumentCommand\eem{m}{e^{\scalebox{.7}{-}}_{#1}}
\NewDocumentCommand\horizontalPart{sm}{\IfBooleanTF{#1}{\underline{#2}}{\bunderline{#2}}}
\NewDocumentCommand\covder{om}{\IfValueTF{#1}{\nabla_{#1}#2}{\nabla#2}}
\NewDocumentCommand\holonomyGroup{mo}{\IfValueTF{#2}{\dot{#1}_{#2}}{\dot{#1}}}
\NewDocumentCommand\restrictedHolonomyGroup{mo}{%%
\IfValueTF{#2}
{\mathring #1_{#2}}
{\mathring #1}
}%%
\NewDocumentCommand\holonomyAlgebra{mo}{\IfValueTF{#2}{\dot{\Lie#1}_{#2}}{\dot{\Lie#1}}}
\NewDocumentCommand\holonomyReduction{mo}{\IfValueTF{#2}{\dot{#1}_{#2}}{\dot{#1}}}
\NewDocumentCommand\area{m}{\operatorname{area}_{#1}}
\let\marginnote\marginpar
\DeclareFontFamily{U}{stixbbit}{}
\DeclareFontShape{U}{stixbbit}{m}{it}{<-> stix-mathbbit}{}
\DeclareRobustCommand{\stixdanger}{%
  {\usefont{U}{stixbbit}{m}{it}\symbol{"F6}}%
}
\NewDocumentCommand\Danger{}%
{%
\rotatebox{45}{\stixdanger}\noindent \emph{Danger:} %
}%

\NewDocumentCommand\VF{s}%
{%
\IfBooleanTF{#1}{\bar{\mathfrak{F}}}{\mathfrak{F}}%
}%

\NewDocumentCommand\OO{D(){0}O{1}}%
{%%
\ensuremath{%%%
\mathscr{O}%
\ifnum\pdfstrcmp{#1}{0}=0{}\else{(#1)}\fi
\ifnum\pdfstrcmp{#2}{1}=0{}\else{
	\ifnum\pdfstrcmp{#2}{0}=0{}\else{^{\oplus #2}}\fi
}\fi
}%%%
}%%

\begin{document}
\newpage
In this book, we explain what Cartan geometries are, aiming at an audience of graduate students familiar with manifolds, Lie groups and differential forms.
The standard reference works on Cartan geometries are Sharpe \cite{Sharpe:1997} and \v{C}ap and Slov\`ak \cite{Cap/Slovak:2009}.
The reader might benefit from studying locally homogeneous structures \cite{Goldman:2010,Goldman:2022,Melnick:2021} and Cartan's famous papers \cite{Cartan:1910,Cartan:68,Cartan:1924,Cartan:136,Cartan:136bis,Cartan:161,Cartan:1938,Cartan:174,Cartan:1992}.

\vskip .5cm 

Thanks to Jacob William Erickson for help with the completeness arguments and for his writings on parabolic geometries. 
Thanks to Dennis The, Boris Kruglikov, Elvind Schneider and Erlend Grong for organizing the Geilo workshop on Cartan geometries and inviting me to speak there.
Thanks to Ian Andersen, Andreas \v{C}ap and Micheal Eastwood for discussions on earlier drafts. 
This publication is based upon work from COST Action CaLISTA, CA21109, supported by COST (European Cooperation in Science and Technology) \url{www.cost.eu}.

\afterpreface
\chapter{Geometries}
\epigraph[author={\'E. Cartan},
source={La méthode du repère mobile, la théorie des groupes continus et les espaces généralisés},translation={It is clear that a large number of generalized geometries are nothing more than geometrical curiosities.}]
{Il est clair qu'un grand nombre de g\'eom\'etries g\'en\'eralis\'ees ne sont jusqu\'a pr\'esent que des curiosit\'es g\'eometriques.}
\begin{marginfigure}[-3.5cm]\includegraphics{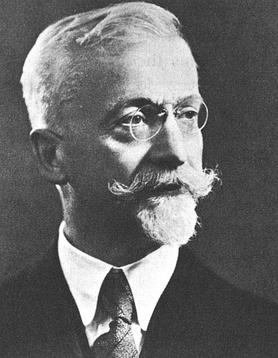}\end{marginfigure}
Riemannian geometry\SubIndex{Riemannian geometry} develops largely by analogy with Euclidean geometry.\SubIndex{Euclidean geometry}
This is natural since Euclidean space has the largest isometry group of any Riemannian manifold.
To be precise, Euclidean space of a given dimension \(n\) has isometry group of dimension \(n(n+1)/2\), as large as that of any Riemannian geometry on any \(n\)-dimensional manifold.
On the other hand, hyperbolic space\SubIndex{hyperbolic space} and the sphere\SubIndex{sphere} both have isometry groups of the same dimension as does Euclidean space.
They are the only Riemannian geometries with isometry groups of such a large dimension, except for quotienting the sphere antipodally to form real projective space,\SubIndex{projective!space!real} and scaling any of these metrics by a positive constant.
We are tempted to consider hyperbolic space, the sphere and the real projective space as ``equal'' models of Riemannian geometry.
Riemannian geometry is a ``lumpy'' cousin of these constant curvature geometries.
These models have always been the most fruitful examples in the study of Riemannian geometry.

Similarly, we study many geometric structures on manifolds by starting with some Platonic ``model'', with exceptionally large symmetry group, and using a combination of geometry and group theory to learn its secrets.
We consider any other ``model'' with equally large symmetry group as being equally valuable.
We learn about geometric structures by analogy with our models.
\[
(X,G)=(\text{model},\text{symmetry group})
\]
\[
\begin{array}{cllll}
\toprule
&\text{Riemannian}&\text{Lorentzian}&\text{Conformal}\SubIndex{conformal!geometry}&\text{K\"ahler}\SubIndex{Kaehler geometry@K\"ahler geometry}\\ 
\cmidrule(r){2-2}\cmidrule(lr){3-3}\cmidrule(lr){4-4}\cmidrule(l){5-5}
X&\R^n&\R^{n-1,1}&S^n&\C^n\\
G&\Orth{n}\ltimes\R^n&\Orth{n-1,1}\ltimes\R^{n-1,1}&\PO{n+1,1}&\Un{n}\ltimes\C^n\\[10pt]
X&S^n&\text{de Sitter}\SubIndex{de~Sitter space}&&\CP{n}\\
G&\Orth{n+1}&\Orth{n,1}&&\PU{n+1}\\[10pt]
X&\mathbb{H}^n&\text{anti-de~Sitter}\SubIndex{anti-de~Sitter space}&&\mathbb{CH}^n\\
G&\Orth{n,1}&\Orth{n-1,2}&&\PU{n,1}\\
\bottomrule
\end{array}
\]
Much as Riemannian manifolds are infinitesimally Euclidean, geometric structures which are infinitesimally analogous in some relevant respects to some model give up their secrets as we follow these analogies.
Locally, our geometric structures exhibit some curvature which makes them different from the model.
Globally, the manifolds that these geometric structures live on could have very different topology from the model.

We will begin by considering the groups themselves, then study homogeneous spaces, and then locally homogeneous geometric structures on manifolds, and finally consider Cartan geometries.
Cartan geometries are the ``lumpy'' geometric structures, perhaps not locally homogeneous, but ``infinitesimally modelled'' on homogeneous spaces, just as Riemannian geometry is ``infinitesimally modelled'' on Euclidean space.
In large part, our aim is to prove that every Cartan geometry shares many global properties of its homogeneous model.

Naturally the student will wonder why we don't expand out everything in coordinates\SubIndex{coordinates} as soon as possible, to make the various objects explicit.
It seems that in any study of differential geometry with a sophisticated local structure, we first take the bird's eye view.
For us, this requires organizing the representation theory of various groups their actions on the various invariants we run in to.
Once the algebra is organized well, we may be able to use it to identify clever choices of coordinates, in which the geometry is expressed clearly enough to start analysis.
At the moment, there are no general methods of analysis applicable to all the different types of Cartan geometries which arise in practice.

\chapter{Lie groups}
A \emph{Lie group}\define{Lie!group} \(G\) is a manifold and a group, so that the group multiplication is a smooth map.
\[
\input{maurer-cartan}
\]
A group has an identity \(1\in G\).
\[
\input{maurer-cartan-1}
\]
The operation of multiplying every element of the group, on the left, by a chosen element \(g\) is \emph{left translation}.
\[
\input{maurer-cartan-6}
\]
The \emph{Lie algebra}\define{Lie!algebra} \(\LieG\) of \(G\) is the tangent space \(T_1 G\):
\[
\input{maurer-cartan-2}
\]
\section{The Maurer--Cartan form}
Picture the tangent space at a point \(g\):
\[
\input{maurer-cartan-3}
\]
Translate \(g\) to \(1\), by left multiplication by \(g^{-1}\), transforming the entire group:
\[
\input{maurer-cartan-4}
\]
This translates \(T_g G\) to \(\LieG:=T_1 G\).
The \emph{Maurer--Cartan \(1\)-form}\define{Maurer--Cartan!form} is the operation we denote by \(\omega_G\) which left translates each vector, at any point \(g\), by \(g^{-1}\).
So \(\omega_G\) is valued in \(\LieG\).
\[
\input{maurer-cartan-5}
\]
\begin{problem}{mc.left.inv}
Prove that the Maurer--Cartan form is left invariant.
Hint: products of left translations are left translations.
\end{problem}
\begin{problem}{a}
What is left translation on the group of affine transformations of the real number line? What is the Maurer--Cartan form?
\end{problem}
\section{Example: square matrices}
Suppose that \(G=\GL{n}\) is the group of invertible real \(n\times n\) matrices.
Since \(\GL{n}\subset\R^{n\times n}\) is an open set of \(n\times n\) matrices, its tangent space at each \(n\times n\) matrix \(x\) is
\[
T_x \GL{n}\cong\R^{n\times n}.
\]
But thinking of \(\GL{n}\) as an abstract manifold, we could also treat each tangent space as a separate vector space, a separate copy of \(\R^{n\times n}\).
Left translation is left multiplication by a given matrix: left translation by an \(n\times n\) matrix \(x\) is the operation \(\LT{x}y=xy\).
The differential of any map is the linear approximation.
Differentiating that operation, which is already linear, gives the same operation:
\[
\LT{x*}y=xy.
\]

Let \(\omega:=\omega_{\GL{n}}\) be the Maurer--Cartan form of \(\GL{n}\).
So \(\omega(y)=x^{-1}y\) for any tangent vector 
\[
y\in T_x\GL{n}.
\]
So then \(\omega\) at each \(n\times n\) matrix \(x\) left translates us back from \(x\) to \(1\),
\[
y\in T_x\GL{n}\xrightarrow{\omega}\LT{x*}^{-1}y=x^{-1}y\in T_1\GL{n}.
\]

Think for a moment about one variable calculus.
In one variable calculus, we write \(x^2\) to mean the function which associates to each real number its square, \(\sin x\) to mean the function which associates to each real number its sine, and so on.
So the expression \(x\) means the identity function, a sophisticated way of thinking of a variable varying over the real number line.
Now we will do the same for matrices in place of real numbers.

We could let the letter \(g\) stand for an arbitrary \(n\times n\) invertible matrix, so a point in  \(\GL{n}\).
But, just as in one variable calculus, we can instead use the letter \(g\) to mean the identity map \(g(x)=x\) for each matrix \(x\); we can say we are thinking about \(g\) as ``every element'' of \(\GL{n}\), not any ``particular element''.
When we differentiate, since \(g(x)=x\), we find \(dg_x=I\), i.e. \(dg_x(y)=y\) for any matrix \(y\).
Left translating by \(x=g(x)\), we can use the popular formal notation
\[
\omega=g^{-1}dg,
\]
by which we mean
\[
\omega_x(y)=g(x)^{-1}dg_x(y)=x^{-1}y.
\]
In this formal notation, thinking of \(g\) as a matrix-valued function on the manifold \(\GL{n}\),
\begin{align*}
(\LT{g_0}^*g)(x)
&=
g(\LT{g_0}x),
\\
&=
g_0x,
\\
&=
g_0g(x),
\\
&=
\LT{g_0}g(x),
\end{align*}
i.e.
\[
\LT{g_0}^*g=\LT{g_0}g.
\]
Take \(d\):
\[
\LT{g_0}^*dg=d\LT{g_0}^*g=d(g_0g)=g_0\,dg=\LT{g_0}dg.
\]
Hence
\[
\LT{g_0}^*(g^{-1}dg)=g^{-1}g_0^{-1}g_0dg=g^{-1}dg.
\]
Similarly, under right action by \(g_0\),
\[
\RT{g_0}^*(g^{-1}dg)=g_0^{-1}g^{-1}dg\, g_0=\Ad_{g_0}^{-1}(g^{-1}dg).
\]
If a Lie group \(G\) admits a representation \(G\to\GL{n}\) with discrete kernel, its Maurer--Cartan form becomes just the pullback of \(g^{-1}\,dg\) to \(G\).
By Ado's theorem (\cite{Hilgert.Neeb:2012} p. 189 Theorem 7.4.1), every Lie group's identity component has a covering space which admits a representation with discrete kernel, so at least near the identity the Maurer--Cartan form of any Lie group arises as \(g^{-1}\,dg\).
Hence the Maurer--Cartan form on any Lie group is often denoted \(g^{-1}dg\) just to pretend it is a group of matrices, and there is no danger in this notation.
\begin{problem}{b}
If \(G\) is the group of \(3\times 3\) matrices of the form
\[
g=
\begin{pmatrix}
1&x&z\\
0&1&y\\
0&0&1
\end{pmatrix}
\]
check that
\[
\omega_G
=
\begin{pmatrix}
0&dx&dz-x\,dy\\
0&0&dy\\
0&0&0
\end{pmatrix}.
\]
\end{problem}
\section{The Lie algebra}
If we pick any element \(A\in\LieG\), we can define a unique vector field \(A_G\) by the rule \(A_G\hook\omega_G=A\), since \(\omega_G\) is a linear isomorphism on each tangent space.
Since \(\omega_G\) is just left translation, \(A_G\) is \emph{left invariant},\define{left invariant} i.e. invariant under left translation.
\begin{problem}{c}
Conversely every left invariant vector field is \(A_G\) for a unique \(A\in\LieG\).
\end{problem}
\begin{problem}{Lie.bracket.invariant}
Prove that the Lie bracket operation on vector fields on a manifold is equivariant under diffeomorphism.
\end{problem}
So if two vector fields are invariant under some diffeomorphism, so is their bracket.
So the brackets of left invariant vector fields are left invariant vector fields.
(Similarly for right invariant vector fields.)
We define the Lie bracket on \(\LieG\) by setting \(C:=\lb{A}{B}\) in \(\LieG\) precisely when \(C_G=\lb*{A_G}{B_G}\) as vector fields.
\begin{lemma}
The Lie bracket \(\lb{A}{B}\) makes \(\LieG\) into a Lie algebra, i.e. it is linear in \(A\) and in \(B\), antisymmetric: \(\lb{B}{A}=-\lb{A}{B}\), and satisfies the Jacobi identity:
\[
\lb{\lb{A}{B}}{C}
\]
sums to zero when we cyclically permute \(A,B,C\).
\end{lemma}
\begin{proof}
These are true for vector fields, hence for the left invariant vector fields.
\end{proof}
\begin{theorem}%
[Lie's Third Theorem \cite{Hilgert.Neeb:2012} p. 334 9.4.11]%
\label{theorem:Lie.3.a}%
\define{theorem!Lie's third}%
\define{Lie's third theorem}%
Every finite dimensional Lie algebra is the Lie algebra of a a connected and simply connected Lie group, unique up to isomorphism.
\end{theorem}
\section{The significance of the Maurer--Cartan form}
\begin{theorem}\label{theorem:open.subset.G}
A diffeomorphism of connected open subsets of \(G\) is left translation by an element of \(G\) just when it preserves the Maurer--Cartan form.
\end{theorem}
\begin{proof}
Compose with a left translation to ensure that our map preserves some point, which we can assume is \(1\in G\).
Our map now commutes with left invariant vector fields, because these are the vector fields on which the Maurer--Cartan form is constant.
So our map is the identity map along the flow lines of the left invariant vector fields.
But left invariant vector fields point in all directions.
\end{proof}
\section{Right translation}
Define \(\Ad_g\) to be the product of left translation by \(g\), and right translation by \(g^{-1}\), so \(\LieG\xrightarrow{\Ad_g}\LieG\).
\begin{problem}{d}
The Maurer--Cartan form transforms under right translation in the \emph{adjoint representation}:\define{adjoint representation} denoting right translation by \(g\in G\) as \(\RT{g}\), \(\RT{g}^*\omega_G=\Ad_g^{-1}\omega_G\).
\end{problem}
\begin{problem}{e}
Denoting left translation by \(g\in G\) as \(\LT{g}\), prove that \(\LT{g*}^{-1}=\LT{g^{-1}*}\).
\end{problem}
\begin{lemma}\label{lemma:invert.Maurer.Cartan}
For any Lie group \(G\), denote the inverse map by
\[
g\in G\xmapsto{\iota}g^{-1}\in G.
\]
Then 
\[
\iota^*\omega_G=-\Ad_g\omega_G.
\]
\end{lemma}
\begin{proof}
Replacing \(g\) by \(g^{-1}\), left translation becomes right translation, and so the left translation to the identity \(\omega_G\) becomes right translation to the identity, which is left translation followed by adjoint action.
But the differential of \(g\mapsto g^{-1}\) at the identity is \(-1\),  so \(\omega_G\) transforms to \(-\Ad_g\omega_G\).
\end{proof}
\begin{problem}{f}
Let \(\iota(g):=g^{-1}\).
Prove that
\[
\iota'(g)=-\RT{g*}^{-1}\LT{g*}^{-1}\colon T_g G\to T_{g^{-1}}G.
\]
\end{problem}
\begin{problem}{right.invariant}
Denote by \(\vc{A}\) the left invariant vector field on a Lie group \(G\) associated to an element \(A\in\LieG\), and by \(\rvc{A}\) the right invariant vector field.
Prove that
\begin{align*}
\lb{\vc{A}}{\vc{B}}&=\vc*{\lb{A}{B}},\\
\lb{\rvc{A}}{\rvc{B}}&=-\rvc*{\lb{A}{B}}.
\end{align*}
\end{problem}
\begin{answer}{right.invariant}
The first equation is our definition of Lie bracket on \(\LieG\).
The second equation follows from the first by applying \(g\mapsto g^{-1}\).
\end{answer}
\section{The structure equations}
\begin{lemma}
The left invariant Maurer--Cartan form \(\omega:=\omega_G\) on any Lie group \(G\) satisfies the \emph{Maurer--Cartan structure equations}\define{Maurer--Cartan!structure equations}\define{structure equations!Maurer--Cartan}
\[
d\omega(A_G,B_G)=-\lb{A}{B}.
\]
\end{lemma}
\begin{proof}
We want to use the Cartan formula
\[
d\omega(v,w)=\LieDer_v(w\hook\omega)-\LieDer_w(v\hook\omega)-\omega(\lb{v}{w}),
\]
which holds for any vector fields \(v,w\) on \(G\).
We will choose \(v,w\) to be left invariant, i.e. take any \(A,B\in\LieG\) and let \(v:=A_G\), \(w:=B_G\).
Then \(v\hook\omega=A\) is constant so
\begin{align*}
d\omega(v,w)
&=
\LieDer_v(w\hook\omega)-\LieDer_w(v\hook\omega)-\omega(\lb{v}{w}),
\\
&=
\LieDer_v B-\LieDer_w A-\omega(\lb{A}{B}_G),
\\
&=0-0-\lb{A}{B},\\
&=-\lb{A}{B}.
\end{align*}
\end{proof}
\section{Algebra valued forms}
If we have any \(1\)-forms \(\alpha,\beta\) valued in some algebra \(A\), it is traditional to define \(\alpha\wedge\beta\) to mean the \(2\)-form valued in \(A\):
\[
(\alpha\wedge\beta)(v,w):=\alpha(v)\cdot\beta(w)-\alpha(w)\cdot\beta(v),
\]
where \(\cdot\) is multiplication in \(A\).
If \(A\) is a Lie algebra, say write it as \(\LieG\), we write \(\alpha\wedge\beta\) as \(\lb{\alpha}{\beta}\), so
\[
\lb{\alpha}{\beta}(v,w)=\lb{\alpha(v)}{\beta(w)}-\lb{\alpha(w)}{\beta(v)}.
\]
In particular, for a Maurer--Cartan form \(\omega=\omega_G\),
\begin{align*}
\lb{\omega}{\omega}(v,w)
&=\lb{\omega(v)}{\omega(w)}-\lb{\omega(w)}{\omega(v)},
\\
&=2\lb{\omega(v)}{\omega(w)},
\end{align*}
which yields many strange factors of \(2\) or \(1/2\) from here on.
In particular, if \(v:=A_G\) and \(w:=B_G\) for some \(A,B\in\LieG\) then 
\[
\lb{\omega}{\omega}(v,w)=2\lb{A}{B}.
\]
We strive to avoid more complicated expressions than this.
\section{The structure equations again}
\begin{lemma}\label{lemma:mc}
The left invariant Maurer--Cartan form \(\omega:=\omega_G\) on any Lie group \(G\) satisfies the \emph{Maurer--Cartan structure equations}\define{Maurer--Cartan!structure equations}\define{structure equations!Maurer--Cartan}
\[
d\omega+\frac{1}{2}\lb{\omega}{\omega}=0.
\]
\end{lemma}
\begin{proof}
\[
\lb{A}{B}=\frac{1}{2}\lb{\omega}{\omega}(A_G,B_G)
\]
\end{proof}
In a basis of \(\LieG\), we can write \(\omega\) as having components \(\omega^i\) and write these structure equations as
\[
d\omega^i+c^i_{jk}\omega^j\wedge\omega^k=0,
\] 
for some constants \(c^i_{jk}\), the \emph{structure constants};\define{structure constants}
we also call these equations the \emph{Maurer--Cartan structure equations}.\define{Maurer--Cartan!structure equations}\define{structure equations!Maurer--Cartan}
In any basis \(e_1,\dots,e_n\) of a Lie algebra \(\LieG\), we can write the Lie brackets uniquely as
\[
\lb{e_i}{e_j}=c^k_{ij}e_k.
\]
\begin{problem}{struc.con}
Prove that these \(c^k_{ij}\) are the structure constants.
\end{problem}
\begin{example}
Take the group\SubIndex{affine group} \(G\) of all affine transformations of the plane which preserve the Euclidean metric up to some constant factor.
In complex notation, each looks like \(z\mapsto az+b\) for some complex numbers \(a\ne 0,b\).
We can write them as matrices
\[
g=
\begin{pmatrix}
a&b\\
0&1
\end{pmatrix}.
\]
Hence the Maurer--Cartan form is
\[
\omega=g^{-1}dg=
\begin{pmatrix}
a&b\\
0&1
\end{pmatrix}^{-1}
\begin{pmatrix}
da&db\\
0&0
\end{pmatrix}
=
\begin{pmatrix}
\frac{da}{a}&\frac{db}{a}\\
0&0
\end{pmatrix}.
\]
Writing this as
\[
\omega=
\begin{pmatrix}
\gamma&\sigma\\
0&0
\end{pmatrix}
\]
compute
\begin{align*}
d\omega
&=
-\frac{1}{2}\lb{\omega}{\omega},\\
&=-\omega\wedge\omega,\\
&=
-
\begin{pmatrix}
\gamma&\sigma\\
0&0
\end{pmatrix}
\wedge
\begin{pmatrix}
\gamma&\sigma\\
0&0
\end{pmatrix},
\\
&=
-
\begin{pmatrix}
0&\gamma\wedge\sigma\\
0&0
\end{pmatrix}.
\end{align*}
So
\begin{align*}
d\sigma&=-\gamma\wedge\sigma,\\
d\gamma&=0.
\end{align*}
\end{example}
\begin{problem}{translate.mult}
Suppose that \(G\) is a Lie group, say with multiplication \(x,y\mapsto xy\).
For any fixed \(g\in G\), let
\[
x*_gy:=xg^{-1}y.
\]
Prove that this operation, the \emph{left translated multiplication}, is the multiplication of a Lie group, with identity \(g\), and that left translation \(x\mapsto g^{-1} x\) identifies the Lie group structure of \(*_g\) with the usual Lie group structure, i.e. that of \(*_1\).
Find the inverse operation and the left translation operation of \(*_g\) in terms of \(*_1\).
Prove that, for suitable isomorphisms of their Lie algebras, all \(*_g\) have the same Maurer--Cartan forms.
\end{problem}
\begin{answer}{translate.mult}
If \(\varphi(x)=g^{-1} x\), check that \(\varphi(x *_g y)=\varphi(x)\varphi(y)\), a Lie group isomorphism, so clearly this is a Lie group with identity \(g\).
The inverse operation is \(\iota^g(x):=gx^{-1}g\).
The left translation is \(\ell^g_xy:=x *_g y=xg^{-1}y=\ell_{xg^{-1}}y\) i.e. \(\ell^g_x=\ell_{xg^{-1}}\).
The Maurer--Cartan form, on a tangent vector \(v\in T_x G\), is
\begin{align*}
\omega^g(v)&:=(\ell^g_{\iota^g(x)})'(x)v,
\\
&=(\ell_{gx^{-1}})'(x)v,
\\
&=\ell_g'(1)(\ell_{x^{-1}})(x)v,
\\
&=
\ell_g'(1)\omega(v).
\end{align*}
So the linear isomorphism \(\ell_g'(1)\) taking \(\LieG\) to the Lie algebra \(\LieG^g:=T_g G\) of \(*_g\) identifies the Maurer--Cartan forms.
\end{answer}
\section{Recovering a Lie group}
When is a differential form the Maurer--Cartan form of a Lie group, and in how many different ways?

A \emph{Maurer--Cartan system} on a manifold \(M\) is a collection of linear independent \(1\)-forms \(\omega^i\) spanning its cotangent spaces, so that every \(\omega^i\) has exterior derivative \(d\omega^i\) a constant coefficient linear combination of wedge products \(\omega^j\wedge\omega^k\).
The \emph{constant vector fields} of a Maurer--Cartan system are the vector fields on which these \(\omega^i\) are constant.

A \emph{Maurer--Cartan form} on a manifold \(M\) is a \(1\)-form \(\omega\) valued in Lie algebra \(\LieG\) so that \(\omega\) is a linear isomorphism of each tangent space of \(M\) to \(\LieG\) and so that
\[
0=d\omega+\frac{1}{2}\lb{\omega}{\omega}.
\]

\begin{lemma}
By picking a basis of the Lie algebra, every Maurer--Cartan form gives rise to a Maurer--Cartan system.
Every Maurer--Cartan system arises in this way from a Maurer--Cartan form valued in a Lie algebra with a basis, uniquely up to basis preserving Lie algebra isomorphism.
\end{lemma}
\begin{proof}
Given a Maurer--Cartan form \(\omega\), valued in a Lie algebra \(\LieG\), take a basis \(e_1,\dots,e_n\) of \(\LieG\), and write \(\omega=\omega^ie_i\).
By lemma~\vref{lemma:mc}, these \(\omega^i\) are a Maurer--Cartan system.

Given a Maurer--Cartan system \(\set{\omega^i}\), suppose
\[
d\omega^i=-c^i_{jk}\omega^j\wedge\omega^k,
\]
for constants \(c^i_{jk}\).
We can assume that \(c^i_{jk}=-c^i_{kj}\).
Let \(n:=\dim M\).
Let \(\LieG:=\R^n\) with standard basis \(e_1,\dots,e_n\) and with Lie bracket
\[
\lb{e_i}{e_j}:=c^k_{ij}e_k.
\]
Take exterior derivative of the equations
\[
d\omega^i=-c^i_{jk}\omega^j\wedge\omega^k,
\]
to find that these \(c^i_{jk}\) satisfy the Jacobi identity.
So \(\omega:=\omega^ie_i\) is a Maurer--Cartan form.
\end{proof}

Since differential forms pullback by smooth maps, any Maurer--Cartan system pulls back by any local diffeomorphism to a Maurer--Cartan system.
\begin{theorem}\label{theorem:recover.Lie.group}
Take a connected manifold \(M\) with a Maurer--Cartan system.
There is a Lie group \(G\) with the constant coefficients of the Maurer--Cartan system as its structure constants.
The pullback of the Maurer--Cartan system to the universal covering space \(\tilde{M}\to M\)  is also the pullback from a local diffeomorphism \(\tilde{M}\xrightarrow{\delta}G\):
\[
\begin{tikzcd}
\tilde{M}\arrow[r,"\delta"]\arrow[d]&G\\
M
\end{tikzcd}
\]
the \emph{developing map}.\define{developing map}
The developing map is equivariant for a unique group morphism \(\pi_1(M)\xrightarrow{h} G\), the \emph{holonomy morphism}.\define{holonomy!morphism}
The pair \((\delta,h)\) of developing map and holonomy morphism are unique up to replacing by \((g\delta,\Ad_g h)\) for any \(g\in G\).
 \end{theorem}
Theorem~\vref{thm:flat} has this theorem as a trivial consequence, and the proof is almost the same, but we give a proof here as a preview of the proof of that theorem.
\begin{proof}
By Lie's third theorem (theorem~\vref{theorem:Lie.3.a}), there is a Lie group \(G\) with the structure constants of the Maurer--Cartan system on \(M\) as its structure constants.
On \(M\times G\), the Pfaffian system \(\set{\omega^i_M-\omega^i_G}\) arises from a foliation, by the Frobenius theorem (theorem~\vref{theorem:Frobenius}), and every leaf has the same dimension as the kernel of \(\omega_M-\omega_G\), i.e. the dimension of \(M\).
Denote by \(\tilde{M}\to M\) the universal covering space of \(M\).
Let \(\tilde{M}\) be a leaf, i.e. a maximal integral manifold of the foliation, through some point \((m_0,1)\in M\times G\). 
The inclusion \(\tilde{M}\subseteq M\times G\) composes with the obvious projection maps \(M\times G\to M, G\) to give maps \(\tilde{M}\xrightarrow{p}M\) and \(\tilde{M}\xrightarrow{\delta}G\).
By theorem~\vref{theorem:orbit.maps}, these maps are smooth submersions, hence local diffeomorphisms.

Take any constant vector field \(A_M\) on \(M\) associated to a vector \(A\in\LieG\).
The vector field \(A_{\tilde{M}}:=(A_M,A_G)\) on \(M\times G\) is tangent to \(\tilde{M}\), since \(\omega_M=\omega_G\) on \(\tilde{M}\).
If \(A_M\) has a flow line \(m(t)\) on \(M\) with \(m(0)=m_0\), then \(A_{\tilde{M}}\) has associated flow line through any point \((m_0,g_0)\) given by
\[
\dot{g}(t)\hook\omega_G=A,
\]
so \(g(t)=g_0e^{tA}\).
So as long as we can solve for the flow of \(A_M\), we can solve for the flow of \(A_{\tilde{M}}\).
In the terminology of theorem~\vref{theorem:orbit.maps}, the vector fields \(A_M\) and \(A_{\tilde{M}}\) are \emph{equicomplete} for \(\tilde{M}\to M\), i.e. have flows defined on the same time intervals.

These vector fields \(A_{\tilde{M}}\), \(A_M\) span the tangent spaces of \(\tilde{M}\) and \(M\).
By theorem~\vref{theorem:orbit.maps}, \(\tilde{M}\xrightarrow{p}M\) is a fiber bundle map.
But it is also a local diffeomorphism.
So it is a covering map with \(\delta^*\omega_G=\omega_{\tilde{M}}:=p^*\omega_M\).

So far the construction is unique up to the choice of leaf.
You pick the leaf through some point \((m_0,1)\), and I prefer \((m_1,1)\). 
But we have seen that every leaf maps to \(M\) by a covering map.
So some point of my leaf maps to \(m_0\), say the point \((m_0,g)\).
After left translation of my leaf by \(g^{-1}\), my leaf passes through \((m_0,1)\).
But there is a unique leaf through each point, so my leaf is now your leaf.
In other words, all leaves are permuted by left translations in \(G\).

We can replace \(\tilde{M}\), if we like, by the universal covering space with base point \(m_0\in M\), and still have covering maps to \(M,G\), pulling back the Maurer--Cartan systems to agree.

The deck transformations of \(\tilde{M}\) preserve
\[
\omega^i_{\tilde{M}}:=p^*\omega^i_M, 
\]
since deck transformation act trivially on \(M\).
But the map \(\tilde{M}\xrightarrow{\delta}G\) might not be invariant under the deck transformations.
Locally, each deck transformation \(\gamma\in\fundamentalGroup{M}\) acts by an automorphism of the Maurer-Cartan system, replacing the initial point \((m_0,1)\) by some other point \((m_0,g)\).
As we saw, it thus acts as a left translation by \(g\).
Let \(h(\gamma):=g\).

Take any developing map \(\tilde{M}\xrightarrow{\delta}G\).
Define a map \(\tilde{M}\to M\times G\) by \((p,\delta)\).
Since \(p\) is a covering map, \((p,\delta)\) is a covering map to its image, which is a leaf.
So every developing map arises as above.
\end{proof}
A Maurer--Cartan system is \emph{complete} if the constant vector fields are complete, i.e. their flows are defined for all time.
\begin{theorem}
Take a connected manifold \(M\) with a Maurer--Cartan system.
The following are equivalent:
\begin{itemize}
\item
The Maurer--Cartan system is complete.
\item
Every developing map \(\tilde{M}\xrightarrow{\delta}G\) is a covering map to a path component of \(G\).
\item
The universal covering space manifold \(\tilde{M}\) is a Lie group, for some Lie group structure, so that the Maurer--Cartan system arises from its Maurer--Cartan form. 
\end{itemize}
Suppose that these occur.
This Lie group structure on \(\tilde{M}\) is unique up to left translating its multiplication: it is the Lie group structure of the universal covering group of the identity component of \(G\).
\end{theorem}
\begin{proof}
The Maurer--Cartan system on \(M\) is complete just when its pullback to \(\tilde{M}\) is complete.
So we can assume that \(M=\tilde{M}\) is simply connected.
Since \(M\) is connected, the image of the developing map lies in a single component of \(G\); after left translation, we can assume that it is the identity component.
By theorem~\vref{theorem:orbit.maps}, the developing map \(M\xrightarrow{\delta}G\) is a fiber bundle map to a path component of \(G\).
Being a local diffeomorphism, it is a covering map.
Replace \(G\) by the universal covering space of its identity component.
Since \(G\) is connected and simply connected, it is a diffeomorphism.
\end{proof}

\begin{corollary}
Take a connected manifold \(M\) with a complete Maurer--Cartan system.
Let \(\pi:=\fundamentalGroup{M}\).
Let \(G\) be the connected and simply connected Lie group with the same structure constants.
Let \(\pi\xrightarrow{h}G\) be a holonomy morphism associated to a developing map as above.
Then \(\pi\) has image a discrete subgroup \(\Gamma\subseteq G\) and \(M=\Gamma\backslash G\) is the quotient by left translations.
This description of \(M\) is unique up to conjugation in \(G\).
Conversely, every quotient \(M:=\Gamma\backslash G\) of a connected and simply connected Lie group is a connected manifold, and the Maurer--Cartan form of \(G\) descends to a complete Maurer--Cartan system on \(M\).
\end{corollary}
\begin{corollary}
On a connected and simply connected manifold, a Maurer--Cartan system arises from a Lie group structure if and only if it is complete.
\end{corollary}
\begin{corollary}
On a compact, connected and simply connected manifold, every Maurer--Cartan system arises from a Lie group structure, unique up to left translating the multiplication.
\end{corollary}

\section{Point set topology detour}
Most authors require second countability\SubIndex{second countability} as part of the definition of a manifold (i.e. existence of a countable basis of open sets), but not all \cite{Helgason:1978}.
Even without the hypothesis that manifolds be paracompact, Lie groups are paracompact \cite{Helgason:1978} p. 88.
Hence a Lie group has countably many components just when it is second countable.

\chapter{Physical analogies and Maurer--Cartan forms}
\begin{center}
\input{skate}
\end{center}
\section{Example: skating}
Skate on ice.
Your skate sits at a point \(z\) of the plane \(X=\C\), pointing in the direction of an angle \(\theta\).
You can only slide forward on a skate: \(\dot{z}=ue^{i\theta}\) where \(u\) is your speed.
You control your speed by pushing backward with the other skate (on the other foot).
You lean to one side or to other to control how your angle changes: \(\dot{\theta}=\kappa\), as you move forward, so that \(\kappa\) is the curvature of the path you follow along the ice.
On the manifold \(G=\C \times S^1\) with coordinates \((z,\theta)\) we have vector fields
\begin{align*}
X(z,\theta) &= e^{i\theta}\partial_z, \\
Y(x,y,\theta) &= e^{i\theta}\partial_z+\partial_{\theta}.
\end{align*}
In other words, the differential equations
\[
\frac{dz}{dt}=e^{i\theta}, \frac{d\theta}{dt}=0\text{ or }1
\]
are the equations of the flow lines of \(X\) and \(Y\).
The vector field \(X\) represents skating without leaning over, while \(Y\) has you leaning to your left.
There is some choice here about how much curvature you generate in \(Y\) when you lean over; we set the coefficient of \(\partial_{\theta}\) to \(1\) for simplicity of notation.
Skate by following the flow of some linear combination of \(X,Y\).
\begin{problem}{i}
Compute the flows of \(X\) and \(Y\).
Draw pictures of skates in the plane to explain \(X\) and \(Y\).
\end{problem}
The bracket is
\[
Z:=\lb{X}{Y} = ie^{i\theta} \partial_z,
\]
a third, linearly independent, vector field.
Geometrically, \(Z=\lb{X}{Y}\) fixes the angle of the skate, and slides the skate along the ice, in the direction perpendicular to the blade.
This is not physically possible motion for a skate.
Switching rapidly between flowing along \(X\) and \(Y\), we can approximate the flow of \(Z=\lb{X}{Y}\) by moving back and forth along the flows of \(X\) and \(Y\) small amounts.
Indeed, if we draw out the motion \(e^{-\sqrt{t}Y}e^{-\sqrt{t}X}e^{\sqrt{t}Y}e^{\sqrt{t}X}\) explicitly, we see the skate ``parallel parking''.
We can thereby move between any points of the plane \(X\): we can skate from anywhere to anywhere, with any starting and ending angles.

We can think of \(G=\C\rtimes S^1\) as the group of orientation preserving rigid motions of the plane \(X=\C\).
Write elements of \(G\) as matrices
\[
\begin{pmatrix}
e^{i\theta}&z\\
0&1
\end{pmatrix}.
\]
The left invariant Maurer--Cartan \(1\)-form is 
\begin{align*}
\omega_G
&=g^{-1}\,dg,
\\
&=
\begin{pmatrix}
e^{i\theta}&z\\
0&1
\end{pmatrix}^{-1}
d
\begin{pmatrix}
e^{i\theta}&z\\
0&1
\end{pmatrix},
\\
&=
\begin{pmatrix}
id\theta&e^{-i\theta}dz\\
0&0
\end{pmatrix}.
\end{align*}
The left invariant vector fields of \(G\) are the vector fields so that \(\omega_G\) is constant, i.e. \(\dot\theta\), \(e^{-i\theta}\dot{z}\) are constant, so spanned by our \(X,Y,Z\) as above.
In terms of these, we see how to direct the skater to move around.

The right invariant vector fields of \(G\) are the vector fields so that \(dg\,g^{-1}\) is constant, i.e. \(\dot\theta\), \(\dot{z}-i\dot\theta z\) constant, so the vector fields are spanned by: \(\dot\theta=0\), \(\dot{z}\) constant, i.e. translations, and \(\dot\theta=\kappa_0\) constant, \(\dot{z}=i\kappa_0 z\), rotation of the plane around the origin.
The flows of the right invariant vector fields act precisely as the isometries of the plane.
The right invariant vector fields of \(G\) have flows given by the left translations of \(G\), so they also act on the space \(G\) of configurations of the skater, preserving \(X,Y,Z\), i.e. preserving the geometry of the skater.
\section{Example: a spaceship}
Recall that an \emph{affine space}\define{affine space} \(X\) is a set acted on transitively by a vector space \(V\), its \emph{translations}, with some (hence every) point in \(X\) having stabilizer \(\set{0}\subseteq V\).
Picking any point of \(X\) suffices to identify \(V\) with \(X\), but we resist the temptation to pick a point of \(X\).
Every tangent space of \(X\) is thus identified with \(V\).
Every \(V\)-invariant Riemannian metric on \(X\) is identified with an inner product on \(V\); pick one.
Let \(G\) the group of rigid motions of \(X\), i.e. preserving the Riemannian metric.
Since the Riemannian metric is translation invariant \(V\subseteq G\).
The action of \(G\) preserves the translations, i.e. they form a closed normal subgroup, with quotient group \(H:=G/V\) the group of orthogonal transformations of \(V\).
So \(G\) is an extension
\[
1\to V\to G\to H\to 1,
\]
and we write the morphism \(G\to H\) as \(g\mapsto\bar{g}\in H\).
If we pick a reference point of \(X\), we can identify \(H\) with the subgroup of \(G\) fixing that point, splitting the exact sequence.

Let \(X\) be a \(3\)-dimensional affine space with a translation invariant Riemannian metric.
\emph{Warning}: In this section, \(e_1, e_2, e_3\) denotes an \emph{arbitrary} orthonormal basis  of \(X\), \emph{not} necessarily the standard basis of \(\R^3\). 
We write \(e\) to mean an orthonormal basis
\[
e=(e_1,e_2,e_3).
\]
\subsection{The orthonormal frame bundle}
An \emph{orthonormal frame}\define{orthonormal!frame} is a choice of a point \(x \in X\) and an orthonormal basis \(e\) of \(T_x X\cong V\).
\begin{center}
\begin{tikzpicture}
\newcommand{\rrrr}{.1}
\newcommand{\RRRR}{.75}
\newcommand{\arw}[3]
{
    \draw[-latex] 
    ({\rrrr*cos(#2)},{\rrrr*sin(#2)}) 
    -- 
    ({\RRRR*cos(#2)},{\RRRR*sin(#2)}) 
    node[black,#3]{\(e_{#1}\)};
}
\fill (0,0) circle (1pt) node[above right]{\(x\)};
\arw{1}{110}{above}
\arw{2}{230}{below left}
\arw{3}{-10}{below right}
%\draw[-latex] (0,.07071) -- (0,.7071) node[black,above]{\(e_1\)};
%\draw[-latex] (-.1,-.1) -- (-.5,-.5) node[black,below left]{\(e_2\)};
%\draw[-latex] ({\rrrr*cos(-30)},{\rrrr*sin(-30)}) -- ({\RRRR*cos(-30)},{\RRRR*sin(-30)}) node[black,below right]{\(e_3\)};
\end{tikzpicture}
\end{center}
The group \(G\) acts on orthonormal frames by
\[
g(x,e):=(gx,\bar{g}e),
\]
where
\[
\bar{g}e=\bar{g}(e_1,e_2,e_3):=(\bar{g}e_1,\bar{g}e_2,\bar{g}e_3).
\]
\begin{problem}{j}
Let \(\G\) be the set of all orthonormal frames, called the \emph{orthonormal frame bundle} of \(X\).
Prove that this is a group action of \(G\) on \(\G\).
\end{problem}
\begin{problem}{k}
Any two orthonormal frames are brought one to the other by a unique rigid motion \(g\in G\), i.e. \(G\) acts freely and transitively on \(\G\) with trivial stabilizer.
\end{problem}
\begin{problem}{smooth.action}
Prove that the action of \(G\) on \(\G\) is smooth for a unique smooth structure on \(\G\).
\end{problem}
The map
\[
(x,e)\in\G\xmapsto{\pi}x\in X
\]
is \(G\)-equivariant, i.e. commutes with the \(G\)-action.
Picture  the \(G\)-action as picking up all of \(X\) and moving its points and frames by rigid motion.
\begin{problem}{l}
The group \(\Orth{3}\) acts on \(\G\) on the right by
\[
(x,e_1,e_2,e_3)h:=(x,h^j_1e_j,h^j_2e_j,h^j_3e_j),
\]
for \(h\in\Orth{3}\), commuting with the \(G\)-action.
\end{problem}
\begin{center}
\begin{tikzpicture}
\newcommand{\rrrr}{.1}
\newcommand{\RRRR}{.75}
\newcommand{\clr}{black}
\newcommand{\arw}[3]
{
    \draw[-latex,\clr] 
    ({\rrrr*cos(#2)},{\rrrr*sin(#2)}) 
    -- 
    ({\RRRR*cos(#2)},{\RRRR*sin(#2)});
}
\fill (0,0) circle (1pt);
\arw{1}{110}{above}
\arw{2}{230}{below left}
\arw{3}{-10}{below right}
\renewcommand{\clr}{gray!75}
\arw{1}{140}{above left}
\arw{2}{260}{below}
\arw{3}{20}{right}
\end{tikzpicture}
\end{center}
The right \(\Orth{3}\)-action fixes every point of \(X\), but spins every orthonormal frame around, while the orthonormal frame sits at its unmoving point.
\begin{problem}{m}
Prove that the map \(\G\xrightarrow{\pi}X\) is a principal right \(\Orth{3}\)-bundle (called the \emph{orthonormal frame bundle})\define{orthonormal!frame bundle}\define{bundle!orthonormal frame}.
\end{problem}
It is always convenient to picture any principal bundle as if its base manifold is a Euclidean plane or a \(3\)-dimensional Euclidean space, and each point of its total space is an orthonormal frame, as Cartan did \cite{Cartan:1926}.
\subsection{Flying the spaceship}
You are the pilot of a spaceship.
I tell you, by radio, how to fly, in terms of the basis you carry with you as you fly \cite{Erickson}.
\[
\begin{array}{p{5cm}}
\input{pilots-frame}\\
\tiny
Creative Commons Attribution-Share Alike 4.0 International license.
Attribution: Acdx, cmglee
\end{array}
\]
You are at a point \(x\in X\) in a \(3\)-dimensional affine space \(X\), with a translation invariant Riemannian metric.
You are holding an orthonormal frame \(e\).
I tell you to move forward, as you measure ``forward'' in the abstract vector space \(X\), but using your orthonormal frame \(e\).
You move ``forward'', i.e. in your \(e_1\) direction: \(\dot{x}=e_1\), \(\dot{e}_1=0\).
In other words, you don't rotate your frame, but slide your base point \(x\) along the direction of the first leg of the frame.
\begin{center}
\begin{tikzpicture}
\node at (2.5,0) {\(\dots\)};
\foreach \i in {30,40,...,70}{
  \fill [fill={gray!\i!white}] ({.1*\i},0) circle (1pt);% node[above];%{\(x\)};
  \draw[-latex,draw={gray!\i!white},fill={gray!\i!white}] ({.1+.1*\i},0) -- ({.5+.1*\i},0);% node[below left];%{\(e_1\)};
  \draw[-latex,draw={gray!\i!white},fill={gray!\i!white}] ({.1*\i},.1) -- ({.1*\i},.5);% node[below right];%{\(e_2\)};
  \draw[-latex,draw={gray!\i!white},fill={gray!\i!white}] ({.07+.1*\i},-.07) -- ({.07+.1*\i+.3},{-.2});% node[below right];%{\(e_2\)};
  }
\node at (8,0) {\(\dots\)};
\end{tikzpicture}
\end{center}
\par\noindent{}%
I can describe motions in other directions, say by asking you to move \(\dot{x}=a_i e_i\), \(\dot{e}_i=0\), for constants \(a_1,a_2,a_3\), and perhaps you can do that using some sort of left-right or up-down thrusters.
In an airplane, you can only move forward.

To instruct you to roll, pitch or yaw, or to move in some other direction than straight ahead, I need a richer language of instructions.
\[
\begin{array}{p{5cm}}
\includegraphics[width=4cm]{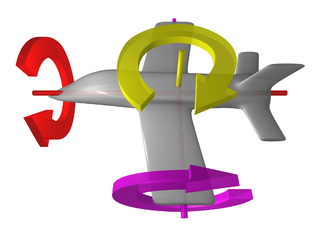}\\
\tiny
Creative Commons Attribution-Share Alike 3.0 Unported license. Attribution: ZeroOne
\end{array}
\]

If I want you to turn to your left, I want \(\dot{x}=0\), \(\dot{e}_1=e_2\).
To keep your basis orthonormal, differentiating the equation \(e_1\cdot e_2=0\), we have to agree that \(\dot{e}_2=-e_1\).
We want \(\dot{e}_3=0\), i.e. you don't move up or down.

Adding up such vector fields, for any \(v\in\R^3\) and \(B\in\LieSO{3}\), we let 
\[
A=
\begin{pmatrix}
B&v\\
0&0
\end{pmatrix}.
\]
Define a vector field \(A_{\G}(x,e)\) by
\[
A_{\G}(x,e)=(\dot{x},\dot{e}),
\]
where \(\dot{x}=v_1e_1+v_2e_2+v_3e_3\), and \(\dot{e}_i=B_{ji}e_j\).
For \(A=B\in\LieH\), the flow of \(A_{\G}\) rotates every frame, fixing the point of \(X\) where it sits.
For \(A=v\in\R^3\), the flow of \(A_{\G}\) moves the point, carrying the frame in parallel.
\begin{problem}{inv.under.rigid}
Prove that every \(A_{\G}\) is \(G\)-invariant.
\end{problem}
\begin{answer}{inv.under.rigid}
Under the action of any \(g\in G\) on \(\G\), for \(A=w\in\R^3\), and
\begin{align*}
(g_*A_{\G})(x,e)
&=
\left.\frac{d}{dt}\right|_{t=0}ge^{tA_{\G}}(x,e),
\\
&=
\left.\frac{d}{dt}\right|_{t=0}g(x+tw_ie_i),e),
\\
&=
\left.\frac{d}{dt}\right|_{t=0}(gx+tw_i\bar{g}e_i,\bar{g}e),
\\
&=
\left.\frac{d}{dt}\right|_{t=0}e^{tA_{\G}}(gx,\bar{g}e),
\\
&=
\left.\frac{d}{dt}\right|_{t=0}e^{tA_{\G}}g(x,e),
\\
&=
A_{\G}(g(x,e)).
\end{align*}
Similarly for \(A=B\in\LieSO{3}\), if \(U(t):=\exp(tB)\),
\begin{align*}
(g_*A_{\G})(x,e)
&=
\left.\frac{d}{dt}\right|_{t=0}ge^{tA_{\G}}(x,e),
\\
&=
\left.\frac{d}{dt}\right|_{t=0}g(x,\set{U(t)^j_ie_j}),
\\
&=
\left.\frac{d}{dt}\right|_{t=0}(gx,\set{U(t)^j_i\bar{g}e_j}),
\\
&=
\left.\frac{d}{dt}\right|_{t=0}e^{tA_{\G}}(gx,\bar{g}e),
\\
&=
\left.\frac{d}{dt}\right|_{t=0}e^{tA_{\G}}g(x,e),
\\
&=
A_{\G}(g(x,e)).
\end{align*}
\end{answer}
\begin{problem}{all.inv}
Prove that all \(G\)-invariant vector fields on \(\G\) have the form \(A_{\G}\).
\end{problem}
Pick an orthonormal frame \((x_0,e_0)\in\G\).
\begin{problem}{G.to.slash.G}
The map
\[
g\mapsto (x,e):=g(x_0,e_0)
\]
taking a rigid motion to its associated orthonormal frame is a \(G\)-equivariant diffeomorphism.
\end{problem}
so we can think of \(G\) as if it were the orthonormal frame bundle \(\G\):
\[
G\cong\G\to X.
\]
Hence the notation \(\G\) to look like \(G\).
\subsection{Picking an orthonormal frame}
The choice of orthonormal frame \((x_0,e_0)\) gives an isometric linear isomorphism \(X\cong\R^3\).
Using this, write elements \(g\in G\) as \(gx=hx+v\), and so the elements of the Lie algebra \(\LieG\) of \(G\) as \(Ax=Bx+v\) for \(B\in\LieH\), \(v\in\R^3\).
Then the left invariant vector fields are precisely these \(A_{\G}=A_G\).

The flows of left invariant vector fields are right translations in \(G\).
They don't act on \(X\) since the map \(G=\G\to X\) is equivariant under the left \(G\)-action, not under the right action.
The left invariant vector fields associated to infinitesimal rotations are my ``instructions'' to you to pitch, roll and yaw, while those associated to infinitesimal translations are my ``instructions'' to you to move in a particular direction relative to the basis you carry with you.
The flows of the right invariant vector fields are left translations in \(G\), so they act on \(G\) and \(X\), and \(G\to X\) is equivariant, i.e. matches up the actions.
So the right invariant vector fields generate the symmetries of the geometry, while the left invariant vector fields generate the geometry itself: the linear motions and the frame rotations.
\section{The Maurer--Cartan form of the space ship}
\subsection{The soldering forms}
If \(v\) is a tangent vector on the orthonormal frame bundle \(\G\), we can write \(v\) as \((\dot{x},\dot{e})\), an infinitesimal motion \(\dot{x}\) of the point \(x\), and an infinitesimal rotation \(\dot{e}\) of the frame \(e\).
As above, write out any tangent vector \(\dot{x}\) in terms of the basis \(e_1, e_2, e_3\), say as \(\dot{x} = a_1 e_1 + a_2 e_2 + a_3 e_3\).
The \emph{soldering forms}\define{form!soldering}\define{soldering forms} are the \(1\)-forms \(\sigma_1, \sigma_2, \sigma_3\) on \(\G\) given by \((\dot{x},\dot{e})\hook\sigma_i = a_i\).
So the soldering forms measure, as we move a frame, how the base point of the frame moves, as measured in the frame itself as a basis.
In other words, the soldering forms are \(\sigma_i(A_{\G})=v_i\) for the vector fields \(A_{\G}\) above.

The \emph{identity function} on \(X\), which we write as \(x\), is defined as \(x(y)=y\) for any point \(y\in X\).
Of course \(dx(v)=v\) for any vector \(v\), i.e. \(dx=I\) is the identity matrix.
We define \(x\) also on \(\G\) by \(x(y,e)=y\).
We can write the \(1\)-forms \(\sigma_1, \sigma_2, \sigma_3\) on \(\G\) as \(\sigma_1=e_1 \cdot dx, \sigma_2=e_2 \cdot dx, \sigma_3=e_3 \cdot dx\).
On our vector 
\[
v=(\dot{x},\dot{e}_1,\dot{e}_2,\dot{e}_3),
\] 
we have \(v\hook\sigma_i=e_i\cdot\dot{x}=a_i\).

\subsection{The connection forms}
When we move a frame \((x,e)\in\G\), the soldering forms measure the motion of the underlying point \(x\).
We want to measure the rotation of the vectors \(e_1, e_2, e_3\).
Infinitesimal rotations are complicated.
Write the inner product on \(\R^3\) as \(x \cdot y=x_1 y_1 + x_2 y_2 + x_3 y_3\). 
So
\[
e_i \cdot e_j=
\begin{cases}
1 & \text{if \(i=j\)}, \\
0 & \text{if \(i\ne j\)}.
\end{cases}
\]
If we rotate an orthonormal basis \(e=(e_1, e_2, e_3)\) through a family of orthonormal bases 
\[
e(t)=(e_1(t), e_2(t), e_3(t)), 
\]
along some curve \(x(t)\), these still have the same constant values of \(e_i(t) \cdot e_j(t)\) at every time \(t\).
Differentiate: \(0=\dot{e}_i(t) \cdot e_j(t) + e_i(t) \cdot \dot{e}_j(t)\).
Therefore we can write any infinitesimal rotation of frame as \(\dot{e}_j = \sum_i B_{ij} e_i\) for an antisymmetric \(3 \times 3\) matrix \(B=(B_{ij})\).
The quantity \(B_{ij}\) measures how quickly \(e_j\) is moving toward \(e_i\).

The \emph{Levi-Civita connection forms}\define{connection forms!Levi-Civita}\define{Levi-Civita connection!forms} are the \(1\)-forms \(\gamma_{ij}=e_i \cdot de_j\), i.e. \(v \hook \gamma_{ij} = B_{ij}\): so \(\gamma_{ij}\) measures the tendency of \(e_j\) to move toward \(e_i\) as the frame moves.
In particular, \(0=\gamma_{ij}+\gamma_{ji}\).
These \(\sigma_i\) and \(\gamma_{ij}\) are defined on \(\G=G\), \emph{not} on \(\R^3\), because they depend on \(x\) and \(e\).
If we move a frame, we said it moves by a velocity vector
\[
v=(\dot{x},\dot{e})
\]
with \(\dot{e}_i = \sum_j B_{ji} e_j\) for an antisymmetric matrix \(B_{ij}\), so
\[
v\hook\gamma_{ij} = v\hook e_i\cdot de_j = e_i\cdot\dot{e}_j = B_{ij}.
\]
So the soldering \(1\)-forms \(\sigma_i\) restrict to any ``moving frame'' \((x(t),e(t))\) to describe how the velocity of the moving point \(x(t)\) is expressed in the moving frame \(e(t)\), while the \(1\)-forms \(\gamma_{ij}\) describe how the infinitesimal rotation of the moving frame \(e(t)\) is expressed at each moment in the moving frame \(e(t)\).

Write our soldering forms as \(\sigma\), thought of as a column of \(1\)-forms
\[
\sigma
=
\begin{pmatrix}
\sigma_1 \\
\sigma_2 \\
\sigma_3
\end{pmatrix}
\]
and our connection \(1\)-forms as
\[
\gamma=
\begin{pmatrix}
\gamma_{11} & \gamma_{12} & \gamma_{13} \\
\gamma_{21} & \gamma_{22} & \gamma_{23} \\
\gamma_{31} & \gamma_{32} & \gamma_{33} 
\end{pmatrix}
=
\begin{pmatrix}
0 & \gamma_{12} & \gamma_{13} \\
-\gamma_{12} & 0 & \gamma_{23} \\
-\gamma_{13} & -\gamma_{23} & 0
\end{pmatrix}
\]
an antisymmetric matrix of \(1\)-forms, the \emph{Levi-Civita connection}.

\begin{problem}{moving.frame:derive.structure.equations}
Prove that the soldering and Levi-Civita connection forms satisfying the \emph{structure equations of Euclidean space}:\define{structure equations!of Euclidean space}
\begin{align*}
d \sigma_i &= - \sum_j \gamma_{ij} \wedge \sigma_j, \\
d \gamma_{ij} &= -\sum_k \gamma_{ik} \wedge \gamma_{kj}.
\end{align*}
\end{problem}
\begin{answer}{moving.frame:derive.structure.equations}
The proof requires some unwinding of notation: the expression \(\sigma_i=e_i \cdot dx\) means that \(\sigma_i = \sum_j e_{ji} dx_j\), which allows us to unwind the following formal steps:
\begin{align*}
d\sigma_i &= d(e_i \cdot dx),\\
&= de_i \wedge dx, \\
&= \sum_j (e_j \cdot de_i) \wedge (e_j \cdot dx), \\
&= \sum_j \gamma_{ji} \wedge \sigma_j.
\end{align*}
Similarly, the expression \(\gamma_{ij}=e_i \cdot de_j\) means that \(\gamma_{ij} = \sum_k e_{ki} de_{kj}\), so:
\begin{align*}
d\gamma_{ij} &= d(e_i \cdot de_j), \\
&= de_i \wedge de_j, \\
&= \sum_k (e_k \cdot de_i) \wedge (e_k \cdot de_j), \\
&= \sum_k \gamma_{ki} \wedge \gamma_{kj}.
\end{align*}
\end{answer}
Write each element \(g\in G\) as a matrix
\[
g=
\begin{pmatrix}
h&v\\
0&1
\end{pmatrix}.
\]
\begin{problem}{find.m.c}
Check that the identification \(G\cong\G\) identifies
\[
\omega_G=g^{-1}\,dg=
\begin{pmatrix}
\gamma&\sigma\\
0&0
\end{pmatrix}.
\]
\end{problem}
\begin{answer}{find.m.c}
The standard basis at the origin is represented as a matrix in \(G\) by the identity matrix.
The element 
\[
g=
\begin{pmatrix}
h&v\\
0&1
\end{pmatrix}
\]
acts on an orthonormal frame \((x,e)\) to give the orthonormal frame \((hx+v,he)\).
So the map \(G\to\G\) is
\[
g\in G\mapsto g(x_0,e_0)=(v,h).
\]
In other words, if we think of each orthonormal frame \(e\) as a matrix whose columns are \(e_1,e_2,e_3\), and write \(h\) as \(e\) and \(v\) as \(x\), we have
\[
g=
\begin{pmatrix}
e&x\\
0&1
\end{pmatrix}
\in G
\mapsto
(x,e)\in\G.
\]
Hence the Maurer--Cartan form is
\begin{align*}
\omega_G
&=
g^{-1}dg,\\
&=
\begin{pmatrix}
e&x\\
0&1
\end{pmatrix}^{-1}
d
\begin{pmatrix}
e&x\\
0&1
\end{pmatrix},
\\
&=
\begin{pmatrix}
\transpose{e}&-\transpose{e}x\\
0&1
\end{pmatrix}^{-1}
\begin{pmatrix}
de&dx\\
0&0
\end{pmatrix},
\\
&=
\begin{pmatrix}
\transpose{e}de&\transpose{e}dx\\
0&0
\end{pmatrix}.
\end{align*}
\end{answer}
\begin{problem}{moving.frame:structure.group.action}
For \(h\in H\) any \(3 \times 3\) orthogonal matrix, write \(r_h(x,e)\) to mean
\[
r_h(x,e_1,e_2,e_3)=(x,e_1,e_2,e_3)h:=(x,h^j_1e_j,h^j_2e_j,h^j_3e_j)
\]
and \(\transpose{h}\) for the transpose of \(h\).
Thinking of \(e=(e_1,e_2,e_3)\) as a matrix with columns the vectors \(e_1,e_2,e_3\), note that \((x,e)h=(x,eh)\), where \(eh\) means matrix multiplication.
Prove that \(r_h^* \sigma=\transpose{h}\sigma\) and \(r_h^* \gamma = \transpose{h}\gamma h\).
Expanding out, this means \(r_h^*\sigma_i = \sum_j h_{ji} \sigma_j\) and \(r_h^*\gamma_{ij} = \sum_{k\ell} h_{ki} \gamma_{k\ell} h_{\ell j}\).
\end{problem}
\begin{answer}{moving.frame:structure.group.action}
Here are two proofs:
\begin{enumerate}
\item
If we think of \(x\) and \(e\) as functions on \(\G\), then \(r_h^* x = x\), \(r_h^* e = eh\).
Hence \(r_h^* dx=dx\) and \(r_h^* de = (de)h\).
So \(r_h^* \sigma = r_h^* (\transpose{e}\, dx)= \transpose{(eh)} dx = \transpose{h} \transpose{e} dx = \transpose{h} \sigma\) and \(r_h^*\gamma = r_h^* (\transpose{e} de) = \transpose{(eh)} d(eh)=\transpose{h} \transpose{e} de \, h=\transpose{h} \gamma h\).
\item
The action is \(r_h(x,e)=(x,eh)\), where \((eh)_i = \sum_j h_{ji} e_j\).
Hence
\[
r_h'(x,e)(\dot{x},\dot{e})=(\dot{x},\sum_j h_{ji} \dot{e}_j).
\]
\begin{align*}
(r_h^* \sigma)_{(x,e)}(\dot{x},\dot{e})
&=
\sigma_{r_h(x,e)}r_h'(x,e)(\dot{x},\dot{e}),
\\
&=
\sigma_{(x,eh)}(\dot{x},\sum_j h_{ji} \dot{e}_j),
\\
&=
(eh) \cdot \dot{x},
\\
&=\sum_j h_{ji} e_j \cdot \dot{x},
\\
&=\sum_j h_{ji} \sigma_{(x,e)}(\dot{x},\dot{e}).
\end{align*}
\end{enumerate}
\end{answer}

\chapter{Homogeneous spaces}
\section{Definition}\label{page:G.mod.H.trouble}
A \emph{prehomogeneous space}\define{homogeneous!space} \((X,G)\) is a manifold \(X\) acted on smoothly and transitively by a Lie group \(G\).
\begin{problem}{h}
If \((X,G)\) is a prehomogeneous space and \(x(t)\in X\) and \(g(t)\in G\) are differentiable paths, what is the product rule to differentiate \(g(t)x(t)\)?
\end{problem}
A Lie group action of a Lie group \(G\) on a manifold \(X\) is 
\begin{itemize}
\item
\emph{infinitesimally locally transitive}\define{locally transitive!infinitesimally}\define{infinitesimally locally transitive} if the vector fields of the induced Lie algebra action span every tangent space of \(X\)
\item
\emph{locally transitive}\define{locally transitive} if, for every point \(x\in X\) and open set \(U\subseteq G\), \(Ux\subseteq X\) is open.
\item
\emph{component transitive}\define{component transitive} if there is a point of \(X\) whose \(G\)-orbit contains a point on every component of \(X\).
\end{itemize}
A \emph{homogeneous space}\define{homogeneous space} \((X,G)\) is a locally transitive and component transitive Lie group action.
Some authors \cite{Sharpe:1997} call a homogeneous space a \emph{Klein geometry}\define{Klein geometry}\define{geometry!Klein} following \cite{Klein:2008}.
\begin{example}
If a Lie group has at most countably many components, then every orbit of any Lie group action of that Lie group is a homogeneous space \cite{Hilgert.Neeb:2012} p.~367 Proposition 10.1.14.
\end{example}
\begin{problem}{open.orb}
Prove that a Lie group action is infinitesimally locally transitive if and only if it is locally transitive.
\end{problem}
\begin{problem}{homo.sp.def}
Take a prehomogeneous space \((X,G)\).
For some point \(x_0\in X\), let \(H:=G^{x_0}\).
Prove that \(gH\in G/H\mapsto gx_0\in X\) is a bijective immersion.
Prove that it is a diffeomorphism if and only if \((X,G)\) is a homogeneous space.
\end{problem} 
\begin{problem}{bad.G.over.H}
Give an example of a prehomogeneous space which is not homogeneous.
\end{problem}
\begin{answer}{bad.G.over.H}
Take \(G\) to be the group of real numbers with the discrete topology, and \(X\) the real number line with the usual topology.
So \(G\) is zero dimensional, and \(X\) is \(1\)-dimensional.
The action is the usual action of real numbers on real numbers, by addition.
The stabilizer of any real number \(x_0\in X\) is \(H=\set{1}\).
The quotient \(G/H\) is \(G\), of dimension zero, while \(X\) has dimension \(1\).
\end{answer}
Take a prehomogeneous space \((X,G)\).
The stabilizer subgroup \(H:=G^{x_0}\subseteq G\) of any point \(x_0\in X\) is closed.
Conversely, every closed subgroup \(H\subseteq G\) of a Lie group is a smooth Lie subgroup and \((G/H,G)\) is a homogeneous space \cite{Mimura/Toda:1991} p. 44.

\begin{problem}{Lie.G.over.H}
Suppose that a Lie group \(G\) acts smoothly and component transitively on a manifold \(X\).
Prove that the following are equivalent:
\begin{enumerate}
\item
\((X,G)\) is a homogeneous space i.e. the action is infinitesimally locally transitive
\item
some subgroup of \(G\), generated by countably many components of \(G\), acts transitively on some component of \(X\)
\item
every orbit in \(X\) of the identity component of \(G\) contains an open subset of \(X\)
\item
the action is locally transitive
\item
some subgroup of \(G\), generated by countably many components of \(G\), acts locally transitively on some component of \(X\)
\item
the identity component of \(G\) acts transitively on some component of \(X\)
\item
the identity component of \(G\) acts transitively on every component of \(X\).
\end{enumerate}
\end{problem}
\section{Example: projective space}
Consider projective space \((X,G)=(\Proj{n},\PGL{n+1})\), over the real or complex numbers, with \(n\ge 2\).
\begin{theorem}%
[Hilbert \cite{Hilbert1971} pp. 79--81]%
\label{theorem:Hilbert}%
\define{theorem!Hilbert}%
\define{Hilbert's theorem}
Over the real or complex numbers, suppose that \(U, U'\subset\Proj{n \ge 2}\) are connected open sets, and \(U\xrightarrow{\varphi}U'\) is a bijection and that \(\varphi\) maps the intersection of \(U\) with any projective line to the intersection of \(U'\) with some projective line.
If we work over the complex numbers, add the hypothesis that \(\varphi\) is a homeomorphism and preserves orientations of lines.
If we work over the real numbers, we need not even assume that \(\varphi\) is continuous.
Then \(\varphi\) is the restriction of a unique projective transformation in \(\PGL{n+1}\).
\end{theorem}
Intuitively, this theorem says that projective geometry is the (local \emph{and} global) geometry of lines in projective space.
The beautiful proof takes us too far afield.
The uniqueness ensures that projective transformations are effective.
We will describe the local geometry of any homogeneous space, but not so elegantly. 
\section{The bundle}
Take a prehomogeneous space \((X,G)\).
Pick any point \(x_0\in X\); denote the stabilizer as
\[
H:=G^{x_0}:=\set{h\in G|hx_0=x_0}.
\]
Map
\[
\begin{tikzcd}
g\in G\ar[d]\\
gx_0\in X.
\end{tikzcd}
\]
If \(h\in H\), \(hx_0=x_0\) so \((gh)x_0=g(hx_0)=gx_0\), i.e. our map is invariant under right \(H\)-action, so drops to a map \(G/H\to X\).
\begin{problem}{GoverH.X}
Prove that this map \(G/H\to X\) is a diffeomorphism if and only if \((X,G)\) is a homogeneous space \cite{Mimura/Toda:1991} p. 44. 
\end{problem}
\begin{problem}{n}
Prove that \(G\to X\) is a principal right \(H\)-bundle if and only if \((X,G)\) is a homogeneous space.
\end{problem}
\begin{example}
The space \(X\) of all lines in the plane is homogeneous under the action of the group \(G\) of isometries of the plane.
Let \(x_0\) be the horizontal axis.
Then the stabilizer \(H=G^{x_0}\) is the group of translations along the horizontal axis and their products with reflection in the axis.
\end{example}
\begin{example}
Return to the space ship: take \(X:=\R^3\) Euclidean space and \(G\) the group of rigid motions and \(\G\) the orthonormal frame bundle.
Then the bundle \(G\to X\) is identified with the bundle \(\G\to X\) by
\[
g\in G\mapsto g(x_0,e_0)\in\G
\]
where \(x_0=0\) is the origin and \(e_0\) is the standard basis of \(T_{x_0}X=\R^3\).
By analogy, we always think of any Lie group \(G\) acting on any homogeneous space \(X\) as if it were something like a bundle of frames.
\end{example}
\section{Example: the projective line}
\begin{problem}{o}
Prove that the set of all smooth quadric curves in the complex projective plane is a homogeneous space under the group of projective linear transformations.
What happens in the real projective plane?
\end{problem}
Working over the field \(k=\R\) or \(\C\), recall \cite{Needham:2023} p. 137 that the group \(G=\PGL{2}\) (with coefficients in \(k\)) acts on the manifold \(X=\Proj{1}=k\sqcup\infty\) as linear fractional transformations, i.e. each invertible matrix
\[
g
=
\begin{pmatrix}
a&b\\
c&d
\end{pmatrix}
\]
acts on each point \(x\in k\sqcup\infty\) by
\[
gx=\frac{ax+b}{cx+d}.
\]

Since \(g\) and \(\lambda g\) have the same action for any nonzero \(\lambda\in k\), this action is defined on \(G=\PGL{2}=\GL{2}/k^{\times}\).
By such a rescaling \(\lambda g\), we can get any element of \(\GL{2}\) to have determinant \(\pm 1\) if \(k=\R\) or \(1\) if \(k=\C\); let \(\hat{G}\) be the group of matrices with these determinants, so that \(\hat{G}\to G=\PGL{2}\) is a finite covering Lie group morphism, hence the identity map on Lie algebras.
In practice, we carry out all calculations in the group \(\hat{G}\), and then quotient down to \(G=\PGL{2}\).
\begin{problem}{p}
Explain why the Lie algebra \(\LieG=\hat\LieG=\LieSL{2}\) is the collection of traceless \(2 \times 2\) matrices.
Show that it has as a basis
\[
X:=
\begin{pmatrix}
0&1\\
0&0
\end{pmatrix},
H:=
\begin{pmatrix}
1&0\\
0&-1
\end{pmatrix},
Y:=
\begin{pmatrix}
0 & 0\\
1 & 0
\end{pmatrix}.
\]
Compute the Lie brackets of these basis elements.
Compute the exponentials
\[
e^{tX},e^{tH},e^{tY}
\] 
as matrices and as linear fractional transformations.
Differentiate to find the vector fields of the Lie algebra action on the projective line.
Find the left invariant and the right vector fields on \(\hat{G}\) in terms of matrices.
Define the map
\[
g=
\begin{pmatrix}
a&b\\
c&d
\end{pmatrix}
\xmapsto{\pi}
x=\frac{b}{d}.
\]
Show that the right invariant vector fields project under \(\pi\) to vector fields on \(X\) as above, while the left invariant vector fields do not project to any vector fields under \(\pi\).
\end{problem}
\begin{answer}{p}
\(\lb{X}{Y}=H\), \(\lb{H}{X}=2X\), \(\lb{H}{Y}=-2Y\)
Exponentiate:
\[
e^{tX} = 
\begin{pmatrix}
1&t\\
0&1
\end{pmatrix},
e^{tH} = 
\begin{pmatrix}
e^t&0\\
0&e^{-t}
\end{pmatrix},
e^{tY} = 
\begin{pmatrix}
1&0\\
t&1
\end{pmatrix}.
\]
Hence, as linear fractional transformations,
\[
e^{tX}x=x+t, e^{tH}x=e^{2t}x, e^{tY}x=\frac{x}{tx+1}.
\]
Differentiate at \(t=0\) to find the vector fields:
\[
X=\partial_x, H=2x\partial_x, Y=-x^2\partial_x.
\]
The Lie algebra \(\LieSL{2,\R{}}\) arises as the quadratic vector fields 
\[
aX+bH+cY=
\begin{pmatrix}
b&a\\
c&-b
\end{pmatrix}
\mapsto
aX+bH+cY=(a+2bx-cx^2)\partial_x.
\]
We can think of \(2\times 2\) matrices as vectors in \(k^4\), so we can identify tangent vectors to \(\hat{G}\) with matrices.
The left invariant vector fields on \(\hat{G}\) are, at each point
\[
g=
\begin{pmatrix}
a&b\\
c&d
\end{pmatrix},
\]
\begin{align*}
X_G(g)&=gX=
\begin{pmatrix}
0&a\\
0&c
\end{pmatrix},
\\
Y_G(g)&=gY=
\begin{pmatrix}
0&a\\
0&c
\end{pmatrix},
\\
H_G(g)&=gH=
\begin{pmatrix}
a&-b\\
c&-d
\end{pmatrix}.
\end{align*}
\end{answer}
\section{Strong and effective}
The \emph{kernel}\define{kernel} of a homogeneous space \((X,G)\) is the set of elements of \(G\) fixing every point of \(X\).
A homogeneous space is \emph{effective}\define{effective homogeneous space}\define{homogeneous!space!effective} if only \(1\in G\) fixes every point of \(X\), i.e. its kernel is \(\set{1}\), \emph{strong}\define{strong homogeneous space}\define{homogeneous!space!strong} if, for any component of \(X\), any element of \(G\) which fixes every point of that component fixes every point of \(X\).
\begin{example}
A pair of circles form a strong effective homogeneous space under the group of rotations of the two by the same angle together with the permutation of the two.
They form an effective but not strong homogeneous space under the group generated by all pairs of rotations of the two, together with that permutation.
\end{example}
\begin{example}
Ineffective homogeneous spaces are somewhat unnatural, but we will find them difficult to avoid in quotient constructions.
The group of invertible linear transformations of a vector space acts via projective transformations on the associated projective space, but the rescalings act trivially, so the action is not effective.
\end{example}
\section{Disconnected homogeneous spaces}
\par\noindent{}
\begin{example}
A pair of spheres is a homogeneous space under simultaneous rotation of the spheres by perhaps different rotations, or reflections of the spheres, or permuting of the spheres.
Each sphere is already homogeneous under just the rotations and reflections.
\end{example}
\begin{theorem}\label{theorem:connected}
Every component \(X_0\subseteq X\) of a homogeneous space \((X,G)\) is a homogeneous space under the subgroup \(G_0\subseteq G\) leaving it invariant.
Each point of \(X_0\) has the same \(G_0\)-stabilizer subgroup as its \(G\)-stabilizer subgroup.
\end{theorem}
\begin{proof}
If \(X\) is not connected, take any one component \(X_0\subseteq X\) and let \(G_0\subseteq G\) be the set of elements of \(G\) preserving \(X_0\).
Clearly \(G_0\subseteq G\) is a closed subgroup and contains the identity component of \(G\), so is a union of components of \(G\) and shares Lie algebra with \(G\).
Pick a point \(x_0\in X\).
For any point \(x_0\in X_0\), the stabilizer \(H:=G^{x_0}\) fixes \(x_0\), so \(x_0\) stays on the same component, so \(H\subseteq G_0\).
The group \(G_0\) has orbit \(G_0x_0=G_0H x_0\) with preimage \(G_0 H=G_0\) a union of components, and \(G_0\) is both open and closed in \(G\) and right \(H\)-invariant, so \(G_0x_0\subseteq X\) is both open and closed, hence a union of components.
In particular \(X_0\) is an orbit of \(G_0\).
\end{proof}
\begin{problem}{hom.spaces:components}
Take a union
\[
X':=\bigcup_{a\in A} X_a
\] 
of components \(X_a\subseteq X\) of a homogeneous space \((X,G)\).
Let \(G'\subseteq G\) be the set of elements preserving \(X'\).
Prove that \((X',G')\) is a homogeneous space if and only if \(G'\) has an orbit containing a point on each \(X_a\).
\end{problem}
\section{Analyticity}\label{section:analytic}
Every topological group locally homeomorphic to Euclidean space admits a real analytic structure as a real analytic Lie group, unique up to real analytic Lie group isomorphism \cite{Montgomery.Zippin:1974}.
Every closed subgroup is a closed embedded real analytic Lie subgroup \cite{Mimura/Toda:1991} p. 44.
Consequently, every homogeneous space admits an invariant real analytic structure, unique up to equivariant isomorphism.
By analyticity, any element of \(G\) acting trivially on a nondiscrete subset of \(X\) acts trivially on every component of \(X\) containing a limit point of that subset.
In particular, every connected homogeneous space \((X,G)\) is strong.

\section{Universal covering homogeneous spaces}
Take a connected homogeneous space \((X,G)\).
Let \(\tilde{X}\xrightarrow{\pi} X\) be its universal covering space \cite{Munkres:2000} p. 498 corollary 82.2.
Let \(\tilde{G}\) be the set of all pairs \((\varphi,g)\) where \(g\in G\), \(\tilde{X}\xrightarrow{\varphi}\tilde{X}\) is a continuous map lifting \(g\), i.e. \(\pi(\varphi(x))=g\pi(x)\) for all \(x\in\tilde{X}\).
Pick a point \(\tilde{x}_0\in\tilde{X}\) mapping to some \(x_0\in X\).
By standard covering space theory (\cite{Munkres:2000} p. 478 lemma 79.1), for each element \(\tilde{y}\in\tilde{X}\), if \(\pi(\tilde{y})=gx_0\), then there is a unique continuous map \(\tilde{X}\xrightarrow{\varphi}\tilde{X}\) lifting \(g\) so that \(\varphi(\tilde{x}_0)=\tilde{y}\).
Since \(g\) acts as a diffeomorphism, \(\varphi\) is a local diffeomorphism.
Applying the same construction to \(g^{-1}\), we find that \(\varphi\) is a diffeomorphism.
Hence \((\varphi,g)\in\tilde{G}\).
So we have an exact sequence of groups
\[
1\to\pi_1(X)\to\tilde{G}\to G\to 1,
\]
making \(\tilde{G}\) a covering space of \(G\).
There is a unique Lie group structure on \(\tilde{G}\) for which \(\tilde{G}\to G\) is a Lie group morphism and a covering map \cite{Hilgert.Neeb:2012} p. 342 Corollary 9.5.14.
The homogeneous space \((\tilde{X},\tilde{G})\) is the \emph{universal covering homogeneous space}\define{universal covering homogeneous space}\define{homogeneous!space!universal covering} of \((X,G)\).
The stabilizer \(\tilde{G}^{\tilde{x}_0}\subseteq\tilde{G}\) always fixes the ``sheet'' of \(\tilde{x}_0\), so intersects \(\pi_1(X)\) trivially.
Hence this stabilizer is \(H\): \(H\to G\to X\) and \(H\to\tilde{G}\to\tilde{X}\).
\begin{example}
Real projective space is homogeneous under projective transformations, and is connected, but is not simply connected, being covered \(2\)-to-\(1\) by the sphere.
\[
\begin{tikzcd}[background color=exampleBackgroundColour]
S^n \arrow[r,equal] & \tilde{X} \arrow[d] & \tilde{G} \arrow[l] & \GL{n+1}/\R^+ \arrow[l,equal] \\
\RP{n} \arrow[r,equal] & X & G \arrow[l] & \GL{n+1}/\R^{\times} \arrow[l,equal] \\
\end{tikzcd}
\]
If instead we take \(G=\PSL{n+1}\) then \(\tilde{G}=\SL{n+1}\).
\end{example}
\section{Quotient notation}
If \(G\) acts on manifolds \(X,Y\), let \(\amal{X}{G}{Y}\) be the quotient by the diagonal action.
\section{Homogeneous vector bundles}
A \emph{homogeneous vector bundle}\define{homogeneous!vector bundle}\define{vector!bundle!homogeneous} on a homogeneous space \((X,G)\) is a vector bundle on \(X\) for which the action of \(G\) on \(X\) lifts to vector bundle automorphisms.
Pick a point \(x_0\in X\) and let \(H:=G^{x_0}\).
Every \(H\)-module \(V\) gives an \emph{associated vector bundle}\define{associated!vector bundle}\define{vector!bundle!associated} \(\vb{V}:=\amal{G}{H}{V}\).
(See problem~\vref{problem:a.vb.fiber} to see why \(\vb{V}\to X\) is a vector bundle.)
\begin{problem}{q}
Every homogeneous vector bundle is an associated vector bundle and vice versa.
\end{problem}
\begin{answer}{q}
Take a homogeneous vector bundle \(W\to X\).
Since \(G\) acts on \(W\) as vector bundle isomorphisms, \(H:=G^{x_0}\subseteq G\) acts as vector bundle isomorphisms fixing \(x_0\), so \(V:=W_{x_0}\) is an \(H\)-module:
\[
H\xrightarrow{\rho_V}\GL{V}.
\]
Map
\[
(g,v)\in G\times V\to gv\in W.
\]
This map is invariant under the right \(H\)-action
\[
(g,v)h:=(gh,\rho_V(h)^{-1}v)
\]
so descends to a map
\[
\vb{V}:=\amal{G}{H}{V}\to W.
\]
Define a left action of \(G\) on \(G\times V\):
\[
g'(g,v):=(g'g,v),
\]
and note that this commutes with the \(H\)-action, so descends to a left \(G\)-action on \(\vb{V}\).
The map \(\vb{V}\to W\) is a linear isomorphism on the fibers
\[
\vb{V}_{x_0}\to W_{x_0},
\] 
as
\[
\vb{V}_{x_0}:=\amal{H}{H}{V}=V=W_{x_0}.
\]
By \(G\)-equivariance, all fibers are mapped by linear isomorphism.
So the differential of \(\vb{V}\to W\) has no kernel on vertical vectors.
The map \(\vb{V}\to W\) projects to the identity map \(X\to X\), so its differential has trivial kernel.
By dimension count, \(\vb{V}\to W\) is a local diffeomorphism.
Since \(\vb{V}\to W\) is bijective on a fiber, and \(G\)-invariant, it is bijective.
So it is a vector bundle isomorphism.
\end{answer}
\begin{problem}{q.1}
The sections of the associated vector bundle \(\vb{V}=\amal{G}{H}{V}\) are identified with the \(H\)-equivariant functions \(G\to V\).
\end{problem}
\begin{answer}{q.1}
Take a section \(s\) of \(\vb{V}\).
So at each point \(x\in X\), \(s(x)\in\vb{V}_x\), i.e. \(s(x)=(g,v)H\) for some \(g\in G\) and \(v\in V\) with \(gx_0=x\).
Picking any  \(g\in G\), let \(x:=gx_0\), and then \(s(x)=(g',v')H\), say.
But \(g'x_0=x=gx_0\) so \(g'=gh\) for a unique \(h\in H\), so \(s(x)=(g,v)H\) where \(v:=hv'\).
So to each point \(g\in G\), we can associate a unique vector \(v=v(g)\) so that \(s(gx_0)=(g,v(g))H\).
This \(v(g)\) is smooth by local triviality of \(g\in G\mapsto gx_0\in X\).
Moreover \(v(gh)=\rho_V(h)^{-1}v(g)\) by definition.

Conversely, if \(G\xrightarrow{v}V\) is \(H\)-equivariant, we let
\[
s(g):=(g,v(g))H\in\vb{V}_{gx_0},
\]
and note that \(s(gh)=s(g)\) for all \(g\in G\), so \(s\) descends to a map \(X\to\vb{V}\), which is clearly a section.
\end{answer}
\begin{problem}{q.2}
An exact sequence of \(H\)-modules gives rise to and arises from an exact sequence of homogeneous vector bundles.
In particular, homogeneous vector bundles are isomorphic just when their fibers over a point are isomorphic modules of the stabilizer of that point.
\end{problem}
\begin{problem}{hom.vb.pullback.iso}
Define
\[
g\in G\xmapsto{\pi}\pi(g):=gx_0\in X.
\]
Prove that \(\pi^*\vb{V}\) is \(G\)-equivariantly isomorphic to the trivial bundle \(G\times V\) over \(G\).
\end{problem}
\begin{answer}{hom.vb.pullback.iso}
Each element \((g,v)\in G\times V\) quotients to an element \((g,v)H\in\vb{V}:=\amal{G}{H}{V}\) in the fiber \(\vb{V}_x\) where \(x=\pi(g)=gx_0\).
Each element \((g,w)\in\pi^*\vb{V}\) is a pair of \(g\in G\) and \(w=(g_1,v_1)H\in\vb{V}_x\) where
\[
x=\pi(g)=gx_0.
\]
Hence \(\pi(g_1)=g_1x_0=x\).
But then \(g_1=gh\) for a unique \(h\in H\).
So \(w=(g,v)\) where \(v=\rho_V(h)v_1\), and \(v\) is uniquely determined.
Map \((g,v)\to(g,w)\) by this recipe: \(w=(g,v)H\).
This is a smooth map, a linear isomorphism on fibers, and is \(G\)-invariant, hence an isomorphism of vector bundles.
\end{answer}
\begin{proposition}[\cite{Sharpe:1997} p. 188, theorem 3.15]\label{prop:TM.homog}
If \(V:=\LieG/\LieH\), then \(\vb{V}=\vbTM=TX\) 
\end{proposition}
\begin{proof}
Pick a point \(x_0\in X\) and let \(H:=G^{x_0}\).
Left translation gives a commutative diagram of \(H\)-modules:
\[
\begin{tikzcd}
&0\arrow[d] & 0 \arrow[d]\\
0\arrow[r]&T_1 (gH)\arrow[d] \arrow[r,"\omega"] & \LieH \arrow[d]\arrow[r]&0\\
0\arrow[r]&T_1 G\arrow[d,"\pi'(p_0)"] \arrow[r,"\omega"] & \LieG \arrow[d]\arrow[r]&0\\
0\arrow[r]&T_{x_0} X \arrow[r] \arrow[d]& \LieG/\LieH \arrow[d]\arrow[r]&0\\
&0 & 0. \\
\end{tikzcd}
\]
Apply problem~\vref{problem:q.2}.
\end{proof}
\begin{example}
The cotangent bundle is \((\vbTM)^*=\vbh^{\perp}\subseteq\vbg^*\), and similarly for the various tensor bundles.
\end{example}

\section{Homogeneous principal bundles}\label{sec:homog.principal.bundles}
A \emph{homogeneous principal bundle} is a principal bundle \(\Bun\to X\) with a Lie group \(G\) acting smoothly on \(\Bun\) and \(X\) by principal bundle automorphisms, and acting transitive on \(X\) so that \((X,G)\) is a homogeneous space.
\Danger{} the definition does not require that the Lie group \(G\) act transitively on \(\Bun\).
Take a Lie group \(G\) and a closed Lie subgroup \(H\subseteq G\) and a Lie group \(F\) with a Lie group morphism which we denote \(h\in H\mapsto \bar{h}\in F\).
Let
\[
E:=\amal{G}{H}{F}:=(G\times F)/H,
\]
where we quotient out by the right \(H\)-action
\[
(g,f)h=(gh,\bar h^{-1}f).
\]
The left \(G\)-action
\[
g(g_0,f_0):=(gg_0,f_0)
\]
and the right \(F\)-action
\[
(g_0,f_0)f:=(g_0,f_0f)
\]
clearly commute with that \(H\)-action so survive to the quotient.
The map \((g,f)\mapsto gx_0\in X\) is invariant under the \(H\)-action, so descends to the map \((g,f)H\in E\to gx_0\in X\).
\begin{theorem}\label{theorem:build.homo.princ.bun}
Take Lie groups \(G,F,H\) with \(H\subseteq G\) closed and a Lie group morphism \(H\to F\).
Let \(E:=\amal{G}{H}{F}\), \(X:=G/H\).
There is a unique smooth structure on \(E\) so that the quotient map 
\[
G\times F\to E=\amal{G}{H}{F}
\]
is a smooth principal right \(H\)-bundle.
The \(F\)-action and \(G\)-action descend to smooth actions on \(E\), making the quotient map \(E\to X\) a principal right \(F\)-bundle, invariant under the \(G\)-action, so a homogeneous right principal bundle.
\end{theorem}
\begin{problem}{invariant.bundles}
Prove it.
\end{problem}
\begin{answer}{invariant.bundles}
Since \(H\subseteq G\) is a closed subgroup, the right \(H\)-action on \(G\) is free and proper, so the right \(H\)-action on \(G\times F\) is free and proper, hence the quotient space \(\Bun:=\amal{G}{H}{F}\) is a smooth manifold and the quotient map
\[
(g,f)\in G\times F\mapsto (g,f)H\in\amal{G}{H}{F}
\]
is a smooth principal right \(H\)-bundle, for a unique smooth structure on \(\Bun\).
The map \(\Bun\to G/H\) is the quotient of the smooth map 
\[
(g,f)\in G\times F\mapsto gH\in G/H
\]
so is smooth.
The \(G\)-action and the \(F\)-action on \(\Bun\) arise from smooth actions on \(G\times F\), so are smooth.
The action of \(F\) on \(\Bun\) is easily seen to be free.

Let us check that \(F\) acts properly on \(\Bun\).
Suppose that in \(\Bun\) we have points
\begin{align*}
\bar{p}_i&\to\bar{p},\\
\bar{p}'_i&=\bar{p}_if_i'\to\bar{p}'
\end{align*}
for some elements \(f_i'\in F\).
To prove properness of the \(F\)-action on \(\Bun\), we have to prove that some infinite subsequence of these \(f_i'\) converges.

Lift these up to sequences \(p_i,p'_i\), perhaps not converging, in \(G\times F\).
By local triviality of the \(H\)-action, we can alter these \(p_i\) by elements of \(H\) to ensure that they converge, say \(p_i\to p\).
By the same argument, there are elements \(h_i\in H\) so that \(p_i'h_i\) converges, say \(p_i'h_i\to p'\).
Write these points of \(G\times F\) as pairs
\begin{align*}
p_i=(g_i,f_i)&\to p=(g,f), \\
p_i'h_i=p_if_i'h_i=(g_ih_i,\bar h_i^{-1}f_if_i')&\to p'=(g',f').
\end{align*}
Take the first components of these pairs: \(g_i\to g\) and \(g_ih_i\to g'\).
So \(g_i^{-1}\to g^{-1}\) and so \(h_i=g_i^{-1}g_ih_i\to g^{-1}g'\).
So then
\[
\bar h_i\to\overline{g^{-1}g'},
\]
and so
\[
f_if_i'\to \overline{g^{-1}g'}f',
\]
and \(f_i^{-1}\to f^{-1}\) so
\[
f_i'\to f^{-1}\overline{g^{-1}g'}f'.
\]
So \(F\) acts freely and properly on \(\Bun\).

Therefore \(F\to\Bun\to\Bun/F\) is a principal right \(F\)-bundle.
Since \(G\times F\) acts smoothly and transitively on \(\Bun=\amal{G}{H}{F}\), \(G\) acts smoothly and transitively on \(\Bun/F\).
Take the identity components \(G_0\subseteq G\), \(F_0\subseteq F\).
Then \(G_0\times F_0\) is the identity component of \(G\times F\), so acts transitively on a component of \(\Bun=\amal{G}{H}{F}\), and so \(G_0\) acts transitively on a component of \(\Bun/F\), so \(\Bun/F\) is a homogeneous space of \(G\).

The map \(\Bun\to X:=G/H\) is clearly \(F\)-invariant, so descends to a smooth map
\[
\Bun/F\to X,
\]
clearly \(G\)-equivariant.
The point \(\bar{p}_0F\in\Bun/F\) has stabilizer \(H\), so this is an isomorphism of \(G\)-homogeneous spaces.
\end{answer}
\begin{theorem}\label{theorem:build.homo.princ.bun.2}
Every homogeneous right principal bundle arises as in theorem~\vref{theorem:build.homo.princ.bun}.
\end{theorem}
\begin{problem}{invariant.bundles.2}
Prove it.
\end{problem}
\begin{answer}{invariant.bundles.2}
Suppose that we have a principal right $F$-bundle $\Bun\xrightarrow{\pi}X$, with an action of a Lie group $G$ by principal bundle automorphisms, acting transitively on \(X\), so that \((X,G)\) is a homogeneous space.
Since \(G\) acts by principal bundle automorphisms, it commutes with the right \(F\)-action on \(\Bun\).
Pick a point \(p_0\in\Bun\) and let \(x_0\in X\) be its projection to \(X\).
Let \(H:=G^{x_0}\).
Map
\[
h\in H\mapsto hp_0\in\Bun.
\]
Since \(h\in H=G^{x_0}\), \(\pi(hp_0)=h\pi(p_0)=hx_0=x_0\), so \(hp_0\in\Bun_{x_0}\).
Since \(F\) acts freely and transitively on the fibers of \(\Bun\to X\), \(hp_0=p_0f\) for a unique element \(f\in F\), which we write as \(f=\bar h\), for a unique map \(H\to F\). 
Using a local trivialization of \(\Bun\to M\), we see that \(f\in F\) is a smooth function of \(p_0f\), hence of \(h\): \(H\to F\) is smooth. 

Note that
\[
p_0=h^{-1}(hp_0)=h^{-1}(p_0\bar h)=(h^{-1}p_0)\bar h=p_0\overline{h^{-1}}\bar h,
\]
so
\[
\overline{h^{-1}}=\bar h^{-1}.
\]
Composing, if \(f_0=\bar h_0\) and \(f_1=\bar h_1\),
\[
h_1(h_0p_0)=h_1(p_0f_0)=(h_1p_0)f_0=(p_0f_1)f_0=p_0(f_1f_0).
\]
So \(H\to F\) is a group morphism, hence a smooth Lie group morphism.

Map
\[
(g,f)\in G\times F\mapsto gp_0f\in\Bun.
\]
This map is onto, with fiber over \(p_0\) precisely the set of pairs \((g,f)\) so that \(gp_0=p_0f\), i.e. \(g\in G^{x_0}=H\) and \(f=\bar h\).
If \(h\in H\), 
\[
(g,f)h=(gh,\bar h^{-1}f)\mapsto ghp_0\bar h^{-1}f=ghh^{-1}p_0f=gp_0f.
\]
So the map is \(H\)-invariant, so descends to a smooth map \(\amal{G}{H}{F}\to\Bun\).
The original map is \(G\)-equivariant and \(F\)-equivariant, and these actions commute with the \(H\)-action, so these equivariances survive the quotient.

We want to see that \(\Bun\) is a homogeneous space of \(G\times F\), under the left action
\[
(g,f)p=gpf^{-1}.
\]
If we prove this then clearly the stabilizer is precisely \(H\), embedded by \((h,\bar h)\), and so the map \(\amal{G}{H}{F}\to\Bun\) is an isomorphism of homogeneous spaces, so a \(G\)-equivariant \(F\)-bundle isomorphism.

Since \((X,G)\) is a homogeneous space, we can suppose that some subgroup \(G_0\subseteq G\) generated by countably many components of \(G\) acts transitively on some component \(X_0\subseteq X\).
After perhaps conjugating \(G_0\), we can assume that \(X_0\) is the component of \(x_0\).
Let \(\Bun_0\subseteq\Bun\) be the component of \(p_0\in\Bun\).
Let \(F_0\subseteq F\) be the identity component.
Then the group \(G_0\times F_0\) acts on \(\Bun\), preserving \(\Bun_0\).
Since the action is transitive on \(X_0\), it acts on \(\Bun_0\) taking any fiber to any other fiber, with orbits transverse to every fiber.
Locally trivializing \(\Bun\to X\) near \(x_0\), we see that the action of \(F_0\) is locally transitive on every fiber of \(\Bun\to X\).
Hence the action of \(G_0\times F_0\) has an open orbit in \(\Bun_0\) through every point.
Since \(\Bun_0\) is connected, it is an orbit of \(G_0\times F_0\).
Hence \(\Bun\) is a homogeneous space of \(G\times F\).
\end{answer}
\begin{theorem}
Take a homogeneous right principal bundle \(\Bun\to X\).
As above write it as \(\Bun=\amal{G}{H}{F}\) for Lie groups \(G,F,H\) with \(H\subseteq G\) closed and a Lie group morphism \(H\to F\), denoted \(h\mapsto \bar h\).
The bundle is equivariantly trivial if and only if \(\bar h=1\) for all \(h\in H\).
\end{theorem}
\begin{problem}{invariantly.trivial}
Prove it.
\end{problem}
\begin{answer}{invariantly.trivial}
Take a \(G\)-equivariant trivialization
\[
\Bun\to X\times F,
\]
and then follow the construction above to give an isomorphism \(\amal{G}{H}{F}=\Bun\) associated to a choice of point \(p_0\in\Bun\) and its image \(x_0\in X\) with \(H=G^{x_0}\).
Compose with \(G\times F\to \amal{G}{H}{F}\) to get a map
\[
G\times F\to (G/H)\times F
\]
which is \(G\)-equivariant, \(F\)-equivariant and \(H\)-invariant.
Note that 
\[
(1,1)\in G\times F\mapsto p_0\in\Bun\mapsto x_0\in X
\]
so
\[
(1,1)\in G\times F\mapsto (x_0,f_0)\in X\times F,
\]
for some \(f_0\in F\).
Applying the \(G\)-equivariance and the \(F\)-equivariance,
\[
(g,f)=g(1,1)f\in G\times F\mapsto g(x_0,f_0)f=(gx_0,f_0f)\in X\times F.
\]
By \(H\)-invariance,
\[
(g,f), (gh,\bar h^{-1}f)
\]
map to the same place, i.e.
\[
f_0\bar h^{-1}f=f_0f,
\]
so \(\bar h=1\) for all \(h\in H\).
Reverse steps to find that if \(\bar h=1\) for all \(h\in H\), then \(\Bun\to M\) is \(G\)-invariantly trivial.
\end{answer}
\begin{theorem}\label{thm:homog.principal.bundle.isos}
Take two homogeneous right principal bundles \(\Bun,\Bun'\to X\).
As above, by picking a point \(x_0\in X\) and points \(p_0\in\Bun\), \(p_0'\in\Bun'\), write them as \(\Bun=\amal{G}{\alpha}{F}\), \(\Bun'=\amal{G}{\alpha'}{F}\) for Lie groups \(G,F,H\) with \(H=G^{x_0}\subseteq G\) and Lie group morphisms \(H\xrightarrow{\alpha,\alpha'}F\)
There is an equivariant bundle isomorphism \(\Bun\to\Bun'\) if and only if there is an element \(f_0\in F\) so that
\[
\alpha'(h)=f_0\alpha(h)f_0^{-1},
\]
for all \(h\in H\).
\end{theorem}
\begin{problem}{equivar.equival}
Prove it.
\end{problem}

\chapter{Locally homogeneous structures}
\section{Definition}
Take a homogeneous space \((X,G)\) and a manifold \(M\) with \(\dim M=\dim X\).
An \emph{\(X\)-chart}\define{X-chart@$X$-chart} on \(M\) is a diffeomorphism
\[
\text{open }\subseteq M\to\text{open }\subseteq X.
\]
\[
\input{xg-chart}
\]
Two \(X\)-charts are \emph{\(G\)-compatible}\define{G-compatible@$G$-compatible} if, on any connected open set where both are defined, they agree up to action of an element of \(G\) (a unique element if \((X,G)\) is strong and effective).
\[
\input{xg-struc}
\]
An \emph{\((X,G)\)-atlas}\define{(X,G)-atlas@$(X,G)$-atlas} is a collection of \(G\)-compatible \(X\)-charts whose domains cover \(M\).
An \emph{\((X,G)\)-structure}\define{(X,G)-structure@$(X,G)$-structure} is a maximal \((X,G)\)-atlas, also called a \emph{locally homogeneous structure modelled on \((X,G)\)}.\define{locally homogeneous structure}
If \(G\) does not act effectively on \(X\), then an \((X,G)\)-structure is precisely an \((X,G/K)\)-structure, \(K\subseteq G\) the kernel of \((X,G)\).
\begin{marginfigure}
\includegraphics[width=2cm]{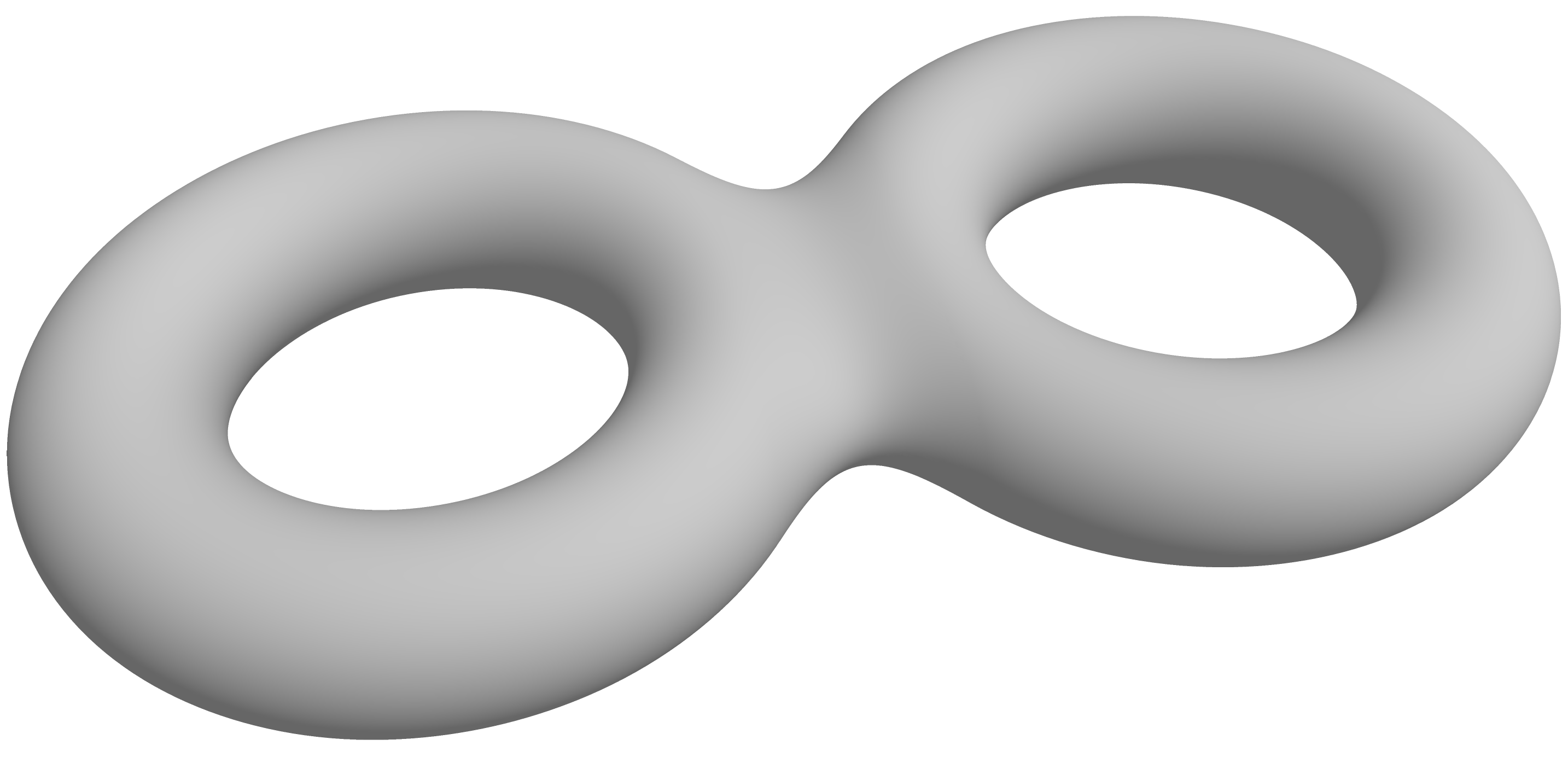}
\end{marginfigure}
\begin{example}
By Poincar\'e's uniformization theorem \cite{Saint-Gervais:2010}, every connected surface admits a complete Riemannian metric of constant curvature everywhere equal to \(k_0=1\), \(0\) or \(-1\).
Any two such, with the same constant \(k_0\), are locally isometric.
Such a metric is an \((X,G)\)-structure where \(X\) is the real projective plane, Euclidean plane, or hyperbolic plane and \(G\) its isometry group.
To see this, note that every local isometry of connected open sets of \(X\) arises from a global isometry of \(X\).
\end{example}
\section{Example: 2\textsuperscript{nd} order ordinary differential equations}\label{example:linear.ode}
Take any 2\textsuperscript{nd} order complex analytic linear ordinary differential equation
\[
0=w''+pw'+qw
\]
with \(w=w(z)\), \(p=p(z)\), \(q=q(z)\) defined for \(z\) in some region of the complex plane.
We can imagine this $z$ to be merely a local coordinate on a Riemann surface $M$.
To make sense of the linearity of the differential equation, imagine that $w$ is the fiber coordinate in a local trivialization of a holomorphic line bundle $L\to M$.
Imagine that our differential equation extends to be defined on every local trivialization of $L$ in every local coordinate domain on $M$.

The local solutions of the differential equation, in some disk around a point for example, form a \(2\)-dimensional complex vector space \(V\).
But there might not be so many global solutions on all of $M$: a local solution might, as we travel around a loop on $M$, return to have a different value.
For example,
\[
w''+(1/4z^2)w=0
\]
has solution $w=\sqrt{z}$, not globally defined on the complex plane punctured at the origin, but defined on the universal covering space.
The universal covering space $\tilde{M}\to M$ has the pullback differential equation, for sections of the pullback line bundle, with a \(2\)-dimensional vector space $V$ of solutions.

For each point $z=z_0\in\tilde{M}$, there is global solution \(w(z)\) vanishing at that point: \(w(z_0)=0\), with nonzero first derivative $w'(z_0)\ne 0$, i.e. not vanishing everywhere.
That solution is uniquely determined up to constant scaling.
The line in $V$ spanned by that solution is a point of the projective line $X:=\mathbb{P}V=(V-\set{0})/\mathbb{C}^{\times}$.

The \emph{developing map}\SubIndex{developing map} maps each point $z=z_0\in\tilde{M}$ to that solution, i.e. to that point of $X$.
We leave the reader to puzzle out why the developing map is a local biholomorphism of Riemann surfaces.
The deck transformations of $\tilde{M}$ over $M$ are, by definition, automorphisms of the ordinary differential equation on $\tilde{M}$, since it is pulled back from $M$.
So deck transformations act on the space of solutions by linear transformation, hence on $X$ by linear fractional transformation, i.e. projective transformation.

In other words, $M$ bears a holomorphic projective structure, i.e. a holomorphic $(X,G)$-structure, for $X$ the complex projective line and $G$ the group of linear fractional transformations.
In fact, by construction, we can take $G$ to be the group of unimodular linear transformations acting on the global solutions of the differential equation, so we can see that any holomorphic projective structure arising in this way lifts to some such structure.

In particular, rather surprisingly, the cross ratio of any \(4\) points of \(\tilde{M}\) is well defined, a global invariant of any holomorphic 2nd order ordinary differential equation.

Conversely, any holomorphic $(X,G)$-structure has an associated line bundle, the $\mathcal{O}(1)$-bundle on $X=\mathbb{P}^1_{\mathbb{C}}$, i.e. the line bundle whose global sections are the complex linear functions on \(V\).

Note that this really uses that $G$ is the group of unimodular linear transformations of \(V\), not the linear fractional transformations of \(X\), since the projective linear group does \emph{not} act on this line bundle.
This is an important example to explain why we often need to allow our homogeneous spaces \((X,G)\) to have group \(G\) not acting effectively.

The line bundle $\mathcal{O}(1)$ has a \(2\)-dimensional space of global sections, which we take to be our solutions of our differential equation: the differential equation is just \(w''=0\) on projective linear charts on \(X\).
Local \((X,G)\)-charts identify the differential equation on \(M\) with this one, and the transition maps glue them together.

Take any holomorphic line bundle \(L\to M\) on a Riemann surface \(M\), and a \(2\)-dimensional complex vector space $V$ of global holomorphic sections of the pullback \(L\to\tilde{M}\) on the universal covering space \(\tilde{M}\to M\).
Take two sections from \(V\); in local coordinates, write them as holomorphic functions \(w_1,w_2\).
They are linearly dependent as elements of \(V\) just when they satisfy a constant coefficient linear relation 
\[
0=\alpha_1w_1(z)+\alpha_2w_2(z)
\]
for all \(z\) in some open set, with \(\alpha_1,\alpha_2\) not both zero.
\begin{lemma}
Holomorphic functions \(w_1,w_2\) are linearly independent just when the \emph{Wronskian}\define{Wronskian}
\[
\det
\begin{pmatrix}
w_1&w_2\\
w_1'&w_2'
\end{pmatrix}
\]
vanishes only on a discrete set of points.
\end{lemma}
\begin{proof}
A nonzero linear relation makes one function a constant multiple of the other: zero Wronskian.
Suppose they are linearly independent.
If they vanish to the same order, rescale to get the same lowest order term, and replace one of them by their difference.
So they vanish to different orders.
The lowest order term of the Wronskian is the Wronskian of the lowest order terms: not zero.
\end{proof}
Define a differential equation
\[
0=\det
\begin{pmatrix}
w&w'&w''\\
w_1&w_1'&w_1''\\
w_2&w_2'&w_2''
\end{pmatrix},
\]
which clearly has a nonzero coefficient in front of \(w''\), so is a nontrivial second order ordinary differential equation satisfied by all elements of \(V\), after perhaps deleting a discrete set from our Riemann surface.

For examples with explicit computation of the action of deck transformations, see the beautiful \cite{Kuga:1993}.

\epigraph[author={Galileo Galilei},
source={Dialogo sopra i due massimi sistemi del mondo tolemaico e copernicano},translation={Shut yourself up with some friend in the main cabin below decks on some large ship, \ldots throwing something to your friend, you need to throw it no more strongly in one direction than another, the distances being equal; jumping with your feet together, you pass equal spaces in every direction.}]
{Riserratevi con qualche amico nella maggiore stanza che sia sotto coverta di alcun gran navilio, \ldots gettando all’amico alcuna cosa, non più gagliardamente la dovrete gettare verso quella parte che verso questa, quando le lontananze sieno eguali; e saltando voi, come si dice, a piè giunti, eguali spazii passerete verso tutte le parti.}
\begin{marginfigure}[-7.25cm]\includegraphics[width=3cm]{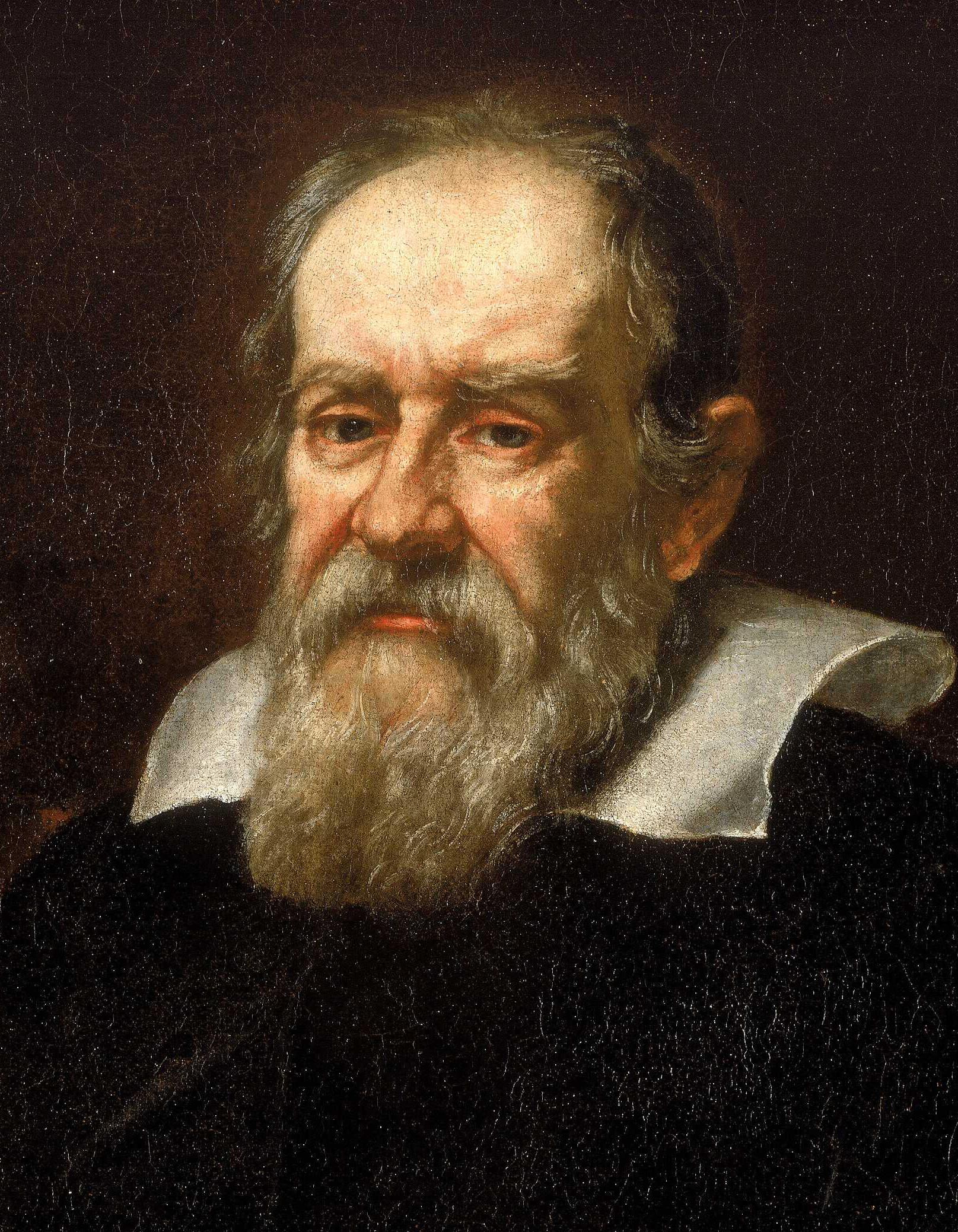}\end{marginfigure}

\begin{problem}{x}
Recall that affine space is the homogeneous space
\[
(X,G)=(\R^n,\GL{n}\ltimes\R^n),
\]
i.e. Euclidean space acted on by translations and invertible linear transformations.
An \emph{affine structure}\define{affine structure} is an \((X,G)\)-structure modelled on affine space.
Prove:
\end{problem}
\begin{theorem}[Galilei's Principal of Relativity]%
\define{Galilei's principal of relativity}
There is an equivalence of categories between (i) pairs of affine structures on a curve and (ii) a pair consisting of an affine structure on that curve and a \(1\)-form.
The equivalence maps a pair of affine structures, with local affine charts \(x,y\), to the pair consisting of the affine structure of \(x\) and the \(1\)-form
\[
\frac{y''(x)}{y'(x)}dx.
\]
\end{theorem}
\begin{problem}{y}
How does Galilei's theorem work in all dimensions?
\end{problem}
\begin{example}
Look back at the biholomorphisms of the unit ball \(X\subseteq\C^n\)~\vpageref{example:unit.ball}.
Each biholomorphism has numerator a linear map, denominator a linear function, so is a complex projective transformation extending uniquely to \(\Proj{n}\).
So any smooth real \((X,G)\)-structure, where \(G\) is the biholomorphism group of the unit ball, determines a holomorphic projective structure, i.e. a structure modelled on \((\Proj{n},\PSL{n+1})\).
\end{example}
There are some very beautiful examples of affine structures \cite{Sullivan/Thurston:1983}.
\section{Kleinian groups}
Suppose that a group \(\Gamma\) acts by homeomorphisms on a topological space \(X\).
A \emph{\(\Gamma\)-house}\define{Γ-house@\(\Gamma\)-house} or \emph{house}\define{house} of a point \(x\in X\) is an open set \(U\subseteq X\) which does not intersection any of its translates \(gU\), for any \(g\in\Gamma\) unless \(g\) fixes every point of \(X\) \cite{Siegel:1943}.
The \emph{free regular set}\define{free regular set} \(\Omega=\Omega_{\Gamma,X}\) of the group action is the set of points which lie in houses.
\begin{problem}{house.covering}
If \(\Omega\) is not empty then prove that \(\Omega\to\bar\Omega:=\Gamma\backslash\Omega\) is a covering space and \(\Gamma\) acts on \(\Omega\) as a group of deck transformations.
\end{problem}
\begin{problem}{house.covering.2}
If \(\Omega\) is not empty and \(X\) is a manifold and \(\Gamma\) acts on \(X\) by diffeomorphisms, then prove that the quotient map \(\Omega\to\bar\Omega\) is a local diffeomorphism for a unique smooth structure on \(\bar\Omega\).
\end{problem}
\begin{example}
Take a homogeneous space \((X,G)\).
A subgroup of \(G\) is a \emph{Kleinian group}\define{Kleinian!group} if its free regular set is not empty \cite{Maskit:1988} p. 15.
The most important examples of locally homogeneous structures are the \emph{Klein manifolds}:\define{Klein!manifold} quotients \(\Gamma\backslash\Omega\) of free regular sets of Kleinian groups.
\end{example}
\begin{problem}{z}
If the free regular set is not empty, prove that \(\Gamma\subseteq G\) is discrete.
\end{problem}
\marginnote{\vspace{1cm}%
\begin{tikzpicture}
\draw[white,thick,inner color=white,outer color=gray!40] (0,0) circle (1cm);
\draw[white,thick,inner color=white,outer color=gray!40] (0,0) circle (.5cm);
\draw[white,thick,inner color=white,outer color=gray!40] (0,0) circle (.25cm);
\draw[white,thick,inner color=white,outer color=gray!40] (0,0) circle (.125cm);
\draw[white,thick,inner color=white,outer color=gray!40] (0,0) circle (.0625cm);
\end{tikzpicture}}
\begin{example}
Take the group \(G=\R^n\ltimes\GL{n}\) of affine transformations of affine space \(X=\R^n\).
Take any affine transformation \(g\) of \(X\), and let \(\Gamma\subseteq G\) be the subgroup generated by \(g\).
If the linear part of \(g\) has spectrum in the unit disk then \(g\) has a unique fixed point \(x_0\in X\) and \(\Gamma\) has free regular set \(\Omega:=X-\set{x_0}\) with quotient the \emph{Hopf manifold}\define{Hopf manifold} \(M:=\Gamma\backslash\Omega\).
Without loss of generality, \(x_0=0\) and \(g\) is linear.
In the special case where \(g\) is a scalar multiple of the identity \(g=\lambda I\), the Hopf manifold is the quotient of the annulus \(|\lambda|\le|x|\le 1\) gluing the inner and outer spheres.
Starting from any affine transformation \(g\) whose linear part has spectrum in the unit disk, we can deform smoothly through such affine transformations until we get to one which is a multiple of the identity.
As we deform our affine transformation \(g\), keeping the spectrum in the unit disk, the Hopf manifold does not change diffeomorphism type, by Ehresmann's theorem (corollary~\vref{corollary:Ehresmann.thm}).
Hence the Hopf manifold is a compact manifold, diffeomorphic to a product of a circle and a sphere.
The Hopf manifold is acted on by the group of all affine transformations commuting with \(g\) (modulo \(g\) itself) as symmetries of its \((X,G)\)-structure.
In particular, if \(g=\lambda I\) is a constant rescaling, then the Hopf manifold is acted on by \(\GL{n}/\left<g\right>\), transitively, so these special Hopf manifolds are homogeneous.
There is one very special case: \(n=1\), we glue the nonzero real numbers together by a dilation, giving a homogeneous affine structure on the circle, not equivalent to the usual affine structure which arises by quotienting \(\R/\Z\).
\end{example}
\begin{example}
Picture the unit disk in the complex plane.
Picture the complex projective line as the Riemann sphere.
Picture the complex plane as the complex affine line, i.e. the Riemann sphere punctured at infinity.
Hence we embed the unit disk into the complex plane, and the complex plane into the Riemann sphere.
A \emph{Riemann surface}\define{Riemann surface} is a connected complex manifold of complex dimension 1.
By Poincare's uniformization theorem \cite{Saint-Gervais:2010}, every Riemann surface $M$ has universal covering space $\tilde{M}$ the disk, the complex plane, or the complex projective line (i.e. the Riemann sphere). 
In any case, we see $\tilde{M}$ sitting in the complex projective line as an open set.
The quotient map $\tilde{M}\to M$ is the quotient by deck transformations $\tilde{M}\to\tilde{M}$.
Look at all of the biholomorphisms of the disk, of the plane, and of the Riemann sphere: each is a linear fractional transformation of the Riemann sphere, i.e. a complex projective linear transformation of the complex projective line.
Hence each Riemann surface $M$ has a \emph{canonical} holomorphic projective structure, i.e. an $(X,G)$-structure where $X$ is the complex projective line, and $G$ the complex projective linear transformations of $X$.
This structure is holomorphic in the sense that all of the charts of the structure are holomorphic maps of Riemann surfaces.
\end{example}
\marginnote{\vspace{1cm}\includegraphics[width=2cm]{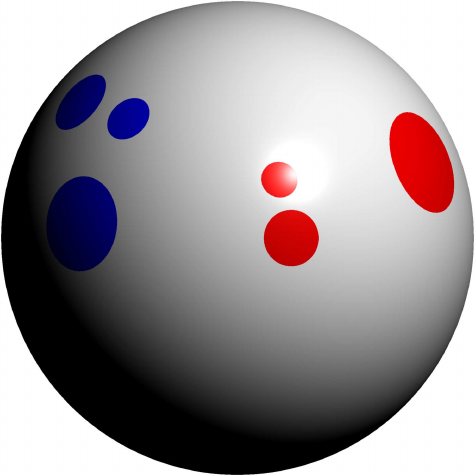}}
\begin{example}
Think of the sphere as the Riemann sphere.
Pick a disk on the sphere.
Inversion in the boundary circle of our disk is conjugate holomorphic.
Compose that inversion with a reflection across an equator through the center of that disk.
The composition is a linear fractional transformation, interchanging the interior of the disk with its exterior.

Fix a positive integer $g$. On the sphere, draw $2g$  disjoint disks, half of which we draw in blue, the other half in red.
Associate bijectively red disks to blue ones.
For each blue disk, pick a linear fractional transformation of the Riemann sphere taking that blue disk to the complement of its associated red disk.
The group generated by these linear fractional transformations is a discrete group called a \emph{classical Schottky group}.\define{Schottky group!classical}\define{classical Schottky group}
Each classical Schottky group has free regular set including the complement of the disks, so the quotient of that complement is a compact Riemann surface of genus $g$ with a holomorphic projective structure \cite{Maskit:1988} p. 82, J20, p. 168, F8.
There is a family of such groups: for nearby classical Schottky groups, we can allow ourselves to use the same blue disks.
Any linear fractional transformation near a generator carries the blue disk of that generator to some disk, which we colour red, and hence we construct a new classical Schottky group.
Each linear fractional transformation depends on \(3\) complex parameters.
Clearly we get an isomorphic classical Schottky group if we conjugate all of the generators by the same linear fractional transformation.
So the space of classical Schottky groups on \(2g\) disks depends on \(3g-3\) complex parameters, up to isomorphism.
\end{example}
\begin{example}
Returning to the Siegel upper half space, the group \(\Gamma\) of integer symplectic matrices acts as change of symplectic basis of the homology of a compact Riemann surface, hence the periods of a compact Riemann surface form an element of \(\Gamma\backslash X\) \cite{Arbarello.Cornalba.Griffiths.Harris:1985} p. 22, \cite{Griffiths.Harris:1978} p. 333, \cite{Siegel:3} chapter 6.
\end{example}
\begin{example}
The only groups of Lorentz transformations of de~Sitter space which act with the entire de~Sitter space as free regular set are finite \cite{Calabi.Markus:1962}.
There are no compact mannifolds with Lorentz metric locally isometric to de~Sitter space, i.e. with de~Sitter space as universal covering space \cite{Klingler:1996}.
The compact manifolds with Lorentz metric locally isometric to anti-de~Sitter space are not yet classified.
\end{example}
\begin{example}
Suppose that a group \(\Gamma\) acts by homeomorphisms on a topological space \(X\).
Let \(K\subseteq\Gamma\) be the elements that fix every point of \(X\).
A \emph{\(\Gamma\)-shack}\define{Γ-shack@$\Gamma$-shack} or \emph{shack}\define{shack} of a point \(x\in X\) is an open set \(U\subseteq X\) intersecting only finitely many of its translates \(gU\) for \(g\in\Gamma/K\).
It might be helpful to consider the set \(\Omega'\supseteq\Omega\) of points \(x\in X\) lying in shacks, so that \(\Gamma\backslash\Omega'\) is an orbifold with locally homogeneous structure, the true object of study in the field of Kleinian groups, apparently.
(I don't know how to define an orbifold Cartan geometry, but maybe we should.)
\end{example}
\begin{center}\small
\begin{tabular}{p{5.5cm}p{5.5cm}}
\includegraphics[width=5.4cm]{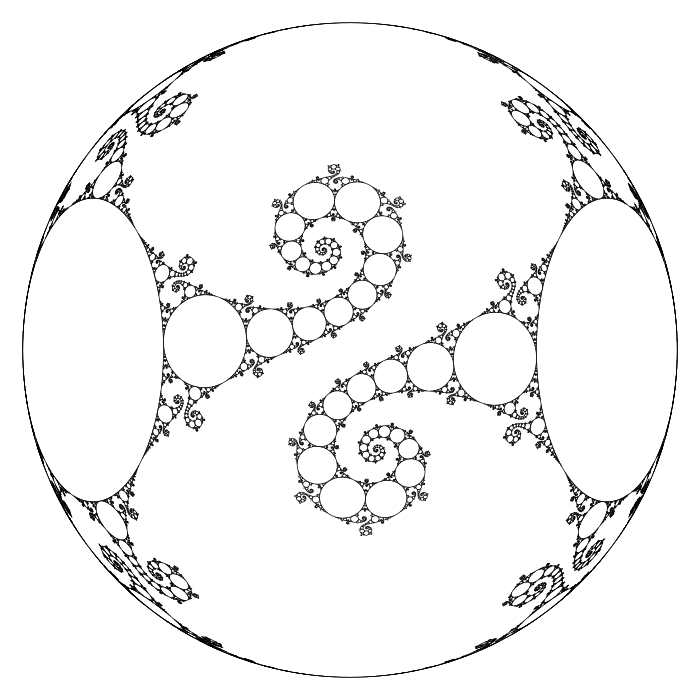}
&\includegraphics[width=5.4cm]{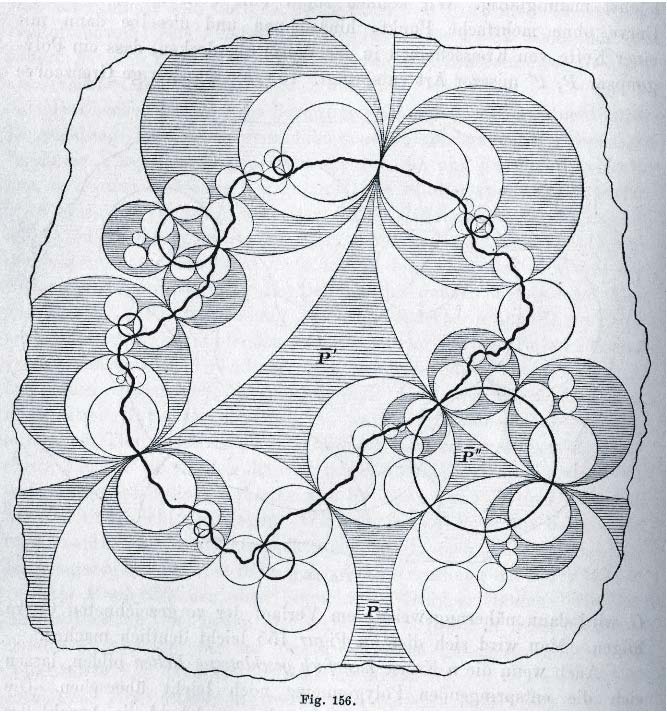}\\
\begin{tabular}{p{5.4cm}}
Adam majewski, CC BY-SA 4.0 \\ 
\url{https://creativecommons.org/licenses/by-sa/4.0}, \\
via Wikimedia Commons
\end{tabular}& 
\begin{tabular}{p{5.4cm}}
R. Fricke, F. Klein, \textbf{Lectures on the theory of automorphic functions}. Vol. 1. 
\cite{Fricke.Klein:2017}
\end{tabular}
\end{tabular}
\end{center}
\begin{example}
A Klein manifold \(\Gamma\backslash\Omega\) is \emph{complete}\define{complete!Klein manifold}\define{Klein!manifold!complete} if \(\Omega=X\), i.e. \(\Gamma\backslash X\) where \(\Gamma\subseteq G\) is any discrete subgroup acting on \(X\) as deck transformations. \Danger{}Sharpe calls these, surprisingly, \emph{locally Klein geometries}\define{locally Klein geometry} \cite{Sharpe:1997} p. 154 definition 3.10.
Goldman calls them \emph{tame geometries} \cite{Goldman:2022} p. 113.
Borel calls them \emph{Clifford--Klein space forms}\define{Clifford-Klein space form@Clifford--Klein space form} \cite{Borel:1963}, which has become the standard term.
But I think the term \emph{complete Klein manifold} indicates more clearly that we quotient the complete space \(X\).
A Lie group \(G\) is \emph{linearly reductive}\define{Lie!group!linearly reductive}\define{linearly reductive!Lie group} if every finite dimensional representation is a direct sum of irreducible finite dimensional representations.
If \(G\) is a linearly reductive Lie group and \(X\) has a \(G\)-invariant Riemannian metric then there is a compact complete Klein manifold \(\Gamma\backslash X\) \cite{Borel:1963}. 
For more general \((X,G)\), it is not known when a compact complete Klein manifold modelled on \((X,G)\) exists.
\end{example}
\begin{example}
Again by Poincar\'e's\SubIndex{uniformization theorem}\SubIndex{theorem!uniformization} uniformization theorem \cite{Saint-Gervais:2010}, every connected surface with a Riemannian metric, after a conformal rescaling, has a complete metric of constant curvature.
It is thus a complete Klein manifold, for a geometry modelled on the real projective plane, Euclidean plane or hyperbolic plane, with its usual constant curvature metric.
\end{example}
\begin{example}
An \emph{elliptic curve}\define{elliptic curve} is a quotient \(E=\C/\Lambda\) where \(\Lambda\subseteq\C\) can be any lattice (i.e. any discrete abelian subgroup under addition), hence a homogeneous space.
Any biholomorphism between elliptic curves arises from a complex affine linear map \(\C\to\C\) matching lattices.
We can always arrange by isomorphism that \(\Lambda\) is generated by \(1,\tau\) for some complex number \(\tau\) with positive imaginary part.
But then \(\tau\) is only determined up to linear fractional transformation of the upper half plane \(\mathbb{CH}^1\) with integer coefficients.
Hence the moduli space of elliptic curves is the complete Klein manifold \(\Gamma\backslash\mathbb{CH}^1\) where \(\Gamma\) is the set of integer coefficient linear fractional transformations.
Similarly, the space of cubic surfaces is a complete Klein manifold \(\Gamma\backslash\mathbb{CH}^4\) where \(\mathbb{CH}^n\) is the complex hyperbolic space of complex dimension \(n\)  \cite{Allcock.Carlson.Toledo:2011}. (Note: \(\mathbb{CH}^n\) is equivariantly biholomorphic to the unit ball in \(\C^n\) \cite{Rudin:2008}.)
Morever, the space of cubic threefolds, after some blow up, is a complete Klein manifold \(\Gamma\backslash\mathbb{CH}^{10}\) \cite{Allcock.Carlson.Toledo:2011}.
\end{example}
\begin{problem}{aa}
Returning to the theory of 2\textsuperscript{nd} order complex analytic linear ordinary differential equations \vpageref{example:linear.ode}: which of the associated geometries is a complete Klein geometry?
\end{problem}
\section{Replacing the model by a connected model}
Suppose that \((X,G)\) is a homogeneous space and \(X_0\subseteq X\) is a component of \(X\).
Let \(G_0\subseteq G\) be the subgroup preserving \(X_0\).
Every \((X_0,G_0)\)-structure induces an \((X,G)\)-structure: (i) take the same charts, but (ii) add more charts by letting every element of \(G\) translate every chart, and (iii) add even more: add in every \(X\)-chart which locally arises in this way.
\begin{theorem}
Suppose that \((X,G)\)is a homogeneous space and \(X_0\subseteq X\) is a component of \(X\).
Let \(G_0\subseteq G\) be the subgroup preserving \(X_0\).
Then every \((X,G)\)-structure is induced by a unique \((X_0,G_0)\)-structure and vice versa.
\end{theorem}
\begin{proof}
Every \((X,G)\)-chart locally maps into a single component of \(X\).
Translate by an element of \(G\) to arrange that this component is \(X_0\).
\end{proof}
Note that every connected homogeneous space is strong, as we saw~\vpageref{section:analytic}.
So we can safely assume that our models are strong and connected.
\section{Developing maps}
To see how close a locally homogeneous structure is to a Klein geometry, consider the following theorem, which we will prove later (see theorem~\vref{thm:flat}):
\begin{theorem}\label{thm:flat.first.appearance}
Take a homogeneous space \((X,G)\).
Take an \((X,G)\)-structure on a connected manifold \(M\).
Its pullback to its universal covering space \(\tilde{M}\to M\) is also the pullback from a local diffeomorphism \(\tilde{M}\xrightarrow{\delta}X\):
\[
\begin{tikzcd}
\tilde{M}\arrow[r,"\delta"]\arrow[d]&X\\
M
\end{tikzcd}
\]
the \emph{developing map}.\define{developing map}

Suppose that \(G\) has kernel \(K\subseteq G\).
Let \(\bar{G}:=G/K\), so \((X,\bar{G})\) is the associated effective homogeneous space.
The developing map is equivariant for a unique group morphism \(\pi_1(M)\xrightarrow{h} \bar{G}\), the \emph{holonomy morphism}.\define{holonomy!morphism}
The pair \((\delta,h)\) of developing map and holonomy morphism are unique up to replacing by \((g\delta,\Ad_g h)\) for any \(g\in\bar{G}\).
\end{theorem}
By theorem~\vref{theorem:connected}, we can replace \(X\) by any of its connected components, so assume that \(X\) is connected and hence \((X,G)\) is strong and \((X,\bar{G})\) is strong effective.
\[
\begin{array}{p{6cm}}
\includegraphics[width=6cm]{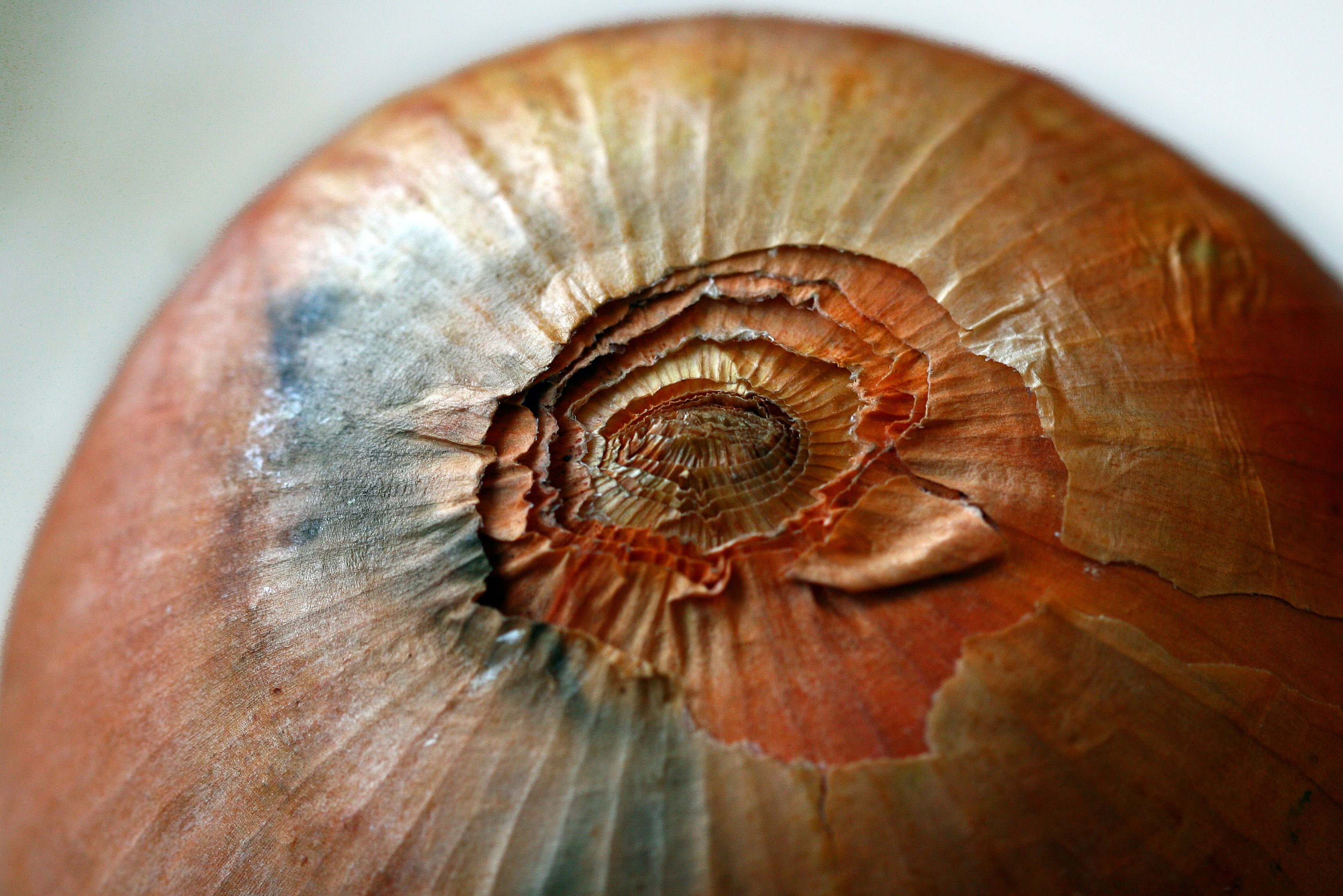}\\
\tiny{Muffet, CC BY 2.0, via Wikimedia Commons}\\ 
\tiny{\verb!<https://creativecommons.org/licenses/by/2.0>!}
\end{array}
\]
\begin{example}
If \(M\) is compact, connected and simply connected then \(M=\tilde{M}\) is compact, so \(\tilde{M}\to X\) is a local diffeomorphism from a compact manifold, so a covering map to its image.
If \(X\) is connected, clearly \(M=\tilde{X}\) is the universal covering space, and there is a unique \((X,G)\)-structure on \(M\) up to isomorphism.
For example, if \(X\) is the sphere, for any group \(G\) acting transitively on the sphere, the unique \((X,G)\)-structure on any compact, connected and simply connected manifold \(M\) is the one on \(X\).
\end{example}
\begin{problem}{ab}
Use this theorem to find all \((X,G)\)-structures on \(X\) where \(X=\RP{n}\) and \(G=\PGL{n+1,\R}\).
\end{problem}
\begin{example}
Let \((X,G)\) be the unit sphere \(X\) with \(G\) the isometry group of the usual round metric.
Let \(M\) be \(X\) punctured at a nonempty finite set of points.
The \emph{onion}\define{onion} \(\tilde{M}\) is the universal covering space of the punctured sphere, with pullback locally homogeneous geometry modelled on 
The developing map \(\tilde{M}\to X\) is just the lift of the inclusion map \(M\to X\).
Therefore the holonomy is trivial: \(\pi_1(M)\xrightarrow{1}\set{1}\subset G\).
\end{example}
\begin{example}
Recall how one proves the classification of curves, i.e. \(1\)-dimensional manifolds.
Assume, without loss of generality, that our curve \(M\) is connected.
By a partition of unity argument, it bears a Riemannian metric.
Because it is \(1\)-dimensional, the curvature vanishes.
Let \((X,G)\) be the real number line with its isometry group \(G\).
The usual proof now proceeds by constructing a developing map \(\tilde{M}\to X\), from the universal covering space \(\tilde{M}\) of the curve \(M\), to the line \(X\) with its standard flat metric.

The image \(\delta(\tilde{M})\subseteq X\) is a connected open set, hence an interval.
The holonomy morphism preserves that interval; it can't reverse the direction of the interval, since it acts without fixed points on that interval.
So the holonomy morphism acts by orientation preserving isometries, i.e. translations.
Throw out \(G\) and replace it with its subgroup of translations, which we now call \(G\), with the same developing map.
The unit speed translation vector field on \(X\) is holonomy invariant so defined on \(M\).
Maximal flow lines of any vector field are equal or disjoint as subsets.
The open sets in \(M\) (or in \(\tilde{M}\)) which are the flow lines of that vector field are thus equal or disjoint.
But \(M\) and \(\tilde{M}\) are both connected, so there is only one flow line.
Thinking of points of \(X\) as numbers, the developing map increases at unit speed along its flow, so is injective.
So the developing map takes \(\tilde{M}\) diffeomorphically to an interval.
If the holonomy is trivial, \(M\) is an interval.

If the holonomy is not trivial, \(\tilde{M}\) is an interval containing all translates of itself by some translation, so \(\tilde{M}=X\).
So \(M\) is the quotient by a discrete group of translations, so a circle.
\end{example}
\begin{example}\label{example:TranslationAndAffine}
Every elliptic curve $E=\C{}/\Lambda$  has a translation structure, i.e. an $(X,G)$-structure where 
$X=G=\C{}$.
The developing map is the identity map \(X\to X\) and the associated holonomy morphism is the inclusion map \(\Lambda\to\C\). 
A holomorphic translation structure is equivalent to a holomorphic nowhere vanishing \(1\)-form.
The automorphisms of any translation structure are the automorphisms of the \(1\)-form
\end{example}
\begin{example}\label{example:elliptic.curve.affine}
Every elliptic curve $E=\C{}/\Lambda$ has a  complex affine structure, i.e. an $(X,G)$-structure where $G=\C^{\times} \rtimes \C{}$ is the group of affine linear transformations \(z\mapsto az+b\), and \(X=\C\).
The developing map is the identity
\(
\mapto[\delta]{z \in \C{}}{z \in \C{}}
\)
and the holonomy morphism is
\[
\mapto[h]{\lambda \in \Lambda}{\left(1,\lambda\right) \in G}.
\]
There are more affine structures: if the elliptic curve is $E=\C{}/\Lambda$, then we can produce, for every $c \in \C^{\times}$, an affine structure on $E$ via the developing map 
\(
\mapto[\delta]{z \in \C{}}{e^{cz} \in \C^{\times}},
\)
with holonomy morphism
\[
\mapto[h]{\lambda \in \Lambda}{\left(e^{c\lambda},0\right) \in G}.
\]
Since the developing map of an affine structure is unique precisely up to affine transformations, we could just as well use the map
\[
\mapto{z}{\frac{e^{cz}-1}{c}},
\]
and think of \mapto{z}{z} as the limiting case $c\to 0$.
To forget the choice of affine coordinate $z$, we can identify the moduli space of affine structures
on an elliptic curve $E$ with  the space of holomorphic \(1\)-forms on $E$: each \(1\)-form is expressed
in terms of an affine chart as $c \, dz$. 
The developing map is a map to \((X',G'):=(\C^{\times},\C^{\times})\).
The holomorphic \(1\)-form \(d\zeta/\zeta\) on \(X'=\C^{\times}_{\zeta}\) is \(G'\)-invariant, and pulls back to \(c\,dz\).
Hence we compute the isomorphism class of the affine structure by  this induced holomorphic \(1\)-form, and take the zero \(1\)-form as the limit \(c\to 0\).
It is well known \cite{Loray/MarinPerez:2009} that every holomorphic affine structure arises uniquely in this way.
The automorphism group of the affine structure with \(c=0\), hence linear developing map, is the biholomorphism group of the elliptic curve $C$.
The automorphism group of any of the affine structures with \(c\ne 0\), i.e. nonlinear developing map,  is the group of translations.
\end{example}
\begin{example}
A \emph{complex projective structure}\define{complex projective structure}\define{projective!structure!complex} on a surface is an $(X,G)$-structure where $X=\Proj{1}$ and $G=\PSL{2,\C{}}$.
Every affine structure imposes a projective structure via the obvious inclusion $\C{}=\Proj{1} \setminus \infty \to \Proj{1}$.
Among the various affine structures on elliptic curves given above, each induces a distinct projective structure, except for pairs $c$ and $-c$, which are clearly related by a linear fractional
transformation 
\[
e^{-cz} = \frac{1}{e^{cz}}.
\]
It turns out that these are all of the projective structures: we can identify the moduli space of
projective structures on any elliptic curve with the space of quadratic
differentials $c^2 \, dz^2$, as we will see.
For $c \ne 0$, the developing map is
\(
\mapto[\delta]{z \in \C{}}{e^{cz} \in \Proj{1}}
\)
and the holonomy morphism is
\[
\mapto[h]{\lambda \in \Lambda}%
{
\begin{bmatrix} 
e^{c\lambda/2} & 0\\
0 & e^{-c\lambda/2}
\end{bmatrix}
\in G}.
\]
For $c=0$, the developing map is
\(
\mapto[\delta]{z \in \C{}}{z \in \Proj{1}}
\)
and the holonomy morphism is
\[
\mapto[h]{\lambda \in \Lambda}%
{
\begin{bmatrix} 
1 & \lambda \\
0 & 1
\end{bmatrix}
\in G}.
\]

The automorphism group of the projective structure at $c=0$ is the biholomorphism group of the elliptic curve.
At $c \ne 0$, the automorphism group is $\Z{}_2 \ltimes \C{}$, where $\Z{}_2$ acts on $\C{}/\Lambda$ by $z \mapsto -z$, and $\C{}$ acts by translations.
\end{example}
\begin{example}
Clearly $\Proj{1}$ has an obvious complex projective structure: as the model itself.
The developing map is
\(
\mapto[\delta]{z \in \Proj{1}}{z \in \Proj{1}}
\)
and the holonomy morphism is
\[
\mapto[h]{1}%
{
\begin{bmatrix} 
1 & 0\\
0 & 1
\end{bmatrix}
\in G}.
\]
The biholomorphism group is $G=\PSL{2,\C{}}$.
\end{example}
\begin{theorem}\label{theorem:morphism.match}
Suppose that \((X,G)\to(X',G')\) is a morphism of homogeneous spaces.
The map \(X\to X'\) is a covering map just when it is both onto and somewhere a local diffeomorphism, by equivariance.
Suppose that this occurs and that \(X\) is connected.
(Hence \(X\) and \(X'\) are strong and connected.)
Suppose that both homogeneous spaces are effective, and that \(G\to G'\) is onto and that every deck transformation of \(X\to X'\) arises from an element of the kernel of \(G\to G'\).
Then every \((X',G')\)-structure arises from a unique \((X,G)\)-structure and vice versa.
\end{theorem}
\begin{proof}
Pick points \(x_0\in X\) mapping to \(x_0'\in X'\).
Any \((X,G)\)-structure determines an \((X',G')\)-structure by composing each \(X\)-chart with \(X\to X'\).
On the other hand, the developing map of any \((X',G')\)-structure has a unique lift to a map to \(X\), a local diffeomorphism, lifting \(x_0'\) to \(x_0\).
The effect of any element of the fundamental group is then to move the map to \(X'\) by some element of \(G'\), and to move the lifted map to \(X\) by some diffeomorphism which lifts that element of \(G'\).
So we need that every diffeomorphism \(X\to X\) which lifts a diffeomorphism \(g'\in G'\) of \(X'\to X'\) is the action of some \(g\in G\).
Since \(G\to G'\) is onto, lift \(g'\in G'\) to some \(g\in G\).
Applying \(g^{-1}\), it matches the diffeomorphism \(X\to X\) up to covering the identity, i.e. up to a deck transformation.
There is an element of the kernel of \(G\to G'\) which carries out that deck transformation.
It is unique since \((X,G)\) is effective.
\end{proof}
\begin{example}\label{example:real.projective.lift}
Let \(X\) be the unit sphere of dimension \(n\): \(X=(\R^{n+1}-0)/\R^+\).
Let \(G\) be the group of linear isomorphisms of \(\R^{n+1}\) of determinant \(\pm 1\).
Let \(X':=\RP{n}\) be the real projective space of dimension \(n\), i.e. the set of all real lines through the origin on \(\R^{n+1}\), and \(G':=\PGL{n+1}\) be the set of all real projective transformations of \(X'\), i.e. \(G':=G/\set{\pm I}\).
Then \(X\to X'\) is a \(2-1\) covering map equivariant for \(G\to G'\), which is also a \(2-1\) covering map.
The deck transformations of \(X'\to X\) are \(\pm I\), the kernel of \(G\to G'\).
\end{example}
\section{Cohomological obstructions}
Suppose that \((X,G)\) is a homogeneous space and that \(\xi_X\) is a closed \(G\)-invariant differential form on \(X\).
If \(M\) has an \((X,G)\)-structure, each chart in the atlas identifies \(\xi_X\) with a closed differential form \(\xi_M\).
These closed forms agree on overlaps, by \(G\)-invariance.
So \(\xi_M\) is globally defined on \(M\).
If \(\xi_X\) has a power which is a volume form, then the same is true of \(\xi_M\), so if \(M\) is compact, \(\xi_M\) is not zero in cohomology.
Vanishing of the real cohomology of \(M\) in that degree ensures that \(M\) has no \((X,G)\)-structure.
\begin{example}\label{example:aff.R}
For example, take \((X,G)=(\CP{n},\PU{n+1})\), so that \(G\) preserves the Fubini--Study symplectic form on \(X\); hence every compact manifold \(M\) with an \((X,G)\)-structure is symplectic.
\end{example}
\begin{example}
The group \(G:=\Aff{\R}\) of affine transformations of the real number line has Maurer--Cartan form
\[
\omega=
\begin{pmatrix}
\alpha&\beta\\
0&0
\end{pmatrix}
\]
so that \(0=d\alpha=d\beta-\alpha\wedge\beta\), as the reader can check.
Take \(X:=G\), and look for an \((X,G)\)-structure on a compact oriented surface \(M\).
Take \(\alpha\wedge\beta\) as an area form, and then
\begin{align*}
0&<
\int_M\alpha\wedge\beta,
\\
&=
\int_M d\beta,
\\&=
\int_{\partial M}\beta,
\\
&=0,
\end{align*}
(If \(M\) is not orientable, pullback the structure to its oriented double cover.)
Hence there is no \((X,G)\)-structure on any compact surface.
\end{example}
\begin{problem}{ac}\label{exercise:univ.cover.replace}
Suppose that \((\tilde{X},\tilde{G})\) is the universal covering homogeneous space of a homogeneous space \((X,G)\). 
Prove that any \((X,G)\)-structure is induced by an \((\tilde{X},\tilde{G})\)-structure.
If \((X,G)\) is connected and effective, prove that \((\tilde{X},\tilde{G})\) is connected, strong and effective.
\end{problem}
\section{Local data}
How can we describe any locally homogeneous structure using only ``local data'' (tensor fields on \(M\) or on principal bundles over \(M\))?
It is convenient to think of a principal bundle over a manifold as having ``very small'' fibers, acted on by a ``very small'' group, so that the principal bundle looks like a fuzzy version of the manifold.
This picture fits with the picture of a principal bundle being glued together by transition maps, which are local data, and also fits with the idea of tangent vectors being infinitesimal motions.
\subsection{The bundle of a locally homogeneous structure}\label{subsection:G.over.H.bundle}
Take a homogeneous space \((X,G)\) and pick a point \(x_0\in X\) and let \(H:=G^{x_0}\).
The map
\[
g\in G\mapsto gx0\in X,
\]
a principal \(H\)-bundle.
It gives a bundle over every open set of \(X\): just the open set of \(G\) which maps to that open set of \(X\).
Each \(X\)-chart identifies an open set in \(M\) with one in \(X\), 
\[
\input{xg-chart}
\]
so determines a bundle over that open set in \(M\): just the open set in \(G\):
\[
\input{xg-bundle}
\]
Take a \(G\)-compatible chart.
\[
\input{xg-struc}
\]
On the overlap of the charts, the two charts are matched up locally by some element \(g\in G\), i.e. by the left action, which matches up those open sets in \(G\) and the Maurer--Cartan forms on those open sets in \(G\).
If \((X,G)\) is strong and effective, \(g\) is unique.
Glue these open sets in \(G\) together over the overlap.
Put together all of the charts: you get a principal \(H\)-bundle \(\G\to M\), locally made from the one over \(X\) in each chart:
\[
\input{xg-bundle-2}
\]
Each gluing is by the left action of elements of \(G\).
The Maurer--Cartan form \(\omega_G\) is left invariant.
So we glue together open sets in \(G\) with Maurer--Cartan forms matching.
So this bundle \(\G\to M\) has a \(1\)-form \(\omega_{\G}\) on it, the \emph{Cartan connection},\define{Cartan!connection}\define{connection!Cartan} locally identified by charts with \(\omega_G\).
The left invariant vector fields on \(G\) get copied onto the bundle \(\G\) when we glue it together; these are the \emph{constant vector fields}\define{constant vector field} of the locally homogeneous structure, as \(\omega_{\G}\) is constant on each of them.
\begin{theorem}\label{theorem:XG.geom.to.struc}
Given a locally homogeneous structure on a manifold \(M\), with strong effective model, the bundle \(\G\to M\) with its Cartan connection \(\omega_{\G}\) determines the locally homogeneous structure on \(M\) uniquely up to isomorphism.
\end{theorem}
This is an obvious corollary of theorem~\vref{thm:flat.first.appearance}.
Theorem~\ref{thm:flat.first.appearance} follows from theorem~\vref{thm:flat}, which itself follows from theorem~\vref{theorem:integrate.Maurer}, but we give a proof here of theorem~\ref{thm:flat.first.appearance} as a preview of the proofs of those theorems.
On a first reading, the reading might skip this long proof.
\begin{proof}
We will provide three independent proofs.
Our first proof relies explicitly on homotopies of paths.
Denote the bundle \(\G\to M\).
Every point of \(M\) lies in an open set \(U\subseteq M\) over which the preimage \(\G_U\subseteq\G\) is identified with an \(H\)-invariant open set \(G_U\subseteq G\), uniquely up to left translation.
Cover \(M\) in such open sets, say \(M_a\subseteq M\), and pick developing maps \(\G_a:=\G_{M_a}\xrightarrow{\delta_a}G\).
We can assume, by perhaps refining our open cover, that each \(M_a\) is connected and that each overlap \(M_{ab}:=M_a\cap M_b\) is connected.
By strong effectivity of \((X,G)\), a unique element \(g_{ab}\) over each overlap \(M_{ab}\) left translates one developing map into the other:
\[
\delta_b=g_{ab}\delta_a.
\]
Again by strong effectivity, on any triple overlap \(M_{abc}:=M_{ab}\cap M_c\), 
\[
g_{ca}g_{bc}g_{ab}=1.
\]
We can alter all of the developing maps at once by one left translation \(g_a\) above each open set \(M_a\).

Fix a point \(m_0\in M\), say in some open set \(M_{a_0}\) and pick a developing map 
\[
\G_{a_0}\xrightarrow{\delta_{a_0}}G.
\]
Take a continuous path \(m(t)\in M\), \(t\in I:=[a,b]\) with \(m(a)=m_0\) and let \(m_1:=m(b)\).
Cover the interval \(I\) by the inverse images \(m^{-1}M_a\) of the various \(M_a\).
Each point of \(I\) lies a connected open set lying in one such preimage.
By compactness of \(I\), we can take a finite subcover.
Partition our interval into finitely many closed subintervals \(I_j\) on which the path lies in a single \(M_{a_j}\), with each interval ending where the next begins, say at a point \(t_j\).
(One open set \(M_a\) might appear several times in this sequence.)
We have picked a developing map \(\delta_{a_0}\) on the first interval.
For each interval \(I_j\), there is a unique \(g_{a_j}\) to arrange that \(\delta_{a_j}=g_{a_j}\delta_{a_{j-1}}\) over an open set in \(M\) containing \(m(t_j)\).
Hence we can arrange a sequence \(\delta_j:=g_{a_j}^{-1}\delta_{a_j}\) of developing maps of open sets \(M_{a_j}\), each containing the image of \(I_j\), each agreeing with the next along our path.

Refine the open cover: pick smaller open sets inside these.
By the same argument, at each step, the patching is uniquely determined by having each developing map agree on some overlap with the previous one.
So the patching together agrees with the previous patching together step by step, i.e. the result is the same for all sufficiently fine open covers.

Patching, we have defined a map
\[
I\times\G\supseteq\text{open}\xrightarrow{\delta}G
\]
defined on an \(H\)-invariant open set of \(I\times\G\),  containing the points \((t,p)\) with \(\pi(p)=m(t)\).
Our map \(\delta\) is locally constant in \(t\).
For each fixed \(t\), \(\delta\) maps an open set in \(\G\) to \(G\) over an open subset of \(M\) containing \(m(t)\).
This map depends only on the path \(m(t)\) and the \((X,G)\)-structure, not on the open cover.

Take a homotopy \(m_s(t)\) of our path \(m(t)\), keeping the same end points.
For each value \(s=s_0\), cover the path \(m_{s_0}(t)\) by open sets as above, constructing some map \(\delta_{s_0}(t,p)\).
For small changes in \(s\), the path \(m_s(t)\) stays in one of these open covers, and therefore \(\delta_s(t,p)\) stays the same, i.e. is locally constant in \(s\).
So \(\delta_s(b,p)\) is constant in \(s\) for any \(p\in\G_{m_1}\).
In other words, \(\delta(b,p)\), as computed along some path \(m(t)\), for any \(p\in\G_{m_1}\), depends only on the homotopy class of the path \(m(t)\).

Each point \(\tilde{m}_1\in\tilde{M}\) of the universal covering space \(\tilde{M}\xrightarrow{\pi}M\) is the end point \(\tilde{m}_1=\tilde{m}(b)\) of a path \(\tilde{m}(t)\) lifted from \(M\), starting at the base point \(\tilde{m}(a)=\tilde{m}_0\), uniquely modulo homotopy.
So the value \(\delta(b,p)\) as computed along the path \(m(t)\) is defined by the choice of the point \(\tilde{m}_1\in\tilde{M }\) and of the point \(p\in\G_{m_1}\) over the image point \(m_1:=\pi(\tilde{m}_1)\in M\).
In other words, \(\delta\) is defined precisely on the pullback of \(\G\to M\) to the universal covering space \(\tilde{M}\to M\).

Each point of \(\tilde{M}\) is the end point of a path lifted from \(M\), starting at \(m_0\), uniquely modulo homotopy.
Each deck transformation of \(\tilde{M}\to M\) is the composition of such a path with a loop from \(m_0\) to \(m_0\), uniquely up to homotopy.
The deck transformation takes our map \(\delta\) to another such \(\delta\), giving by construction such a \(\delta\) for the loop.
By uniqueness of \(\delta\) up to left translation, the deck transformations act by elements of \(G\), a holonomy morphism.
The theorem is proved.

Our second proof uses developing maps of Maurer--Cartan forms.
Write the components of \(\G\) as \(\G_i\).
Pick a point \(m_0\in M\).
Since \(M\) is connected, each component of \(\G\) has a point mapping to \(m_0\).
Pick one point \(p_i\in\G_i\) mapping to \(m_0\) from each component \(\G_i\).
In particular, pick one such point and denote it \(p_0\in\G_0\).
Each \(p_i\) maps to \(m_0\), so all lie in the same fiber of the principal bundle \(H\to\G\to M\).
Hence every \(p_i\) is uniquely written as \(p_i=p_0h_i\) for uniquely \(h_i\in H\).
Let \(\tilde\G_i\to\G_i\) be the universal covering map with a chosen point \(\tilde{p}_i\in\tilde\G_i\) mapping to \(p_i\in G_i\).
By theorem~\vref{theorem:recover.Lie.group}, there is a developing map for the Maurer--Cartan form:
\[
\begin{tikzcd}
\tilde\G_i\arrow[d]\arrow[r]&G\\
\G_i,
\end{tikzcd}
\]
and a holonomy morphism \(\fundamentalGroup{\G_i}\xrightarrow{h}G\).
These developing maps are determined up to left translation.
So there is a unique choice to get \(\tilde{p}_i\mapsto h_i\).
Let \(\tilde\G\) be the disjoint union of the \(\tilde\G_i\), giving a map
\[
\begin{tikzcd}
\tilde\G\arrow[d]\arrow[r,"\tilde\delta"]&G\\
\G.
\end{tikzcd}
\]

By hypothesis, \(M\) is covered by open sets \(M_a\subseteq M\) over which the preimage \(\G_a\subseteq\G\) admits a developing map \(\G_a\xrightarrow{\delta_a}G\) for the Maurer--Cartan form.
Take the preimage \(\tilde\G_a\subseteq\tilde\G\).
The composition
\[
\tilde\G_a\xrightarrow{\pi}\G_a\xrightarrow{\delta_a} G
\]
is a developing map for the Maurer--Cartan form, so agrees with \(\tilde\delta\) up to \(G\)-action on each component.
In particular, if \(m_0\in M_a\) then we can pick \(\delta_a(p_0)=1\) and then \(\delta_a(p_i)=h_i\) by right \(H\)- equivariance.
So then  the composition has \(\delta_a\circ\pi=\tilde\delta\) on the points \(p_i\) and therefore on \(\tilde\G_a\).

Suppose that we have arranged \(\delta_b\circ\pi=\tilde\delta\) on some \(\tilde\G_b\) and that \(M_{bc}\) is not empty.
Compose \(\delta_c\) with a left translation to arrange \(\delta_c\circ\pi=\tilde\delta\) over \(M_{bc}\).
By uniqueness of developing maps, \(\delta_c\circ\pi=\tilde\delta\) over \(M_c\).
Repeating this process, we have picked out the needed \(g_a\) so we can arrange \(\delta_a\circ\pi=\tilde\delta\) for all \(M_a\).

Since \(\tilde\delta\) descends on these open sets \(M_a\) to be defined on \(H\)-bundles, \(\tilde\delta\) doesn't change value under the deck transformation of any loop in \(H\).
As we move to any overlapping open set \(M_b\), this invariance holds on the overlap \(\tilde\G_a\cap\tilde\G_b\), so on all of \(\tilde\G_b\) by uniqueness of the developing map and holonomy morphism.

Since \((X,G)\) is strongly effective, \(h\) is uniquely defined by its interaction with \(\tilde\delta\), so \(h\) is trivial on \(\fundamentalGroup{H}\), so descends to a morphism \(\fundamentalGroup{M}\xrightarrow{h}G\) by the homotopy sequence
\[
\begin{tikzcd}
\dots\arrow[r]&\fundamentalGroup{H}\arrow[r]&\fundamentalGroup{\G}\arrow[r]&\fundamentalGroup{M}\arrow[r]&\dots
\end{tikzcd}
\]
By \(h\)-equivariance, \(\tilde\delta\) is defined on the pullback of \(\G\) to \(\tilde{M}\); call it \(\delta\).
So \(\delta\) agrees with each \(\delta_a\) up to left translation hence is a developing map of the locally homogeneous structure.
Moreover, \(\delta\) is \(H\)-equivariant on that pullback bundle, since the \(H\)-action descends over each of these open sets \(M_a\) to make the developing map \(H\)-equivariant.
The same argument gives the holonomy morphism.

Our third proof is in the spirit of Cartan geometries.
On the manifold \(\G\times G\), the set \(V\) of tangent vectors to \(\G\times G\) satisfying \(\omega_{\G}=\omega_G\) is a subbundle of the tangent bundle, invariant under diagonal right \(H\)-action by lemma~\vref{lemma:bracket}, and invariant under left \(G\)-action on \(G\), because \(\omega_G\) is left invariant. 
It satisfies the hypothesis of the Frobenius theorem, since both \(\omega=\omega_{\G}\) and \(\omega=\omega_G\) satisfy
\[
0=d\omega+\frac{1}{2}\lb{\omega}{\omega}.
\]
So \(\G\times G\) is foliated by leaves, i.e. maximal connected submanifolds whose tangent space at each point \((p,g)\) is the linear subspace \(V_{(p,g)}\), so of dimension equal to the dimension of \(G\).

An \emph{\(H\)-folium} is a minimal \(H\)-invariant nonempty set which contains every \(V\)-leaf through any of its points.
Since the \(H\)-orbits are tangent to \(V\), and \(V\) is \(H\)-invariant, by the equivariant Frobenius theorem (corollary~\vref{corollary:equivariant.Frobenius}), \(\G\times G\) is foliated by \(H\)-folios of dimension equal to that of \(G\).
Moreover, the \(H\)-folios are the preimages of the leaves of a unique foliation on \(\amal{\G}{H}{G}\).

Let \(\G'\) be one such \(H\)-folio.
Let \(\omega\) be the pullback of \(\omega_{\G}\) to \(\G'\); because \(\G'\) is \(V\)-tangent, \(\omega_{\G}=\omega_G\) on \(\G'\), i.e. the pullbacks agree.
The constant vector fields of \(G\) are the left invariant vector fields, so complete.
The constant vector fields of \(\G'\), i.e. vector fields on which \(\omega\) is constant, are precisely the vector fields \(A_{\G'}=(A_{\G},A_G)\), which span \(V\).
The flow of each constant vector field \(A_{\G'}\) through any point of \(\G'\) is therefore defined for precisely as long as that of \(A_{\G}\) is defined on the image point in \(\G\). 
By theorem~\vref{cor:orbitMapEquiv} the map \(\G'\to\G\) is a fiber bundle mapping as is the quotient \(M':=\G'/H\to M\).
Since each is a local diffeomorphism, each fiber bundle map is a covering map.
The quotient \(M'\) is a leaf of the foliation of \(\amal{\G}{H}{G}\).
Lift the bundle \(\G'\to M'\) to the universal covering space \(\tilde{M}\):
\[
\begin{tikzcd}
\G\arrow[d]&\tilde\G\arrow[d]\arrow[r]\arrow[l]&G\arrow[d]\\
M&\tilde{M}\arrow[l]\arrow[r]&X.
\end{tikzcd}
%\begin{tikzcd}
%&\tilde\G\arrow[dl]\arrow[dr]\arrow[rrr,bend left]&&&\tilde{M}\arrow[dl]\arrow[dr]&\\
%\G\arrow[rrr,bend right]&&G\arrow[rrr,bend right,crossing over]&M&&X.
%\end{tikzcd}
\]
The deck transformations of \(\tilde{M}\) act as bundle automorphisms on \(\tilde\G\), since the bundle is pulled back from \(M\).
These transformations preserve \(\omega\) and the \(H\)-action, as these are also pulled back.
But the map to \(G\) might not be invariant under the deck transformations.
Locally, each deck transformation \(\gamma\in\pi_1(M)\) acts by an automorphism, so some element \(g=h(\gamma)\in G\).
By connectivity of \(M\), this element is unique.

The construction is unique up to the choice of leaf \(\G'\).
If we pick the leaf through some point \((p_0,g_0)\), changing our choice to \((p_0,gg_0)\) alters the developing map and holonomy morphism as described.
But every leaf maps to \(p_0\), since it covers \(\G\).
The leaves in \(\amal{\G}{H}{G}\) have preimages in \(\G\times G\) which are \(H\)-invariant, so \(H\)-bundles, on which \(\omega_{\G}=\omega_G\).
\end{proof}
\begin{corollary}
Given a locally homogeneous structure on a manifold \(M\), with strong effective model,
the bundle \(\G\to M\) with the association of a constant vector field on \(\G\) to each element of the Lie algebra \(\LieG\) of \(G\) determines the locally homogeneous structure on \(M\) uniquely up to isomorphism.
\end{corollary}
\begin{theorem}
Given a locally homogeneous structure, with strong effective connected and simply connected model, its constant vector fields are complete if and only if the locally homogeneous structure is a complete Klein geometry.
\end{theorem}
This follows from theorem~\vref{thm:develop.submanifold}; also see corollary~\vref{cor:vast.flat}.
\begin{example}
Take a Lie group \(G\) and let \(X:=G\), with \(G\) acting on \(X\) by left translation.
Then an \((X,G)\)-structure on a manifold \(M\) determines a differential form \(\omega\) on \(M\) which is locally identified with the  Maurer--Cartan form of \(G\), by theorem~\vref{theorem:open.subset.G}.
By theorem~\vref{theorem:recover.Lie.group},  
\end{example}

\chapter{Submanifolds of homogeneous spaces}
\section{Subgeometries of immersions}
We can easily generalize the concept of locally homogeneous structure: we will capture the local geometry of an immersed submanifold of a homogeneous space, expressed purely as local data on the submanifold.
Take a homogeneous space \((X,G)\).
Take a point \(x_0\in X\) and let \(H:=G^{x_0}\).
Take an immersion \(M\xrightarrow{\delta} X\) from a manifold \(M\).
Define a principal \(H\)-bundle \(\G_M:=\delta^*G\):
\[
\begin{tikzcd}
H\arrow{r}&\G_M\arrow{d}\arrow[r,"\Delta"]&G\arrow[d]\\
&M\arrow[r,"\delta"]&X.
\end{tikzcd}
\]
The points of \(\G_M\) are pairs \((m,g)\) so that \(m\in M\) and \(g\in G\) and \(gx_0=\delta(m)\), and \(\Delta\) is
\[
(m,g)\in\G_M\xmapsto{\Delta}g\in G.
\]
Let \(\omega_M:=\Delta^*\omega_G\) be the pullback of the Maurer--Cartan form.
This bundle \(\G_M\to M\) with this form \(\omega_M\) is the \emph{subgeometry} associated to the immersion.

The subgeometry of a composition of immersions
\[
\begin{tikzcd}
M_0\arrow[r]& M_1\arrow[r]&X
\end{tikzcd}
\]
is thus the pullback bundle with the pullback form, the pullback subgeometry.
So a subgeometry is a local geometric structure on a manifold.
The subgeometry of an immersion is invariant under composition of the immersion with left translation by any constant element \(g\in G\):
\[
\begin{tikzcd}
H\arrow{r}&\G_M\arrow{d}\arrow[r,"g\Delta"]&G\arrow[d]\\
&M\arrow[r,"g\delta"]&X.
\end{tikzcd}
\]

\begin{problem}{local.picture:properties}
Prove that \(\omega:=\omega_M\) is injective on tangent spaces of \(\G_M\) and that 
\[
0=d\omega+\frac{1}{2}\lb{\omega}{\omega}.
\]
\end{problem}

\section{Subgeometries in general}
Once again, our goal is to capture the local geometry of an immersed submanifold of a homogeneous space.
We will soon prove that the subgeometry defined above is precisely that local geometry.
But to make sense of this approach, 
\begin{itemize}
\item
we will define subgeometries abstractly, as geometric structures on manifolds
\item
then we will find differential invariants which determine whether one of these abstract subgeometries arises locally as above from an immersion in a homogeneous space
\item
and finally, we will identify the topological obstruction to turn ``locally'' into ``globally''.
\end{itemize}
\subsection{Definition}
Abstractly, an \((X,G)\)-\emph{subgeometry} on a manifold \(M\) is a principal right \(H\)-bundle
\[
\begin{tikzcd}
H\arrow{r}&\G_M\arrow{d}\\
&M
\end{tikzcd}
\]
and a \(\LieG\)-valued \(1\)-form \(\omega=\omega_M\) on \(\G_M\)
\begin{itemize}
\item
\(H\)-equivariant, i.e. transforming in the adjoint representation under right \(H\)-action and
\item
 injective on tangent spaces of \(\G_M\) and 
 \item
 with \(A_{\G_M}\hook\omega=A\) for all \(A\in\LieH\).
 \end{itemize}
\begin{example}
Clearly the subgeometry of an immersion to a homogeneous space is a subgeometry in this sense.
\end{example}
\subsection{Morphisms}
The \emph{pullback} of a subgeometry \(\G_{M'}\to M'\) via an immersion \(M\to M'\) is the pullback bundle with the pullback differential form.
A \emph{morphism} of \((X,G)\)-subgeometries \(\G_M\to M\), \(\G_{M'}\to M'\) is an \(H\)-equivariant immersion of the bundles matching up the forms.
Clearly it induces an immersion \(M\to M'\).
An \emph{isomorphism} is a morphism which is also a diffeomorphism, hence has an inverse morphism.
\subsection{Curvature}
The \emph{curvature} of a subgeometry is
\[
\Omega:=d\omega+\frac{1}{2}\lb{\omega}{\omega}.
\]
The subgeometry is \emph{flat} if \(\Omega=0\).
We will prove soon that any subgeometry is flat just when it is locally the subgeometry of an immersion.
In the process, we characterize flat subgeometries globally, not just locally (see theorem~\vref{theorem:integrate.Maurer}).

\section{Equivalent immersions}
\begin{lemma}\label{lemma:picture.uniqueness}
Take two immersions
\[
M_0\xrightarrow{\delta_0}X\xleftarrow{\delta_1}M_1
\]
of manifolds to the same homogeneous space \((X,G)\).
Suppose that \(M_0\) is connected.
Denote the pullback bundles as \(\G_0:=\delta_0^*G\), \(\G_1:=\delta_1^*G\), and the associated bundle maps as
\[
\begin{tikzcd}
\G_0\arrow[r,"\Delta_0"]&G&\G_1\arrow[l,"\Delta_1"'].
\end{tikzcd}
\]
There is a morphism
\[
\begin{tikzcd}
\G_0\arrow[d]\arrow[r]&\G_1\arrow[d]\\
M_0\arrow[r]&M_1
\end{tikzcd}
\]
precisely when there is an element \(g\in G\) extending the diagram to:
\[
\begin{tikzcd}
&G
\arrow[dd]
&\\
\G_0\arrow[dd]\arrow[rr,crossing over]\arrow[ur,"g\Delta_0"]
&&\G_1\arrow[dd]\arrow[ul,"\Delta_1"']
\\
&X&\\
M_0
\arrow[ur,"g\delta_0"]\arrow[rr]
&&M_1
\arrow[ul,"\delta_1"']
\end{tikzcd}
\]
This element \(g\in G\) is then unique.
\end{lemma}
\begin{proof}
The result only depends on the pullbacks to \(M_0\), so we can assume that \(M_1=M_0\), which we denote by \(M\).
Identify the bundles, so both are some \(\G:=\G_M\to M\), so the bundle maps are two \(H\)-equivariant maps
\[
\G\xrightarrow{\Delta_0,\Delta_1}G
\]
with the same \(\omega_M\).
Take a point \(p_0\in\G\).
After left translation by a unique element of \(G\), arrange that \(\Delta_0(p_0)=\Delta_1(p_0)\).
Take any locally absolutely continuous path \(p(t)\in\G\) with \(p(0)=p_0\).
Then the paths \(\Delta_0(p(t)), \Delta_1(p(t))\in G\) have the same Darboux derivatives:
\begin{align*}
\frac{d}{dt}\Delta_i(p(t))\hook\omega_G
&=
(\Delta_i'(p(t))p'(t))\hook\omega_G,
\\
&=
p'(t)\hook\Delta_i^*\omega_G,
\\
&=
p'(t)\hook\omega_M
\end{align*}
and start at the same point, so are the same path by the uniqueness part of theorem~\vref{theorem:Lie.equations}.
Hence \(\Delta_0=\Delta_1\) on the path component of \(p_0\) in \(\G\).
By \(H\)-equivariance, \(\Delta_0=\Delta_1\) above a path component of \(M\).
Because \(M\) is connected, \(\Delta_0=\Delta_1\).
Quotient by \(H\)-action to find \(\delta_0=\delta_1\).
\end{proof}

\section{Multivalued immersions and developing pairs}
\subsection{Covering immersions}
Take a connected manifold \(M\) and let \(\pi:=\fundamentalGroup{M}\) and let \(\tilde{M}\to M\) be the universal covering space.
A flat \((X,G)\)-subgeometry on \(M\), in this generality, might arise not from an immersion \(M\to X\) but instead perhaps from a \emph{covering immersion},\define{covering immersion} i.e. from an immersion \(\tilde{M}\to X\) of the universal covering space \(\tilde{M}\to M\), but with a lift of the action of \(\fundamentalGroup{M}\) to an action by automorphisms of the associated flat subgeometry.
\subsection{Developing pairs}
We will need a more concrete description of covering immersions.
Take a homogeneous space \((X,G)\), pick a point \(x_0\in X\) and let \(H:=G^{x_0}\).
An \((X,G)\)-\emph{developing pair}\define{developing!pair} \((\delta,h)\) on a manifold \(M\) is an immersion 
\[
\begin{tikzcd}
\tilde{M}\arrow[d]\arrow[r,"\delta"]&X\\
M,
\end{tikzcd}
\]
called the \emph{developing map}\define{developing!map}\define{map!developing} and a group morphism \(\pi\xrightarrow{h}G\), called the \emph{holonomy morphism},  so that 
\[
\delta\circ\gamma=h(\gamma)\delta
\]
for all \(\gamma\in\pi:=\fundamentalGroup{M}\).
\begin{lemma}
Every covering immersion has a developing pair, uniquely defined, and conversely every developing pair arises from a unique covering immersion.
\end{lemma}
Henceforth, we forget about covering immersions, and make frequent use of developing pairs.
\begin{proof}
Start with a developing pair \((\delta,h)\).
Since 
\[
\delta\circ\gamma=h(\gamma)\delta
\]
we see that
\[
gx_0=\delta(\tilde{m}) \text{ if and only if } h(\gamma)gx_0=\delta(\gamma\tilde{m}).
\]
Let \(\pi:=\fundamentalGroup{M}\).
Define an action of \(\pi\) on \(\tilde{M}\times G\) by the product action
\[
\gamma(\tilde{m},g):=(\gamma\tilde{m},h(\gamma)g)
\]
for \(\gamma\in\pi\).
This \(\pi\)-action leaves \(\G_{\tilde{M}}=\delta^*G\subseteq\tilde{M}\times G\) invariant, since \((\tilde{m},g)\in\G_{\tilde{M}}\) just when 
\[
\delta(\tilde{m})=gx_0,
\]
and this occurs just when
\[
\delta(\gamma\tilde{m})=h(\gamma)\delta(\tilde{m})=h(\gamma)gx_0,
\]
i.e. just when 
\[
\gamma(\tilde{m},g)\in\G_{\tilde{M}}.
\]
This \(\pi\)-action commutes with the \(H\)-action
\[
(\tilde{m},g)h=(\tilde{m},gh).
\]
It acts preserving \(\omega_{\tilde{M}}\), as \(\omega_{\tilde{M}}\) is pulled back from \(G\), and is invariant under left translation by \(h(\gamma)\).
Hence the immersion \(\delta\) is equivariant for this action, a covering immersion.

Conversely, take a covering immersion, i.e. an immersion 
\[
\begin{tikzcd}
\tilde{M}\arrow[d]\arrow[r,"\delta"]&X\\
M,
\end{tikzcd}
\]
which is equivariant for an action of \(\pi:=\fundamentalGroup{M}\) by bundle automorphisms of \(\G_{\tilde{M}}\) preserving \(\omega:=\omega_{\tilde{M}}\).
Denote the bundle map by
\[
\begin{tikzcd}
\G_{\tilde{M}}\arrow[d]\arrow[r,"\Delta"]&G\arrow[d]\\
\tilde{M}\arrow[d]\arrow[r,"\delta"]&X\\
M,
\end{tikzcd}
\]
i.e. \(\Delta(\tilde{m},g):=g\).
Pick some \(\gamma\in\pi\).
Since \(\gamma\) acts by automorphisms of the subgeometry, by lemma~\vref{lemma:picture.uniqueness},  \(\Delta\circ\gamma=h(\gamma)\Delta\) for a unique element \(h(\gamma)\in G\).
By existence and uniqueness of \(h(\gamma)\) for each \(\gamma\in\pi\), \(\pi\xrightarrow{h}G\) is a group morphism.
\end{proof}
Given a developing pair, the action of \(\pi\) on \(\G_{\tilde{M}}\) is free and proper, with quotient subgeometry \(\G_M:=\pi\backslash\G_{\tilde{M}}\to M\), the \emph{associated} flat subgeometry of the developing pair.
\begin{problem}{developing.pair:kernel.zero}
If \(M\) and \(X\) have the same dimension and \((X,G)\) is strong effective, prove that the developing map of any developing pair determines the holonomy morphism.
\end{problem}
\subsection{Action on developing pairs}
Take a homogeneous space \((X,G)\) and a manifold \(M\).
The group \(G\) acts on the space of \((X,G)\)-developing pairs on \(M\) by
\[
g(\delta,h):=(g\delta,Ad_g h).
\]

\section{Integration}
We want to integrate the local data of a subgeometry to the global data of a developing pair.
\begin{theorem}\label{theorem:integrate.Maurer}
Take a homogeneous space \((X,G)\) and a connected manifold \(M\).
Any \((X,G)\)-subgeometry on \(M\) is flat just when it is isomorphic to the \((X,G)\)-subgeometry of a developing pair.
Two \((X,G)\)-developing pairs on \(M\) have isomorphic associated subgeometries just when they agree up to \(G\)-action.
\end{theorem}
Intuitively, this theorem says that the moduli space of flat subgeometries is the set of developing pairs (i.e. of covering immersions) modulo \(G\)-action.
\begin{proof}
We start with an \((X,G)\)-subgeometry \(\G_M\to M\).
Let \(\G:=\G_M\).
Let \(Z:=\G\times G\).
Denote points of \(Z=\G\times G\) as \((p,g)\) for \(p\in\G\) and \(g\in G\).
Get \(G\) to act on \(Z\) by
\[
g'(p,g):=(p,g'g).
\]
Get \(H\times H\) to act on \(Z\) on the right by 
\[
(p,g),(h,h'):=(ph,gh').
\]
Clearly the quotient is \(M\times X\).
Get \(H\) to act on on \(Z\) on right by the diagonal action
\[
(p,g)h:=(ph,gh).
\]
The action is free and proper, so \(\bar{Z}:=Z/H=\amal{\G}{H}{G}\) is a smooth manifold and the quotient map \(Z\to\bar{Z}\) is a principal right \(H\)-bundle.
The map \((p,g)\in Z\mapsto (m,gx_0)\in M\times X\) is invariant under the right \(H\)-action \((p,g)h=(ph,gh)\), so descends to a smooth map \(\bar{Z}\to M\times X\).
Since \(Z\to M\times X\) is a principal right \(H\times H\)-bundle, \(\bar{Z}\to M\times X\) is a fiber bundle with fibers diffeomorphic to \(H\), perhaps not a principal bundle.

Let \(V_Z\subseteq TZ\) be the set of tangent vectors \(v\in T_z Z\) on which \(\omega_M=\omega_G\).
Pick a point \(z=(p,g)\in Z=\G\times G\).
Take any vector \(u\in T_p\G\).
Let \(A:=u\hook\omega_M\) and let \(w:=\LT{g*}A\in T_g G\).
Clearly \(v:=(u,w)\in V_Z\), and conversely every vector \(v\in V_Z\) arises uniquely as a such pair \(v=(u,w)\).
Hence the map
\[
(g,u)\in G\times T\G\mapsto v=(u,w)\in V_Z
\]
is a diffeomorphism, and a linear isomorphism on fibers \(\set{g}\times T_p\G\mapsto V_{Z,(p,g)}\).
So \(V_Z\cong G\times T\G\) is a vector subbundle of \(TZ\).

The tangent spaces of the \(H\)-orbits in \(Z\) are precisely the span of the vectors \(v=(u,w)\) for \(u\hook\omega_M=w\hook\omega_G\in\LieH\).
The \(H\)-orbits in \(Z\) are \(V_Z\)-tangent and the \(\LieH\)-action is by \(V_Z\)-tangent vector fields.

Under \(H\)-action on \(Z\), 
\[
r_h^*(\omega_M-\omega_G)=\Ad_h^{-1}\omega_M-\Ad_h^{-1}\omega_G=\Ad_h^{-1}(\omega_M-\omega_G).
\]
So \(V_Z\) is \(H\)-invariant.

For \(\omega=\omega_M\) or \(\omega=\omega_G\), we have
\[
0=d\omega+\frac{1}{2}\lb{\omega}{\omega}.
\]
So the difference \(\xi:=\omega_M-\omega_G\) satisfies \(d\xi=0\) on \(V_Z\).
Hence \(V_Z\subseteq TZ\) is bracket closed.

By the Frobenius theorem~\vref{theorem:Frobenius}, \(V_Z\) is the tangent bundle \(V_Z=T\mathscr{F}_Z\) of a unique foliation \(\mathscr{F}_Z\) on \(Z\).
Since \(V_Z\) is \(H\)-invariant, so is the foliation \(\mathscr{F}_Z\).

Recall that the \emph{\(H\)-folios} of \(\mathscr{F}_Z\) are the minimal \(H\)-invariant unions of \(\mathscr{F}_Z\)-leaves in \(Z\).
Note that \(H\) acts freely and properly on \(Z\).
By the equivariant Frobenius theorem (corollary~\vref{corollary:equivariant.Frobenius}), \(V_Z\) projects to a subbundle \(V_{\bar{Z}}\subseteq T\bar{Z}\), also bracket closed.
Hence \(V_{\bar{Z}}\) is the tangent bundle of a unique foliation \(\mathscr{F}_{\bar{Z}}\) on \(\bar{Z}\).
Moreover, the equivariant Frobenius theorem tells us further that each leaf \(M'\subseteq\bar{Z}\) is precisely the \(H\)-quotient of a unique \(H\)-folio \(\G'\subset Z\) of \(\mathscr{F}_Z\), and \(\G'\subseteq Z\) is the preimage of \(M'\subseteq\bar{Z}\).
Since \(H\) acts freely and properly on \(Z\), preserving each \(H\)-folio, each \(H\)-folio \(\G'\subseteq Z\) is a principal right \(H\)-bundle over the leaf \(M':=\G'/H\subset\bar{Z}\).

Let \(M'\subseteq\bar{Z}\) be one such leaf and take the pullback bundle
\[
\begin{tikzcd}
\G'\arrow{r}\arrow{d}&Z\arrow{d}\\
M'\arrow{r}&\bar{Z}.
\end{tikzcd}
\]
So \(\G'\) is the \(H\)-folio in \(\G\times G\) which is the preimage of \(M'\subseteq Z\), and \(\G'\to M'\) is a principal right \(H\)-bundle.
Let \(\omega\) be the pullback of \(\omega_G\) to \(\G'\).
Because \(\G'\) is \(V_Z\)-tangent, \(\omega_M=\omega_G\) on \(\G'\), i.e. the pullbacks agree.

Take any smooth vector field \(u\) on \(\G\).
Construct from it, as above, a \(V_Z\)-tangent vector field \(v=(u,w)\), i.e. at each point \((p,g)\in\G\times G\), we let \(w\in T_g G\) be the left \(g\)-translate \(w:=\LT{g*}A\) of \(A:=u\hook\omega_M\).
Take a flow line of \(u\) in \(\G\), say \(p(t)\in\G\), say with \(p(0)=p_0\).
Pick a point \(g_0\in G\).
Let \(A(t):=\dot{p}(t)\hook\omega_M\).
We can construct a unique smooth curve \(g(t)\in G\) by solving the Lie equation
\[
\dot{g}(t)=A(t)
\]
with \(g(0)=g_0\), by theorem~\vref{theorem:Lie.equations}.
The curve \((p(t),g(t))\in\G\times G\) is a flow line of \(v\), and this flow line exists for as long as \(p(t)\) does.
In other words, the obvious projection map \(\G\times G\to\G\) is equicomplete taking \(V_Z\) to \(T\G\).

By theorem~\vref{theorem:orbit.maps} the map \(\G'\to\G\) is a fiber bundle mapping.
By theorem~\vref{cor:orbitMapEquiv}, the quotient \(M':=\G'/H\to M\) is also a fiber bundle mapping.
Since each is a local diffeomorphism, each is a covering map.
Every point \((p_0,g_0)\) lies in a unique \(H\)-folio.
The left \(G\)-action on \(Z\):
\[
g(p_0,g_0):=(p_0,gg_0)
\]
preserves \(V_Z\) and commutes with the \(H\)-action, hence permutes the \(H\)-folios in \(Z\) and so permutes the leaves in \(\bar{Z}\).

Each \(H\)-folio \(\G'\) is a covering space of \(\G\), so maps onto \(\G\).
So for any point \(p_0\in\G\) and any \(H\)-folio \(\G'\subseteq Z\), there is a point \(g_0\in G\) so that \((p_0,g_0)\in\G'\).
Take any two \(H\)-folios \(\G',\G''\subseteq Z\).
Find points \((p_0,g')\in\G'\), \((p_0,g'')\in\G''\) in them.
Let \(g:=g''(g')^{-1}\).
Then \(\G''\) and \(\LT{g}\G'\) are \(H\)-folios in \(Z\) passing through the same point \((p_0,g'')\).
So \(\G''=\LT{g}\G'\), i.e. the \(H\)-folios in \(Z\) are transitively permuted by the left \(G\)-action.
Hence the leaves in \(\bar{Z}\) are transitively permuted by the left \(G\)-action on \(\bar{Z}\).

Since \(M'\to M\) is a covering map, \(M'\) and \(M\) share the same universal covering space
\[
\tilde{M}\to M'\to M.
\]
Lift the bundle \(\G'\to M'\) to \(\tilde{M}\):
\[
\begin{tikzcd}
\G\arrow[d]&\tilde\G\arrow[d]\arrow[r]\arrow[l]&G\arrow[d]\\
M&\tilde{M}\arrow[l]\arrow[r]&X.
\end{tikzcd}
%\begin{tikzcd}
%&\tilde\G\arrow[dl]\arrow[dr]\arrow[rrr,bend left]&&&\tilde{M}\arrow[dl]\arrow[dr]&\\
%\G\arrow[rrr,bend right]&&G\arrow[rrr,bend right,crossing over]&M&&X.
%\end{tikzcd}
\]
The deck transformations of \(\tilde{M}\) act as bundle automorphisms on \(\tilde\G\), since the bundle is pulled back from \(M\).
These transformations preserve \(\omega\) and the \(H\)-action, as these are also pulled back from \(M\).
But the map to \(G\) might not be invariant under the deck transformations.

The map \(\tilde\G\to\G'\to Z\) is a covering map to an \(H\)-folio, and remains such under deck transformations.
So it locally diffeomorphically identifies \(\tilde\G\) with that \(H\)-folio.
The deck transformations preserve \(\omega\) and the right \(H\)-action, so take it to another \(H\)-folio, so acting by left translation by  some constant element \(g\in G\).
This element is thus locally constant on \(\tilde\G\).
Since the left \(G\)-action commutes with the right \(H\)-action, this element \(g\in G\) is locally constant on \(\tilde{M}\).
Since \(M\) is connected, so is \(\tilde{M}\), so \(g\in G\) is constant, depending only on \(\gamma\in\pi\); denote it \(h(\gamma)\).
By uniqueness of \(H\)-folio up the left \(G\)-action, \(\gamma\mapsto h(\gamma)\) is a group morphism.

We have a developing map \(\tilde{M}\xrightarrow{\delta}M\) given by quotienting \(\tilde{G}\to X\), and a holonomy morphism.
The construction of the developing map and holonomy morphism is determined by the choice of \(M'\subseteq Z\), i.e.  the choice of folio \(\G'\subseteq\G\times G\).

If we pick the folio \(\G'\) which passes through some point \((p_0,g_0)\), changing our choice to \((p_0,gg_0)\) alters the developing map and holonomy morphism as described.
But every folio maps to \(p_0\), since every folio covers \(\G\).
\end{proof}

\chapter{The moving frame}
The reader may skip this chapter without loss of continuity.
We will consider the differential invariants of submanifolds of homogeneous spaces.
\section{Example: the Frenet--Serret frame}
Consider an immersed curve \(C\xrightarrow{\varphi} X\) in Euclidean space \(X:=\R^3\).
For each point \(p\in C\), there is a translation taking 
\[
x=
\begin{pmatrix}
x_1\\
x_2\\
x_3
\end{pmatrix}
=\varphi(p)
\]
to the origin.
Rotate the tangent line into the direction of the \(x_1\)-axis, so that the curve becomes locally the graph \(x_2=x_2(x_1),x_3=x_3(x_1)\) with
\[
x_2=\frac{1}{2}k_2 x_1^2+\dots, x_3=\frac{1}{2}k_3 x_1^2+\dots.
\]
Rotate the normal plane to get the vector \((k_2,k_3)\) into the form \((k,0)\):
\[
x_2=\frac{1}{2}kx_1^2+\dots, x_3=\frac{1}{6}\tau x_1^3+\dots.
\]
These expressions \(k,\tau\) are differential invariants; it turns out that \(k\) is the curvature and \(\tau/k\) is the torsion in the standard theory of curves.

We can view this computation of a Taylor series as a construction of a submanifold inside the rigid motion group \(G:=\Orth{3}\ltimes\R^3\) of \(X\).
Note that \(H:=\Orth{3}\) is the stabilizer of the origin \(0\in X\).
Write elements of \(G\) as matrices
\[
g=
\begin{pmatrix}
h&x\\
0&1
\end{pmatrix}, h\in\Orth{3}, x\in\R^3,
\]
acting on points \(y\in X\) by
\[
gy:=hy+x.
\]
Take a point \(p\in C\).
It has image \(x:=\varphi(p)\in X\).
This point is the image of the origin under translation by \(x\), or indeed by any element of \(G\) of the form
\[
g=
\begin{pmatrix}
h&x\\
0&1
\end{pmatrix}.
\]
Let \(\G\subseteq C\times G\) be the set of all pairs \((p,g)\) for \(p\in C\) and
\[
g=
\begin{pmatrix}
h&x\\
0&1
\end{pmatrix}
\]
with \(x=\varphi(p)\).
We have two obvious maps: \((p,g)\in\G\mapsto g\in G\) and \((p,g)\in\G\mapsto p\in C\); the second map is a principal right \(H\)-bundle:
\[
\begin{tikzcd}
\G\arrow[r]\arrow[d]&G\arrow[d]\\
C\arrow[r]&X.
\end{tikzcd}
\]
Note that we can construct this diagram, i.e. this bundle \(H\to\G\to C\) and map \(\G\to G\), for any map \(C\to X\), which need not be an immersed curve.
The construct takes zero derivatives at each point.

If \(C\) is an immersed curve in \(X\) then we can do better.
First, we note that \(C\to X\) is an immersion if and only if \(\G\to G\) is an immersion.
The curve \(C\) has a tangent line at each point, with one of two possible unit tangent vectors.
(We resist the temptation to orient the curve \(C\), i.e. prefer one unit tangent vector to the other.)
Define a submanifold \(\G'\subset\G\): we let \(\G'\) be the set of pairs \((p,g)\in\G\) so that \(g\) takes the \(x_1\)-axis to the tangent line:
\[
g
=
\begin{pmatrix}
h&x\\
0&1
\end{pmatrix}
\]
with \(x=\varphi(p)\) and 
\[
h
\begin{pmatrix}
1\\
0\\ 
0
\end{pmatrix}
\in T_p C.
\]
So then \(h\) is determined up replacing by \(hh'\) where \(h'\) is an orthogonal transformation with preserving \((1,0,0)\) up to \(\pm\).
So \(h'\) can orthogonally transform arbitrarily in perpendicular directions to \((1,0,0)\).
So \(h'\in H'\) where \(H'\subset H\) is the stabilizerof the origin and the\(x_1\)-axis, i.e. the set of orthogonal matrices of the form
\[
\begin{pmatrix}
\pm 1&0&0\\
0&a&b\\
0&c&d
\end{pmatrix}
\]
where
\[
\begin{pmatrix}
a&b\\
c&d
\end{pmatrix}\in\Orth{2}
\]
so we can write
\[
\begin{pmatrix}
a&b\\
c&d
\end{pmatrix}
=
\begin{pmatrix}
\cos\theta&\mp\sin\theta\\
\sin\theta&\pm\cos\theta
\end{pmatrix}.
\]
So \(H'\to\G'\to C\) is a principal right \(H'\)-bundle.
Let \(X':=G/H'\) be the set of pairs \((x,\ell)\) where \(x\in X\) and \(\ell\subseteq T_x X\) is a tangent line.
Map \((p,g)\in\G'\mapsto g\in G\), and then quotient by \(H'\) to get a diagram
\[
\begin{tikzcd}
\G'\arrow[r]\arrow[d]&G\arrow[d]\\
C\arrow[r]&X'
\end{tikzcd}
\]
which simply maps each point of \(C\) to its tangent line.
\section{Reduction}
Take smooth actions of Lie groups \(G,G'\) on manifolds \(Q,Q'\).
An \emph{action morphism} \(Q'\to Q\) is a smooth map equivariant for a Lie group morphism \(G'\to G\).
\begin{example}
Every morphism of homogeneous spaces is an action morphism.
\end{example}
\begin{example}
Let \(V\) be a finite dimensional real or complex vector space and \(G=\GL{V}\) its group of real or complex linear automorphisms.
Let \(Q\) be the set of all linear maps \(V\to V\) with only simple eigenvalues.
Fix a basis of \(V\) and let \(Q'\) be the set of all linear maps \(V\to V\) which are diagonal in that basis, with only simple eigenvalues.
Let \(G'\) be the group of linear maps permuting those basis elements up to rescaling, so \(G'\) has dimension equal to the dimension of \(V\).
\end{example}
As in this example, we always think of the elements of \(Q'\) as the elements of \(Q\) brought to some ``normal form''.
Given an action morphism \(Q'\to Q\) over \(G'\to G\), define a \(G'\)-action on \(G\times Q'\):
\[
h(g,q):=(g\Phi(h)^{-1},hq)
\]
and a map, the \emph{extrusion map},
\[
(g,q)\in G\times Q'\mapsto g\varphi(q)\in Q.
\]
\[
\includegraphics[width=10cm]{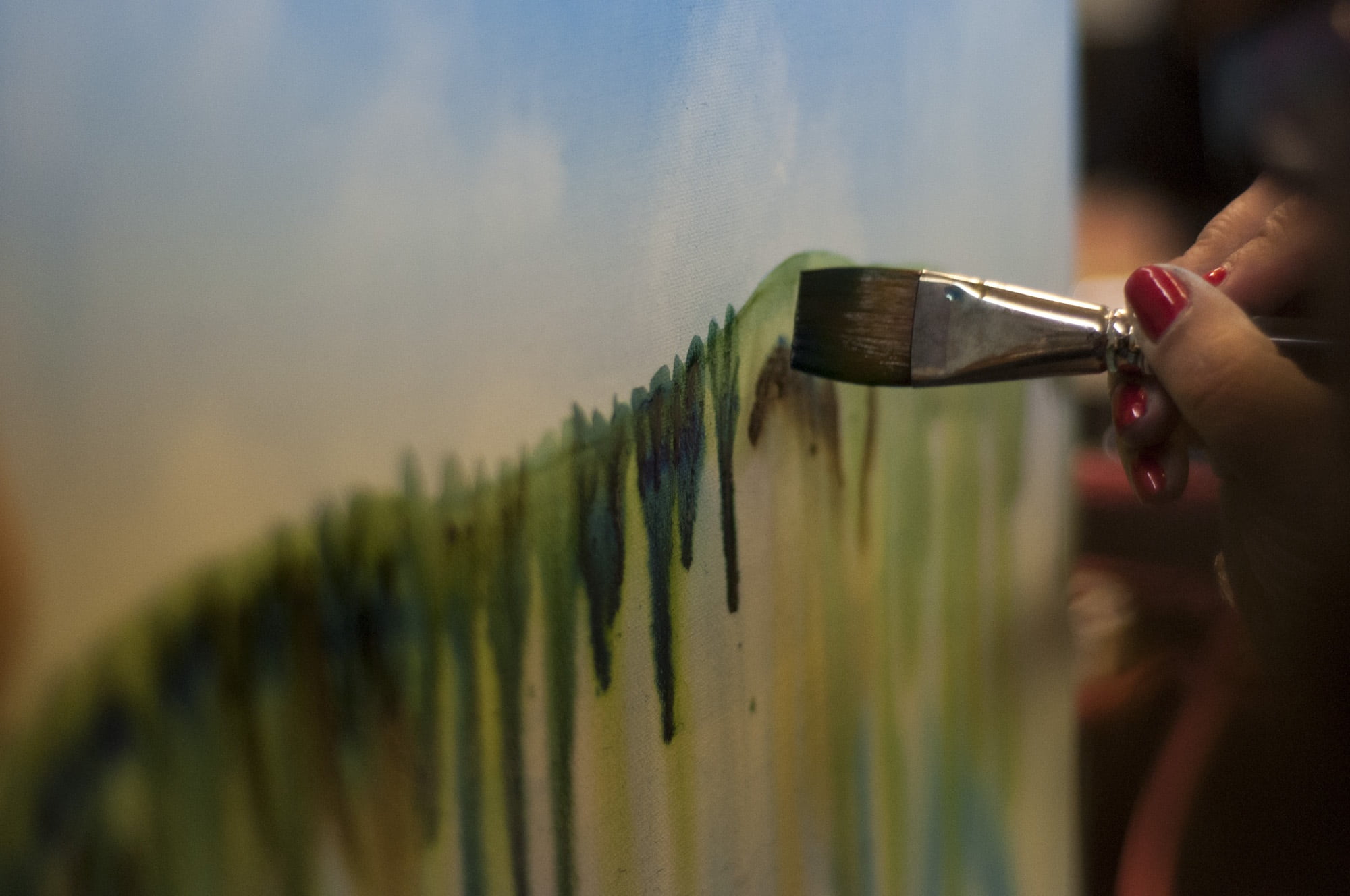}
\]
\begin{example}
If \(Q\) is a wall, \(G\) is the real line acting by moving points up and down the wall, \(Q'\) is a horizontal line on the wall along which we trace a paint brush, and the extrusion map extrudes the paint up and down the wall.
\end{example}
The extrusion map is clearly invariant under the \(G'\)-action.
A \emph{reduction} is an action morphism for which the extrusion map is a principal \(G'\)-bundle under this \(G'\)-action.
Intuitively, reductions are rare and useful, so it helps to have a tool for identifying them, when we suspect we have one.
See \cite{Guillemin.Pollack.2010} chapter 2 for the theory of transverse maps.
\begin{lemma}\label{lemma:extrusion.submersion}
Take an action morphism \(Q'\to Q, G'\to G\).
The associated extrusion map is a submersion just when \(Q'\to Q\) is transverse to all \(G\)-orbits in \(Q\).
\end{lemma}
\begin{proof}
At each point \((1,q')\in G\times Q'\), the extrusion map maps the tangent space to \(G\times Q'\) to the sum of the tangent space to the \(G\)-orbit and the image under the differential of \(Q'\to Q\) of \(T_{q'} Q'\).
So the extrusion map is a submersion at every point of \(1\times Q'\subset G\times Q'\) just when \(Q'\to Q\) is transverse to all \(G\)-orbits.
Translation by any \(g\in G\) permutes tangent spaces to \(G\)-orbits, and so takes any linear subspace transverse to that tangent space to another also transverse.
\end{proof}
\begin{proposition}\label{prop:prereduction}
Take a Lie group morphism \(G'\to G\) with closed image.
An action morphism \(Q'\to Q\) over \(G'\to G\) is a reduction if and only if all of the following hold:
\begin{itemize}
\item
every element of \(Q\) is, up to \(G\)-action, the image of an element of \(Q'\), uniquely precisely up to \(G'\)-action and
\item
the kernel of \(G'\to G\) acts freely and properly on \(Q'\) and
\item
the map \(Q'\to Q\) is transverse to the \(G\)-orbits in \(Q\).
\end{itemize}
\end{proposition}
\begin{proof}
Reductions have the first two of these properties, and the third by lemma~\vref{lemma:extrusion.submersion}.
Suppose conversely that we have an action morphism with those three properties.

Freedom of \(G'\) acting on \(G\times Q'\) is precisely \emph{relative freedom} of \(Q'\to Q\), i.e. 
\[
h'(g,q')=(g,q')
\] 
for some \(h'\in G'\), \(g\in G\) and \(q'\in Q\) just when \(\Phi(h')=1\) and \(h'q'=q'\).
So \(G'\) acts freely on \(G\times Q'\) just when any element of the kernel of \(G'\to G\) acts freely on \(Q'\).

Properness of \(G'\) acting on \(G\times Q'\) is precisely \emph{relative properness} of \(Q'\to Q\), i.e. given  a convergent sequence
\[
(g_1,q'_1),
(g_2,q'_2),
\dots
\to
(g,q'),
\]
and a sequence \(h_1',h_2'\in G'\) so that
\[
h_1'(g_1,q'_1),
h_2'(g_2,q'_2),
\dots
\] 
converges, say to \((\bar{g},\bar{q}')\), properness states that
\[
h_1',h_2',\dots
\] 
is convergent as well, perhaps after taking an infinite subsequence.
Let \(h_i:=\Phi(h_i')\in G\).
This expands out to asking that if
\begin{align*}
q'_1,q'_2\to q'&\in Q',\\
h_1',h_2',\dots,&\in G',\\
h_1'q'_1,h_2'q'_2,\dots\to\bar{q}'&\in Q',\\
h_1,h_2,\dots\to h&\in G,
\end{align*}
then 
\[
h_1',h_2',\dots\in G'
\]
converges, perhaps after taking an infinite subsequence.
Pick any \(h_1',h_2',\dots\) in the kernel of \(G'\to G\): the kernel of \(G'\to G\) acts properly on \(Q'\).

Suppose that \(G'\to G\) has closed image \(H\subseteq G\), and that the kernel of \(G'\to G\) acts properly on \(Q'\). 
Again we pick some sequences
\begin{align*}
q'_1,q'_2\to q'&\in Q',\\
h_1',h_2',\dots,&\in G',\\
h_1'q'_1,h_2'q'_2,\dots\to\bar{q}'&\in Q',\\
h_1,h_2,\dots,\to h&\in G,
\end{align*}
where \(h_i:=\Phi(h_i')\).
Since \(h_i\in H\) and \(H\subseteq G\) is closed, \(h\in H\), so \(h\) is the image of some \(h'\in G'\).

Being a surjective Lie group morphism, \(G'\to H\) is a fiber bundle mapping, with fibers diffeomorphic to the kernel
\[
1\to K\to G'\to H\to 1.
\]
Take a local section \(s\) of that fiber bundle mapping \(G'\to H\), defined on some open set \(U\subseteq H\) containing \(h\in H\), say
\[
U\xrightarrow{s}G'.
\]
Translate to ensure that \(s(h)=h'\).
Then \(s_i:=s(h_i)\to s(h)=h'\).
Each \(s_i\) projects to \(h_i\) so agrees with \(h_i'\) up to translation by some element of the kernel \(K\), say
\[
h_i'=k_is_i.
\]
Define a new symbol \(q_i'\) to represent the old \(s_iq_i'\), and a new \(h_i\) to represent the old \(k_i\in K\).
The values \(h_i'q'_i\) are unchanged, as is convergence of \(q_i'\).
So it suffices to suppose that \(h_i'=k_i\in K\).
But the kernel acts properly, so \(k_1,k_2,\dots\)  converges, perhaps after taking an infinite subsequence.

The \(G'\)-action is free and proper, so the quotient \(G\times^{G'}Q'\) is a smooth manifold and
\[
G'\to G\times Q'\to G\times^{G'}Q'
\]
is a smooth principal bundle.
The extrusion map, call it \(\alpha\), is a submersion, as we have seen, and \(G'\)-invariant; hence  \(\alpha\) descends to a smooth submersion
\[
G\times^{G'} Q'\xrightarrow{\bar\alpha} Q.
\]
The existence of \((g,q')\) ensures that \(\alpha\) is onto, hence so is \(\bar\alpha\).
We has assumed that \((g,q')\) is determined uniquely up to replacing by \((gh^{-1},h'q')\) for an arbitrary \(h'\in G'\) and \(h=\Phi(h')\).
This ensures that \(\alpha\) is injective up to \(G'\)-action, hence \(\bar\alpha\) is injective, so bijective.
So \(\bar\alpha\) is a bijective submersion, hence a diffeomorphism by the inverse function theorem.
\end{proof}
\begin{problem}{reduction:homog}
Prove that a morphism \((X',G')\xrightarrow{\varphi,\Phi} (X,G)\) of homogeneous spaces is a reduction just when
\begin{itemize}
\item
the morphism is surjective on stabilizers \((G')^{x_0'}\to G^{x_0}\) for some \(x_0'\in X'\) mapping to some \(x_0\in X\) and
\item
the kernel of \(G'\to G\) acts freely and properly on \(X'\).
\end{itemize}
\end{problem}
\begin{answer}{reduction:homog}
The whole of \(X\) is a single \(G\)-orbit, so the map \(X'\to X\) is transverse to \(G\)-orbits.
Pick \(x_0'\in X'\) mapping to \(x_0\in X\).
The fiber of the extrusion map over a point \(x=g_0x_0\in X\) consists of pairs  \((g,x')=(g,h'x_0')\), say, with
\[
ghx_0=g_0x_0
\]
where \(h=\Phi(h')\), i.e. 
\[
(g,x')=(g_0h_0h^{-1},h'\ell'x_0),
\]
for any \(h'\in G'\), \(h_0\in G^{x_0}\), \(\ell'\in (G')^{x_0'}\).
Note that varying \(\ell'\) doesn't move this point of \(G\times X'\).
To ensure that the stalks of the extrusion map are precisely the \(G'\)-orbits in \(G\times X'\), we need thus precisely that we can always pick \(\ell'\) to have image \(\ell\in G\) agreeing with \(h_0^{-1}\), i.e. we need precisely that \(G^{x_0}\) is the image of \((G')^{x_0'}\).
\end{answer}
\begin{example}
For any homogeneous space \((X,G)\), if we set \(X':=\set{x_0}\) a single point of \(X\) and \(G':=G^{x_0}\), then \((X',G')\to (X,G)\) is a morphism of homogeneous spaces and a reduction.
\end{example}
\begin{example}
Take a vector space \(V=\R^{p+q}\) with an inner product of signature \(p,q\), with \(p\ge 1\).
Let \(G:=\Orth{V}=\Orth{p,q}\) be the orthogonal group of \(V\).
Let \(Q\subset V\) be the vectors on which the inner product is positive.
Pick a vector \(v_0\in V\) with \(\ip{v_0}{v_0}=1\).
Let \(Q'\subset Q\) be the ray of vectors \(v=tv_0\) for \(t>0\) and \(G'\subset G\) the orthogonal transformations fixing \(v_0\).
By proposition~\vref{prop:prereduction}, \(Q'\) is a \(G'\)-reduction.
Clearly we can apply this to the study of positive definite curves in pseudo--Riemannian manifolds.
\end{example}
\begin{proposition}\label{proposition:reduction.bundle}
If \(Q'\to Q\) is a \(G'\)-reduction of a \(G\)-action, and if \(G\to Q\to M\) is a principal bundle then \(G'\to Q'\to M\) is also a  principal bundle over the same manifold \(M=Q/G=Q'/G'\).
\end{proposition}
\begin{proof}
We want to prove that \(G'\) acts freely on \(Q'\).
Take some \(k'\in G'\) fixing some point \(q'\in Q'\).
Let \(k\in G\) be its image.
So \(k\) fixes the point \(q\in Q\) in the image of \(q'\in Q'\), so \(k=1\).
Hence \(k'\) is in the kernel of \(G'\to G\).
Since \(Q'\to Q\) is a reduction, \(k=1\).
So \(G'\) acts freely on \(Q'\).

We want to show that \(G'\) acts properly on \(Q'\).
Take some sequences
\begin{align*}
g_1',g_2',\dots&\in G',\\
q_1',q_2',\dots &\to 	q'\in Q',\\
g_1'q_1',g_2'q_2',\dots &\to\bar{q}'\in Q'
\end{align*}
and write their images as
\begin{align*}
g_1,g_2,\dots&\in G,\\
q_1,q_2,\dots &\to q\in Q,\\
g_1q_1,g_2q_2,\dots &\to\bar{q}\in Q.
\end{align*}
Because \(G\to Q\to M\) is a principal bundle, after perhaps taking an infinite subsequence, \(g_1,g_2,\dots\) converges, say to \(g\in G\).
By relative properness of reductions (as in the proof of proposition~\vref{prop:prereduction}), \(g_1',g_2',\dots\) converges, perhaps after taking an infinite subsequence.
So \(G'\to P'\to M':=P'/G'\) is a fiber bundle mapping.

The map \(Q'\to Q\to M:=Q/G\) is a composition of smooth maps, so smooth, and \(G'\)-invariant, so drops to a unique smooth map \(M':=Q'/G'\to M\).
The derivative 
\[
T_{q'} Q'\to T_q Q
\]
maps tranverse to the \(G\)-orbits, i.e. fibers of \(Q\to M\), so the composition
\[
T_{q'}Q'\to T_q Q\to T_m M
\]
is onto.
Therefore \(Q'\to M\) is a submersion.
It drops to \(M'\to M\), so \(M'\to M\) is also a submersion.
Take two points \(q_0',q_1'\in Q'\to M\) mapping to the same point \(m_0\in M\).
So the images \(q_0,q_1\in Q\) of these points lie in the same \(G\)-orbit, say 
\[
q_1=gq_0.
\]
Hence \((g,q_0')\) and \((1,q_1')\) are in the same fiber of the extrusion map.
Because \(Q'\to Q\) is a reduction, there is a unique \(h'\in G'\) so that
\[
h'(g,q_0')=(1,q_1'),
\]
i.e. if \(h=\Phi(h')\), 
\[
(gh^{-1},h'q_0')=(1,q_1').
\]
Hence \(q_1'=h'q_0'\), i.e. these represent the same point of \(M'\).
So \(M'\to M\) is an injection.
So \(M'\to M\) is a submersive injection, hence a local diffeomorphism by the inverse function theorem.

Take a point \(m_0\in M\).
Since \(P\to M\) is a principal bundle, \(m_0\) is the image of some point \(p_0\in P\), and is the image precisely of those points in \(Gp_0\).
Let \(q_0\in Q\) be the image of \(p_0\in P\).
The map
\[
(g,q')\in G\times Q'\mapsto gq\in Q,
\]
where \(q\in Q\) is the image of \(q'\in Q'\), is a principal \(G'\)-bundle.
Therefore we can write \(q_0\) as \(q_0=gq\) for some \(q'\in Q'\) mapping to \(q\in Q\) and some \(g\in G\).
If we let \(p:=g^{-1}p_0\), we find \(p\) and \(q\) mapping to \(q'\).
So the point
\[
p'=(p,q',q)
\]
belongs to \(P'\), and maps to \(m_0\in M\).
The map \(P'\to M\) descends to the map \(M'\to M\), so \(M'\to M\) is onto.
So \(M'\to M\) is a submersive bijection, hence a diffeomorphism.
\end{proof}
\begin{theorem}\label{theorem:fiber.product.reduction}
Suppose that \(Q'\to Q\) is a \(G'\)-reduction for some Lie group morphism \(G'\to G\).
Suppose that \(G\) also acts smoothly on a manifold \(P\) with a \(G\)-equivariant smooth map \(P\to Q\).
Let \(P'\) be the \emph{fiber product},\define{fiber product} i.e. the points of \(P'\) are the tuples
\[
p'=(p,q')\in P\times Q'
\]
for which \(p\) and \(q'\) map to the same point \(q\in Q\), with obvious maps
\[
\begin{tikzcd}
P'\arrow[d]\arrow[r]&P\arrow[d]\\
Q'\arrow[r]&Q
\end{tikzcd}
\]
with \(G'\)-action
\[
g'(p,q'):=(gp,g'q')
\]
where \(g'\in G'\) maps to \(g\in G\).
Then \(P'\to P\) is a \(G'\)-reduction.
If the \(G'\)-action on \(P\times Q'\) is free and proper then \(G'\to P'\to M':=P'/G'\) is a principal bundle.
\end{theorem}
\begin{proof}
Since \(P\to Q\) is \(G\)-equivariant, it maps \(G\)-orbits onto \(G\)-orbits by submersions.
Since \(Q'\to Q\) is transverse to \(G\)-orbits and \(P\to Q\) maps \(G\)-orbits to \(G\)-orbits by submersions, \(Q'\to Q\) and \(P\to Q\) are transverse.
Denote our maps
\[
P\xrightarrow{\psi} Q, Q'\xrightarrow{\varphi} Q \text{ and } G'\xrightarrow{\Phi}G.
\]
By transversality, the fiber product \(P'\) is a smooth closed embedded submanifold of \(P\times Q'\) \cite{Guillemin.Pollack.2010} chapter 2, with tangent space \(T_{p'}P'\) at each point 
\[
p'=(p,q')
\]
given by the set of vectors \((v,w')\) for \(v\in T_p P\) and \(w'\in T_{q'} Q\) so that
\[
\psi'(p)v=\varphi'(q')w'\in T_q Q.
\]
The map
\[
p'=(p,q')\in P'\xmapsto{\beta}p\in P
\]
 is an action morphism.

Take an element \(p\in P\) and let \(q:=\psi(p)\in Q\).
Since \(Q'\to Q\) is a \(G'\)-reduction, there is some \(q'\in Q'\) and \(g\in G\) so that 
\[
g\varphi(q')=q.
\]
This pair \((g,q')\) is unique up to replacing by \((gh^{-1},h'q')\) for any \(h'\in G'\) with image \(h\in G\).
So both \(g^{-1}p\) and \(q'\) map to \(g^{-1}q\):
\[
p':=(g^{-1}p,q')\in P'.
\]
Hence
\[
g\beta(p')=p.
\]
In other words, every element of \(P\) is the image of an element of \(P'\), up to \(G\)-action.

The points \((g,p')=(g,p,q')\in P'\) mapping to a given \(p_0\in P\) consist of points with \(gp_0=p\).
Let \(q_0\in Q\) be the image of \(p_0\in P\), so \(gq_0=q\) and \(q'\) maps to \(q\).
So \(q'\) is uniquely determined up to replacing by \(h'q'\), for any \(h'\in G'\), and then we replace
\[
(g,p')=(g,p,q') \mapsto (gh^{-1},h'p')=(gh^{-1},hp,h'q').
\]
So every element of \(P\) is the image of an element of \(P'\), up to \(G\)-action, uniquely up to \(G'\)-action.

Take any \(k\) in the kernel \(K\) of \(G'\to G\).
Suppose that \(kp'=p'\) for some \(p'\in P'\).
Write out \(p'=(p,q')\) as so
\[
kp'=(p,kq')=p'=(p,q'),
\]
and then \(kq'=q'\), so \(k=1\): \(K\) acts freely on \(P'\).

Take a sequence \(k_1,k_2,\dots \in K\) and a convergent sequence \(p_1',p_2',\dots p'\in P'\).
Suppose that 
\[
k_1p_1',k_2p_2',\dots \to\bar{p}'\in P'.
\]
Write out \(p'_i\) as
\[
p'_i=(p_i,q'_i).
\]
Expand out:
\begin{align*}
q_1',q_2',\dots &\to 	q',\\
k_1q_1',k_2q_2',\dots &\to\bar{q}'.
\end{align*}
But \(K\) acts properly on \(Q'\), so some infinite subsequence of \(k_1,k_2,\dots\) converges.
So \(K\) acts properly on \(P'\).

We want to prove that \(P'\to P\) is tranverse to the \(G\)-orbits in \(P\).
Pick a point \(p'=(p,q')\in P'\).
Let \(q:=\psi(p)\).
Take some vector \(v\in T_p P\).
Map it to a vector \(w\in T_q Q\) by \(\psi'(p)v=w\).
By transversality, we can write
\[
w=A(q)+\varphi'(q')w',
\]
for some vector \(A\in\LieG\) and some vector \(w'\in T_{q'}Q'\).
So then
\[
v-A(p)\mapsto w-A(q)=\varphi'(q')w'.
\]
So then
\[
v':=(v-A(p),w')\in T_{p'}P'\mapsto v-A(p)\in T_p P,
\]
hence
\[
v=A(p)+\beta'(p')v':
\]
\(P'\xrightarrow{\beta}P\) is transverse to the \(G\)-orbits in \(P\).

Since \(P'\subseteq P\times Q'\) is a closed embedded submanifold, if \(G'\) acts freely (or properly) on \(P\times Q'\) then it does so on \(P'\).
\end{proof}
\section{The Gauss map}
Recall that if \(V\) is a vector space, \(\Gr{k}{V}\) is the \emph{Grassmannian}: the set of \(k\)-dimensional linear subspaces of \(V\).
Take a manifold \(M\) and an integer \(0\le k\le\dim M\).
The \emph{Grassmann bundle}\define{Grassmann bundle}\define{bundle!Grassmann} \(\Gr{k}{M}\) is the set of pairs \((m,V)\) where \(m\in M\) and \(V\subseteq T_m M\) is a linear subspace of dimension \(k\).
The Grassman bundle \(\Gr{k}{M}\to M\) is a fiber bundle with fibers \(\Gr{k}{T_m M}\subseteq\Gr{k}{M}\).
The \emph{Gauss map}\define{Gauss map} of an immersed \(k\)-dimensional submanifold \(S\xrightarrow{\iota_S}M\) is the map
\[
s\in S\xmapsto{\Gau}T_s S\in\Gr{k}{M}
\]
which lifts \(S\to M\):
\[
\begin{tikzcd}
&\Gr{k}{M}\arrow[d]\\
S\arrow[ur,"\Gau"]\arrow[r]&M.
\end{tikzcd}
\]
\begin{problem}{cn}
Write out the Gauss map in coordinates.
\end{problem}
\section{The Grassmann bundle of a homogeneous space}
Henceforth in this chapter, take a homogeneous space \((X,G)\), a point \(x_0\in X\), and let \(H:=G^{x_0}\).
Define
\[
g\in G\xrightarrow{\pi}gx_0\in X.
\]
The map
\[
(g,V)\in G\times\Gr*{k}{\LieG/\LieH}\mapsto \LT{g*}V\subseteq T_{gx_0} X
\]
is invariant under the right \(H\)-action
\[
(g,V)h:=(gh,\Ad_h^{-1}V)
\]
so drops to a map
\[
\amal{G}{H}{\Gr*{k}{\LieG/\LieH}}\to\Gr{k}{X}.
\]
\begin{problem}{grass.map}
Prove that this map is a diffeomorphism on the fibers over \(x_0\).
\end{problem}
It is \(G\)-equivariant under
\[
g'(g,V):=(g'g,V),
\]
so an isomorphism of fiber bundles
\[
\begin{tikzcd}
\amal{G}{H}{\Gr*{k}{\LieG/\LieH}}\arrow[rr]\arrow[dr]&&\Gr{k}{X}\arrow[dl]\\
&X&
\end{tikzcd}
\]

We can invert this map trivially by inverse left translation: take any linear subspace \(T\subseteq T_x X\), pick any \(g\in G\) with \(gx_0=x\), and then map \(T\) to
\[
V:=\LT{g*}^{-1}T,
\]
so that
\[
T=\LT{g*}V.
\]

Consider the preimage \(T':=\pi_*^{-1}T\subseteq TG\).
At each point \(g\in G\), \(\omega\) is left translation taking \(T_g G\to\LieG=T_1 G\).
So \(\omega\) applied to \(T'\subseteq T_g G\) left translates \(T'\) to \(V':=\LT{g*}^{-1}T'\subset\LieG\).
The fiber of \(T_g G\to T_x X\) is the kernel of \(\pi'\).
Quotienting by these,  \(\omega\) takes \(T\) to \(V+\LieH\subset\LieG/\LieH\).
Our map
\[
\amal{G}{H}{\Gr*{k}{\LieG/\LieH}}\to\Gr{k}{X}
\]
is left translation by \(g\), while \(\omega\) is left translation by \(g^{-1}\), so \(\omega\) gives the inverse map
\[
\amal{G}{H}{\Gr*{k}{\LieG/\LieH}}\leftarrow\Gr{k}{X}.
\]

\section{The soldering form}
The \emph{soldering form}\define{soldering form} of \((X,G)\) is the \(H\)-equivariant \(1\)-form \(\sigma\) valued in \(\LieG/\LieH\) given by
\[
\sigma=\omega+\LieH\in\nForms{1}{G}\otimes^H(\LieG/\LieH).
\]
The soldering form vanishes precisely on the left invariant vector fields from \(\LieH\), i.e. on the fibers of \(G\to X\).
Applied to our linear subspace \(T\),  the soldering form is the same bundle isomorphism again
\[
T=T'/\ker\pi'(g)\xrightarrow{\sigma}\LieG/\LieH
\]
given precisely by left translation to \(1\in G\).
So in our bundle isomorphism, we  map
\[
T\in\Gr{k}{X}\to(g,V)H\in\amal{G}{H}{\Gr*{k}{\LieG/\LieH}}
\]
where \(\sigma(T')/=V\) and we can take any \(g\in G\) here with \(gx_0=x\) and calculate \(\sigma\) at this \(g\in G\).

\section{The Gauss map of a submanifold in a homogeneous space}
Take any immersed submanifold \(S\xrightarrow{\delta} X\).
Take any point \(g\in G\) so that \(x:=\delta(s)=gx_0\).
As above, we identify
\[
T_s S\in\Gr{k}{X}\to(g,V)H\in\amal{\G_S}{H}{\Gr*{k}{\LieG/\LieH}}
\]
where \(\sigma(T')=V\) and \((s,g)\in\G_S\).
In other words, the Gauss map of \(S\) is
\[
\Gau_{(s,g)}:=g\sigma(T_p\G_S).
\]
By our identification 
\[
\Gr{k}{X}=\amal{G}{H}{(\LieG/\LieH)},
\]
the Gauss map is identified with the image of the soldering form on tangent spaces of \(\G_S\):
\[
p\in\G_S\mapsto\sigma(T_p\G_S)\in\Gr*{k}{\LieG/\LieH}.
\]
Denote \(\G_S\to S\) by \(\G_S\xrightarrow{\pi}S\).
The soldering form \(\sigma\) trivializes the pullback Grassmann bundle:
\begin{align*}
\pi^*\Gr{k}{X}&=\G_S\times\Gr*{k}{\LieG/\LieH},\\
\intertext{hence}
\delta^*\Gr{k}{X}&=\amal{\G_S}{H}{\Gr*{k}{\LieG/\LieH}}
\end{align*}
so the Gauss map lifts to an \(H\)-equivariant commutative diagram
\[
\begin{tikzcd}[cramped]
\G_S\arrow[d]\arrow[r]&
\G_S\times\Gr*{k}{\LieG/\LieH}\arrow[d]\\
S\arrow[r]&\Gr{k}{X}.
\end{tikzcd}
\]

\section{The ``Gauss map'' of a subgeometry}
If \(\G_S\to S\) is merely a subgeometry, perhaps not flat, so perhaps not locally arising from an immersion into \(X\), we still have a soldering form 
\[
\sigma:=\omega+\LieH.
\]
We can take the soldering form's image
\[
p\in\G_S\mapsto\sigma(T_p\G_S)\in\Gr*{k}{\LieG/\LieH}.
\]
Note that \(\sigma\) vanishes precisely on the fibers of \(\G_S\to S\), so has rank \(k:=\dim S\), so 
\[
\sigma(T_p\G_S)\subset\LieG/\LieH
\]
is a \(k\)-dimensional linear subspace.
We can define the \emph{Gauss map}\define{Gauss map} to be the soldering form image, a section of the bundle
\[
\begin{tikzcd}
\amal{\G_S}{H}{\Gr*{k}{\LieG/\LieH}}\arrow[d]\\
S.
\end{tikzcd}
\]

\section{Equations on submanifolds}
A \emph{invariant first order differential equation on \(k\)-dimensional submanifolds}\define{invariant first order differential equation on submanifolds} (or \emph{equation}\define{equation} for short) of a homogeneous space \((X,G)\) is a immersion \(E\xrightarrow{\iota_E}\Gr*{k}{\LieG/\LieH}\) equivariant for a smooth \(H\)-action on a manifold \(E\).
It follows that 
\[
\begin{tikzcd}[cramped]
\amal{G}{H}{E}\arrow[rr]\arrow[dr]&&\Gr{k}{X}\arrow[dl]\\
&X&
\end{tikzcd}
\]
is an immersion of fiber bundles.
An \emph{\(E\)-geometry} is a subgeometry whose Gauss map 
\[
S\xrightarrow{\Gau}\amal{\G_S}{H}{\Gr*{k}{\LieG/\LieH}}
\]
lifts to a map
\[
S\to\amal{\G_S}{H}{E}\to\amal{\G_S}{H}{\Gr*{k}{\LieG/\LieH}}.
\]
(For simplicity of notation, also write the map
\[
S\to\amal{\G_S}{H}{E}
\]
as \(\Gau\).)
If \(E\) is an embedded submanifold of \(\Gr*{k}{\LieG/\LieH}\) then an \(E\)-geometry is precisely a subgeometry whose Gauss map has image in the embedded submanifold 
\[
\amal{\G_S}{H}{E}\subseteq\amal{\G_S}{H}{\Gr*{k}{\LieG/\LieH}}.
\]
An \emph{\(E\)-submanifold} is an immersion \(S\to X\) whose associated subgeometry is an \(E\)-geometry.
Hence the Gauss map lifts to an \(H\)-equivariant commutative diagram
\[
\begin{tikzcd}[cramped]
\G_S\arrow[d]\arrow[r]&\G_S\times E\arrow[d]\arrow[r]&\G_S\times{\Gr*{k}{\LieG/\LieH}}\arrow[d]\\
S\arrow[r]&\amal{\G_S}{H}{E}\arrow[r]&\Gr{k}{X}.
\end{tikzcd}
\]
If \(E\) is an embedded submanifold of \(\Gr*{k}{\LieG/\LieH}\) then an \(E\)-submanifold is precisely an immersion \(S\to X\) whose Gauss map has image in the embedded submanifold 
\[
\amal{G}{H}{E}\subseteq\amal{G}{H}{\Gr*{k}{\LieG/\LieH}}=\Gr{k}{X}.
\]
\section{Reduction}
Take a homogeneous space \(X\), and a point \(x_0\in X\) and let \(H:=G^{x_0}\).
For some integer \(k\), suppose that \(E\to\Gr*{k}{\LieG/\LieH}\) is an equation, and that
\[
\G_S\to S
\]
is an \(E\)-subgeometry.

Since our equation \(E\) is \(H\)-invariant, it may map to many elements of \(\Gr*{k}{\LieG/\LieH}\) which are equivalent under \(H\)-action.
We want to pick out one element (or maybe just a discrete set of elements) from each \(H\)-orbit, which we think of as an element of \(E\) in ``normal form''. 
Suppose that \(H'\to H\) is a Lie group morphism and \(E'\to E\) is an \(H'\)-reduction.
Let \(\G'_S\) be the fiber product
\[
\begin{tikzcd}[cramped]
\G'_S\arrow[d]\arrow[r]&\G_S\arrow[d]\\
E'\arrow[r]&E.
\end{tikzcd}
\]
By theorem~\vref{theorem:fiber.product.reduction} and proposition~\vref{proposition:reduction.bundle}, \(H'\to \G'_S\to S\)  is a principal bundle and \(\G_S\cong\G'_S\times^{H'}H\) as bundles over \(S\).
The subgeometry \(\G_S\to S\) has differential form \(\omega\), which we can pullback to \(\G'_S\).
Let \(X':=G/H'\).
Then \(\G'_S\to S\) is an \((X',G')\)-subgeometry, but not modelled on \((X',G)\), flat if \(\G_S\to S\) is flat.
If \(\G_S\to S\) is associated to an immersion \(S\to X\) then that immersion lifts to an immersion
\[
\begin{tikzcd}[cramped]
{}&X'\arrow[d]\\
S\arrow[r]\arrow[ur]&X
\end{tikzcd}
\]
which is an \(E'\)-submanifold with associated bundle \(\G'_S\to S\) and with Maurer--Cartan form pulling back to \(\omega\).
The lift is defined by taking the composition
\[
\G'_S\to\G_S\to G
\]
and quotienting by \(H'\).
\section{The soldering form after reduction}
The soldering form on \(\G'_S\) is now valued in \(\LieG/\LieH'\), so the Gauss map associated to it is valued in \(\Gr*{k}{\LieG/\LieH'}\).
But the soldering form on \(\G'_S\) projects by \(\LieG/\LieH'\to\LieG/\LieH\) to agree with the soldering form on \(\G_S\).
Let \(F:=E\).
A linear subspace \(P'\subseteq\LieG/\LieH'\) of dimension \(k\) is a \emph{nondegenerate} if \(P'\) is taken by the linear projection 
\[
\LieG/\LieH'\to\LieG/\LieH
\]
to a linear subspace \(P\subseteq\LieG/\LieH\) of the same dimension.
Let \(U\subseteq\Gr*{k}{\LieG/\LieH'}\) be the set of nondegenerate subspaces, a Zariski open set.
Let 
\[
F'\subseteq E'\times U
\] 
be the set of pairs \((e',P')\) so that,  if we let \(e\in E\) be the image of  \(e'\in E'\), and \(P\in\Gr*{k}{\LieG/\LieH}\) the image of \(P\), our equation \(E\) maps \(e\in E\mapsto P\in\Gr*{k}{\LieG/\LieH}\).
If \(E\) is an embedded submanifold of \(\Gr*{k}{\LieG/\LieH}\) and \(E'\) is an embedded submanifold of \(E\), then \(E=F\) but \(F'\) is the preimage of \(E'\) inside \(U\).

On \(\G_S\), the Gauss map is valued in \(F=E\).
On \(\G'_S\), we have the fiber product maps
\[
\begin{tikzcd}[cramped]
\G'_S\arrow[d]\arrow[r]&\G_S\arrow[d]\\
E'\arrow[r]&E.
\end{tikzcd}
\]
we define the \emph{Gauss map} to be valued in \(F'\), by taking
\[
p'\in\G'_S\mapsto (e'(p'),\sigma(T_{p'} P')),
\]
with \(e'(p)\in E'\) the image from the fiber product maps of \(\G'_S\).
So subsequent reductions will occur by picking not a reduction of the equation \(E'\) but a reduction of the equation \(F'\), which contains more information.

The automorphisms of the \((X,G)\)-subgeometry on \(S\) are precisely the automorphisms of the \((X',G)\)-subgeometry:
\[
\Aut[S]=\Aut[S']
\]
so \(\dim\Aut[S]\le\dim S+\dim\LieH'\).
\section{Example: the Frenet--Serret frame}
\begin{theorem}
A connected immersed curve in a Riemannian \(3\)-manifold has automorphism group of dimension at most \(2\), and if equal to \(2\), the curve is a geodesic.
A nongeodesic connected immersed curve has automorphism group of dimension at most \(1\), and if equal to \(1\), the curve has constant curvature and torsion.
\end{theorem}
Note that the automorphisms here are those of the subgeometry, and may include among them the isometries of the ambient \(3\)-manifold, but need not.
First consider curves in \(X=\R^3\) under the action of the rigid motion group \(G=\Orth{3}\ltimes\R^3\), with stabilizer the orthogonal group \(H=\Orth{3}\).
Write elements of \(G\) as matrices
\[
g=
\begin{pmatrix}
h&x\\
0&1
\end{pmatrix}, h\in\Orth{3}, x\in\R^3.
\]
So the Lie algebra \(\LieG\) consists of the matrices
\[
A=
\begin{pmatrix}
B&x\\
0&0
\end{pmatrix}, B\in\SO{3}, x\in\R^3.
\]
The quotient \(\LieG/\LieH\) is identified with the choices of vector \(x\).
The action of \(h\in H\) on \(x\in\R^3\) comes from
\[
\begin{pmatrix}
h&0\\
0&1
\end{pmatrix}
\begin{pmatrix}
0&x\\
0&0
\end{pmatrix}
\begin{pmatrix}
h&0\\
0&1
\end{pmatrix}^{-1}
=
\begin{pmatrix}
0&hx\\
0&0
\end{pmatrix},
\]
so just orthogonal transformation of the vector \(x\).
Consider all possible curves, we set our equation to be
\[
E=\Gr*{1}{\LieG/\LieH}=\RP{2}.
\]
The action of \(H\) on \(E\) is rotation of the real projective plane, which acts transitively: any two lines through the origin rotate into one another.
Since \((E,H)\) is homogeneous, we can take as reduction any one point, say
\[
E':=\set{e_1}
\]
where
\[
e_1=
\begin{bmatrix}
1\\
0\\
0
\end{bmatrix}
\in\RP{2}.
\]
The stabilizer \(H'\subset H\) is the set of orthogonal matrices of the form
\[
\begin{pmatrix}
\pm 1&0&0\\
0&a&b\\
0&c&d
\end{pmatrix}
\]
where
\[
\begin{pmatrix}
a&b\\
c&d
\end{pmatrix}\in\Orth{2}
\]
so we can write
\[
\begin{pmatrix}
a&b\\
c&d
\end{pmatrix}
=
\begin{pmatrix}
\cos\theta&\mp\sin\theta\\
\sin\theta&\pm\cos\theta
\end{pmatrix}.
\]
Write each element of \(\LieG\) as
\[
\begin{pmatrix}
0&-a&-b&x\\
a&0&-c&y\\
b&c&0&z\\
0&0&0&0\\
\end{pmatrix}.
\]
So the elements of \(\LieH'\) are those of the form
\[
\begin{pmatrix}
0&0&0&x\\
0&0&-c&y\\
0&c&0&z\\
0&0&0&0\\
\end{pmatrix}.
\]
So the quotient \(\LieG/\LieH'\) is the space of matrices of the form
\[
\begin{pmatrix}
0&-a&-b&x\\
a&0&*&y\\
b&*&0&z\\
0&0&0&0\\
\end{pmatrix}
\]
and this projects to \(\LieG/\LieH\) to be
\[
\begin{pmatrix}
0&*&*&x\\
*&0&*&y\\
*&*&0&z\\
0&0&0&0\\
\end{pmatrix}.
\]
Therefore the preimage \(F'\subset\Gr*{1}{\LieG/\LieH'}\) of \(E'\) is the set of lines each spanned by some nonzero matrix of the form
\[
\begin{pmatrix}
0&-a&-b&1\\
a&0&*&0\\
b&*&0&0\\
0&0&0&0\\
\end{pmatrix}.
\]
Under the \(H'\)-action, this transforms according to
\[
\begin{pmatrix}
a\\
b
\end{pmatrix}
\mapsto
\begin{pmatrix}
\cos\theta&(-1)^{\delta}\sin\theta\\
\sin\theta&-(-1)^{\delta}\cos\theta
\end{pmatrix}.
\]
There are now going to be two equations we need to consider: 
\begin{align*}
(0=a=b)&=E_1\subset F',\\
(0\ne(a,b))&=E_2\subset F'.
\end{align*}
We split up our study into those curves which satisfy \(E_1\) everywhere, and those which satisfy \(E_2\) everywhere.
Since \(E_1\) is a single point, it admits no reduction.
For \(E_2\)-curves, take the subset \(E_2'\) with:
\[
a>0, b=0
\]
which is an \(H''\)-reduction for the group \(H''\) of matrices of the form
\[
\begin{pmatrix}
\pm 1&0&0\\
0&1&0\\
0&0&\pm 1
\end{pmatrix}.
\]

Let's apply this to curves in any Riemannian \(3\)-manifold \(M\).
If \(S\) is a curve and \(S'\) its lift for the first reduction (note that \(S'=S\) as a curve), and \(S'\) maps by its Gauss map to some elements of \(\LieG/\LieH'\) which have nonzero coefficient of \((a,b)\), then we have a second reduction \(S''\) (and again note that \(S''=S\) as a curve), which has structure group \(H''=\set{\pm 1}^3\), and so \(\G_{S''}\to S\) is at most an \(8\)-to-\(1\) covering.
Hence a curve, not necessarily oriented, with nonzero curvature, in a Riemannian \(3\)-manifold, not necessarily oriented, has subgeometry with automorphism group at most of dimension \(1\), with stabilizer a subgroup of \(H''=\set{\pm 1}^2\).
So for a curve in a Riemannian \(3\)-manifold, if the dimension of the automorphism group exceeds \(1\), then the curve is a geodesic: as we will see,\(E_1\)-curves are precisely geodesics.

Consider this story more explicitly in the Cartan geometry of a Riemannian \(3\)-manifold.
We have Cartan connection
\[
\omega
=
\begin{pmatrix}
0&-\gamma^2_1&-\gamma^3_1&\sigma^1\\
\gamma^2_1&0&-\gamma^3_2&\sigma^2\\
\gamma^3_1&\gamma^3_2&0&\sigma^3\\
0&0&0&0
\end{pmatrix}.
\]
For any curve \(S\), \(\G_S\) has the same \(1\)-forms, but has \(\sigma^1,\sigma^2,\sigma^3\) vanishing on the fibers, so of rank \(1\).
So two of these are some multiples of the third.
Note that \(\sigma\) transforms under the \(H\)-action in the obvious representation
\[
\RT{h}^*\sigma=h^{-1}\sigma.
\]
Hence we can find points of \(\G_S\) where 
\[
0=\sigma^2=\sigma^3
\]
cuts out the relations between the components of \(\sigma\).
This occurs on an \(H'\)-subbundle \(\G_{S'}\), since \(H'\) stabilizes this equation on components of \(\sigma\).
Hence on \(\G_{S'}\), 
\[
0=\sigma^2=\sigma^3.
\]
But now the structure group is reduced to \(H'\), so the three differential forms
\[
\gamma^2_1,\gamma^3_1,\gamma^3_2
\]
are no longer linearly independent.
Indeed, looking at the matrices in \(H'\), we see that \(\gamma^2_1,\gamma^3_1\) vanish on the fibers, so multiples of \(\sigma^1\), say
\[
\begin{pmatrix}
\gamma^2_1\\
\gamma^3_1
\end{pmatrix}
=
\begin{pmatrix}
k^2\\
k^3
\end{pmatrix}
\sigma^1.
\]
Writing the frame on \(M\) dual to \(\sigma^1,\sigma^2,\sigma^3\) as \(e_1,e_2,e_3\), the reader can puzzle out why \(k_2e_2+k_3e_3\) descends to become the curvature \(1\)-form of the space curve \(S\).
For an \(E_1\)-curve, those multiples vanish everywhere, i.e. our equations are now
\[
0=\sigma^2=\sigma^3=\gamma^2_1=\gamma^2_3,
\]
giving us
\[
\omega
=
\begin{pmatrix}
0&0&0&\sigma^1\\
0&0&-\gamma^3_2&0\\
0&\gamma^3_2&0&0\\
0&0&0&0
\end{pmatrix},
\]
which is the equation of a geodesic.
Note that it has only \(2\) independent \(1\)-forms, so automorphism group of dimension at most \(2\).

On the other hand, consider any \(E'_2\)-curve.
 Under \(H'\)-action we can move
\[
(k_2,k_3)\ne(0,0)
\]
by orthogonal linear transformations as we did \((a,b)\) above, so we can arrange 
\[
k_2>0, k_3=0
\]
using our \(H''\)-slice, writing \(k_2\) as \(k\) for simplicity.
We reduce to
\[
0=\sigma^2=\sigma^3=\gamma^2_1-k\sigma^1=\gamma^3_1,
\]
leaving \(\gamma^2_3\) now vanishing on the fibers (which are \(8\) points as \(H''\) is a group of \(4\) elements).
So now
\[
\gamma^2_3=t\sigma^1,
\] 
for a unique function \(t\) on \(S''\):
\[
\omega
=
\begin{pmatrix}
0&0&0&1\\
k&0&-t&0\\
0&t&0&0\\
0&0&0&0
\end{pmatrix}\sigma^1.
\]
The quantities
\[
\begin{tikzcd}[cramped]
&\R^+\\
\G_{S''}\arrow[ur,"k"]\arrow[dr,"t"]&\\
&\R
\end{tikzcd}
\]
are the \emph{curvature}\define{curvature!space curve} and \emph{torsion}.\define{torsion!space curve}
\begin{problem}{co}
How do \(k\) and \(t\) transform under the \(H''\)-action on \(\G_{S''}\)?
\end{problem}
Since \(d\sigma^1=0\), at least locally \(\sigma^1=ds\) for a unique function \(s\) on \(S\).
Hence our eight copies of our curve are integral curves of
\[
\omega
=
\begin{pmatrix}
0&0&0&1\\
k&0&-t&0\\
0&t&0&0\\
0&0&0&0
\end{pmatrix}ds.
\]
Such integral curves exist locally, clearly.
We can do better.
\begin{problem}{cp}
For any bounded smooth functions \(\R\xrightarrow{k,t}\R\) with \(k>0\), prove that there is a curve in \(3\)-dimensional Euclidean space \(\R^3\) with curvature \(k\) and torsion \(t\).
\end{problem}
Any complete Riemannian \(3\)-manifold has a development of such a curve, hence also has a curve with the same curvature and torsion. 
\subsection{Example: projective connections}
Take \((X,G)=(\Proj{n},\PGL{n+1})\) over the real or complex numbers.
Consider \(k\)-dimensional submanifolds of a manifold with projective connection.
It is convenient to write points of \(X\) as nonzero vectors with entries defined up to rescaling, divided into vectors of size \(1,k,n-k\):
\[
x=
\begin{bmatrix}
x^0\\
x^a\\
x^A
\end{bmatrix}.
\]
Write elements of \(G\) as invertible matrices with entries defined up to rescaling.
It is convenient to split each matrix into blocks as
\[
g=
\begin{bmatrix}
g^0_0 & g^0_b & g^0_B \\
g^a_0 & g^a_b & g^a_B \\
g^A_0 & g^A_b & g^A_B
\end{bmatrix}
\]
in blocks of sizes
\[
\begin{bNiceMatrix}[margin,first-row,first-col]
 & 1 & k & n-k \\
1 & \ & \ & \ \\
k & \ & \ & \ \\
n-k & \ & \ & \ \\
\end{bNiceMatrix}
\]
The subgroup \(H\) consists of the matrices
\[
h=
\begin{bmatrix}
h^0_0 & h^0_b & h^0_B \\
0 & h^a_b & h^a_B \\
0 & h^A_b & h^A_B
\end{bmatrix}.
\]
The quotient \(\LieG/\LieH\) consist of the matrices
\[
\left[
\begin{NiceMatrix}
\Block[fill=gray!20]{1-3} \ & \ & \ \\
A^a_0 & \Block[fill=gray!20]{1-2} \ & \ \\
A^A_0 & \Block[fill=gray!20]{1-2} \ & \
\end{NiceMatrix}
\right],
\]
which we write as
\[
\left[
\begin{NiceMatrix}
\Block[fill=gray!20]{1-3} \ & \ & \ \\
x^a & \Block[fill=gray!20]{1-2} \ & \ \\
x^A & \Block[fill=gray!20]{1-2} \ & \
\end{NiceMatrix}
\right],
\]
and which we think of as
\[
\LieG/\LieH\cong\R^n.
\]
The action of \(H\) on this quotient is
\[
h
\begin{pmatrix}
x^a\\
x^A
\end{pmatrix}
=
\frac{1}{h^0_0}
\begin{pmatrix}
h^a_b & h^a_B \\
h^A_b & h^A_B
\end{pmatrix}
\begin{pmatrix}
x^a\\
x^A
\end{pmatrix}.
\]
These are arbitrary linear transformations: the projective transformations fixing a point of projective space act on its tangent space by arbitrary linear transformations.

We can see this geometrically: the linear transformations acting on a vector space \(V\) act on its projectivization, with rescalings acting trivially.
The tangent space to the projective space is
\[
T_{x_0}X=T_{[v]}\Proj{}V=[v]^*\otimes(V/[v]),
\]
acted on by any invertible linear transformation that fixes the vector \(v\) up to rescaling, so can induced any invertible linear transformation on the quotient \(V/[v]\).

We take as equation just the entire Grassmannian \(E=\Gr*{k}{\LieG/\LieH}\), so making no hypothesis on our submanifold of projective space.
Since \(E\) is homogeneous, we can take as reduction a single point
\[
E':=\set{\operatorname{span}(e_0,e_1,\dots,e_k)/[e_0]},
\]
in the standard basis \(e_0,\dots,e_n\), i.e. so that \(x^A=0\), or we could write this as \(A^A_0=0\).
The group \(H'\subset H\) stabilizing this linear subspace is the group of invertible matrices of the form
\[
h=
\begin{bmatrix}
h^0_0 & h^0_b & h^0_B \\
0 & h^a_b & h^a_B \\
0 & 0 & h^A_B
\end{bmatrix}
\]
So the quotient \(\LieG/\LieH'\) consists of the matrices of the form
\[
\left[
\begin{NiceMatrix}
\Block[fill=gray!20]{1-3} \ & \ & \ \\
A^a_0 & \Block[fill=gray!20]{1-2} \ & \ \\
A^A_0 & A^A_b & \Block[fill=gray!20]{1-1} \  
\end{NiceMatrix}
\right].
\]
The preimage of the subspace \(E'\) is cut out by \(A^A_0=0\):
\[
A=
\left[
\begin{NiceMatrix}
\Block[fill=gray!20]{1-3} \ & \ & \ \\
A^a_0 & \Block[fill=gray!20]{1-2} \ & \ \\
0 & A^A_b & \Block[fill=gray!20]{1-1} \  
\end{NiceMatrix}
\right].
\]
Under the adjoint action, \(H'\) acts on these matrices by
\[
hAh^{-1}
=
\left[
\begin{NiceMatrix}
\Block[fill=gray!20]{1-3} \ & \ & \ \\
h^a_c A^c_0/h^0_0 & \Block[fill=gray!20]{1-2} \ & \ \\
0 & h^A_CA^C_d (h^{-1})d_b & \Block[fill=gray!20]{1-1} \  
\end{NiceMatrix}
\right].
\]
So in the Grassmannian, our equation \(F'\) consists of the linear subspaces of such matrices which project to have \(A^a_0\) coefficients arbitrary, i.e. project to a \(k\)-dimensional subspace in \(\LieG/\LieH\), with \(A^A_0=0\).
Hence on any such subspace, the coefficients \(A^A_b\) are linear functions of these \(A^a_0\), i.e.
\[
A^A_b=a^A_{bc}A^c_0,
\]
for unique coefficients \(a^A_{bc}\).
Hence each point of \(F'\) is a choice of such coefficients \(a^A_{bc}\) for all indices \(A,b,c\).

In fact, we can see that we can do a bit better.
Take a projective connection \(\omega\) on a manifold \(M\), which we write as
\[
\omega=
\begin{bmatrix}
\omega^0_0 & \omega^0_b & \omega^0_B \\
\omega^a_0 & \omega^a_b & \omega^a_B \\
\omega^A_0 & \omega^A_b & \omega^A_B
\end{bmatrix}
\]
as above, the soldering form is
\[
\sigma=
\begin{pmatrix}
\omega^a_0\\
\omega^A_0
\end{pmatrix}.
\]
Take an immersed submanifold \(S\) of \(M\), our reduction \(S'\) to our slice has its soldering form in the span of \(e_1,\dots,e_k\), so has \(0=\omega^A_0\).
Differentiate to find
\begin{align*}
0&=d\omega^A_0,
\\
&=
-\omega^A_b\wedge\omega^b_0
+\frac{1}{2}k^A_{0ab}\omega^a_0\wedge\omega^b_0
\end{align*}
hence
\[
\omega^A_b=a^A_{ab}\omega^b_0,
\]
where we can write
\[
a^A_{ab}=-\frac{1}{2}k^A_{ab}+s^A_{ab},
\]
and compute that
\[
s^A_{ab}=s^A_{ba}.
\]
We leave the reader to check that, if \(e_1,\dots,e_n\) is the basis of \(T_m M\) dual to the coframing given by \(\omega^1_0,\dots,\omega^n_0\), then
\[
\shapeOp:=s^A_{ab}\omega^a_0\otimes\omega^b_0e_A
\]
is the pullback to \(\G_{S'}\) of unique tensor, \emph{the shape operator}\define{shape operator} \(\shapeOp\), which is a symmetric \(2\)-tensor on \(S\) valued in the normal bundle \(T_s M/T_s S\).

Consider the special case of \(S\) a hypersurface in \(M\), i.e. \(k=\dim S=n-1\), so the shape operator is a quadratic form valued in the normal line bundle, and the capital indices take only the value \(n\). 
Suppose that the shape operator is everywhere positive definite.
Then it imposes a conformal Riemannian metric on the hypersurface \(S\).
Under the action of \(H'\), we can arrange as slice that we require \(a^A_{ab}\) to be the identity matrix.
This condition is preserved by the subgroup \(H''\subset H'\) of matrices
\[
h=
\begin{bmatrix}
h^0_0 & h^0_b & h^0_B \\
0 & h^a_b & h^a_n \\
0 & 0 & h^n_n
\end{bmatrix}
\]
with 
\[
h^d_bh^d_c=h^n_n\delta_{bc}.
\]
So then on \(S''\),
\[
\omega^n_a=\omega^a_0,
\]
which we differentiate to find
\[
\omega^a_b+\omega^b_a-\delta_{ab}(\omega^0_0+\omega^n_n)=
(s_{abc}+\frac{1}{2}(k^n_{abc}-k^a_{0bc}))\omega^c_0
\]
for some \(s_{abc}\) symmetric in all lower indices.
We leave the reader to check that 
\[
s_{abc}\omega^a_0\otimes\omega^b_0\otimes\omega^c_0\otimes e_n
\]
is a cubic form valued in the normal bundle.
For example, if \(S\) is a surface, and the cubic form is not everywhere zero, the automorphism group has to preserve the zero lines of the cubic form, so has finite stabilizer of each point, hence dimension at most \(2\), so Lie algebra either abelian or the unique nonabelian Lie algebra of dimension \(2\).

\chapter{Connections}
\[
\includegraphics[width=4cm]{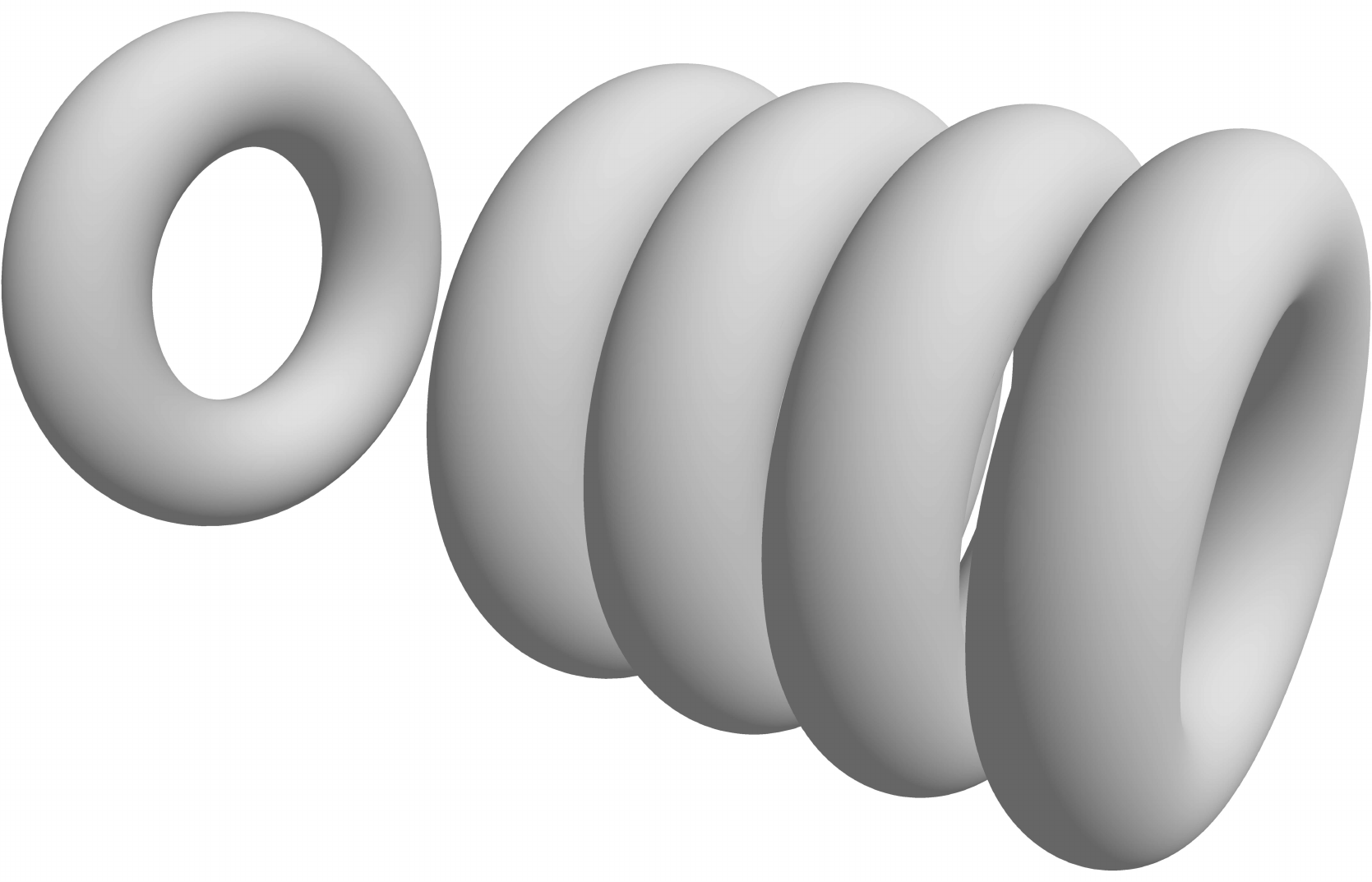}
\]
\section{Definition}
Each fiber \(\Bun_{m_0}\) of a right principal \(G\)-bundle \(\Bun\to M\), say over a point \(m_0\in M\), is homogeneous under \(G\), with trivial stabilizer.
Pick a point \(p_0\in \Bun_{m_0}\) and map \(g\in G\to p_0g \in \Bun_{m_0}\), a diffeomorphism.
If we replace \(p_0\) by a different point of \(E_{m_0}\), we change the diffeomorphism by a left translation.
So the left invariant Maurer--Cartan form \(\omega_G\) is well defined on each fiber \(\Bun_{m_0}\), independent of choice of point \(p_0\).
But it is defined only on tangent vectors of \(\Bun\) which are tangent to the fibers \(\Bun_{m_0}\) (the \emph{vertical vectors}).\define{vertical!vector}
A \emph{connection}\define{connection} on a right principal \(G\)-bundle \(\Bun\to M\) is a \(1\)-form \(\omega_{\Bun}\) on \(\Bun\), valued in \(\LieG\), equal to the left invariant Maurer--Cartan form on the fibers, transforming in the adjoint representation under the right \(G\)-action:
\[
\RT{g}^*\omega_{\Bun}=\Ad_g^{-1}\omega_{\Bun}.
\]
\begin{problem}{left.translaty}
Translate the definitions from right principal bundles above into definitions for left principal bundles.
\end{problem}
\begin{answer}{left.translaty}
Each fiber \(\Bun_{m_0}\) of a left principal \(G\)-bundle \(\Bun\to M\), say over a point \(m_0\in M\), is homogeneous under \(G\), with trivial stabilizer.
Pick a point \(p_0\in \Bun_{m_0}\) and map \(g\in G\to gp_0 \in \Bun_{m_0}\), a diffeomorphism.
If we replace \(p_0\) by a different point of \(E_{m_0}\), we change the diffeomorphism by a right translation.
So the right invariant Maurer--Cartan form \(\rho_G=-\Ad_g\omega_G\) is well defined on each fiber \(\Bun_{m_0}\), independent of choice of point \(p_0\).
But it is defined only on tangent vectors of \(\Bun\) which are tangent to the fibers \(\Bun_{m_0}\) (the \emph{vertical vectors}).\define{vertical!vector}
A \emph{connection}\define{connection} on a left principal \(G\)-bundle \(\Bun\to M\) is a \(1\)-form \(\rho_{\Bun}\) on \(\Bun\), valued in \(\LieG\), equal to the right invariant Maurer--Cartan form on the fibers, transforming in the adjoint representation under the left \(G\)-action:
\[
\LT{g}^*\omega_{\Bun}=\Ad_g\omega_{\Bun}.
\]
\end{answer}
\begin{example}
If \(\Bun\to M\) is a trivial right principal bundle, \(\Bun=M\times G\),  \(\omega_{\Bun}=\omega_G\) is the \emph{standard flat connection}.\define{connection!standard flat}\define{standard flat connection}
\end{example}
\begin{problem}{std.flat.aut}
What are the automorphisms of a trivial principal bundle?
Which preserve the standard flat connection?
Prove that the automorphisms preserving the standard flat connection act transitively.
If the base manifold has positive dimension, prove that the automorphism group contains abelian Lie groups of all dimensions as subgroups.
\end{problem}
A connection is \emph{flat}\define{connection!flat}\define{flat connection} if it is locally isomorphic to the standard flat connection.
\begin{example}
Let us consider what connections look like in local coordinates.
Since we work locally, we can suppose that \(\Bun\to M\) is a trivial right principal bundle, \(\Bun=M\times G\).
We can further suppose that \(M\) is an open set in Euclidean space \(\R^n\).
Take any connection \(\omega\) on \(\Bun\), and let \(\omega_G\) be the standard flat connection.
For simplicity of notation, denote \(\omega_G\) also as \(g^{-1}\,dg\).
Returning to the definition of connection, \(\omega-\omega_G\) vanishes on the fibers of \(\Bun\to M\) , i.e. each fiber \(\set{m_0}\times G\).
Both \(\omega\) and \(\omega_G\) transform in the adjoint representation under right \(G\)-action, hence their difference does:
\[
r_g^*(\omega-\omega_G)=\Ad_g^{-1}(\omega-\omega_G).
\]
Therefore the differential form
\[
\eta:=\Ad_g(\omega-\omega_G)
\]
is \(G\)-invariant and vanishes on the fibers.
So on \(M\times G\), in any local coordinates \(x^1,\dots,x^n\) on \(M\), \(\eta\) is a linear combination of \(dx^1,\dots,dx^n\).
The \(1\)-form \(\eta\) is \(G\)-invariant, so its coefficients don't depend on \(g\in G\), i.e. \(\eta\) is a differential form on \(M\):
\[
\eta=\Gamma_i(x)dx^i.
\]
Solve for \(\omega\):
\[
\omega=g^{-1}\,dg+\Ad_g^{-1}\Gamma_i(x)dx^i,
\]
where \(M\xrightarrow{\Gamma_i}\LieG\) are smooth functions.
\end{example}
\begin{problem}{construct.conn}
Prove that, for any open set \(M\) in Euclidean space, and for any smooth functions \(M\xrightarrow{\Gamma_i}\LieG\), this expression
\[
\omega=g^{-1}\,dg+\Ad_g^{-1}\Gamma_i(x)dx^i,
\]
is a connection on the trivial bundle over \(M\).
\end{problem}
So we have the general description of a connection in local coordinates.
\begin{problem}{left.conn}
What is the analogous expression of connections on left principal bundles?
\end{problem}
\begin{problem}{af}
Prove that the bundle automorphisms of the trivial right principal bundle are precisely the maps
\[
(m,g)\mapsto(m,h(m)g)
\]
where \(M\xrightarrow{h}G\) is any smooth map.
Prove that the standard flat connection then transforms from \(\omega_G\) to
\[
\omega_G+\Ad_g^{-1}h^*\omega_G.
\]
Prove that the automorphisms that preserve the standard flat connection are precisely those with \(h\) locally constant.
\end{problem}
\begin{problem}{ag}
Use a partition of unity subordinate to a cover by open sets on which the bundle is trivial to prove that every principal bundle admits a connection, and that this connection is unique up to adding a section of \(T^*M \otimes (\amal{\Bun}{G}{\LieG})\).
\end{problem}
\begin{answer}{ag}
Take a principal bundle \(G\to\Bun\to M\).
Cover \(M\) with open sets \(M_a\subseteq M\) and take \(\Bun_a\subseteq\Bun\) to be the preimage of \(M_a\).
Suppose that each \(G\to\Bun_a\to M_a\) admits a connection \(\omega_a\).
Take a partition of unity \(\set{f_a}\) subordinate to the open cover \(\set{M_a}\).
Let \(\eta_a:=f_a\omega_a\) on \(\Bun_a\).
Extend \(\eta_a\) to a \(1\)-form on \(\Bun\) valued in \(\LieG\), by asking it to be zero outside \(\Bun_a\).
Note \(\eta_a\) is a smooth \(1\)-form, since the support of \(f_a\) is compact and contained in \(M_a\).
Moreover, \(\eta_a\) is \(G\)-equivariant.
Let \(\omega=\sum_a\eta_a\).
So \(\omega\) is \(G\)-equivariant, and agrees with the Maurer--Cartan form on each fiber as the \(f_a\) sum to \(1\).
So \(\omega\) is a connection on \(\Bun\).
\end{answer}
\section{Horizontal and vertical}
Take a connection \(\omega_{\Bun}\) on a bundle \(\Bun\to M\).
Note that the action of \(G\) on \(\Bun\) moves points of \(\Bun\) with vertical velocity: for a right principal \(G\)-bundle, we let
\[
A_{\Bun}(p):=\left.\frac{d}{dt}\right|_{t=0}pe^{tA},
\]
giving a linear injection
\[
A\in\LieG\mapsto A_{\Bun}(p)
\]
with image the vertical vectors at \(p\in\Bun\).
A \emph{horizontal vector}\define{horizontal!vector} of \(\omega_{\Bun}\) is a tangent vector \(v\) to \(\Bun\) on which \(0=v\hook\omega_{\Bun}\).
The \emph{horizontal space}\define{horizontal!space} (and the \emph{vertical space})\define{vertical!space} at a point \(p\in\Bun\) is the set of horizontal vectors (vertical vectors) in \(T_p\Bun\).
These linear subspaces are carried invariantly by \(G\)-action, since \(\omega_{\Bun}\) transforms by the adjoint representation.
The horizontal space is a complement to the vertical space, since \(\omega_{\Bun}\) at each point \(p\in\Bun\) maps isomorphically the vertical vectors to the Lie algebra.
The bundle map \(\Bun\xrightarrow{\pi}M\) therefore has derivative \(\pi'(p)\) identifying the horizontal space at any point \(p\in\Bun\) with the tangent space \(T_m M\) at the underlying point \(m:=\pi(p)\).
The inverse of that linear isomorphism takes each tangent vector in \(T_m M\) to its unique horizontal lift,\define{horizontal lift}\define{lift!horizontal} the horizontal tangent vector in \(T_p\Bun\) which projects to it.
The horizontal lifting is \(G\)-invariant.

\section{Lifting curves}
Take a principal bundle \(G\to\Bun \to M\) with a connection.
A submanifold of \(\Bun\) is \emph{horizontal}\define{horizontal!submanifold} if its tangent vectors are horizontal.
\begin{theorem}\label{theorem:horizontal.lift}
Take a principal bundle \(\Bun\xrightarrow{\pi}M\) with a connection \(\omega_{\Bun}\).
Pick points \(p_0\in\Bun\) and \(m_0\in M\), so that \(m_0:=\pi(p_0)\).
Take a locally absolutely continuous curve \(m(t)\) in \(M\), defined on an interval of the real number line,  with \(m(t_0)=m_0\) at some time \(t_0\) in that interval.
Then the curve has a unique horizontal lift \(p(t)\in\Bun\), i.e. a locally absolutely continuous curve with \(\pi(p(t))=m(t)\), so that \(p(t_0)=p_0\) and \(\dot{p}\hook\omega_{\Bun}=0\).
\end{theorem}
\begin{proof}
Suppose first that \(\Bun\to M\) is trivial: \(\Bun=M\times G\).
As above, define a \(1\)-form \(\eta\) by
\[
\eta:=\Ad_g(\omega_{\Bun}-\omega_G), 
\]
so we can write the connection as
\[
\omega_{\Bun}=\omega_G+\Ad_g^{-1}\eta.
\]
Check that \(\eta\) vanishes on vertical vectors for \(\Bun\to M\) and that \(\eta\) is \(G\)-invariant. 
By theorem~\vref{theorem:basic}, \(\eta\) is the pullback of a unique \(1\)-form, which we also denote by \(\eta\), on \(M\).
 
The horizontal vectors are those on which \(\omega_{\Bun}=0\), i.e. on which
\[
\omega_G=-\Ad_g^{-1}\eta.
\]
Take the map \(g\mapsto g^{-1}\), applied to \(\Bun=M\times G\) as \((p,g)\mapsto (p,g^{-1})\).
By lemma~\vref{lemma:invert.Maurer.Cartan}, under this map \(\omega_G\) becomes \(-\Ad_g\omega_G\).
So our equation becomes
\[
\omega_G=\eta,
\]
a Lie equation along any curve.
The lift exists and is unique by theorem~\vref{theorem:Lie.equations}.
Local existence and uniqueness implies global existence and uniqueness.
\end{proof}

\section{Horizontal bundle}
The union of the vertical spaces of a principal bundle is the \emph{vertical bundle}.
The union of the horizontal spaces is the \emph{horizontal bundle} of the connection.
\begin{lemma}\label{lemma:horizontal.bundle}
Take a principal bundle \(G\to\Bun\to M\).
The horizontal bundle of every connection is a \(G\)-invariant subundle of \(T\Bun\) complementary to the vertical bundle.
Every \(G\)-invariant subundle of the \(T\Bun\) complementary to the vertical bundle
is the horizontal bundle of a unique connection.
\end{lemma}
\begin{problem}{connection.bundle}
Prove it.
\end{problem}
\begin{answer}{connection.bundle}
Take a connection \(\omega\); it transforms in the adjoint representation, so the subspace \((0=\omega)\) is carried to itself.
In local coordinates on \(M\), 
\[
\omega=g^{-1}\,dg+\Ad_g^{-1}\Gamma_i(x)dx^i.
\]
The horizontal bundle is the set of tangent vectors \((\dot{x},A_G(g))\) so that
\[
0=\omega(
	\dot{x},
	A_G(g)
	)=A+\Ad_g^{-1}\Gamma_i(x)\dot{x}^i,
\]
which is a linear equation determining \(A\) from \(\dot{x}\), smoothly in parameters \((x,g)\).
So the horizontal bundle is a smooth vector subbundle of \(T\Bun\).

On the other hand, take some \(G\)-invariant vector subbundle \(H\subseteq T\Bun\) complementary to the vertical.
In local coordinates as above, since \(H\) is complementary to the vertical, it is given at each point \(p\in\Bun\) by some linear equation which determines the vertical part of a vector \((\dot{x},A_G)\in H_p\) in terms of the part \(\dot{x}\), which can be chosen arbitrarily.
So at each point \((x,1)\), let \(\Gamma_i(x)\) be defined by
\[
(\partial_i,-\Gamma_i(x))\in H,
\]
where 
\[
\partial_i:=\frac{\partial}{\partial x^i}.
\]
Then let
\[
\omega(\dot{x},A_G):=A+\Ad_g^{-1}\Gamma_i(x)dx^i.
\]
By the existence and uniqueness of \(\omega\) with the given horizontal space at each point \((x,g)\), and by \(G\)-equivariance, \(\omega\) is determined uniquely locally and hence globally.
\end{answer}
\begin{theorem}\label{theorem:induced.connection}
Take a Lie group morphism \(G\to G'\) and a right principal \(G\)-bundle \(G\to\Bun\to M\).
Under the bundle morphism \(\Bun\to\Bun'\) to the associated \(G'\)-bundle, every connection on \(\Bun\) has horizontal spaces mapped to the horizontal spaces of a unique connection on \(\Bun'\), the \emph{induced connection}.\define{connection!induced}\define{induced!connection}
\end{theorem}
\begin{problem}{induced.connection}
Prove it.
\end{problem}
\begin{answer}{induced.connection}
Geometrically: the map \(\Bun\to\Bun'\) projects to the identity map \(M\to M\), 
so the bundle morphism takes tangent vectors to \(\Bun\) with nonzero projection to \(TM\) to tangent vectors to \(\Bun'\) with nonzero projection to \(TM\).
So the horizontal bundle of any connection on \(\Bun\) maps to an invariant subbundle, complementary to the vertical, on \(\Bun'\).
By lemma~\vref{lemma:horizontal.bundle}, this is the horizontal bundle of a unique connection.
We give a more algebraic proof, directly in terms of connections.
Recall that \(\Bun'=\amal{\Bun}{G}{G'}\) is the quotient of \(\Bun\times G'\) by the right \(G\)-action
\[
(p,g')g=(pg,g^{-1}g').
\]
Recall that on \(\Bun\), 
\[
A_{\Bun}(p):=\left.\frac{d}{dt}\right|_{t=0}pe^{tA}
\]
for \(A\in\LieG\).
The tangent space to the fiber at \((p,g)\) is therefore the set of vectors of the form
\[
(A_{\Bun},-\RT{g*}A).
\]
Denote the connection on \(\Bun\) by \(\omega\).
On \(\Bun\times  G'\), take the \(1\)-form
\[
\omega'_{(p,g')}:=\omega_{G,g'}+\Ad_{g'}^{-1}\rho\omega_p.
\]
This form vanishes precisely on tangent vectors of the form \((u,w)\) where \(u\in T_p\Bun\) and \(w\hook\omega_{G',g'}=-\Ad_{g'}^{-1}\rho\omega(u)\).
So if \(A:=\omega(u)\) and \(B:=w\hook\omega_G\), 
\[
B=-\Ad_{g'}^{-1}A=-\LT{g'*}^{-1}\RT{g'*}A
\]
so that \(w=\LT{g'*}B=-\RT{g'*}A\).
Hence
\[
(u,w)=(\horizontalPart{u},0)+(A_{\Bun},-\RT{g'*}A)
\]
is the sum of a horizontal vector for the connection and a tangent vector to the fibers of the map to \(\Bun_G\).
Hence \(\omega'\) descends to \(\Bun'\), and vanishes on tangent vectors to \(\Bun'\) precisely when they are the images of horizontal vectors for \(\Bun\).
On \(\Bun\times G'\), on tangent vector fields generating the right \(G'\)-action, \(\omega'=\omega_{G'}\), so the same holds on \(\Bun'\).
Hence \(\omega'\) is a connection form on \(\Bun'\) with horizontal space precisely the image of the horizontal space of \(\omega\) on \(\Bun\).
If there is another such connection \(\omega''\) on \(\Bun'\), the difference \(\omega'-\omega''\) vanishes on the vertical and the horizontal vectors, so vanishes.
\end{answer}

\section{Quotients}
Many important examples in differential geometry can be expressed as abstract quotients of simple explicit examples.
\begin{problem}{quotient.conn}
Suppose that \(G\to\Bun\to M\) is a right principal bundle.
Suppose that \(H\) is a Lie group acting smoothly on \(\Bun\) and on \(M\), commuting with the map \(\Bun\to M\) and with the \(G\)-action.
Suppose that \(H\) acts freely and properly on \(M\).
Prove \(\bar\Bun:=H\backslash\Bun\to\bar{M}:=H\backslash M\) is a smooth principal right \(G\)-bundle, for unique smooth structures so that the maps
\[
\begin{tikzcd}
\Bun\arrow[r]\arrow[d]&\bar\Bun\arrow[d]\\
M\arrow[r]&\bar{M}
\end{tikzcd}
\]
are principal bundle maps.
Prove that a connection on \(G\to\Bun\to M\) arises from a connection on \(G\to\bar\Bun\to\bar{M}\) if and only if both (i) it is \(H\)-invariant and (ii) 
the \(H\)-fibers in \(\Bun\) are horizontal.
When this happens, prove that this connection on \(G\to\bar\Bun\to\bar{M}\) is unique.
\end{problem}
\begin{answer}{quotient.conn}
All but the last statement is just problem~\vref{problem:biquot}.
Clearly any connection on \(\bar\Bun\to\bar{M}\) pulls back to an invariant differential form on \(\Bun\), equal to the Maurer--Cartan form on each fiber, hence a connection.
But being pulled back from \(\bar\Bun\), it vanishes on the fibers of \(\Bun\to\bar\Bun\).

Any invariant connection on \(\Bun\to M\) vanishing on the fibers of \(\Bun\to\bar\Bun\) is locally pulled back from a form on \(\bar\Bun\) by theorem~\vref{theorem:basic}.
This form agrees with the Maurer--Cartan form on the fibers of \(\bar\Bun\to\bar{M}\), because they are matched up with the fibers of \(\Bun\to M\).

The quotient map takes each horizontal space to a linear subspace in the tangent space of \(\bar\Bun\).
This linear subspace is defined independently of choice of point of \(\Bun\), by invariance.
\end{answer}
\begin{example}
Suppose that \(H\to\Bun\to M\) is a principal bundle and that \(\Gamma\) is a group of diffeomorphisms of \(\Bun\), commuting with \(H\), hence acting on \(M\).
Let \(M'\subseteq M\) be the free regular set of the \(\Gamma\)-action, \(\Bun'\subseteq\Bun\) its preimage.
Suppose that \(M'\) is not empty.
Then \(\Bun'\) lies in the free regular set of the \(\Gamma\)-action on \(\Bun\).
Let \(\bar\Bun:=\Gamma\backslash\Bun'\) and \(\bar{M}:=\Gamma\backslash M'\).
From problem~\vref{problem:quotient.conn}, there is a unique smooth structure and \(H\)-action on \(\bar\Bun\) for which \(\Bun'\to\bar\Bun\) is a smooth \(H\)-equivariant covering and \(H\to\bar\Bun\to\bar{M}\) is a smooth principal bundle.
Every \(\Gamma\)-invariant connection on \(\Bun\) is pulled back from a unique connection on \(\bar\Bun\), and every connection on \(\bar\Bun\) pulls back to a \(\Gamma\)-invariant connection on \(\Bun\).
\end{example}
\section{Example: the hyperbolic plane}
The \emph{hyperbolic plane} \(X\) is the upper half plane, i.e. the set of complex numbers \(z=x+iy\) with \(y>0\).
It is acted on by the group \(G=\PSL{2,\R}\) whose elements are matrices
\[
g=
\begin{bmatrix}
a&b\\
c&d
\end{bmatrix}
\]
with \(ad-bc=1\), defined up to \(\pm\).
Each element
\[
g=
\begin{bmatrix}
1&b\\
0&1
\end{bmatrix}
\]
translates \(z=x+iy\mapsto (x+b)+iy\).
Each element
\[
g=
\begin{bmatrix}
a&0\\
0&a^{-1}
\end{bmatrix}
\]
rescales \(z\mapsto a^2z\).
So \(G\) acts transitively on \(X\).
The stabilizer of \(z=i\) is the group of rotations
\[
h=
\begin{bmatrix}
a&b\\
-b&a
\end{bmatrix}
\]
with \(a^2+b^2=1\).
It acts on the tangent space \(T_i X\) by differentiating:
\begin{align*}
\begin{bmatrix}
a&b\\
-b&a
\end{bmatrix}
(i+\dot{z})
&=
\frac{a(i+\dot{z})+b}{-b(i+\dot{z})+a},
\\
&=
\frac{ai+b+a\dot{z}}{a-bi-b\dot{z}},
\\
&=
\frac{(a+bi)(ai+b+a\dot{z})}{(a+bi)(a-bi-b\dot{z})},
\\
&=
\frac{i+(a+bi)a\dot{z}}{1-(a+bi)b\dot{z}},
\\
&=
(i+(a+bi)a\dot{z})(1+(a+bi)b\dot{z}+\dots),
\\
&=
i+(a+bi)^2\dot{z}+\dots,
\end{align*}
So the action of \(H\) on \(T_i X\) is
\[
\begin{bmatrix}
a&b\\
-b&a
\end{bmatrix}
\dot{z}
=
(a+bi)^2\dot{z}.
\]
Here \(a+bi\) can be any unit complex number.
The stabilizer thus preserves precisely (i) the orientation and (ii) the Euclidean metric on the tangent space \(T_i X\).
Move that metric and orientation around \(X\) using \(G\): the hyperbolic plane \(X\) is an oriented Riemannian surface and \(G\) is a group of orientation preserving isometries.
Since \(G\) acts transitively on \(X\), this metric is complete with constant curvature.
Note that \(H\) acts transitively, strongly and effectively on the oriented orthonormal bases of \(T_i X\) in that \(H\)-invariant metric.
So \(G\) acts transitively, strongly and effectively on the oriented orthonormal bases of \(X\) in that \(G\)-invariant metric.
Hence \(G\to X\) is the bundle of oriented orthonormal bases of the hyperbolic plane. 
The Lie algebra \(\LieH\) of \(H\) is the span of the matrix
\[
\begin{pmatrix}
0&1\\
-1&0
\end{pmatrix}.
\]
In \(G\), note that \(ad-bc=1\) forces
\[
a\,dd+d\,da-b\,dc-c\,db=0,
\]
which we write as
\[
a\,dd-c\,db=-(d\,da-b\,dc).
\]
The Maurer--Cartan form is
\begin{align*}
\omega&:=
g^{-1}dg,\\
&=
\begin{pmatrix}
d&-b\\
-c&a
\end{pmatrix}
\begin{pmatrix}
da&db\\
dc&dd
\end{pmatrix},
\\
&=
\begin{pmatrix}
d\,da-b\,dc&d\,db-b\,dd\\
-c\,da+a\,dc&-c\,db+a\,dd
\end{pmatrix},
\\
&=
\begin{pmatrix}
d\,da-b\,dc&d\,db-b\,dd\\
-c\,da+a\,dc&-(d\,da-b\,dc)
\end{pmatrix}
\end{align*}
On \(H\), we have \(c=-b\) and \(d=a\), so this becomes
\[
h^{-1}dh
=
(a\,db-b\,da)
\begin{pmatrix}
0&1\\
-1&0
\end{pmatrix}.
\]

The transpose operation on \(\LieG=\LieSL{2,\R}\) is \(H\)-invariant.
So we can \(H\)-equivariantly split any element \(A\in\LieG\) into its antisymmetric and symmetric parts:
\[
A=K+P
\]
where \(K:=(A-A^t)/2\in\LieH\) and \(P:=(A+A^t)/2\).

To get a connection on the oriented orthonormal frame bundle of the hyperbolic plane, we need an \(H\)-equivariant \(1\)-form on \(G\) valued in \(\LieH\) equal to the Maurer--Cartan form of \(H\) on fibers of \(G\to X\).
Apply to \(\omega=g^{-1}\,dg\) the \(H\)-invariant projection \(\LieG\to\LieH\), i.e. the \(K\)-projection:
\[
\frac{\gamma}{2}
\begin{pmatrix}
0&1\\
-1&0
\end{pmatrix}
\]
where
\[
\gamma=d\,db-b\,dd+c\,da-a\,dc.
\]
We will see that this connection is the Levi--Civita connection.

The \emph{soldering form} is the quotient \(g^{-1}dg+\LieH\), which we can identify with the \(P\)-projection of \(\omega=g^{-1}\,dg\):
\[
\frac{1}{2}
\begin{pmatrix}
\sigma_1&\sigma_2\\
\sigma_2&-\sigma_1
\end{pmatrix}
\]
where
\begin{align*}
\sigma_1&=2(d\,da-b\,dc),\\
\sigma_2&=d\,db-b\,dd+a\,dc-c\,da\\
\end{align*}
Denoting \(g^{-1}\,dg\) as \(\omega\),
\begin{align*}
d\omega
&=-\omega\wedge\omega,
\\
&=
\frac{1}{2}
d
\begin{pmatrix}
\sigma_1&\sigma_2+\gamma\\
\sigma_2-\gamma&-\sigma_1
\end{pmatrix},
\\
&=
-
\frac{1}{4}
\begin{pmatrix}
\sigma_1&\sigma_2+\gamma\\
\sigma_2-\gamma&-\sigma_1
\end{pmatrix}
\wedge
\begin{pmatrix}
\sigma_1&\sigma_2+\gamma\\
\sigma_2-\gamma&-\sigma_1
\end{pmatrix},
\end{align*}
so that
\begin{align*}
d\sigma_1&=-\gamma\wedge\sigma_2,\\
d\sigma_2&=\gamma\wedge\sigma_1,\\
d\gamma&=-\sigma_1\wedge\sigma_2,
\end{align*}
which we can think of as the Maurer--Cartan structure equations of \(G=\PSL{2,\R}\), or as the structure equations of the hyperbolic plane.
It is convenient to let \(\sigma:=\sigma_1+i\sigma_2\), giving structure equations
\begin{align*}
d\sigma&=i\gamma\wedge\sigma,\\
d\gamma&=-\frac{i}{2}\sigma\wedge\bar\sigma.
\end{align*}

The Killing form on \(\LieG=\LieSL{2,\R}\) is \(4\) times trace squared:
\begin{align*}
\operatorname{tr}
\begin{pmatrix}
\sigma_1&\sigma_2+\gamma\\
\sigma_2-\gamma&-\sigma_1
\end{pmatrix}
\begin{pmatrix}
\sigma_1&\sigma_2+\gamma\\
\sigma_2-\gamma&-\sigma_1
\end{pmatrix}
&=
\operatorname{tr}
\begin{pmatrix}
\sigma_1^2+\sigma_2^2-\gamma^2&0\\
0&\sigma_1^2+\sigma_2^2-\gamma^2
\end{pmatrix},
\\
&=
2(\sigma_1^2+\sigma_2^2-\gamma^2).
\end{align*}
So on \(H\), \(-2\gamma^2\) is the Killing form, while on the horizontal space, the Killing form is \(2(\sigma_1^2+\sigma_2^2)\).
In particular, \(\sigma_1^2+\sigma_2^2\) is an invariant quadratic form defined on \(G\), vanishing precisely on the fibers of \(G\to X\), hence descending to a \(G\)-invariant Riemannian metric on the hyperbolic plane \(X\).
Note that the \(G\)-invariant metric on \(X\) is unique, up to a positive constant scaling.
We will see that this metric has curvature \(-1\).

Consider the quotient \(\bar{X}:=\Gamma\backslash X'\) of the free regular set \(X'\subseteq X\) of a subgroup \(\Gamma\subset G=\PSL{2,\R}\).
The associated bundle is \(\bar\G:=\Gamma\backslash G'\), where \(G'\subseteq G\) is the preimage of \(X'\subseteq X\), equipped with the form \(\omega\) above.
By the lemma, the connection \(\gamma\) descends.
In fact, this is already clear in this example since \(\gamma\) is the invariantly defined \(K\)-part of \(\omega\), and \(\omega\) itself descends by left invariance.

We will see that every oriented complete surface of constant curvature \(-1\) is covered by the hyperbolic plane \(X\), and hence arises as such a quotient \(\bar{X}=\Gamma\backslash X\), uniquely up to \(G\)-conjugation of \(\Gamma\), where \(\Gamma\subset G\) is a discrete group acting freely and properly.
Hence the Levi--Civita connection of any such surface pulls back to \(\gamma\).
This concrete computation therefore provides an explicit description of the pullback \((X,G)\)-structure and the Levi--Civita connection of any complete oriented surface of constant curvature \(-1\).

\section{Connections as splittings}
Suppose that \(G\to\Bun\xrightarrow{\pi}M\) is a principal right \(G\)-bundle.
\begin{problem}{Atiyah.class}
Prove that the vertical bundle is trivial as a bundle over \(\Bun\), isomorphic to the trivial bundle \(\underline{\LieG}:=\Bun\times\LieG\).
\end{problem}
\begin{answer}{Atiyah.class}
Map
\[
(p,A)\in\Bun\times\LieG\xmapsto{\Phi}A_{\Bun}(p)\in T_p\Bun.
\]
Locally trivializing \(\Bun\to M\), \(A_{\Bun}\) becomes the left invariant vector field \(A_G\) on \(G\) at each fiber.
Hence \(\Phi\) is a linear isomorphism to the vertical bundle \(\ker\pi'\).
\end{answer}
\begin{problem}{Atiyah.class.2}
Prove that the sequence of vector bundles
\[
\begin{tikzcd}
0\arrow[r]&\underline{\LieG}\arrow[r]\arrow[dr]&T\Bun\arrow[r]\arrow[d]&\pi^*TM\arrow[r]\arrow[dl]&0\\
&&\Bun&&
\end{tikzcd}
\]
is exact, where the map \(T\Bun\to\pi^*TM\) maps each tangent vector to \(\Bun\) to its image tangent vector by \(\pi_*\).
\end{problem}
\begin{answer}{Atiyah.class.2}
Applying the previous problem, \(\ker\pi'\cong\underline{\LieG}\).
We have an obvious exact sequence on each fiber: if \(m=\pi(p)\),
\[
\begin{tikzcd}
0\arrow[r]&\ker\pi'(p)\arrow[r]&T_p\Bun\arrow[r,"\pi'"]&(\pi^TM)_p=T_m M\to 0
\end{tikzcd}
\]
\end{answer}
\begin{problem}{Atiyah.class.3}
Prove that every connection splits the exact sequence: the horizontal space sits inside \(T\Bun\) as a complement to the vertical, hence isomorphic to \(\pi^*TM\).
\end{problem}
\begin{answer}{Atiyah.class.3}
Take a connection \(\omega\).
The map
\[
v\in T_p\Bun\mapsto (\omega(v),\pi'(p)v)
\]
is (again by looking at a local trivialization) an isomorphism of vector bundles
\[
T\Bun\to\underline{\LieG}\oplus\pi^*TM,
\]
so we split the exact sequence by composing the inverse of this isomorphism with the inclusion
\[
\pi^*TM\to\underline{\LieG}\oplus\pi^*TM.
\]
\end{answer}
\begin{problem}{Atiyah.class.4}
Prove that every \(G\)-invariant splitting arises uniquely from a connection.
\end{problem}
\begin{answer}{Atiyah.class.4}
Given a \(G\)-equivariant splitting
\[
\pi^*TM\xrightarrow{s}T\Bun,
\]
each vector \(v\in T_p \Bun\) has associated horizontal lift 
\[
\bar{v}:=s(\pi'(p)v).
\]
Hence \(v-\bar{v}\) is vertical, so agrees with \(A_{\Bun}(p)\) for precisely one \(A\in\LieG\).
We define a connection \(\omega\) by
\[
\omega(v):=A.
\]
\end{answer}
\begin{problem}{Atiyah.class.5}
Prove that the quotient of the vertical bundle by the \(G\)-action is a vector bundle over \(M\), the \emph{adjoint bundle}\define{adjoint bundle} \(\ad_{\Bun}\). Hint: problem~\vref{problem:quotient.v.b}.
\end{problem}
\begin{answer}{Atiyah.class.5}
The action is free and proper by vector bundle isomorphisms.
We can see all of this easily in a local trivialization.
\end{answer}
\begin{problem}{Atiyah.class.6}
Prove that the quotient of \(T\Bun\) by the \(G\)-action is a vector bundle over \(M\); denote it \(T\Bun/G\).
\end{problem}
\begin{answer}{Atiyah.class.6}
The action is free and proper by vector bundle isomorphisms.
We can see all of this easily in a local trivialization.
\end{answer}
\begin{problem}{Atiyah.class.7}
Prove that the exact sequence quotients by \(G\)-action to an exact sequence
\[
\begin{tikzcd}
0\arrow[r]&\ad_{\Bun}\arrow[r]\arrow[dr]&T\Bun/G\arrow[r]\arrow[d]&TM\arrow[r]\arrow[dl]&0\\
&&M&&
\end{tikzcd}
\]
\end{problem}
\begin{answer}{Atiyah.class.7}
Again, by local triviality, we get this sequence is exact in each fiber over open subsets of \(M\) on which \(\Bun\to M\) is locally trivial, so it is enough to suppose \(\Bun\to M\) is locally trivial, and then it is easy to check.
\end{answer}
\begin{problem}{Atiyah.class.8}
Prove that every connection also splits this exact sequence, by the image of the horizontal bundle.
\end{problem}
\begin{answer}{Atiyah.class.8}
Again check by local triviality.
\end{answer}
\begin{problem}{Atiyah.class.9}
Conversely, prove that every splitting of this sequence of vector bundles on \(M\) arises uniquely from a connection on \(\Bun\).
\end{problem}
\begin{answer}{Atiyah.class.9}
Again check by local triviality.
\end{answer}
\begin{example}
Take complex Euclidean space \(V:=\C^{n+1}\) for some \(n\ge 0\), and let \(M\) be the associated complex projective space \(M:=\Proj{n}\).
Let \(\Bun\subset V\) be a compact connected real smooth hypersurface, transverse to every real line through the origin (for example: the unit sphere).
Map \(\Bun\xrightarrow{\pi}M\) by taking each point \(p\in\Bun\) to the line \(m:=\C p\) that it spans.
The map is a submersion by transversality.
The preimage of a point of \(M\) is the intersection of \(\Bun\) with a complex line through the origin, a smooth curve by transversality, compact since \(\Bun\) is compact.
Since the fibers of the map are all compact, the map is a fiber bundle map, by corollary~\vref{corollary:Ehresmann.thm}.
Since complex projective space is connected and simply connected, the exact sequence in homotopy tells us that the fibers are connected, so a fiber bundle \(G\to\Bun\to M\) where \(G=S^1\).
Each fiber of \(\Bun\to M\) is a connected compact curve (hence diffeomorphic to a circle) in a complex line through the origin, transverse to the real rays through the origin, hence (by connectivity and winding number) containing a unique point on each real ray.
We get the unit circle \(G:=S^1\subset\C\) to act on \(\Bun\) by usual action scaling complex vectors, and hence acting on the set of real rays.
Since the \(G\)-action is free, and \(G\) is compact, \(G\to\Bun\to M\) is a principal right \(G\)-bundle.
Each subspace \(H_p:=(ip)^{\perp}\cap T_p\Bun\subset V\) is a complement to the \(G\)-orbit through \(p\), smoothly varying with \(p\), so a connection bundle, so the horizontal spaces of a unique connection.
In the special case when \(\Bun\) is the unit sphere, this is the Hopf fibration, with its standard connection, invariant under the unitary group.
\end{example}
\section{Connections on homogeneous principal bundles}
Take Lie groups \(G,F,H\) with \(H\subseteq G\) closed.
To each a Lie group morphism \(h\in H\to\bar{h}\in F\), we have an associated principal \(F\)-bundle \(\Bun:=\amal{G}{H}{F}\).
A \emph{die}\define{die} is an \(H\)-equivariant linear map \(A\in\LieG\to\bar{A}\in\LieF\) extending the associated Lie algebra morphism \(\LieH\to\LieF\).
Consider the following recipe: each die has associated \(1\)-form
\[
\omega:=\omega_F+\Ad_f^{-1}\bar\omega_G
\]
on \(G\times F\).
We conjugate the die by replacing the Lie group morphism by
\[
h\mapsto\Ad_{f_0} \bar{h},
\]
and the die \(\bar{\,}\) by \(\Ad_{f_0}\bar{\,}\), for some \(f_0\in F\).
\begin{theorem}[Wang \cite{Wang1958}]\label{theorem:invariant.connection}
Take a homogeneous right principal bundle \(\Bun\to X\).
As~\vpageref{sec:homog.principal.bundles}, write it as \(\Bun=\amal{G}{H}{F}\), \(X=G/H\) for Lie groups \(G,F,H\) with \(H\subseteq G\) closed and a Lie group morphism \(H\to F\).
The recipe above associates to each die a \(1\)-form.
There is a unique \(G\)-invariant connection on \(\Bun\to X\) which pulls back to that \(1\)-form via
\[
G\times F\to\amal{G}{H}{F}=\Bun.
\]
Conversely, every invariant connection on any homogeneous principal bundle arises by this recipe, uniquely up to conjugating the die.
\end{theorem}
\begin{problem}{invariant.connections}
Prove it.
\end{problem}
\begin{answer}{invariant.connections}
Take a \(G\)-invariant connection \(\omega=\omega_\Bun\) on \(\Bun\).
Pullback by
\[
G\times F\to\Bun=\amal{G}{H}{F}.
\]
Also denote the pullback by \(\omega\).
So \(\omega\) is a \(G\)-invariant, \(H\)-invariant \(1\)-form on \(G\times F\), transforming in the adjoint representation under \(F\)-action.
Take any vector \(v\in T_x (G\times F)\), at some point \(x=(g_0,f_0)\).
We can write that vector as
\[
v=(\LT{g_0*}A,\RT{f_0*}B),
\]
for unique \(A\in\LieG\) and \(B\in\LieF\).
By right \(F\)-equivariance in the adjoint representation and left \(G\)-invariance,
\[
v\hook\omega=\Ad_{f_0}^{-1}((A,0)\hook\omega+B).
\]
Let \(\beta A:=A\hook\omega\), so \(\LieG\xrightarrow{\beta}\LieF\), so 
\[
v\hook\gamma=\Ad_{f_0}^{-1}(B+\beta A).
\]
The \(H\)-action is 
\[
(g_0,f_0)h=(g_0h,a(h)^{-1}f_0).
\]
Consider a one parameter flow in \(H\), say \(e^{tA}\) for some \(A\in\LieH\).
This path in \(H\) moves the point \((g_0,f_0)\) along a path in \(G\times F\):
\[
(g_0e^{tA},e^{-ta A}f_0).
\]
Since it this path arises from a path in \(H\), it lies inside a fiber of
\[
G\times F\to\Bun=\amal{G}{H}{F}.
\]
Since \(\omega\) is pulled back from \(\Bun=\amal{G}{H}{F}\), \(\gamma\) vanishes on the fibers of  \(G\times F\to\Bun\), so vanishes on the velocity vector of this path.
That vector, at \(t=0\), is
\[
v=(\LT{g_0*}A,-\RT{f_0*}a A).
\]
So
\[
0=v\hook\omega=\Ad_{f_0}^{-1}(-a A+\beta A)
\]
vanishes for all \(A\in\LieH\) at every point, i.e. \(\beta=a\) on \(\LieH\).
So we can extend \(a\) to be defined on all of \(\LieG\) by setting \(a:=\beta\).
Under the \(H\)-action, \(\omega\) is invariant, hence \(a\) is invariant.

Conversely, start with any \(H\)-invariant \(\LieG\xrightarrow{\beta}\LieF\) extending \(a\).
Take any point \(x=(g_0,f_0)\in G\times F\) and any vector \(v\in T_x(G\times F)\).
Write \(v\) as
\[
v=(\LT{g_0*}A,\RT{f_0*}B),
\]
for unique \(A\in\LieG\) and \(B\in\LieF\).
Define a \(1\)-form \(\omega\) on \(G\times F\), valued in \(\LieF\), by
\[
v\hook\omega:=\Ad_{f_0}^{-1}(B+\beta A).
\]
Reverse our steps to see that, since \(\beta\) extends \(a\), this form vanishes on the fibers of 
\[
G\times F\to\Bun=\amal{G}{H}{F},
\]
and since \(\beta\) is \(H\)-invariant, so is \(\omega\), so \(\omega\) is semibasic and invariant.
By lemma~\vref{lemma:G.inv.basic} and theorem~\vref{theorem:basic} (using both that \(\omega\) is semibasic and invariant), \(\omega\) is basic, i.e. pulled back from \(\Bun\).
\end{answer}
If we prefer to write the formal notation \(g^{-1}dg\) in place of \(\omega_G\), we can write
\[
\omega_\Bun=f^{-1}df+\Ad_f^{-1}\overline{g^{-1}dg}.
\]
For each \(A\in\LieG\), the vector field
\[
v_A(g,f):=(\LT{g*}A,\RT{f*}B)),
\]
on \(G\times F\), where \(B=-\bar{A}\), projects to a horizontal vector field on \(\Bun=\amal{G}{H}{F}\).
The projection is zero if and only if \(A\in\LieH\).
These vector fields span the horizontal vector fields at every point.
\section{Vector bundle connections}\label{section:vb.connections}
Take a vector bundle \(\vb{V}\to M\) on a manifold \(M\).
A \emph{vector bundle connection}\define{vector!bundle!connection}\define{connection!vector bundle} \(\nabla\) on \(\vb{V}\) is a linear map associating to each smooth section \(s\) of \(\vb{V}\) a smooth section \(\nabla s\) of \(\vb{V}\otimes T^*M\), so that \(\nabla(fs)=f\nabla s+df\,s\) for any smooth function \(f\) on \(M\).

Take a principal bundle \(G\to\Bun\to M\) and a \(G\)-module \(G\xrightarrow{\rho}V\).
Recall that the associated vector bundle \(\vb{V}:=\amal{\Bun}{G}{V}\) has sections identified with the \(G\)-equivariant maps \(\Bun\xrightarrow{s}V\), i.e.
\[
r_g^*s=\rho(g)^{-1}s.
\]
Take a connection \(\omega=\omega_{\Bun}\) on \(\Bun\to M\).
\begin{problem}{compute.nabla}
The \(1\)-form \(\nabla s\) defined by \(ds+\rho(\omega)s\) is a \(1\)-form on \(\Bun\), valued in \(V\), but vanishes on the fibers.
\end{problem}
\begin{answer}{compute.nabla}
For any \(A\in\LieG\),
\begin{align*}
\LieDer_{A_{\Bun}}s
&=
\left.\frac{d}{dt}\right|_{t=0}r_{e^{tA}}^*s,
\\
&=
\left.\frac{d}{dt}\right|_{t=0}\rho(e^{-tA})s,
\\
&=
-\rho(A)s,
\\
&=
-\rho(A_{\Bun}\hook\omega)s.
\end{align*}
\end{answer}
Since \(\nabla s\) vanishes on the fibers, and is \(G\)-equivariant, we can identify it with a section of \(T^*M\otimes\vb{V}\) by problem~\vref{problem:ck}, so \(\nabla s\) is the \emph{covariant derivative}\define{covariant derivative}\define{derivative!covariant} of \(s\).
Any connection on any principal bundle, by this construction, induces a vector bundle connection on every associated vector bundle.
\section{The frame bundle of a vector bundle}
Take a vector space \(V_0\) and a vector bundle \(V\to M\) whose fibers are of the same dimension as \(V_0\).
The \(V_0\)-valued \emph{frame bundle}\define{frame bundle} \(\framebundle{V}\) of \(V\to M\) is the set of pairs \((m,u)\) where \(m\in M\) and \(V_m\xrightarrow{u}\R^N\).
The \emph{projection map}\define{projection map} is
\[
(m,u)\in\framebundle{V}\mapsto m\in M.
\]
The group \(\GL{V_0}\) acts on \(\framebundle{V}\) by
\[
(m,u)g:=(m,g^{-1}u).
\]
\begin{problem}{vb.frame.bundle}
Prove that \(\framebundle{V}\xrightarrow{\pi}M\) is a principal right \(\GL{V_0}\)-bundle.
\end{problem}
Clearly this generalizes the definition of frame bundle of a manifold: our old frame bundle which we denoted \(\framebundle{M}\) is really the frame bundle of \(TM\), so should perhaps be denoted \(\framebundle{TM}\).
We will see in problem~\vref{problem:Koszul.to.Ehresmann.II} that every vector bundle connection on any vector bundle is associated to a unique connection on its frame bundle.
%\begin{problem}{Koszul.to.Ehresmann}
%Prove that a vector bundle connection on a vector bundle is associated to a unique connection on its frame bundle.
%\end{problem}
%\begin{answer}{Koszul.to.Ehresmann}
%Conversely, given a connection \(\nabla\) on a vector bundle, and a tangent vector
%\[
%v\in T_{(m,u)}\framebundle{\vb{V}},
%\]
%pick a basis of local sections \(s_1,\dots,s_N\) defined on some open set of \(M\).
%Each becomes some \(V\)-valued function \(f_1,\dots,f_N\) on \(\framebundle{\vb{V}}\).
%Each has covariant derivative \(\nabla s_i\) a section of \(T^*M\otimes\vb{V}\), so becoming some \(V\)-valued semibasic \(1\)-form \(\nabla f_i\) on \(\framebundle{\vb{V}}\).
%We have to uniquely define \(\omega\) by 
%\[
%\omega f_i = \nabla f_i-df_i.
%\]
%This defines the components of \(\omega\) in that basis, at each point.
%We can change basis, from \(s_i\) to some \(s'_i\) say, and this will change these \(f_i\) to some new \(f'_i\), say.
%But then \(s'_i=g_i^js_j\) for a unique invertible matrix \(g\in\GL{N}\) at each point of \(M\), and hence 
%\[
%\nabla s'_i=g_i^j\nabla s_j+s_j dg^j_i,
%\]
%giving
%\[
%\nabla f'_i=g_i^j\nabla f_j+f_j dg^j_i,
%\]
%and so changes \(\omega\) to \(\omega'\) by
%\begin{align*}
%\omega' f_i'&=
%\nabla f_i'-df_i',
%\\
%&=
%g_i^j\nabla f_j+f_j dg^j_i-g^j_i df_j-f_j dg^j_i,
%\\
%&=
%g_i^j(\nabla f_j-df_j),
%\\
%&=
%\omega g_i^jf_j,
%\\
%&=
%\omega f_i',
%\end{align*}
%i.e. \(\omega'=\omega\).
%So \(\omega\) is well defined.
%\end{answer}

\chapter{Reductive homogeneous spaces}\label{reductive.homogeneous.spaces}
\section{Definition}
A homogeneous space \((X,G)\) with stabilizer \(H:=G^{x_0}\) is \emph{reductive}\define{homogeneous!space!reductive}\define{reductive!homogeneous space} if \(\LieG\) splits as an \(H\)-module, say \(\LieG=\LieH\oplus𝔪\) where \(𝔪\) is some \(H\)-module \cite{Nomizu:1954}.
(Even though \(𝔪\) is \emph{not} necessarily a Lie subalgebra, it is still traditional to write it with a fraktur letter.)
\begin{theorem}
A homogeneous space \((X,G)\) is reductive precisely when the bundle \(G\to X\) has a \(G\)-invariant connection.
Every \(H\)-invariant linear subspace of \(\LieG\) complementary to \(\LieH\) is the horizontal space at \(1\in G\) of a unique \(G\)-invariant connection on the right \(H\)-bundle \(G\to X\).
Conversely, every \(G\)-invariant connection on the bundle \(G\to X\) arises in this way uniquely.
\end{theorem}
\begin{proof}
Take a homogeneous space \((X,G)\), say \(X=G/H\).
Take a \(G\)-invariant connection on the bundle \(G\to X\).
So at each point \(g\in G\), we have some linear subspace \(𝔪_g\subseteq T_g G\), the horizontal space of the connection, which is complementary to the vertical \(T_g (gH)\), and \(G\)-invariance is precisely that \(\LT{g_0*}𝔪_{g_1}=𝔪_{g_0 g_1}\), while being a connection is precisely \(H\)-invariance, i.e. \(\RT{h*}𝔪_g=𝔪_{gh}\). 
So at the point \(1\in G\), this is a linear subspace \(𝔪:=𝔪_1\subseteq \LieG=T_1 G\), complementary to \(\LieH\), and invariant under adjoint \(H\)-action.
Reversing our steps, we get the other direction.
\end{proof}
\begin{theorem}
Take a homogeneous space \((X,G)\), pick a point \(x_0\in X\), and let \(H:=G^{x_0}\subseteq G\).
Any \(G\)-invariant connection on \(g\in G\mapsto gx_0\in X\), with horizontal space \(𝔪\subseteq\LieG\) at \(1\in G\), is flat if and only if \(𝔪\subseteq\LieG\) is an ideal, in which case the immersed Lie subgroup \(N:=e^W\subseteq G\) is a normal subgroup, complementary to \(H\), so that \(n\in N\mapsto nx_0\in X\) is a local diffeomorphism and \(N\)-equivariant, so a covering map to a component of \(X\).
\end{theorem}
\begin{proof}
A connection is flat just when its horizontal spaces are bracket closed.
\end{proof}
\begin{example}
If \(G:=H\ltimes N\) is a semidirect product of Lie groups, then the tangent spaces of the left translates of \(N\) form a flat connection on \(G\to X:=G/H\).
Note that \(N\cong G/H\) is an \(N\)-equivariant diffeomorphism.
\end{example}
\begin{example}
Take \(G:=T^n\) any torus, write \(n=p+q\), \(p,q\ge 1\), let \(H:=T^p\subset G\) any \(p\)-dimensional subtorus, and let \(N:=\R^q\subset G\) be any linearly immersed but not embedded subgroup.
Since \(G\) is abelian, all subgroups of \(G\) are normal.
So \(N\) is a normal subgroup complementary to \(H\) but \(N\to G/H=T^q\) is an infinite covering map, so \(G\) is not a semidirect product.
The \(H\)-invariant linear subspace \(W=N\subseteq\R^{p+q}=\LieG\) is the horizontal space of a unique flat connection.
\end{example}
\begin{problem}{ah}
Prove that any homogeneous space \((X,G)\) is reductive just when \(G\) acts freely on the frame bundle of \(X\).
\end{problem}
In particular, the natural invariant connection of a reductive homogeneous space, associated to the choice of a particular splitting \(\LieG=\LieH\oplus𝔪\), has geodesics precisely the images in the frame bundle of the \(1\)-parameter subgroups generated by elements of \(𝔪\).
\begin{example}
The quotient \(X=G/H\) with \(G:=\SL{n,\R}\) and \(H:=\SO{n}\) has \(H\)-invariant splitting
\[
\LieG=\LieSL{n,\R}=\LieH\oplus𝔪=\LieSO{n}\oplus\Sym^2{\R^n}_0.
\]
The elements of \(𝔪=\Sym^2{\R^n}_0\) are the traceless symmetric matrices.
By Cartan decomposition, every element of \(G\) is uniquely decomposed into a product:
\[
(x,h)\in𝔪\times H\mapsto e^x h\in G=e^{𝔪} H,
\]
is an \(H\)-equivariant diffeomorphism, so \(X\cong𝔪\).
The elements of \(𝔪\), the traceless symmetric matrices, generate the \(1\)-parameter subgroups \(e^{tx}\in X\).
Up to conjugate, we can ensure that \(A\) is diagonal, and the geodesics are obvious: exponential growth or decay of eigenvalues.
\end{example}
\begin{example}
Suppose that \((X,G)\) is a reductive homogeneous space, with stabilizer \(H\subseteq G\), with Lie algebra decomposed \(\LieG=\LieH\oplus𝔪\) as \(H\)-modules. 
Take a Lie group morphism \(H\xrightarrow{a}F\).
The associated Lie algebra morphism \(\LieH\xrightarrow{a}\LieF\) has obvious extension to an \(H\)-equivariant linear map \(\LieG\xrightarrow{a}\LieF\) given by being zero on \(𝔪\).
Hence every homogeneous principal bundle and every homogeneous vector bundle on \(X\) has a canonical invariant connection.
We can split any \(H\)-equivariant linear map \(\LieG\xrightarrow{a}\LieF\) into its components, say \(\LieH\xrightarrow{\alpha}\LieF\) and \(𝔪\xrightarrow{\beta}\LieF\).
But then \(\alpha\) is the Lie algebra morphism induced by the fixed Lie group morphism \(H\to F\).
So the choice of invariant connection on any homogeneous principal bundle is precisely the choice of \(H\)-module morphism \(𝔪\xrightarrow{\beta}\LieF\), so parameterized by the vector space \((𝔪^*\otimes\LieF)^H\).
\end{example}
\begin{problem}{more.conns.reductive.conns}
Given a reductive homogeneous space \((X,G)\), with a fixed splitting \(\LieG=\LieH\oplus\redComplement\) as above, explain how to identify the \(G\)-invariant affine connections with
\[
(\redComplement^*\otimes\redComplement^*\otimes\redComplement)^H.
\]
\end{problem}

\subsection{Affine form}
Given any reductive homogeneous space \((X,G)\), say with splitting \(\LieG=\LieH\oplus𝔪\) as \(H\)-modules, define a new homogeneous space \((X',G')\), the \emph{affine form} of \((X,G)\), with \(X'=𝔪\) and with \(G'=H\ltimes𝔪\), where \(𝔪\) is now an abelian Lie subalgebra.
Every Cartan geometry modelled on \(X=G/H\) is also a Cartan geometry modelled on \(X'=G'/H\), as the Lie algebra \(\LieG'=𝔪\ltimes\LieH=\LieG\) is the same, as is the structure group \(H\).
So the spirals of such a Cartan geometry are modelled on the spirals either of the model \(X=G/H\) or of the affine form \(X'=G'/H=𝔪\).
\begin{example}
Consider the spirals of the affine form.
Denote elements of the group \(G'\) as matrices
\[
g=
\begin{pmatrix}
h&x\\
0&1
\end{pmatrix}
\]
for \(h\in H\) and \(x\in 𝔪=X'\).
The spirals associated to 
\[
B=
\begin{pmatrix}
A&v\\
0&0
\end{pmatrix}
\]
are the integral curves of the differential equation \(\dot{g}=gB\), i.e.
\begin{align*}
\dot{h}&=hA,\\
\dot{x}&=\rho_{𝔪}(h)v.
\end{align*}
So \(h(t)=h_0e^{tA}\) and so we have only to solve
\[
\dot{x}=\rho_{𝔪}(h_0e^{tA})v,
\]
for constant \(h_0,A,v\).
If there is some vector \(w\) so that \(\rho_{𝔪}(A)w=v\), then 
\[
x(t)=x_0+\rho_{𝔪}(h_0e^{tA})w.
\]
(If there is no such \(w\), the integration is more complicated.)
Replacing \(A\) by \(h_0Ah_0^{-1}\) and \(w\) by \(\rho_{𝔪}(h_0)w\), we can arrange that \(h_0=I\) without changing the resulting class of spiral curves in the affine form:
\[
x(t)=x_0+\rho_{𝔪}(e^{tA})w.
\]
The geodesics in the affine form are thus precisely the curves \(x(t)=x_0+tv_0\).

For example, in affine space \(𝔪=ℝ^n\), \(H=\GL{n,ℝ{}}\), and the generic \(A\) is diagonalizable and invertible with some real and some complex eigenvalues, making \(x(t)\) turn around in the directions of complex eigenvalues and decay or grow in the directions of real eigenvalues, hence the name ``spiral''.
\end{example}
\begin{example}
Suppose that \((X,G)\) is a quotient \(X=G/H\) where \(G\) has a Cartan decomposition, in the sense that the map
\[
(h,x)\in H\times 𝔪\mapsto e^xh\in G
\]
is a diffeomorphism, where \(\LieG=\LieH\oplus 𝔪\) as \(H\)-modules.
In particular, \((X,G)\) is a reductive homogeneous space, but also \(X=𝔪\).
Computing the spirals in the affine form as above:
\begin{align*}
\dot{h}&=hA,\\
\dot{x}&=\Ad_h v,
\end{align*}
we easily find that the map
\[
(h,x)\in H\times 𝔪\mapsto e^xh\in G
\]
takes spirals to spirals.
For example, consider the spirals in \(X=G/H=\SL{n,\R}/\SO{n}\).
Each arises from a \(1\)-parameter subgroups \(h=e^{tA}\in H=\SO{n}\); these are easy to classify up to conjugacy.
Then we have to integrate the equation \(\dot{x}=\Ad_h v\) for a symmetric matrix \(x(t)\).
\end{example}
\section{Homogeneous reductive geometries with given symmetries}\label{section:reductive.homogeneous}
As a special case of the construction of homogeneous Cartan geometries (see chapter~\ref{chapter:homogeneous.examples}), suppose that \((X,G)\) is a reductive homogeneous space, say
\[
\LieG=\LieH\oplus𝔪,
\]
and \((X',G')\) is another homogeneous space, and we wish to find all \(G'\)-invariant \((X,G)\)-geometries on \(X'\).
From chapter~\ref{chapter:homogeneous.examples}, we see that the geometry arises from a die, i.e. a Lie group morphism \(H'\to H\) and an \(H'\)-equivariant linear map
\[
A\in\LieG'\mapsto\bar{A}\in\LieG
\]
equal to the induced Lie algebra morphism on \(\LieH'\), and yielding a linear isomorphism
\[
\LieG'/\LieH'=𝔪'\to\LieG/\LieH=𝔪.
\]
\begin{lemma}
Every homogeneous reductive Cartan geometry is on a reductive homogeneous space.
\end{lemma}
\begin{proof}
The splitting 
\[
\LieG=\LieH\oplus𝔪,
\]
splits the die into, say, \(\bar{A}=a_{\LieH}(A)+a_{𝔪}(A)\).
Let \(𝔪'\) be the kernel of \(a_{\LieH}\).
By equivariance of the die, \(𝔪'\) is \(H'\)-invariant, so that \((X',G')\) is also a reductive homogeneous space.
\end{proof}
We now suppose that \((X',G')\) is a reductive homogeneous space with an \(H'\)-equivariant splitting
\[
\LieG'=\LieH'\oplus 𝔪',
\]
fixed independently of the choice of \(G'\)-invariant reductive Cartan geometry.
The Cartan geometry then has a die, but its restriction to \(\LieH'\) is the Lie algebra morphism, so it consists of two \(H'\)-equivariant linear maps
\[
𝔪\xleftarrow{ι}𝔪'\xrightarrow{\Gamma}𝔥,
\]
so that \(ι\) is a linear isomorphism.
Split \(\omega_{G'}\) by the splitting
\[
\LieG'=\LieH'\oplus 𝔪',
\]
as a connection form and a soldering form:
\[
\omega_{G'}=\gamma+\sigma.
\]
In particular, if \((X',G')\) is in affine form, i.e. if \(G'=H'\ltimes{𝔪'}\), say writing each element \(g'\in G'\) as \(g'=(h',x')\), then we can write
\begin{align*}
\gamma&=\omega_{H'},\\
\sigma&=\rho_{𝔪'}(h')^{-1}dx'.
\end{align*}
By the recipe from chapter~\ref{chapter:homogeneous.examples},
\begin{align*}
\omega_{\G}
&=
\Ad_h^{-1}\bar\omega_{G'}+\omega_H,
\\
&=
(\Ad_h^{-1}\Gamma\sigma+\Ad_h^{-1}\bar\gamma+\omega_H,
\rho_{𝔪}(h)^{-1}{ι}\sigma).
\end{align*}
Putting this together, we get a recipe: take
\begin{itemize}
\item
reductive homogeneous spaces \((X,G)\) and \((X',G')\) with \(\dim X=\dim X'\), say \(X=G/H\), \(X'=G'/H'\) and
\item
\(H'\)-invariant and \(H\)-invariant splittings
\begin{align*}
\LieG'&=\LieH'\oplus 𝔪',\\
\LieG&=\LieH\oplus 𝔪
\end{align*}
\item
a Lie group morphism \(k\in H'\to\bar{k}\in H\) and
\item
an \(H'\)-equivariant linear map \(𝔪'\xrightarrow{\Gamma}\LieH\)
\item
an \(H'\)-equivariant linear isomorphism \(𝔪'\xrightarrow{ι}𝔪\)
\end{itemize}
On \(G'\times H\), applying the splitting \(\LieG'=\LieH'\oplus 𝔪'\) to construct connection and soldering form:
\[
\omega_{G'}=(\gamma_{G'},\sigma_{G'}),
\]
and denoting points of \(H\) as \(h\in H\), let
\[
\omega_\G=(\gamma_\G,\sigma_\G)
\]
where
\begin{align*}
\gamma_\G&=\Ad_h^{-1}\Gamma\sigma_{G'}+\Ad_h^{-1}\bar\gamma_{G'}+\omega_H,\\
\sigma_\G&=\rho_{𝔪}(h)^{-1}{ι}\sigma_{G'}.
\end{align*}
\begin{theorem}
The recipe above generates a homogeneous reductive Cartan geometry: the form \(\omega_{\G}\) is the pullback from \(\G:=\amal{G'}{H'}{H}\) of a unique \((X,G)\)-Cartan connection.
Every homogeneous reductive Cartan geometry arises this way uniquely up to isomorphism.
\end{theorem}
Denote the right invariant Maurer--Cartan form on the Lie group \(H\) by \(\rightMC[H]\).
The geodesic equations of any homogeneous reductive Cartan geometry are \(0=\gamma_\G\), \(A\,dt=\sigma_\G\), for any constant \(A\in\redComplement\), i.e.
\begin{align*}
\sigma_{G'}&=ι^{-1}\rho_{𝔪}(h)A\,dt,\\
\rightMC[H]+\bar\gamma_{G'}&=-\Gamma ι^{-1}\rho_{𝔪}(h)A\,dt.
\end{align*}
These are equations of a curve through each point of \(\G\), but in \(G'\times H\) they are equations of a principal \(H'\)-bundle over each of those curves. 
In particular, every geodesic on \(\G\) lifts to a curve on which \(\gamma_{G'}=0\), a horizontal lift of the geodesic, satisfying
\begin{align*}
\sigma_{G'}&=ι^{-1}\rho_{𝔪}(h)A\,dt,\\
\gamma_{G'}&=0,\\
\rightMC[H]&=-\Gamma ι^{-1}\rho_{𝔪}(h)A\,dt.
\end{align*}
We see that the final equation, defined inside \(H\), determines a unique curve \(h(t)\in H\) through each point of \(H\), and then the two previous equations are Lie equations turning this curve into a curve in \(G'\times H\).
In particular, if \(\Gamma ι^{-1}\rho_{𝔪}\) is \(H\)-equivariant, then this is also a Lie equation in \(H\), so the reductive Cartan geometry is geodesically complete.

\chapter{Cartan geometries}\label{section:Cartan.geometries}
\section{Definition}
Take a locally homogeneous structure with bundle \(\G={\G}_H\to M\).
Recall from \vref{subsection:G.over.H.bundle} the associated form \(\omega_{\G}\) on \(\G\).
This form \(\omega_{\G}\) is \emph{not} a connection on \(\G\to M\): it is valued in the wrong Lie algebra: \(\LieG\), not \(\LieH\).
But it \emph{is} a flat connection on \(\G_G\), as we will see.

An \((X,G)\)-\emph{geometry}\define{Cartan!geometry}\define{geometry!Cartan},\label{def:Cartan.geometry} also called a \emph{Cartan geometry} modelled on \((X,G)\), on a manifold \(M\) is a principal \(H\)-bundle \(\G=\G_H\to M\), and a connection \(\omega=\omega_{\G_G}\) on the associated \(G\)-bundle \(\G_G\), called the \emph{Cartan connection},\define{Cartan!connection}\define{connection!Cartan} so that, just as for the model, \(\omega\) is a linear isomorphism of each tangent space
\[
\begin{tikzcd}
0\arrow[r]&
T_p \G_H\arrow[r,"\omega"]
&\LieG\arrow[r]&0.
\end{tikzcd}
\]
The group \(H\) is the \emph{structure group}\define{structure group} of the Cartan geometry.
\begin{problem}{ai}
Prove that this condition of being a linear isomorphism is equivalent to: \(\omega\ne 0\) on nonzero tangent vectors to the \(H\)-bundle \(\G=\G_H\), i.e. nonzero tangent vectors to \(\G\) are not horizontal.
\end{problem}
An \emph{isomorphism}\define{isomorphism!of Cartan geometries}\define{Cartan!geometry!isomorphism} (automorphism) of Cartan geometries is a bundle isomorphism (automorphism) preserving the Cartan connection.
\begin{example}
Pick a Cartan geometry \(\G\xrightarrow{\pi}M\) and an open set \(U\subseteq M\).
The \emph{induced Cartan geometry}\define{induced!Cartan geometry}\define{Cartan!geometry!induced} on \(U\) is \(\left.\G\right|_U:=\pi^{-1}U\to U\) with the Cartan connection restricted to that open set \(\left.\G\right|_U\subseteq \G\).
\end{example}
\begin{example}
Similarly if \(M\xrightarrow{\varphi}M'\) is a local diffeomorphism, each Cartan geometry \(\G'\to M'\) pulls back to the \emph{pullback Cartan geometry}\define{pullback!Cartan geometry}\define{Cartan!geometry!pullback} \(\G\to M\) defined by taking \(\G\) the \emph{pullback bundle},\define{pullback!bundle} i.e. 
\[
\G:=\set{(m,p')|m \in M, p'\in \G'_{m'}, m':=\varphi(m)},
\]
with obvious bundle map
\[
p=(m,p')\in \G\mapsto p'\in \G',
\]
and Cartan connection the pullback by this bundle map.
\end{example}
\begin{example}
If \(X'\subseteq X\) is a path component and \(G'\subseteq G\) is the subgroup preserving \(X'\) then any \((X,G)\)-geometry is precisely an \((X',G')\)-geometry, so we can assume, without loss of generality, that the model \(X\) of any Cartan geometry is connected.
\end{example}
\begin{example}
If some local diffeomorphism \(X\to X'\) is equivariant for a Lie group morphism \(G\xrightarrow{\Phi}G'\), every \((X,G)\)-geometry induces an \((X',G')\)-geometry 
with \(\G':=\amal{\G}{H}{H'}\) and with \(\omega':=\Phi'(1)\omega\).
So a Riemannian geometry induces a conformal geometry, and also a projective geometry, and so on.
\end{example}
\begin{problem}{aj}
Suppose that \((\tilde{X},\tilde{G})\to(X,G)\) is the universal covering homogeneous space.
Prove that any \((X,G)\)-geometry is precisely a \((\tilde{X},\tilde{G})\)-geometry.
\end{problem}
\begin{example}
Take finitely many Cartan geometries \(H_i\to\G_i\to M_i\), say with models \((X_i,G_i)\).
The \emph{product}\define{product!of Cartan geometries}\define{Cartan!geometry!product} has model \((\prod_i X_i,\prod_i G_i)\), structure group \(\prod_i H_i\), total space \(\prod_i\G_i\), base space \(\prod_i M_i\).
The Cartan connection of the product is the sum of the Cartan connections.
\end{example}
\section{Discrete quotients}
\begin{theorem}
Suppose that \(\G\to M\) is a Cartan geometry and that \(\Gamma\) is a group of diffeomorphisms of \(\G\), each commuting with the structure group action and preserving the Cartan connection.
Let \(M'\subseteq M\) be the free regular set of the \(\Gamma\)-action, \(\G'\subseteq\G\) its preimage.
Suppose that \(M'\) is not empty.
Let \(\bar\G:=\Gamma\backslash\G'\) and \(\bar{M}:=\Gamma\backslash M'\).
Then \(\G'\to M'\) is the pullback Cartan geometry of a unique Cartan geometry on \(\bar\G\to\bar{M}\).
\end{theorem}
\begin{proof}
Apply problem~\vref{problem:quotient.conn} to make the quotients \(\bar\G\to\bar{M}\) and the connection on the quotient \(\bar\G_G=\overline{\G_G}\).
On each sheet of the covering, the tangent spaces of \(\G_H\) and the horizontal spaces of the connection are identified with those of \(\bar\G_H\) and the horizontal spaces of the connection on the quotient, so the connection remains a Cartan connection.
\end{proof}
\section{Coordinate expressions}
Take a homogeneous space \((X,G)\) and a point \(x_0\in X\) and let \(H:=G^{x_0}\).
Take a manifold \(M\) of the same dimension as \(X\).
A \emph{Cartan gauge}\define{Cartan gauge} \cite{Sharpe:2002} p. 174 on \(M\) is a \(1\)-form \(\eta\) on \(M\), valued in \(\LieG\), so that \(\eta+\LieH\) is a linear isomorphism of each tangent space of \(M\) to \(\LieG/\LieH\).

Take an \((X,G)\)-geometry \(\G\to M\) on a trivial bundle \(\G=M\times H\) (which is always the case after we perhaps replace \(M\) by a sufficiently small open subset of \(M\) around any given point).
Denote points of \(M\times H\) by \((m,h)\).
Since the Cartan connection \(\omega_{\G}\) is \(H\)-equivariant, and transforms in the adjoint \(H\)-representation, the form
\[
\eta:=\Ad_h(\omega_{\G}-\omega_H)
\]
is \(H\)-invariant on \(M\times H\) and vanishes on each fiber \(\set{m_0}\times H\).
So on \(M\times G\), in any local coordinates \(x^1,\dots,x^n\) on \(M\), \(\eta\) is a linear combination of \(dx^1,\dots,dx^n\).
But it doesn't transform under the \(H\)-action, so its coefficients don't depend on \(h\in H\), i.e. \(\eta\) is a differential form on \(M\).
Replacing \(M\) by an open subset of \(M\), in which we have coordinates \(x^1,\dots,x^n\),
\[
\eta=A_i(x)dx^i
\]
for functions
\[
M\xrightarrow{A_i}\LieG.
\]
Solve for \(\omega_{\G}\):
\[
\omega_{\G}=\omega_H+\Ad_h^{-1}(A_i(x)dx^i).
\]

For \(\omega_{\G}\) to be a Cartan connection, we need
\[
v\in T_p\G\mapsto v\hook\omega_{\G}\in\LieG
\]
to be a linear isomorphism, or equivalently to be nonzero for nonzero \(v\).
Write out our tangent vector as
\[
v=(v_M,A_H)
\]
so that
\[
v\hook\omega_{\G}=A+\Ad_h^{-1}(v_M\hook\eta).
\]
So \(\omega_{\G}\) is a Cartan connection precisely when 
\[
v_M\hook\eta\ne \Ad_h A,
\]
for any \(A\in\LieH\), \(v_M\in T_m M\), unless \(v_M=0\) and \(A=0\).
Equivalently, we need that
\[
v_M\mapsto v_M\hook\eta+\LieH\in\LieG/\LieH
\]
is a linear isomorphism.
In other words, \(\omega_{\G}\) is a Cartan connection precisely when \(\eta\) is a Cartan gauge.
So \(\omega_{\G}\) is a Cartan connection precisely when the matrix valued functions
\[
A_1,A_2,\dots,A_n\colon M\to\LieG
\] 
when projected to \(\LieG/\LieH\), span \(\LieG/\LieH\) at every point.
\begin{problem}{ak}
Recall that we can recover any principal \(H\)-bundle \(\G\) as the quotient
\[
\G\cong\left(\bigsqcup_a M_a\times H\right)/\sim
\]
where
\[
(m,h)\in M_a\times H\sim(m,h_{ab}h)\in M_b\times H,
\]
for some open cover of \(M\) by open sets \(M_a\) and principal bundle transition maps 
\[
M_a\cap M_b\xrightarrow{h_{ab}}H.
\]
Show that a Cartan connection on this bundle is precisely some Cartan gauge \(\eta_a\) on each \(M_a\) so that
\[
\eta_b(x)=h_{ab}^*\omega_G+\Ad_{h_{ab}(x)}^{-1}\eta_a(x).
\]
\end{problem}
We won't mention Cartan gauges again.
\section{Another definition}
We prove that our definition above of Cartan geometry corresponds to another definition, which is frequently encountered.

Take a homogeneous space \((X,G)\) and a point \(x_0\in X\) and let \(H:=G^{x_0}\).
Take a manifold \(M\).
Many authors define an \((X,G)\)-Cartan connection on \(M\) to be precisely an \((X,G)\)-subgeometry, but with the added hypothesis that \(\dim M=\dim X\).
We will see that this is equivalent to the definition~\vpageref{def:Cartan.geometry}.

First some motivation.
Recall that a locally homogeneous space modelled on \((X,G)\) (at least if \((X,G)\) is strong and effective) is an \(H\)-bundle and a Cartan connection on the \(H\)-bundle, identifying all of the tangent spaces of the \(H\)-bundle with \(\LieG\).
Let us generalize.
Suppose that \(\G_H\to M\) is a principal right \(H\)-bundle and that \(\dim M=\dim X\).
A \emph{preconnection}\define{Cartan!preconnection} is a \(\LieG\)-valued \(1\)-form \(\omega\) on \(\G_H\) transforming in the adjoint representation under the \(H\)-action, and agreeing with the Maurer--Cartan form on the fibers of \(\G_H\to M\).
Note that a preconnection  is a subgeometry just when it is an injective map on all tangent spaces of \(\G_H\) to \(\LieG\).

Let us expand out this definition.
Denote the right \(H\)-action on the bundle as \(\RT{h}p=ph\).
For any \(A\in\LieH\), let
\[
A_{\G}=\left.\frac{d}{dt}\right|_{t=0}\RT{e^{tA}}.
\]
A preconnection is precisely a \(\LieG\)-valued \(1\)-form \(\omega\in\Omega^1_{\G_H}\otimes\LieG\) satisfying
\[
\RT{h}^*\omega=\Ad_h^{-1}\omega, \ h\in H,
\]
and
\[
A_{\G}\hook\omega=A, \ A\in\LieH.
\]
\begin{problem}{preconn}
Prove that any preconnection on a bundle \(H\to\G\to M\) is the pullback of a unique connection on \(\G_G\).
Conversely, prove that every connection on \(\G_G\) pulls back to a preconnection.
If \(\dim M=\dim X\) then prove that  a connection on \(\G_G\) is a Cartan geometry just when the preconnection it pulls back to is a subgeometry, i.e. is injective on all tangent spaces of \(\G_H\).
\end{problem}
Henceforth we will not distinguish a Cartan connection (in our sense above) from the associated preconnection on a manifold \(M\) with \(\dim M=\dim X\).
\begin{answer}{preconn}
Suppose that \(\omega\) is a preconnection.
Define a \(\LieG\)-valued \(1\)-form 
\[
\omega'\in\Omega^1_{\G_H\times G}\otimes\LieG
\]
by the rule: at each point \((p_0,g_0)\in\G_H\times G\), let
\[
\omega':=\Ad_{g_0}^{-1}\omega+\omega_G.
\]
Clearly \(\omega'\) pulls back by \(p\in\G_H\mapsto(p,1)\in\G_H\times G\) to become \(\omega\).
We need to prove that \(\omega'\) is the pullback of a unique connection on \(\G_G\).
Denote the right \(H\)-action on \(\G_H\times G\) by
\[
\rho_h(p_0,g_0)=(p_0h,h^{-1}g_0),
\]
and the right \(G\)-action by
\[
\RT{g}(p_0,g_0)=(p_0,g_0g).
\]
Note that \(H\) acts freely and properly on \(\G_H\), hence on \(\G_H\times G\), so \(\G_G\) is a smooth manifold and 
\[
H\to\G_H\times G\to\G_G
\]
is a principal right \(H\)-bundle.
The reader can check that \(G\) acts freely and properly on \(\G_G\), which is therefore a principal right \(G\)-bundle over \(M\).
Under these actions
\begin{align*}
\rho_h^*\omega&=\Ad_h^{-1}\omega,\\
\rho_h^*\omega_G&=\omega_G,\\
\RT{g}^*\omega&=\omega,\\
\RT{g}^*\omega_G&=\Ad_g^{-1}\omega_G.
\end{align*}
Therefore 
\begin{align*}
\rho_h^*\omega'&=\omega',\\
\RT{g}^*\omega'&=\Ad_g^{-1}\omega'.
\end{align*}
Denote the vector fields by
\begin{align*}
A_H&=\left.\frac{d}{dt}\right|_{t=0}\rho_{e^{tA}},\\
A_G&=\left.\frac{d}{dt}\right|_{t=0}\RT{e^{tA}}.
\end{align*}
Check that at any point \((p_0,g_0)\in\G_H\times G\),
\begin{align*}
A_H\hook\omega&=A,\\
A_H\hook\omega_G&=-\Ad_{g_0}^{-1}A,\\
A_G\hook\omega&=0,\\
A_G\hook\omega_G&=A.
\end{align*}
Therefore
\begin{align*}
A_H\hook\omega'&=0,\\
A_G\hook\omega'&=A.
\end{align*}
So \(\omega'\) vanishes on the \(H\)-orbits in \(\G_H\times G\), and is constant under the \(H\)-action.
We leave the reader to prove that this occurs precisely when \(\omega'\) is the pullback of a differential form \(\omega''\) from \(\G_G\), and this form is then unique, perhaps using appendix~\ref{appendix:semibasic}.
It follows by uniqueness that \(\RT{g}^*\omega''=\Ad_g^{-1}\omega''\) and that \(A_G\hook\omega''=A\).
Hence \(\omega''\) is a connection.

Suppose there are two connections on \(\G_G\) pulling back to the same \(\omega\) on \(\G_H\).
Their difference \(\delta\) pulls back to zero on \(\G_H=\G_H\times\set{1}\).
On \(\G_H\times G\), \(\delta=0\) on the \(G\)-orbits, since both connections have \(A_G\hook{}=A\).
On \(\G_H\times G\), \(\delta=0\) on the \(H\)-orbits and is \(H\)-invariant, since \(\delta\) is pulled back from \(\G_G=\amal{\G}{H}{G}\).
By \(H\)-invariance, \(\delta=0\) in all directions of \(\G_H\times G\).
\end{answer}
\section{Flat Cartan geometries}
A Cartan geometry is \emph{flat}\define{flat!Cartan geometry}\define{Cartan!geometry!flat} if its connection is flat, or in other words if it is flat as a subgeometry.
\begin{example}
The principal \(H\)-bundle \(g \in G \mapsto gx_0 \in X\) is a flat Cartan geometry on \(X\), called the \emph{model geometry},\define{model geometry}\define{geometry!model} with the Maurer--Cartan \(1\)-form as Cartan connection, as we will see.
\end{example}
We restate theorem~\vref{theorem:XG.geom.to.struc}.
\begin{theorem}
Take a strong effective homogeneous space \((X,G)\).
Every \((X,G)\)-structure determines a flat \((X,G)\)-Cartan geometry as above.
Conversely, every flat \((X,G)\)-Cartan geometry arises from a unique \((X,G)\)-structure. 
\end{theorem}
This follows from theorem~\vref{thm:flat}.
\section{Holomorphic Cartan geometries}
A \emph{holomorphic Cartan geometry}\define{holomorphic Cartan geometry}\define{Cartan!geometry!holomorphic} is a Cartan geometry \(\G\to M\) modelled on a complex homogeneous space, with \(\G\) and \(M\) complex manifolds, \(\G\to M\) a holomorphic principal bundle, and the Cartan connection a holomorphic \(1\)-form.
\begin{example}
Recall from~\vpageref{example:linear.ode} that every holomorphic linear second order ordinary differential equation gives rise to a holomorphic Cartan geometry, on some Riemann surface, modelled on the complex projective line.
Conversely, suppose we take any holomorphic Cartan geometry with that model.
Since the geometry is on a Riemann surface, i.e. a complex manifold of complex dimension one, and the Cartan geometry is holomorphic, the Cartan connection has curvature zero.
We will see that this implies that the Cartan geometry is flat, and hence it is a locally homogeneous structure.
As we saw in our study of holomorphic linear second order ordinary differential equations, every locally homogeneous structure modelled on the complex projective line arises from a holomorphic linear second order ordinary differential equation.
Hence we have an isomorphism of categories, i.e. nothing but a change of terminology: a holomorphic Cartan geometry, on a connected manifold, modelled on the complex projective line, is precisely a holomorphic linear second order ordinary differential equation.
\end{example}

\section{Pseudo-riemannian geometry}\label{section:pseudo.Riem}
Fix a finite dimensional vector space \(V\) with a definite quadratic form.
Let \(X:=V\), which we think of as an affine space with \(V\) as its group of translations.
(We resist the temptation to use the fancy letter \(\redComplement\) here in place of \(V\), because it is easier somehow to think of a vector space as \(V\), but much of this discussion generalizes to geometries modelled on reductive homogeneous spaces; we will discuss reductive geometries \vpageref{reductive.homogeneous.spaces}.)
Let \(H:=\Orth{V}\), the orthogonal group of that form, and let \(G:=H\ltimes V\) the group of rigid motions of \(X\).
Suppose that \(M\) is a manifold equipped with a pseudo-Riemannian geometry, of the same signature as the form on \(X\).
The \emph{orthonormal frame bundle}\define{frame bundle!orthonormal}\define{orthonormal!frame bundle} \(\G=\G_H\) is the set of pairs \((m,u)\) where \(m\in M\) and \(T_m M\xrightarrow{u}V\) is a linear isometry
\begin{problem}{al}
Give \(\G\) the smooth structure of an \(H\)-bundle under the right \(H\)-action
\[
(m,u)h:=(m,h^{-1}u).
\]
Hint: Gram--Schmidt orthogonalization is a smooth map taking bases to orthonormal bases.
\end{problem}
(The right action turns out to be more convenient for us than the more obvious left action.)
We also denote this action as \(\RT{h}(m,u)=(m,u)h\), so \(\G\xrightarrow{\RT{h}}\G\).
Denote the Lie algebra of \(H\) by \(\LieH\).
For any \(A\in\LieH\), let \(A_{\G}\) be the vector field on \(\G\) which gives the infinitesimal Lie algebra action:
\[
A_{\G}(m,u):=\left.\frac{d}{dt}\right|_{t=0}(m,e^{-tA}u).
\]
Denote the bundle projection mapping by
\[
(m,u)\in \G\xmapsto{\pi} m\in M.
\]
As we move a frame, we move both its underlying point \(m\) and its linear isomorphism \(u\).
The \emph{soldering form}\define{soldering form} measures the velocity \(\dot{m}\) of the motion of the point, but as measured in the frame \(u\).
In other words, the \emph{soldering form} \(\sigma\) is the \(V\)-valued \(1\)-form defined on tangent vectors \(v\in T_{(m,u)}\G\) by
\[
v\hook\sigma_{(m,u)}:= u(\pi'(m,u)v).
\]
\begin{problem}{am}\(\RT{h}^*\sigma=h^{-1}\sigma\), for any \(h\in H\), and hence \(\LieDer_{A_{\G}}\sigma=-A\sigma\), for any vector \(A\in\LieH\).
\end{problem}
A connection on \(\G\) is a \(1\)-form \(\gamma\) valued in \(\LieH\) so that \(A_{\G}\hook\gamma=A\) for all \(A\in\LieH\) and so that
\[
\RT{h}^*\gamma=\Ad_h^{-1}\gamma
\]
for all \(h\in H\).

What does it mean? A connection should measure the motion \(\dot{u}\) of the frame \(u\), ``ignoring'' the motion \(\dot{m}\) of the point \(m\).
But this does not have an unambiguous interpretation because at different times, \(u\) is an isomorphism between different vector spaces \(T_m M\xrightarrow{u}V\).
The choice of a connection precisely breaks that ambiguity.

The equation \(A_{\G}\hook\gamma=A\) means that the connection, faced with a rotating frame \(e^{-tA}u\) at a fixed point \(m\), will of course tell us that \(u\) identifies fixed vectors in \(V\) with moving vectors in \(T_m M\) rotating with angular velocity \(A\).
The connection thus equals the left invariant Maurer--Cartan form on the fibers.
Adjoint invariance of \(\gamma\) matches precisely how the left invariant Maurer--Cartan form on \(H\) transforms under right action.
In any basis \(\set{e_i}\) of \(V\), \(\gamma\) becomes a matrix whose entry \(\gamma^i_j\) represents the tendency of a moving frame to rotate the vector corresponding to \(e_i\) in the direction of the vector corresponding to \(e_j\).
In a Riemannian geometry, with an orthonormal basis, just as the \(e_j\) vector turns toward the \(e_i\) vector, at the same rate the \(e_i\) vector turns away, hence the antisymmetry of \(\gamma\).
\begin{problem}{an}For any connection \(\gamma\), \(A_{\G}\hook(d\sigma+\gamma\wedge\sigma)=0\) for any \(A\in\LieH\).
\end{problem}
\begin{answer}{an}
Note that
\begin{align*}
A_{\G}\hook\sigma&=0,\\
A_{\G}\hook\gamma&=A.
\end{align*}
so that
\[
A_{\G}\hook(\gamma\wedge\sigma)=A\sigma.
\]
\begin{align*}
-A\sigma
&=
\LieDer_{A_{\G}}\sigma,
\\
&=
A_{\G}\hook d\sigma+d(A_{\G}\hook\sigma),
\\
&=
A_{\G}\hook d\sigma+d(0),
\\
&=
A_{\G}\hook d\sigma.
\end{align*}
So
\[
0=A_{\G}\hook(d\sigma+\gamma\wedge\sigma).
\]
\end{answer}
\begin{problem}{conn.iso}
Prove that at each point \(p\in\G\),
\[
T_p\G\xrightarrow{(\sigma,\gamma)}V\oplus\LieH
\]
is a linear isomorphism, and that the kernel of \(\sigma\) is precisely the set of vertical vectors.
\end{problem}
\begin{answer}{conn.iso}
By definition,
\[
v\hook\sigma=u(\pi'(m,u)v).
\]
Since \(\G\to M\) is a bundle map, it is a submersion, so \(\pi'(m,u)\) is onto, and \(u\) is a linear isomorphism, so \(\sigma\) is onto \(V\).
Its kernel is precisely the vectors \(v\) with \(0=u(\pi'(m,u)v)\).
Since \(u\) is a linear isomorphism, this kernel is precisely the kernel of \(\pi'(m,u)\), i.e. precisely the vertical vectors.
Since \(\G\to M\) is a principal bundle, each vertical vector has the form \(v=A_{\G}\), for a unique \(A\in\LieH\).
So \(\gamma\) is a linear isomorphism of the vertical vectors to \(\LieH\).
Hence, at each point \(p\in\G\),
\[
T_p\G\xrightarrow{(\sigma,\gamma)}V\oplus\LieH
\]
is a linear isomorphism.
\end{answer}
\begin{problem}{conn.components}
In any basis of \(V\), and associated basis of \(\LieH\), the components of \(\sigma,\gamma\) put together form a basis of the \(1\)-forms on \(\G\).
The components of \(\sigma\) form a basis of the semibasic \(1\)-forms.
\end{problem}
\begin{problem}{dsigma}
Prove that
\[
d\sigma+\gamma\wedge\sigma=\frac{1}{2}t\sigma\wedge\sigma,
\]
for a unique \(H\)-equivariant map \(\G\xrightarrow{t}V\otimes\Lm{2}{V}^*\), the \emph{torsion}\define{torsion!of an affine connection} of the connection \(\gamma\).
\end{problem}
\begin{answer}{dsigma}
Let \(\Sigma:=d\sigma+\gamma\wedge\sigma\).
We have seen that \(A_{\G}\hook\Sigma=0\) for all \(A\in\LieH\).
Every vertical vector \(v\) has the form \(v=A_{\G}\) for a unique \(A\in\LieH\), since \(\G\to M\) is a principal right \(H\)-bundle.
In other words, \(\Sigma\) is a semibasic \(2\)-form.
Take a basis \(\set{e_i}\) of \(V\) and a basis \(\set{A_{\alpha}}\) of \(\LieH\).
In these bases, expand out \(\sigma,\gamma\) into their components, a basis of \(1\)-forms on \(\G\), by problem~\vref{problem:conn.components}: 
\[
\sigma=\sigma^ie_i, \gamma=\gamma^{\alpha}A_{\alpha}.
\]
Write out the components of \(\Sigma\) in this basis, say:
\[
\Sigma^i=t^i_{jk}\sigma^j\wedge\sigma^k+a^i_{\alpha j}\gamma^{\alpha}\wedge\sigma^j+
b^i_{\alpha\beta}\gamma^{\alpha}\wedge\gamma^{\beta}.
\]
We can assume that \(t^i_{jk}+t^i_{kj}=0\) and that \(b^i_{\alpha\beta}+b^i_{\beta\alpha}=0\).
Plug in each vertical vector \(v=A_{\alpha,\G}\):
\[
0=v\hook\Sigma^i=a^i_{\alpha j}\sigma^j+b^i_{\alpha\beta}\gamma^{\beta}.
\]
So \(0=a=b\).
Similarly, every semibasic form is a unique wedge product of components of \(\sigma\).
\end{answer}
The fundamental lemma of Riemannian geometry \cite{ONeill:1983} p. 61 Theorem 11 tells us that there is a unique torsion free connection \(\gamma\) on \(\G\), i.e. so that
\[
0=d\sigma+\gamma\wedge\sigma.
\]
\begin{problem}{fund.lemma.riem.geom}
Prove the fundamtal lemma of Riemannian geometry, assuming only that every principal bundle admits a connection.
\end{problem}
\begin{answer}{fund.lemma.riem.geom}
Take a connection \(\gamma\).
Any other connection \(\gamma'\) agrees with \(\gamma\) on the fibers of the frame bundle, so disagrees by something \(\gamma'-\gamma\) vanishing on the fibers, so depending at each point only on the projection to the base, i.e. on \(\pi'\).
But then \(\sigma\) has the same kernel as \(\pi'\), so
\[
\gamma'=\gamma+q\sigma,
\]
where \(\G\xrightarrow{q}V^*\otimes\LieH\) is \(H\)-equivariant, since \(\gamma,\gamma'\) are.
Conversely, any \(H\)-equivariant \(q\) gives a new connection
\[
\gamma'=\gamma+q\sigma.
\]
The torsion of \(\gamma'\) is 
\[
d\sigma-\gamma'\wedge\sigma=\frac{1}{2}t'\sigma\wedge\sigma,
\]
but plug in \(\gamma'=\gamma+q\sigma\) to see that
\[
t'=t-\delta_q
\]
where 
\[
(\delta_q)^i_{jk}=\frac{1}{2}(q^i_{jk}-q^i_{kj}).
\]
At this point, algebra steps in, which we leave to the reader, to prove that there is a unique \(q\) which will make \(\gamma'\) torsion free.
\end{answer}
Write elements of \(G\) as matrices
\[
g=
\begin{pmatrix}
h&v\\
0&1
\end{pmatrix}
\]
so that \(h\in H=\Orth{V}\) and \(v\in V\).
\begin{problem}{ao}
For any connection \(\gamma\) on \(\G\to M\), the \(1\)-form
\[
\omega
=
\begin{pmatrix}
\gamma&\sigma\\
0&0
\end{pmatrix}
\]
is an \((X,G)\)-Cartan connection on \(\G\to M\).
Conversely, every \((X,G)\)-Cartan connection is carried by a unique bundle isomorphism to one obtained this way.
\end{problem}
Consequently there are \((X,G)\)-Cartan connections with torsion, corresponding to choices of connection \(\gamma\) with torsion, so the set of \((X,G)\)-Cartan geometries is a bit larger than the set of pseudo-Riemannian manifolds of the given signature.

\section{The method of equivalence}
A host of theorems demonstrate an isomorphism of categories between some type of geometric structure (e.g.: a Riemannian metric) and a Cartan geometry with a certain model (e.g. Euclidean space) and some condition on its curvature (e.g. torsion freedom); see \v{C}ap and Slov\`ak \cite{Cap/Slovak:2009}.
Choices need to be made but, modulo those choices, there is a general procedure for carrying out these constructions, known as \emph{Cartan's method of equivalence}.\SubIndex{method of equivalence}
In Cartan's hands \cite{Cartan:1910,Cartan:136,Cartan:136bis,Cartan:161}, the method is famously opaque:
\begin{quotation}
We must confess that we find most of these papers extremely rough going and we certainly cannot follow all the arguments in detail. The best procedure is to guess at the theorems, then prove them, then go back to Cartan.
\par\noindent\qquad\rule[.35em]{1cm}{.2pt} I.~M.~Singer, S.~Sternberg \cite{Singer.Sternberg:1965}
\end{quotation}
\begin{quotation}
It is this problem with which Cartan deals in the present book, and in some way he reduces the second influence, the choice of parameters, to the choice of the frame. I did not quite understand how he does this in general, though in the examples he gives the procedure is clear. \dots

All of the author's books, the present one not excepted, are highly stimulating,
full of original viewpoints, and profuse in interesting geometric details. Cartan is undoubtedly the greatest living master of differential geometry. \dots

\dots \ Nevertheless, I must admit that I found the book, like most of Cartan's papers, hard reading.
\par\noindent\qquad\rule[.35em]{1cm}{.2pt} Hermann Weyl \cite{Weyl1938}
\end{quotation}
Bryant \&{} Griffiths \cite{Bryant.Griffiths:1995} and Gardner \cite{Gardner:1989} make the method clear, but leave the reader to normalize curvature (``the choice of the parameters'').
Generalizing work of Tanaka, \v{C}ap and Slov\`ak \cite{Cap/Slovak:2009} naturally normalize curvature of a huge collection of geometric structures to associate to each one a unique Cartan geometry from which it can be recovered, and characterize the Cartan geometries that arise from their construction.
Researchers are solving some very important problems using their approach \cite{Kruglikov.The:2017}.
The method of equivalence does not always yield a Cartan geometry.
We will not explain the method of equivalence.
\section{Constant vector fields}\label{section:constant.vec.fields}
Take a homogeneous space \((X,G)\) and an \((X,G)\)-Cartan geometry \(\G\to M\).
A vector field \(Z\) on \(\G\) is \emph{constant}\define{constant vector field} if \(Z\hook\omega\) is constant.
Since \(\omega\) is a linear isomorphism \(T_p \G\to\LieG\), each \(A\in\LieG\) has a unique constant vector field \(Z\) defined by \(Z\hook\omega=A\); we denote this vector field as \(A_{\G}\).
We denote the time one flow of a vector field \(Z\) by \(e^Z\), but let \(\fl{A}:=e^{A_{\G}}\) to save ink.
\begin{example}
On the model \((X,G)\), the constant vector fields are precisely the left invariant vector fields on \(G\).
\end{example}
\begin{lemma}\label{lemma:adj.H}
For any \(h\in H\) and \(A\in\LieG\), \(\RT{h*}(A_{\G})=(\Ad_h^{-1}A)_{\G}\).
\end{lemma}
\begin{problem}{lemma:adj.H.problem}
Compute out the left hand side to give the proof.
\end{problem}
\begin{answer}{lemma:adj.H.problem}
\begin{align*}
\RT{h*}(A_{\G})\hook\omega
&=
A_{\G}\hook \RT{h}^*\omega,
\\
&=
A_{\G}\hook \Ad_h^{-1}\omega,\\
&=
\Ad_h^{-1}(A_{\G}\hook\omega),\\
&=
\Ad_h^{-1}A,\\
&=
(\Ad_h^{-1}A)_{\G}\hook\omega.
\end{align*}
\end{answer}
\begin{problem}{ap}
The Cartan geometry is determined by its \(H\)-bundle structure and the linear map taking \(\LieG\) to the constant vector fields.
\end{problem}
Each flow line of any constant vector field in \(\G\) has image a curve in \(M\), sometimes called a \emph{spiral}\define{spiral} (or a \emph{generalized circle}\define{generalized circle} \cite{Sharpe:1997} p.210 or a \emph{canonical curve}\define{canonical curve} \cite{Cap/Slovak:2009} p. 112).
\begin{example}
In affine space \(X=\R^n\), \(G=\GL{n}\ltimes\R^n\), the generic \(1\)-parameter subgroup of \(G\) is, up to adjoint action, a linear rescaling of axes by some \(e^{t\lambda}\) for some real \(\lambda\), and some complex rescaling of some complex planes by some \(e^{t\lambda}\) for some complex \(\lambda\).
Hence geometrically, the generic vector draws out a curve which has various real components just growing or shrinking exponentially, and various complex components growing or shrinking exponentially while spinning around in a complex plane:
\end{example}
\[
\begin{tikzpicture}[scale=.3]
\begin{axis}[hide axis,colormap={bluered}{color=(red) color=(blue)}] 
 \addplot[domain=0:4,samples=401,smooth,very thick,mesh,point meta=y] 
 	({1.5*(x+0.4*x*x)+x*sin(x*360)},{-x*cos(360*x)});
\end{axis}
\end{tikzpicture}
\]
\begin{example}
In any homogeneous space, the differential invariants of any spiral are constant, as the symmetries commute with the constant vector fields, so symmetries bring any point of a spiral to any other point.
In particular, in the Euclidean plane with its flat Riemannian geometry, the spirals are the circles and the straight lines.
In \(3\)-dimensional Euclidean space, the spirals are the curves of constant curvature and torsion, i.e. helices, circles and straight lines:
\end{example}
\[
\begin{tikzpicture}[scale=.3]
\begin{axis}[hide axis,colormap={bluered}{color=(red) color=(blue)}] 
 \addplot[domain=0:4,samples=401,smooth,very thick,mesh,point meta=y] 
 	({x+cos(360*x)},{.1*x+sin(360*x)});
\end{axis}
\end{tikzpicture}
\]
Spirals in a homogeneous space \((X,G)\) are invariant both under the action of \(G\) and under affine changes of parameter, since the constant vector fields are left invariant, and rescaling a constant vector field gives a constant vector field.
\begin{problem}{Minkowski.spirals}
What are the spirals in the Minkowski plane, i.e. the translation invariant Lorentz geometry on the plane?
\end{problem}
\section{Dropping}
Take an equivariant map \(X\to X'\) of homogeneous \(G\)-spaces.
Take a point \(x_0\in X\) and its image \(x_0'\in X'\).
Let \(H:=G^{x_0}\) and \(H':=G^{x_0'}\).
so \(H\subseteq H'\subseteq G\) are closed subgroups.
Take an \((X',G)\)-geometry \(\G \to M'\).
Let \(M := \G/H\); then \(\G \to M\) is an \((X,G)\)-geometry with the same Cartan connection as the original \((X',G)\)-geometry, called the \emph{lift}\define{lift!of Cartan geometries}\define{Cartan!geometry!lift} or the \emph{twistor correspondence space}\define{twistor correspondence space} of \(\G \to M'\), a fiber bundle
\[
\begin{tikzcd}[cramped]
G\arrow[d]&&\G\arrow[d]\\
X=G/H\arrow[d]&&M=\G/H\arrow[d]\\
X'=G/H'&&M'=\G/H'.
\end{tikzcd}
\]
\begin{example}
Suppose that \(G\) is the group of isometries of a translation invariant pseudo-Riemannian metric \(g_0\) on a finite dimensional vector space \(X'\), say of dimension \(n\).
We have seen that every \(n\)-dimensional manifold \(M'\) with a pseudo-Riemannian metric \(g\) of the same signature as \(g_0\) has an \((X',G)\)-structure.
So \(X'=G/H'\) where \(H'\subseteq G\) is the orthogonal group of \(g_0\), i.e. the stabilizer of the origin \(0\in X'\).
Let \(H\subseteq H'\) be the subgroup fixing a vector \(v_0\in T_0 X'\).
Then \(M\to M'\) is the set of pairs \((m,v)\) of point \(m\in M'\) and tangent vector \(v\in T_m M'\) with \(g(v,v)=g_0(v_0,v_0)\).
So if \(v_0\) is a nonzero \(g_0\)-null vector, then \(M\) is the bundle of nonzero null vectors, while if \(v_0\) is a unit vector, then \(M\) is the unit tangent bundle.
\end{example}
\begin{example}
Denote by \(\C^{n+1,1}\) a complex vector space with a Lorentz Hermitian inner product.
Let \(X\) be the set of real null lines, and \(X'\) the set of complex null lines, in \(\C^{n+1,1}\).
Let \(G:=\SU{n+1,1}\).
Each real null line spans a unique complex null line: \(X\to X'\).
The manifold \(X\) has a \(G\)-invariant conformal structure, since \(G\subseteq\SO{2n+2,2}\).
The boundary \(M'=\partial D\) of a pseudoconvex domain in \(\C^{n+1}\) bears an \((X',G)\)-geometry, a CR-geometry \cite{Jacobowitz1990}.
Each CR geometry lifts to an \((X,G)\)-geometry \(M\to M'\), a conformal structure on a circle bundle.
\end{example}
Conversely, a given \((X,G)\)-geometry \emph{drops}\define{drop!of Cartan geometries}\define{Cartan!geometry!drop} to a given \((X',G)\)-geometry if it is isomorphic to the \((X,G)\)-lift of that \((X',G)\)-geometry.
\begin{example}
If \(X\) is a connected homogeneous \(G\)-space and \(*\) is a point with trivial \(G\)-action, then an \((X,G)\)-geometry on a connected manifold drops to a \((*,G)\)-geometry just when the geometry is isomorphic to its model.
\end{example}
More generally, if a geometry on some manifold \(M\) drops to some manifold \(M'\), then we can recover the manifold \(M\) and the original geometry on \(M\) directly from the geometry on \(M'\).
\begin{example}
If a \(2n+2,2\)-signature conformal structure drops to a CR-geometry, we can study that special type of conformal geometry using several variable complex analysis.
\end{example}
When we drop, an effective model could become ineffective; we need to allow ineffective models.
\section{Twistors}
Take \(G\)-equivariant maps of homogeneous spaces
\[
\begin{tikzcd}
&X'\arrow[dl]\arrow[dr]&\\
X&&X''
\end{tikzcd}
\]
Lift an \((X,G)\)-geometry to an \((X',G)\)-geometry; if it drops to an \((X'',G)\)-geometry, the geometries are \emph{twistor transforms}.\define{twistor transform}
The complex analytic Penrose twistor transform is a particular case \cite{Baston/Eastwood:1989}.

\section{From modules to vector bundles}
Take a homogeneous space \((X,G)\), a point \(x_0\in X\), and let \(H:=G^{x_0}\).
Take an \((X,G)\)-Cartan geometry \(\G \to M\).
For any \(H\)-module \(V\), the \emph{associated vector bundle}\define{associated!vector bundle} is \(\vb{V}:=\amal{\G}{H}{V}\to M\); its sections are the \(H\)-equivariant maps \(\G \to V\).
\begin{problem}{aq}
Prove that \(\vb{V}\) is a vector bundle, with those maps canonically identified with its sections.
\end{problem}
We use the same symbol \(\vb{V}\) for the associated vector bundle on the model \(X\) or on \(M\).
Consider the obvious linear projection mapping \(\LieG\xrightarrow{\pi}\LieG/\LieH=V\).
The \emph{soldering form}\define{soldering form} \(\sigma\) is the \(V\)-valued \(1\)-form \(\sigma:=\pi\circ\omega\).
\begin{proposition}[\cite{Sharpe:1997} p. 188, theorem 3.15]\label{prop:TM}
If \(V:=\LieG/\LieH\), then \(\vb{V}=\vbTM=TM\) 
\end{proposition}
\begin{proof}
Denote \(\G\to M\) as \(\G\xrightarrow{\pi}M\).
Pick \(p_0\in \G\) and let \(m_0:=\pi(p_0)\).
The commutative diagram of linear maps
\[
\begin{tikzcd}
&0\arrow[d] & 0 \arrow[d]\\
0\arrow[r]&T_{p_0} (\G_{m_0})\arrow[d] \arrow[r,"\omega"] & \LieH \arrow[d]\arrow[r]&0\\
0\arrow[r]&T_{p_0} \G\arrow[d,"\pi'(p_0)"] \arrow[r,"\omega"] & \LieG \arrow[d]\arrow[r]&0\\
0\arrow[r]&T_{m_0} M \arrow[r] \arrow[d]& \LieG/\LieH \arrow[d]\arrow[r]&0\\
&0 & 0 \\
\end{tikzcd}
\]
gives a commutative diagram of vector bundles
\[
\begin{tikzcd}
&0\arrow[d] & 0 \arrow[d]\\
0\arrow[r]&\ker\pi'\arrow[d] \arrow[r,"\omega"] &\pi^* \vbh \arrow[d]\arrow[r]&0\\
0\arrow[r]&T\G\arrow[d,"\pi'(p_0)"] \arrow[r,"\omega"] & \pi^*\vbg \arrow[d]\arrow[r]&0\\
0\arrow[r]&\pi^*TM \arrow[r] \arrow[d]& \pi^*(\vbg/\vbh) \arrow[d]\arrow[r]&0\\
&0 & 0 \\
\end{tikzcd}
\]
The bundle \(\vbh\) is the bundle \(\ad_{\G}\) which we saw in problem~\vref{problem:Atiyah.class}, where we saw that the exact sequence quotients by \(H\)-action:
\[
\begin{tikzcd}
&0\arrow[d] & 0 \arrow[d]\\
0\arrow[r]&\ker\pi'/H\arrow[d] \arrow[r,"\omega"] &\vbh \arrow[d]\arrow[r]&0\\
0\arrow[r]&T\G/H\arrow[d,"\pi'(p_0)"] \arrow[r,"\omega"] & \vbg \arrow[d]\arrow[r]&0\\
0\arrow[r]&TM \arrow[r] \arrow[d]& \vbg/\vbh\arrow[d]\arrow[r]&0\\
&0 & 0 \\
\end{tikzcd}
\]
\end{proof}
\begin{example}
The cotangent bundle is \((\vbTM)^*=\vbh^{\perp}\subseteq\vbg^*\), and similarly for the various tensor bundles.
\end{example}
\begin{example}
Every \(G\)-invariant tensor \(\tau\) on \(X\) has an \emph{associated tensor}\define{associated!tensor} \(\tau\) on the manifold \(M\) of any \((X,G)\)-geometry \(\G\xrightarrow{\pi}M\).
To define it, we take again \(V:=\LieG/\LieH\), so \(\vb{V}=\vbTM=TM\).
Let \(\tau_0\) be the value of \(\tau\) in the tangent space \(T_{x_0} X=V\).
So \(\tau_0\) is an \(H\)-invariant tensor in the vector space \(V\).
It therefore defines a tensor in each fiber of the associated vector bundle \(\vb{V}=TM\).
\end{example}

\chapter{Space forms}
We work out in detail the long example of the de~Sitter and anti-de~Sitter spaces, and the null quadric as a conformal pseudo-Riemannian geometry.
\section{Orthogonal linear algebra}
An \emph{inner product space}\define{inner product space} is a vector space with an inner product; we assume it is real and finite dimensional.
An \emph{orthogonal transformation}\define{orthogonal!transformation} is a linear isomorphism of inner product spaces preserving inner products.
A \emph{conformal transformation}\define{conformal!transformation} is a linear isomorphism of inner product spaces preserving inner products up to a positive constant factor.
Suppose that \(V\) is a finite dimensional inner product space.
Denote by \(\CO{V}\) the \emph{conformal group}:\define{conformal!group} the group of conformal transformations, which contains the \emph{orthogonal group}\define{orthogonal!group} \(\Orth{V}\): the group of orthogonal transformations.
Vectors in an inner product space are \emph{perpendicular}\define{perpendicular} if their inner product vanishes.
A vector \(x\) is \emph{unit}\define{unit vector} if \(\ip{x}{x}=\pm 1\).
A basis of \(V\) is \emph{orthonormal}\define{orthonormal!basis} if its elements are mutually perpendicular unit vectors.
Denote by \(\R^{p,q}\) the vector space \(\R^{p+q}\) with the inner product
\[
dx_1^2+dx_2^2+\dots+dx_p^2-dy_1^2-\dots-dy_q^2.
\]
Let \(\eep 1,\dots,\eep p\) be the standard basis in the \(x\)-variables, \(\eem1,\dots,\eem q\) the standard basis in the \(y\)-variables.
\begin{theorem}%
[Sylvester's Theorem of Inertia \cite{Garling:2011}]%
\define{theorem!Sylvester's theorem of inertia}\define{Sylvester!theorem of inertia}
Take a finite dimensional vector space \(V\) with inner product.
Any isometry between subspaces of isometric inner product spaces extends to an isometry of the spaces.
Hence every orthonormal basis for a subspace of \(V\) extends to an orthonormal basis of \(V\).
Hence \(V\) is isomorphic to \(\R^{p,q}\) for unique \(p,q\).
Hence the pair \(p,q\) is the unique isometry invariant of finite dimensional inner product spaces.
Any two orthonormal bases can be taken to one another, and hence to this standard basis, by a unique isomorphism of inner product spaces.
\end{theorem}
If \(p=q\), we say that \(V\) has \emph{split signature}.\define{split signature}
A vector \(v\in V\) is \emph{space-like}\define{vector!space-like}\define{space-like vector} if \(\ip{v}{v}>0\), \emph{time-like}\define{vector!time-like}\define{time-like vector} if \(\ip{v}{v}<0\), \emph{null}\define{vector!null}\define{null!vector} or \emph{light-like}\define{vector!light-like}\define{light-like vector} if \(\ip{v}{v}=0\).
A linear subspace of \(V\) is \emph{space-like}, \emph{time-like}, \emph{null} if all its nonzero vectors are.
A \emph{space and time splitting}\define{space and time splitting} is a perpendicular direct sum of a maximal dimensional space-like linear subspace, and a maximal dimensional time-like linear subspace.
For example, the span \(V^+\) of these \(\eep i\) with  the span \(V^-\) of these \(\eem i\) is a space and time splitting.

Every conformal linear transformation \(g\) scales inner products by some factor \(\lambda>0\).
Replacing \(g\) by a constant positive multiple, we can arrange that \(\lambda=1\).
We can scale a space-like or time-like vector to unit, equivariant under orthogonal transformations.
We can take any orthonormal basis to any other.
Any nonnull vector scales to unit, and then lies in some orthonormal basis.
\section{Affine quadric hypersurfaces}
So we can take any unit space-like vector to any other, any unit time-like vector to any other, and any nonzero null vector to any nonzero null vector, by orthogonal transformation.
We let \(Y_+\) be the set of space-like unit vectors, similarly \(Y_-\) the set of time-like unit vectors and \(Y_0\) the \emph{light cone},\define{light cone} i.e. the set of all nonzero null vectors.
Note that \(Y_+, Y_-\) are affine quadric hypersurfaces, while \(Y_0\) is a quadratic cone punctured at the origin, and \(\Orth{V}\) acts transitively on each of \(Y_+,Y_-,Y_0\),  so they are smooth submanifolds.
Note that \(\CO{V}\) acts only on \(Y_0\), not on \(Y_+\) or \(Y_-\), but it acts transitively on the open sets \(U_+,U_-\subset V\) of space-like and time-like vectors with \(\partial U_+=\partial U_-=Y_0\).
Each element of \(U_+\) becomes an element of \(Y_+\) under a unique positive rescaling: \(U_+\cong \R^+\times Y_+\) and similarly \(U_-\).
\begin{example}
If \(q=0\), \(V=\R^p\) with Euclidean inner product, \(Y_+=S^{p-1}\) is the unit sphere, \(Y_0\) and \(Y_-\) are empty.
\end{example}
\begin{example}
If \(p,q=2,1\), \(V=\R^{2,1}\) with Lorentzian inner product \(dx^2+dy^2-dz^2\), \(Y_+\) is the connected hyperboloid \(x^2+y^2=z^2+1\), \(Y_0\) the cone \(x^2+y^2=z^2\) and \(Y_-\) the disconnected hyperboloid \(x^2+y^2+1=z^2\).
\end{example}
\[
\includegraphics[width=3cm]{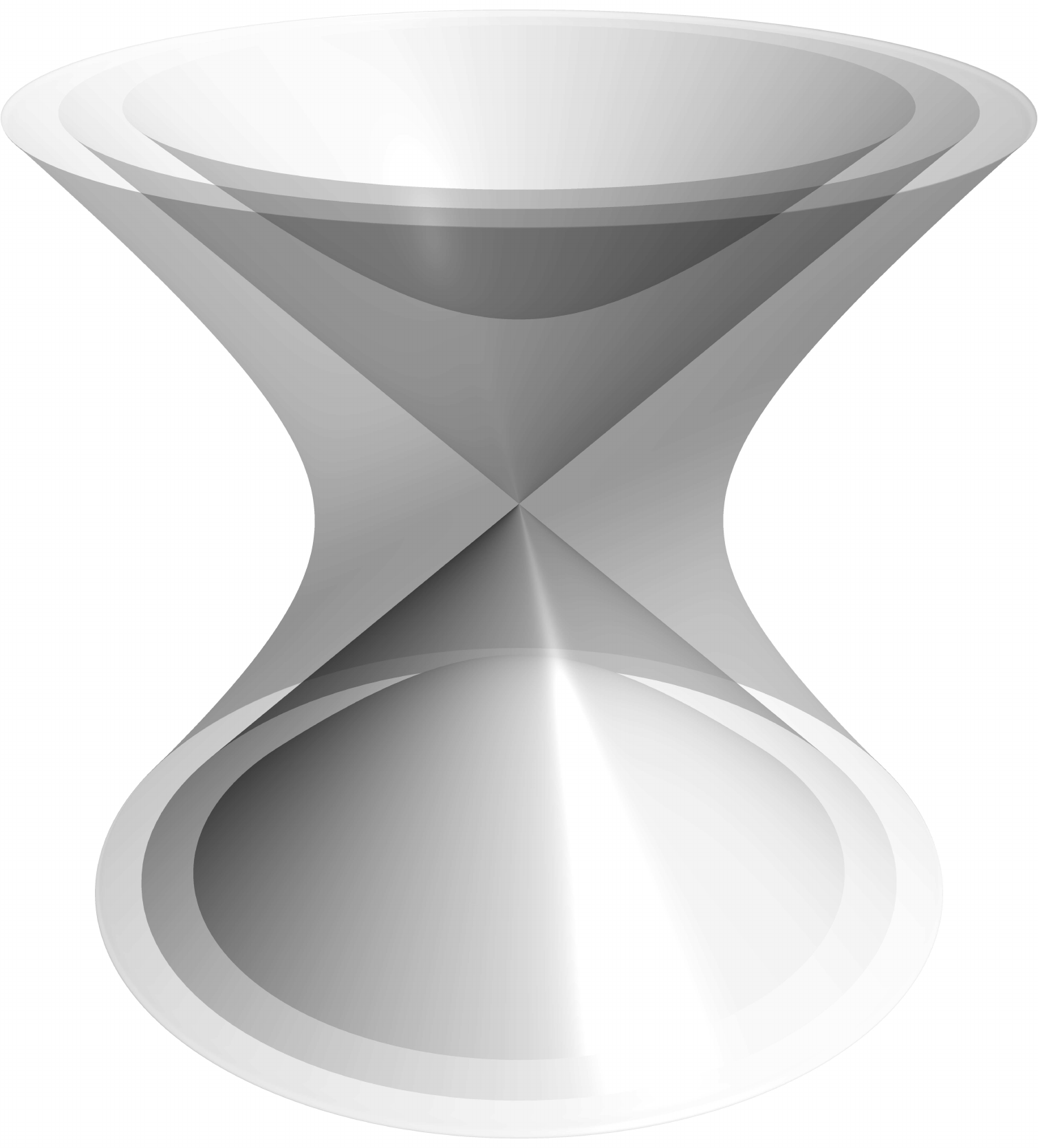}
\]
\section{Projective quadric hypersurfaces}
Let \(X_+\) be the set of space-like lines, \(X_-\) the set of time-like lines, \(X_0\) the set of null lines, called the \emph{null quadric}:\define{null!quadric}
\[
X_+,X_-,X_0\subset\Proj{}_V=\RP{p+q-1}.
\]
Clearly \(\CO{V}\) acts on lines through the origin, i.e. on \(\Proj{}V\), with orbits \(X_+,X_-,X_0,\set{0}\).
Note that every smooth quadratic hypersurface in real projective space arises as \(X_0\) for some quadratic form on \(V\).
A linear transformation in \(\CO{V}\) which acts trivially on \(X_+\) preserves all space-like lines, so has all space-like vectors as eigenvectors, hence more than a basis of eigenvectors, so is a multiple of the identity.
The same for time-like and even for null lines, so \(X_+,X_-,X_0\) are acted on trivially precisely by the elements \(\lambda I\in\CO{V}\), i.e. \(\R^{\times}\). 
The \emph{projective orthogonal group}\define{projective!orthogonal group} is \(G:=\CO{V}/\R^{\times}=\Orth{V}\!/\!\pm\! I\).
Clearly \(X_+,X_-\subset\Proj{}_V\) are open sets.
Being the orbit of Lie group action of \(G\) on \(\Proj{}_V\),  \(X_0\subset\Proj{}_V\) is a smooth submanifold, closed since it is locally the zero locus of a quadratic form in affine charts.
\begin{lemma}\label{lemma:projective.conformal}
An inner product space not of split signature has a nonzero null vector if and only if the group of projective transformations preserving its null quadric  is precisely the projective orthogonal group \(G\), i.e. the group of projective transformations arising from conformal (or orthogonal) transformations of the inner product space.
An inner product space of split signature has group of projective transformations preserving its null quadric precisely the group \(G\ltimes\set{I,T}\), where \(G\) is the projective orthogonal group and \(T\) is an involutive projective linear transformation interchanging the two components of the complement of the null quadric.
\end{lemma}
\begin{proof}
Every element of \(G\) arises, by definition, from a conformal linear transformation.
Conversely, suppose that \(g\) is a projective transformation preserving the null quadric, and bears a nonzero null vector, i.e. the null quadric is not empty.
So \(g\) is the image in the projective group of a linear transformation, which we also call \(g\), preserving the null vectors (i.e. the preimage of the null quadric).
So \(g\) takes space-like to space-like, or else swaps space-like and time-like.

Suppose that \(g\) swaps space-like and time-like.
Hence the signature is split.
Compose \(g\) with a linear transformation interchanging the space-like and time-like elements of an orthonormal basis, so preserving the null quadric.
So we can arrange that \(g\) takes space-like to space-like.

Pick some space-like \(x\) and time-like \(y\); rescale so that \(x+y\) is null.
Then \(g(x+y)\) is also null, so \(\ip{gx}{gx}=\ip{gy}{gy}\). 
Hence the scale factor \(\ip{gx}{gx}/\ip{x}{x}\) for space-like and time-like \(x,y\) is the same.
Switching which space-like we use, but keeping the same time-like, and vice versa, the scale factor is the same for all vectors.
Composing with a rescaling, \(g\) becomes orthogonal.
\end{proof}
\section{Bundles and product decompositions}
The obvious surjections are \(\Orth{V}\)-equivariant fiber bundles
\[
\pm 1\to Y_+\to X_+, \pm 1\to Y_-\to X_-, \R^{\times}\to Y_0\to X_0
\]
while
\[
\R^{\times}\to U_+\to X_+, \R^{\times}\to U_-\to X_-, \R^{\times}\to Y_0\to X_0
\]
are \(\CO{V}\)-equivariant fiber bundles.

To see the diffeomorphism types of these manifolds, pick a space time splitting \(V=V^+\oplus V^-\).
Let \(|x|=\sqrt{|\ip{x}{x}|}\) for \(x\in V^+\) or \(x\in V^-\).
Let \(S_+\subset V^+\), \(S_-\subset V^-\) be the unit spheres \(|x|=1\).
Then each nonzero vector is uniquely \(x+y\) with \(x\in V^+\), \(y\in V^-\) and being null expands out to \(|x|=|y|\).
Scale to have \(1=|x|=|y|\); call it a \emph{unit} null vector.
So the set of unit null vectors is precisely \(S_+\times S_-\subset Y_0\).
\[
\includegraphics[width=3cm]{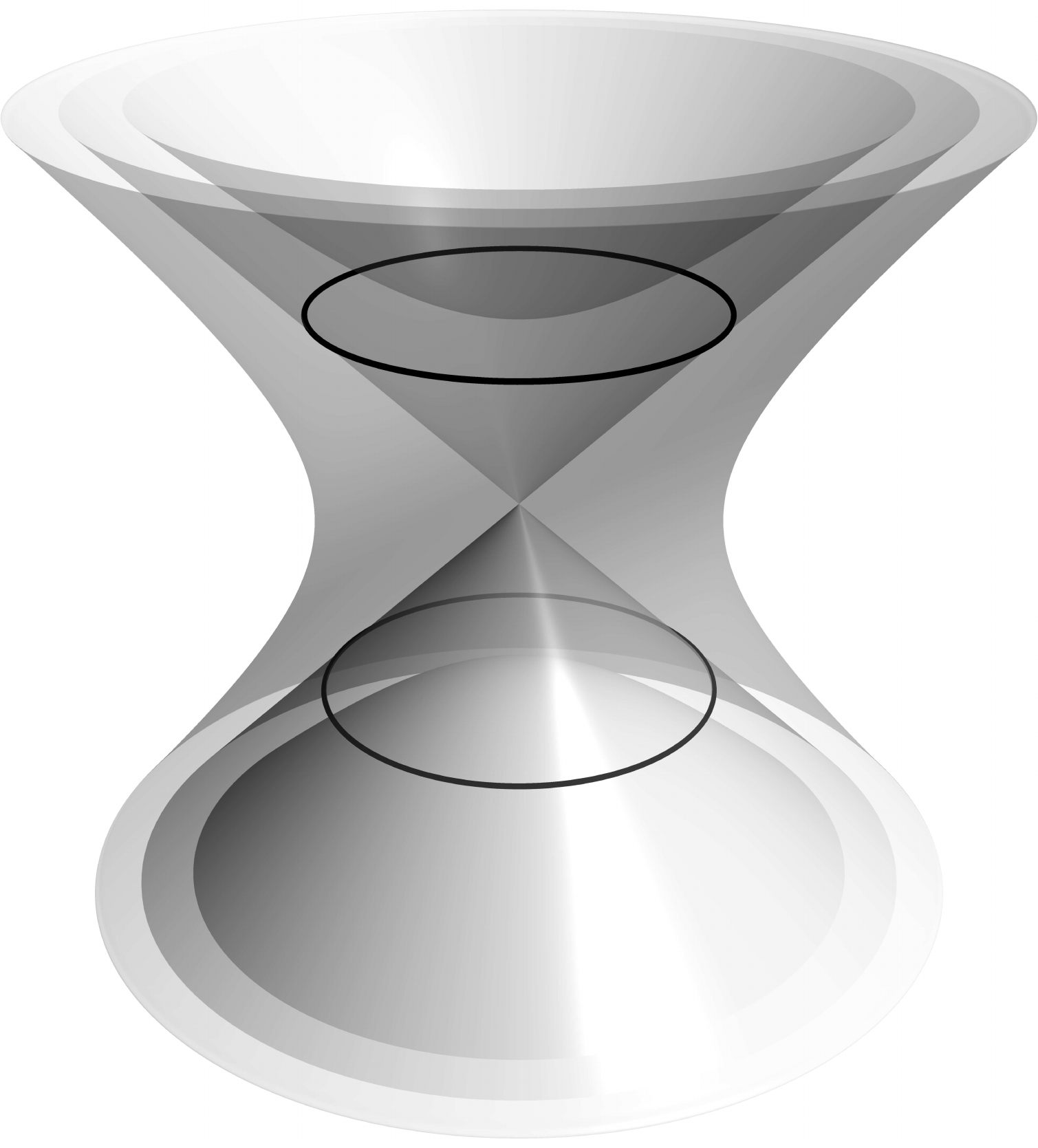}
\]
Every nonzero null vector \((x,y)\) scales to  a unit null vector:
\[
Y_0\cong\R^{+}\times S_+\times S_-.
\]
Hence a covering map
\[
\pm 1 \to S_+\times S_- \to X_0,
\]
taking \((x,y)\) to the line spanned by \(x+y\), so the null quadric \(X_0\) is the quotient of the unit null vectors by the double antipodal map \(x+y\mapsto (-x)+(-y)\).
Note that \(S_+\times S_-\subset V\) is not invariant under \(\Orth{V}\), because its depends on the choice of splitting \(V=V^+\oplus V^-\).

Every element of \(Y_+\) splits uniquely into a sum \(x+y\) with \(x\in V^+, y\in V^-\) with  \(|x|^2=1+|y|^2>0\).
Map
\[
x+y\in Y_+\mapsto (x/|x|,y)\in S_+\times V^-,
\]
and invert:
\[
(u,y)\in S_+\times V^-\mapsto x+y, x=\sqrt{1+|y|^2}u,
\]
hence
\[
Y_+\cong S_+\times V^-.
\]
Similarly
\[
Y_-\cong V^+\times S_-.
\]
Similarly
\[
x+y\in U_+\mapsto (|x|,x/|x|,y)\in \R^+\times S_+\times V^-
\]
is a diffeomorphism.

Every Lie group with finitely many components smoothly retracts to its maximal compact subgroup, and so every homogeneous space of that group retracts to the orbit of that maximal compact subgroup.
So \(\Orth{p,q}\) retracts to \(\Orth{p}\times\Orth{q}\).
In particular, \(\Orth{p,q}\) has one component if \(p=q=0\), two components if \(p=0\) or if \(q=0\) but not both, and otherwise has \(4\) components.
As above, each of \(Y_+,Y_-,Y_0\) retracts \(\Orth{p}\times\Orth{q}\)-equivariantly to the orbits of \(\eep 1,\eem 1,\eep 1+\eem 1\) respectively, i.e. to \(S^{p-1},S^{q-1},S^{p-1}\times S^{q-1}\).
If \(p=1\), \(S^{p-1}=\pm 1\), while if \(p=0\), \(S^{p-1}\) is empty.
\section{Tangent spaces}
Differentiate constant inner product to find that, if \(Y=Y_+,Y_-\) or \(Y_0\) then
\[
T_x Y = x^{\perp}.
\]
Denote by \(\eep{2\dots p}\) the list of vectors \(\eep 2,\dots, \eep p\), etc.
The various tangent spaces and stabilizers \(H:=G^{x_0}\) are:
\[
\begin{array}{ccccc}
\toprule
X&x_0&T_{x_0} X=\operatorname{span}&\text{signature of }X&H\\
\cmidrule(r){1-1}\cmidrule(lr){2-2}\cmidrule(lr){3-3}\cmidrule(lr){4-4}\cmidrule(l){5-5}
Y_+&\eep1&\eep{2\dots p},\eem{1\dots q}&p-1,q&\Orth{p-1,q}\\
Y_-&\eem1&\eep{1\dots p},\eem{2\dots q}&p,q-1&\Orth{p,q-1}\\
Y_0&\eep1+\eem1&\eep1-\eem 1,\eep{2\dots p},\eem{2\dots q}&p-1,q-1&\R^{\times}\Orth{p-1,q-1}\\
%U_0&\eep1+\eem1&\eep{2\dots p},\eem{2\dots q}&p-1,q-1&\\
\bottomrule
\end{array}
\]
On each \(\Orth{p,q}\)-orbit, by \(\Orth{p,q}\)-invariance of the inner product, the signature of the restriction of the inner product to the orbit is the same at all points.
Note that on \(Y_0\), the inner product has a null direction, precisely the rescaling direction.
Under the rescalings, the inner product rescales, so \(X_0\) has only a conformal metric of \(p-1,q-1\) signature.
When we quotient \(Y_+,Y_-\) by \(\pm I\), we swap the choice of preimage of the resulting point in \(X_+,X_-\), but we apply an orthogonal transformation which identifies the tangent spaces at the two points, so the resulting metric is defined on \(X_+,X_-\).
So \(X_+\) has a smooth metric of signature \(p-1,q\), while \(X_-\) has a smooth metric of signature \(p,q-1\), and \(X_0\) has a smooth conformal structure of signature \(p-1,q-1\).
The orthogonal group acts as stabilizers for \(X_+,X_-\), so these metrics are invariant under automorphisms which act transitively on points and orthonormal frames, i.e. on the orthonormal frame bundles.
The stabilizer acting on \(X_+,X_-\), as subgroups  of \(\Orth{V}/\pm 1\), are the isomorphic images of \(\Orth{v^{\perp}}\) with \(v=\eep1,\eem1\) respectively, as we can lift each transformation uniquely to have it fix \(v\).
\section{The Levi--Civita connection}
Each point \(x\in Y_{\pm}\) has normal space the span of \(x\), and hence an orthogonal decomposition into tangent and normal spaces.
Given a curve \(x(t)\in Y_{\pm}\), any vector field \(v(t)\in V\) defined ``along the curve'' decomposes into tangent and normal parts.
The derivative \(\dot{v}(t)\) in the usual Euclidean sense splits into tangent and normal parts as well.
If \(v\) is everywhere tangent (normal), define \(\nabla_{\dot{x}} v\) to be the tangent (normal) part of \(\dot{v}\).
Clearly this definition is invariant under the orthogonal group.
\begin{problem}{w}
This \(\nabla\) operator has value at any time \(t=t_0\) depending only on the vector \(\dot{x}(t_0)\), the tangent vector \(v(t_0)\) and its first derivative \(\dot{v}(t_0)\).
It thus extends uniquely to a connection on the tangent (normal) bundle, and descends to \(X_{\pm}\).
It is torsion free, and compatible with the metric on \(Y_{\pm}\) given by the inner product on \(V\).
Consequently, it is the Levi--Civita connection.
\end{problem}
By invariance under the orthogonal group, this connection is the Levi--Civita connection on \(X_{\pm}\) too.
Take a tangent vector \(v\) to \(Y_+\) at some point \(x\in Y_+\).
Let \(P\) be the span of \(x,v\) in \(V\).
By orthogonal transformation we can arrange that \(x=\eep1\).
If \(v\) is space-like or time-like, we can arrange by orthogonal transformation that \(v\) is a multiple of \(\eep2\) or \(\eem1\).
Hence we can construct the linear reflection map fixing \(x\) and \(v\) and changing the signs of all other vectors in an orthonormal basis, i.e. a reflection fixing all vectors in \(P\) changing signs on all vectors in \(P^{\perp}\).
The curvature of the curve \(Y_+\cap P\) is thus invariant, but lies in \(P^{\perp}\), so vanishes: a geodesic.
Similarly if \(v\) is null, we can arrange that \(v=\eep2+\eem2\).
Linearly transforming the basis \(\eep2,\eem2\) by 
\[
\begin{pmatrix}
\cosh t&\sinh t\\
\sinh t&\cosh t
\end{pmatrix}
\]
scales \(v\) by \(e^t\) and the other null line in the span of \(\eep2,\eem2\) by \(e^{-t}\).
Hence the curvature is invariant, so vanishes: a geodesic.
Similarly this works for \(Y_-\).
This is only one geodesic in each direction.
So we have all of the geodesics: a curve on \(Y_{\pm}\) is a geodesic just when, at some point, it is the intersection with a plane containing the normal vector at that point.
Mapping to projective space, the geodesics on \(X_{\pm}\) are its intersection with projective lines.
\[
\includegraphics[width=3cm]{hyperboloid}
\]
\begin{example}
In \(\R^{2,1}\), the surface \(Y_-\) is the hyperboloid of \(2\) sheets: \(x^2+y^2+1=z^2\), with a Riemannian metric invariant under the orthogonal group.
The vertical plane through \(\eep1\) intersects the cone at two points, so the geodesic is the hypebola in between.
All the geodesics in \(Y_-\) are carried to one another by isometries.

The surface \(Y_+\) is the hyperboloid of \(1\) sheet: \(x^2+y^2=z^2+1\).
Every point \(x\in Y_+\) can be brought to the point \(\eem1\) by orthogonal transformation.
The horizontal plane through \(\eem1\) and \(\eem2\) intersects \(Y_+\) along a time-like geodesic circle.
The plane through \(\eem1\) and either \(\eep1+\eem1\) or \(-\eep1+\eem1\) intersects \(Y_+\) along a null geodesic: the lines \(x=\pm 1\).
The plane through \(\eem1\) and \(\eep1\) intersects \(Y_-\) on a space-like geodesic hyperbola.
All pointed geodesics are brought to one of these by orthogonal transformation.
\end{example}
\begin{theorem}
The group of pseudo-Riemannian isometries of de~Sitter and anti-de~Sitter spaces \(X_+, X_-\) is precisely the projective orthogonal group.
\end{theorem}
\begin{proof}
The projective orthogonal group acts transitively.
Consider the stabilizer of a point.
Looking at our description of tangent spaces above, the stabilizer of a point acts transitively on orthogonal frames at that point, i.e. as the orthogonal group of that tangent space.
Take any isometry.
Composing it with an element of the projective orthogonal group, we can arrange that it fixes a given point, and a frame at that point, i.e. acts trivially on the tangent space at that point.
It commutes with the exponential map, so preserves an open set.
By analyticity and connectivity, this isometry is the identity map.
\end{proof}
\section{Stabilizer of the null quadric}
The description of the stabilizer of a point of the null quadric \(X_0\) is more complicated.
Change the variables of our standard example to get the axes of the first and last variables to be null, so \(V=\R\oplus V'\oplus\R\) with \(V'\) of signature \(p-1,q-1\).
If the signature is split, so is that of \(V'\); take a conformal transformation \(T'\in\CO{V'}\) changing sign of the inner product, and let
\[
T=
\begin{pmatrix}
1&0&0\\
0&T'&0\\
0&0&-1
\end{pmatrix}
\]
so that \(T\) changes the sign of the inner product, and fixes the standard basis vector \(e_1\).
Let \(P_O\subseteq\Orth{V}\) be the group of orthogonal transformations which stabilizer the span of the null vector \(e_1\) in these coordinates.
Let \(P_C\subseteq\CO{V}\) be the group of conformal transformations which stabilizer the span of the null vector \(e_1\) in these coordinates.
\begin{lemma}
In this basis, the  group \(P_O\) consists precisely of the matrices of the form
\[
g=\lambda
\begin{pmatrix}
a&a v^* h&\frac{a}{2}\ip{v}{v}\\
0&h&v\\
0&0&a^{-1}
\end{pmatrix}
\]
where \(a\in\R^{\times}\), \(v\in V'\), \(h\in\Orth{V'}\), \(v^*\in (V')^*\) the dual  in the inner product.
The group \(P_C\) consists precisely of the matrices  of the form \(\lambda g\), \(\lambda>0\) and \(g\in P_O\).
\end{lemma}
\begin{proof}
After a positive scaling, these conformal transformations are orthogonal transformations.
So assume that \(g\) is orthogonal and preserves the line spanned by the null vector \(e_1\) (in the standard basis) and so preserves the perpendicular space of that null vector, i.e. the span of the entire standard basis with the last vector omitted.
Hence \(g\) has the form
\[
g=
\begin{pmatrix}
a&\xi&t\\
0&h&v\\
0&0&b
\end{pmatrix}
\]
for some \(a,b\in\R^{\times}\), \(h\in\GL{V'}\), \(v\in V'\), \(\xi\in(V')^*\).

Write the inner product on \(V'\) as 
\[
\ip{v}{w}=\ip{v}{\eta w}_{\text{Euclid}},
\]
for a symmetric matrix \(\eta\) of signature \(p-1,q-1\).
Let
\[
S:=
\begin{pmatrix}
0&0&-1\\
0&\eta&0\\
-1&0&0
\end{pmatrix}.
\]
Plug in to \(g^tSg=S\) to find exactly the form specified.
\end{proof}
Our goal is to produce homogeneous pseudo--Riemannian manifolds and conformal manifolds, so now we will change the indices \(p,q\).
Our null quadrics \(X_0\) are conformal of signature \(p-1,q-1\), so if we just add one to \(p,q\), we can denote \(X_0\) as
\[
\begin{tikzcd}
\pm 1\arrow[r]&S^p\times S^q\arrow[r]&X^{p,q}_0=\CO{p+1,q+1}/P_C.
\end{tikzcd}
\]
To make a pseudo--Riemannian manifold of signature \(p,q\), we can choose \(X_+\) with \(V=\R^{p+1,q}\) or \(X_-\) with \(V=\R^{p,q-1}\).
So we find \emph{de~Sitter space}\define{de~Sitter space} and \emph{anti-de~Sitter space}:\define{anti-de~Sitter space}
\[
\begin{tikzcd}
\pm1\arrow[r]&Y_+^{p,q}=S^p\times\R^q=\Orth{p+1,q}/\Orth{p,q}\arrow[d]\\&X_+^{p,q}=(S^p\times\R^q)/\pm 1=(\Orth{p+1,q}/\pm 1)/\Orth{p,q}\\
\pm1\arrow[r]&Y_-^{p,q}=\R^p\times S^q=\Orth{p,q+1}/\Orth{p,q}\arrow[d]\\&X_-^{p,q}=(\R^p\times S^q)/\pm 1=(\Orth{p,q+1}/\pm 1)/\Orth{p,q}
\end{tikzcd}
\]
Swapping the sign of the quadratic form on \(V\),
\begin{align*}
X_+^{p,q}&\cong X_-^{q,p},\\
X_0^{p,q}&\cong X_0^{q,p}.
\end{align*}

\chapter{Homogeneous geometries with given symmetries}\label{chapter:homogeneous.examples}
Let us try to find all homogeneous Cartan geometries with a given model on a given homogeneous space.
This is quite different to section~\ref{section:homog.aff.surf}, where we look for all homogenous Cartan geometries with a given model, on all possible homogeneous spaces.
\begin{example}
The flat torus has a flat Riemannian geometry invariant under translations, hence a homogeneous projective connection.
Similarly the conformal geometry of a flat torus of dimension \(3\) or more is homogeneous.
\end{example}
Take homogeneous spaces \((X,G)\) and \((X',G')\) with \(\dim X=\dim X'\), say \(X=G/H\), \(X'=G'/H'\).
We want to build \((X,G)\)-geometries on \(X'\), invariant under \(G'\). 

A \emph{die}\define{die} is a pair of
\begin{itemize}
\item
a Lie group morphism \(k\in H'\to\bar{k}\in H\) and
\item
a linear map \(A\in\LieG'\mapsto\bar{A}\in\LieG\) so that
\begin{itemize}[label=$\diamond$]
\item
the die is \(H'\)-equivariant: \(\overline{\Ad}_k=\Ad_{\bar{k}}\bar{\,}\) for \(k\in H'\) and
\item
the die restricted to  \(\LieH'\) is the induced Lie algebra morphism \(\LieH'\to\LieH\) and
\item
\(A\in\LieG'/\LieH'\mapsto\bar{A}\in\LieG/\LieH\) is a linear isomorphism.
\end{itemize}
\end{itemize}
Consider the following recipe: given a die, let
\[
\omega_\G:=\Ad_h^{-1}\bar\omega_{G'}+\omega_H
\]
on \(G'\times H\).
We conjugate the die by replacing the Lie group morphism by
\[
k\mapsto \Ad_{h_0} \bar{k},
\]
and \(\bar{\,}\) by \(\Ad_{h_0}\bar{\,}\), for some \(h_0\in H\).
Let \(\G:=\amal{G'}{H'}{H}\).
\begin{theorem}[Kobayashi \cite{Kobayashi:1956},Hammerl \cite{Hammerl:2007}]
Every die generates a homogeneous Cartan geometry, by the recipe:  \(\omega_{\G}\) descends to a Cartan connection on \(\G\) for an \((X,G)\)-geometry on \(X'\), invariant under \(G'\).
Every homogeneous Cartan geometry arises from a die, uniquely up to conjugating the die.
\end{theorem}
\begin{proof}
Suppose that \(H\to\G\to X'\) is a \(G'\)-invariant Cartan geometry.
By theorem~\vref{theorem:build.homo.princ.bun.2}, \(\G=\amal{G'}{H'}{H}\) for a Lie group morphism \(H'\to H\).
By theorem~\vref{thm:homog.principal.bundle.isos}, this Lie group morphism is uniquely determined up to conjugation in \(H\).
To determine a Cartan geometry, we need an invariant connection on \(\G_G=\amal{G'}{H'}{G}\).
By theorem~\vref{theorem:invariant.connection}, this is precisely determined by an extension of our Lie group morphism to an \(H'\)-equivariant linear map \(\LieG'\xrightarrow{a}\LieG\), as in our recipe.
To give a Cartan geometry, this connection cannot vanish on any \(\G\)-tangent vectors.
This is precisely that
\[
C+\Ad_h^{-1}\bar B\ne 0
\]
for all \(B\in\LieG'\) and \(C\in\LieH\) unless \(B\in\LieH'\) and \(C=-\bar{B}\).
Mod out \(\LieH\): this is precisely that \(\LieG'/\LieH'\to\LieG/\LieH\) is a linear isomorphism.
The construction of the connection is unique up to conjugation.
\end{proof}

Homogeneous Cartan geometries need not be complete: the torus and the plane have homogeneous flat incomplete projective connections.
\section{Example: projective connections on the sphere}\label{section:proj.conn.sphere}
\begin{theorem}\label{thm:inv.proj.conn.sphere}
Up to isomorphism, there is a \(1\)-parameter family of projective connections on the sphere invariant under rotation, with die
\[
\begin{pmatrix}
0 & -c^t\\
c & d
\end{pmatrix}
\in\LieSO{n+1}
\mapsto
\begin{pmatrix}
0&-\beta c^t\\
c&d
\end{pmatrix}
\in\LieSL{n+1}
\]
for an arbitrary real constant \(\beta\).
\end{theorem}
\begin{proof}
Consider the projection connections on the sphere invariant under rotations:
\begin{align*}
(X',G')&=(S^n,\SO{n+1}),\\
(X,G)&=(\RP{n},\PGL{n+1}).
\end{align*}
The stabilizer subgroups are 
\begin{align*}
H'&=\SO{n},\\
H&=
\set{
\begin{bmatrix}
a&b\\
0&d
\end{bmatrix}
}.
\end{align*}
A die is an \(H\)-equivariant linear map
\(
\LieG'\to\LieG,
\)
which restricts to a Lie algebra morphism \(\LieH'\to\LieH\), from a given Lie group morphism, and the induced linear map
\[
\LieG'/\LieH'\to\LieG/\LieH,
\]
is an \(H'\)-equivariant linear isomorphism.
But the die is determined only up to conjugation inside \(H\).

Since \(H\) has maximal compact subgroup the orthogonal group, and this subgroup is unique up to conjugation, \(H'\to H\) has image inside the obvious orthogonal group.
Since 
\[
\LieG'/\LieH'\to\LieG/\LieH,
\]
is an \(H'\)-equivariant linear isomorphism, and \(H'\) acts on \(\LieG'/\LieH'\) faithfully, \(H'\to H\) is injective.
All automorphisms of the special orthogonal Lie algebra are inner, so we can arrange by conjugation in \(H\) that \(H'\to H\) is the standard inclusion.
We extend the die to \(\LieG'=\LieSO{n+1}=\LieSO{n}\oplus\R^n\), so that it is a linear isomorphism on the quotient
\[
\LieG'/\LieH'=\R^n\to\R^n=\LieG/\LieH.
\]
Split 
\[
\LieG'=\LieH'\oplus\R^n
\]
as \(H'\)-modules corresponding to writing each element of \(\LieG'\) as
\[
\begin{pmatrix}
0 & -c^t\\
c & d
\end{pmatrix}
\]
with \(d\in\LieH'\), \(c\in\R^n\).
Writing every element of \(\LieG\) as
\[
\begin{pmatrix}
a&b\\
c&d
\end{pmatrix}
\]
in blocks of sizes
\[
\begin{bNiceMatrix}[margin,first-row,first-col]
 & 1 & n \\
1 & \ & \ \\
n & \ & \ \\
\end{bNiceMatrix}
\]
Split these, over \(H'\), into irreducibles
\[
\begin{pmatrix}
-nt&b\\
c&tI+d+e
\end{pmatrix}
\]
where \(d\) is symmetric and traceless, \(e\) antisymmetric.

Note that \(\R^n,\Lm{2}{\R^n}\) and \(\SymZ{2}{\R^n}\) are irreducible \(H'\)-modules.
By Schur's lemma, an invariant linear map between irreducible modules is zero unless the modules are isomorphic.
This determines all of the maps except \(b=b(c)\); \(H'\)-equivariance (since \(H'\) contains reflections in all hyperplanes) easily shows that \(b=\beta c\) for a unique real constant \(\beta\).
From this observation, the map is easily seen to be, up to \(H'\)-equivariant isomorphism,
\[
\begin{pmatrix}
0 & -c^t\\
c & d
\end{pmatrix}
\mapsto
\begin{pmatrix}
0&-\beta c^t\\
c&d
\end{pmatrix}
\]
for some real constant \(\beta\).
\end{proof}

\section{Remarks}
The \emph{codimension}\define{codimension} of a submanifold is the difference between the dimension of the ambient manifold and the dimension of the submanifold.
The \emph{homogeneity}\define{homogeneity} of a Lie group action is the maximal orbit dimension; the \emph{cohomogeneity}\define{cohomogeneity} of a Lie group action is the codimension of such an orbit. 
Cohomogeneity \(1\) Cartan geometries on connected manifolds are probably described by some ``linear algebra data'' as above, but satisfying some ordinary differential equations.
To the best of my knowledge, nobody has found out how to express those equations.

A Cartan geometry is \emph{quasihomogeneous}\define{quasihomogeneous!Cartan geometry}\define{Cartan!geometry!quasihomogeneous} if it is homogeneous on a dense open set.
The only result I know in this direction is the classification of germs of real analytic torsion-free quasihomogeneous affine connections on surfaces \cite{Dumitrescu:2014}.
To the best of my knowledge, the classification of quasihomogeneous Cartan geometries with a given model is unknown for every model.

\chapter{Homogeneous spaces again}
\section{Morphisms}
A \emph{morphism}\define{morphism!of homogeneous spaces}\define{homogeneous space!morphism} \((X,G)\xrightarrow{(\varphi,\Phi)}(X',G')\) of homogeneous spaces  is a smooth map \(X \xrightarrow{\varphi}X'\) equivariant for a Lie group morphism \(G\xrightarrow{\Phi}G'\).
\begin{problem}{when.morphism}
Take a Lie group morphism \(G\xrightarrow{\Phi}G'\) and a pair of points  \(x\in X\) and \(x'\in X'\).
Prove that \(\Phi(G^x)\subseteq (G')^{x'}\) if and only if \(\Phi\) arises from a morphism of homogeneous spaces taking \(x\) to \(x'\), and then the morphism is determined by \(\varphi(gx)=\Phi(g)x'\).
\end{problem}
\begin{example}
Let \(X\) be the unit sphere of dimension \(n\): \(X=(\R^{n+1}-\set{0})/\R^+\), let \(G\) be the group of linear isomorphisms of \(\R^{n+1}\) of determinant \(\pm 1\).
Let \(X':=\RP{n}\) be the real projective space of dimension \(n\), i.e. the set of all real lines through the origin on \(\R^{n+1}\), and \(G':=\PGL{n+1}\) be the set of all real projective transformations of \(X'\), i.e. \(G':=G/\set{\pm I}\).
Then \(X\to X'\) is a \(2-1\) covering map equivariant for \(G\to G'\), which is also a \(2-1\) covering map.
\end{example}
The \emph{product}\define{product!of homogeneous spaces}\define{homogeneous space!product} of homogeneous spaces 
\[
(X_1,G_1), (X_2,G_2), \dots, (X_k,G_k)
\]
is \((X,G)\) where
\[
X=\prod_i X_i, G=\prod_i G_i,
\]
with the obvious product action.
\section{Homogeneous Riemannian manifolds}
Every homogeneous Riemannian manifold has constant injectivity radius, convexity radius, and so on.
As you travel along a geodesic, wherever you are, your geodesic keeps going up to that injectivity radius: every homogeneous Riemannian manifold is complete.
\begin{example}
Euclidean space, hyperbolic space, the sphere, the flat torus, the real, complex and quaternionic projective spaces and Grassmannians, and their products, with their standard metrics, are homogeneous.
\end{example}
A homogeneous space \((X,G)\) is \emph{proper}\define{proper!homogeneous space}\define{homogeneous!space!proper} if 
\[
(g,x)\in G\times X\mapsto (gx,x)\in X\times X
\]
is a proper map.
Equivalently, if a sequence \(x_1,x_2,\dots\in X\) converges and a sequence
\[
g_1x_1,g_2x_2,\dots\in X
\]
converges then some subsequence of \(g_1,g_2,\dots\in G\) converges.
\begin{problem}{r}
For a strong effective homogeneous space \((X,G)\), the following are equivalent:
\begin{itemize}
\item
\((X,G)\) is proper
\item
\(G\) is a closed subgroup of the isometry group of a \(G\)-invariant metric imposing the usual topology on \(X\)
\item
the stabilizer subgroup of some point is compact
\item
the stabilizer subgroup of every point is compact 
\item
\(G\) is a closed subgroup of the isometry group of a \(G\)-invariant Riemannian metric on \(X\).
\end{itemize}
\end{problem}
\begin{example}\label{example:unit.ball}
The unit ball 
\[
X:=\set{z\in\C^n\,|\,|z|=1},
\]
is homogeneous under the group \(G\) of projective linear transformations of the form
\[
z\in X\mapsto\frac{a-\pi_az-\sqrt{1-|a|^2}\pi^{\perp}_az}{1-\left<z,a\right>}\in X
\]
where \(\left<z,a\right>\) is the usual Hermitian inner product in \(\C^n\), \(a\in X\), \(\pi_a\) is orthogonal projection to the complex span of \(a\), and \(\pi^{\perp}_a=I-\pi_a\) is orthogonal projection to \(a^{\perp}\) \cite{Rudin:2008}.
The subgroup \(H\subseteq G\) fixing \(0\in X\) is precisely the group of unitary linear transformations \(H=\Un{n}\).
\end{example}
\begin{example}
For any integer \(g\ge 1\), the \emph{Siegel upper half space}\define{Siegel upper half space} \(X\) is the set of all \(g\times g\) complex matrices with positive definite imaginary part \cite{Klingen:1990} p. 1, \cite{Siegel:3} chapter 6.
The Siegel upper half space arises naturally: the period matrix of any compact Riemann surface, in a symplectic basis of the homology, is such a matrix.
The group of \(G\) of real linear symplectic matrices (the \emph{symplectic group})\SubIndex{symplectic group} acts transitively on \(X\) by linear fractional transformations: if 
\[
g=
\begin{pmatrix}
A&B\\
C&D
\end{pmatrix}
\]
then
\[
Z\mapsto gZ:=(AZ+B)(CZ+D)^{-1}.
\]
The stabilizer of the symmetric matrix \(X=\sqrt{-1}I\) is isomorphic to the unitary group, compact, so a proper homogeneous space.
\end{example}
\begin{example}
Affine space \((X,G)=(\R^n,\GL{n}\ltimes\R^n)\) has stabilizer \(G^0=\GL{n}\), not compact: an improper strong effective homogeneous space.
The action of the stabilizer on the tangent space is by arbitrary linear maps, some of which have eigenvalues not unit complex numbers, so the action does not preserve a Riemannian metric.
\end{example}
\begin{example}
The torus \(X=\R^2/\Z^2\) is acted on by the group \(G=\GL{2,\Z{}}\ltimes\R^2\) of affine transformations with integer coefficient linear part: each \(g\in G\) acts on the plane as
\[
v\mapsto Av+b,
\]
so if add an integer vector \(z\) to \(v\),
\[
v+z\mapsto (Av+b)+Az.
\]
The stabilizer \(H=\GL{2,\Z{}}\) contains
\[
\begin{pmatrix}
\pm 1&n\\
0&\pm 1
\end{pmatrix}
\]
for any integer \(n\), so \(H\) has noncompact closure, and \(X\) is an improper strong effective compact homogeneous space with no \(G\)-invariant Riemannian metric.
\end{example}
\begin{example}
Pick an integer \(n\ge 3\) and let \(X\) be the set of pairs \((u,v)\) of perpendicular unit vectors in \(\C^n\).
Pick real numbers \(\alpha,\beta\), neither a rational multiple of the other. 
Let \(G:=\SU{n}\times\R\) be the group of transformations of \(X\) generate by (i) \(\SU{n}\) acting on \(\C^n\) and (ii) maps \((u,v)\mapsto (e^{i\alpha t}u,e^{i\beta t}v)\) for \(t\in\R\).
Each pair \((u,v)\in X\) is a point \((u,v)\in\C^{2n}\), making \(X\subseteq\C^{2n}\) a compact submanifold.
Already \(\SU{n}\) acts transitively.
The Hermitian metric on \(\C^{2n}\) pulls back to a  \(G\)-invariant metric on \(X\).
But \(G\) acts improperly on \(X\), not being closed in the isometry group of \(X\), and also not having compact stabilizers.
\end{example}
\begin{problem}{s}
Prove that, for any proper homogeneous space \((X,G)\) with \(n:=\dim X\), the space of \(G\)-invariant metrics is finite dimensional, of dimension at most \(n(n+1)/2\).
For which \((X,G)\) does it reach this bound?
\end{problem}
\section{Holomorphy}
A \emph{biholomorphism}\define{biholomorphism} of a complex manifold is a diffeomorphism which is holomorphic with holomorphic inverse.
A \emph{homogeneous complex space}\define{homogeneous!space!complex} is a homogeneous space \((X,G)\) with \(X\) a complex manifold and \(G\) acting by biholomorphism of \(X\).
A \emph{complex homogeneous space}\define{complex homogeneous space} is a homogeneous complex space so that \(G\) is a complex Lie group acting holomorphically on \(X\).
\begin{problem}{t}
Under its biholomorphism group, is the Riemann sphere a complex homogeneous space? The complex plane \(\C\)? The unit disk in the complex plane? An annulus in the complex plane? \(\C^2\)?
\end{problem}
\section{The frame bundle}\label{section:frame.bundle}
The \emph{frame bundle} of a manifold is the set of all linear isomorphisms of tangent spaces of the manifold with some fixed vector space. 
In other words: take a manifold \(M\) and a vector space \(V_0\) of the same dimension as \(M\).
The \(V_0\)-valued \emph{frame bundle}\define{frame bundle} \(\framebundle{M}\) of \(M\) is the set of all pairs \((m,u)\) where \(m\in M\) and \(T_m M\xrightarrow{u}V_0\) is linear isomorphism.
The \emph{projection map}\define{projection map} \(\framebundle{M}\xrightarrow{\pi}M\) is the map \(\pi(m,u)=m\).
The group \(\GL{V_0}\) acts on \(\framebundle{M}\) by the right action
\[
(m,u)g:=(m,g^{-1}\circ u).
\]
\begin{problem}{frame.bundle.bundle}
Prove that the frame bundle \(\framebundle{M}\xrightarrow{\pi}M\) is a principal right \(\GL{V_0}\)-bundle.
\end{problem}
\section{Effective homogeneous spaces}\label{subsection:effective.homog}
The \emph{kernel}\define{homogeneous!space!kernel}\define{kernel!of a homogeneous space} \(K\subseteq G\) of a homogeneous space \((X,G)\) is the set of elements of \(G\) fixing every element of \(X\), so \((X,G/K)\) is effective, the \emph{associated effective}\define{associated!effective!homogeneous space}\define{homogeneous!space!associated effective} homogeneous space.
Since \(K\) fixes all points of \(X\), \(K\) lies in the stabilizer \(G^{x_0}\) of any point of \(X\).
Lying in all the other stabilizers \(G^{gx_0}=gG^{x_0}g^{-1}\), \(K\) is normal in \(G^{x_0}\).
Reversing this argument, any subgroup of \(G^{x_0}\) normal in \(G\) fixes all points of \(X\), so \(K\) is the maximal subgroup of \(H\) normal in \(G\).
Being the intersection of the stabilizer subgroups, which are closed subgroups of \(G\), \(K\) is closed in \(G\) and in \(H\).
If \(K\) is discrete in \(G\), \((X,G)\) is \emph{almost effective};\define{almost effective homogeneous space}\define{homogeneous!space!almost effective} if \(K=\set{1}\), \((X,G)\) is \emph{effective}.\SubIndex{effective homogeneous space}\SubIndex{homogeneous space!effective}
\begin{problem}{u}
Prove that a homogeneous space with finite kernel is proper if and only if its associated effective homogeneous space is proper.
\end{problem}
\section{Example: the kernel of real projective space}
Let \((X,G)=(\RP{n},\SL{n+1})\) so \(H\) is the group of matrices with block sizes:
\[
\begin{pNiceMatrix}[margin,first-row,first-col]
 & 1 & n \\
1 & \Block[fill=gray!40]{1-2}{} \ & \ \\
n & 0 & \Block[fill=gray!40]{1-1}{} \ \\
\end{pNiceMatrix}
\]
i.e. preserving the point 
\[
x_0:=
\begin{bmatrix}
1\\
0\\
\vdots\\
0
\end{bmatrix}\in X
\]
The elements of \(H\) act on elements of \(\LieG\) in the adjoint action as
\[
\begin{pmatrix}
a&b\\
0&d
\end{pmatrix}
\begin{pmatrix}
A&B\\
C&D
\end{pmatrix}
\begin{pmatrix}
a&b\\
0&d
\end{pmatrix}^{-1}.
\]
\begin{problem}{v}
For any homogeneous space \((X,G)\) explain why the \(H\)-module isomorphism \(\LieG\cong T_1 G\) descends to an \(H\)-module isomorphism \(\LieG/\LieH\cong T_{x_0} X\).
\end{problem}
The tangent space \(T_{x_0}X=\LieG/\LieH\) is acted on then by
\[
\begin{pmatrix}
a&b\\
0&d
\end{pmatrix}
\begin{pmatrix}
*&*\\
C&*
\end{pmatrix}
\begin{pmatrix}
a&b\\
0&d
\end{pmatrix}^{-1}
=
\begin{pmatrix}
*&*\\
dCa^{-1}&*
\end{pmatrix}
\]
So the elements of \(H\) which act trivially on \(T_{x_0} X\) are precisely the subgroup \(H_1\subseteq H\) of matrices with \(dCa^{-1}=C\) for all vectors \(C\in\R^n\), i.e. those of the form
\[
\begin{pmatrix}
a&*\\
0&aI
\end{pmatrix},
\]
with unit determinant, i.e. \(a^{n+1}=1\).
This subgroup \(H_1\subseteq G\) is precisely the subgroup preserving the point \(x_0\in X\) and the standard basis for \(T_{x_0}X\), so \(X_1:=G/H_1\) is the frame bundle of \(X\), i.e. the set of choices of point \(x\in X\) and linear isomorphism \(T_x X\xrightarrow{u}\R^n\).
To act trivially on \(X\), an element of \(G\) must lie inside \(H_1\), but must also act trivially on the frame bundle of \(X\), i.e. on \(X_1\).
By the same reasoning, it must act trivially on \(T_{u_0} X_1\) where \(u_0\) is the standard basis of \(\R^n\) as a frame on \(X\).
But \(T_{u_0}X_1=\LieG/\LieH_1\), and we compute that \(H_1\) acts on \(\LieG/\LieH_1\) by
\[
\begin{pmatrix}
a&b\\
0&aI
\end{pmatrix}
\begin{pmatrix}
A&*\\
C&D
\end{pmatrix}
\begin{pmatrix}
a&b\\
0&aI
\end{pmatrix}^{-1}
=
\begin{pmatrix}
A+a^{-1}bC&*\\
C&-a^{-1}Cb+D
\end{pmatrix}
\]
Let \(H_2\) be the elements of \(H_1\) which act trivially, i.e. 
\[
0=a^{-1}bC=a^{-1}Cb
\]
for all \(C\), and hence \(b=0\), so \(H_2\) consists precisely in the group of matrices
\[
\begin{pmatrix}
a&0\\
0&aI
\end{pmatrix}
\]
with determinant \(1\), i.e. \(a^{n+1}=1\).
This is clearly the kernel, as these are precisely the linear transformations trivial as projective transformations.
\section{Finding the kernel}
The kernel is the smallest in a nested sequence of subgroups.
Suppose that \(G\) is a Lie group.
To each linear subspace \(\LieL\subseteq\LieG\) and closed subgroup \(H\subseteq G\), associate the closed subgroup \(H_{\LieL}\) of elements \(h\in H\) for which 
\[
(\Ad_h-I)\LieG\subseteq\LieL,
\]
or equivalently \(\LieL\) is \(H_{\LieL}\)-invariant and \(\Ad_h\) is trivial on \(\LieG/\LieL\) for \(h\in H_{\ell}\).
\begin{lemma}\label{lemma:eff.hom}
If \(\LieL\subseteq\LieH\) is an \(H\)-invariant linear subspace then \(H_{\LieL}\subseteq H\) is a closed normal subgroup.
\end{lemma}
\begin{proof}
For any \(h\in H\) and \(a\in H_{\LieL}\),
\begin{align*}
\Ad_{ha h}-I
&=
\Ad_h\Ad_{a}\Ad_h^{-1}-I,
\\
&=
\Ad_h(\Ad_{a}-I)\Ad_h^{-1},
\end{align*}
takes \(\LieG\) to \(\Ad_h \LieL=\LieL\).
\end{proof}
Let
\[
H_0:=H, \ H_1:=H_{\LieH_0},\ H_2:=H_{\LieH_1}, \dots.
\]
Since \(H_i\subseteq H\) is normal, \(\LieH_i\subseteq\LieH\) is an ideal, so \(H_{i+1}\subseteq H\) is a closed normal subgroup.
\begin{lemma}\label{lemma:kernel}
Suppose that \((X,G)\) is a strong homogeneous space.
Pick a point \(x_0\in X\) and let \(H:=G^{x_0}\).
Then the kernel of \((X,G)\) is
\[
\bigcap_i H_i.
\]
\end{lemma}
\begin{proof}
Every element \(h\) of the kernel of \((X,G)\) acts trivially on \(X\), hence preserves \(\LieH\) and acts trivially on \(\LieG/\LieH=T_{x_0} X\), so lies in \(H_1=H_{\LieH_0}\).
So \(X_1:=G/H_1\) is identified with a \(G\)-invariant set of linear isomorphisms \(T_x X\to\LieG/\LieH\).
So the kernel acts trivially on \(X_1\). 
By induction, the kernel of \((X,G)\) lies in this intersection.

Denote by \(G^0\) the identity component of \(G\), i.e. the path component containing the identity element.
Let \(G'\subseteq G\) be the subgroup generated by \(G^0\) and \(H\).
Let \(G''\) be the union of all components of \(G\) of the form \(hG^0\) for \(h\in H\).
Claim: \(G'=G''\).
Proof: Clearly \(H,G^0\subseteq G''\).
Each component \(hG^0\) lies in \(G'\), so \(G''\subseteq G'\).
For any \(h,h'\in H\), \(hG^0h'G^0\) is a connected subset of \(G\) containing \(hh'G^0\), which is a component of \(G\), so
\[
hG^0h'G^0\subseteq hh'G^0.
\]
The reverse inclusion is also clear.
Hence \(G''\) is a subgroup of \(G\).
So finally \(G'=G''\).

Suppose that \(N\) is a subgroup of \(H\) normal in the group \(G'\)
Then for \(g\in G'\) and \(n\in N\),
\[
g^{-1}ng\in N 
\]
so
\[
ng\in gN
\]
so
\[
n(gH)=(ng)H=gNH=gH.
\]
Hence \(N\) fixes every point of \(X':=G'/H\subseteq G/H=X\), a union of components of \(X\).
But \((X,G)\) is strong, so \(N\) lies in the kernel of \((X,G)\).
So the kernel of \((X,G)\) contains every subgroup of \(H\) normal in \(G'\).

Let \(K:=\bigcap_i H_i\); being a closed subgroup, \(K\) is also an embedded Lie subgroup \cite{Mimura/Toda:1991} p. 44.
The Lie algebra \(\LieK\) of \(K\) is the smallest of the nested \(\LieH_i\).
For any \(A\in\LieG\), \(h\in K\),
\[
\Ad_h A-A\in\LieK.
\]
For \(A\in\LieG/\LieK\), 
\[
\Ad_hA=A.
\]
If we exponentiate in \(G/K\),
\[
\Ad_h e^{tA}=e^{tA}.
\]
If \(g\) is in the identity component \((G/K)^0\subseteq G/K\), write \(g\) as a product of elements of the form \(e^{tA}\) to find that
\[
hgh^{-1}g^{-1}=1.
\]
For \(g\) in the identity component \(G^0\subseteq G\),
\[
hgh^{-1}g^{-1}\in K.
\]
So \(K\) is normal in \(G^0\) and in \(H\), so normal in \(G'\).
In particular, \(K\) lies in the kernel of \((X,G)\).
\end{proof}

\section{Jets and homogeneous spaces}
\subsection{Jets}
Take manifolds \(P,Q\), a point \(p_0\in P\) and an integer \(r\ge 0\).
Two \(C^r\) maps from neighborhoods of \(p_0\) to \(Q\) are \emph{\(r\)-jet equivalent}\define{jet!equivalence} if they both map \(p_0\) to the same point \(q_0\) and, in some coordinates on \(P\) near \(p_0\) and on \(Q\) near \(q_0\), they agree up to and including order \(r\) terms in a Taylor series.
The equivalence classes are \(r\)-\emph{jets}.\define{jet}
By the chain rule, \(r\)-equivalence is independent of coordinates.
\begin{example}
For one function of one variable, say
\[
y=a_0+a_1x+\dots,
\]
we have Taylor coefficients
\[
a_r:=\frac{1}{r!}\left.\frac{d^ry}{dx^r}\right|_{x=0},
\]
which we can denote \(a_r(0)\) to indicate that it arises at \(x=0\).
These coefficients are related by
\begin{align*}
a_{r+1}(0)
&=
\frac{1}{(r+1)!}
\left.\frac{d}{dx}\right|_{x=0}\left.\frac{d^r}{dt^r}\right|_{t=x}y(t),
\\
&=
\frac{1}{r+1}
\left.\frac{d}{dx}\right|_{x=0}a_r(x),
\end{align*}
so that the \((r+1)\)-jet is identified with the \(1\)-jet of the \(r\)-jet.
\end{example}
The same holds, with more notation, in any number of variables: every \((r+1)\)-jet is a \(1\)-jet of an \(r\)-jet.
A family of \(r\)-jets can vary with arbitrary \(1\)-jet, depending on a parameter, but arises as an \((r+1)\)-jet just when it satisfies symmetry in partial derivatives, by the Poincar\'e lemma.
\subsection{Homogeneous spaces}
\begin{lemma}
Take a homogeneous space \((X,G)\) and a point \(x_0\in X\); let \(H:=G^{x_0}\).
For each \(k=0,1,2\dots\), the subgroup \(H_k\subset G\) above is precisely the set of elements of \(G\) which act on \(X\) at \(x_0\in X\) with the same \(k\)-jet as the identity element, i.e. the identity Taylor series to order \(k\), in any coordinates around \(x_0\).
\end{lemma}
\begin{proof}
For \(k=0\), \(H_0=H\) is precisely the stabilizer of \(x_0\).

For \(k=1\), \(H_1=H_{\LieH_0}\) is the set of elements \(h\in H\) so that \(\LieH\subset\LieG\) is \(h\)-invariant and \(\Ad_h\) is trivial on \(\LieG/\LieH=T_{x_0}X\), i.e. \(h\) acts on \(X\) fixing \(x_0\) and is trivial on \(T_{x_0} X\).

Each \(g\in G\) gives us an element \(x:=g(x_0)\) and a linear map 
\[
T_x X\xrightarrow{(g^{-1})'(x)}T_{x_0}X=\LieG/\LieH.
\]
Changing \(g\) to \(gh\) preserves this point just when \(h\in H\):
\[
ghx_0=gx_0 \text{ iff } hx_0=x_0.
\]
In addition it preserves this linear map just when \(h\in H_1\).
Hence we map
\[
gH_1\in X_1:=G/H_1\mapsto(g^{-1})'(x)\in T^*_x X\otimes(\LieG/\LieH)
\]
by a \(G\)-equivariant injective map.
Since it is an equivariant injection of homogeneous spaces, it is an immersion.

Suppose that the result is true for \(H_0,H_1,\dots,H_{k-1}\).
For each \(p\le k\), associate to each \(g\in G\) the \(p\)-jet \(j^pg(x_0)\) of \(g\) as a map of \(X\), so a map
\[
G\xrightarrow{j^p}X^{(p)},
\]
to the set \(X^{(p)}\) of \(p\)-jets of maps from \(X\) to \(X\) based at \(x_0\).
By hypothesis, for \(p<k\), \(H_p\) is precisely the subgroup leaving this map invariant:
\[
j^p(gh)=j^pg.
\]
Hence the map drops to a \(G\)-equivariant injective immersion
\[
G/H_p\xrightarrow{j^p}X^{(p)}.
\]
By definition, \(H_k\) is the subgroup of \(H_{k-1}\) acting trivially on the tangent space of \(X_{k-1}:=G/H_{k-1}\) at the identity jet.
By thinking of \(H_{k-1}\) as the group \(H\) above, we have \(H_k\) precisely the elements of \(H_{k-1}\) acting trivially on the image of \(G/H_{k-1}\) in the \(1\)-jets of the \((k-1)\)-jets.
But these are precisely the \(k\)-jets, as the \(k\)-jets sit inside the space of \(1\)-jets of the \((k-1)\)-jets \(G\)-invariantly.
\end{proof}
\begin{example}\label{page:proj.jets}
For projective space \((X,G)=(\RP{n},\PGL{n+1})\) (note that we are using a different group \(G\) here, but essentially the same), we saw that \(G/H=X\) is projective space, \(G/H_1\) is the frame bundle of projective space, i.e. \(X^{(1)}\), so any linear isomorphism of tangent spaces of \(X\) arises as the first derivative of a projective transformation.
Then \(H_2\) is trivial, so \(G/H_2=G=\PGL{n+1}\subset X^{(2)}\) is the set of \(2\)-jets of projective transformations: a projective transformation is determined by its \(2\)-jet.
Moreover, the set of \(2\)-jets with given \(1\)-jet is \(H_1/H_2=H_1\), the group of projective transformations of the form
\[
\begin{bmatrix}
1&b\\
0&I
\end{bmatrix},
\]
defined up to multiplying by \((n+1)\)-roots of unity.
So the space of \(2\)-jets of projective transformations with given \(1\)-jet has dimension \(n\).
\end{example}

\chapter{Curvature}
\section{Curvature}
Let us examine flatness more carefully.
Suppose that \(G\to\Bun\to M\) and \(G\to\Bun'\to M'\) are principal right bundles with connections \(\omega,\omega'\), with the same structure group \(G\), and that \(\dim M=\dim M'\).
Pick points \(p\in\Bun\) and \(p'\in\Bun'\), above points \(m\in M\) and \(m'\in M'\).
Any linear isomorphism \(T_m M\xrightarrow{\varphi}T_{m'}M'\) identifies the horizontal spaces of the connections at \(p,p'\), since we saw that these are identified by projection with \(T_m M\) and \(T_{m'} M'\).
Identify the vertical vectors at \(p,p'\) by the obvious linear isomorphism
\[
A_{\Bun}(p)\leftarrow A\in\LieG\to A_{\Bun'}(p').
\]
We have a linear isomorphism \(T_p\Bun\to T_{p'}\Bun'\) identifying \(\omega_p\cong\omega'_{p'}\).
Hence there is no zero order invariant of a connection.
We guess that we should look at a first order invariant, say at \(d\omega\).
Our experience with Lie groups points to a more convenient invariant: since \(\omega\) agrees with the Maurer--Cartan form on the fibers, the \emph{curvature}\define{curvature} \cite{Chern:1989} p. 6
\[
\Omega:=d\omega+\frac{1}{2}\lb{\omega}{\omega}
\]
vanishes on the fibers.
Stronger: the curvature is semibasic.
As a \(2\)-form, the curvature transforms in the adjoint representation:
\[
r_g^*\Omega=\Ad_g^{-1}\Omega.
\]
Hence is a section of the \emph{curvature bundle}
\[
\Lm{2}{T^*M}\otimes(\amal{\Bun}{G}{\LieG})
\]
which is a vector bundle on \(M\), not on \(\Bun\).
If \(u,v\) are vector fields on \(M\), say with horizontal lifts \(\horizontalPart{u},\horizontalPart{v}\), we can also write \(\Omega(u,v)\) to mean \(\Omega(\horizontalPart{u},\horizontalPart{v})\).
\begin{problem}{curvature.values}
Prove that, for any point of the base manifold of a principal bundle, and any element of the curvature bundle at that point, there is a connection on the bundle which acheives that value of curvature at that point.
\end{problem}
\begin{answer}{curvature.values}
Using a partition of unity, it suffices to prove the result locally.
So we can assume the bundle is trivial and that the manifold is a ball in Euclidean space.
Take any element \(F_{ij}dx^i\wedge dx^j\) of the curvature bundle.
Let 
\[
\omega_{\Bun}=\omega_G+\Ad_g^{-1}(F_{ij}x^idx^j).
\]
\end{answer}
\section{Flatness}
\begin{lemma}
A principal bundle is flat just when its curvature vanishes.
Take any principal bundle \(G\to\Bun\to M\) with connection on a connected manifold \(M\).
Let \(\pi:=\pi_1(M)\).
Then the connection is flat just when the bundle is the quotient bundle
\[
\Bun=\amal{\tilde{M}}{\pi}{G}
\]
associated to a group morphism \(\pi\to G\), and the connection is the quotient of the standard flat connection on \(\tilde{M}\times G\).
\end{lemma}
\begin{proof}
Suppose that \(\Bun,\Bun'\to M\) are principal right \(G\)-bundles with connections \(\omega,\omega'\).
Let \(Z:=\Bun\times_M \Bun'\), i.e. the set of pairs \((p,p')\in\Bun\times\Bun'\) projecting to the same point of \(M\).
Let \(V\subset TZ\) be the set of tangent vectors \((v,v')\in TZ\) projecting to the same tangent vector in \(M\) and with
\[
v\hook\omega=v'\hook\omega'.
\]
So \(V\) is the kernel of \(\omega-\omega'\).
Now suppose that both connections \(\omega,\omega'\) are flat.
Then \(V\) is bracket closed precisely because
\[
d(\omega-\omega')=-\frac{1}{2}\lb{\omega}{\omega}+\frac{1}{2}\lb{\omega'}{\omega'}=0\quad\operatorname{mod}\omega-\omega'.
\]
Each map
\[
\begin{tikzcd}
&(p,p')\in Z\arrow[dl]\arrow[dr]&\\
p\in\Bun&&p'\in\Bun'
\end{tikzcd}
\]
has fiber given by a fiber of \(\Bun'\) or \(\Bun\) respectively.
The linear constraints \(\omega=\omega'\) on \(V\) constrain the leaves of the associated folation \(F\), with \(V=TF\), to be complementary to the fibers of both maps.
So the leaves project locally diffeomorphically to both \(\Bun\) and \(\Bun'\). 

Apply the equivariant Frobenius theorem (corollary~\vref{corollary:equivariant.Frobenius}) to \(F\).
The \(G\)-folios project \(G\)-equivariantly and locally diffeomorphically to \(\Bun\) and to \(\Bun'\), and quotient to the identity map on \(M\).

The vector fields \((A_{\Bun},A_{\Bun'})\) are complete.
For any complete vector field \(v\) on \(M\), its horizontal lifts \(\horizontalPart{v},\horizontalPart{v}'\) on \(\Bun,\Bun'\) are complete by theorem~\vref{theorem:horizontal.lift}.
They give a vector field \((\horizontalPart{v},\horizontalPart{v}')\) on \(Z\), also complete.
Together, these \((A_{\Bun},A_{\Bun'})\) and these \((\horizontalPart{v},\horizontalPart{v}')\) span the tangent spaces of \(F\) at every point. 

By Blumenthal's theorem (theorem~\vref{cor:orbitMapEquiv}) each \(G\)-folio \(\Bun''\subseteq Z\) maps to \(\Bun\) and to \(\Bun'\) by fiber bundle maps.
But they have equal dimension, so these are covering maps.
So \(\Bun''\) is a covering space of \(\Bun\) and of \(\Bun'\), and quotients by \(G\)-action to covering spaces
\[
\begin{tikzcd}
&M''\arrow[dl]\arrow[dr]&\\
M&&M'.
\end{tikzcd}
\]
So the pullback bundles are identified on \(M''\), with their connections.
Hence the pullback bundles are identified on the universal covering space \(\tilde{M}\), with their connections.

Now suppose that we pick \(\Bun'\to M\) to be the trivial bundle with the standard flat connection.
So \(\Bun\) pulls back to a trivial bundle on \(\tilde{M}\) with the standard flat connection.
The deck transformations of \(\tilde{M}\to M\) preserve the bundle and the connection, and so act by some bundle automorphism of the standard flat connection on the trivial bundle: \((m,g)\mapsto (\gamma(m),\rho(\gamma)g)\) where \(\rho(\gamma)\in G\).
\end{proof}
\section{Adjoint bundle notation}
Take a right principal bundle \(G\to\Bun\to M\) with connection \(\omega=\omega_\Bun\).
There are two classes of vertical vector fields which we want to consider: (a) those with constant value under the connection and (b) those which which are invariant under the structure group of the bundle.

The first class: every element \(A\in\LieG\) gives rise to a vertical vector field \(A_\Bun\) uniquely determined by the equation \(A=A_\Bun\hook\omega\) or by the equation
\[
A_\Bun(p):=\left.\frac{d}{dt}\right|_{t=0}pe^{tA}.
\]
Such a vector field satisfies
\[
r_{g*}A_\Bun=(\Ad_g^{-1} A)_\Bun.
\]

We can generalize this.
Take any function \(\Bun\xrightarrow{A}\LieG\) and define a unique vertical vector field \(A_\Bun\) by 
\[
A(p)=A_\Bun(p)\hook\omega_p
\]
for any \(p\in\Bun\).
We could also write this as
\[
A_\Bun(p):=\left.\frac{d}{dt}\right|_{t=0}pe^{tA(p)}.
\]
Conversely, every vertical vector field \(v\) arises uniquely as \(v=A_\Bun\) where \(A=v\hook\omega\).
(Some authors use the notation \(\omega^{-1}A\), instead of \(A_\Bun\), for this vector field.)

The second class we want to consider are those vertical vector fields \(A_\Bun\) which are invariant under the structure group:
\[
r_{g*}A_\Bun=A_\Bun,
\]
for all \(g\in G\), or in other words, \(A\) transforms in the adjoint representation:
\[
r_g^*A=\Ad_g^{-1}A.
\]
In other words, these are the sections \(A\) of the bundle
\[
\ad_\Bun:=\amal{\Bun}{G}{\LieG}\to M.
\]
They are identified by the connection with vector fields \(A_\Bun\) which are sections of \(\ad_\Bun\subseteq T\Bun/G\).
Henceforth, we use the notation \(A_\Bun\) only as applied to these two classes of vertical vector fields.
\section{Curvature and Lie bracket}
\begin{lemma}\label{lemma:curv.brack}
Suppose that \(G\to\Bun\to M\) is a principal bundle with connection \(\omega=\omega_\Bun\).
Given a vector field \(u\) on \(M\), let \(\horizontalPart{u}\) be its horizontal lift.
Take two vector fields \(u,v\) on \(M\).
Let \(A:=\Omega(\horizontalPart{u},\horizontalPart{v})\), so \(\Bun\xrightarrow{A}\LieG\).
Then
\[
A_\Bun=\horizontalPart*{\lb{u}{v}}-\lb{\horizontalPart{u}}{\horizontalPart{v}}.
\]
\end{lemma}
\begin{problem}{curvature:hor.bracket}
Prove it.
\end{problem}
\begin{answer}{curvature:hor.bracket}
Since \(\horizontalPart{u},\horizontalPart{v}\) project to \(u,v\), \(\lb{\horizontalPart{u}}{\horizontalPart{v}}\) projects to \(\lb{u}{v}\), as does \(\horizontalPart*{\lb{u}{v}}\), so \(\horizontalPart*{\lb{u}{v}}-\lb{\horizontalPart{u}}{\horizontalPart{v}}\) is vertical.
By the Cartan lemma:
\[
d\omega(\horizontalPart{u},\horizontalPart{v})
=
\LieDer_{\horizontalPart{u}}(\horizontalPart{v}\hook\omega)
-
\LieDer_{\horizontalPart{v}}(\horizontalPart{u}\hook\omega)
-
\lb{\horizontalPart{u}}{\horizontalPart{v}}\hook\omega.
\]
By definition of horizontal, \(\horizontalPart{u}\hook\omega=0\), so
\begin{align*}
A
&=
\Omega(\horizontalPart{u},\horizontalPart{v}),
\\
&=
d\omega(\horizontalPart{u},\horizontalPart{v})+\lb{\horizontalPart{u}\hook\omega}{\horizontalPart{v}\hook\omega},
\\
&=
-\lb{\horizontalPart{u}}{\horizontalPart{v}}\hook\omega,
\\
&=
(\horizontalPart*{\lb{u}{v}}-\lb{\horizontalPart{u}}{\horizontalPart{v}})\hook\omega,
\\
&=
A_\Bun\hook\omega.
\end{align*}
\end{answer}
\begin{lemma}\label{lemma:hor.lift.brackets}
Take a right principal bundle \(G\to\Bun\to M\) with a connection \(\omega\) with curvature \(\Omega\).
Suppose that \(\horizontalPart{u},\horizontalPart{v}\) are the horizontal lifts of vector fields \(u,v\) on \(M\) and 
Then these are also smooth functions on \(\Bun\) transforming in the adjoint representation, hence have associated vertical vector fields \(\Omega(\horizontalPart{u},\horizontalPart{v})_\Bun,dA(u')_\Bun\).
The Lie brackets on \(\Bun\) are
\begin{align*}
\lb{\horizontalPart{u}}{\horizontalPart{v}}&=\horizontalPart*{\lb{u}{v}}-\Omega(\horizontalPart{u},\horizontalPart{v})_\Bun,\\
\lb{\horizontalPart{u}}{A_\Bun}&=dA(\horizontalPart{u})_\Bun,\\
\lb{A_\Bun}{B_\Bun}&=\lb{A}{B}_\Bun.
\end{align*}
\end{lemma}
\begin{problem}{curvature:hor.bracket.2}
Prove it.
\end{problem}
\begin{answer}{curvature:hor.bracket.2}
The first bracket is the one we computed in lemma~\vref{lemma:curv.brack}.
Consider the second.
By definition of horizontal, \(\horizontalPart{u}\hook\omega=0\).
The vector fields \(\horizontalPart{u},A_\Bun\) project to \(u,0\), so their bracket projects to \(0\), i.e. is vertical, say \(\lb{\horizontalPart{u}}{A_\Bun}=B_\Bun\).
By the Cartan lemma:
\begin{align*}
d\omega(\horizontalPart{u},A_\Bun)
&=
\LieDer_{\horizontalPart{u}}(A_\Bun\hook\omega)
-
\LieDer_{A_\Bun}(\horizontalPart{u}\hook\omega)
-
\lb{\horizontalPart{u}}{A_\Bun}\hook\omega,
\\
&=
dA(\horizontalPart{u})
-
\lb{\horizontalPart{u}}{A_\Bun}\hook\omega,
\\
&=
dA(\horizontalPart{u})
-
B.
\end{align*}
Because the curvature vanishes on vertical vector fields, 
\begin{align*}
0
&=
\Omega(\horizontalPart{u},A_\Bun),
\\
&=
d\omega(\horizontalPart{u},A_\Bun)+\lb{\horizontalPart{u}\hook\omega}{A_\Bun\hook\omega},
\\
&=
d\omega(\horizontalPart{u},A_\Bun)+\lb{0}{A_\Bun\hook\omega},
\\
&=
d\omega(\horizontalPart{u},A_\Bun),
\\
&=
dA(\horizontalPart{u})
-
B.
\end{align*}

In a local trivialization, say if \(\Bun=M\times G\), we see that \(A_\Bun,B_\Bun\) are vector fields on \(G\) whose coefficients depend on points of \(M\), which we think of as a parameter.
Brackets are computed for each value of that parameter, and then the parameter plugged in.
Hence the bracket doesn't notice the parameter.
So we can compute on \(G\) directly.
\end{answer}
\begin{example}
Take any Lie algebra \(\LieL\) of vector fields on \(M\), and construct their horizontal lifts.
The brackets above describe a Lie algebra of vector fields on \(\Bun\), generated by the horizontal lifts of \(\LieL\), perhaps infinite dimensional, projecting to \(\LieL\).
\end{example}
\section{Curvature of homogeneous connections}
Recall our notation for homogeneous principal bundles.
Take Lie groups \(G,F,H\) with \(H\subseteq G\) closed.
A \emph{die}\define{die} is a pair of
\begin{itemize}
\item
a Lie group morphism \(h\in H\to\bar{h}\in F\) and
\item
an \(H\)-equivariant linear map \(A\in\LieG\mapsto\bar{A}\in\LieF\) extending the associated Lie algebra morphism.
\end{itemize}
Consider the following recipe: given a die, the associated \(1\)-form is
\[
\gamma:=\omega_F+\Ad_f^{-1}\bar\omega_G
\]
on \(G\times F\).
We conjugate the die by replacing the Lie group morphism by replacing \(\bar{\,}\) by \(\Ad_{f_0} \bar{h}\), for some \(f_0\in F\).
Recall that every homogeneous principal bundle with invariant connection has the form \(\Bun:=\amal{G}{H}{F}\) with a unique connection form pulling back to \(\gamma\), for a die uniquely determined up to conjugation.
Denote by \(\delta\in\LieF\otimes\Lm*{2}{\LieG}^*\) the expression
\[
\delta(A,B):=\lb{\bar{A}}{\bar{B}}-\overline{\lb{A}{B}},
\]
for \(A,B\in\LieG\), i.e. the obstruction to the die being a Lie algebra morphism \(\LieG\to\LieF\).
Note that \(\delta\) vanishes if either \(A\) or \(B\) lie in \(\LieH\), since the die is \(H\)-equivariant, so \(\delta\in\LieF\otimes\Lm*{2}{\LieG/\LieH}^*\).
\begin{theorem}[Wang \cite{Wang1958}]%
\label{thm:hom.bundle.curv}%
Take a homogeneous right principal bundle \(F\to\Bun\to X\).
As above, it can be represented, uniquely up to isomorphism, as \(X=G/H\) and \(\Bun=\amal{G}{H}{F}\) for a Lie group morphism \(H\xrightarrow{a}F\).
Take a \(G\)-invariant connection.
As above, the connection form is represented uniquely as pulling back to \(G\times F\) to become the \(1\)-form
\[
\gamma=\omega_F+\Ad_f^{-1}\bar\omega_G,
\]
for a die \(\LieG\to\LieF\).
The curvature form pulls back to \(G\times F\) to become the \(2\)-form
\[
d\gamma+\frac{1}{2}\lb{\gamma}{\gamma}=\frac{1}{2}\Ad_f^{-1}\delta(\sigma,\sigma),
\]
where \(\sigma=\sigma_G=\omega_G+\LieH\) is the soldering form on \(G\).
In particular, the connection is flat just when the die is a Lie algebra morphism.
\end{theorem}
\begin{problem}{more.connections:compute.hom.curv}
Prove it.
\end{problem}
\begin{answer}{more.connections:compute.hom.curv}
Since the connection form pulls back to \(\gamma\), the curvature pulls back to
\[
\Omega:=d\gamma+\frac{1}{2}\lb{\gamma}{\gamma}.
\]
Take \(A,A'\in\LieG\), \(B,B'\in\LieF\), and associate to these the vector fields
\begin{align*}
v(g,f)&:=(\LT{g*}A,\RT{f*}B),\\
v'(g,f)&:=(\LT{g*}A',\RT{f*}B').
\end{align*}

Denote by \(\vc{A}\) the left invariant vector field associated to a Lie algebra element \(A\in\LieG\), and by \(\rvc{A}\) the right invariant vector field.
Recall (from problem~\vref{problem:right.invariant}) that
\begin{align*}
\lb{\vc{A}}{\vc{A}'}&=\vc*{\lb{A}{A'}},\\
\lb{\rvc{B}}{\rvc{B'}}&=-\rvc*{\lb{B}{B'}},\\
\end{align*}
Our vector fields are
\begin{align*}
v&=\vc{A}+\rvc{B},\\
v'&=\vc{A}'+\rvc{B}',
\end{align*}
so
\[
\gamma(v)=\Ad_f^{-1}(B-\bar A).
\]
The vector fields on \(G\) commute with those on \(F\), so
\[
\lb{v}{v'}=\vc*{\lb{A}{A'}}-\rvc*{\lb{B}{B'}}.
\]
\begin{align*}
\gamma(\lb{v}{v'})
&=
\gamma(\vc*{\lb{A}{A'}})-\gamma(\rvc*{\lb{B}{B'}}),
\\
&=
-\Ad_f^{-1}\overline{\lb{A}{A'}}-\Ad_f^{-1}\lb{B}{B'}.
\end{align*}
We can also compute
\begin{align*}
\frac{1}{2}\lb{\gamma}{\gamma}(v,v')
&=
\lb{\gamma(v)}{\gamma(v')},
\\
&=\lb*{\Ad_f^{-1}(B-\bar A)}{\Ad_f^{-1}(B'-\bar A')},
\\
&=
\Ad_f^{-1}\lb*{B-\bar A}{B'-\bar A'}
\end{align*}
Again
\[
v\hook\gamma=\Ad_f^{-1}(B-\bar A).
\]
Note that the right hand side depends on \(f\in F\), but is constant in \(G\), so 
\[
\LieDer_{v'}(v\hook\gamma)=\LieDer_{\rvc{B}'}(v\hook\gamma).
\]
Recall that the flow of the right invariant vector field \(\rvc{B}'\) through a point \(f\in F\) is the path
\[
e^{t\rvc{B'}}f=e^{tB'}f,
\]
while left invariant vector fields have flows
\[
e^{t\vc{A}}g=ge^{tA}.
\]
(Hence our choice of directions for the arrows on \(\vec{A},\rvc{B}\).)
If we let \(C:=B+\bar A\), then 
\begin{align*}
\LieDer_{v'}(v\hook\gamma)
&=
\LieDer_{\rvc{B}'}(\Ad_f^{-1}C),
\\
&=
\left.
\frac{d}{dt}
\right|_{t=0}
\Ad_{e^{tB'}f}^{-1}C,
\\
&=
\left.
\frac{d}{dt}
\right|_{t=0}
\Ad_f^{-1}\Ad_{e^{tB'}}^{-1}C,
\\
&=
\Ad_f^{-1}\lb{C}{B'},
\\
&=
\Ad_f^{-1}\lb*{B+\bar A}{B'},
\\
&=
\Ad_f^{-1}\lb{B}{B'}+\Ad_f^{-1}\lb*{\bar A}{B'}.
\end{align*}
We can now compute
\begin{align*}
d\gamma(v,v')
&=
\LieDer_v(\gamma(v'))
-
\LieDer_{v'}(\gamma(v))
-\gamma(\lb{v}{v'}),
\\
&=\Ad_f^{-1}\lb{B'}{B}+\Ad_f^{-1}\lb*{\bar A'}{B}\\
&\phantom{=}-\Ad_f^{-1}\lb{B}{B'}+\Ad_f^{-1}\lb*{\bar A}{B'}\\
&\phantom{=}-\Ad_f^{-1}\overline{\lb{A}{A'}}+\Ad_f^{-1}\lb{B}{B'}).
\end{align*}
So
\begin{align*}
\Ad_f d\gamma(v,v')
&=\lb{B'}{B}+\lb*{\bar A'}{B}\\
&\phantom{=}-\lb{B}{B'}-\lb*{\bar A}{B'}\\
&\phantom{=}-\overline{\lb{A}{A'}}+\lb{B}{B'},
\\
&=
-\overline{\lb{A}{A'}}
+\lb*{\bar A'}{B}
-\lb*{\bar A}{B'}
-\lb{B}{B'}
\end{align*}
Finally,
\begin{align*}
\Ad_f\Omega(v,v')
&=
\Ad_f d\gamma
+\frac{1}{2}\Ad_f \lb{\gamma}{\gamma}(v,v'),
\\
&=
-\overline{\lb{A}{A'}}
+\lb*{\bar A'}{B}
-\lb*{\bar A}{B'}
-\lb{B}{B'}
+\lb*{B+\bar A}{B'+\bar A'},
\\
&=
\lb*{\bar A}{\bar A'}-\overline{\lb{A}{A'}}.
\end{align*}
\end{answer}
\begin{example}
As above, for a reductive homogeneous space \((X,G)\), a point \(x_0\in X\) with stabilizer \(H\), and a Lie group morphism \(H\xrightarrow{a}F\), we have a homogeneous principal bundle \(E:=\amal{G}{H}{F}\), and there is a splitting \(\LieG=\LieH\oplus𝔪\) of \(H\)-modules.
To each \(H\)-equivariant linear map \(𝔪\xrightarrow{\beta}\LieF\), as above we generate an invariant connection on \(E\).
We can then compute that the curvature is
\[
\Omega=d\gamma+\frac{1}{2}\lb{\gamma}{\gamma}=\frac{1}{2}\delta_{\beta}(\sigma,\sigma)
\]
\end{example}

\chapter{Curvature in Cartan geometries}
\section{Connections}
Take a homogeneous space \((X,G)\) with \(X=G/H\).
Given any finite dimensional \(H\)-module \(H\xrightarrow{\rho_V}\GL{V}\), denote its differential as \(\LieH\xrightarrow{\rho_V}V^* \otimes V\).
Let \(\vb{V}:=\amal{G}{H}{V}\) be the associated vector bundle \(\vb{V} \to X\).
Let \(A:=\Conn{V}\) be the set of all \(G\)-invariant connections on \(\vb{V}\).
By theorem~\vref{theorem:invariant.connection}, we identify \(A\) with the set of all \(H\)-equivariant linear maps \(\LieG\xrightarrow{\varphi}V^* \otimes V\) for which \(\left.\varphi\right|_{\LieH}=\rho_V\).
If \((X,G)\) is reductive, say as usual \(\LieG=\LieH\oplus\redComplement\), then as usual we can pick \(\varphi=0\) on \(\redComplement\), a canonical connection.

If \(\G\to M\) is an \((X,G\))-geometry with Cartan connection \(\omega\), every element \(\varphi\in A\) gives a \(1\)-form \(\varphi\circ\omega\) which is the pullback of a unique connection \(1\)-form on the associated principal bundle, and hence of a connection on \(\vb{V} \to M\), with curvature
\[
\frac{1}{2}\lb{\phi\omega}{\phi\omega}
-\frac{1}{2}\phi\lb{\omega}{\omega}+
\frac{1}{2}(\phi k)\sigma\wedge\sigma.
\]
If \(\rho_V=0\) on \(\LieH\), we can pick \(\phi:= 0\): a canonical flat connection.
More generally, if we let \(\vb{A}:= \amal{\G}{H}{A}\to M\), then the \(C^{\infty}\) sections of \(\vb{A}\to M\) are precisely the \(C^{\infty}\) connections on \(\vb{V} \to M\).

If \(V\) is a \(G\)-module, the vector bundle \(\vb{V}\to M\) is a \emph{tractor bundle}\define{tractor!bundle}\define{bundle!tractor} and has a canonical connection on \(\vb{V} \to M\), the one induced by the Cartan connection, as a connection on \(\G_G\to M\), called the \emph{tractor connection},\define{tractor!connection}\define{connection!tractor} with \(\varphi\) the differential of the \(G\)-module map \(G\to\GL{V}\).

\section{Curvature}
\begin{example}
The \emph{curvature module}\define{curvature!module} of a homogeneous space \((X,G)\) with \(X=G/H\) is the \(H\)-module \(V:=\LieG \otimes \Lm{2}{(\LieG/\LieH)}^*\).
The \emph{curvature bundle}\define{curvature!bundle}\define{bundle!curvature} of an \((X,G)\)-geometry is the vector bundle \(\vb{V}\) associated to the curvature module \(V\).
\end{example}
The \emph{curvature}\define{curvature!of a Cartan geometry}\define{Cartan!geometry!curvature} of an \((X,G)\)-geometry, with Cartan connection \(\omega\), is the curvature of \(\omega\) as a connection on \(\G_G\), i.e. the section \(k\) of \(\vb{V}\) so that 
\[
d\omega + \frac{1}{2} \lb{\omega}{\omega}=\frac{1}{2}k\sigma\wedge\sigma.
\]
This is the usual definition of curvature of a connection \cite{Chern:1989} p. 6; the usual proof shows that flatness of the connection is equivalent to vanishing of the curvature.

\begin{example}
Any lift or drop has the same curvature.
\end{example}
\begin{example}
Curvature of a pullback is pulled back.
\end{example}
\begin{example}
Any pseudo-Riemannian metric has Cartan connection
\[
\omega
=
\begin{pmatrix}
\gamma&\sigma\\
0&0
\end{pmatrix},
\]
with curvature
\begin{align*}
d\omega
+\frac{1}{2}\lb{\omega}{\omega}
&=
\begin{pmatrix}
d\gamma+\frac{1}{2}\lb{\gamma}{\gamma}&d\sigma+\gamma\wedge\sigma\\
0&0
\end{pmatrix},
\\
&=
\begin{pmatrix}
\frac{1}{2}R\sigma\wedge\sigma&0\\
0&0
\end{pmatrix}
\end{align*}
where \(R\) is the Riemann curvature tensor.
\end{example}
\begin{example}
Penrose's original construction of twistor theory, in the setting of Riemannian \(4\)-manifolds, starts by assuming we have a spin structure, i.e. an \((X',G)\)-geometry, with \(X'=\R^4\) and \(G=\Spin{4}\ltimes\R^4\), so \(H'=\Spin{4}=\SU{2}\times\SU{2}\), acting on the typical tangent space of \(X'\)
\[
T_{x'_0}X'=\LieG/\LieH'=\R^4,
\]
as left and right multiplications by quaternions: \(H'=\SUL\times\SUR\).
This representation of \(H'\) is a real representation, i.e. does not preserve either a complex vector space structure, or a quaternionic vector space structure, on \(\LieG/\LieH'\).
Consider the subgroup \(H=\SUL\times\UnR\subset H'\) with \(\UnR\) the right multiplication by complex numbers, and let \(X=G/H\).
So \(X\to X'\) is a sphere bundle with fibers \(H'/H=\SUR/\UnR=S^2\).
The action of \(H\) on the typical tangent space
\[
T_{x_0} X=\LieG/\LieH=(\LieSUR/\LieUnR)\oplus\R^4=\C\oplus\C^2,
\]
is that of left multiplication by \(\UnR\) on the first factor, and left multiplication by \(\SUL\) and complex scaling by \(\UnR\) on the second.
So this is a complex representation.
Hence if \(M'\) is a Riemannian \(4\)-manifold with a spin structure \(\G\to M'\), then \(M=\G/H\) has an almost complex structure.
Naturally, there is some condition on the curvature of the Cartan connection (which is the \emph{same} Cartan connection for both \(\G\to M\) and for \(\G\to M'\)) under which this is a complex manifold.
So some condition on the curvature of a Riemannian manifold \(M'\) makes the sphere bundle \(M\to M'\) into a complex \(3\)-manifold.
We can make this explicit: for \(\G\to M'\)
\[
d\omega+\frac{1}{2}\lb{\omega}{\omega}=\frac{1}{2}k\sigma\wedge\sigma,
\]
which expand out in terms of \(\sigma\) and \(\gamma\) to become
\begin{align*}
d\sigma&=-{\ell}_{\gamma^{\ell}}\wedge\sigma+\RT{\gamma^r}\wedge\sigma,\\
d\gamma^{\ell}&=-\frac{1}{2}\lb{\gamma^{\ell}}{\gamma^{\ell}}+R^{\ell}\sigma\wedge\sigma,\\
d\gamma^r&=-\frac{1}{2}\lb{\gamma^r}{\gamma^r}+R^r\sigma\wedge\sigma,\\
\end{align*}
in terms of left and right multiplication by quaternions, and various linear combinations \(R^{\ell},R^r\) of the entries of the Riemann curvature tensor.
Think of \(\R^4\) as quaternions, so \(\sigma=\sigma_0+i\sigma_1+j\sigma_2+k\sigma_3\).
The soldering forms of the Cartan geometry \(\G\to M\) are \(\sigma_0+i\sigma_1,\sigma_2-i\sigma_3,\gamma_{r2}-i\gamma_{r3}\); write these as \(\sigma^1,\sigma^2,\sigma^3\).
Compute
\[
d
\begin{pmatrix}
\sigma^1\\
\sigma^2\\
\sigma^3
\end{pmatrix}
=
-
\begin{pmatrix}
i(\gamma^{\ell}_1+\gamma^r_1)&-\gamma^{\ell}_2-i\gamma^{\ell}_3&0\\
\gamma^{\ell}_2-i\gamma^{\ell}_3)&i(\gamma^{\ell}_1-\gamma^r_1&0\\
0&0&2i\gamma^R_1
\end{pmatrix}
\wedge
\begin{pmatrix}
\sigma^1\\
\sigma^2\\
\sigma^3
\end{pmatrix}
+
\begin{pmatrix}
\sigma^3\wedge\bar\sigma^1\\
-\sigma^3\wedge\bar\sigma^2\\
\tau
\end{pmatrix}
\]
where
\[
\tau=a\sigma^1\wedge\sigma^2+b\bar{\sigma}^1\wedge\bar{\sigma}^2+c_{pq}\sigma^p\wedge\bar\sigma^q
\]
for some complex numbers \(a,b,c_{pq}\), \(p,q=1,2\).
By the Newlander--Nirenberg theorem \cite{Demailly2012} p. 399 theorem 11.8, \cite{Kazdan} p. 73, the almost complex structure is complex just when there are no \((0,2)\) terms in the exterior derivatives of the \((1,0)\)-forms, i.e. just when \(b=0\).
A computation with constant coefficient combinations of curvature components shows that \(b\) is a component of the anti-self-dual Weyl curvature.
Vanishing of \(b\), by \(H\)-invariance, is equivalent to vanishing of some irreducible components of the curvature, hence the anti-self-dual Weyl curvature.
Moreover, \(H\subset H'\) is the largest group for which \(\LieG/\LieH\) is a complex representation.
There is a natural generalization of this method to relate Cartan geometries to complex geometry, but never explored, as far as I know.
\end{example}
\begin{example}
Take \(G=\Aff{\R}\), the group of affine transformations of the real number line, and \(X=G\).
We saw~\vpageref{example:aff.R} that no compact surface \(M\) has a flat \((X,G)\)-geometry.
Recall that \(G\) has Maurer--Cartan form
\[
\omega=
\begin{pmatrix}
\alpha&\beta\\
0&0
\end{pmatrix}
\]
so that \(0=d\alpha=d\beta-\alpha\wedge\beta\).
For any \((X,G)\)-Cartan geometry, the curvature consists of \(k,\ell\) in
\begin{align*}
d\alpha&=k \alpha\wedge\beta,\\
d\beta&=(1+\ell)\alpha\wedge\beta.
\end{align*}
There is no \((X,G)\)-Cartan geometry on any connected compact surface \(M\) with curvature \(k\) nowhere zero or with \(\ell\) nowhere equal to \(-1\), since we can replace \(M\) with its oriented covering, to ensure \(M\) is oriented, and then
\begin{align*}
\int_M(1+\ell)\alpha\wedge\beta
&=
\int_M d\beta,
\\&=
\int_{\partial M}\beta,
\\
&=0,
\end{align*}
and similarly with \(d\alpha\).
\end{example}
\begin{lemma}[\cite{Sharpe:1997} p. 187 corollary 3.10]\label{lemma:bracket}
The curvature of any Cartan geometry \(\G\to M\) satisfies the \emph{curvature deformed bracket}\define{curvature!deformed bracket}
\[
\lb{A_{\G}}{B_{\G}}\hook\omega=\lb{A}{B}-k(A,B).
\]
\end{lemma}
\begin{proof}
By the Cartan lemma:
\[
d\omega(A_{\G},B_{\G})
=
\LieDer_{A_{\G}}(B_{\G} \hook \omega)
-
\LieDer_{B_{\G}}(A_{\G} \hook \omega)
-
\lb{A_{\G}}{B_{\G}}\hook \omega.
\]
Since \(A_{\G}\hook\omega=A\) by definition,
\[
d\omega(A_{\G},B_{\G})
=
-
\lb{A_{\G}}{B_{\G}} \hook \omega.
\]
Expand out the definition of curvature:
\begin{align*}
k(A,B)
&=
(A_{\G},B_{\G})\hook (d\omega+ \frac{1}{2}\lb{\omega}{\omega}),
\\
&= -\lb{A_{\G}}{B_{\G}} \hook \omega + \lb{A}{B}.
\end{align*}
\end{proof}
\section{Curvature of homogeneous Cartan geometries}
For a homogeneous Cartan geometry \(\G\to X'\), on a homogeneous space \((X',G')\), modelled on another homogeneous space \((X,G)\), there is a soldering form \(\sigma_\G\) on \(\G:=\amal{G'}{H'}{H}\to X'\) and another \(\sigma_{G'}\) on \(G'\).
It will be convenient to write the die as \(\LieG'\xrightarrow{a}\LieG\), rather than as \(A\mapsto\bar{A}\).
Recall
\[
\delta_a(A,B):=\lb*{a A}{a B}-a\lb{A}{B},
\]
for \(A,B\in\LieG'\), i.e. the obstruction to \(a\) being a Lie algebra morphism \(\LieG'\to\LieG\).
The curvature of the Cartan geometry is
\[
d\omega_\G+\frac{1}{2}\lb{\omega_\G}{\omega_\G}
=
\Ad_h^{-1}\delta_a(\sigma_{G'},\sigma_{G'}),
\]
where \(\LieG'\xrightarrow{a}\LieG\) is the die.
We want to write the curvature in terms of \(\sigma_\G\).
As usual, if \(V\) is an \(H\)-module, denote by
\[
H\xrightarrow{\rho_V}\GL{V}
\]
the associated \(H\)-representation.
Write out
\[
\omega_\G=\omega_H+\Ad_h^{-1}a\omega_{G'},
\]
so quotienting out by \(\LieH\),
\[
\sigma_\G=\rho_{\LieG/\LieH}(h)^{-1}\bar{a}\sigma_{G'},
\]
where \(\LieG'/\LieH'\xrightarrow{\bar{a}}\LieG/\LieH\) is the induced linear isomorphism:
\[
\sigma_{G'}
=\bar{a}^{-1}\rho_{\LieG/\LieH}(h)\sigma_\G,
\]
Let 
\[
k_a(A,B):=\delta_a(\bar{a}^{-1}A,\bar{a}^{-1}B)
\]
for \(A,B\in\LieG/\LieH\).
So \(k_a\in K:=\Lm*{2}{\LieG/\LieH}\otimes\LieG\), the curvature module for \((X,G)\)-geometries.
So the curvature as a Cartan geometry is
\[
d\omega_\G+\frac{1}{2}\lb{\omega_\G}{\omega_\G}
=
(\rho_K(h)^{-1}
k_a)(\sigma_\G,\sigma_\G).
\]

\begin{example}\label{ex:conn.proj.conn.sphere}
On page~\vpageref{section:proj.conn.sphere}, we saw that the projective connections on the sphere which are invariant under rotations are given by the die
\[
\begin{pmatrix}
0 & -c^t\\
c & C
\end{pmatrix}
\in\LieSO{n+1}
\mapsto
\begin{pmatrix}
0&-\beta c^t\\
c&C
\end{pmatrix}
\in\LieSL{n+1}
\]
for an arbitrary real constant \(\beta\) (with \(c\in\R{n}\) and \(C\in\LieSO{n}\), valued in the Lie algebra of the group \(H\) of projective linear transformations fixing the point
\[
[1:0:\dots:0]\in\RP{n}.
\]
Denoting the Cartan connection as \(\omega\), the curvature of the projective connection on \(\G=\amal{\SO{n+1}}{\SO{n}}{H}\) is
\[
d\omega+\frac{1}{2}\lb{\omega}{\omega}=\frac{1}{2}\Ad_h^{-1}\delta(\sigma,\sigma),
\]
where if
\[
\delta
\left(
\begin{pmatrix}
0&-a^t\\
a^t&A
\end{pmatrix},
\begin{pmatrix}
0&-b^t\\
b^t&B
\end{pmatrix}
\right)
=
\beta
\begin{pmatrix}
0&b^tA-a^tB\\
0&ba^t-ab^t
\end{pmatrix}.
\]
In particular, there is a unique flat invariant projective connection: \(\beta=0\).
Since every flat projective connection on a simply connected compact manifold is precisely the pullback \(S^n\to\RP{n}\) by the antipodal covering map, this is the unique flat projective connection on the sphere, up to isomorphism.

Since the resulting curvature is valued inside
\[
\LieH\otimes\Lm*{2}{\LieG/\LieH}
\subseteq
\LieG\otimes\Lm*{2}{\LieG/\LieH}
\]
the projective connection is torsion-free, for any \(\beta\).
\end{example}

\section{Curvature and small parallelograms}\label{subsec:small.parallelograms}
Recall that, for any two vector fields \(X,Y\) on a manifold \(M\),
\[
e^{sX}e^{tY}=e^{sX+tY+st\lb{X}{Y}/2+\dots}
\]
for small \(s,t\), and so
\[
e^{-tY}e^{-sX}e^{tY}e^{sX}=e^{st\lb{X}{Y}+\dots}.
\]
In particular, for any Lie group \(G\) and \(A,B\in\LieG\), since the left invariant vector fields generate the right action, we can rewrite this as
\[
e^Ae^Be^{-A}e^{-B}=e^{\lb{A}{B}+\dots}.
\]
Take the exponent of the right hand side as a definition: let
\[
\lambda_{AB}:=\log\left(e^Ae^Be^{-A}e^{-B}\right)\in\LieG,
\]
for \(A,B\) small enough.
(We could expand this out using the Campbell--Baker--Hausdorff--Dynkin formula.)
By definition, in \(G\),
\[
e^{\lambda_{AB}}e^Ae^Be^{-A}e^{-B}=1,
\]
and
\[
\lambda_{AB}=\lb{A}{B}+\dots
\]
Take a Cartan geometry \(\G\to M\) with model \((X,G)\).
In \(\G\),
\[
\fl{-\lambda_{AB}}\fl{-B}\fl{-A}\fl{B}\fl{A}p=\fl{k_{\G}(p)(A,B)+\dots}
\]
Consider the path in \(\G\) given by following first the flow of \(A_{\G}\) for time \(1\), then \(B_{\G}\), then \(-A_{\G}\), then \(-B_{\G}\), then \(C_{\G}\) where \(C=-\lambda_{AB}\):
\[
\begin{tikzpicture}[decoration={
    markings,
    mark=at position 0.5 with {\arrow{latex}}}
    ] 
\coordinate (p) at (0,0) {};
\coordinate (p') at (0,.3) {};
\coordinate (p'') at (-.3,0) {};
\coordinate (q) at (-.1,1) {};
\coordinate (s) at (1,1.2) {};
\coordinate (r) at (.9,.3) {};
\draw[postaction={decorate}] (p) to[bend left] (r);
\draw[postaction={decorate}] (r) to[bend left,-latex] (s);
\draw[postaction={decorate}] (s) to[bend right,-latex] (q);
\draw[postaction={decorate}] (q) to[bend right] (p');
\draw[postaction={decorate}] (p') to[bend right] (p'');
\fill[black,draw=white] (p) circle (1pt);
\fill[black,draw=white] (q) circle (1pt);
\fill[black,draw=white] (r) circle (1pt);
\fill[black,draw=white] (s) circle (1pt);
\fill[black,draw=white] (p') circle (1pt);
\fill[black,draw=white] (p'') circle (1pt);
%\begin{scope}[scale=.9]
%\coordinate (p) at (0,0) {};
%\coordinate (p') at (0,.3) {};
%\coordinate (p'') at (-.3,0) {};
%\coordinate (q) at (-.1,1) {};
%\coordinate (s) at (1,1.2) {};
%\coordinate (r) at (.9,.3) {};
%\draw[postaction={decorate}] (p) to[bend left] (r);
%\draw[postaction={decorate}] (r) to[bend left,-latex] (s);
%\draw[postaction={decorate}] (s) to[bend right,-latex] (q);
%\draw[postaction={decorate}] (q) to[bend right] (p');
%\draw[postaction={decorate}] (p') to[bend right] (p'');
%\fill[gray,draw=white] (p) circle (1pt);
%\fill[gray,draw=white] (q) circle (1pt);
%\fill[gray,draw=white] (r) circle (1pt);
%\fill[gray,draw=white] (s) circle (1pt);
%\fill[gray,draw=white] (p') circle (1pt);
%\fill[gray,draw=white] (p'') circle (1pt);
%\end{scope}
\end{tikzpicture}
\]
In the model, this closes up, giving a ``pentagonal near parallelogram''.
So the curvature of a Cartan geometry is the second order obstruction to the closing up of this figure in \(\G\) which, by design, closes up in the model.
If the Cartan geometry is flat, then the local isomorphism with the model ensures that it closes up in \(\G\) as well.
\section{Local dropping}
Take a homogeneous space \((X,G)\) and an \((X,G)\)-Cartan geometry \(\G\to M\).
The \emph{local expected structure algebra}\define{local expected structure algebra} \(\LieH_0\subset\LieG\) of \(\G\) is the set of all \(A\in\LieG\) so that \(A\hook k=0\), i.e. \(A\) is null for the curvature, at every point of \(\G\).
Note that \(\LieH\subseteq\LieH_0\).
\begin{lemma}\label{lemma:H.nought.H.invariant}
The linear subspace \(\LieH_0\subseteq\LieG\) is \(H\)-invariant.
\end{lemma}
\begin{proof}
Let \(\Omega:=d\omega+\frac{1}{2}\lb{\omega}{\omega}\).
Note that \(A\hook k=0\) just when \(A_{\G}\hook\Omega=0\).
For any \(h\in H\), \(\RT{h}^*\omega=\Ad_h^{-1}\omega\), so expand out to find that, for any \(h\in H\),
\[
\RT{h}^*\Omega=\Ad_h^{-1}\Omega.
\]
\begin{problem}{H.nought.H.inv}
Apply lemma~\vref{lemma:adj.H} to finish the proof.
\end{problem}
\begin{answer}{H.nought.H.inv}
For any \(A\in\LieG\) and \(h\in H\), the following are equivalent by lemma~\vref{lemma:adj.H}:
\begin{align*}
0&=A\hook k,\\
0&=A_{\G}\hook\Omega,\\
0&=\Ad_h^{-1}(A_{\G}\hook\Omega),\\
0&=A_{\G}\hook\Ad_h^{-1}\Omega,\\
0&=A_{\G}\hook \RT{h}^*\Omega,\\
0&=(\RT{h*}(A_{\G}))\hook\Omega,\\
0&=(\Ad_h^{-1}A)_{\G}\hook\Omega,\\
0&=(\Ad_h^{-1}A)\hook k.
\end{align*}
\end{answer}
\end{proof}
\begin{lemma}\label{lemma:exp.A}
Take a vector \(A\in\LieH_0\).
Let \(h:=\fl{tA}\) and \(\bar{h}:=e^{tA}\).
At all points of \(\G\) where \(h\) is defined,
\[
h^*\omega=\Ad_{\bar{h}}^{-1}\omega, h^*k=\Ad_{\bar{h}}^{-1}k\circ\Lm{2}{\Ad_{\bar{h}}}.
\]
\end{lemma}
\begin{problem}{lemma:exp.A.problem}
Prove it.
\end{problem}
\begin{answer}{lemma:exp.A.problem}
\begin{align*}
\LieDer_{A_{\G}}\omega
&=A_{\G}\hook d\omega+d(A_{\G}\hook\omega),
\\
&=A_{\G}\hook d\omega+d(A_{\G}\hook\omega),
\\
&=A_{\G}\hook(-\frac{1}{2}\lb{\omega}{\omega}+\frac{1}{2}k\sigma\wedge\sigma)+dA,
\\
&=-\lb{A}{\omega}+k(A,\sigma),
\\
&=-\lb{A}{\omega}.
\end{align*}
So if
\[
\omega_t:=\Ad_{e^{tA}}\fl{tA*}\omega
\]
\begin{align*}
\frac{d}{dt}\omega_t&=
\lb{A}{\omega_t}-\lb{A}{\omega_t},
\\
&=0.
\end{align*}
Integrate to find the law for transformation of \(\omega\) under the flow; then differentiate to find the law for \(k\).
\end{answer}
\begin{lemma}
For any \(A,B\in\LieH_0\), if \(h:=\fl{tA}\), then \(\RT{h*}(B_{\G})=(\Ad_h^{-1}B)_{\G}\).
\end{lemma}
\begin{proof}
The same proof as for lemma~\vref{lemma:adj.H}.
\end{proof}
\begin{lemma}\label{lemma:H.nought.H.nought.invariant}
The linear subspace \(\LieH_0\subseteq\LieG\) is \(e^{\LieH_0}\)-invariant.
\end{lemma}
\begin{proof}
The same proof as for lemma~\vref{lemma:H.nought.H.invariant}.
\end{proof}
\begin{lemma}
The local expected structure algebra \(\LieH_0\subset\LieG\) is a Lie subalgebra under the Lie bracket of \(\LieG\).
\end{lemma}
\begin{proof}
The definition of \(\LieH_0\) consists of linear equations on elements of \(\LieG\), so \(\LieH_0\subseteq\LieG\) is a linear subspace.
Differentiate lemma~\vref{lemma:H.nought.H.nought.invariant}: it is closed under bracket.
\end{proof}
The \emph{local expected structure group}\define{local expected structure group} \(H_0\subseteq G\) is the smallest immersed Lie subgroup containing \(H\) and invariant under the flows of the local expected structure algebra, i.e. whose Lie algebra contains \(\LieH_0\).
Since \(\LieH_0\) is \(H\)-invariant, \(H_0\) has Lie algebra \(\LieH_0\).
\begin{proposition}
Suppose that \(H'\subseteq H_0\) is a subgroup of the local expected structure group, and that \(H'\subseteq G\) is closed.
Let \(X':=G/H'\).
Then \(\G\) is covered by open sets \(U_a\), each of which embeds as an open set \(U_a\subset\G_a\) into the total space of an \((X',G)\)-Cartan geometry \(\G_a\to M_a\), with Cartan connection \(\omega_a\) pulling back to \(U_a\) to equal \(\omega\), the Cartan connection of \(\G\to M\).
\end{proposition}
This follows from lemma~\vref{lemma:local.dropping}.
An extreme case: if \(H=\set{1}\), we have only a manifold \(\G\) with a \(1\)-form \(\omega\) valued in \(\LieG\), and we ask when it arises, at least locally, from a Cartan geometry with model \((X,G)\), \(X=G/H\).
By this proprosition, that occurs just when the local expected structure algebra contains the Lie algebra of \(H\).
\section{Global dropping}
The \emph{expected structure algebra}\define{expected structure algebra} is the subset \(\LieH'\subseteq\LieH_0\) of complete vector fields in the local expected structure algebra.
(We will define the \emph{complete cone}\define{complete!cone} to be the elements of \(\LieG\) whose associated constant vector fields are complete; the expected structure algebra is the intersection of the complete cone and the local expected structure algebra.)
Note that \(\LieH\subseteq\LieH'\).
It is clear that \(\LieH'\) is a cone, i.e. closed under scaling and containing the origin, but a priori it is perhaps not clear that \(\LieH'\) is a Lie subalgebra, so we will prove this.
\begin{lemma}
The expected structure algebra \(\LieH'\subset\LieG\) is closed under the Lie bracket of \(\LieG\).
\end{lemma}
\begin{proof}
By lemma~\vref{lemma:exp.A}, if \(A\in\LieH'\) and \(B\in\LieG\) then
\[
\fl{sA}_*B_{\G}=(\Ad_{e^{-sA}}B)_{\G}.
\]
Any diffeomorphism preserves and reflects completeness.
For \(A,B\in\LieH'\), differentiate through this family of complete vector fields, to get a complete vector field \(\lb{A}{B}_{\G}\).
\end{proof}
Consider pairs \((h,\bar{h})\) for which \(\G\xrightarrow{h}\G\) is a diffeomorphism and \(\bar{h}\in G\) with
\[
h^{-1*}\omega=\Ad_{\bar{h}}^{-1}\omega.
\]
Get such pairs to act on \(\G\times G\) on the right by
\[
(p,g)(h,\bar{h}):=(h^{-1}(p),g\bar{h}).
\]
The \emph{expected structure group}\define{expected structure group} \(H'\) is the subgroup of these pairs generated by 
\begin{itemize}
\item
pairs \((h,\bar{h})=(\RT{\bar{h}}^{-1},\bar{h})\) for \(\bar{h}\in H\) and
\item
pairs \((\fl{-tA},e^{tA})\) for \(A\in\LieH'\)
\end{itemize}
Map \((h,\bar{h})\in H'\mapsto \bar{h}\in G\).
The \emph{expected model}\define{expected model}\define{model!expected} is \((X',G)\) with \(X':=G/H'\).
Our aim is to prove:
\begin{theorem}
Take a homogeneous space \((X,G)\) and an \((X,G)\)-Cartan geometry \(\G\to M\).
The Cartan geometry drops to its expected model just when the expected structure group \(H'\) is injectively mapped \(H'\to G\) with closed image, and \(H'\) acts properly on \(\G\).
\end{theorem}
\begin{lemma}
The expected structure algebra is a Lie subalgebra.
\end{lemma}
\begin{proof}
We only need to prove that it is closed under addition.
On the manifold \(\G\times G\), consider the vector fields \(A_{\G}+A_G\), for \(A\in\LieH_0\).
Since curvature vanishes on \(\LieH_0\), the equivariant Frobenius theorem (corollary~\vref{corollary:equivariant.Frobenius}) shows that \(\G\times G\) is foliated by \(H\)-folios of these vector fields.
Each \(H\)-folio projects by local diffeomorphism to \(H_0\subset G\) and to a leaf of the \(\LieH_0\) constant vector fields in \(\G\).
By the equivariant orbit theorem (theorem~\vref{theorem:orbit.maps}), the same is true for the vector fields \(A_{\G}+A_G\) for \(A\in\LieH'\).
But these vector fields are complete, so the \(H\)-folio projects by a fiber bundle map, again by theorem~\vref{theorem:orbit.maps}.
Since the projection is a local diffeomorphism, and a fiber bundle map, it is a covering map.
Covering an immersed Lie subgroup of \(G\), matching up constant vector fields, all constant vector fields are complete.
\end{proof}
\begin{theorem}%
[Palais \cite{Palais:1957} p. 58 theorem XI]%
\label{theorem:Palais}%
\define{Palais theorem}
\define{theorem!Palais}
If a collection of smooth complete vector fields on a connected manifold generates a finite dimensional Lie algebra of vector fields, then there is a connected Lie
group acting smoothly and faithfully inducing that Lie algebra action, unique up to unique isomorphism matching up the action.
\end{theorem}
Hence the expected natural structure algebra generates flows forming a Lie group acting on \(\G\times G\) smoothly and faithfully, and sitting inside the expected structure group.
\begin{theorem}
Take a homogeneous space \((X,G)\) and an \((X,G)\)-Cartan geometry \(\G\to M\).
The expected structure group is a Lie group, acting faithfully on \(\G\times G\), with Lie algebra the expected structure algebra, splitting the Lie group inclusion
\[
H\to G
\]
into Lie group morphisms
\[
H\to H'\to G
\]
with associated Lie algebra morphisms the inclusions
\[
\LieH\to\LieH'\to\LieG
\]
and with every component of \(H'\) containing elements of \(H\).
If \(\G\) is connected then \(H'\) acts freely on \(\G\times G\).
If \(H'\to G\) has finite kernel and closed image, then \(H'\) acts properly on \(\G\times G\).
\end{theorem}
\begin{proof}
Since the expected natural structure algebra \(\LieH'\) is \(H\)-invariant, so is the group \(H_0'\) generated by the flows of \(\LieH'\).
So \(H_0'\subseteq H'\) is a normal subgroup.
This group \(H_0'\) is a Lie group acting smoothly by Palais's theorem.
Take \(H'\)-translates of the open sets of \(H_0'\) as open sets of \(H'\).
Because \(H'_0\subseteq H'\) is normal, acted on smoothly by \(H\), this is well defined, making \(H'\) a Lie group with \(H_0'\) an open subgroup.
Writing each element of \(H'\) as a product of elements of \(H_0'\) and elements of \(H\), \(H'\) acts by diffeomorphisms on \(\G\), and these depend smoothly on the Lie group structure, so \(H'\) is Lie group acting smoothly on \(\G\).
By definition, \(H'\) is a group of diffeomorphisms of \(\G\times G\), so it acts faithfully on \(\G\times G\).
The map \(H'\to G\) is a Lie group morphism.

Any element \((h,\bar{h})\in H'\) fixing a point \((p,g)\) has \(\bar{h}=1\) so  \(h^{-1*}\omega=\omega\) with \(h(p)=p\), so \(h=1\) on the component of \(p\) in \(\G\).
So if \(\G\) is connected then \(h\) is an automorphism fixing a point, so the identity.
So \(H'\) acts freely.

Suppose that \(H'\to G\) has finite kernel and closed image.
The image acts freely and properly on \(G\), being a closed Lie subgroup, hence \(H'\) acts properly on \(\G\times G\).
\end{proof}
Suppose that \(H'\to G\) is injective with closed image.
The image is therefore a closed subgroup isomorphic to \(H'\) and \(H'\to G\) is a closed embedding.
\begin{example}
Take a parabolic subgroup \(H\subset G\) where \(G\) is the split real form of a reductive complex linear algebraic group.
Then every Lie subgroup of \(G\) containing \(H\) is also a parabolic subgroup, therefore a closed connected subgroup.
So if we take the expected structure group \(H\to H'\to G\) of some parabolic geometry, the image of \(H'\to G\) is a closed parabolic subgroup and \(H'\) is a covering of that group.
\end{example}

\section{Torsion}
The \emph{soldering form}\define{soldering form} \(\sigma\) of a Cartan geometry is composition of the Cartan connection \(\omega\) (valued in \(\LieG\)) with the linear projection \(\LieG\to\LieG/\LieH\); we can write this as \(\sigma=\omega+\LieH\).
The \emph{torsion}\define{torsion} of a Cartan geometry is the projection of the curvature from the curvature \(H\)-module \(\LieG \otimes \Lm{2}{(\LieG/\LieH)}^*\) to \((\LieG/\LieH) \otimes \Lm{2}{(\LieG/\LieH)}^*\).
\subsection{Torsion in reductive geometries}
A Cartan geometry is a \emph{reductive geometry}\define{reductive!geometry}\define{Cartan!geometry!reductive} if its model is reductive.
Writing the model Lie algebra as \(\LieG=\LieH\oplus\redComplement\), we can identify \(\LieG/\LieH=\redComplement\) so identify the soldering form \(\sigma\) with the \(\redComplement\)-part of \(\omega\).
The \(\LieH\)-part of \(\omega\) is then a connection \(\gamma\) on \(\G=\G_H\), and determines an affine connection on \(M\).
Hence every reductive geometry imposes an affine connection.
The torsion of the reductive geometry is then identified with \(d\sigma+\gamma\wedge\sigma\), i.e. with the torsion of the affine connection.
\begin{problem}{reductive.spiral}
Prove that the flow line of any constant vector field in any reductive geometry becomes horizontal after suitable right action by the structure group, in which case the value of the Cartan connection changes from a constant to an exponential horizontal curve in the model. 
Hint: use lemma~\vref{lemma:reparam}.
\end{problem}
\subsection{Torsion and small parallelograms}
Recall our \(5\)-sided near parallelogram from section~\vref{subsec:small.parallelograms}:
\[
\begin{tikzpicture}[decoration={
    markings,
    mark=at position 0.5 with {\arrow{latex}}}
    ] 
\coordinate (p) at (0,0) {};
\coordinate (p') at (0,.3) {};
\coordinate (p'') at (-.3,0) {};
\coordinate (q) at (-.1,1) {};
\coordinate (s) at (1,1.2) {};
\coordinate (r) at (.9,.3) {};
\draw[postaction={decorate}] (p) to[bend left] (r);
\draw[postaction={decorate}] (r) to[bend left,-latex] (s);
\draw[postaction={decorate}] (s) to[bend right,-latex] (q);
\draw[postaction={decorate}] (q) to[bend right] (p');
\draw[postaction={decorate}] (p') to[bend right] (p'');
\fill[black,draw=white] (p) circle (1pt);
\fill[black,draw=white] (q) circle (1pt);
\fill[black,draw=white] (r) circle (1pt);
\fill[black,draw=white] (s) circle (1pt);
\fill[black,draw=white] (p') circle (1pt);
\fill[black,draw=white] (p'') circle (1pt);
%\begin{scope}[scale=.9]
%\coordinate (p) at (0,0) {};
%\coordinate (p') at (0,.3) {};
%\coordinate (p'') at (-.3,0) {};
%\coordinate (q) at (-.1,1) {};
%\coordinate (s) at (1,1.2) {};
%\coordinate (r) at (.9,.3) {};
%\draw[postaction={decorate}] (p) to[bend left] (r);
%\draw[postaction={decorate}] (r) to[bend left,-latex] (s);
%\draw[postaction={decorate}] (s) to[bend right,-latex] (q);
%\draw[postaction={decorate}] (q) to[bend right] (p');
%\draw[postaction={decorate}] (p') to[bend right] (p'');
%\fill[gray,draw=white] (p) circle (1pt);
%\fill[gray,draw=white] (q) circle (1pt);
%\fill[gray,draw=white] (r) circle (1pt);
%\fill[gray,draw=white] (s) circle (1pt);
%\fill[gray,draw=white] (p') circle (1pt);
%\fill[gray,draw=white] (p'') circle (1pt);
%\end{scope}
\end{tikzpicture}
\]
associated to any \(A,B\in\LieG\) small enough.
Imagine scaling the figure, i.e. replacing \(A,B\) by \(tA,tB\).
\[
\begin{tikzpicture}[decoration={
    markings,
    mark=at position 0.5 with {\arrow{latex}}}
    ] 
\begin{scope}[scale=.60,draw=gray!25]
\coordinate (p) at (0,0) {};
\coordinate (p') at (0,.3) {};
\coordinate (p'') at (-.3,0) {};
\coordinate (q) at (-.1,1) {};
\coordinate (s) at (1,1.2) {};
\coordinate (r) at (.9,.3) {};
\draw[postaction={decorate}] (p) to[bend left] (r);
\draw[postaction={decorate}] (r) to[bend left,-latex] (s);
\draw[postaction={decorate}] (s) to[bend right,-latex] (q);
\draw[postaction={decorate}] (q) to[bend right] (p');
\draw[postaction={decorate}] (p') to[bend right] (p'');
\fill[gray!50,draw=white] (p) circle (1pt);
\fill[gray!50,draw=white] (q) circle (1pt);
\fill[gray!50,draw=white] (r) circle (1pt);
\fill[gray!50,draw=white] (s) circle (1pt);
\fill[gray!50,draw=white] (p') circle (1pt);
\fill[gray!50,draw=white] (p'') circle (1pt);
\end{scope}
\begin{scope}[scale=.70,draw=gray!50]
\coordinate (p) at (0,0) {};
\coordinate (p') at (0,.3) {};
\coordinate (p'') at (-.3,0) {};
\coordinate (q) at (-.1,1) {};
\coordinate (s) at (1,1.2) {};
\coordinate (r) at (.9,.3) {};
\draw[postaction={decorate}] (p) to[bend left] (r);
\draw[postaction={decorate}] (r) to[bend left,-latex] (s);
\draw[postaction={decorate}] (s) to[bend right,-latex] (q);
\draw[postaction={decorate}] (q) to[bend right] (p');
\draw[postaction={decorate}] (p') to[bend right] (p'');
\fill[gray!50,draw=white] (p) circle (1pt);
\fill[gray!50,draw=white] (q) circle (1pt);
\fill[gray!50,draw=white] (r) circle (1pt);
\fill[gray!50,draw=white] (s) circle (1pt);
\fill[gray!50,draw=white] (p') circle (1pt);
\fill[gray!50,draw=white] (p'') circle (1pt);
\end{scope}
\begin{scope}[scale=.85,draw=gray]
\coordinate (p) at (0,0) {};
\coordinate (p') at (0,.3) {};
\coordinate (p'') at (-.3,0) {};
\coordinate (q) at (-.1,1) {};
\coordinate (s) at (1,1.2) {};
\coordinate (r) at (.9,.3) {};
\draw[postaction={decorate}] (p) to[bend left] (r);
\draw[postaction={decorate}] (r) to[bend left,-latex] (s);
\draw[postaction={decorate}] (s) to[bend right,-latex] (q);
\draw[postaction={decorate}] (q) to[bend right] (p');
\draw[postaction={decorate}] (p') to[bend right] (p'');
\fill[gray,draw=white] (p) circle (1pt);
\fill[gray,draw=white] (q) circle (1pt);
\fill[gray,draw=white] (r) circle (1pt);
\fill[gray,draw=white] (s) circle (1pt);
\fill[gray,draw=white] (p') circle (1pt);
\fill[gray,draw=white] (p'') circle (1pt);
\end{scope}
\coordinate (p) at (0,0) {};
\coordinate (p') at (0,.3) {};
\coordinate (p'') at (-.3,0) {};
\coordinate (q) at (-.1,1) {};
\coordinate (s) at (1,1.2) {};
\coordinate (r) at (.9,.3) {};
\draw[postaction={decorate}] (p) to[bend left] (r);
\draw[postaction={decorate}] (r) to[bend left,-latex] (s);
\draw[postaction={decorate}] (s) to[bend right,-latex] (q);
\draw[postaction={decorate}] (q) to[bend right] (p');
\draw[postaction={decorate}] (p') to[bend right] (p'');
\fill[black,draw=white] (p) circle (1pt);
\fill[black,draw=white] (q) circle (1pt);
\fill[black,draw=white] (r) circle (1pt);
\fill[black,draw=white] (s) circle (1pt);
\fill[black,draw=white] (p') circle (1pt);
\fill[black,draw=white] (p'') circle (1pt);

\end{tikzpicture}
\]
Suppose it starts at some point \(p_0\in\G\), not changing with \(t\), and ends at a point \(p(t)\).
As we scale, \(p(0)=p_0\) and \(p'(0)=0\).
So the end point \(p(t)\) moves, near \(t=0\), with zero velocity and with acceleration
\[
p''(0)\hook\omega=2k(p)(A,B).
\]
Consider the image \(m(t)=\pi(p(t))\in M\) of \(p(t)\in\G\).
Recall from proposition~\vref{prop:TM} that \(TM\) is identified with the vector bundle associated to \(\LieG/\LieH\).
So \(m'(0)=0\) while the acceleration \(m''(0)\) is represented by \(2k(p)(A,B)+\LieH\), the torsion.
So our near parallelogram, closed in the model, has image in \(M\) closing up, up to second order, and up to third order just at points where the torsion vanishes.
\begin{example}
Suppose that \(\lb{A}{B}=0\) in the model.
Our near parallelogram is a parallelogram in \(G\).
The abelian Lie algebra generated by \(A,B\) maps by exponential map to an abelian Lie subgroup, a quotient by a discrete abelian subgroup.
All parallelograms in the Lie algebra exponentiate to close up in the Lie group.
But then any curvature of the Cartan geometry still prevents them, for small \(A,B\), closing up in \(\G\).
\end{example}
\begin{example}
For any reductive geometry, i.e. modelled on a reductive homogeneous space \((X,G)\) with \(\LieG=\LieH\oplus\redComplement\) as \(H\)-modules, we can pick to use only elements of \(\redComplement\) to form our near parallelograms.
We can replace the model with a new model with \(X=\redComplement\) and with \(G=H\ltimes\redComplement\), where \(\redComplement\) is now an abelian Lie subalgebra, the \emph{affine form}\define{affine form} of the origin model.
In the model, the parallelogram is the horizontal lift (via the affine connection, which is flat) of a parallelogram in \(X\).
In any \((X,G)\)-Cartan geometry \(\G\to M\), the four sides of the near parallelogram project to geodesics in \(M\).
\end{example}
This recovers the well known pictures of torsion \cite{Misner/Thorne/Wheeler:1973} p. 250 box 10.2C, \cite{Penrose:2005} p. 306.

\section{Invariant metrics}
If \(H\) is compact, every \((X,G)\)-geometry with \(X=G/H\) has an invariant Riemannian metric from each \(G\)-invariant metric on \(X\).

If \(H\) is not compact?
The group of components of a Lie group \(H\) has a maximal finite subgroup, unique up to conjugacy, just when \(H\) contains a maximal compact subgroup \(H_c\subseteq H\), unique up to conjugacy \cite{Hilgert.Neeb:2012} p. 531 Theorem 14.1.3.

If the group of components of \(H\) does not have a maximal finite subgroup? 
Take \(H_c\) to be any compact subgroup \(H_c\subseteq H\) so that the identity component of \(H_c\) is a maximal compact subgroup of the identity component of \(H\).

Let \(X_c:=G/H_c\), a homogeneous Riemannian manifold.
Every \((X,G)\)-geometry \(\G\to M\) lifts to a \((X_c,G)\)-geometry \(\G\to M_c=\G/H_c\), a fiber bundle \(H/H_c\to M_c\to M\), the \emph{scaffold}.\define{scaffold}

Since \(H_c\) is compact, it has a biinvariant smooth measure \(dV\), scaled by a positive constant to have unit total volume \cite{Hilgert.Neeb:2012} p. 407 section 10.4, from which we can average
\[
\left<x,y\right>^{H_c}:=\int_{h \in H_c} \left<hx,hy\right>dV
\]
to take any positive definite inner product \(\left<x,y\right>\) on any \(H_c\)-module, and define an \(H_c\)-invariant positive definite inner product.
Hence take an \(H_c\)-invariant positive definite inner product \(\left<,\right>\) on \(\LieG\).
The Riemannian metric \(\left<\omega,\omega\right>\) on \(\G\) drops to a Riemannian metric on \(M_c\).
The automorphisms of any \((X,G)\)-geometry act as isometries of the scaffold.
The scaffold \(H/H_c\to M_c\to M\) has fibers diffeomorphic to Euclidean space, with a complete Riemannian metric of nonpositive curvature on each fiber.
\begin{example}
In general relativity, \((X,G)=(\R^{n,1},\Orth{n,1}\ltimes\R^{n,1})\), \(H=\Orth{n,1}\), \(H_c=\Orth{n}\times\Orth{1}\) is the stabilizer of a time line and a space hyperplane, so the scaffold of a relativistic space-time \(M=M^{n,1}\) is the bundle \(M_c\to M\) of all choices of splitting of tangent spaces of \(M\) into time line and space hyperplane.
\end{example}
It seems reasonable to say that the Hodge theory of a Cartan geometry is the Hodge theory of its scaffold.
\begin{lemma}\label{lemma:completeness.shared}
For any Cartan geometry \(\G\to M\), in the Riemannian metrics defined above, the scaffold \(M_c\) is a complete Riemannian manifold just when \(\G\) is a complete Riemannian manifold. 
\end{lemma}
\begin{proof}
The bundle map \(H_c\to\G\to M_c\) is proper and nondecreasing in the Riemannian metrics.
On the perpendicular space to the fibers of \(\G\to M_c\) (the ``horizontal space''), the Riemannian metric is an isometry.
Suppose that \(M_c\) is complete.
Given a Cauchy sequence \(p_1,p_2,\dots\in\G\), its image \(m_1,m_2,\dots\) in \(M_c\) is a Cauchy sequence, so convergent, say \(m_i\to m\).
Pick a compact ball \(\bar{B}\) of positive radius around \(m\).
All but a finite number of these \(m_i\) stay inside \(\bar{B}\).
The preimage of \(\bar{B}\) in \(\G\) is compact, so some subsequence of these \(p_i\) converges, hence the original sequence converges.

Suppose that \(\G\) is complete.
Take a Cauchy sequence \(m_i\in M_c\).
Pick a path from each \(m_i\) to \(m_{i+1}\), of lengths approaching the infimum length.
Lift the paths to horizontal paths in \(\G\), say with ends \(p_i\), so \(p_i\) maps to \(m_i\).
So these \(p_i\) form a Cauchy sequence, so converge \(p_i\to p\), so \(p\) maps to a limit for the \(m_i\).
\end{proof}

\section{Affine connections}
More generally, consider any reductive homogeneous space \((X,G)\), with stabilizer \(H=G^{x_0}\) for some \(x_0\in X\), so \(\LieG=\LieH\oplus\redComplement\) as an \(H\)-module.
A \emph{soldering form}\define{soldering form} \(\sigma\) on an \(H\)-bundle \(\G\) is an \(H\)-equivariant \(1\)-form \(\sigma\in\Omega^1_{\G}\otimes\redComplement\) vanishing precisely on the vertical vectors, i.e. the tangent spaces of the fibers of \(\G\to M\).
A \emph{frame bundle}\define{frame bundle} for a manifold \(M\) is a principal \(H\)-bundle \(\G\to M\) with a soldering form.

Because the Lie algebra splits \(\LieG=\LieH\oplus\redComplement\) as an \(H\)-module, every \((X,G)\)-geometry \(\G\to M\) has Cartan connection splitting accordingly:
\[
\omega=(\gamma,\sigma),
\]
where \(\sigma\) is a soldering form, and \(\gamma\) is a connection.
Conversely, any \(H\)-bundle with soldering form and connection becomes an \((X,G)\)-geometry in precisely this way.

Affine space is \((X,G):=(V,\GL{\redComplement}\ltimes\redComplement)\) with stabilizer \(H:=\GL{\redComplement}\).
With affine space as model, the \emph{standard frame bundle}\define{standard frame bundle}\define{frame bundle!standard} of \(M\) is the set of pairs \((m,u)\) of \(m\in M\) and linear isomorphism \(T_m M\xrightarrow{u}\redComplement\), with right \(H\)-action
\[
(m,u)h:=(m,h^{-1}u)
\]
and bundle map 
\[
(m,u)\xmapsto{\pi}m.
\]
(Clearly this generalizes the definition from page~\vpageref{section:frame.bundle} which we called, and often continue to call, the \emph{frame bundle}.)
On the standard frame bundle, we define a soldering form \(\sigma\) by
\[
\sigma_{(m,u)}v=u(\pi'(m,u)v)
\]
for \(m\in M\), \((m,u)\in\G\), and \(v\in T_{(m,u)}\G\).
As we will see in section~\vref{section:soldering}, every frame bundle modelled on affine space is canonically isomorphic to the standard frame bundle.

As we have seen \vpageref{section:vb.connections}, every affine connection, i.e. vector bundle connection on the tangent bundle, determines a connection on the standard frame bundle and vice versa.
Hence Cartan geometries modelled on affine space are identified with affine connections
The torsion of the Cartan geometry is precisely the torsion of the affine connection.
\begin{example}
By theorem~\vref{theorem:invariant.connection}, any homogeneous space \(X=G/H\) admits a \(G\)-invariant affine connection, i.e. connection on its tangent bundle, hence on its frame bundle, just when the adjoint representation of \(G\), restricted to \(H\) and then differentiated at \(1\in H\):
\[
\LieH\xrightarrow{\alpha}\LieGL{\LieG/\LieH}
\]
extends to an \(H\)-module morphism
\[
\LieG\xrightarrow{\beta}\LieGL{\LieG/\LieH}.
\]
Each extension yields a different invariant affine connection.
\end{example}

\section{Projective connections}
A \emph{projective connection}\define{projective!connection}\define{connection!projective} is a Cartan geometry modelled on \((X,G):=(\Proj{n},\PGL{n+1})\).
Let us consider how we might run into such an object in the woods.
How many affine connections have the same geodesics? 
Two affine connections, on a manifold of dimension \(2\) or more, are \emph{projectively equivalent}\define{affine connection!projectively equivalent}\define{projectively!equivalent!affine connections} if they have the same geodesics (up to reparameterization) and the same torsion.
A \emph{projective structure}\define{projective!structure} on a manifold \(M\) is a maximal covering of \(M\) by open sets, with an affine connection on each open set, any two projectively equivalent.
\begin{example}
The standard Euclidean metric on the unit ball has straight lines as its geodesics, but so does the hyperbolic metric, in the Betrami--Klein model, so they have the same projective structure.
\end{example}
A projective connection is \emph{normal}\define{normal!projective connection}\define{projective!connection!normal} if it satisfies some linear equations on its curvature, which we won't make precise.
We also won't define its geodesics, but the reader can guess what they  should be, in terms of the Cartan connection, by analogy with the model \(X=\Proj{n}\).
Nor will we prove:
\begin{theorem}[Cartan \cite{Cartan:1924}, Kobayashi \cite{Kobayashi:1995}, Kobayashi \& Nagano \cite{KobayashiNagano:1964}]
To any torsion-free projective structure there is associated a unique normal projective connection with the same geodesics, as unparameterized curves.
Conversely, every normal projective connection is induced by a unique projective structure.
\end{theorem}
The interested reader can pursue \cite{Molzon/Mortensen:1996}, \cite{Spivak:1979b} p. 273.
Using partitions of unity, any projective structure arises from an affine connection.
It is not known whether every real analytic projective structure is induced by a real analytic torsion-free affine connection.

For example, take a Riemannian metric on a surface.
The Cartan geometry of the Riemannian metric is modelled on \((X,G)=(\R^2,\SO{2}\ltimes\R^2)\).
It is convenient to identify \(X=\R^2=\C\) and \(\SO{2}=\Un{1}\).
The Cartan connection is thus valued in \(\LieUn{1}\ltimes\C\), so we can write it split into 
\[
\sigma\in\C,\gamma\in\R
\]
so that \(i\gamma\in\LieUn{1}\) and 
\[
d\sigma=i\gamma\wedge\sigma
\]
and the curvature is
\[
d\gamma=G\sigma^1\wedge\sigma^2
\]
This function \(G\) is precisely the Gauss curvature.
Differentiate the Gauss curvature by writing, on the frame bundle,
\[
dG=G_1\sigma^1+G_2\sigma^2
\]
for unique functions \(G_1,G_2\).

We won't prove: the bundle of orthonormal frames is a subbundle of the bundle of the projective connection.
Nor will we prove: on that subbundle, the projective connection curvature is the tensor
\[
(G_2\sigma^1-G_1\sigma^2)\otimes(\sigma^1\wedge\sigma^2)
\]
defined on the surface, a section of \(T^*M\otimes\Lm{2}{T^*M}\), vanishing exactly at the critical points of Gauss curvature.
\begin{theorem}
The normal projective connection of a Riemannian metric on a surface is flat just when the Gauss curvature is constant on every connected component of the surface.
\end{theorem}
A \emph{projective automorphism}\define{projective!automorphism} of a Riemannian manifold is an automorphism of its associated normal projective connection, i.e. of its projective structure, i.e. of its geodesics, as a set of unparameterized curves.
It is \emph{projectively homogeneous}\define{projectively!homogeneous} if these automorphisms act transitively on the Riemannian manifold.
\begin{theorem}
The projectively homogeneous compact connected surfaces with Riemannian metric are the sphere, the real projective plane, the torus, and the Klein bottle, the first two with their standard round metrics up to constant rescaling, and the latter two with any flat metric.
Any projectively homogeneous noncompact surface with Riemannian metric is diffeomorphic to a plane, cylinder or M\"obius strip; either it has
\begin{itemize}
\item
constant curvature or 
\item
its Gauss curvature has no critical points and its projective automorphisms preserve the foliation by gradient lines of the Gauss curvature.
\end{itemize}
\end{theorem}
\begin{proof}
A connected surface admits a transitive Lie group action if and only if it is diffeomorphic to one of
\begin{itemize}
\item
the compact:
\begin{itemize}
\item
the sphere,
\item
the real projective plane,
\item
the torus,
\item
the Klein bottle,
\end{itemize}
\item
or the noncompact:
\begin{itemize}
\item
the plane,
\item
the cylinder,
\item
the M\"obius strip;
\end{itemize}
\end{itemize}
\cite{Cartan1930} p. 29, \cite{Chevalley:1955.3} p. 351, \cite{Komrakov.Churyumov.Doubrov.1993},
\cite{Mostow:1950} p. 606.
The Gauss curvature on a compact surface has a maximum and minimum, hence a critical point, so our \(3\)-tensor vanishes at some point.
By homogeneity, every point is critical for Gauss curvature, hence Gauss curvature is constant.
For a noncompact surface, our \(3\)-tensor
\[
(G_2\omega^1-G_1\omega^2)\otimes(\omega^1\wedge\omega^2)
\]
has a well defined kernel: the lines \(G_2\omega^1=G_1\omega^2\), i.e. the zeroes of 
\[
*dG=G_2\omega_1-G_1\omega^2,
\]
i.e. the perpendicular lines to \(0=dG\), i.e. the gradient lines of \(G\).
\end{proof}

\chapter{Development and rolling}
\epigraph[author={\`E. Cartan}, source={La méthode du repère mobile, la théorie des groupes continus et les espaces généralisés},translation={On any such surface, the theory of curves will be identically the same as in the plane and nothing will change for application of the method of the moving frame.}]{Sur une telle surface, la th\'eorie des courbes sera identiquement la m\^eme que dans le plan et rien ne sera chang\'e \'a l'application de la m\'ethode du r\'epere mobile.}
\section{Definition}
\begin{problem}{ar}
Suppose that \(P\xrightarrow{\pi}Q\) is a smooth submersion.
A \emph{submersion connection}\define{submersion connection}\define{connection!submersion} is a vector subbundle \(V\subseteq TP\) splitting the tangent bundle into a direct sum \(TP=V\oplus\ker\pi'\).
Take a locally absolutely continuous curve \(q(t)\in Q\) defined on an interval \(I\subseteq\R\) and time \(t_0\in I\) and point \(p_0\in P\) with \(\pi(p_0)=q(t_0)\).
A \emph{lift}\define{lift!of curve} of \(q(t)\) through \(p_0\) is a locally absolutely continuous curve \(p(t)\in P\) for \(t\) in some interval \(J\subseteq I\) containing \(t_0\) with \(p(t_0)=p_0\) and \(\pi(p(t))=q(t)\) for all \(t\in I\).
Prove that there is a unique maximal lift, i.e. with maximal interval \(J\).
Give an example where \(J\subsetneq I\).
Prove that the lift of a reparametrization is the reparametrization of the lift.
\end{problem}
Take a homogeneous space \((X,G)\) and two \((X,G)\)-geometries \(\G\to M\) and \(\G'\to M'\), with Cartan connections \(\omega,\omega'\).
Consider when a curve in \(M\) ``looks like'' a curve in \(M'\).
On \(\G\times \G'\), there is a vector subbundle of the tangent bundle: the vectors on which \(\omega=\omega'\), spanned at each point by the vectors \(A_{\G}+A_{\G'}\), sums of the associated constant vector fields.
Intuitively, these are the tangent directions which ``look the same''.
This subbundle is \(H\)-invariant.
By the equivariant Frobenius theorem (corollary~\vref{corollary:equivariant.Frobenius}), this subbundle projects to a subbundle \(V\) of the tangent bundle of the fiber bundle \(\amal{\G}{H}{\G'}\to M\times M'\), of rank equal to the dimension of \(X\), and nowhere tangent to the vertical for either map
\[\label{page:double.submersion}
\begin{tikzcd}
&\amal{\G}{H}{\G'}\arrow[d]&\\
&M\times M'\arrow[dl]\arrow[dr]&\\
M&&M'.
\end{tikzcd}
\]
Hence \(V\) is a submersion connection for both maps.
Every lift of a  curve in \(M\) consists of a curve in \(\amal{\G}{H}{\G'}\), projecting to a curve in \(M'\), its \emph{development}\define{development} or \emph{developing curve}.\define{developing curve}
Say that \emph{curves in \(M\) develop to curves in \(M'\)}\define{curves develop} if this lift \(\ell(t)\) is globally defined, for as long as \(m(t)\) is defined, and hence \(m'(t)\) is defined then too.
\begin{problem}{tilting}
Take a Cartan geometry \(\G\to M\) and absolutely continuous paths \(p(t)\in\G\), \(h(t)\in H\).
Let \(q(t):=p(t)h(t)\). 
Prove the \emph{product rule}\define{product rule} of Cartan geometries:
\[
\dot{q}\hook\omega_{\G}=\dot{h}\hook\omega_H+\Ad_h^{-1}(\dot{p}\hook\omega_{\G}).
\]
\end{problem}
\begin{answer}{tilting}
Let 
\begin{align*}
q(t)&:=p(t)h(t),\\
A(t)&:=\dot{p}(t)\hook\omega,\\
B(t)&:=\dot{h}(t)\hook\omega_H.
\end{align*}
Compute
\begin{align*}
\dot{q}(t)
&=
\left.\frac{d}{ds}\right|_{s=0}
p(t+s)h(t+s),
\\
&=
\left.\frac{d}{ds}\right|_{s=0}
p(t)h(t+s)
+
\left.\frac{d}{ds}\right|_{s=0}
p(t+s)h(t),
\\
&=
\left.\frac{d}{ds}\right|_{s=0}
p(t)h(t)e^{sB(t)}
+
\RT{h(t)*}
\left.\frac{d}{ds}\right|_{s=0}
p(t+s),
\\
&=
B(t)_{\G}(q(t))
+
\RT{h(t)*}\dot{p}(t).
\end{align*}
Hence 
\[
\dot{q}\hook\omega=B(t)+\Ad_{h(t)}^{-1}A(t).
\]
By the same calculation,
\[
\frac{d}{dt}g(t)h(t)=B(t)+\Ad_{h(t)}^{-1}A(t).
\]
\end{answer}
\begin{lemma}\label{lemma:reparam}
Take a Cartan geometry \(\G\to M\) and a point \(p_0\in\G\) with image \(m_0\in M\).
Take an absolutely continuous path \(m(t)\in M\), \(a\le t\le b\).
It lifts to an absolutely continuous path \(p(t)\in\G\) starting at \(p_0\), i.e. \(p(a)=p_0\) and \(p(t)\) projects to \(m(t)\).
Develop that path: there is a unique absolutely continuous path \(g(t)\in G\) with \(g(a)=1\) and \(\dot{g}\hook\omega_G=\dot{p}\hook\omega_{\G}\).
Any other lift of \(m(t)\) has the form \(p(t)h(t)\) for some absolutely continuous path \(h(t)\in H\).
The development of \(p(t)h(t)\) is \(g(t)h(t)\).
If we reparameterize \(m(t)\), we reparameterize \(p(t)\) by the same reparameterization.
\end{lemma}
The same is true for a subgeometry in place of a Cartan geometry, with the same proof.
\begin{proof}
Since \(\G\to M\) is a fiber bundle, every path in \(M\) lifts to a path in \(\G\), locally by picking any smooth local trivialization, and then globally by patching together.
The development \(g(t)\) exists and is unique, being the solution of a Lie equation, by theorem~\vref{theorem:Lie.equations}.
Reparameterization clearly solves the Lie equation; apply existence and uniqueness of solution.
Apply the product rule in \(\G\) and \(G\) to see that tilting by the same \(h(t)\) stills solve the Lie equation; apply existence and uniqueness of solution.
\end{proof}
\begin{lemma}\label{lemma:to.the.model}
Curves in any Cartan geometry develop to the model.
\end{lemma}
Again the same is true for a subgeometry in place of a Cartan geometry, with the same proof.
\begin{proof}
Take a curve \(m(t)\in M\).
Lift to a curve \(p(t)\in \G\).
Let \(g(t)\) solve the Lie equation \(\omega_G=\dot{p}\hook\omega_{\G}\); see appendix~\vref{appendix:Lie.equations}.
Project \(g(t)\) to a curve \(x(t)\in X\).
\end{proof}
\section{Developing to the model is parallel transport}
In the special case when \(\G'\to M'\) is the model \(G\to X\), we want to see why development is essentially just parallel transport on the associated \(G\)-bundle \(\G_G\).
As when we defined development above, we get \(H\) to act on the right on \(\G\times G\):
\[
(p,g)h:=(ph,gh),
\]
denoting the quotient as \(\amal{\G}{H}{G}\), a principal left \(G\)-bundle.
Let \(\pi\) be the projection map
\[
(p,g)\in G\times\G\xmapsto{\pi}(p,g)H\in\amal{\G}{H}{G}.
\]

Recall that \(\G_G\) is the right \(G\)-bundle associated to \(\G_H=\G\), i.e. the quotient of \(\G\times G\) by a different \(H\)-action:
\[
(p,g)h:=(ph,h^{-1}g).
\]
Let \(\pi'\) be the projection map
\[
(p,g)\in G\times\G\xmapsto{\pi'}(p,g)H\in\G_G.
\]
Hence the map
\[
(p,g)\xmapsto{\iota}(p,g^{-1})
\]
descends to a diffeomorphism
\[
\amal{\G}{H}{G}\xrightarrow{\iota}\G_G,
\]
making commute
\[
\begin{tikzcd}
\G\times G\arrow[d,"\pi"]\arrow[r,"\iota"]&\G\times G\arrow[d,"\pi'"]\\
\amal{\G}{H}{G}\arrow[r,"\iota"]&\G_G
\end{tikzcd}
\]
and intertwining the left \(G\)-action with the right one by
\[
\iota\circ\LT{g}=\RT{g}^{-1}\circ\iota.
\]
By definition, the Cartan connection \(\omega_{\G}\) on \(\G\) is the pullback to \(\G=\G_H\subset\G_G\) of a connection form \(\omega_{\G_G}\) on \(\G_G\).
\begin{lemma}\label{lemma:forms.on.product}
\begin{align*}
\pi^*\iota^*\omega_{\G_G}&=\omega_G+\Ad_g^{-1}\omega_{\G},\\
\iota^*\pi^{\prime *}\omega_{\G_G}&=\Ad_g(\omega_{\G}-\omega_G).
\end{align*}
\end{lemma}
\begin{problem}{invert.geometry}
Prove it.
\end{problem}
\begin{answer}{invert.geometry}
Along \(\G\times\set{1}\), \(\pi'(p,1)=p\), so \(\pi^{\prime *}\omega_{\G_G}=\omega_{\G}\).
On that fiber, since \(\omega_{\G_G}\) is a connection, for any \(A\in\LieG\),
\[
A_{\G_G}\hook\omega_{\G_G}=A,
\]
so
\[
A_G\hook\pi^{\prime *}\omega_{\G_G}=A.
\]
Under left \(H\)-action, \(\pi'\LT{h}=\pi'\), so
\begin{align*}
\LT{h}^*\pi^{\prime *}\omega_{\G_G}
&=
(\pi^{\prime}\LT{h})^*\omega_{\G_G},
\\
&=
\pi^{\prime *}\omega_{\G_G}
\end{align*}
i.e. \(\pi^{\prime *}\omega_{\G_G}\) is \(H\)-invariant.
Under right \(G\)-action, 
\begin{align*}
\RT{g}^*\pi^{\prime *}\omega_{\G_G}
&=
\pi^{\prime *}\RT{g}^*\omega_{\G_G},
\\
&=
\pi^{\prime *}\Ad_g^{-1}\omega_{\G_G},
\\
&=
\Ad_g^{-1}\pi^{\prime *}\omega_{\G_G},
\end{align*}
i.e. \(\pi^{\prime *}\omega_{\G_G}\) transforms in the adjoint.
Hence, at each point \((p,g)\),
\[
\pi^{\prime *}\omega_{\G_G}=\omega_G+\Ad_g^{-1}\omega_{\G}.
\]
By lemma~\vref{lemma:invert.Maurer.Cartan},
\[
\iota^*\omega_G=-\Ad_g\omega_G.
\]
Hence, intertwining the actions,
\begin{align*}
\iota^*\pi^{\prime *}\omega_{\G_G}
&=
\iota^*\pi^{\prime *}(\omega_G+\Ad_g^{-1}\omega_{\G}),
\\
&=-\Ad_g\omega_G+\Ad_g\omega_{\G},
\\
&=
\Ad_g(\omega_{\G}-\omega_G).
\end{align*}
\end{answer}
\begin{theorem}
Take a Cartan geometry \(\G\to M\) with model \((X,G)\).
The map
\[
\amal{\G}{H}{G}\xrightarrow{\iota}\G_G
\]
identifies the development lift of each curve from \(M\) to the horizontal lift of that curve for the Cartan connection, as a connection on the principal bundle \(\G_G\).
\end{theorem}
\begin{proof}
Applying the map \(\iota\), the vectors \(A_{\G}+A_{\G'}\) are the vectors on which \(\iota^*\pi^{\prime *}\omega_{\G_G}\) vanishes, so are mapped by \(\iota\) precisely to the horizontal vectors on the right principal \(G\)-bundle \(\G_G\).
\end{proof}
\section{Rolling is development}
Let \((X,G)\) a \(3\)-dimensional homogeneous space with a \(G\)-invariant Riemannian metric on \(X\) so that \(G\) acts transitively on the orthonormal frames of \(X\).
\begin{problem}{as}
Prove that \(X\) has constant curvature under that metric.
Prove that \(X\) is isometric to Euclidean space, a sphere, a real projective space, or hyperbolic space. Hint: \cite{Wolf:1967}.
What are all of the possibilities of \((X,G)\)?
\end{problem}
Fix two immersed surfaces \(M,\otM\) in \(X\).
\[
\begin{tikzpicture}
\coordinate (ball-center) at (2,.25);
\shadedraw (-1.5,0) -- (-.5,.5) -- (1,0) -- (.5,-.75) -- cycle;
\shadedraw[ball color=white] (ball-center) circle (.5cm);
\end{tikzpicture}
\]
We want to roll \(\otM\) on \(M\) 
\[
\begin{tikzpicture}
\coordinate (ball-center) at (2,.25);
\shadedraw (-1.5,0) -- (-.5,.5) -- (1,0) -- (.5,-.75) -- cycle;
\coordinate (moved-ball-center) at (0,.25);
\shadedraw[ball color=white] (moved-ball-center) circle (.5cm);
\shadedraw[ball color=white,opacity=0] (ball-center) circle (.5cm);
\end{tikzpicture}
\]
along some path \(x(t)\in M\):
\[
\begin{tikzpicture}
\coordinate (ball-center) at (2,.25);
\shadedraw[ball color=white,opacity=0] (ball-center) circle (.5cm);
\shadedraw (-1.5,0) -- (-.5,.5) -- (1,0) -- (.5,-.75) -- cycle;
\draw[midarrow] plot [smooth] coordinates {(-.8,.1) (-.5,.15) (.2,-.1) (.3,-.3) (.7,-.1)};
\coordinate (moved-ball-center) at (0,.25);
\shadedraw[ball color=white,opacity=.6] (moved-ball-center) circle (.5cm);
\end{tikzpicture}
\]
Model the rolling by applying a time-varying rigid motion \(g(t)\) to the fixed surface \(\otM\) to take it to touch \(M\):
\[
\begin{tikzpicture}
\coordinate (ball-center) at (2,.25);
\shadedraw[ball color=white] (ball-center) circle (.5cm);
\shadedraw (-1.5,0) -- (-.5,.5) -- (1,0) -- (.5,-.75) -- cycle;
\draw[midarrow] plot [smooth] coordinates {(-.8,.1) (-.5,.15) (.2,-.1) (.3,-.3) (.7,-.1)};
\coordinate (moved-ball-center) at (0,.25);
\shadedraw[ball color=white,opacity=.6] (moved-ball-center) circle (.5cm);
\path[-stealth,bend right,very thick,line cap=round](ball-center) edge node[black,midway,above]{\(g(t)\)} (moved-ball-center);
\end{tikzpicture}
\]
and to take some path \(\otx(t)\in\otM\) to \(g(t)\otx(t)=x(t)\), and \(g(t)_*T_{\otx(t)}\otM=T_x M\), the surfaces touch and remain tangent.
But to say that this rigid motion models rolling, we need physics to encode that the rolling surface is not slipping or twisting as it rolls.

Not slipping means: at the instant that a chosen point of the rolling surface touches the fixed surface, it does not slide along the fixed surface, so its velocity at the point of contact is zero.
Not twisting means: at the instant that a chosen point of the rolling surface touches the fixed surface, bringing a vector tangent to the rolling surface to coincide with one on the fixed surface, the end point of that moving tangent vector is not being twisted tangentially, i.e. the end point has velocity  normal to both surfaces.

First, let's derive a no slip equation.
Picture a particular point \(\otx_0\in\otM\), not varying over time, and watch where it goes when we roll: along the path \(g(t)\otx_0\in X\) in space.
This path is perhaps never a point of \(M\).
Even if \(g(t)\otx_0\) sometimes touches \(M\), perhaps it does not remain in contact with \(M\) for any interval of time.
The physical condition for rolling without slipping is thus that, if a chosen point \(\otx_0\in\otM\) happens to be one of the points of the contact curve at some time \(t=t_0\), i.e. 
\[
\otx_0=\otx(t_0),
\]
then the point \(\otx_0\) is carried by the rigid motion \(g(t)\) on a path in \(X\)  which has zero velocity at that time \(t=t_0\).
So no slipping means precisely that, for any time \(t_0\),
\[
\left.\frac{d}{dt}\right|_{t=t_0}g(t)\otx(t_0)=0.
\]
By the obvious product rule for Lie group actions, we can write this in terms of the path \(\otx(t)\) on which we roll as
\[
\frac{d}{dt}g(t)\otx(t)=g(t)_*\dot\otx(t).
\]
\begin{problem}{product.rule}
What is the product rule for Lie group actions?
Prove it.
\end{problem}
\begin{answer}{product.rule}
Suppose that a Lie group \(G\) acts smoothly on a manifold \(M\).
Write the associated Lie algebra action as: for any \(A\in\LieG\) has associated vector field \(A_M\) on \(M\), so that for any point \(m\in M\):
\[
\left.\frac{d}{dt}e^{tA}m\right|_{t=0}=A_M(m).
\]
Recall that 
\[
\lb*{A_M}{B_M}=-\lb{A}{B}_M.
\]
Take a path \(g(t)\in G\) and a path \(m(t)\in M\).
Consider the path \(g(t)m(t)\).
Let \(A(t):=\dot{g}(t)\hook\omega_G\).
Our product rule is
\[
\frac{d}{dt}g(t)m(t)=(A(t))_M(m(t))+g(t)_*\frac{dm}{dt}.
\]
\end{answer}
In other words, no slip means exactly that the rigid motion identifies the velocity of  the paths on the fixed surfaces \(\otM\) and \(M\).
(To me, this is not obvious.)

To write the no twist equation, we need the Levi--Civita connection of \(X\); remember that it measures twisting of frames.
As in section~\vref{section:pseudo.Riem}, denote the frame bundle of \(X\) as \(\G_X\xrightarrow{\pi}X\), but write the points of the frame bundle as
\[
p=(x,e_1,e_2,e_3)
\]
where \(e_1,e_2,e_3\) can be any orthonormal basis of \(T_x X\).
\begin{problem}{at}
What is the action of \(\Orth{3}\) on \(\G_X\)?
What natural bundle isomorphism identifies \(\G_X\) with the bundle we called \(\G\) in section~\vref{section:pseudo.Riem}?
\end{problem}
Denote the Cartan connection of the Riemannian metric on \(X\)  as soldering form \(\sigma\) and Levi-Civita connection \(\gamma\).
Write the three component \(1\)-forms of the soldering \(\sigma\) as 
\[
\sigma_1,\sigma_2,\sigma_3,
\]
so that, for any tangent vector \(v\in T_p\G_X\), 
\[
\pi'(p)v=\sigma_i(v)e_i.
\]
Similarly, \(\gamma\) is valued in \(\LieSO{3}\), say \(\gamma=(\gamma_{ij})\), so that \(\gamma_{12}\), applied to a moving frame
\[
p(t)=(x(t),e_1(t),e_2(t),e_3(t))\in\G_X
\] 
measures the angular velocity with which \(e_2(t)\) turns in the direction of \(e_1(t)\). 
We denote \(\gamma,\sigma\) as \(\otgamma,\otsigma\) when we work above points of \(\otM\).

A frame 
\[
p=(x,e_1,e_2,e_3)
\]
is \emph{adapted}\define{frame!adapted}\define{adapted!frame} to the surface \(M\) if \(x\in M\) and \(e_1,e_2\in T_x M\).
So our no twist equation is precisely that, if \(\otp_0\) is a frame adapted to \(\otM\) at a point \(\otx_0\), and we let \(p(t):=g(t)\otp_0\), and if, at some time \(t=t_0\), \(\otx_0=\otx(t)\) then
\[
\dot{p}(t_0)\hook\gamma_{12}=0.
\]
If instead we replace the fixed frame \(\otp_0\) with a moving adapted frame \(\otp(t)\), and let \(p(t):=g(t)\otp(t)\), again the product rule gives
\begin{align*}
\dot{p}\hook\gamma_{12}
&=\left(g(t)_*\dot\otp\right)\hook\gamma_{12},
\\
&=\dot\otp\hook g(t)^*\gamma_{12},
\\
&=\dot\otp\hook\otgamma_{12}.
\end{align*}

Let \(\G_M\to M\) be the \emph{adapted frame bundle}\define{adapted!frame bundle} of the surface \(M\), i.e. the set of all adapted frames.
The subgroup 
\[
H_{\text{surface}}:=\Orth{2}\times\set{\pm 1}\subset\Orth{3}
\]
acts on \(\G_M\) by orthogonal transformation of \(e_1,e_2\) vectors in the tangent plane of the surface, and reflection of \(e_3\) in the normal line to the surface.
\begin{problem}{au}
Prove that \(\G_M\to M\) is a smooth \(H_{\text{surface}}\) bundle smoothly immersed in \(\G_X\).
Prove that \(\sigma_3=0\) on \(\G_M\) and that
\[
\begin{pmatrix}
\gamma_{13}\\
\gamma_{23}
\end{pmatrix}
=
\begin{pmatrix}
a_{11}&a_{12}\\
a_{21}&a_{22}
\end{pmatrix}
\begin{pmatrix}
\sigma_1\\
\sigma_2
\end{pmatrix}
\]
for unique functions \(a_{ij}=a_{ji}\) on \(\G_M\).
How do these \(a_{ij}\) transform under \(H_{\text{surface}}\)-action?
\end{problem}

Let \(I\subset\R\) be the interval of time \(t\) during which the rolling occurs.
Let \(\G_I\) denote the set of all 
\[
(t,e_1,e_2,e_3)
\]
so that \(t\in I\) and 
\[
(x(t),e_1,e_2,e_3)\in\G_M,
\]
i.e. the pullback bundle
\[
\begin{tikzcd}
\G_I\arrow[r]\arrow[d]&\G_M\arrow[d]\\
I\arrow[r]&M.
\end{tikzcd}
\]
So \(\G_I\to I\) is a smooth principal \(H_{\text{surface}}\)-bundle smoothly immersed in \(\G_M\).
Similarly for \(\otM\)
\[
\begin{tikzcd}
\otG_I\arrow[r,"x"]\arrow[d]&\G_{\otM}\arrow[d]\\
I\arrow[r]&\otM.
\end{tikzcd}
\]
At each point \(t\in I\), the rigid motion \(g(t)\) identifies these bundles
\[
\begin{tikzcd}
\otG_I\arrow[rr,"g"]\arrow[dr]&&\G_I\arrow[dl]\\
&I&
\end{tikzcd}
\]
This bundle isomorphism is \(H_{\text{surface}}\)-equivariant, since rotating or reflecting an adapted frame on \(\otM\) rotates or reflects the resulting frame on \(M\) via \(g(t)\).
So \(g\) is a principal bundle isomorphism.

The forms \(\gamma,\sigma\) on \(\G_M\) pullback to \(\G_I\).
Since \(\sigma\) is semibasic on \(\G_X\), it is also semibasic on \(\G_I\).
Take the tangent vector \(\partial_t\) tangent to \(I\).
At some point of \(\otG_I\), pick a tangent vector \(\otv\in T_{\otp}\otG_I\) mapping to \(\partial_t\).
Let \(v\in T_p\G_I\) be its image under the bundle isomorphism.
Compute
\begin{align*}
v\hook\sigma_i
&=
\left<e_i,\pi_*(x(t),e_1,e_2,e_3)x_*v\right>,
\\
&=
\left<e_i,\dot{x}(t)\right>,
\\
&=
\left<g(t)_*\ote_i,g(t)_*\dot\otx(t)\right>,
\\
&=
\left<\ote_i,\dot\otx(t)\right>,
\\
&=
\otv\hook\otsigma_i.
\end{align*}
Hence our bundle isomorphism yields \(g^*\sigma=\otsigma\).
Similarly \(g^*\gamma_{12}=\otgamma_{12}\).

We can interpret these equations geometrically.
The constraint \(g^*\sigma=\otsigma\) says that velocity vectors of the curves are matched up by the rigid motion.
Suppose that one of the curves in an immersed curve, so that both then are.
We can reparameterize one to have unit speed, and then they both do.
Along any unit speed curve on \(M\), there is bundle of orthonormal frames with first leg \(e_1\) tangent to the curve, and second leg normal, given by \(\sigma_2=0\), and then \(\gamma_{12}=\kappa\sigma_1\) is the geodesic curvature of the curve.
So the constraint \(g^*\gamma_{12}=\otsigma\) says that the curves have the same geodesic curvature.
In particular, geodesics roll only along geodesics, and curves of constant curvature only along curves of the same constant curvature.

So the graph of the bundle isomorphism is an integral surface of the foliation
\[
0=\sigma-\otsigma=\gamma_{12}-\otgamma_{12}
\]
on \(\G_{\otM}\times\G_M\).
This integral surface is invariant under diagonal \(H_{\text{surface}}\).
That diagonal action preserves this foliation.
By the equivariant Frobenius theorem, the quotient \(5\)-manifold
\[
\begin{tikzcd}
&\amal{\G_M}{H_{\text{surface}}}{\G_{\otM}}\arrow[d]&\\
&M\times\otM\arrow[dl]\arrow[dr]&\\
M&&\otM
\end{tikzcd}
\]
is foliated by curves, which are the integral curves of this quotient system, giving two submersion connections.
Each curve represents the rolling of one surface on the other, determining which points and frames are carried to which, and so uniquely determining the required rigid motion at each time.
Once we pick a curve \(x(t)\in M\), it lifts to a unique integral curve of this foliation through each point \((p_0,\otp_0)H_{\text{surface}}\), i.e. a choice of linear isometry of the tangent spaces, at the points \(x_0,\otx_0\).

The surface \(M\) inherits a Riemannian metric, which is a Cartan geometry given by its soldering form \(\sigma_1,\sigma_2\) and Levi-Civita connection \(\gamma_{12}\); similarly for \(\otM\).
The equations of rolling above are precisely the equations of development of our two surfaces \(M\) and \(\otM\)  as Cartan geometries; our diagram above is the diagram~\vpageref{page:double.submersion}.
The equations of rolling are thus invariant under isometry of the surfaces, as abstract surfaces with Riemannian metrics, independent of the choice of isometric immersions of the surfaces into \(X\).
(To me, this is not obvious.)
Imagine isometrically immersing the same surfaces into different constant curvature geometries \(X\), or bending each of the surfaces without stretching or squishing, and we uncover the same rolling motions.
\begin{problem}{av}
Take an integral curve
\[
(p(t),\otp(t))H_{\text{surface}}.
\] 
Suppose that, for each time \(t\), some \(g(t)\in G\) exists taking \(\otp(t)\) to \(p(t)\).
Prove that then \(g(t)\) is uniquely and smoothly determined by the curve \(x(t)\) and linear isometry of tangent spaces at \(x_0,\otx_0\). Hint: you might use theorem~\vref{thm:strong}.
\end{problem}
If we repeat the story with \((X,G)\) any homogeneous Riemannian \(3\)-manifold, we might not be able to bring \(\otM\) tangent to \(M\) at a specified point by action of an element of \(G\).
\begin{problem}{aw}
If \((X,G)\) is a homogeneous Riemannian \(3\)-manifold and \(G\) acts transitively on tangent planes of \(X\), prove that \((X,G)\) has constant curvature.
\end{problem}
Even still, the notion of development remains well defined purely in terms of the Riemannian geometry of the two surfaces, without any ambient \(3\)-manifold \(X\).

Note that, in a more realistic model of rolling, one surface could get ``stuck'' in the other; this relies on understanding how the surfaces are immersed into \(X\).

Suppose that \((X,G)\) is a constant curvature homogeneous Riemannian \(3\)-manifold.
Consider \(g\) as a map
\[
(p,\otp)\in\G_X\times\G_X\to g\in G
\]
defined by \(gp=\otp\).
Identify \(\G_X=G\), so this map becomes
\[
(p,\otp)\in G\times G\mapsto g=\otp p^{-1}\in G.
\]
Using the convenient notation \(g^{-1}dg\) rather than \(\omega_G\), we can write
\[
g^{-1}\,dg=\Ad_p(\otp^{-1}d\otp-p^{-1}dp).
\]
The equations of development along \(\G_M\times\G_{\otM}\) are the equations that say that \(g^{-1}dg\) has the form 
\[
g^{-1}dg
=
\Ad_p
\begin{pmatrix}
0&0&\otgamma_{13}-\gamma_{13}&0\\
0&0&\otgamma_{23}-\gamma_{23}&0\\
\gamma_{13}-\otgamma_{13}&\gamma_{23}-\otgamma_{23}&0\\
0&0&0&0
\end{pmatrix}.
\]
The shape operators of the surfaces are defined by
\begin{align*}
\gamma_{i3}&=a_{ij}\sigma_j,\\
\otgamma_{i3}&=\ota_{ij}\otsigma_j,\\
\end{align*}
At any one moment, we can left translate to arrange that \(p=I\in G\).
The rolling equations say precisely that
\[
g^{-1}dg
=
\begin{pmatrix}
0&0&(a_{1j}-\ota_{1j})\sigma_j&0\\
0&0&(a_{2j}-\ota_{2j})\sigma_j&0\\
(\ota_{1j}-a_{1j})\sigma_j&(\ota_{2j}-a_{2j})\sigma_j&0\\
0&0&0&0
\end{pmatrix}.
\]
The infinitesimal rotation of frame caused by the rolling is trivial in the tangent directions, but controlled in the normal directions by the difference of shape operators.
We recover
\begin{theorem}[Nomizu \cite{Nomizu:1978}]
The rolling of a moving surface along a path on a fixed surface, inside any constant curvature Riemannian \(3\)-manifold, is carried out by a curve in the rigid motion group of the \(3\)-manifold.
This curve is an immersed curve except precisely at moments when the surfaces have the same value of their shape operators on the tangent vector to the path.
\end{theorem}

\chapter{Developability}
\section{Definition}
A Cartan geometry is \emph{developable}\define{Cartan!geometry!developable}\define{developable!Cartan geometry} if curves in any Cartan geometry with the same model develop to it.
\begin{example}
Let \((X,G)\) be the real number line \(X=\R\) with \(G=\R^{\times}\ltimes\R\) the group of affine transformations.
Take an \((X,G)\)-geometry \(\G\to M\) on a connected curve \(M\).
Suppose first that \(M\) is diffeomorphic to \(\R\).
The developing map \(M\to X\) makes \(M\) an open interval of \(X=\R\), uniquely determined up to affine transformation, giving the classification of \((X,G)\)-geometries on any connected and simply connected curve.
If \(M\) is developable, we can develop \(X\to M\), and also \(M\to X\), hence compose 
\[
X\to M\to X,
\]
giving a development \(X\to X\), so an element of \(G\), so an affine transformation, so \(X\) with the standard \((X,G)\)-geometry is the unique developable \((X,G)\)-geometry, up to isomorphism, on any connected and simply connected curve.

If \(M\) is instead diffeomorphic to \(S^1\), then by the same argument \(\tilde{M}\) is developable if and only if \(\tilde{M}=X\) with the standard structure, so \(M\) is the quotient by some affine transformation with no fixed points, so a translation.
Hence the unique developable \((X,G)\)-geometry on any connected but not simply connected curve is the standard \(\R/\Z\), up to isomorphism.
\end{example}
\begin{example}
A Cartan geometry on a manifold is developable just when the induced Cartan geometry on each component of the manifold is developable.
\end{example}
\begin{example}
A Cartan geometry on a connected manifold is developable just when the pullback Cartan geometry on any covering space is developable.
\end{example}
\begin{example}
Even a homogeneous flat Cartan geometry on a compact manifold need not be developable.
The model projective connection on projective space pulls back to affine space, to a projective connection which is not developable, since we can't develop the projective lines, which are compact, to the affine lines, which are not.
Affine space is a covering space of the torus, and the projective connection quotients to the torus, to give again a projective connection which is not developable, since developments lift and project in covering maps.
\end{example}
\begin{problem}{ax}
Prove that a Cartan geometry is developable if and only if there is a point in the model so that every curve in the model, starting at that point, develops to the Cartan geometry.
\end{problem}
Similarly, say that \emph{curves in \(\G\) develop to curves in \(\G'\)}\define{curves develop} if, for any absolutely continuous curve \(p(t)\in\G\), defined on some interval containing a time \(t_0\), and any point \(p_0'\in\G'\), there is an absolutely continuous curve \(p'(t)\in\G'\) with the same Darboux derivative
\[
\frac{dp}{dt}\hook\omega=\frac{dp'}{dt}\hook\omega'.
\]
\begin{lemma}
Take Cartan geometries \(\G\to M\) and \(\G'\to M'\) with the same model.
Curves in \(M\) develop to curves in \(M'\) just when curves in \(\G\) develop to curves in \(\G'\). 
\end{lemma}
\begin{proof}
If curves in \(\G\) develop to curves in \(\G'\), take any curve \(m(t)\in M\) and pick some absolutely continuous lift \(p(t)\) of it, and then construct such a \(p'(t)\) and then project to \(m'(t)\).
Conversely, suppose that curves in \(M\) develop to curves in \(M'\).
Take an absolutely continuous curve \(p(t)\in\G\), defined on some interval containing a time \(t_0\), and any point \(p_0'\in\G'\).
Denote its Darboux derivative by
\[
A(t):=\frac{dp}{dt}\hook\omega.
\]
Project \(p(t)\) to a curve \(m(t)\in M\).
Develop as above: we construct a curve \(\ell(t)\in\amal{\G}{H}{\G'}\) through the point \(\ell_0:=(p(t_0),p_0')H\).
Such a curve lifts to a curve in \(\G\times\G'\), tangent to the preimage of \(V\), uniquely up to \(H\)-action, and uniquely if we in addition specify that it map to \(p(t)\in\G\).
The preimage of \(V\) consists of the vectors of the form \(A_{\G}+A_{\G'}\).
\end{proof}
\begin{corollary}
A Cartan geometry is developable just when any of its lifts or drops are developable.
\end{corollary}
\begin{theorem}%
[Clifton I \cite{Clifton:1966}]%
\label{theorem:Clifton}%
\define{theorem!Clifton I}%
\define{Clifton!theorem I}
The following are equivalent conditions on a Cartan geometry:
\begin{itemize}
\item
The Cartan geometry is developable, i.e. curves in any Cartan geometry with the same model develop to it.
\item
Curves in the model develop to it.
\item
Smooth curves in the model develop to it.
\item
The scaffold is a complete Riemannian manifold.
\item
The bundle is a complete Riemannian manifold.
\end{itemize}
\end{theorem}
\begin{proof}
Use the Riemannian metric \(\left<\omega,\omega\right>\) on \(\G\), for some \(H_c\)-invariant metric on \(\LieG\).
The map \(\G\to M_c=\G/H_c\) is proper, so completeness is preserved and reflected.

We can develop to the model by lemma~\vref{lemma:to.the.model}.
If we can develop from any geometry with the same model, we can develop from the model.
If we can develop from the model \((X,G)\) to \(M\), then we can develop from any Cartan geometry \(M'\), first to \(X\), and then to \(M\).

Suppose that the metric on \(\G\) is complete.
When we develop a smooth curve from \(G\), our velocity stays bounded in any relatively compact time interval, so we stay inside some compact metric ball, in which we can bound from below the amount of time we can continue to solve the ordinary differential equation of development.
So we can develop curves from \(G\) to \(\G\), so from \(X\) to \(M\), so we can develop to \(M\) from any \((X,G)\)-Cartan geometry.

Suppose that all curves develop from the model to \(M\).
Pick a constant speed geodesic on \(\G\), defined on an open interval.
Take its development to \(G\), by solving the Lie equation.
The Darboux derivative is a constant length vector, smoothly varying on that interval.
Extend it to remain bounded and thus extend to an absolutely continuous curve defined for all time.
Develop it back to \(\G\): it extends our geodesic, and it continues to locally minimize distance up to the ends of our interval, so is still geodesic at those points, and so admits a smooth extension.
The set of points we can reach on curves of a given length from a given point is thus compact, and by the Hopf--Rinow theorem \cite{Chavel:2006} \S I7, \cite{Gromov:2007} p. 9, \cite{Petersen:2016} p. 137 theorem 16 the Riemannian metric is complete.

Suppose that smooth curves develop from the model to \(M\).
Pick a constant speed geodesic on \(\G\), defined on an open interval.
Take its development to \(G\), by solving the Lie equation.
Since the development in \(G\) has bounded velocity, it is Lipschitz.
Its velocity bound keeps it, near each end of our interval, inside some small closed ball.
This ball is compact because the metric on \(G\) is complete.
The radius of this ball gets smaller as we approach the end of our interval, so we approach the intersection point of these various compact balls, a single point.
We uniquely extend the curve to be continuous on the closure of our interval, asking it to reach that point.
The curve remains Lipschitz, precisely because it stays in these balls.

Our curve is smooth in the interior of that interval.
Reparameterize it to have velocity vanishing at the end points, to all orders; it is smooth.
Develop it back to \(\G\), extended, and then undo the reparameterization, to get a Lipschitz extension.
Apply the Hopf--Rinow theorem \cite{Chavel:2006} \S I7, \cite{Gromov:2007} p. 9, \cite{Petersen:2016} p. 137 theorem 16 as we did before: the metric is complete.

By lemma~\vref{lemma:completeness.shared}, the scaffold is a complete Riemannian manifold just when \(\G\) is.
\end{proof}
\begin{example}
A Riemannian manifold is its own scaffold, hence is developable just when complete as a metric space, i.e. just when geodesically complete.
\end{example}
\begin{proposition}
A Cartan geometry is developable precisely when all absolutely continuous curves in any Cartan geometry with the same model develop to it.
\end{proposition}
\begin{proof}
Take a homogeneous space \((X,G)\).
Denote by \(\LieG\) the Lie algebra of \(G\).
Pick a positive definite inner product \(g_0\) on \(\LieG\).
Take two \((X,G)\)-geometries \(\G\to M\), \(\G'\to M'\), say with Cartan connections \(\omega,\omega'\).
Let \(g:=\omega^*g_0\), \(g':={\omega'}^*g_0\).
Take an absolutely continous path \(p(t)\in\G\).
By theorem~\vref{thm:reparam}, \(p(t)=P(\rho(t))\) for some unit speed path \(P(t)\).
The geometry \(\G'\to M'\) is developable just when its metric \(g'\) is complete.
The equation \(\dot P\hook\omega=\dot Q\hook\omega'\) has a unique local solution with given initial condition, by absolute continuity.
By compactness of the geodesic balls, the solution extends to a compact interval, hence extends to a unique global solution.
The curve \(q(t)=Q(\rho(t))\) is the development of \(p(t)\).
\end{proof}
\section{Uniform limits}
Pick a Lie group \(G\) with Lie algebra \(\LieG\).
The \emph{velocity space} \(\mathscr{V}_{\LieG}\) is the set of tuples \((t_0,g_0,I,A)\) with \(I\subseteq\R\) a closed interval, \(t_0\in I\), and \(I\xrightarrow{A}\LieG\) locally integrable and \(g_0\in G\).
Topologize by the obvious product topology.

Take an \((X,G)\)-geometry \(\G\to M\).
The \emph{path space} \(\mathscr{P}_\G\) is the set of pairs \((t_0,I,p)\) with \(I\subseteq\R\) a closed interval and \(p\) and absolutely continuous path defined on \(I\) and \(t_0\in I\).
Topologize with compact open \(C^0\) topology on the path and the topology of integrable functions on its velocity computed via \(\dot p\hook\omega\).
\begin{corollary}
Applying the Cartan connection to the velocity of a path,
\[
\mathscr{P}_\G\to\mathscr{V}_{\LieG}
\]
is a homeomorphism to an open subset of the velocity space.
\end{corollary}
\begin{proof}
See theorem~\vref{thm:continuity}.
\end{proof}
\begin{corollary}
The operation of development from one Cartan geometry to another, when the geometries have the same model, is defined on an open set in the path space of the source, and is a homeomorphism to an open set in the path space of the target.
\end{corollary}
\begin{theorem}\label{thm:abs.cont.int.curves}
Take a manifold \(M\) and a smooth vector subbundle \(V\subseteq TM\).
Take a closed interval \(I\subseteq\R\).
Suppose that \(p_1,p_2,\dots\) are absolutely continuous \(V\)-tangent curves defined on \(I\), and that they converge uniformly (but perhaps with no control on their derivatives) to an absolutely continuous curve \(p\), and that their velocities are uniformly bounded in \(L^1\).
Then \(p\) is \(V\)-tangent.
\end{theorem}
\begin{proof}
Apply theorem~\vref{thm:lift.continuity}.
\end{proof}
\begin{corollary}\label{cor:dev.uniform.limit}
Take a homogeneous space \((X,G)\) and \((X,G)\)-geometries \(\G\to M\), \(\G'\to M'\).
Take a closed interval \(I\subseteq\R\), a time \(t_0\in I\).
Take a sequence \(p_1,p_2,\dots\) of absolutely continuous paths in \(\G\).
Suppose that these paths converge uniformly (but perhaps with no control on derivatives) to an absolutely continuous path \(p\) in \(\G\).

Take their developments \(p'_1,p'_2,\dots\) and \(p'\) through some points of \(\G'\), so that \(p'_j(t_0)\to p'(t_0)\).
Suppose that \(p'\) is defined on \(I\).
Then so are all but finitely many of \(p_1',p_2',\dots\) and the paths \(p'_1,p'_2,\dots\) converge uniformly to \(p'\).
Moreover
\[
\dot{p}_j\hook\omega\to\dot{p}\hook\omega\text{ in \(L^1\)},
\]
if and only if
\[
\dot{p}'_j\hook\omega\to\dot{p}'\hook\omega\text{ in \(L^1\)}.
\]
\end{corollary}
\begin{proof}
Apply theorem~\vref{thm:lift.continuity} to the manifold \(\G\times\G'\) and the subbundle \((\omega=\omega')\).
\end{proof}
\begin{theorem}
Take a homogeneous space \((X,G)\) and an \((X,G)\)-geometry \(\G\to M\) with Cartan connection \(\omega\).
Take a closed interval \(I\subseteq\R\), a time \(t_0\in I\).
Suppose that \(p\) is an absolutely continuous path in \(\G\).
Let \(A:=\dot p\hook \omega\).
Suppose that \(I\xrightarrow{A_j}\LieG\) is a sequence of \(L^1\) maps with
\[
\int_{t_0}^t A_j\to\int_{t_0}^t A,
\]
uniformly in \(t\).
Pick points \(x_1,x_2,\dots\to p(t_0)\).
Then, for all but finitely many \(A_j\), there is a unique absolutely continuous path \(I\xrightarrow{p_j}\G\) with \(p_j(t_0)=x_j\) and \(A_j=\dot p_j\hook\omega\).
Moreover \(p_j\to p\) uniformly.
\end{theorem}
\begin{proof}
We can consider \(\G\) with its Cartan connection \(\omega\) to be a Cartan geometry modelled on \((X',G'):=(\LieG,\LieG)\), with the Lie algebra \(G'=\LieG\) acting on itself \(X'=\LieG\) by translation.
The curves 
\begin{align*}
\bunderline{p}(t)&:=\int_{t_0}^t A(s)\,ds,\\
\bunderline{p}_j(t)&:=\int_{t_0}^t A_j(s)\,ds
\end{align*}
develop to \(p,p_j\), and converge uniformly.
Apply corollary~\vref{cor:dev.uniform.limit}.
\end{proof}
\begin{example}
Take a homogeneous space \((X,G)\).
Denote by \(\LieG\) the Lie algebra of \(G\).
Take \(B,C\in\LieG\).
For each positive integer \(n\), define some function \(A_n(t)\) with \(A_n=B\) half the time and \(C\) the other half, switching rapidly, say
\[
A_n(t):=
\begin{cases}
B, &\lceil 2^nt\rceil \text{ odd},\\
C, &\text{otherwise}.
\end{cases}
\]
Let \(A:=(B+C)/2\).
Then 
\[
\int_{t_0}^t A_n(s) \, ds\to \int A\,ds=(t-t_0)A,
\]
uniformly.
So by our theorem, we obtain the \emph{Trotter product formula}: for any \((X,G)\)-geometry \(\G\to M\),
\[
(e^{t C/2^n}e^{t B/2^n})^{2^n}\to e^{t(B+C)/2},
\]
wherever the right hand side is defined, which ensures that the left hand side is also defined for large enough \(n\).
Careful: where the left hand side is defined, even for every \(n\), the right hand side might not be.
\end{example}
\section{Developability and induction}
\begin{theorem}
Suppose that \((X,G)\to(X',G')\) is a morphism of homogeneous spaces and that \(X\to X'\) is a local diffeomorphism and \(X'\) is connected.
The following are equivalent:
\begin{itemize}
\item
The map \(X\to X'\) is onto.
\item
Every \((X,G)\)-geometry is developable if and only if its  induced \((X',G')\)-geometry is developable.
\item
There is an \((X,G)\)-geometry whose  induced \((X',G')\)-geometry is developable.
\end{itemize}
\end{theorem}
\begin{proof}
Suppose that \(X\to X'\) is onto.
It is a covering map, by homogeneity, so paths on \(X'\) lift to paths on \(X\).
Every \((X,G)\)-geometry which is developable can develop that lifted path, so its induced \((X',G')\)-geometry is developable.
Conversely, if the induced \((X',G')\)-geometry is developable then we can use it to take a path in \(X\), map it to \(X'\), develop the path from \(X'\), and then see that the development is a development of the path from \(X\) because the pullback bundles of the geometries are identified with those in the model.

Suppose that every \((X,G)\)-geometry is developable if and only if its  induced \((X',G')\)-geometry is developable.
Taking the model, we see that there is an \((X,G)\)-geometry whose  induced \((X',G')\)-geometry is developable.

Suppose that there is an \((X,G)\)-geometry whose  induced \((X',G')\)-geometry is developable.
Given a point of \(X'\) in the image of \(X\to X'\) and another point of \(X'\), connect the two by a path, develop it to the \((X,G)\)-geometry, and then develop it back to \(X\), to see that \(X\to X'\) is onto.
\end{proof}
\begin{example}
Take \((X,G)\to(X',G')\) the obvious morphism from Euclidean space \((X,G)=(\R^n,\Orth{n}\ltimes\R^n)\) to affine space \((X',G')=(\R^n,\GL{n}\ltimes\R^n)\).
Every \((X,G)\)-geometry is its own scaffold, so is developable just when complete in its metric.
So the induced affine geometry of its affine connection is developable.
\end{example}
\begin{example}
The same argument works for any reductive geometry: it is developable just when its induced affine connection is developable.
\end{example}
\begin{example}
Take \((X,G)\to(X',G')\) the obvious morphism from affine space to projective space.
Since \(X\to X'\) is \emph{not} onto, every affine connection induces a nondevelopable projective connection.
\end{example}
\section{Automorphisms and developability}
\begin{theorem}
If the automorphisms of a Cartan geometry \(\G\to M\) act locally transitively on the scaffold then the Cartan geometry is developable.
\end{theorem}
\begin{proof}
The injectivity radius of the Riemannian metric on the scaffold is locally constant.
So as we travel along a geodesic, it continues to exist forward or backward in time for a never diminishing amount of time, i.e. every geodesic is complete.
\end{proof}
\begin{corollary}
If the automorphisms of a Cartan geometry \(\G\to M\) act locally transitively on the total space \(\G\) then the Cartan geometry is developable.
\end{corollary}
We can generalize this slightly.
A topological space is \emph{pathwise compact}\define{pathwise compact} if every one of its path components lies in a compact set.
A group of homeomorphisms of a topological space acts \emph{pathwise cocompactly}\define{pathwise cocompact} if the quotient space (i.e. the space of group orbits) is pathwise compact in the quotient topology.
\begin{theorem}
If the automorphism group of a Cartan geometry acts pathwise cocompactly on the scaffold then the Cartan geometry is developable.
\end{theorem}
\begin{proof}
As above, the injectivity radius of the Riemannian metric on the scaffold is constant along orbits of the automorphism group.
In any Riemannian manifold, the injectivity radius of a point varies continuously with the choice of point, i.e. injectivity radius is a continuous function \cite{Klingenberg:1995} p. 131 Proposition 2.1.10.
Every continuous function invariant under automorphism descends to a continuous function on the quotient space.
(The quotient space is a manifold by theorem~\vref{thm:strong}, but we don't need to know that for this proof. 
In particular, pathwise cocompactness of the automorphism group action is precisely the compactness of all path components of the quotient space.)
Every geodesic is a continuous path, so maps to a continuous path in a single path component of the quotient space.
In some compact set containing that path component, the injectivity radius function has a minimum.
Since injectivity radius is positive, that minimum is positive.
So as we travel along a geodesic, the geodesic continues to exist forward or backward in time for a never diminishing amount of time, i.e. every geodesic is complete.
\end{proof}
\begin{example}
Some noncompact Riemannian manifolds have compact quotient spaces by automorphisms, so are developable.
\end{example}
{\centering\tiny
\begin{tabular}{p{5cm}p{5cm}}
\includegraphics[width=5cm]{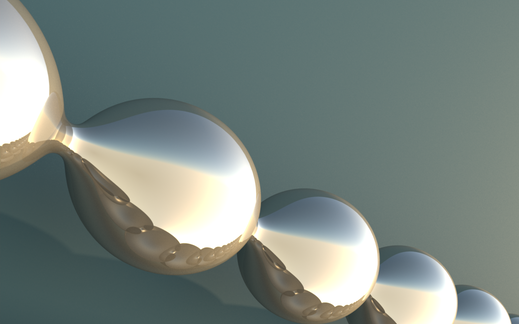} & \includegraphics[width=5cm]{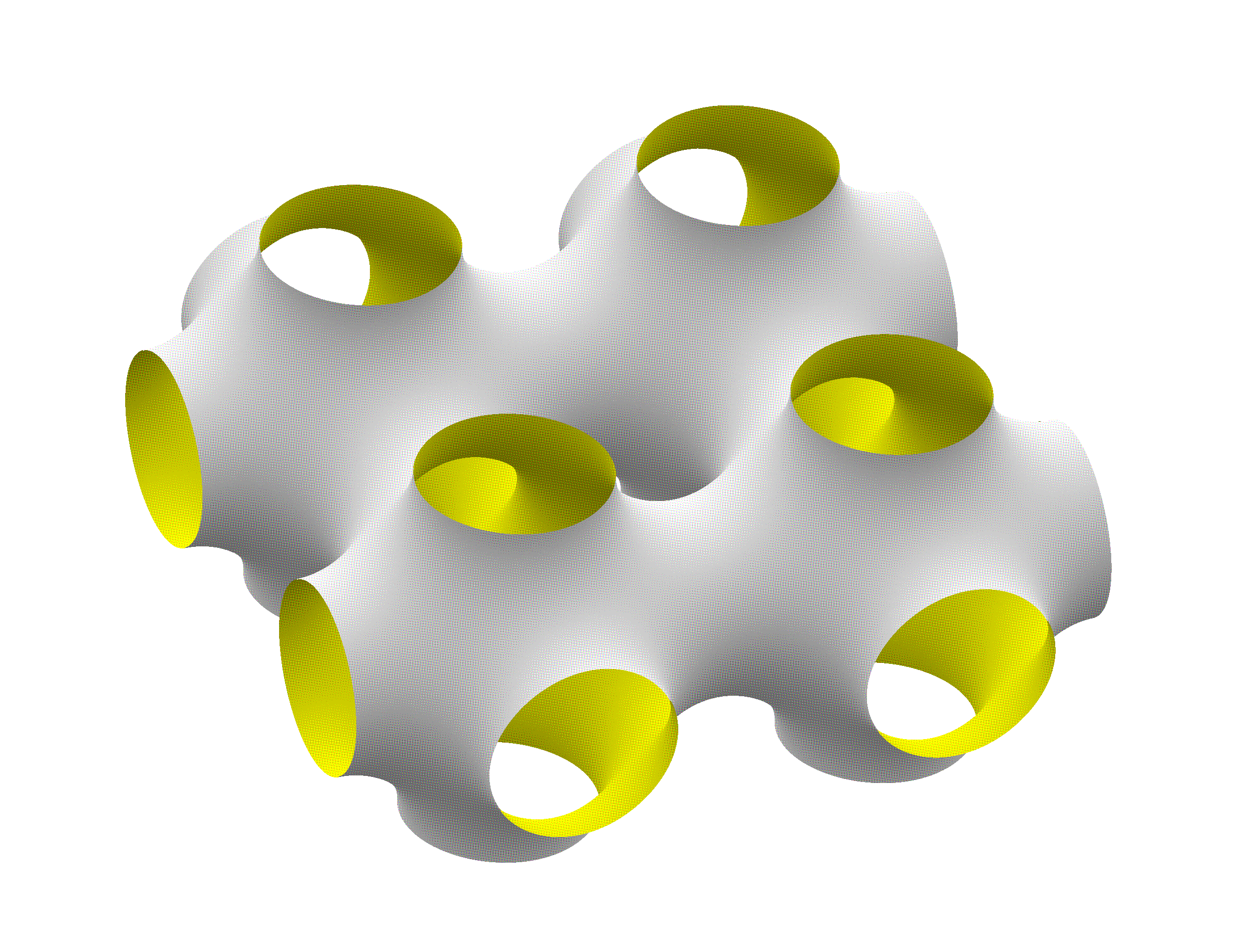} \\
GeometrieWerkstatt CC BY-NC-SA-3.0 &
By Anders Sandberg \\ & Own work \\ & CC BY-SA 3.0 \\
& {\verb!commons.wikimedia.org/w/index.php?curid=20737176!}
\end{tabular}
\par
}
\begin{example}
The automorphism groups of de~Sitter and anti-de~Sitter space act transitively on the orthonormal frame bundle, which is the total space of their Lorentz geometries, i.e. Cartan geometries modelled on Minkowski space.
Hence they are developable.
\end{example}
\section{Developing and the kernel}
If \(K\) is the kernel of a homogeneous space \((X,G)\), let \(\bar{G}:=G/K\), and denote the Lie algebra of \(K\) as \(\LieK\).
If \(\G\to M\) is an \((X,G)\)-geometry, let \(\bar\G:=\G/K\).
The \(1\)-form $\bar\omega:=\omega + \LieK$ is invariant under $K$-action, and vanishes on the fibers of $\G\to\bar\G$, so is defined on $\bar\G$, a Cartan connection.
With this Cartan connection, \(\bar\G\) is the \emph{associated effective}\define{associated!effective!Cartan geometry}\define{Cartan!geometry!associated effective} \((X,\bar{G})\)-geometry.
An \((X,\bar{G})\)-geometry lifts to an \((X,G)\)-geometry just when a certain cocycle in \(\cohomology{1}{M,\nForms{1}{M} \otimes (\amal{\G}{H}{\LieK})}\) vanishes, and it is then unique up to \(\cohomology{0}{M,\nForms{1}{M} \otimes (\amal{\G}{H}{\LieK})}\).
\begin{theorem}
A Cartan geometry is developable just when its associated effective Cartan geometry is developable.
\end{theorem}
The proof is as for theorem~\vref{thm:eff.complete}.
\begin{example}
A spin Lorentz manifold is developable just when its underlying Lorentz metric is developable.
\end{example}
\section{Developing submanifolds}
Take two Cartan geometries \(\G\to M\), \(\G'\to M'\) with the same model \((X,G)\), \(X=G/H\).
Take two immersions
\[
\begin{tikzcd}
M&S\arrow[l,"\iota"']\arrow[r,"\iota'"]&M'
\end{tikzcd}
\]
On the fiber bundle
\[
\begin{tikzcd}
\amal{\iota^*\G}{H}{\iota^{\prime *}\G'}
\arrow[d]\\
S
\end{tikzcd}
\]
we have the pullback \(V\) of the submersion connection above.
The immersions are \emph{developments}\define{development} of one another if this fiber bundle is foliated by leaves of \(V\), all projecting diffeomorphically to \(S\).
Note that, if both have the same constant curvature, then the equivariant Frobenius theorem (corollary~\vref{corollary:equivariant.Frobenius}) ensures that the fiber bundle is foliated, but does not ensure that the leaves project diffeomorphically to \(S\).
On the other hand, the 

A \emph{pancake}\define{pancake} is an immersion \(S\xrightarrow{\iota}M\) from a connected manifold \(S\) on which the curvature of the Cartan geometry vanishes, so 
\[
\begin{tikzcd}
\iota^*\G_G\arrow[d]&G\arrow[l]\\
S
\end{tikzcd}
\]
is a \(G\)-bundle with flat connection \(\omega_{\G}\).
In other words, a pancake is precisely a submanifold for which \(\iota^*\G_H\to S\) is a flat subgeometry.
%Hence by the equivariant Frobenius theorem (corollary~\vref{corollary:equivariant.Frobenius}), every pancake locally develops to the model.
\begin{theorem}\label{thm:develop.submanifold}
Take a Cartan geometry \(\G\to M\).
Every pancake \(S\xrightarrow{\iota}M\) has a development from its universal covering space
\[
\begin{tikzcd}
\tilde{S}\arrow[d,"\pi"]\arrow[r,"\delta"]&X\\
S
\end{tikzcd}
\]
called its \emph{developing map}\define{developing map} \(\delta\), equivariant for a unique group morphism \(\pi_1(S)\xrightarrow{h}G\), its \emph{holonomy morphism},\define{holonomy!morphism} so that the developing map lifts to a smooth immersion
\[
\begin{tikzcd}
(\iota\circ\pi)^*\G\arrow[r]\arrow[d]&\delta^*G\arrow[d]\\
\tilde{S}\arrow[r,"\iota\circ\pi"]&X
\end{tikzcd}
\]
equivariant for the holonomy morphism and the \(H\)-action, and matching \(\omega_{\G}\) and \(\omega_G\).
The pair \((\delta,h)\) of developing map and holonomy morphism are unique up to replacing them by \((g\delta,\Ad_g h)\).
\end{theorem}
\begin{proof}
Let \(\G_S:=\iota^*\G\).
Vanishing of the curvature on \(\G_S\) ensures that the set \(V\) of tangent vectors to \(\G_S\times G\) satisfying \(\omega_{\G}=\omega_G\) is a subbundle of the tangent bundle, invariant under diagonal right \(H\)-action by lemma~\vref{lemma:bracket}, and invariant under left \(G\)-action on \(G\), because \(\omega_G\) is left invariant.
Since the \(H\)-orbits are tangent to \(V\), and \(V\) is \(H\)-invariant, by the equivariant Frobenius theorem (corollary~\vref{corollary:equivariant.Frobenius}), \(\G_S\times G\) is foliated by \(H\)-folios of dimension equal to that of \(\G_S\).
Let \(\G'\) be one such \(H\)-folio.
Let \(\omega\) be the pullback of \(\omega_{\G}\) to \(\G'\); because \(\G'\) is \(V\)-tangent, \(\omega_{\G}=\omega_G\) on \(\G'\), i.e. the pullbacks agree.
By theorem~\vref{cor:orbitMapEquiv} the map \(\G'\to\G\) is a fiber bundle mapping as is the quotient \(S':=\G'/H\to S\).
Since each is a local diffeomorphism, each is a covering map.
Lift the bundle \(\G'\to S'\) to the universal covering space \(\tilde{S}\):
\[
\begin{tikzcd}
\G_S\arrow[d]&\tilde\G\arrow[d]\arrow[r]\arrow[l]&G\arrow[d]\\
S&\tilde{S}\arrow[l]\arrow[r]&X.
\end{tikzcd}
%\begin{tikzcd}
%&\tilde\G\arrow[dl]\arrow[dr]\arrow[rrr,bend left]&&&\tilde{S}\arrow[dl]\arrow[dr]&\\
%\G_S\arrow[rrr,bend right]&&G\arrow[rrr,bend right,crossing over]&S&&X.
%\end{tikzcd}
\]
The deck transformations of \(\tilde{S}\) act as bundle automorphisms on \(\tilde\G\), since the bundle is pulled back from \(S\).
These transformations preserve \(\omega\) and the \(H\)-action, as these are also pulled back.
But the map to \(G\) might not be invariant under the deck transformations.
Locally, each deck transformation \(\gamma\in\pi_1(S)\) acts by an automorphism, so some element \(g=h(\gamma)\in G\).
By connectivity of \(S\), this element is unique.

The construction is unique up to the choice of leaf \(\G'\).
If we pick the leaf through some point \((p_0,g_0)\), changing our choice to \((p_0,gg_0)\) alters the developing map and holonomy morphism as described.
But every leaf maps to \(p_0\), since every leaf covers \(\G\).
The leaves in \(\amal{\G_S}{H}{G}\) have preimages in \(\G_S\times G\) which are \(H\)-invariant, so \(H\)-bundles, on which \(\omega_{\G}=\omega_G\).
\end{proof}
\begin{corollary}
Take a connected manifold \(M\) with a locally homogeneous structure modelled on a strong effective homogeneous space \((X,G)\).
The locally homogeneous structure is pulled back by the developing map \(\tilde{M}\to X\), which is uniquely determined up to left \(G\)-translation.
\end{corollary}
\begin{example}
Take as model \((X,G)\) the unit sphere \(X=S^n\subset\R^{n+1}\) with its isometry group \(G=\Orth{n+1}\).
Every flat \((X,G)\)-geometry is precisely a Riemannian metric locally isometric to the sphere, but it also a pancake, i.e. by inclusion map into itself.
On a compact connected manifold \(M\), such a geometry is a complete Riemannian metric, so complete as a Cartan geometry, so has universal covering space \(\tilde{M}\) with developing map \(\tilde{M}\to X\), a complete Klein manifold, hence compact and connected, so \(\tilde{M}=X\).
Hence \(M=\Gamma\backslash X\) is a quotient of the sphere by a finite group, with \(\Gamma\subset\Orth{n+1}\) some finite group acting without fixed points.
Conversely, every finite subgroup \(\Gamma\subseteq\Orth{n+1}\) fixing no point of the unit sphere arises in this way uniquely; these are classified \cite{Allcock:2018}, \cite{Wolf:1967}.
\end{example}

\chapter{Completeness}
A Cartan geometry is \emph{complete}\define{complete!Cartan geometry}\define{Cartan!geometry!complete} if its constant vector fields are complete, i.e. have flows defined for all time.
\begin{example}
In the model, the constant vector fields are the generators of the \(1\)-parameter subgroups, so the model is complete.
\end{example}
\begin{problem}{ay}
Pulling back a Cartan geometry by a covering map preserves and reflects completeness.
\end{problem}
\begin{example}
The model geometry \((X,G)=(\RP{n},\PGL{n})\) lifts to the \(n\)-sphere and drops to every quotient \(\Gamma\backslash S^n\) by any finite group \(\Gamma\subset\SO{n+1}\) of rotations acting freely, remaining complete.
\end{example}
\[
\begin{array}{p{6cm}}
\includegraphics[width=6cm]{Layers_of_an_onion.jpg}\\
\tiny{Muffet, CC BY 2.0, via Wikimedia Commons}\\ 
\tiny{\verb!<https://creativecommons.org/licenses/by/2.0>!}
\end{array}
\]
\begin{example}
Take the unit sphere with its usual round metric, and puncture it at a finite set of at least two points; the \emph{onion}\SubIndex{onion} is the universal covering space of the punctured sphere, with the pullback Riemannian metric, hence the flat pullback Cartan geometry modelled on the sphere.
The onion is not complete, because the punctured sphere is not complete as a metric space, hence as a surface with Riemannian metric.
The metric space completion of the onion attaches one point for each puncture, not a smooth surface.
So this smooth flat Cartan geometry does not extend to a complete geometry.
It does not even embed into any larger Cartan geometry.
\end{example}
\begin{example}
Affine space \((X,G)=(\R^n,\GL{n}\ltimes\R^n)\) sits as an open set inside projective space \((X',G')=(\RP{n},\PGL{n})\), equivariantly for an obvious morphism \(G\to G'\).
This ensures that every manifold with an affine connection, i.e. an \((X,G)\)-geometry, has an induced \((X',G')\)-geometry, its \emph{projective connection}.\define{projective!connection!induced}
Even if the manifold is compact and the affine connection is complete (even if it is the Levi--Civita connection of a Riemannian metric) the projective connection may be incomplete.
For example, any flat torus has a flat complete homogeneous affine connection, but has no flat complete projective connection, as we will see.
\end{example}
\begin{problem}{az}
Take the plane \(\R^2\) with Lorentz metric \(dx \,dy\).
Let \(G\) be the group of affine Lorentz transformations of \(\R^2\) which preserve the line \(y=0\).
Let \(X\) be the complement of that line.
Check that \((X,G)\) is a homogeneous incomplete Lorentz surface
\end{problem}
\begin{example}
A \emph{Killing field}\define{Killing field} of a pseudo-Riemannian metric is a vector field whose flow preserves the metric.
It therefore arises from a unique infinitesimal automorphism of the Cartan geometry of the Riemannian metric, as the model of pseudo-Riemannian geometry has trivial kernel.
It is easy to check \cite{Arnold:1989} p. 88 that if \(v\) is a Killing field of \(g\) then \(g(v,\dot{x})\) is constant along any geodesic \(x(t)\).
On the plane  with coordinates \(x,y\), the Lorentz metric \(g=dx\,dy+\frac{1}{2}f(x)dy^2\) has Killing field \(\partial_{y}\) hence integral of motion \(I=\dot{x}+f(x)\dot{y}\) \cite{Guediri1995}.
Suppose that \(f(0)\ne 0\) but that \(f(x)=0\) for some values of \(x\), both positive and negative.
We want to see why this Lorentz metric on the torus is incomplete.
Factor the metric
\[
g=(dx+\frac{1}{2}f(x)dy)dy.
\]
So the null geodesics are the curves 
\[
\dot{x}+\frac{1}{2}f(x)\dot{y}=0
\]
together with the curves
\[
\dot{y}=0.
\]
The null geodesics of this second family, \(\dot{y}=0\), are just linear motion in \(x\) (since \(I\) is constant) with fixed \(y\).
Travel along any of the null geodesics of the first family:
\[
\dot{x}+\frac{1}{2}f(x)\dot{y}=0
\]
starting with some nonzero velocity \(\dot{y}\ne 0\) and at some point where \(f(x)\ne 0\).
Adding and subtracting \(g\) from \(I\) yields 
\[
\dot{x}=-I, \ f(x)\dot{y}=2I.
\]
At the start, \(f(x)\dot{y}\ne 0\) so \(I\ne 0\), and \(I\) stays constant.
Thus \(x(t)\) moves linearly, with nonzero velocity, and so \(f(x)\) reaches zero in finite time.
At time \(t=0\), \(f(x)\dot{y}=2I\ne 0\), and \(f(x)\dot{y}\) is constant while \(f(x)\to 0\), so \(\dot{y}\to\pm\infty\) in finite time.
Take \(f\) to be periodic with period \(1\) and the Lorentz metric descends to the torus \(M=\R^2/\Z^2\), hence an incomplete Lorentz metric on the torus, even with a circle symmetry group.
\end{example}
\begin{example}
The \emph{Clifton--Pohl torus}\define{Clifton--Pohl!torus} is the quotient of the punctured plane \(\R^2-0\) by \((x,y)\mapsto(\lambda x,\lambda y)\) for a fixed \(\lambda\ne\pm 1\) with the Lorentz metric
\[
\frac{dx \,dy}{x^2+y^2}
\]
Changing the sign of both \(x\) and \(y\) preserves the metric, so the fixed points of this transformation are totally geodesic, i.e. the axes are null geodesics.
The Clifton--Pohl torus is incomplete because the null geodesics \((\pm 1/(1-t),0)\) along this axis cannot be extended, and similarly along the other axis.
\end{example}
\Danger{} there are two inequivalent definitions of \emph{complete} Cartan geometry.
Sharpe's book \cite{Sharpe:1997} has become the standard reference on Cartan geometries, so I think it best to follow his definition, as we have done.
The older definition \cite{Ehresmann:1936,Ehresmann:1938,Ehresmann:1951} is what we have called \emph{developable}.\SubIndex{developable!Cartan geometry}\SubIndex{Cartan!geometry!developable}
Ehresmann \cite{Ehresmann:1938} also uses the term \emph{normal}\SubIndex{normal!Cartan geometry}\SubIndex{Cartan!geometry!normal} for what we call \emph{developable}; the term \emph{normal} today has a different meaning in standard use for a condition on the curvature of certain types of Cartan geometry \cite{Cap/Slovak:2009}.
Kobayashi's fundamental theorem \cite{Kobayashi:1954} Th\'eor\`eme 1 (which he attempted to prove in \cite{Kobayashi:1957} p. 172 (D)) states that these are equivalent.
Clifton showed that this theorem is wrong \cite{Clifton:1966}.
Clifton's paper seems to have gone unnoticed for many years.
Kobayashi's approach is to apply the Trotter product formula to constant vector fields, to try to approximate a time-varying vector field.
The mistake is to assume convergence.
\begin{theorem}%
[Clifton II \cite{Clifton:1966}]%
\label{thm:dev.implies.cmplt}%
\define{theorem!Clifton II}%
\define{Clifton!theorem II}
Developable implies complete.
\end{theorem}
\begin{proof}
Developable is equivalent to completeness of the metric on \(\G\) by Clifton's theorem~\vref{theorem:Clifton}.
The constant vector fields are of constant length in that metric.
As we move for a time \(<t\) along a flow of a constant vector field \(A_{\G}\), we stay inside the closed ball of radius \(t|A|\).
This ball is compact, since the metric is complete.
In any compact set, we have a positive lower bound on how much longer the flow continues.
And so the flow continues: the Cartan geometry is complete.
\end{proof}
\begin{problem}{ba}
Do complete Cartan geometries form a closed set among Cartan geometries? An open set? (In what topology?)
\end{problem}
\begin{problem}{bb}
Do developable Cartan geometries form a closed set among Cartan geometries? An open set? (In what topology?)
\end{problem}
\begin{example}
Pick an invertible linear isomorphism \(V\xrightarrow{g}V\) of a finite dimensional vector space \(V\), with all eigenvalues in the unit disk. 
The quotient of \(M:=(V-0)/(x\sim gx)\) is the Hopf manifold, with its flat affine structure.
The structure is incomplete because, for any \(x_0\in V-0\), the geodesic \(t\mapsto tx_0\) maps to an incomplete geodesic in \(M\).
But \(M\) is compact, so the geodesic remains in the compact set \(M\): it doesn't wander off to infinity in finite time.
For simplicity, suppose that \(x_0\) is an eigenvector of \(g\), say with eigenvalue \(\lambda>0\).
So \(t\mapsto tx_0\) is a geodesic loop which passes through the image point of \(x_0\) at times \(t=\lambda^n\), \(n\in\Z\).
In the picture up in \(V\), as \(t\) goes from \(1\) to \(\lambda^{-n}\), \(tx_0\) moves to a different sheet of \(V-0\to M\).
We can move that sheet back to the other by applying \(g^n\), which affects our velocity by making it \(\lambda^{-n}x_0\) as measured on the original sheet.
So we travel along our geodesic loop more and more slowly as time increases.
In the other direction, as \(t\to 0\) we travel around our loop faster and faster.
The velocities with which this geodesic loop passes through the image point of \(x_0\) in \(M\) are unbounded, growing to infinity in finite time.
\end{example}
Recall that the projection to \(M\) of a flow line in \(\G\) of a constant vector field of a Cartan geometry \(\G\to M\) is a \emph{spiral}.\SubIndex{spiral}
\begin{theorem}\label{theorem:spirals}
Take a Cartan geometry \(\G\to M\).
Take a constant vector field in that geometry.
Consider its flow line through some point of \(\G\).
Suppose that \(I\subseteq\R\) is the largest interval on which that flow line is defined.
Consider the spiral, i.e. the image in \(M\) of the flow line.
Then either 
\begin{itemize}
\item
\(I=\R\), the flow line is complete or
\item
the spiral does not extend continously to \(\bar{I}\) or 
\item
the spiral has velocity leaving every compact set in \(TM\) and the velocity is not locally integrable on any open subinterval \(J\subseteq I\) for which \(\bar{J}\) shares an endpoint with \(\bar{I}\).
\end{itemize}
\end{theorem}
\begin{proof}
Suppose that the model is \((X,G)\) with \(X=G/H\).
Suppose that the flow line is defined on \(I\subsetneq\R\).
For \(t\) near one end of \(I\), say \(t\) near \(T\), continuously extend the projection \(m(t)\in M\) of the flow line to \(t=T\).
Near \(m(T)\), take a local trivialization of \(\G\to M\).
In our trivialization, writing points as \(p=(m,h)\).
Write out the Cartan connection as
\[
\omega_{\G}=\omega_H+\Ad_h^{-1}\eta
\]
where \(\eta\) is a \(\LieG\)-valued \(1\)-form on \(M\) defined near \(m(T)\), as in section~\vref{section:Cartan.geometries}.
The differential equation of the flow of the constant vector field is
\[
A=\dot{h}\hook\omega_H+\Ad_h^{-1}(\dot{m}\hook\eta).
\]
In the group \(G\), let
\[
g:=e^{tA}h^{-1}.
\]
Claim: \(\dot{g}\hook\omega_G=\dot{m}\hook\eta\).
If this is true, then \(g\) satisfies a smooth Lie equation, so has solution for time \(t<T\); see appendix~\vref{appendix:Lie.equations}.
So \(g\) extends locally absolutely continuously to \(\bar{I}\) just when \(\dot{m}\hook\eta\) remains locally integrable, and then \(h=e^{tA}g^{-1}\) also extends locally absolutely continuously to \(\bar{I}\).
By existence and uniqueness of constant vector field flows, the flow of the constant vector field through the limit point in \(\G\) extends our flow line, and so the flow extends to a larger open interval.
So \(\dot{m}(t)\) is locally integrable on \(0\le t<T\) just when smooth beyond time \(t=T\).

We recall that
\[
\frac{d}{dt}h^{-1}=-\RT{h*}^{-1}\LT{h*}^{-1}\dot{h}.
\]
Hence we differentiate
\begin{align*}
\dot{g}
&=
\RT{h*}^{-1}\LT{e^{tA}*}A
+
\LT{e^{tA}*}\frac{d}{dt}h^{-1},
\\
&=
\RT{h*}^{-1}\LT{e^{tA}*}(\dot{h}\hook\omega_H+\Ad_h^{-1}(\dot{m}\hook\eta))
+
\LT{e^{tA}*}(-\RT{h*}^{-1}\LT{h*}^{-1}\dot{h}),
\\
&=
\RT{h*}^{-1}\LT{e^{tA}*}\left(
\LT{h*}^{-1}\dot{h}
+\LT{h*}^{-1}
\RT{h*}
(\dot{m}\hook\eta)
\right)
+
\LT{e^{tA}*}(-\RT{h*}^{-1}\LT{h*}^{-1}\dot{h}),
\\
&=
\LT{g*}(\dot{m}\hook\eta),
\end{align*}
hence our Lie equation.
\end{proof}
\begin{theorem}\label{thm:eff.complete}
A Cartan geometry is complete just when its associated effective Cartan geometry is complete.
\end{theorem}
\begin{proof}
They have the same spirals.
\end{proof}
\begin{problem}{bc}
If a Riemannian manifold has induced projective connection complete, is it a complete Riemannian manifold?
\end{problem}
\section{Example: complete undevelopable}
The \emph{Clifton can}\define{Clifton!can} is an example of a complete Cartan geometry which is not developable \cite{Clifton:1966}.
Take the half cylinder \(M\subseteq\R^3\) cut out by \(x^2+y^2=1\), \(z>0\).
\[
\includegraphics{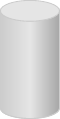}
\]
In cylindrical coordinates, 
\[
x=r\cos\theta,y=r\sin\theta,z
\]
take the usual basis of vector fields \(\partial_z,\partial_{\theta}\).
\[
\includegraphics{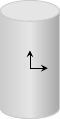}
\]
Pick a smooth real-valued function \(\varphi(z)\) with \(|\varphi(z)|\to\infty\) as \(z\to 0\).
For example, Clifton takes
\[
\varphi(z)=\frac{1}{z}.
\]
At each point of the cylinder, consider the two vector fields given by rotating the usual basis by angle \(\varphi(z)\).
\[
\includegraphics{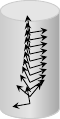}
\]
The angle turns rapidly as we approach the bottom of the cylinder, and the vector fields do not approach a limit.
Pick either of the two vector fields.
That vector field is a unit vector field, so every one of its flow lines either survives for all time or eventually approaches \(z=0\); we want to see why it survives for all time.
Our vector field is invariant under rotation of the cylinder.
Its angle to the vertical, \(\varphi(z)\), turns faster and faster as \(z\to 0\).
So our vector field becomes horizontal for infinitely many heights \(z\), arbitrarily close to \(z=0\).
At each of these heights \(z\), our vector field is tangent to a circle \(z=\text{constant}\).
If our flow lines heads down toward \(z=0\), eventually it hits one of these circles; but they are invariant, so it never hits any of them or is always on one of them.
Hence every flow line is bounded below by the next invariant circle, so the vector fields are complete.

Let \((X,G)=(\R^2,\R^2)\), so an \((X,G)\)-structure is exactly a pair of vector fields, and we have a Cartan geometry on the half cylinder.
Completeness of the Cartan geometry is precisely the completeness of every constant coefficient combination of the two vector fields; such combinations are complete by the same argument.
The scaffold has the usual metric of the half cylinder, since our vector fields are orthonormal for that metric.
So the scaffold is not a complete metric space.
So this Cartan geometry on the half cylinder is \emph{not} developable: 
The pullback of this Cartan geometry to the plane, as universal covering space of the cylinder, is also complete but not developable.
In this argument, we can replace the circle with any compact manifold \(S\) with trivial tangent bundle, to get undevelopable complete Cartan geometries on products \(M:=\R^+\times S\) and on their covering spaces.

\chapter{Inextendable geometries}
\section{The complete cone}
Take a Cartan geometry \(\G\to M\) with model \((X,G)\).
Let \(\CompleteCone=\CompleteCone[M]\subseteq\LieG\) be the set of \(A\in\LieG\) for which the associated constant vector field \(A_{\G}\) is complete.
In particular, a Cartan geometry is complete just when \(\CompleteCone=\LieG\).
\begin{example}
Take \((X,G)=(\R^2,\R^2)\) the plane under translation.
Then an \((X,G)\)-geometry is precisely a trivial bundle \(\G=M\to M\), with a Cartan connection being precisely a \(1\)-form \(\omega\) valued in \(\R^2\), i.e. a pair of \(1\)-forms, linearly independent everywhere.
So consider
\begin{align*}
\omega^1&=\frac{dx}{1+y^2},\\
\omega^2&=\frac{dy}{1+x^2}.
\end{align*}
The constant vector fields associated to \((a,b)\in\LieG\) are
\[
(a,b)_{\G}=a(1+y^2)\partial_x+b(1+x^2)\partial_y.
\]
For instance 
\begin{align*}
(1,0)_{\G}&=(1+y^2)\partial_x,\\
(0,1)_{\G}&=(1+x^2)\partial_y,\\
\end{align*}
are complete, but
\[
(1,1)_{\G}=(1+y^2)\partial_x+(1+x^2)\partial_y
\]
is not, having flow through the origin 
\[
(x(t),y(t))=(\tan t,\tan t).
\]
So \(\CompleteCone\) is not necessarily closed under addition.
\end{example}
\begin{example}
The model translation structure on Euclidean space becomes incomplete on each solid cylinder, with \(\CompleteCone\) a one dimensional vector subspace.
\end{example}
\begin{example}
On the plane, take a ``lump'': an open set containing exactly one line, the line \(y=0\), and intersecting all other lines in a relatively compact set.
\[
\begin{tikzpicture}[domain=-4:4]
%  \draw[very thin,color=gray] (-0.1,-1.1) grid (3.9,3.9);
%
%  \draw[->] (-0.2,0) -- (4.2,0) node[right] {$x$};
%  \draw[->] (0,-1.2) -- (0,4.2) node[above] {$f(x)$};
%
%  \draw[color=red]    plot (\x,\x)             node[right] {$f(x) =x$};
%  % \x r means to convert '\x' from degrees to _r_adians:
%  \draw[color=blue]   plot (\x,{sin(\x r)})    node[right] {$f(x) = \sin x$};
\draw[->] (-4,0) -- (4.1,0);
  \fill[color=blue!20,draw=blue,samples=100] plot (\x,{exp(-\x*\x)});
  \fill[color=blue!20,draw=blue,samples=100] plot (\x,{-exp(-\x*\x)});
  % node[right] {$f(x) = \frac{1}{20} \mathrm e^x$};
\end{tikzpicture}
\]
Take the translation structure on the plane with Cartan connection
\[
\omega=\frac{(dx,dy)}{1+r^2h},
\]
where \(r^2=x^2+y^2\) and \(h=h(x,y)\ge 0\) is equal to one near \(y=0\) and equal to zero outside the lump.
The constant vector fields are
\[
v_{\G}=(1+r^2h)v,
\]
for each vector \(v=(a,b)\in\R^2\).
The flow of this vector field agrees with that of a linear vector field, once our point moves outside of the lump.
It points in the same direction as that linear vector field, so follows along a straight line.
If \(b\ne 0\), our vector field disagrees with the linear vector field only on a compact set on that line, where it is therefore bounded from above and below, so complete.
Since the vector field flows at faster than constant rate, every flow line leaves the lump in finite time, unless \(b=0\).
If \((a,b)=(0,0)\), we stay put, so complete flow.
But if \(a\ne 0\) and \(b=0\) then the flow is incomplete.
So \(\CompleteCone=\set{(a,b)|a\ne 0\text{ or }(a,b)=(0,0)}\):
\[
\begin{tikzpicture}
\fill[gray!50] (-1,.05) rectangle (1,1);
\fill[gray!50] (-1,-.05) rectangle (1,-1);
\draw[dashed] (-1,0) -- (1,0);
\fill[gray!50] (-.05,-.05) rectangle (.05,.05);
\draw[gray] (0,-1) -- (0,1);
\end{tikzpicture}
\]
\end{example}
\begin{theorem}\label{thm:complete.cone.brackets}
Take a Cartan geometry \(\G\to M\).
The complete cone \(\CompleteCone=\CompleteCone[M]\) is a cone, invariant under adjoint \(H\)-action, closed under brackets with elements of \(\LieH\), and containing \(\LieH\).
Suppose that \(\LieG'\subseteq\LieG\) is a linear subspace on which the curvature of \(\G\) is constant, forming a Lie algebra under curvature deformed bracket, and that, as a Lie algebra, \(\LieG'\) is generated by elements of \(\CompleteCone\).
Then \(\LieG'\subseteq\CompleteCone\).
\end{theorem}
\begin{proof}
Take a Cartan geometry \(\G\to M\) with model \((X,G)\), \(X=G/H\).
By definition it is clear that the complete cone is a cone, i.e. contains the origin and  is closed under rescaling by real constants.
The vector fields \(A_{\G}\) for \(A\in\LieH\) generate the \(H\)-action, so have complete flows.
By definition of constant vector fields, if \(A\in\LieG\),
\[
\RT{h*}(A_{\G})=(\Ad_h^{-1}A)_{\G}
\]
so that adjoint \(H\)-action preserves the complete cone.
Let \(h_s:=e^{sB}\), pick a point \(p_0\in\G\) and let
\[
\varphi(s,t):=e^{t\Ad_{h_s}^{-1} A_{\G}}=e^{\RT{h_s*}tA_{\G}}p_0,
\]
which is clearly defined for all \(s,t\).
Differentiating in \(s\), 
\[
\lb{\LieH}{\CompleteCone}\subset\CompleteCone.
\]
The Lie bracket of elements of \(\LieG'\) is identified by \(A\mapsto A_{\G}\) with the Lie bracket of the vector fields \(A_{\G}\), so the Jacobi identity for the curvature deformed bracket is the Jacobi identity for vector fields.
Apply Palais's theorem (theorem~\vref{theorem:Palais}) to construct a Lie group with Lie algebra \(\LieG'\) acting on \(\G\).
\end{proof}
\begin{lemma}
Take a Cartan geometry \(\G\to M\) with model \((X,G)\) and stabilizer \(H:=G^{x_0}\).
The complete cone is invariant under the expected structure group and under Lie brackets with the expected structure algebra.
\end{lemma}
\begin{proof}
For \((h,\bar{h})\in H'\),
\[
h^{-1*}\omega=\Ad_{\bar{h}}^{-1}\omega,
\]
so
\[
h^{-1}_*A_{\G}=(\Ad_{\bar{h}}^{-1}A)_{\G},
\]
with flows intertwined by \(h\).
\end{proof}
\section{Example: schichttorten}
Suppose that \(M\) is an oriented surface with Riemannian metric and that \(v\) is a nowhere vanishing Killing field, i.e. a vector field whose flow preserves the Riemannian metric.
The flow also preserves the orientation, by connectivity.
Moreover, it preserves the vector field: every vector field is invariant under its own flow.
There is a unique positively oriented orthonormal coframing \(\xi_1,\xi_2\) with \(\xi_1(v)=0\) and \(\xi_2(v)>0\).
Since the coframing is unique, it is invariant under the flow of \(v\):
\begin{align*}
0
&=
\LieDer_v \xi_1,
\\
&=
v\hook d\xi_1+d(v\hook\xi_1),
\\
&=
v\hook d\xi_1
\end{align*}
so \(d\xi_1=0\), locally \(\xi_1\) is exact; let us suppose that \(\xi_1\) is globally exact, say \(\xi_1=dz\).
Since \(v\hook dz=0\), \(z\) is constant on the flow lines of \(v\).

Let \(\lambda:=\log\xi_2(v)\); since \(\xi_1,\xi_2\) are invariants of \(v\) and the metric,  \(\lambda\) is also invariant, hence constant along the flow lines of \(v\), so \(d\lambda=\lambda'\,dz\) for some function \(\lambda'\).
Locally \(\lambda=\lambda(z)\); let us suppose that \(\lambda=\lambda(z)\) globally.

Since \(d\xi_2\) is a \(2\)-form on a surface, it is a multiple of \(dz\wedge \xi_2\), say write \(d\xi_2=f\,dz\wedge\xi_2\).
By \(v\)-invariance,
\begin{align*}
0
&=
\LieDer_v \xi_2,
\\
&=
v\hook d\xi_2+d(v\hook\xi_2),
\\
&=
v\hook f\,dz\wedge\xi_2+d(v\hook\xi_2),
\\
&=
-fe^{\lambda}\, dz+de^{\lambda},
\\
&=
-fe^{\lambda}\, dz+e^{\lambda}\lambda'\,dz.
\end{align*}
So \(f=\lambda'\), i.e.
\[
d\xi_2=\lambda'\,dz\wedge\xi_2.
\]
Hence
\begin{align*}
d(e^{-\lambda}\xi_2)
&=-\lambda'e^{-\lambda}dz\wedge\xi_2+e^{-\lambda}\lambda'dz\wedge\xi_2,
\\
&=0.
\end{align*}
So locally \(e^{-\lambda}\xi_2\) is exact:
\[
\xi_2=e^{\lambda}d\theta,
\]
for some local variable \(\theta\).
Summing up, so far:
\[
\xi_1=dz, \xi_2=e^{\lambda}d\theta, \lambda=\lambda(z).
\]

In these \(z,\theta\) coordinates, 
\[
0=v\hook\xi_1=v\hook dz, e^{\lambda}=v\hook \xi_2=v\hook e^{\lambda}d\theta,
\]
so
\[
v=\partial_{\theta}.
\]
Of course all constant multiples of \(v\) are also Killing fields, i.e. infinitesimal automorphisms of the oriented Riemannian geometry.
We can extend our surface to allow \(\theta\) to take values on the whole real number line.
We can then quotient our surface by \(\theta\sim\theta+2\pi\), so our surface is now locally isometric to a surface of revolution in \(3\)-dimensional Euclidean space.

The soldering form of the Cartan geometry of the oriented Riemannian geometry is
\[
\sigma=e^{i\varphi}\xi,
\]
where \(\xi=\xi_1+i\xi_2\) and \(\varphi\) represents an angle of rotation of an orthonormal frame from \(\xi\).
We want to use the torsion-free Cartan geometry given by the Levi--Civita connection, i.e. we want
\[
d\sigma=i\gamma\wedge\sigma.
\]
\begin{problem}{schichttorte.Levi}
Compute the Cartan connection in terms of \(\varphi,\lambda,z,d\varphi,dz,d\theta\).
\end{problem}
\begin{answer}{schichttorte.Levi}
If we let \(\xi=\xi_1+i\xi_2=dz+ie^{\lambda}d\theta\),
\[
d\xi=i\lambda'e^{\lambda}dz\wedge d\theta=i\lambda'\xi_1\wedge\xi_2.
\]
But
\[
\xi\wedge\bar\xi=(\xi_1+i\xi_2)\wedge(\xi_1-i\xi_2)=-2i\xi_1\wedge\xi_2,
\]
so
\[
\xi_1\wedge\xi_2=\frac{i}{2}\xi\wedge\bar\xi.
\]
Hence
\[
d\xi=-\frac{\lambda'}{2}\xi\wedge\bar\xi.
\]
Differentiating,
\begin{align*}
0
&=
d\sigma-i\gamma\wedge\sigma,
\\
&=
d(e^{i\varphi}\xi)-i\gamma\wedge\sigma,
\\
&=
i\, d\varphi\wedge\sigma+e^{i\varphi}d\xi-i\gamma\wedge\sigma,
\\
&=
i\, (d\varphi-\gamma)\wedge\sigma+e^{i\varphi}(-\frac{\lambda'}{2}\xi\wedge\bar\xi),
\\
&=
i\, (d\varphi-\gamma)\wedge\sigma-\frac{\lambda'}{2}\sigma\wedge\bar\xi,
\\
&=
i\, (d\varphi-\gamma)\wedge\sigma+\frac{\lambda'}{2}\bar\xi\wedge\sigma,
\\
&=
i\, \left(d\varphi-\gamma-i\frac{\lambda'}{2}\bar\xi\right)\wedge\sigma,
\\
&=
i\, \left(d\varphi-\gamma-i\frac{\lambda'}{2}\bar\xi
+i\frac{\lambda'}{2}\xi
\right)\wedge\sigma,
\end{align*}
hence
\[
0
=
\left(d\varphi-\gamma-i\frac{\lambda'}{2}\bar\xi
+i\frac{\lambda'}{2}\xi
\right)\wedge\xi.
\]
\begin{align*}
-i\frac{\lambda'}{2}\bar\xi
+i\frac{\lambda'}{2}\xi
&=
-i\frac{\lambda'}{2}(dz-i\xi_2)
+i\frac{\lambda'}{2}(dz+i\xi_2),
\\
&=
-i\frac{\lambda'}{2}(-i\xi_2)
+i\frac{\lambda'}{2}(+i\xi_2),
\\
&=
-2\frac{\lambda'}{2}\xi_2,
\\
&=
-\lambda'\xi_2.
\end{align*}
So we get
\[
\left(d\varphi-\gamma-\lambda'\xi_2\right)\wedge\xi=0,
\]
hence
\[
\gamma=d\varphi-\lambda'e^{\lambda}d\theta,
\]
So the soldering form is
\[
\sigma=e^{i\varphi}(dz+ie^{\lambda(z)}d\theta),
\]
and the Levi--Civita connection is
\[
\gamma=d\varphi-\lambda'(z)e^{\lambda(z)}d\theta.
\]
\end{answer}
Hence the Cartan connection is
\[
\omega=
\begin{pmatrix}
i\gamma&\sigma\\
0&0
\end{pmatrix}
\]
a \(1\)-form valued in the Lie algebra of the symmetry group \(G=\Un{1}\ltimes\C\) of the model \(X=\C\).
\begin{problem}{bd}
Compute the curvature.
\end{problem}
The equations of motion along constant vector fields are
\begin{align*}
e^{i\varphi}(dz+ie^{\lambda(z)}d\theta)&=(a_0+ib_0)\,ds,\\
d\varphi-\lambda'(z)e^{\lambda(z)}d\theta&=c_0\,ds,
\end{align*}
with each choice of real constants \(a_0,b_0,c_0\) determining a single constant vector field.

Solving for the differentials,
\begin{align*}
dz+ie^{\lambda}d\theta
&=
e^{-i\varphi}(a_0+ib_0)ds,
\\
&=
(\cos\varphi-i\sin\varphi)(a_0+ib_0)ds,
\\
&=
(a_0\cos\varphi+b_0\sin\varphi)ds
+
i(b_0\cos\varphi-a_0\sin\varphi)ds,
\end{align*}
so
\begin{align*}
dz&=(a_0\cos\varphi+b_0\sin\varphi)ds,\\
d\theta&=e^{-\lambda}(b_0\cos\varphi-a_0\sin\varphi)ds.
\end{align*}
We can always rotate and rescale to arrange that \(a_0+ib_0=1\) or \(0\), so consider \(a_0=1\):
\begin{align*}
dz&=\cos\varphi\,ds,\\
d\theta&=-e^{-\lambda}\sin\varphi\,ds,\\
d\varphi&=(-\lambda'\sin\varphi+c_0)ds.
\end{align*}

The geodesics are precisely those with \(c_0=0\).
For example, we see the \emph{meridians}\SubIndex{meridian} are the geodesics with \(\theta=\theta_0\) constant, \(\varphi=0\) and \(c_0=0\).
On any unit speed curve \(C\) between points with coordinates \((z_0,\theta_0)\) and \((z_1,\theta_1)\), 
\[
z_1-z_0=\int dz\le \int \sqrt{|\xi_1|^2+|\xi_2|^2}=\operatorname{length}_C,
\]
with equality just when \(\xi_2=0\), i.e. just on meridians: among all continuous paths between any two circles \(z=z_0\) and \(z=z_1\), those of minimal length are precisely the meridians.
In particular, each meridian is the shortest path between any of its points.
The \emph{latitudes}\SubIndex{latitude} are the compact connected curves with \(z=z_0\) constant and \(\varphi=\pi/2\) or \(\varphi=3\pi/2\), so \(\cos\varphi=0\), \(\sin\varphi=\pm 1\), and \(c_0=\mp\lambda'(z_0)\).
All of the latitudes are spirals, but a latitude is a geodesic precisely when \(\lambda'(z_0)=0\).

Consider the Clairault integral
\[
I:=e^{\lambda}\sin\varphi.
\]
Along the flow of our constant vector field,
\begin{align*}
dI
&=
\lambda'e^{\lambda}\sin\varphi\,dz
+
e^{\lambda}\cos\varphi\,d\varphi,
\\
&=
\lambda'e^{\lambda}\sin\varphi\cos\varphi\,ds
-
e^{\lambda}\cos\varphi(\lambda'\sin\varphi+c_0)ds,
\\
&=
-c_0e^{\lambda}\cos\varphi\,ds,
\\
&=
-c_0e^{\lambda}dz.
\end{align*}
For any \(t_0,t_1\) during which the flow line is defined,
\[
I(t_1)-I(t_0)=e^{\lambda(z_1)}\sin\varphi(t_1)-e^{\lambda(z_0)}\sin\varphi(t_0)
=
-c_0\int_{z_0}^{z_c} e^{\lambda}dz,
\]
where \(z_0=z(t_0)\), so
\[
\sin\varphi(t_1)=e^{\lambda(z_0)-\lambda(z_1)}\sin\varphi(t_0)
-c_0e^{-\lambda(z_1)}\int_{z_0}^{z_c} e^{\lambda}dz.
\]
For any fixed \(z_0,z_1\), by picking \(c_0\) large enough in absolute value, we can ensure that the right hand side is larger in absolute value than \(1\), for any initial value \(\varphi(t_0)\).
So we can ensure that our constant vector fields, starting at \(z_0\), never reach \(z_1\).
So there are complete spirals on any surface of revolution.
\[
\includegraphics[width=5cm]{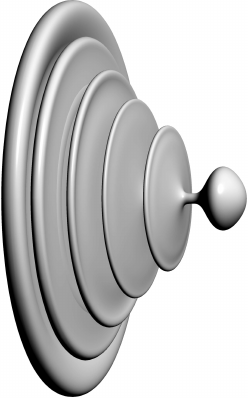}
\]
A \emph{schichttorte}\define{schichttorte} is 
\begin{itemize}
\item
a closed disk \(\bar{D}\) with interior \(D\) and
\item
an orientation and 
\item 
a smooth Riemannian metric on \(D\) and
\item
a vector field \(v\) on \(\bar{D}\), nonzero near the boundary, tangent to the boundary, with flow periodic of period \(2\pi\) through every point near the boundary and
\item
\(v\) is a Killing field near the boundary and
\item
the orthonormal frame \(dz,e^{\lambda}d\theta\) associated to \(v\) has 
\item
\(\lambda(z)\) a smooth function which has large negative values, say \(\lambda(z_n)\to-\infty\), at a sequence of points which approach the boundary and
\item
\(\lambda(z_{2n})=\lambda(z_{2n+1})\) for each integer \(n\) and
\item
\(\lambda\) is very large positive inbetween these \(z_n\), i.e.
\[
\int_{z_{2n}}^{z_{2n+1}}e^\lambda \,dz\to\infty.
\]
\end{itemize}
\begin{theorem}
Schichttorten exist.
Every schichttorte has all spirals defined for all time, except meridians.
Hence, as a Cartan geometry modelled on the oriented Euclidean plane, the complete cone of any schichttorte is precisely the set of matrices
\[
\begin{pmatrix}
ic_0 & a_0+ib_0 \\
0 & 0
\end{pmatrix}
\in \Un{1}\ltimes\C,
\]
so that either (i) \(c_0\ne 0\) or (ii) \(a_0=b_0=c_0=0\).
In particular, every schichttorte is incomplete and geodesically incomplete.
\end{theorem}
\begin{proof}
We can suppose that \(\bar{D}\) is the unit disk and that \(v=\partial_{\theta}\) near the boundary.
We can take \(z=\sqrt{x^2+y^2}\).
We can easily pick some function \(\lambda=\lambda(z)\) with the required properties on some interval \(\delta\le z<1\), and then extend the metric smoothly inside \(|z|\le\delta\) to an arbitrary Riemannian metric.
If there is a constant vector field with a flow line, i.e. a spiral, which reaches all of these \(z_n\), say at increasing times \(t_n\), we find
\[
\sin\varphi(t_{2n+1})-\sin\varphi(t_{2n})
=
-c_0e^{-\lambda(t_{2n})}
\int_{z_{2n}}^{z_{2n+1}}e^{\lambda}dz.
\]
For any \(c_0\ne 0\), the right hand side eventually gets larger, in absolute value, than \(2\), a contradiction.
So any spiral which is not a geodesic does not get outside some compact set \(z\le z_n\), and moves at bounded velocity inside that compact set in our metric, and hence in the Euclidean metric (since Riemannian metrics on compact sets are comparable).
By theorem~\vref{theorem:spirals}, these spirals continue to be defined for all time.

If \(c_0=0\) then for any circles \(z=z_0,z=z_1\),
\[
e^{\lambda(z_1)}\ge e^{\lambda(z_1)}\sin\varphi(z_1)=I(z_1)=I(z_0)=
e^{\lambda(z_0)}\sin\varphi(z_0),
\]
so no geodesic with \(\sin\varphi(z_0)\ne 0\), i.e. not a meridian, can ever enter a region where the circle \(z=z_1\) has \(e^{\lambda(z_1)}\) smaller than \(e^{\lambda(z_0)}\sin\varphi(z_0)\). 
But the values of \(e^{\lambda}\) get arbitrarily small as we go to the boundary, so the nonmeridian geodesics stay inside a compact annulus in our surface: they are complete.
Hence none of the spirals, except the meridians, can reach every \(z_n\).
Therefore they stay in some compact set, away from the singularity, so they are complete.
\end{proof}
\begin{problem}{be}
Prove that, for every compact surface with boundary \(\bar{M}\), every component of which has nonempty boundary, the surface \(M:=\bar{M}-\partial\bar{M}\) has a Riemannian metric whose associated Cartan geometry has the same complete cone as a schichttorte.
Note: the surface might not be oriented, or perhaps not even orientable, so the model for such a Riemannian geometry is \((X,G)=(\R^2,\Orth{2}\ltimes\R^2)\).
Hint: you might need the collaring theorem \cite{Connelly:1971}, \cite{Gauld:2014} section 3.2 pp. 42--44.
\end{problem}
\section{Local isomorphisms}
A \emph{local isomorphism}\define{local isomorphism!of Cartan geometries}\define{Cartan!geometry!local isomorphism} of two Cartan geometries with the same model is a diagram
\[
\begin{tikzcd}
\G\arrow[r,"\Phi"]\arrow[d]&\G'\arrow[d]\\
M\arrow[r,"\varphi"]&M'.
\end{tikzcd}
\]
of \(H\)-equivariant maps matching Cartan connections.
\begin{problem}{bf}
Prove that both \(\Phi\) and \(\varphi\) are local diffeomorphisms.
\end{problem}
\begin{problem}{bg}
What are all of the affine connections on the real line, up to diffeomorphism?
What are all of the affine connections on the circle, up to diffeomorphism?
Which are complete?
Which are developable?
\end{problem}
\begin{lemma}
Every local isomorphism of real analytic Cartan geometries is real analytic.
\end{lemma}
\begin{proof}
The flows of constant vector fields yield local coordinates: their flow charts, as in the Orbit Theorem (theorem~\vref{thm:orbit}).
\end{proof}
\section{Vast geometries}
For any Cartan geometry, denote its curvature as \(k\).
The \emph{complete module}, denoted \(\barCompleteCone\subseteq\LieG\), is the smallest \(H\)-submodule containing the complete cone \(\CompleteCone\) and so that, if \(A,B\in\barCompleteCone\) and \(k(A,B)\) is constant, then \(\lb{A}{B}-k(A,B)\in\barCompleteCone\).
A Cartan geometry  is \emph{vast}\define{vast} if \(\barCompleteCone=\LieG\).
\begin{example}
Every schichttorte is vast, but incomplete.
\end{example}
We will prove the implications
\[
\begin{tikzpicture}
\node (dev) at (0,3) {developable};
\node (cpl) at (0,2) {complete};
\node (vst) at (0,1) {vast};
\node (inx) at (0,0) {inextendable};
\draw[-latex] (dev) -- (cpl);
\draw[-latex] (cpl) -- (vst);
\draw[-latex] (vst) -- (inx);
\end{tikzpicture}
\]
and also prove by examples that none of them reverse.
For reductive geometries, we will prove the same for the larger diagram of implications
\[
\begin{tikzpicture}
\node (dev) at (0,4) {developable};
\node (cpl) at (0,3) {complete};
\node (gc) at (0,2) {geodesically complete};
\node (vst) at (0,1) {vast};
\node (inx) at (0,0) {inextendable};
\draw[-latex] (dev) -- (cpl);
\draw[-latex] (cpl) -- (gc);
\draw[-latex] (gc) -- (vst);
\draw[-latex] (vst) -- (inx);
\end{tikzpicture}
\]
\begin{example}
The Clifton can is complete but not developable so complete does not imply developable, so vast does not imply developable.
The schichtorte is vast but not complete, as is the Clifton plane (see page~\pageref{section:Clifton.2}).
\end{example}
My guess is that \emph{vast} is the most useful of these (and related) concepts, since it is easiest to test for and has the following theorems as consequences, which were previously proven only for those stronger concepts.
\begin{example}
If the complete cone \(\CompleteCone[M]\subseteq\LieG\) projects under \(\LieG\to\LieG/\LieH\) to a spanning set of \(\LieG/\LieH\) then the geometry is vast.
\end{example}
\begin{example}
The standard flat projective connection on hyperbolic space is not vast. 
Take as model \((X,G)\) the real projective plane \(X\) and the group \(G\) of projective linear transformations, and \(X=G/H\) where the structure group \(H\) is the group of projective linear transformations fixing a point \(x_0\in X\) of projective space.
Hyperbolic space \(M\) embeds into projective space as the interior of a ball in affine space.
The boundary of hyperbolic space is then a sphere.
Let \(M\subset X\) be hyperbolic space, i.e. that ball, and \(\G\subseteq G\) its induced flat projective connection: \(\G\) is the open set of points of \(G\) which map to \(M\) under the projection \(g\in G\mapsto gx_0\in X\).
Recall~\vpageref{example:hyperbolic.plane}: infinitesimal isometries of hyperbolic space preserve \(\G\) and \(M\), while flows of nonzero constant vector fields preserve \(\G\) just when they are in \(\LieH\); all other constant vector fields are incomplete on \(\G\).
Let us see this from another point of view.

Take a constant vector field \(A\) on \(G\). 
From lemma~\vref{lemma:adj.H}, as we move up the fiber of the bundle, the projection of \(A\) to a tangent vector in the underlying projective space transforms in the adjoint action of the structure group. 
This action gives rise to arbitrary linear transformations of that projection.
Hence the resulting vector does not stay tangent to any submanifold unless it vanishes.
The constant vector fields on \(\G\) are precisely those of \(G\) restricted to \(\G\).
A constant vector field on \(G\) has flow preserving \(\G\) just when its projection from any point is tangent to the sphere \(\partial M\), so just when the constant vector field is vertical everywhere.
Hence \(\CompleteCone=\LieH\): hyperbolic space is not vast, as a manifold with projective connection.
\end{example}
\begin{example}
The same argument proves that the complement of any subset of projective space is not vast, as a manifold with projective connection.
For instance affine space is not vast.
\end{example}
\begin{example}
The same argument proves that the complement of any subset of the sphere is not vast, as a manifold with conformal connection.
For instance conformal Euclidean space is not vast.
\end{example}
\begin{example}
There are Lorentzian geometries whose null geodesics are complete, but with no complete spacelike or timelike geodesics \cite{Beem.Ehrlich.Easley:1996} p. 203, \cite{Ehrlich:2006} p. 7, \cite{ONeill:1983} chapter 5 example 43.
Take as model for Lorentzian geometry \(X=\R^{3,1}\) with \(G=H\ltimes X\), where \(H\) is the Lorentz group.
So \(\LieG=\LieH\oplus\R^{3,1}\).
The null geodesics are the flows of the constant vector fields \(v_{\G}\) for any nonzero null vector \(v\in\R^{3,1}\subset\LieG\).
Clearly these \(v\) span \(\R^{3,1}=\LieG/\LieH\), so these geometries are vast but geodesically incomplete, so incomplete and so not developable. 
\end{example}
\begin{theorem}
A Cartan geometry is vast just when its associated effective Cartan geometry is vast.
\end{theorem}
The proof is as for theorem~\vref{thm:eff.complete}.
\subsection{Isomorphisms: local to global}
\begin{theorem}
In a vast real analytic Cartan geometry on a connected manifold, every infinitesimal automorphism defined on a connected open set of the total space extends uniquely to a complete infinitesimal automorphism.
Any two infinitesimal automorphisms which agree on an open set agree everywhere.
\end{theorem}
\begin{proof}
Our infinitesimal automorphism commutes with the complete constant vector fields, and with the structure group action, hence wiith the infinitesimal generators of the structure group action.
Hence it commutes with the brackets of these various vector fields.
By analyticity, the infinitesimal automorphism commutes with the (perhaps locally defined) flows of those vector fields, so is globally defined by vastness.
The time for which its flow is defined is invariant under those flows, and under the structure group action, so the same everywhere, so it is complete.
\end{proof}
\begin{theorem}\label{thm:from.vast}
A local isomorphism of Cartan geometries on connected manifolds, from a vast Cartan geometry, is a covering map.
\end{theorem}
\begin{proof}
Take a local isomorphism \(\G\to\G'\).
Let \(\CompleteCone\subseteq\LieG\) be the complete cone of \(\G\) and \(\barCompleteCone\) its associated complete module.
Let \(\CompleteCone[\G],\barCompleteCone[\G]\) be the associated vector fields on \(\G\), and
\(\CompleteCone[\G'],\barCompleteCone[\G']\) be the associated vector fields on \(\G'\) \emph{but} from the complete cone and complete module on \(\G\), not on \(\G'\).
By the orbit theorem (theorem~\vref{thm:orbit}), the orbits of the complete constant vector fields in \(\G\) are preserved by the vector fields of \(\barCompleteCone[\G]\).
Since \(\barCompleteCone=\LieG\), the \(\CompleteCone[\G]\)-orbits are open sets.
Hence they are components of \(\G\).
So the \(H\CompleteCone[\G]\)-orbits in \(\G\) are preimages of components of \(M\).
But \(M\) is connected, so \(\G\) is a single \(H\CompleteCone[\G]\)-orbit.

Each constant vector field on \(\G\) maps to a constant vector field on \(\G'\).
Flow lines map to flow lines.
Since the flow lines of \(\CompleteCone[\G]\) are complete, they map to complete flow lines on \(\G'\).
(We \emph{dont'} assert that the complete cone of \(\G\) maps to that of \(\G'\).)
The \(\CompleteCone[\G]\)-orbit through a point maps to the \(\CompleteCone[\G']\)-orbit through the corresponding point.
So these are open too, throughout the points in the image of the local isomorphism \(\G\to\G'\).

Take any point of \(\G'\) in the closure of the image of \(\G\to\G'\).
The elements of \(\barCompleteCone=\LieG\) are constructed by iterated Lie brackets in the curvature deformed bracket, on elements of \(\CompleteCone\) for which the curvature is constant.
The constant stays the same in the closure of the image of \(\G\to\G'\), by continuity.
So we can build up any constant vector field on \(\G'\), iteratively, by such brackets, and its flow preserves the \(\CompleteCone[\G']\)-orbit.
So that orbit is open, for any point in the closure of the image of \(\G\to\G'\).
But then it overlaps some orbit in the image.
Hence it equals that orbit.
So the \(\CompleteCone[\G]\)-orbits in \(\G'\) extend past every point of the closure of the image of \(\G\to\G'\).
But they consist of the images of the (complete!) flow lines in \(\G\).
So they lie in the image of \(\G\to\G'\).
So each compoent of \(\G\) has image a single component of \(\G'\), which is a single \(\CompleteCone[\G]\)-orbit.
So there is a single \(H\CompleteCone[\G]\)-orbit in \(\G\): all of \(\G\), and the same in \(\G'\).
Apply theorem~\vref{cor:orbitMapEquiv} to see that \(M\to M'\) is a covering map.
\end{proof}
\begin{corollary}\label{cor:vast.flat}
On any connected manifold with any flat Cartan geometry, with strong effective connected model, the following are equivalent:
\begin{itemize}
\item
The geometry is vast.
\item
The geometry is complete.
\item
The geometry is developable.
\item
The geometry has universal covering space the same as the model, with the same geometry.
\item
The manifold is a pancake as a submanifold of itself, and its developing map is a covering map.
\item
The geometry is a complete Klein geometry \(M=\Gamma\backslash \tilde{X}\) of the universal covering homogeneous space of the model.
\end{itemize}
\end{corollary}
\begin{corollary}\label{corollary:analytic.vast}
Suppose that \(M,M'\) are manifolds with real analytic Cartan geometries with the same model.
Suppose that \(M\) is connected, and that the geometry on \(M'\) is vast.
Take a connected open set \(U\subseteq M\).
Suppose that \(U\to M'\) is a real analytic map which is locally an isomorphism of Cartan geometries.
Then this map extends uniquely to a real analytic local isomorphism of Cartan geometries on a connected covering space \(\tilde{M}\to M\),
\[
\begin{tikzcd}
&\tilde{M}\arrow[r]\arrow[d]&M'\\
U\arrow[ur]\arrow[r]&M
\end{tikzcd}
\]
\end{corollary}
This corollary generalizes a result of Ehresmann \cite{Ehresmann:1951} p. 50.
\begin{proof}
Denote the Cartan geometries as \(\G\xrightarrow{\pi}M\), \(\G'\xrightarrow{\pi'}M'\).
Inside the product bundle \(\G\times\G'\), the graph of our isomorphism \(\left.\G\right|_U\to\G'\) is a submanifold tangent to all vector fields \(A_{\G}+A_{\G'}\), \(A\in\LieG\).
Under the flows of any of these vector fields, this graph is carried to a submanifold which, for small time of flow, overlaps the original submanifold on an open set.
By analyticity, the submanifold continues to be tangent to those vector fields.
Let \(\G''\) be the orbit under those vector fields of some point of our submanifold.
So \(\G''\) is again an immersed connected submanifold, of the same dimension, tangent to those vector fields; see theorem~\vref{thm:orbit}. 
Under projections to \(\G\) and \(\G'\), \(A_{\G}+A_{\G'}\mapsto A_{\G},A_{\G'}\), so \(\G''\to\G,\G'\) are local diffeomorphisms.
Clearly \(\G''\) is invariant under diagonal action of the identity component of \(H\).
By the equivariant Frobenius theorem (corollary~\vref{corollary:equivariant.Frobenius}), the union of \(H\)-translates of \(\G''\) is an immersed submanifold, which we also call \(\G''\), \(H\)-invariant and invariant under our flows.
By theorem~\vref{cor:orbitMapEquiv} the map \(\G''\to\G\) is a fiber bundle mapping as is the quotient \(M'':=\G''/H\to M\).
Since each is a local diffeomorphism, it is a covering map.
Our map \(\G''\to G'\) by construction contains the graph of \(\left.\G\right|_U\to\G'\).
Quotient: \(M''\) contains a lift of \(U\).
\end{proof}
\begin{example}
Any local isomorphism of real analytic schichttorten defined on a connected open set extends to a global isomorphism, as they are connected and simply connected and vast.
Note that their Riemannian metrics are not complete.
\end{example}
\begin{example}
Our theorem describes a vast target.
What about a vast source?
For any homogeneous space \((X,G)\), every flat \((X,G)\)-geometry contains an open subset of \(X\), but it is complete just when that inclusion extends to a map from \(X\), even if \(X\) is compact and simply connected, and the geometry is analytic on a compact manifold.
\end{example}
\section{Geodesic completeness}
Take a reductive homogeneous space \((X,G)\), say \(X=G/H\) so 
\[
\LieG=\LieH\oplus V
\]
as \(H\)-modules.
Take any \((X,G)\)-Cartan geometry \(\G\to M\), with Cartan connection 
\[
\omega=(\gamma,\sigma).
\]
Take the vector field \(Z\) on \(\G\times V\) given by
\[
Z(p,v)\hook
\omega
=
(0,v)
\]
and \(Z\hook dv=0\), the \emph{geodesic vector field};\define{geodesic!vector field} its flow is the \emph{geodesic flow}.\define{geodesic!flow}
A reductive Cartan geometry is \emph{geodesically complete}\define{geodesically complete} if its geodesic vector field is complete.
\begin{problem}{bh}
Prove that this vector field \(Z\) descends under the diagonal \(H\)-action to a vector field on \(\amal{\G}{H}{V}=TM\), also called the \emph{geodesic vector field};\define{geodesic!vector field} its flow is also called the \emph{geodesic flow}.\define{geodesic!flow}
\end{problem}
\Danger{} A reductive geometry is geodesically complete just when, if \(\LieG=\LieH\oplus V\), the constant vector fields of \(V\) are complete.
Hence completeness of the Cartan geometry implies completeness of the geodesic flow.
(In particular, for example, the model of a reductive geometry, i.e. any reductive homogeneous space, is complete as a Cartan geometry and hence geodesically complete.) 
The Clifton plane (see page~\pageref{section:Clifton.2}) is an incomplete but geodesically complete reductive geometries, contradicting the generally accepted theorem \cite{Kobayashi:1957} p. 181 Theorem 5.15 claiming that geodesic completeness is equivalent to completeness.

\Danger{} Noting that the constant vector fields of \(\LieH\) are always complete, it would appear that the Trotter product formula could apply, on a reductive geometry with complete geodesic flow, to a constant vector field from \(\LieH\) and one from \(V\), and thereby prove completeness of their sum, and hence completeness of the Cartan geometry \cite{Kobayashi:1957} p. 182.
Keep in mind that two complete vector fields can have incomplete sum, so the Trotter product formula's terms are defined, but might leave every compact set as we take the limit.

\Danger{} Any Cartan geometry is complete just when its scaffold is geodesically complete.
But the scaffold's Riemannian geometry is geodesically complete just when that Riemannian geometry is metrically complete, i.e. just when the Cartan geometry is developable.
Paradoxically, it would then seem that completeness and developability are the same, and this has probably trapped many researchers into thinking so, but the affine connection of the scaffold is \emph{not} necessarily the Levi--Civita connection, so there are two notions of geodesic completeness on the scaffold.
\begin{theorem}
Geodesically complete implies vast.
\end{theorem}
\begin{proof}
Every Cartan geometry has \(\LieH\subseteq\CompleteCone\).
A reductive geometry has \(H\)-module decomposition \(\LieG=\LieH\oplus V\).
It is geodesically complete just when \(V\subseteq\CompleteCone\).
But then \(\CompleteCone\) contains \(H\cup V\), so the Lie algebra \(\barCompleteCone\) contains the sum, i.e. \(\LieG\).
\end{proof}
\begin{theorem}
A local isomorphism of reductive geometries with the same model, on connected manifolds, from a geodesically complete reductive geometry, is a covering map.
\end{theorem}
\begin{proof}
Apply theorem~\vref{thm:from.vast}.
\end{proof}
\begin{corollary}\label{corollary:analytic.complete.2}
Suppose that \(M,M'\) are manifolds with real analytic Cartan geometries with the same model, a reductive homogeneous space.
Suppose that \(M\) is connected, and that the geometry on \(M'\) is geodesically complete.
Take a connected open set \(U\subseteq M\).
Suppose that \(U\to M'\) is a real analytic map which is locally an isomorphism of Cartan geometries.
Then this map extends uniquely to a real analytic local isomorphism of Cartan geometries on a connected covering space \(\tilde{M}\to M\),
\[
\begin{tikzcd}
&\tilde{M}\arrow[r]\arrow[d]&M'\\
U\arrow[ur]\arrow[r]&M
\end{tikzcd}
\]
\end{corollary}
\begin{proof}
Apply corollary~\vref{corollary:analytic.vast}.
\end{proof}
\begin{problem}{bi}
If a homogeneous space \((X,G)\) carries an invariant Riemannian metric, then any \((X,G)\)-geometry is geodesically complete if and only if it is developable.
\end{problem}

\subsection{Developability from completeness}
If we have a Cartan geometry \(\G\to M\) with reductive model \((X,G)\), say \(\LieG=\LieH\oplus V\), then the Cartan connection splits correspondingly \(\omega=(\gamma,\sigma)\) into a connection \(\gamma\) for \(\G\to M\) and the soldering form \(\sigma\).
Denote by \(H\xrightarrow{\rho}\GL{V}\) and \(\LieH\xrightarrow{\rho}\LieGL{V}\) the representation of \(H\) on \(V\).
Every connection on \(\G\to M\) has the form \(\gamma':=\gamma+q\sigma\) for a unique \(H\)-equivariant \(\G\xrightarrow{q}\LieH\otimes V^*\).
The torsion changes to \(t'=t+\delta q\), where
\[
\LieH\otimes V^*\xrightarrow{\delta}V\otimes\Lm{2}{V}^*
\]
is defined by, for \(A\in\LieH\) and \(\xi\in V^*\), if \(q=A\otimes\xi\), by
\[
(\delta q)(v,w)=\rho(A)(v)\xi(w)-\rho(A)(w)\xi(v).
\]
\begin{theorem}
Take a Cartan geometry with curvature \(k\) and with model \((X,G)\) where \(X=G/H\).
Suppose that \(H\) has finitely many components and that \(H_c\subseteq H\) is a maximal compact subgroup.
Suppose that there is a constant \(q_0\in\LieH_c\otimes(\LieG/\LieH_c)^*\) so that
\[
k+\delta q_0\in\LieH_c\otimes\Lm*{2}{\LieG/\LieH_c}.
\]
Then the Levi--Civita connection of the scaffold is the connection induced from the splitting into \(\LieG=\LieH_c\oplus(\LieG/\LieH_c)\) up to adding \(q_0\) times the soldering form.
Hence the Cartan geometry is developable just when it is complete.
\end{theorem}
\begin{proof}
We can assume that \(H_c=H\) since we just work on the scaffold directly.
Split \(\omega=\gamma+\sigma\) invariant under \(H\).
Compute the torsion of the scaffold geometry as the projection of the curvature:
\[
t:=k+\LieH\otimes\Lm*{2}{\LieG/\LieH}.
\]
Hence the torsion of the scaffold arises in the soldering form as
\[
d\sigma+\gamma\wedge\sigma=\frac{1}{2}t\sigma\wedge\sigma.
\]
So the affine connection \(\gamma+q_0\sigma\) has vanishing torsion.
Let \(n:=\dim(\LieG/\LieH)\).
On \(\G\times\Orth{n}\), define \(\gamma':=\omega_H+\Ad_h^{-1}(\gamma+q_0\sigma)\) for \(h\in\Orth{n}\), descending to \(\G':=\amal{\G}{H}{\Orth{n}}\).
Vanishing torsion ensures that \(\gamma'\) has vanishing torsion too, so is the Levi--Civita connection.
Suppose that the Cartan geometry is complete.
We demand that \(q_0\) is constant to ensure that the geodesic flow vector fields of \(\gamma+q_0\sigma\) are constant vector fields, hence are complete too.
Completeness of these geodesic flows of \(\gamma+q_0\sigma\) on \(\G\) implies completeness on \(\amal{\G}{H}{\Orth{n}}\) of the same vector fields, by right invariance of the vector fields, hence completeness of the Riemannian geometry.
\end{proof}
\begin{corollary}
Flat Cartan geometries are developable just when complete and complete just when vast.
\end{corollary}
\section{Extendability}
\begin{example}
On the unit disk, with \(r^2:=x^2+y^2\), the Riemannian metric
\[
ds^2=(1-r^2)(dx^2+dy^2)
\]
extends as a metric to the one point compactification; the resulting metric is not given by a Riemannian metric on the sphere.
The metric is invariant under reflection in any line through the origin, so lines through the origin are geodesics.
The metric has Gauss curvature 
\[
G=\frac{4}{(1-r^2)^3},
\]
so in the orthonormal frame \(\omega_1=\sqrt{1-r^2}dx, \omega_2=\sqrt{1-r^2}dy\), \(dG=G_1\omega_1+G_2\omega_2\) where 
\begin{align*}
G_1&=\frac{24x}{(1-r^2)^{9/2}},\\
G_2&=\frac{24y}{(1-r^2)^{9/2}}.
\end{align*}
The Riemannian metric induces a normal projective connection, with curvature represented by the tensor
\[
(G_2\omega^1-G_2\omega^2)\otimes(\omega^1\wedge\omega^2)
=
\frac{24}{(1-r^2)^3}(y\,dx-x\,dy)\otimes(dx\wedge dy),
\]
so the projective connection does not extend to the points of the unit circle.
If it extended to the one point compactification smoothly, there would be an area form near that point, so an area form in our \(x,y\) coordinates which gives small areas to regions near the unit circle, hence \(f(x,y)dx\wedge dy\) with \(f\to 0\) as we approach the unit circle.
Our tensor, in terms of this area form, becomes a \(1\)-form with large integral around circles.
So the normal projective connection does not extend to the metric space completion.
\end{example}
\begin{example}
Consider the Lorentz metric
\[
g=\frac{dx\,dy}{\cos^2x\sin^2y+\sin^2x\cos^2y},
\]
which is defined on the plane except at the points
\[
(x,y)=\frac{\pi}{2}(k,\ell)
\]
for integers \(k,\ell\) with \(k+\ell\) even \cite{Bavard.Mounoud:2013}.
The \emph{Clifton--Pohl plane}\define{Clifton--Pohl!plane} \(M\) is the plane punctured at those points and with that metric.
The open subset \(U\subset M\) consisting of points \((x,y)\) with \(-\pi/2<x,y<\pi/2\) maps to the plane punctured at the origin by
\[
(x,y)\mapsto(X,Y):=(\tan x,\tan y).
\]
Check that this map pulls back the Lorentz metric
\[
\frac{dX\,dY}{X^2+Y^2}
\]
to \(g\).
The \emph{Clifton--Pohl torus}\SubIndex{Clifton--Pohl torus} is the quotient of the punctured plane by
\[
(X,Y)\mapsto(\lambda X,\lambda Y)
\]
for some real \(\lambda\ne\pm 1\).
This Lorentz metric descends to a Lorentz metric on the Clifton--Pohl torus invariant under the group of all such rescalings, which is a circle action on the Clifton--Pohl torus.
Hence the open set \(U\subset M\) is an infinite covering space of the Clifton--Pohl torus.
Note that \(U\subsetneq M\), i.e. the universal covering space of the Clifton--Pohl torus is extendable, hence incomplete.
Hence the Clifton--Pohl torus is incomplete.
It turns out that the Clifton--Pohl plane is complete.
\end{example}
\begin{example}
Similarly, on any vector space \(V\), take two quadratic forms \(q,Q\), with \(Q\) positive definite and \(q\) nondegenerate, and take some real number \(\lambda>1\).
The pseudo-Riemannian metric \(q(dx)/Q(x)\) for \(x\in V\) descends to the Hopf manifold \(M:=(V-0)/(x\cong\lambda x)\); with this pseudo--Riemannian metric, \(M\) is the \emph{Clifton--Pohl Hopf manifold}.\define{Clifton--Pohl!Hopf manifold}
If \(\dim V=2\) and \(q\) is Lorentian then \(M\) is the Clifton--Pohl torus.
It would be nice to find a complete Lorentzian geometry containing the universal covering space of the Clifton--Pohl Hopf manifold, as \cite{Bavard.Mounoud:2013} did for the Clifton--Pohl torus.
\end{example}
\begin{problem}{cj}
The Clifton--Pohl Hopf manifold is incomplete.
\end{problem}
Take a homogeneous space \((X,G)\), \(X=G/H\), and suppose that \(H\) has maximal compact subgroup \(H_c\subseteq H\).
As usual, we take an \(H_c\)-invariant positive definite inner product \(\left<,\right>\) on \(\LieG\).
The resulting Riemannian metric \(\left<\omega,\omega\right>\) on \(\G\) drops to a Riemannian metric on the scaffold \(M_c\).
Take the metric space completion \(\bar{\G}\) of \(\G\).
Since the metric is \(H_c\)-invariant, \(\G\to M_c\) extends to metric space completions \(\bar{\G}\to\bar{M}_c\).
The action of \(H\) on \(\bar{\G}\) is continuous, but not necessarily isometric.
Let \(\bar{M}:=\bar{\G}/H\), a topological space containing \(M\).
The topological space \(\bbb{M_c}:=\bar{M}_c-M_c\) with \(H\)-action is the \emph{\bboundary{}} \cite{Schmidt:1974,Schmidt:1979,Clarke:1979,Frances:2012}.
(\Danger{} This is not quite the usual definition of \bboundary{}, which is usually \(\bar\G-\G\), but I prefer this definition.)
The image \(\bbb{M}\subseteq\bar{M}\) of \(\bbb{M_c}\subseteq\bar{M}_c\) is also sometimes called the \bboundary{}.
\begin{example}
A Cartan geometry is developable if and only if \(M_c\) is complete, i.e. \(\bar{M}_c=M_c\), i.e. \(\bbb{M_c}\) is empty, i.e. \(\bbb{M}\) is empty.
\end{example}
\begin{example}
If \(M\subseteq X\) is an open set, then \(\bbb{M}_c\) is the boundary of \(M_c\) as a subset of \(X_c\).
So if we puncture \(X\) at a single point, \(M:=X-\set{x_0}\) has \(\bbb{M_c}=H/H_c\).
\end{example}
\begin{example}
Take \((X,G):=(S^n,\PO{n+1,1})\), the model of conformal geometry; \(G\) is the group of linear transformations preserving a Lorentz quadratic form, up to scaling, \(X\) the set of null lines of the quadratic form, \(H\) the subgroup of \(G\) fixing a null line, \(H_c\) the subgroup fixing that null line, and some time-like line, and the plane between them, and a complementary space-like \(n\)-dimensional linear subspace.
So \(X_c\) is the collection of all such.
We leave the reader to check that each fiber \(H/H_c\) is isometric to hyperbolic space.
So if \(M=X-\set{x_0}\) is a punctured sphere, by Ptolemaic projection \(M\) is conformally diffeomorphic to Euclidean space, and \(M_c\) is a trivial bundle of hyperbolic spaces over Euclidean space, \(\bbb{M}_c\) is hyperbolic space, and \(\bbb{M}\) is a point.
\end{example}
\begin{example}
Again take \((X,G):=(S^n,\PO{n+1,1})\).
Let \(M:=V/\Lambda\) a torus, i.e. a quotient of Euclidean space by a lattice.
So \(\tilde{M}=V\subset X=S^n=V\sqcup\set{\infty}\).
Under \(\Lambda\)-action on \(V\), every point gets arbitrarily close to the north pole \(x_0\in S^n\).
So \(\bbb{M_c}\) contains the quotient of the hyperbolic space over \(x_0\), acted on by the free action of \(\Lambda\) by isometries, arbitrarily close to every fiber of \(M_c\to M\).
So \(\bbb{M}\) is nowhere Hausdorff.
\end{example}
An \((X,G)\)-geometry \emph{embeds}\define{embedding!of Cartan geometries}\define{Cartan!geometry!embedding} in another if the bundle embeds equivariantly, matching Cartan connections.
If some \((X,G)\)-geometry \(\G\to M\) embeds into another one \(\G'\to M'\) then the same holds for the associated scaffolds.
If \(M'\) is not the image of \(M\), the Cartan geometry \(\G\to M\) is \emph{extendable}.\define{Cartan!geometry!extendable}\define{extendable!Cartan geometry}
Near each point \(\bar{m}\) in the topological closure \(\bar{M}_c^{\text{top}}\) of \(M_c\subseteq M'_c\), \(M_c'\) is a smooth Riemannian manifold, and so its sufficiently small metric balls are compact. 
Every point of the metric space completion \(\bar{M}_c\) of \(M_c\) within a small enough distance of some point of \(M_c\subseteq M'_c\) lies inside that compact ball, as a point of \(\bar{M}_c^{\text{top}}\): a map \(\bar{M}_c^{\text{top}}\to\bar{M}_c\).
A point \(\bar{m}\in \bbb{M_c}\) is \emph{extendable}\define{extendable point} if it arises as a point of \(\bar{M}_c^{\text{top}}\) in this way.
By \(H\)-invariance, the extendable points of \(\bbb{M_c}\) are \(H\)-orbits above points of \(\bbb{M}\), the \emph{extendable} points of \(M\).
Hence:
\begin{theorem}%
[Frances \cite{Frances:2014}]%
\define{theorem!Frances}
\define{Frances theorem}
Every extendable point of \(\bbb{M}\) lies in a Hausdorff open set in \(\bbb{M}\) over which \(\bbb{M_c}\to\bbb{M}\) is a topological fiber bundle, with every fiber isometric to \(H/H_c\).
\end{theorem}
\begin{theorem}
Take a Cartan geometry \(\G\to M\) and consider its scaffold \(H_c\to\G\to M_c\).
Take some polynomial invariants of \(H_c\) acting on the curvature and its derivatives, as \(H_c\)-modules.
Put these invariants into a single column vector \(M_c\to \R^N\).
If this invariant map is a proper map then the Cartan geometry \(\G\to M\) is not extendable.
\end{theorem}
\begin{proof}
Take an embedding \(M\subseteq M'\) and compute these invariants on \(M_c'\), near some point \(p\) of the boundary of \(M\) in \(M'\), giving a point of \(\R^N\).
Take a compact set around that point; its preimage in \(M_c\) is compact, so doesn't approach \(p\).
\end{proof}
\begin{example}
Trivially the Riemannian metric \(|x| \sum dx_i^2\) on Euclidean space punctured at the origin is conformal to the standard one, so extends conformally but the Riemannian metric does not extend.
\end{example}
\begin{example}
The standard flat affine connection on Euclidean space extends, as a projective connection, to projective space.
On the other hand, every torsion-free affine connection on a manifold of dimension \(4\) or more induces a unique normal projective connection for which the curvature of the projective connection is the Weyl curvature of the original affine connection \cite{KobayashiNagano:1964}.
Take some polynomials in the Weyl tensor components invariant under \(\Orth{n}\).
If these form a proper map from \(M_c\), then the affine connection does not extend projectively.
One can use this technique to prove that various spacetimes in general relativity are projectively complete, via long computations (which clearly we leave to the reader).
\end{example}

\chapter{Flat geometries}
Restating theorem~\vref{thm:develop.submanifold} as a theorem about a manifold rather than a submanifold:
\begin{theorem}\label{thm:flat}
Take a strong effective homogeneous space \((X,G)\), and a connected manifold \(M\) with a flat \((X,G)\)-Cartan geometry.
The geometry on \(M\) is the quotient of the pullback from a local isomorphism
\[
\begin{tikzcd}
\tilde{M}\arrow[r]\arrow[d]&X\\
M
\end{tikzcd}
\]
the \emph{developing map},\define{developing map} from the universal covering space \(\tilde{M}\to M\), equivariant for a group morphism \(\fundamentalGroup{M}\to G\), the \emph{holonomy morphism}.\define{holonomy!morphism}
The pair \((\delta,h)\) of developing map and holonomy morphism are unique up to replacing by \((g\delta,\Ad_g h)\) for any \(g\in G\).
\end{theorem}
\Danger{} Sharpe \cite{Sharpe:1997} has a notion of \emph{geometric orientability};\SubIndex{geometric orientability} ignore it and you get this stronger result.
\begin{example}
The model \(X\), if connected, has developing map its universal covering map, holonomy trivial.
\end{example}
\begin{example}
The sphere, as \(2\)-to-\(1\) cover of real projective space, has a pullback \((\RP{n},\PGL{n+1})\)-structure, with developing map the \(2\)-to-\(1\) covering map, and trivial holonomy.
\end{example}
\begin{example}
Euclidean space which maps conformally diffeomorphically to the punctured sphere, by Ptolemaic projection.
So this must be the developing map.
So any open subset of the sphere conformal to Euclidean space is the sphere punctured once.
\end{example}
\begin{example}
Take a vector space \(V\) with a positive definite inner product and a lattice \(\Lambda\subseteq V\), i.e. the set of integer linear combinations of some basis vectors.
Consider the flat conformal structure on a torus \(M=V/\Lambda\).
Its universal covering space is \(\tilde{M}=V\) with developing map \(V\to S^n=V\cup\set{\infty}\) as seen via Ptolemaic projection.
\end{example}
\begin{example}
We can quotient out the kernel of the holonomy morphism; the developing map is defined as a map \(\hat{M}\to X\), where \(\hat{M}=\tilde{M}/\ker h\).
For a connected model \((X,G)\), the developing map (thus quotiented) is the identity map \(X\to X\).
\end{example}
\begin{example}
Take a homogeneous space \((X,G)\), an open set \(X'\subseteq X\) and a Lie subgroup \(G'\subseteq G\) acting transitively on \(X'\).
Any flat \((X,G)\)-Cartan geometry arises from a flat \((X',G')\)-geometry just when it has holonomy morphism lying in \(G'\) and developing map image in \(X'\).
\end{example}
\begin{corollary}
A Cartan geometry is flat just when it is locally isomorphic to its model.
\end{corollary}
\begin{example}
On a sphere of dimension \(2\) or more, there is no flat affine structure: its developing map would be a local diffeomorphism \(S^n\to\R^n\), and we leave the reader to argue that there is no such.
In the same way, if a connected compact manifold \(M\) of finite fundamental group admits a flat Cartan geometry, the model is compact with the same universal covering space as \(M\): the universal covering homogeneous space, so \(M=\Gamma\backslash\tilde{X}\), some finite subgroup \(\Gamma\subseteq\tilde{G}\).
\end{example}
\begin{problem}{bl}
Which Hopf manifolds are isomorphic?
\end{problem}
\begin{problem}{bm}
Prove that \(M=\R\) admits three affine structures, i.e. Cartan geometries modelled on \((X,G)\) with \(X\) the real number line and \(G\) the group of affine transformations \(x\mapsto ax+b\) of \(X\), \(a\ne 0\).
Hint: \(\R,\R^{>0},(0,1)\).
\end{problem}
\subsection{Curves}\label{subsec:curves}
The classification of connected and simply connected effective homogeneous spaces \((X,G)\) with \(G\) connected is elementary but long for \(\dim X=0,1,2\) \cite{Komrakov.Churyumov.Doubrov.1993}.

Suppose that \((X,G)\) is an effective homogeneous curve with \(X\) connected. 
We saw \vpageref{exercise:univ.cover.replace} that \((X,G)\)-structures are the same if we replace \((X,G)\) by the universal covering homogeneous space, so we can assume that \(X\) is a connected and simply connected curve, i.e. \(X=\R\) up to diffeomorphism.

Consider \((X,G)\)-geometries on \(M=\R\).
We have a developing map \(M\to X\), a local diffeomorphism, so an increasing or decreasing map, so a diffeomorphism to its image, which is an open interval \((a,b)\subset\R\), say \(M=(a,b)\subseteq X=\R\).
The fundamental group is \(\set{1}\), so the holonomy morphism is \(1\mapsto 1\).
Hence the developing map is uniquely determined up to action of an element \(g\in G\).
So the moduli space of \((X,G)\)-geometries on \(M=\R\) is the set of all intervals \(M=(a,b)\subseteq X=\R\), modulo action of \(G\).
The automorphism group of the geometry on \(M\) then consists of the elements of \(G\) acting preserving the interval \(M=(a,b)\subseteq X=\R\), by theorem~\vref{theorem:open.set.G.aut}.
In particular, \(M=X=\R\) if and only the geometry is isomorphic to the model.
If \(b=\infty\) and \(a\) is finite, by transitivity of the \(G\)-action on \(X\), we can arrange \(a=0\), so \(M=(0,\infty)\), and the symmetry group of the Cartan geometry is the subgroup of \(G\) consisting of orientation preserving maps fixing the origin.
Similarly, if \(a=-\infty\) and \(b\) is finite, we can arrange \(b=0\) so \(M=(-\infty,0)\).
But then there is the question whether \(G\) acts preserving orientation of \(X=\R\), i.e. by increasing functions, in which case we can only arrange \(M=(-\infty,0)\), and the automorphism group of the Cartan geometry is the subgroup of \(G\) consisting of increasing maps fixing the origin.
On the other hand, if \(G\) contains an orientation reversing map, we return to the case of \(M=(0,\infty)\).
The classification of \((X,G)\)-geometries on \(M=(a,b)\subset X\) if \(a\) and \(b\) are both finite, is the classification of ordered pairs of distinct points of \(X\), up to \(G\)-action.
\begin{example}
If \(M\) is a connected and simply connected curve and \((X,G)=(\R{},\Euc{\R})\) then \(G\) contains an orientation reversing map \(x\mapsto -x\), so \(M=X\) or \(M=(0,\infty)\) or \(M=(a,b)\) and \(G\) acts on ordered pairs of distinct points in \(X\) preserving their distance, but transitively on pairs of points of the same distance, so \(M=(0,\ell)\) with \(\ell>0\) arbitrary.
Our classification is
\[
\begin{array}{ll}
\toprule
M&\Aut[M]\\
\cmidrule(r){1-1}\cmidrule(l){2-2}
(0,\ell)&\set{1,x\mapsto\ell-x}\\
(0,\infty)&1\\
X&G\\
\bottomrule
\end{array}
\]
\end{example}
\begin{example}
If \(M\) is a connected and simply connected curve and \((X,G)=(\R{},\Aff{\R})\) then \(G\) contains an orientation reversing map \(x\mapsto -x\), so \(M=X\) or \(M=(0,\infty)\) or \(M=(a,b)\).
The group \(G\) acts on ordered pairs of distinct points in \(X\) transitively, so \(M=(-1,1)\).
Our classification is
\[
\begin{array}{ll}
\toprule
M&\Aut[M]\\
\cmidrule(r){1-1}\cmidrule(l){2-2}
(-1,1)&\set{1,x\mapsto-x}\\
(0,\infty)&\set{x\mapsto \lambda x|\lambda>0}\\
X&G\\
\bottomrule
\end{array}
\]
\end{example}
%
%\R&\Aff{\R}&\R&\R&&\Aff{\R}\\
%\R&\Aff{\R}&\R&(0,\infty)&&\R^{+}\\
%\R&\Aff{\R}&\R&(-1,1)&&\pm 1\\
%\R&\Aff{\R}&S^1&\R&(1,1)&\pm 1\ltimes (\R/\Z{})\\ 
%\R&\Aff{\R}&S^1&(0,\infty)&(\lambda,0)&\R^+/\left<\lambda\right>\\ 
%\R&\Euc{\R}&\R&\R&&\Euc{\R}\\ 
%\R&\Euc{\R}&\R&(0,\infty)&&1\\ 
%\R&\Euc{\R}&\R&(0,\ell)&&\set{1,(-1,\ell)}\\
%\R&\Euc{\R}&S^1&\R&(1,\ell)&\pm{1}\ltimes(\R/\ell\Z)\\
%\bottomrule
%\end{array}
%\]
%\[
%M=
%\begin{cases}
%(0,\ell), \ell>0\text{ arbitrary}
%\]
%\end{example}

Similarly, for any homogeneous curve \((X,G)\), if \(X=\R\) and \(M=S^1\), we find that \(\tilde{M}=\R\) has pullback Cartan geometry as above, and we have to quotient by some holonomy generator.
But the holonomy generator acts on \(\tilde{M}\) without fixed point, i.e. as deck transformations of \(\tilde{M}\to M\), so in particular cannot reverse orientation of \(\tilde{M}\).
So the classification is that of \((X,G)\)-structures on connected and simply connected curves \(\tilde{M}\) together with the classification, up to conjugation by an automorphism, of the automorphisms that act on \(\tilde{M}\) without fixed points.
\begin{example}
If \(M\) is a connected but not simply connected curve and \((X,G)=(\R{},\Euc{\R})\) then, as above, \(\tilde{M}=(0,\ell), (0\,\infty)\) or \(X\).
The automorphisms \(g\) acting without fixed point are
\[
\begin{array}{ll}
\toprule
\tilde{M}&g\\
\cmidrule(r){1-1}\cmidrule(l){2-2}
(0,\ell)&\text{none}\\
(0,\infty)&\text{none}\\
X&x\mapsto x+r, r\ne 0\\
\bottomrule
\end{array}
\]
Hence the classification: \(M=X/(x\mapsto x+r)\), some \(r>0\), since we can conjugate by automorphism to get a translation to be by a positive difference.
So \(M\) is a circle of some radius \(r\).
Its automorphism group is
\[
\set{x\mapsto x+b}/(b\sim b+r)\sqcup\set{x\mapsto -x+b}/(b\sim b+r),
\]
a pair of circles.
\end{example}
\begin{example}
If \(M\) is a connected but not simply connected curve and \((X,G)=(\R{},\Aff{\R})\) then, as above, \(\tilde{M}=(-1,1), (0,\infty)\) or \(X\).
The automorphisms \(g\) acting without fixed point are
\[
\begin{array}{ll}
\toprule
\tilde{M}&g\\
\cmidrule(r){1-1}\cmidrule(l){2-2}
(-1,1)&\text{none}\\
(0,\infty)&\set{x\mapsto\lambda x|\lambda>1}\\
X&x\mapsto x+r, r\ne 0\\
\bottomrule
\end{array}
\]
Hence the classification: the \emph{Hopf circle}\define{Hopf circle} \(M=(0,\infty)/(x\mapsto \lambda x)\), some \(\lambda>1\), or \(M=X/(x\mapsto x+2\pi)\) (after conjugation), the \emph{Euclidean circle}.\define{Euclidean circle}
The map \(x\mapsto 2\pi \log x/\log \lambda\) descends to a canonical diffeomorphism from the Hopf circle to the Euclidean circle, identifying their automorphisms, but is not an isomorphism of affine structures, as it does not lift to an affine transformation of \(X\).
\[
\begin{array}{lll}
\toprule
M&\text{holonomy generator}&\Aut[M]\\
\cmidrule(r){1-1}\cmidrule(lr){2-2}\cmidrule(l){3-3}
\multicolumn{2}{l}{\textrm{the Hopf circle:}}\\
\qquad
\begin{array}{l}
\R^+/(x\sim\lambda x)\\ 
\lambda>1
\end{array}&
x\mapsto\lambda x
&
\begin{array}{l}
\set{x\mapsto \mu x}\\
\mu>0
\end{array}
\\[14pt]
\multicolumn{2}{l}{\textrm{the Euclidean circle:}}\\
\qquad\R/2\pi\Z&
x\mapsto x+2\pi&
\begin{array}{l}
\set{x\mapsto\pm x+b}\\
b\in\R/2\pi\\
\end{array}\\
\bottomrule
\end{array}
\]
\end{example}

\subsection{Flat geometries and Klein manifolds}
Recall that we can always safely assume that the model of any Cartan geometry is connected.
\begin{corollary}
Every flat vast Cartan geometry, on a connected manifold, with connected model, is a complete Klein manifold of the associated universal covering homogeneous space.
\end{corollary}
\begin{example}
Take any flat conformal structure on a manifold of dimension \(3\) or more, with infinite fundamental group, for example on a flat torus.
Its developing map is from its universal covering space, which is not compact, so is not a covering space of the sphere.
A conformal structure has model \((X,G)=(S^n,\PO{n+1,1})\), so \(X\) is simply connected and compact, so has no noncompact covering space.
So the original conformal geometry is not complete.
\end{example}
\begin{example}
The developing map of the projective connection on Euclidean space is the obvious embedding \(\R^n\to\RP{n}\) as \(\R^n=\RP{n}-\RP{n-1}\), an affine chart.
So the projective connection is not complete.
Projective space has universal covering \(S^n\to\RP{n}\) if \(n\ge 2\).
So the sphere punctured at a point is diffeomorphic to \(\R^n\) but has a different flat projective than the one on the affine chart, with geodesics closed if and only if they don't pass through the origin.
Similarly we can slice \(\RP{2}\) along various closed intervals of \(\RP{1}\), to get various flat projective connections on the plane.
\end{example}
\begin{example}
On every connected manifold with infinite fundamental group (for example, the flat torus) or noncompact (for example, Euclidean space), every flat projective connection is incomplete.
A flat projective connection on a connected manifold \(M\) with finite fundamental group is \(M=\Gamma\backslash S^n\), for some finite group \(\Gamma\subset\SL{n+1,\R}\) acting freely.
Since \(\SO{n+1}\subset\SL{n+1,\R}\) is a maximal compact subgroup, up to isomorphism, we can arrange that \(\Gamma\subset\SO{n+1}\), i.e. \(M\) is a quotient of the sphere by a finite group of rotations \(\Gamma\) acting freely.
\end{example}
\begin{example}
Consider the holomorphic effective Cartan geometries on \(M=\CP{1}\).
Since \(M\) is one dimensional (as a complex manifold), the curvature of any Cartan geometry on \(M\) vanishes, being a \(2\)-form valued in the adjoint bundle.
So the geometry is flat.
Since \(M\) is simply connected, the geometry arises by pullback of a local biholomorphism \(M\to X\) to the model \((X,G)\).
But since \(M\) is compact, this is a covering map.
Note that \(X\) is also a one dimensional complex manifold, so an oriented surface.
By classification of oriented surfaces, \(X=M=\CP{1}\), so the geometry is that of the model.
By the classification of complex homogeneous spaces in dimension one \cite{McKay2011b}, \(G=\PSL{2}\): there is a unique holomorphic effective Cartan geometry on \(\CP{1}\).
\end{example}
The moduli space of complete flat Cartan geometries, with a given connected model \((X,G)\), is thus the set of discrete subgroups \(\Gamma\subseteq G\) acting freely and properly on \(X\), modulo conjugation.
If we continuously vary, because \(\Gamma\) is discrete, its group structure does not vary, so the moduli space is an open set in the quotient space \(\Hom{\Gamma}{G}/G\), quotienting out by conjugation; see \cite{Goldman:2022} p. 165, \cite{Sikora:2012} for the theory of this quotient space.
\begin{corollary}\label{corollary:compact.flat}
If \((X,G)\) is a proper homogeneous space then every flat \((X,G)\)-geometry on any compact manifold is a complete Klein manifold of the associated universal covering homogeneous space.
\end{corollary}
\begin{proof}
By the Hopf--Rinow theorem \cite{Chavel:2006} \S I7, \cite{Gromov:2007} p. 9, \cite{Petersen:2016} p. 137 theorem 16, the metric on \(M\) is complete.
The stabilizer subgroup \(H:=G^{x_0}\) of any point \(x_0\in X\) is compact.
So \(M=M_c\) is its own scaffold.
So metric completeness of \(M\) implies developability, hence completeness, of the \((X,G)\)-geometry.
\end{proof}
\begin{example}
Consider a compact Riemannian manifold \(M\) locally isometric to \(\CP{n}\) with its usual metric.
It is a complete Klein manifold \(M=\Gamma\backslash\CP{n}\), some \(\Gamma\subset\PU{n+1}\).
Every element of \(\PU{n+1}\) acts on \(\CP{n}\) with a fixed point, by the spectral theorem from elementary linear algebra.
Hence \(\Gamma=\set{1}\) and \(M=\CP{n}\).
\end{example}
\begin{theorem}[Auslander--Markus]%
\define{theorem!Auslander--Markus}%
\define{Auslander--Markus theorem}
A flat affine connection on a connected manifold \(M\) is geodesically complete just when the developing map \(\tilde{M}\to\R^n\) to Euclidean space is an isomorphism of affine connections.
\end{theorem}
\begin{proof}
Clearly \(\R^n\) is simply connected, so only diffeomorphisms are covering maps.
\end{proof}

A group $\Gamma$ \emph{defies}\define{defiance} a group $G$ if every morphism $\Gamma\to G$ has finite image.  
\begin{example}
If $\Gamma$ is finite, or $G$ is finite then $\Gamma$ defies $G$.
\end{example}
\begin{theorem}\label{thm:compactToFlat}
For any connected homogeneous space \((X,G)\), with universal covering homogeneous space \((\tilde{X},\tilde{G})\), any flat \((X,G)\)-Cartan geometry, defined on a compact connected manifold $M$ with fundamental group defying $\tilde{G}$, is a complete Klein \((\tilde{X},\tilde{G})\)-geometry and \(X\) is compact with finite fundamental group.
\end{theorem}
\begin{proof}
We can assume \(X=\tilde{X}\) so \(G=\tilde{G}\). 
Write $M = \Gamma \backslash \tilde{M}$ for some pullback $\tilde{M}\to X$, so  $\Gamma = \fundamentalGroup{M}/\fundamentalGroup{\tilde{M}} \subseteq G$, and since $\fundamentalGroup{M}$ defies $G$, $\Gamma$ is finite, and $\tilde{M}$ is compact. 
The local diffeomorphism $\tilde{M} \to X$ is a covering map to its image. 
The bundles on which the Cartan connections live, say
  \[
  \begin{tikzcd}
    \G \arrow[d] & \tilde{\G} \arrow[l] \arrow[d] \arrow[r] & G \arrow[d] \\
    M        & \tilde{M} \arrow[l] \arrow[r] & X \\
  \end{tikzcd}
  \]
  are all pullbacks via covering maps, so completeness is preserved
  from $X$ and reflected to $M$.
\end{proof}
\begin{problem}{bo}
The conformal group of the sphere we have noted above; what is the conformal group of the standard flat conformal geometry on the real projective space of dimension \(n\ge 3\)?
What are all flat conformal geometries on that real projective space?
\end{problem}
\begin{corollary}
For any homogeneous space \((X,G)\), if $\dim X \ge 4$, then infinitely many compact manifolds of the same dimension as \(X\) bear no flat $(X,G)$-Cartan geometry.
\end{corollary}
\begin{proof}
  Construct manifolds with fundamental group defying $\tilde{G}$, following
  Massey \cite{Massey:1967} pp. 114--116. The fundamental group can be finite or
  infinite, as long as it has no quotient group belonging to $\tilde{G}$. For
  example, the fundamental group could be a free product of finitely presented
  simple groups not belonging to $\tilde{G}$.
\end{proof}

\subsection{Deformation of flat geometries}
The smooth deformations of a principal bundle on a fixed manifold \(M\), parameterized by an interval (or a connected manifold) are all isomorphic \cite{Steenrod:1999} section 11. 
Indeed, by definition, a ``family of bundles'' on \(M\) parameterized by a connected manifold \(S\) is just exactly a single bundle on \(M\times S\).
Picking a smooth connection on that bundle over \(M\times S\), parallel transport along a curve in \(S\) identifies any two of the bundles in the family.
The same is true topologically: continuous deformations of principal bundles are topologically isomorphic.
To see this: the transition maps of the bundle, valued in \(G\), are continuously varying, hence can be approximated by smoothly varying maps, arranged to agree at two points of \(S\) with the original bundles.
(This works for principal bundles, but not for fiber bundles: we can't always smooth homeomorphisms to diffeomorphisms.)
Therefore, in deforming Cartan geometries, we can consider the bundle to be fixed through the deformation.
\begin{theorem}%
[Ehresmann--Thurston--Weil 
\cite{Canary.Epstein.Green:2006} p. 16 Theorem 1.7.1, 
\cite{Goldman:2022} p. 161 Theorem 7.2.1]%
\define{theorem!Ehresmann--Thurston--Weil}
\define{Ehresmann--Thurston--Weil theorem}
Take a compact connected manifold \(M\) and a homogeneous space \((X,G)\).
Give the set of flat \((X,G)\)-geometries on \(M\) the \(C^1\) compact open topology.
Each flat \((X,G)\)-geometry on \(M\) has a neighborhood in which it is the unique flat \((X,G)\)-geometry with its holonomy, up to isomorphism.
Let \(\pi:=\fundamentalGroup{M}\).
The \emph{monodromy map}\define{monodromy map} associating to each flat \((X,G)\)-geometry its holonomy is an open continuous map to \(\operatorname{Hom}(\pi,G)\).
\end{theorem}
\begin{example}
In particular, if the character variety \(\operatorname{Hom}(\pi,G)/\Ad_G\) has positive dimension, then there are flat Cartan geometries with smoothly varying induced flat \(G\)-bundles, hence locally isomorphic but not globally.
For example, if \(M\) is the Riemann sphere punctured at three points, then the representation variety is \(\C^3\), so \(M\) has a \(3\)-parameter family of flat holomorphic  projective connections \cite{Goldman:2009}.
\end{example}
\begin{proof}
For any morphism \(\pi\to G\) of groups, \(\pi:=\fundamentalGroup{M}\), the flat connection on \(M\times G\) descends to a unique flat connection on \(\G_G:=\amal{M}{\pi}{G}\), hence our morphism arises as the holonomy of a flat connection.
Vary the morphism continuously: \(\pi\) is finitely generated, so continuously move the images in \(G\) of the generators, continuous in some parameter \(s\in \R\), to vary the flat connection, as we identify the bundle \(\G_{G,s}\) with a fixed one \(\G_{G,s_0}\).

In the \(C^1\) strong topology, there is an open set of connections on \(\G_G\) which do not vanish on the tangent spaces of \(\G_H\).
On a compact manifold, the \(C^1\) strong topology is just the \(C^1\) compact open topology.
So near any flat \((X,G)\)-geometry, there are flat \((X,G)\)-geometries giving rise to any nearby morphism \(\pi\to G\): our map is locally onto.

Take two flat \((X,G)\)-geometries on a compact connected manifold \(M\), with the same holonomy morphism \(\pi:=\fundamentalGroup{M}\to G\), with two developing maps \(\tilde{M}\to X\).
Let \(B:=\amal{\tilde{M}}{\pi}{X}\).
The fiber bundle \(X\to B\to M\) depends only on the holonomy, so the two geometries share this same fiber bundle.
Quotienting the graphs of the developing maps by \(\pi\) gives two sections of \(B\to M\).
If the two developing maps are \(C^1\) close on some large compact subset of \(\tilde{M}\), then these sections are \(C^1\) close in \(B\).
From another point of view: \(B=\G_G/H\) with section given by \(\G_H/H=M\), so if two geometries have \(H\)-bundles sitting close to one another in the same \(G\)-bundle, they have nearby sections in \(B\).

When \(\pi\) acts on \(\tilde{M}\times X\), the map to \(M\) is invariant.
The map to \(X\) is not, but the fibers transform among one another, so that they become the leaves of a foliation of \(B\).
The sections are transverse both to the map to \(M\) and to the leaves of the foliation, since the developing maps are local diffeomorphisms.

Pick one of the two \((X,G)\)-structures, i.e. one of the two sections.
Locally trivialize the foliation into a product along an open subset of the section.
\[
\begin{tikzpicture}
\foreach\i in {0,...,10}{
	\draw[gray!50] (0,{.1*\i}) -- (1,{.1*\i});
	\draw[gray!50] ({.1*\i},0) -- ({.1*\i},1);
}%
\draw [thick] (0,0) to (1,1);
\end{tikzpicture}
\]
Every \(C^1\) nearby section has a unique nearby point lying along each leaf of the foliation, on perhaps some slightly smaller open set.
\[
\begin{tikzpicture}
\foreach\i in {0,...,10}{
	\draw[gray!50] (0,{.1*\i}) -- (1,{.1*\i});
}%
\draw [thick] (0,.1) to [bend right] (1,.9);
\draw [thick] (0,0) to (1,1);
\end{tikzpicture}
\]
Map one section to the other, by identifying points that lie along the same leaf and inside one of our local trivializations.
For \(C^1\) close enough sections, this is a diffeomorphism which preserves the leaves, so identifies the developing map images in \(X\), in our local product trivializations.
Since this ``sliding along leaves'' is locally unique and smooth, it is globally so, for \(C^1\) close sections, and hence a diffeomorphism of \(M\). 
\end{proof}
\Danger{} This theorem is often misstated as saying that monodromy map takes the moduli space of \((X,G)\)-structures locally homeomorphically to the moduli space of representations of \(\pi\) in \(G\), but the monodromy map can be, for example, a branched double cover; see \cite{Canary.Epstein.Green:2006,Goldman:2022,Kapovich:1989}.
The monodromy map \emph{is} a local homeomorphism in the special case when \((X,G)=(\CP{1},\PSL{2})\) \cite{Earle:1981,Hubbard:1981}; see \cite{Loray/MarinPerez:2009} for more on these geometries.

\chapter{Automorphisms}
\section{Automorphisms of locally homogeneous structures}\label{subsection:X.G.auts}
For the moment, we just consider automorphisms of some simple examples.
\begin{example}
Take a group generated by four generators.
Draw the Cayley graph of the group with those generators, so each vertex has four edges.
Replace each edge with a tube, each vertex with a sphere punctured by a disk centered at each point of a Platonic tetrahedron, all disks of equal radius.
Glue together smoothly, using the same ``smoothing glue'' to put each edge in, preserving tetrahedral symmetry, to generate a surface whose isometry group contains the original group.
Conformally rescale to form a complete constant curvature metric, of curvature \(-1\) since there are at least two edges.
We find a complete hyperbolic surface, so a complete Cartan geometry modelled on the hyperbolic plane, with a very complicated automorphism group, containing our original group.
\end{example}
Take a strong effective homogeneous space \((X,G)\).
Suppose that \(U\subseteq X\) is an open subset.
Pick a point \(x_0\in X\).
Take \(\G\subseteq G\) to be the pullback bundle, i.e. the set of elements \(g\in G\) so that \(gx_0\in U\).
So \(\G\to U\) is a the bundle of a unique locally homogeneous structure: the one for which the inclusion of \(U\subseteq X\) is a chart.
Take \(\Gamma\subseteq G\) to be the subgroup of elements preserving \(U\subseteq X\).

Left action by any element of \(\Gamma\) preserves \(\G\): if \(\gamma\in\Gamma\) and \(g\in\G\) then 
\(\gamma(gx_0)\in U\), so \(\gamma g\in \G\), i.e. \(\Gamma\G=\G\).
Clearly \(\Gamma\subseteq G\) is a closed subgroup, hence a closed Lie subgroup  \cite{Mimura/Toda:1991} p. 44.
Left action preserves the Maurer--Cartan form as well, so \(\Gamma\) acts on the left as a group of automorphisms of the locally homogeneous structure \(U\).
\begin{theorem}\label{theorem:open.set.G.aut}
Take a strong effective homogeneous space \((X,G)\).
Take an connected open set \(U\subseteq X\).
The automorphism group of the induced \((X,G)\)-structure on \(U\) is precisely the subgroup of \(G\) preserving that open subset.
\end{theorem}
\begin{proof}
Any automorphism of \(\G\to U\) preserves the Maurer--Cartan form on the open set \(\G\subseteq G\); apply theorem~\vref{theorem:open.subset.G} to see that the automorphism is given, on each connected component of \(\G\), by a single element of \(G\).
The automorphism commutes with right action of \(H\) on \(\G\), so this element of \(G\) is constant throughout \(\G/H=U\).
\end{proof}

Right action by \(\Gamma\) might not preserve \(\G\), as we will shortly see.
Right action by an element \(g\in G\) transforms the Maurer--Cartan form in the adjoint representation, so preserves the Maurer--Cartan form just when \(g\) lies in the kernel of the adjoint representation, i.e. commutes with the identity component of \(G\).

The left invariant vector fields and the right invariant vector fields of \(G\) restrict to \(\G\).
Right invariant vector fields generate the left action.
So the flows of right invariant vector fields preserve both \(\G\) and the Maurer--Cartan form: they are infinitesimal automorphisms of the locally homogeneous structure.
The infinitesimal automorphisms from the Lie algebra of \(\Gamma\) are complete vector fields.
\begin{theorem}\label{theorem:open.set.G.inf.auts}
Take a strong effective homogeneous space \((X,G)\).
Take a connected open set \(U\subseteq X\).
The infinitesimal automorphisms of the induced \((X,G)\)-structure on \(U\) are precisely the restrictions of the right invariant vector fields on \(G\) to the preimage of \(G\) above \(U\).
Those which are complete are precisely the Lie algebra of the automorphism group.
\end{theorem}
\begin{proof}
Consider the preimage \(\G\subseteq G\) of \(U\subseteq X\).
Take an infinitesimal automorphism, i.e. a vector field on \(\G\) commuting with all left invariant vector fields and with the right \(H\)-action.
The left invariant vector fields generate the right action, so our vector field is invariant under the local right action.
Take the right invariant vector field which agrees with our vector field at some point, and subtract it.
So our vector field vanishes at some point.
But then it vanishes everywhere nearby, by its local right invariance.
The set of points at which our vector field vanishes is a union of components of \(\G\), and \(H\)-invariance, so the preimage in \(\G\) of a union of components of \(U\).
\end{proof}
Let us consider the left invariant vector fields, and the right action, in a particular example.
\begin{example}\label{example:hyperbolic.plane}
To be more specific, suppose that \(X=\Proj{2}\) is the real projective plane and \(G\) the group of projective linear transformations of \(X\).
Take \(U\subset X\) to be a disk inside the affine plane inside \(X\), and \(x_0\) the origin of the disk.
We don't need to know this, but the group \(\Gamma\) turns out to be the group of isometries of the hyperbolic plane metric on \(U\); indeed Hilbert computed the hyperbolic distance between points in terms of the projective plane geometry \cite{Goldman:2022} p. 69, \cite{Papadopoulos.Troyanov.2014}.
Hence an element \(g_0\in G\) satisfies \(g_0\G=\G\) precisely when \(g_0\in\Gamma\).

Take an element \(g_0\in G\) which moves the point \(x_0\), say to a point \(x_1:=g_0 x_0\ne x_0\).
A projective transformation \(g_1\in\G\) can move \(x_1\) anywhere, even outside of the disk \(U\), while moving \(x_0\) inside the disk \(U\), since projective transformations of the projective plane act transitively on pairs of distinct points.
So \(\G g_0=\G\) just when \(g_0\in H\), i.e. just when \(g_0\) fixes \(x_0\).
In other words, right action by an element \(g_0\in G\) is an automorphism just when \(g\) lies in \(H\) and commutes with the identity component of \(G\), hence fixes the \(x_0\)-component of \(X\).
So right translation of \(\G\) by some \(g\in G\) is an automorphism of the \((X,G)\)-structure on the hyperbolic plane just when \(g=1\).

Left invariant vector fields generate the right action.
The flow of a left invariant vector field from the Lie algebra of \(\Gamma\) is a right translation in \(G\), so preserves \(\G\) precisely when it is in \(\LieH\).
It does not preserve the left invariant Maurer--Cartan form unless it is the zero vector field.
\end{example}
\begin{theorem}
Take a strong effective homogeneous space \((X,G)\) and a subgroup \(\Gamma\subseteq G\).
Suppose that \(M\subseteq X\) is a \(\Gamma\)-invariant connected component of the free regular set.
Let \(\bar{M}:=\Gamma\backslash M\).
Let \(N_{\Gamma}G\) be the normalizer of \(\Gamma\) in \(G\).
Then \(N_{\Gamma}G\) leaves the free regular set invariant.
Let \(\Gamma'\subseteq N_{\Gamma}G\) be the set of elements of \(N_{\Gamma} G\) which leave \(M\) invariant.
Then \(\Gamma'\subseteq N_{\Gamma}G\) is a union of components of \(N_{\Gamma} G\).
The automorphism group \(\Aut[\bar{M}]\) of the \((X,G)\)-structure on \(\bar{M}\) is the quotient
\[
1\to\Gamma\to \Gamma'\to\Aut[\bar{M}]\to 1
\]
The identity component of \(\Gamma'\) is the identity component of the centralizer \(Z_{\Gamma} G\) of \(\Gamma\subseteq G\).
The Lie algebra of infinitesimal automorphisms is thus the Lie algebra centralizer
\[
\LieZ_{\Gamma}\LieG:=\set{A\in\LieG|\Ad_{\gamma}A=A\text{ for all }\gamma\in\Gamma}.
\]
The developing map is the composition
\[
\begin{tikzcd}
\tilde{M}\arrow[d]\arrow[dr,bend left]\\
M\arrow[d]\arrow[r]&X\\
\bar{M}
\end{tikzcd}
\]
given by the inclusion \(M\subseteq X\).
The holonomy morphism
\[
\pi_1(\bar{M})\to G
\]
factors through the inclusion maps \(\Gamma\to\Gamma'\to G\):
\[
\pi_1(\bar{M})\to \Gamma\cong\pi_1(\bar{M})/\pi_1(M)\to\Gamma'\to G.
\]
\end{theorem}
\begin{proof}
An open set \(U\subseteq X\) is a \(\Gamma\)-house for a point \(x\in X\) just when \(gU\) is a \(g\Gamma g^{-1}\)-house for the point \(gx\in X\).
Hence for \(g\) normalizing \(\Gamma\), the \(\Gamma\)-free regular set \(\Omega\subseteq X\) is \(g\)-invariant.
Even if \(\Gamma\subset G\) is not a Lie subgroup, \(N_{\Gamma}G\subseteq G\) is, by definition, a closed subgroup, hence an embedded Lie subgroup.
The identity component of \(N_{\Gamma} G\) fixes every component of the free regular set, so fixes \(M\), so belongs to \(\Gamma'\).

Every automorphism \(\varphi\in\Aut[\bar{M}]\) lifts to a diffeomorphism of \(M\), which is an automorphism of the \((X,G)\)-structure on the connected open set \(M\subseteq X\).
By theorem~\vref{theorem:open.set.G.aut}, the automorphism is given by a unique element of \(g\in G\) preserving \(M\subseteq X\).
Any two such lifts agree up to a deck transformation by an element of \(\Gamma\), replacing \(g\) by \(\gamma g\) for some \(\gamma\in\Gamma\).
Let \(\G\subseteq G\) be the preimage of \(M\subseteq X\) under the map
\[
g\in G\mapsto gx_0\in X.
\]
The lifted automorphism \(g\in G\) acts on \(\bar\G:=\Gamma\backslash\G\), i.e. \(g\Gamma g_0=\Gamma gg_0\) for any \(g_0\in\G\), i.e.\(g\Gamma=\Gamma g\), i.e. \(g\) normalizes \(\Gamma\) in \(G\).
The lifted automorphism also fixes \(M\), so belongs to \(\Gamma'\).

Since \(\Gamma'\) is a closed subgroup of \(N_{\Gamma}G\) containing the identity component, it is a union of components.
Hence it has the same Lie algebra.

Claim: \(\Gamma\subseteq G\) is discrete.
Take a convergent sequence \(\gamma_j\to\gamma_{\infty}\in\Gamma\).
We want to prove that it is constant, after finitely many terms.
Replacing it by \(\gamma_j\gamma_{\infty}^{-1}\), we can assume that \(\gamma_{\infty}=1\).
For any house \(U\),  \(U\cap\gamma_jU\) is empty unless \(\gamma_j=1\).
But \(\gamma_j\to 1\).
So\(\gamma_j=1\) after finitely many terms.
Hence \(\Gamma\) is discrete.

Take a smooth path \(g(t)\in N_{\Gamma}G\) with \(g(0)=1\).
So \(g(t)\in\Gamma'\).
For any fixed \(\gamma\in\Gamma\), \(g(t)\gamma g(t)^{-1}\in\Gamma\), a continuous path in \(\Gamma\), so constant.
Plugging in \(t=0\), this constant is \(\gamma\), so \(g(t)\gamma g(t)^{-1}=\gamma\) for all \(t\).
Hence the identity component of \(N_{\Gamma} G\) is the identity component of the centralizer \(Z_{\Gamma} G\).
\end{proof}
\begin{example}
Take affine space \((X,G)=(V,\GL{V}\ltimes V)\) for some finite dimensional vector space \(V\).
Pick a lattice \(\Lambda\subseteq V\).
Consider the affine structure on the torus \(T:=V/\Lambda\).
The free regular set is all of \(X=V\).
The fundamental group of \(T\) is \(\pi_1(T)=\Lambda\).
The normalizer \(N_{\Lambda}G\) is the set of all \(g=(h,v)\in G\) so that \(h\Lambda+v=\Lambda\).
Clearly this contains all elements \((0,v)\), and contains all \((h,0)\) so that \(h\Lambda=\Lambda\).
So \(g\) is the product \(g=(1,v)(h,0)\) of an element of \(\Lambda\) and one of the subgroup \(\GL{\Lambda}\subseteq\GL{V}\) preserving \(\Lambda\).
So the automorphism group of the affine structure on the torus \(T=V/\Lambda\) is the quotient
\[
\Aut[T]=(N_{\Lambda}G)/\Lambda\cong\GL{\Lambda}\ltimes T
\]
\end{example}
\begin{example}
Take projective space \((X,G):=(\Proj{n},\PGL{n+1})\) acted on by its group of projective transformations.
Take a vector space \(V\) of dimension \(n\) and  a lattice \(\Lambda\subset V\).
Consider \(V\) as affine space sitting in projective space \(X\), as the complement of the affine hyperplane at infiniity.
For each line in affine space, say containing some vector \(u\) and also some other vector \(u+v\), suppose that, in some basis of \(\Lambda\), all components of \(v\) are irrational.
We can approximate points of that line arbitrarily closely by elements of \(\Lambda\).
Hence the point of the affine hyperplane striking that line is a limit point of these elements \(\Lambda\), a point of discontinuity of the \(\Lambda\)-action on \(X\).
Therefore the free regular set of \(\Lambda\) on \(X\) is \(V\), which is a connected open set.
So the automorphism group of the projective structure on the torus \(T:=V/\Lambda\) is the quotient
\[
1\to\Lambda\to N_{\Lambda}G\to\Aut[T]\to 1.
\]
The normalizer sits in the subgroup of \(G\) preserving \(V\), which is precisely the group of affine transformations of \(V\).
So the projective automorphism group of the torus \(T\) is precisely the affine automorphism group, 
\[
\Aut[T]=\GL{\Lambda}\ltimes T.
\]
\end{example}
\begin{example}
The group 
\[
\Gamma:=\PSL{2,\Z}\subset G:=\PSL{2,\C}
\]
is generated by the two matrices \cite{Conrad:2024}
\[
\begin{bmatrix}
0&-1\\
1&0
\end{bmatrix},
\begin{bmatrix}
1&1\\
0&1
\end{bmatrix}.
\]
They act on the complex projective line \(X\) by \(z\mapsto -1/z\) and \(z\mapsto z+1\).
On the real number line, their orbit through any real number is dense.
Hence the \(\Gamma\)-free regular set lies in the upper half plane (also known as the \emph{hyperbolic plane})\SubIndex{hyperbolic space}
\[
X_0:=\set{z\in\C|z=x+iy, y>0}.
\]
It turns out that this is the \(\Gamma\)-free regular set.\SubIndex{free regular set}
The subgroup \(G_0\subset G\) preserving \(X_0\) is precisely \(\PSL{2,\R}\).
Hence the automorphism group of the complex projective structure on \(X_0\) is precisely the group \(G_0\), which is the group of orientation preserving isometries of the hyperbolic metric
\[
\frac{dx^2+dy^2}{y^2}
\]
on \(X_0\).
So \(\Gamma\) also preserves the orientation and the hyperbolic metric on \(X_0\).
It turns out that the action of \(\Gamma\) on \(X_0\) is a covering action  \cite{Conrad:2024}.
A suitable fundamental domain is taken by these two matrices to a tiling of \(X_0\):
\[
\includegraphics[width=6cm]{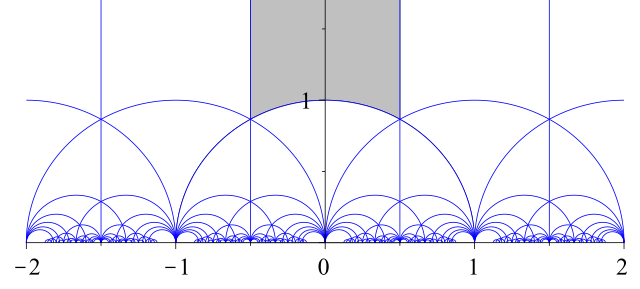}
\]
Take a real \(2\times 2\) traceless matrix \(A\).
To commute with \(\Gamma\), it must commute with those two matrices, which the reader can easily check forces it to vanish.
But suppose more generally that we take some group \(\Gamma\subseteq G\) with free regular set \(X_0\) and so that the only traceless \(2\times 2\) matrix \(A\) commuting with all elements of \(\Gamma\) is \(0\).
Hence the centralizer Lie algebra is trivial.
Consider the Riemann surface \(\bar{M}:=\Gamma\backslash X_0\).
The automorphism group of the complex projective structure on the Riemann surface \(\bar{M}\) is a discrete group of orientation preserving hyperbolic isometries.
If \(\bar{M}\) is compact then this group is finite, as the isometries of a compact Riemannian manifold form a compact Lie group.
\end{example}

\section{Infinitesimal automorphisms}
An \emph{infinitesimal automorphism}\define{infinitesimal automorphism}\define{automorphism!infinitesimal} of a Cartan geometry \(H\to \G\to M\) is an \(H\)-invariant vector field on \(\G\) whose flow preserves the Cartan connection.
\begin{problem}{bp}
This vector field projects to a vector field on \(M\).
\end{problem}
\begin{example}
On the model \((X,G)\)-Cartan geometry, the left invariant vector fields on \(G\) are the constant vector fields, while the right invariant vector fields on \(G\) are the infinitesimal automorphisms.
\end{example}
Clearly every infinitesimal automorphism \(Z\) commutes with every constant vector field.
But then \(Z\) is invariant under the flow of every constant vector field.
Since the constant vector fields point in all directions, \(Z\) is carried by them through an open set of \(\G\).
\begin{proposition}\label{prop:complete.complete}
Every infinitesimal automorphism of a vast Cartan geometry is a complete vector field.
\end{proposition}
Kobayashi \cite{Kobayashi:1957} p. 190 proved this for complete Cartan geometries.
\begin{proof}
At each point, consider the time for which the flow of an infinitesimal automorphism \(Z\) is defined (forward or backward in time).
This time remains constant along the flow of the constant vector fields, so is constant on their orbits.
But if finite, that time diminishes along the flow of \(Z\).
Hence that time is infinite.
\end{proof}
Since \(Z\) also commutes with the \(H\)-action, we can travel around from any point of \(\G\) to any other, assuming \(M\) is connected, determining our vector field \(Z\) at all points from its value near any one point.
Moreover, \(Z\) is determined near our initial point just by taking its value at one point of \(\G\) and keeping it constant along the flows of the constant vector fields and extending by \(H\)-invariance:
\begin{theorem}%
[Amores I \cite{Amores:1979}]%
\label{theorem:Lie.alg}
\define{theorem!Amores I}
\define{Amores!theorem I}
Two infinitesimal automorphisms of a Cartan geometry \(\G\to M\), with \(M\) connected, which agree at one point of \(\G\), agree everywhere.
The Lie algebra of infinitesimal automorphisms of an \((X,G)\)-Cartan geometry has dimension at most that of \(G\).
\end{theorem}
Cartan \cite{Cartan:1910} p. 150, p. 157 knew this; he explicitly computes equations for the infinitesimal automorphisms, and points out that when he has enough \(1\)-forms invariantly defined on \(\G\) to span every cotangent space, then the symmetry group is of dimension at most that of \(\G\), and equal only when the curvature is constant. 
\subsection{Two notions of local}
A \emph{local infinitesimal automorphism}\define{local infinitesimal automorphism}\define{automorphism!local infinitesimal}\define{infinitesimal automorphism!local} of a Cartan geometry \(\G\to M\) is an infinitesimal automorphism of some open set of \(M\).
There is a different definition that gives essentially the same theory: a \emph{microlocal automorphism}\define{microlocal!automorphism}\define{automorphism!microlocal} of a Cartan geometry \(\G\to M\) is a bijection \(U\xrightarrow{\varphi}W\) of open sets \(U,W\subseteq \G\) which commutes with the flows of constant vector fields where defined, so that for any point \(m_0\) and components \(V_0,V_1\) of \(\G_{m_0}\cap U\), there are points \(p_0\in V_0\), \(p_0h\in V_1\) for some \(h\in H\) so that \(\varphi(p_0h)=\varphi(p_0)h\).
(If all intersections with fibers are connected, clearly this holds.)
It follows that such a map is smooth (and real analytic if the geometry is real analytic), since the constant vector fields give coordinates via their flow chart, as in the Orbit Theorem (theorem~\vref{thm:orbit}).
(In particular, if \(U\) has connected intersection with every fiber of \(\G\to M\), then we only need ask that \(\varphi\) commute with the constant vector fields.)
Similarly, a \emph{microlocal infinitesimal automorphism}\define{infinitesimal automorphism!microlocal}\define{automorphism!infinitesimal!microlocal} of a Cartan geometry \(\G\to M\) is a vector field on an open subset \(U\subseteq\G\) whose flow commutes with the flows of constant vector fields where defined and so that for any point \(m_0\) and components \(V_0,V_1\) of \(\G_{m_0}\cap U\), there are points \(p_0\in V_0\), \(p_0h\in V_1\) for some \(h\in H\) so that \(Z(p_0h)=\RT{h}'(p_0)Z(p_0)\).
\begin{lemma}[Microlocal is local]
Take a microlocal automorphism \(\varphi\) of a Cartan geometry \(\G\xrightarrow{\pi}M\) defined on an open subset \(U\subseteq\G\).
It extends to a unique local automorphism on \(\pi^{-1}\pi U\).
Similarly, take a microlocal infinitesimal automorphism \(Z\) defined on an open subset \(U\subseteq\G\).
It extends uniquely to a local infinitesimal automorphism on \(\pi^{-1}\pi U\).
\end{lemma}
After this lemma, we forget about the term \emph{microlocal}, using this lemma to deliberately conflate the two concepts.
\begin{proof}
By definition, a microlocal automorphism commutes with the flows up the fibers of \(\G\to M\), which are copies of \(H\), so is locally a left translation from one fiber to the other, in any local trivialization.
Hence it extends globally to such a translation just when we can get it to agree from one component to the other on which element it translates by.
Hence it extends to be \(H\)-invariant.
Similarly for a microlocal infinitesimal automorphism.
\end{proof}
\subsection{Transitivity}
\begin{corollary}\label{corollary:span.inf.aut}
If the automorphisms of a Cartan geometry \(\G\to M\) permute the components of \(M\) and the (local) infinitesimal automorphisms span the tangent space of some point (a dense set of points) of \(\G\) then the curvature is constant, an \(H\)-invariant element of \(\LieG\otimes\Lm*{2}{\LieG/\LieH}^*\).
\end{corollary}
\begin{example}
For a conformal structure, \((X,G)=(S^n,\PO{n+1,1})\), it is not difficult algebra to check that
\[
\left(\LieG\otimes\Lm*{2}{\LieG/\LieH}^*\right)^H=0.
\]
Infinitesimal automorphism algebra of the same dimension as \(G\) forces flatness.
\end{example}
\begin{example}
For a Riemannian geometry, \((X,G)=(\R^n,\Orth{n}\ltimes\R^n)\), it is not difficult algebra to check that 
\[
\left(\LieG\otimes\Lm*{2}{\LieG/\LieH}^*\right)^H
\]
is \(1\)-dimensional: infinitesimal automorphism algebra of the same dimension as \(G\) forces local isometry to a Euclidean space, sphere or hyperbolic space.
\end{example}
\subsection{Multivalued infinitesimal automorphisms}
Take a Cartan geometry \(H\to\G\to M\).
The infinitesimal automorphism vector fields defined near a point \(p_0\in \G\) form a linear subspace in \(T_{p_0} \G\); the ones which extend globally form a smaller linear subspace, of constant dimension over each component of \(M\).
The \emph{orthodox set}\define{orthodox set} is the set of points of \(\G\) near which this first linear subspace has locally constant dimension, or equivalently has locally maximal dimension.
Clearly the orthodox set is open.
Since the orthodox set is \(H\)-invariant, it is the preimage of an open subset of \(M\), also called the \emph{orthodox set}.
This linear subspace in each tangent space forms a vector subbundle over each component of the orthodox set of \(M\).
Above any component \(U\) of the orthodox set, each vector in that linear subspace locally extends uniquely to a local infinitesimal automorphism.
So this vector bundle has a flat connection.
Hence by continuation it extends to a global infinitesimal automorphism on the universal covering space \(\tilde{U}\to U\).

A Cartan geometry is \emph{orthodox}\define{orthodox!Cartan geometry}\define{Cartan!geometry!orthodox} if the linear subspace is of constant dimension.
(Gromov uses \emph{regular}\SubIndex{Cartan!geometry!regular}\SubIndex{regular!Cartan geometry} for \emph{orthodox}, but the term \emph{regular} is in use already in the theory of Cartan geometries.)
\begin{example}
Any real analytic Cartan geometry, on a connected manifold is orthodox: both linear subspaces are carried invariantly by the constant vector fields.
\end{example}
In a smooth Cartan geometry, the first vector space dimension can drop suddenly as we approach a point of \(\G\).
\begin{example}
Take a flat torus and put a little bump on it, an unorthodox Riemannian geometry.
\end{example}
For each infinitesimal automorphism vector field on \(\G\), its projection to \(M\) has zeroes arising from points of \(\G\) where it is tangent to the vertical, and hence this linear subspace is tangent to the vertical.
So in a real analytic geometry, on any connected manifold, there is a dense open \(H\)-invariant subset of \(\G\) on which the linear projection of that linear subspace achieves maximal rank.
Quotienting by \(H\), there is a dense open subset of \(M\) on which the orbits of the local (or the global) infinitesimal automorphisms form a foliation, by the Orbit Theorem (theorem~\vref{thm:orbit}).
Even away from that subset each orbit is a smooth submanifold of \(M\), as is the case for any collection of vector fields, again by the Orbit Theorem.

Any orthodox Cartan geometry admits a unique minimal normal covering space \(\hat{M}\to M\) for which the local infinitesimal automorphisms are not multivalued.
The Galois group of the covering is the \emph{holonomy group}\define{holonomy!group!of infinitesimal automorphisms} of infinitesimal automorphsms.
\begin{example}
Real projective space inherits its conformal structure from the sphere by antipodal quotient.
But the infinitesimal automorphisms of the conformal geometry of the sphere are larger: picture the sphere as Euclidean space with a point added at infinity.
The dilations of Euclidean space are conformal, but don't commute with the antipodal map on the sphere, i.e. the map \(x\mapsto x/|x|^2\), which generates the holonomy group of infinitesimal automorphisms of real projective space.
\end{example}
\begin{problem}{bq}
Find the infinitesimal automorphisms of the conformal geometry on the sphere and the real projective space of dimension \(3\) or more.
\end{problem}
\begin{problem}{br}
What does this holonomy group have to do with the holonomy of the developing map, if the Cartan geometry is flat?
\end{problem}
\begin{problem}{bs}
Prove that the infinitesimal automorphisms on \(M\) are precisely those on \(\hat{M}\) invariant under the holonomy.
\end{problem}
\begin{theorem}%
[Amores II \cite{Amores:1979}]%
\label{theorem:analytic.extension}%
\define{theorem!Amores II}
\define{Amores!theorem II}
Take a real analytic Cartan geometry \(\G\to M\) on a connected manifold \(M\). 
Take a microlocal infinitesimal automorphism \(Z\) defined on a connected open subset \(U\subseteq\G\).
If \(M\) is simply connected, then \(Z\)  extends uniquely to an infinitesimal automorphism defined over \(M\).
If \(M\) is not simply connected, \(Z\) extends over the universal covering space of \(M\), descending to an infinitesimal automorphism defined over \(M\) just when \(Z\) is invariant under the holonomy of infinitesimal automorphisms.
\end{theorem}
\begin{proof}
We can suppose that \(M\) is simply connected.
Moving by flows of constant vector fields and by \(H\)-action, we cover \(\G\) in open sets on each of which we have defined some vector field.
Analyticity ensures that, since \(Z\) is invariant under all constant vector fields initially (an analytic equation), this remains true as we extend the domain of \(Z\).
Micolocality is similarly preserved by analyticity and permutation of the constant vector fields under \(H\)-action.
Hence each vector field extends uniquely to a local infinitesimal automorphism.
At each step in the process, by invariance under the flows, these local infinitesimal automorphisms agree on overlaps.
Our vector field \(Z\) remains \(H\)-equivariant as we extend it, again by analyticity, so is a section of \(T\G/H\to M\), a vector bundle on a simply connected manifold, so it does not become multivalued.
\end{proof}
\begin{example}
The Clifton plane is simply connected and real analytic, so any local automorphism or infinitesimal automorphism extends uniquely globally, and is real analytic, and the infinitesimal automorphisms are complete vector fields.
\end{example}
For each \(A\in\LieG\), let \(\bar{A}:=A+\LieH\in\LieG/\LieH\).
We can thus write the soldering form of a Cartan connection \(\omega\) as \(\bar\omega\) rather than \(\sigma\).
For any vector field \(Z\) on \(\G\),  let \(A:=Z\hook\omega\).
Conversely for any function \(\G\xrightarrow{A}\LieG\), there is a unique associated vector field \(Z\) on \(\G\) defined by \(Z\hook\omega=A\).
\begin{problem}{infinitesimal.aut.eq	uation}
Prove that, for any function \(\G\xrightarrow{A}\LieG\), the associated vector field \(Z\) with \(Z\hook\omega=A\) preserves the Cartan connection (and hence is \(\LieH\)-invariant) if and only if 
\[
dA = \lb{A}{\omega} + k \bar{A}\bar\omega.
\]
Furthermore, \(Z\) is an infinitesimal automorphism just when \(A\) satisfies this differential equation and \(\RT{h}^*A=\Ad_h^{-1}A\) for at least one \(h\) in each path component of \(H\). 
\end{problem}
\begin{answer}{infinitesimal.aut.eq	uation}
\begin{align*}
\LieDer_Z \omega
&=
d \left(Z \hook \omega \right) + Z \hook d \omega
\\
&=
dA + Z \hook \left( -\frac{1}{2}\lb{\omega}{\omega} + \frac{1}{2}k\sigma\wedge\sigma \right)
\\
&=dA-\lb{A}{\omega} + k \bar{A} \sigma.
\end{align*}
\end{answer}
Note that this tells us once again that \(Z\) is determined by its value at a point: we have a total differential equation for \(A\).
Any $H$-equivariant vector field $Z$ on $\G$ descends to a vector field $\bar{Z}$ on $M$, defined by $\bar{Z}(m) = \pi'(p) Z(p)$ for $p \in \G$ and $m = \pi(p)$. 
Just as $Z$ is associated to $A$ (i.e. $Z \hook \omega = A$), this vector field $\bar{Z}$ is associated to $\bar{A}$.

Take two infinitesimal automorphisms \(V,W\) with associated maps \(\G\xrightarrow{A,B}\LieG\).
Suppose that \(C:=[V,W]\hook\omega\). 
The Cartan equation
\[
d\omega(V,W)=\LieDer_V(W\hook\omega)-\LieDer_W(V\hook\omega)-[V,W]\hook\omega
\]
expands to
\[
C=\lb{A}{B}-k(\bar{A},\bar{B})+\LieDer_V B-\LieDer_W A,
\]
so the Lie algebra of the infinitesimal automorphisms is given by the \emph{curvature deformed bracket},\define{curvature!deformed bracket} with the extra terms we did see for constant vector fields because \(A,B\) are not constant.
From this, we see that the underlying vector field on \(M\) is represented by
\[
\bar{C}=\overline{\lb{A}{B}}-\bar{k}(\bar{A},\bar{B})+\LieDer_V \bar{B}-\LieDer_W \bar{A}.
\]
\begin{problem}{bt}
If an infinitesimal automorphism \(Z\) has associated projection \(\bar{Z}\) and \(\bar{Z}(m_0)=0\) at some point \(m_0\in M\), prove that the linearization of \(\bar{Z}\) at \(m_0\) is the linearization of some right invariant vector field on the homogeneous model.
\end{problem}
\begin{proposition}
On any connected manifold \(M\), the infinitesimal automorphisms of any flat Cartan geometry, say with model \((X,G)\), are identified by the developing map with the holonomy invariant Lie subalgebra of the Lie algebra of \(G\).
\end{proposition}
\begin{proof}
Take an infinitesimal automorphism, lift to the geometry on \(\tilde{M}\).
Locally this is identified by the developing map with an infinitesimal automorphism of the model, some element of \(\LieG\).
Once we match our infinitesimal automorphism with some element of \(\LieG\) near some point of \(\tilde\G\), we continue to match in some open set, by flowing along constant vector fields.
We flow to the entire component of that point.
Extend by \(H\)-invariance to all components.
So the infinitesimal automorphisms form a Lie subalgebra of \(\LieG\), invariant under the fundamental group.
The converse is clear.
\end{proof}
\begin{example}
The conformal geometry of a flat torus of dimension \(3\) or more has developing map identifying Euclidean space with a punctured sphere.
Infinitesimal automorphisms are vector fields on the sphere vanishing at the puncture, and invariant under the cocompact group action of the fundamental group of the torus.
But the model is algebraic, so these vector fields are invariant under the Zariski closure of the cocompact group action, i.e. under the translations of Euclidean space, hence are themselves translations.
The torus has conformal group precisely the torus acting on itself by translation.
\end{example}
\subsection{The unorthodox locus}
In a real analytic Cartan geometry, all points are orthodox.
In a smooth Cartan geometry, there can be infinitely many unorthodox points, although the unorthodox points are nowhere dense.
We expect the set of unorthodox points to be something like a hypersurface. 
At least, we can say that it is never ``small'':
\begin{theorem}[Gromov \cite{Gromov1988} p. 88]
Take a Cartan geometry \(H\to\G\to M\) and an open set \(U\subseteq M\).
Suppose that \(U\) contains a dense subset \(U'\subseteq U\) of orthodox points with trivial holonomy of local infinitesimal automorphisms.
Then every point of \(U\) is orthodox.
\end{theorem}
\begin{proof}
Take a point \(p_0\in\G_{U'}\).
Take an open set \(B\subseteq\G_U\) containing \(p_0\).
We need to prove that \(B\) is orthodox.
We can assume that \(B\) is a coordinate ball.

Claim: \(U'\cap B\subseteq B\) is dense.
Proof: take a point \(p\in B\), and an open subset of \(B\) containing it.
Then \(p\in\bar{U}'=U\), so \(p\) is the limit of a sequence of points of \(U'\).
That sequence eventually enters that open set.

Pick an infinitesimal automorphism defined near \(p_0\).
Take any path \(p(t)\) starting at \(p_0\) and lying in \(B\).
Along that path, we can extend our infinitesimal automorphism to a vector field by solving the differential equation
\[
dA = \lb{A}{\omega} + k \bar{A}\bar\omega
\]
but just along that path.
We can make such paths radially from \(p_0\), since \(B\) is a ball, to extend the vector field smoothly to \(B\).

On \(U'\cap B\), this vector field has to agree with the infinitesimal automorphism, by density.
Hence it is an infinitesimal automorphism on \(B\).
It is determined at each point \(p_1\in B\) by its value \(A(p_1)\in\LieG\).
So the linear subspace of these values in \(\LieG\) has dimension at least as large as that of \(p_0\).
\end{proof}
\begin{corollary}
Take a Cartan geometry \(H\to\G\to M\) and an open set \(U\subseteq M\).
Suppose that the intersection of \(U\) with the unorthodox set of \(M\) sits in a connected hypersurface with nonempty boundary, or in a submanifold of codimension \(2\) or more.
Then that intersection is empty.
\end{corollary}
\begin{proof}
Any higher codimension submanifold sits locally in a hypersurface with nonempty boundary, so it is enough to prove the result for a hypersurface.

Every point of that connected hypersurface maps to any other by a diffeomorphism of \(M\).
So we can assume our point lies in a ball \(B\subseteq U\), and replace \(U\) by \(B\).
We can assume that, in some coordinates, our hypersurface is a half hyperplane, so its complement \(U'\) is a connected and simply connected open subset of \(U\).
\end{proof}
There is a different type of ``singularity'' that we can't get rid of: local infinitesimal automorphisms could be more plentiful at certain points of the underlying manifold \(M\).
\begin{example}
Rotate an ellipsoid of revolution, preserving its metric and orientation.
\[
\includegraphics[width=2cm]{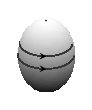}
\]
The two poles are fixed points.
So at these points, a local infinitesimal automorphism is not determined by its zero jet in the tangent space of the ellipsoid.
\end{example}
\subsection{From base to bundle}
Recall that the kernel of \((X,G)\) is the largest normal subgroup of \(G\) contained in the stabilizer \(H:=G^{x_0}\); the kernel has Lie algebra the largest ideal of \(\LieG\) contained in \(\LieH\), trivial just when \((X,G)\) is almost effective.
\begin{proposition}[Sharpe \cite{Sharpe:1997}]
Every infinitesimal automorphism of a Cartan geometry with almost effective model is determined by its projection to the underlying base manifold. 
\end{proposition}
\begin{proof}
If two have the same projection, their difference is an infinitesimal automorphism \(Z\) tangent to the fibers of \(\G\to M\).
Let \(A:=Z\hook\omega\); tangency to the fibers is precisely \(\G\xrightarrow{A}\LieH\), so \(\bar{A}=Z\hook\sigma=0\).
Since \(Z\) is an infinitesimal automorphism,
\[
dA = \lb{A}{\omega} + k\bar{A}\bar{\omega}=\lb{A}{\omega}.
\]
But \(A\) is valued in \(\LieH\), so \(dA\) is a \(1\)-form valued in \(\LieH\), while \(\omega\) is onto \(\LieG\).
Hence \(\lb{A}{\LieG}\subseteq\LieH\); so \(A\) is valued in the subalgebra of \(\LieH\) satisfying this equation.
Repeating the argument, force \(A\) to be valued in successively smaller subalgebras, eventually in the largest ideal of \(\LieG\) contained in \(\LieH\), as in lemma~\vref{lemma:eff.hom}.
\end{proof}

\section{Automorphism groups}\label{section:aut.gps}
An \emph{isomorphism}\define{isomorphism!of Cartan geometries}\define{Cartan!geometry!isomorphism} (automorphism) of Cartan geometries is a smooth bundle isomorphism (automorphism) preserving the Cartan connection.
Equivalently, it is a bijection (not assumed even continuous) of the total space \(\G\) commuting with the structure group action and with the flows of constant vector fields (where these are defined).
\begin{problem}{bu}
Take two finite dimensional real vector spaces \(V,W\) and a lattice \(\Lambda\subseteq V\), i.e. the set of integer linear combinations of some basis vectors.
Let \(T:=V/\Lambda\).
The lattice acts as automorphisms of the affine structure of \(V\times W\), and this structure descends to \(T\times W\).
Find the automorphism group and all of the infinitesimal automorphisms.
Give some examples of orbits of the one parameter subgroups of automorphisms.
For example, take a densely winding irrational rotation inside \(T\) and the exponential of a linear map of \(W\) with both positive and negative eigenvalues:
\[

\begin{tikzpicture}[scale=.3]
\draw[domain=0:1440,smooth,samples=100,black,very thick,variable=\t]
plot ({sin(\t)},0,{cos(\t)});
\draw[domain=0:1440,smooth,samples=150,gray!50,thick,variable=\t]
plot ({sin(\t)},{1/(\t/360+.1)},{cos(\t)});
\end{tikzpicture}

\]
This one parameter subgroup has orbits that foliate the complement \(T\times W-T\times\set{0}\) by curves.
Some of those curves are embedded copies of \(\R\), while others tend in one direction toward all points of \(T\times\set{0}\).
So while the one parameter subgroup is only \(1\)-dimensional, the limit points of one of its orbits can have dimension \(\dim T\), as large as we like.
\end{problem}
\begin{example}
The Riemannian geometry of a countably infinite disjoint union of circles, all of the same length, has automorphism group extending the group of permutations of the circles (an uncountable group) by the product of the group of rotations and reflections of each circle: an infinite dimensional automorphism group with uncountably many components.
If instead all of the circles have different lengths, the automorphism group is still infinite dimensional.
If the automorphism group is to be a Lie group, we must to restrict to manifolds with finitely many components.
Even if we just have two circles, of different radii, the automorphism group is a product, not acting freely on the circles.
\end{example}
\begin{example}
Take a manifold \(M\) with finitely many components.
The automorphism group of any Cartan geometry is the obvious extension of (i) the product of the automorphism groups of each of the components with (ii) the permutations of those components which are isomorphic.
Henceforth we only consider Cartan geometries on connected manifolds.
\end{example}
\begin{example}
There are complete Riemannian geometries on surfaces whose automorphism group is one dimensional with infinitely many components, or is discrete with infinitely many components.
\end{example}
{\centering\tiny
\begin{tabular}{p{5cm}p{5cm}}
\includegraphics[width=5cm]{delaunay.png} & \includegraphics[width=5cm]{Schwarz_P_surface} \\
GeometrieWerkstatt CC BY-NC-SA-3.0 &
By Anders Sandberg \\ & Own work \\ & CC BY-SA 3.0 \\
& {\verb!commons.wikimedia.org/w/index.php?curid=20737176!}
\end{tabular}
\par
}

\bigskip

It is unknown whether the automorphism group of any homogeneous Cartan geometry on any connected and simply connected manifold has finitely many components.
\begin{example}
Rotate an ellipsoid of revolution, preserving its metric and orientation.
\[
\includegraphics[width=2cm]{ellipsoid}
\]
The two poles are fixed points.
The oriented Riemannian geometry is a Cartan geometry modelled on the oriented Euclidean plane.
The automorphisms of the Cartan geometry are precisely the rotations fixing the poles and their composition with the rotation interchanging the poles, turning it upside down.
When that same circle action acts on unit tangent vectors, it has no fixed points: either the tangent vector sits at a point which is not the north or south pole, and the rotation moves the point of the ellipsoid where the tangent vector sits, or the tangent vector sits at the north or south pole, and the circle action rotates the tangent space of the ellipsoid.
Turning the ellipsoid upside down fixes the points of the equator, but doesn't fix any unit tangent vectors, as they reverse direction.
Turn each unit tangent vector by a right angle to give another unit tangent vector, an orthonormal basis.
Hence the unit tangent bundle is identified with the orthonormal frame bundle.
The orthonormal frame bundle is the bundle of the Cartan geometry.
So automorphisms of the Cartan geometry act freely on the bundle, a principal circle action together with the rotation upside down.
Each orbit of the automorphism group in the total space consists of a unit tangent vector at some point, but then moved around by rotating the ellipsoid around its axis, and by turning it upside down: a pair of disjoint circles in the unit tangent bundle.
\end{example}
\begin{problem}{bv}
What is the quotient space by the circle action?
\end{problem}
\begin{problem}{bw}
What are the orbits of the automorphism group of the standard oriented Riemannian metric on the plane, the hyperbolic plane, the sphere, the cylinder?
What are the quotient spaces?
\end{problem}
Consider a product of finitely many Cartan geometries.
Clearly the product of the automorphism groups sits in the automorphism group of the product, as each element commutes with the various constant vector fields and with the structure group.
\begin{problem}{bx}
Prove that the automorphism group of the product of finitely many Cartan geometries is the product of the automorphism groups.
\end{problem}
\begin{example}
On the product of two unit 2-spheres, the Cartan geometry associated to their Riemannian geometries is \emph{not} the product Riemannian geometry.
The Cartan geometry keeps track of the order of the product and has structure group \(\Orth{2}\times\Orth{2}\) while the Riemannian geometry has the permutation of the factors as an isometry, and has structure group \(\Orth{4}\).
\end{example}
\begin{example}
Any group of diffeomorphisms of a manifold, endowed with the discrete topology, acts smoothly as a Lie group, perhaps with uncountably many components.
For instance, the real number line with the discrete topology acts on the real number line with the usual topology, by addition, preserving the translation structure, i.e. the \((X,G)=(\R,\R)\)-structure, acting transitively and smoothly as a Lie group with uncountably many components, but not properly.
So we have to be careful to put the coarsest topology on the automorphism group which will still render it a Lie group acting smoothly.
\end{example}
Our aim in this section is to prove:
\begin{theorem}\label{thm:strong}
Suppose that \(M\) is a connected manifold.
Take a Cartan geometry \(\G\to M\).
The set \(\Aut\) of Cartan geometry automorphisms is a Lie group for a unique Lie group structure for which \(\Aut\to\G\to\G/\Aut\) is a smooth principal bundle.
The Lie algebra of \(\Aut\) is the set of complete infinitesimal automorphisms.
Automorphisms are isometries on the scaffold.
The constant vector fields descend to vector fields on \(\G/\Aut\), spanning every tangent space.
The structure group action descends to an action on \(\G/\Aut\).
\end{theorem}
\begin{corollary}\label{corollary:finite.components}
On any manifold \(M\) with finitely many components, the automorphism group of any Cartan geometry \(\G\to M\) is a Lie group acting smoothly on \(\G\) and \(M\).
\end{corollary}
\begin{corollary}
If \(M\) is a compact manifold with a Cartan geometry modelled on a proper homogeneous space \((X,G)\) then  the automorphism group \(\Aut\) of that Cartan geometry is a compact Lie group.
\end{corollary}
\begin{proof}
Since \(M\) is compact, it has finitely many components.
The subgroup of the automorphism group which preserves all connected components has finite index, so it suffices to prove that this group is compact.
It is the product of the automorphism groups of the connected components.
So we can assume that \(M\) is connected. 
The stabilizer \(H:=G^{x_0}\subseteq G\) of any point \(x_0\in X\) is compact.
The total space \(\G\to M\) of the Cartan geometry has compact base space \(M\) and compact structure group \(H\), so is compact.
So the base of any bundle with total space \(\G\) is compact, and hence the fibers are compact.
The fibers of \(\Aut\to\G\to\G/\Aut\) are therefore compact.
\end{proof}
\begin{example}
The isometry group of any compact Riemannian manifold is a compact Lie group acting smoothly.
If the manifold is connected, then the quotient of the orthonormal frame bundle by that isometry group is a smooth compact manifold.
\end{example}
\subsection{Elementary properties}
\begin{theorem}[Cartan \cite{Cartan:1910,Cartan:161} p. 122]
Any local isomorphism of Cartan geometries 
\[
\begin{tikzcd}
\G\arrow[d]\arrow[r]&\G'\arrow[d]\\
M\arrow[r]&M'
\end{tikzcd}
\]
on a connected manifold \(M\) is determined by how it maps any one point of \(\G\).
\end{theorem}
\begin{proof}
The map \(\G\to\G'\) is \(H\)-equivariant and permutes the constant vector fields, so takes their flows line to one another.
The constant vector fields span every tangent space of \(\G\) and \(\G'\), so their flows move our one point of \(\G\) through a component of \(\G\); \(H\)-equivariance takes us through all other components above a component of \(M\).
\end{proof}
\begin{corollary}\label{corollary:free}
On any connected manifold \(M\), the automorphism group of any Cartan geometry \(\G\to M\) acts freely on \(\G\).
\end{corollary}
\begin{problem}{by}
Prove that the automorphism group of the model \((X,G)\) is \(G\) if \(X\) is connected.
Find an example of a homogeneous space \((X,G)\) for which the automorphism group of the model geometry is not \(G\).
\end{problem}
\begin{problem}{bz}
Suppose that \((X,G)\) is a strong effective homogeneous space.
Prove that every isomorphism between connected open subsets of \(X\) is a unique element of \(G\), so extends to \(X\).
In particular the automorphism group of an open subset of \(X\) is the subgroup of \(G\) preserving that open set.
\end{problem}
\begin{example}
Take \((X,G)=(\R^n\smallsetminus 0,\GL{n})\) and \(M\subset X\) the upper half space \(x_1>0\).
So \(\Aut[M]\) is the group of linear transformations preserving that half space: the matrices 
\[
\begin{pmatrix}
a&0\\
b&c
\end{pmatrix}
\]
with \(a>0\).
\end{example}
\begin{example}
Take \(M=\GL{n}\subset\R^{n\times n}\), with the standard flat affine geometry pulled back from the affine space \(\R^{n\times n}\).
The automorphism group of \(M\) is the group of affine transformations of \(\R^{n\times n}\) preserving \(M\).
Clearly the automorphism group includes \(\GL{n}\) acting by left translation, and also by right translation, and also the map taking any matrix in \(M\) to its transpose.
So \(M\) is a homogeneous space under its automorphism group.
\end{example}
\begin{problem}{ca}
What is the affine automorphism group of this example? (For an answer, see \cite{Dieudonne:1949}).
Prove that, for any discrete subgroup \(\Gamma\subset\GL{n}\), the affine structure descends to \(\Gamma\backslash M\), and find the affine automorphism group of this quotient in terms of the normalizer of \(\Gamma\) in \(\GL{n}\).
\end{problem}
\begin{example}
Affine transformations preserve centers of mass.
Take \((X,G)=(\R^n,\Aff{\R^n})\).
Suppose that some affine structure on a manifold \(M\) has developing map with bounded image \(M'\subset X\).
The automorphism group \(\Aut[M]\) preserves the center of mass of \(M'\) inside \(X\), which we can arrange by isomorphism to be the origin.
Pick a translation invariant volume form \(\Omega\) on \(X\).
For any covectors \(\alpha,\beta\in X^*\), let 
\[
\left<\alpha,\beta\right>:=
(n+2)\frac{\int_{x\in M'} \alpha(x)\beta(x)\Omega}{\int_{M'} \Omega},
\]
a positive definite inner product, so with a dual positive definite inner product on \(X\).
The volume form \(\Omega\) is arbitrary, and unique up to rescaling, which does not change the inner product, the \emph{Binet--Legendre inner product}\define{Binet--Legendre inner product}\define{inner product!Binet--Legendre} \cite{Matveev.Troyanov:2015}.
(The constants are chosen so that if \(M'\) is the unit ball for an inner product, we recover that inner product.)
Hence the automorphism group of \(M\) maps to the orthogonal group, preserving the intersections of \(M'\) with spheres around the origin.
In particular, any homogeneous flat affine geometry has unbounded developing map.
\end{example}
The automorphisms of any Cartan geometry preserve the constant vector fields, so commute with their flows, and commute with the action of the structure group \(H\).
The flows of the constant vector fields need not be defined for all time, but if defined for some time at some point, are defined for the same time throughout the automorphism group orbit through that point.
Hence the orbits of the automorphism group are permuted by those flows, and by the structure group.
Over a connected manifold, every automorphism group orbit is thereby permuted with any other.
\begin{theorem}[Kobayashi \cite{Kobayashi:1957} p. 185 lemma 6.5, \cite{Kobayashi:1995} Theorem 3.2]
The automorphism group orbits of any Cartan geometry \(\G\to M\) are closed.
\end{theorem}
\begin{proof}
Take a point \(p\in\G\) in the closure of an orbit of the automorphism group, say \(g_i p_0\to p\).
Define a map \(\varphi\) on some neighborhood of \(p_0\) by demanding that \(\varphi(p_0)=p\) and that \(\varphi\) commute with the constant vector fields, i.e.
\[
\varphi(\fl{A}p_0)=\fl{A}p;
\]
this uniquely determines \(\varphi\) near \(p_0\), but we still have to see why it is a local automorphism.
Commuting with constant vector fields, convergence of \(g_1p_0,g_2p_0,\dots\to p\) in \(\G\) implies convergence uniformly with all derivatives as maps of \(\G\).
So \(\varphi=\lim g_i\) near \(p_0\). 
Repeating the construction, \(\varphi\) extends smoothly to \(\G\).
By the same construction, define a limit for the \(g_i^{-1}\), so \(\varphi\) is a diffeomorphism.
\end{proof}
\subsection{Orbit tangent spaces}
A \emph{tangent vector}\define{tangent vector!to a subset of a manifold} to a subset \(S\subseteq\G\) at a point \(p_0\) is a vector \(v\in T_{p_0} \G\) so that, for some sequences \(\lambda_i\to\infty\) and \(A_i\to 0\), \(\fl{A_i}p_0\in S\) and \(\lambda_iA_i\to v\hook\omega\); the \emph{tangent space}\define{tangent space!to a subset of a manifold} \(T_{p_0} S\) is the set of tangent vectors to \(S\) at \(p_0\). 
The tangent spaces of each automorphism group orbit are permuted by the automorphism group.

Automorphisms of a Cartan geometry preserve every constant vector field.
So the tangent spaces to each automorphism group orbit \(\mathscr{O}\) consist of the same constant vector fields at all points of \(\mathscr{O}\).
\begin{lemma}\label{lemma:lines}
Take a Cartan geometry \(\G\to M\).
Take a sequence of nonzero tangent vectors \(v_i\in T_{p_0}\G\) with \(v_i\to 0\).
Let \(A_i:=v_i\hook\omega\in\LieG\).
Suppose that infinitely many \(\fl{A_i}p_0\) are in the automorphism orbit \(\mathscr{O}\) of a point \(p_0\in\G\).
After perhaps replacing the \(v_i\) by an infinite subsequence, the lines spanned by the vectors \(v_i\) converge to a line tangent to \(\mathscr{O}\).
A constant vector field is somewhere tangent to an orbit just when it is everywhere tangent to that orbit, which occurs just when its flow preserves that orbit.
\end{lemma}
\begin{proof}
Suppose that \(\fl{A_i}p_0=g_ip_0\), some \(A_i\to 0\), \(g_i\in\Aut\).
Pick \(\lambda_i>0\) so that \(\lambda_iA_i\) stays bounded and stays outside some neighborhood of the origin in \(\LieG\).
In particular, \(\lambda_i\to\infty\).
Replacing with an infinite subsequence, \(\lambda_iA_i\) converges, say \(\lambda_i A_i\to A\in\LieG\).
Pick some real number \(t\).
Round off \(t\lambda_i\): take integers \(n_i\) within a bounded distance of \(t\lambda_i\). 
But \(A_i\to 0\) so \(t\lambda_iA_i-n_iA_i\to 0\).
So \(n_iA_i\to tA\).
But \(\fl{n_iA_i} p_0=g_i^{n_i}p_0\in\mathscr{O}\) while \(\fl{n_iA_i}p_0\to\fl{tA}p_0\).
Since the orbit \(\mathscr{O}\) is closed as a subset of \(\G\), \(\fl{tA}p_0\in\mathscr{O}\) for all \(t\).
In particular, if \(A=v\hook\omega\) for a vector \(v\) tangent to \(\mathscr{O}\) then the flow of \(A\) preserves \(\mathscr{O}\).
\end{proof}
\begin{lemma}\label{lemma:closed.tgts}
Take a Cartan geometry \(\G\to M\).
The tangent spaces of each automorphism orbit \(\mathscr{O}\subseteq\G\) are closed cones, i.e. closed subsets of the tangent spaces of \(\G\), invariant under rescaling by real numbers, and containing the origin.
\end{lemma}
\begin{proof}
Pick a convergent sequence of tangent vectors \(v_i\in T_{p_0}\mathscr{O}\).
Take sequences \(A_{ij}\in\LieG\) with \(A_{ij}\to 0\) and sequences \(\lambda_{ij}\to \infty\) as \(j\to \infty\), so that \(\fl{A_{ij}}p_0\in\mathscr{O}\) and \(\lambda_{ij}A_{ij}\to v_i\hook\omega\) as \(j\to\infty\).
Let \(A_i:=A_{ii}\), \(\lambda_i:=\lambda_{ii}\) and apply lemma~\vref{lemma:lines}.
\end{proof}
Differentiating flows of constant vector fields,
\[
\fl{tB}\fl{tA}
=
\fl{t(A+B)}
+O(t)^2
\]
for any \(A,B\in\LieG\) close enough to zero. (How close to zero we need might vary from point to point around \(\G\)).
Recall the \emph{curvature deformed bracket}\SubIndex{curvature deformed bracket}
\[
A,B\mapsto\lb{A}{B}-k(A,B)
\]
defined by the value of the curvature at each point of \(\G\).
Taking brackets by commutators of flows,
\begin{align*}
\fl{-tB}\fl{-sA}\fl{tB}\fl{sA}
&=
e^{-tB_{\G}}e^{-sA_{\G}}e^{tB_{\G}}e^{sA_{\G}}\\
&=
e^{st\lb{A_{\G}}{B_{\G}}}+O(s,t)^3,
\\
&=
e^{st\lb{A}{B}_{\G}-k(A,B)_{\G}}+O(s,t)^3,
\\
&=
\fl{st(\lb{A}{B}-k(A,B))}
+O(s,t)^3
\end{align*}
for any \(A,B\in\LieG\) close enough to zero, with \(k\) the curvature.
(In the same vein, Melnick \cite{Melnick:2011} discovered a Baker--Campbell--Hausdorff formula.)
\begin{lemma}
Take a Cartan geometry \(\G\to M\).
The tangent spaces of each automorphism group orbit in \(\G\) are linear subspaces of the tangent spaces of \(\G\), and Lie algebras under the curvature deformed bracket.
\end{lemma}
\begin{proof}
Take tangent vectors \(v,v'\in T_{p_0}\mathscr{O}\).
Let \(A:=v\hook\omega\), \(A':=v'\hook\omega\), 
So there are sequences \(\lambda_i,\lambda'_i\to\infty\) and \(A_i,A'_i\to 0\), for which
\[
\fl{A_i}p_0,\fl{A'_i}p_0\in\mathscr{O},
\]
say 
\[
g_ip_0=\fl{A_i}p_0,
g'_ip_0=\fl{A'_i}p_0,
\]
for some \(g_i,g'_i\in\Aut\), and \(\lambda_iA_i\to A\) and \(\lambda'_iA'_i\to A'\).

Pick any sequence \(\sigma_i\to\infty\) of positive numbers.
Replacing \(\lambda_i,\lambda'_i\) by subsequences, we can assume that both sequences grow faster than \(\sigma_i\):
\[
\frac{\lambda_i}{\sigma_i},\frac{\lambda'_i}{\sigma_i}\to\infty.
\]
Let \(n_i,n'_i\) be integers within a bounded distance of
\[
\frac{\lambda_i}{\sigma_i},\frac{\lambda'_i}{\sigma_i}.
\]
So \(\lambda_i/n_i,\lambda'_i/n'_i\sim\sigma_i\to\infty\).
Take new \(A_i, A'_i, \lambda_i, \lambda'_i\) equal to old
\[
n_iA_i, n'_iA'_i, \frac{\lambda_i}{n_i},\frac{\lambda'_i}{n'_i},
\]
so we can assume that \(\lambda_i/\lambda'_i\to 1\).
\begin{align*}
g'_ig_ip_0
&=
g'_i\fl{A_i}p_0,
\\
&=
\fl{A_i}g'_ip_0,
\\
&=
\fl{A_i}\fl{A'_i}p_0,
\\
&=
\fl{A_i+A'_i+\dots}p_0,
\end{align*}
\[
\lambda_i(A_i+A'_i+\dots)=\lambda_iA_i+\frac{\lambda_i}{\lambda'_i}\lambda'_iA'_i+\dots\to A+A'.
\]
By lemma~\vref{lemma:lines}, the lines spanned by \(A_i+A'_i\) converge to a tangent line to the orbit, the span of \(A+A'\).
By the same argument, using the brackets by commutators of flows, \(\lb{A}{A'}-k(A,A')\) is also in the tangent space.
\end{proof}
\subsection{Slicing}
Take a point \(p_0\in\G\) and let \(\mathscr{O}\) denote its automorphism group orbit.
A \emph{slice}\define{slice!of group action} at \(p_0\) is a smooth embedding
\[
U\xrightarrow{\varphi}\G
\]
of an open set \(U\subset V\) of a finite dimensional vector space \(V\), so that \(0\in U\), \(\varphi(0)=p_0\), \(\varphi^{-1}\mathscr{O}=\set{0}\), and
\[
\varphi'(0)V \oplus T_{p_0} \mathscr{O}=T_p \G.
\]
\[
\includegraphics[width=4cm]{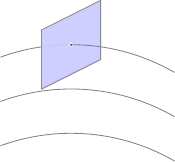}
\]
\begin{lemma}\label{lemma:outside.O}
Take a Cartan geometry \(\G\to M\).
Let \(\mathscr{O}\) be the automorphism orbit of a point \(p_0\in\G\).
Take any linear subspace \(V\subseteq T_{p_0}\G\) complementary to \(T_{p_0}\mathscr{O}\).
There is an open set \(U\subseteq V\) containing the origin so \(v\in U\mapsto\fl{A}p_0\in \G\), with \(A:=v\hook\omega\), is a slice.
\end{lemma}
\begin{proof}
Let \(V':=V\hook\omega\subseteq\LieG\), a linear subspace complementary to \(T':=T_{p_0}\mathfrak{O}\hook\omega\).
If there is a sequence of elements \(A_i\to 0\) with \(A_i\in V'\) and with \(\fl{A_i}p_0\in\mathscr{O}\), then as above we can find a convergent subsequence so that the line spanned by \(A_i\) converges in \(T_{p_0} \G\) to the line spanned by some \(A\ne 0\).
Since \(V'\subseteq\LieG\) is a linear subspace, it is a closed subset, so \(A_i\in V'\) implies \(A\in V'\).
By definition of tangent spaces, \(A\in T'\), so \(A\) is in \(T'\) and in \(V'\), a contradiction since \(V'\) is a linear subspace complementary to \(T'\).

So there is no such sequence, i.e. there is an open set \(U'\) around the origin in \(V'\) in which no point \(A\) has \(\fl{A}p_0\in\mathscr{O}\).
Take the associated set \(U\subseteq V\), i.e. with \(U\hook\omega=U'\).
If we make \(U\) smaller, we can arrange that \(A\in U\mapsto\fl{A}p_0\in\G\) is defined and an embedding.
\end{proof}
\subsection{Splitting}
\[
\includegraphics[width=4cm]{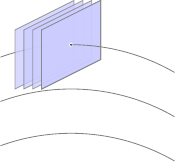}
\]
\begin{lemma}\label{lemma:splitG}
Take a Cartan geometry \(\G\to M\).
Let \(\mathscr{O}\) be the automorphism orbit of a point \(p_0\in\G\).
Take any linear subspace \(V\subseteq T_{p_0}\G\) complementary to \(T_{p_0}\mathscr{O}\).
Pick open sets \(U\subseteq V\) and \(W\subseteq T_{p_0}\mathscr{O}\) containing the origin.
If the open sets are small enough then the map
\[
(v,w)\in V\times T_{p_0}\mathscr{O}\xmapsto{\varphi}\fl{A}\fl{B}p_0\in \G,
\]
(where \(A:=v\hook\omega\) and \(B:=w\hook\omega\)) is a diffeomorphism to an open subset of \(\G\) and a slice for any fixed \(w\).
\end{lemma}
\begin{proof}
By the inverse function theorem, we can pick small enough open sets \(U,W\) to ensure that \(\varphi\) is a diffeomorphism to its image, an open set \(U_{\G}\subseteq\G\).
For any \(w\), \(\varphi(0,w)\in\mathscr{O}\), by lemma~\vref{lemma:lines}.
We can pick \(U\) small enough to ensure that \(v\mapsto\varphi(v,0)\) is a slice, by lemma~\vref{lemma:outside.O}, so stays outside \(\mathscr{O}\) except at \(v=0\).
Let \(U':=U\hook\omega\subseteq\LieG\) and \(W':=W\hook\omega\subseteq\LieG\).
If \(\varphi(v,w)\in\mathscr{O}\) then, setting \(A:=v\hook\omega\), \(B:=w\hook\omega\), we have
\[
\fl{A}\fl{B}p_0=gp_0,
\]
for some automorphism \(g\), so
\[
\fl{-A}p_0=g^{-1}\fl{B}p_0,
\]
lies in \(\mathscr{O}\), so \(A=0\), so \(v=0\).
\end{proof}
\begin{corollary}\label{corollary:aut.orbit.closed.embedded}
Every automorphism group orbit in \(\G\) of any Cartan geometry \(\G\to M\) is a closed embedded submanifold.
\end{corollary}
\begin{corollary}
Take a connected manifold \(M\) and a Cartan geometry \(\G\to M\).
Every automorphism group orbit in \(\G\) is a closed embedded submanifold.
Endow the automorphism group with the smooth structure of any of its orbits, by the map \(g\in\Aut{}\mapsto gp_0\in\G\).
This smooth structure makes the automorphism group a Lie group.
All of these Lie group structures, given by different choice of point \(p_0\), are identified by constant vector field flows and action of the structure group.
The automorphism group acts smoothly on \(\G\) and \(M\).
\end{corollary}
\begin{proof}
We take the smooth structure as an orbit in \(\G\).
Take \(g,h\in\Aut\) and some point \(p_0\in\G\).
We want to prove that \(gh\) is a smooth function of \(g,h\), i.e. that \(ghp_0\) is a smooth function of \(gp_0\), \(hp_0\).
We need only vary \(g,h\) by flows \(\fl{A},\fl{B}\), since these flows are local coordinates \(A,B\) on the orbit near each point i.e. prove that 
\[
(\fl{A}g)(\fl{B}h)p_0=\fl{A}g\fl{B}hp_0
\]
depends smoothly on \(A,B\), which is obvious.
The same computation demonstrates the smoothness of the action.
\end{proof}
\begin{example}The biholomorphism group of the unit ball \(B\subseteq\C^n\) is well known \cite{Goldman:1999} p. 68, \cite{Rudin:2008} p. 25--28 to be the subgroup \(\mathbb{PU}_{n,1}\subset\mathbb{PGL}_{n+1}\).
(The simplest example: the unit disk in \(\C\) has automorphism group consisting of M\"obius transformations.)
This biholomorphism group is also precisely the group of oriented preserving isometries of the standard metric on the ball, the \emph{hyperbolic metric}.\define{hyperbolic metric!Kobayashi}\define{Kobayashi hyperbolic metric}
It is also the group of automorphisms of the pullback to the ball of the standard flat holomorphic projective connection on \(\C^n\subset\Proj{n}\).
But \(\mathbb{PU}_{n,1}\) is \emph{not} a complex Lie group.
It preserves the complex structure on the ball, i.e. acts by biholomorphisms, but is not itself a complex manifold.
So the automorphisms of a holomorphic Cartan geometry form a real Lie group acting smoothly, but perhaps not a complex Lie group, or perhaps not acting holomorphically.
\end{example}
\subsection{Properness}
Recall that a Lie group \(G\) acting smoothly on a manifold \(M\) acts \emph{properly}\define{proper!group action}\define{group action!proper} just when 
\[
(g,m)\in G\times M\mapsto (gm,m)\in M\times M
\]
is a proper map.
Equivalently, for any sequence of \(g_i\in G\) and \(m_i\in M\), if \(m_i\) converges and \(g_im_i\) converges, some infinite subsequence of the \(g_i\) converges \cite{Abraham/Marsden:1978} p. 264.
\begin{example}
A densely winding subgroup of the torus does not act properly on the torus.
\end{example}
\begin{example}
Every action of any compact group is proper.
\end{example}
\begin{example}
If \(k=\R, \C\) or \(\Quat\), the action of \(G=k^{\times}\) on \(k^{n+1}\) is not proper: take \(g_i=2^i\), \(m_i=2^{-i}v\) for some \(v\ne 0\).
But the action on \(M:=k^{n+1}-0\) is proper, with quotient \(G\backslash M=\Proj{n}\). 
\end{example}
\begin{problem}{cb}
If \(G\) is a Lie group and \(H\subseteq G\) is a Lie subgroup, the right action of \(H\) on \(G\) is proper just when \(H\subseteq G\) is closed, and then the quotient is the homogeneous space \(X=G/H\).
\end{problem}
\begin{corollary}
The automorphism group of a Cartan geometry \(\G\to M\) on a connected manifold \(M\) acts freely and properly on \(\G\).
\end{corollary}
\begin{proof}
The topology of the automorphism group is that of the orbit through any point \(p\), i.e. \(g_ip\to gp\) just when \(g_i\to g\).
After taking a subsequence, \(p_i=\fl{A_i}\fl{B_i}p\) for some \(A_i,B_i\to 0\) as in lemma~\vref{lemma:splitG}.
Similarly, after taking another subsequence, \(g_ip_i=\fl{C_i}\fl{D_i}q\) for some \(C_i,D_i\to 0\):
\[
g_i p=\fl{-B_i}\fl{-A_i}\fl{C_i}\fl{D_i}q
\]
approaches \(q\), by continuity of the flows of the constant vector fields.
Since orbits are closed, \(q=gp\) for some \(g\in\G\), so \(g_ip\to gp\) so \(g_i\to g\).
\end{proof}
The quotient by a free and proper action is a smooth manifold, with quotient map a principal bundle map \cite{Duistermaat.Kolk:2000} p.53 theorem 1.11.4, proving the main theorem of this section: theorem~\vref{thm:strong}. (In fact, the main part of the proof that the quotient map is a principal bundle map is constructing a local splitting as we did in lemma~\vref{lemma:splitG}.)
\begin{lemma}\label{lemma:stabilzer.closed}
Take a Cartan geometry \(\G\to M\).
The stabilizer \(\Aut^{m_0}\subseteq\Aut\) of a point \(m_0\in M\) sits inside \(H\) as a closed Lie subgroup.
\end{lemma}
\begin{proof}
We can assume that \(M\) is connected without loss of generality.
Let \(G':=\Aut\) and \(H':=\Aut^{m_0}\).
Pick a point \(p_0\in\G\) mapping to \(m_0\in M\).
Embed \(g\in G'\mapsto gp_0\in\G\).
Each \(k\in H'\) moves \(p_0\) to a point of the fiber over \(m_0\), on which \(H\) acts transitively, so
\[
kp_0=p_0\bark
\]
for some \(\bark\in H\).
But since \(H'p_0=(G'p_0)\cap(p_0H)\), the image of \(k\mapsto\bark\) is a closed embedded submanifold.
\end{proof}
\begin{problem}{cc}
Take a connected manifold with a Cartan geometry.
Prove that the topology on the automorphism group determined by identifying it with an orbit in \(\G\) is the topology of pointwise convergence, but is also the compact-open topology on it as a collection of maps of \(\G\), and is also the topology of uniform convergence on compact sets with all derivatives.
\end{problem}
Sometimes we can allow manifolds with infinitely many components.
\begin{theorem}\label{thm:strongly.effective}
Suppose that \(\Gamma\) is a group of automorphisms of a Cartan geometry \(\G\to M\), closed in the topology of uniform convergence on compact sets with all derivatives.
Suppose that, for every component \(M'\subseteq M\), the only element of \(\Gamma\) that fixes every point of \(\G\) lying above \(M'\) is the identity.
Then \(\Gamma\) is a Lie group acting smoothly on \(\G\) and \(M\) and 
\[
\Gamma\to\G\to\Gamma\backslash\G
\]
is a principal bundle.
There is a finite set of smooth functions on \(\G\), invariant under \(\Gamma\), which distinguish \(\Gamma\)-orbits.
\end{theorem}
\begin{proof}
Suppose that \(\Gamma\) fixes a point of \(\G\).
Commuting with the \(H\)-action and the constant vector field, \(\Gamma\) fixes every element of \(\G\) above some component \(M'\subseteq M\), so is the identity.
Hence \(\Gamma\) acts freely.
The rest of the proof is identical.
The topology on \(\Gamma\) is almost irrelevant by the previous exercise.
As the quotient \(\Gamma\backslash\G\) is a manifold, it admits an embedding into Euclidean space, by the Whitney embedding theorem.
Take the coordinate functions of such an embedding as our finite set of smooth functions.
\end{proof}
\begin{problem}{cd}
Prove that, on any compact manifold, any Cartan geometry with compact structure group has compact automorphism group.
\end{problem}
\begin{problem}{ce}
Suppose that \(\G\to M\) is a Cartan geometry with model \((X,G)\).
Suppose that \(\Aut=\Aut[\G]\) acts transitively on \(\G\), so is identified, by choice of any point \(p_0\in\G\) with \(\G\) itself by \(g\in\Aut\mapsto gp_0\in\G\).
Let let \(\LieA\) be the Lie algebra of \(\Aut\).
Prove that the Cartan connection pulls back to a \(\omega_{\G}=\lambda \omega_G\), for a constant linear map \(\lambda\in\LieA^*\otimes\LieG\).
\end{problem}
\subsection{Developability}
\begin{theorem}
Take a Cartan geometry \(\G\to M\) on a manifold \(M\), say with model \((X,G)\).
Take a positive definite quadratic form on the Lie algebra \(\LieG\) of \(G\).
Pull this back by the Cartan connection to generate a Riemannian metric on \(\G\).
Let \(\bar\G:=\Aut\backslash\G\).
This metric induces a Riemannian metric on \(\bar\G\) by lifting each tangent vector up to a tangent vector perpendicular to the associated \(\Aut\)-orbit.
As a Riemannian manifold, \(\bar\G\) is complete in this metric precisely when \(\G\to M\) is developable.
\end{theorem}
\begin{proof}
Suppose that \(\bar\G\) is complete as a Riemannian manifold.
Take a Cauchy sequence in a path component of \(\G\).
The map \(\G\to\bar\G\) contracts, so its image in \(\bar\G\) converges.
Trivialize the bundle \(\Aut\to\G\to\bar\G\) around the limit point, say in an open set \(U\subseteq\G\) we have \(U=\bar{U}\times\Aut\) where \(\bar{U}\subseteq\bar\G\) is open.
Then \(\omega\) is a direct sum of a parallelism on \(\bar\G\) and the left invariant Maurer--Cartan form on \(\Aut\), with a metric which is a product metric, so our Cauchy sequence converges.

Suppose that \(\G\to M\) is developable.
So the Riemannian metric on \(\G\) is complete by Clifton's theorem (theorem~\vref{theorem:Clifton}).
Take a geodesic in \(\bar\G\) defined on an interval \(a<t<b\), parameterized by arclength.
It lifts to a horizontal path in \(\G\), of the same length, so extends to \(t=a\) and \(t=b\).
\end{proof}
\subsection{Decorations}
A \emph{decoration}\define{decoration} of a Cartan geometry \(\G\to M\) is a smooth map of \(\G\) to some manifold.
A \emph{decorated Cartan geometry}\define{Cartan!geometry!decorated} is a Cartan geometry equipped with a collection, perhaps infinite, of decorations.
\begin{example}
The Cartan connection identifies any vector field \(Z\) on \(M\) with an \(H\)-equivariant map \(A:=Z\hook\omega\colon\G\to\LieG/\LieH\), a decoration.
Similarly, any tensor field is a decoration, and any finite dimensional Lie algebra action.
\end{example}
Clearly a decorated Cartan geometry has automorphism group a closed subgroup of the automorphism group when the decoration is removed.
Hence \(\Aut[M]\to\G\to\Aut[M]\backslash\G\) is still a principal bundle.

\subsection{The Lie algebra}
The complete infinitesimal symmetries of the Cartan geometry flow to produce automorphisms, in the identity component of the automorphism group in this Lie group structure, as they preserve the automorphism group orbits, so flow smoothly along them.
On the other hand, take any right invariant vector field for that Lie group structure, on any one orbit.
The flows of constant vector fields and the action of the structure group move it around \(\G\) to produce a global infinitesimal symmetry, complete because it is a right invariant vector field for the Lie group structure on each orbit.
Hence the Lie algebra of this Lie group structure is the collection of complete infinitesimal symmetries, and the identity component of the automorphism group is generated by the complete infinitesimal symmetries.
\begin{corollary}
Take a Cartan geometry on a connected manifold.
The constant vector fields which are tangent to an automorphism group orbit are the left invariant vector fields on the automorphism group for the Lie group structure on each automorphism orbit.
In particular, the Lie bracket of the automorphism group is isomorphic to the curvature deformed bracket.
\end{corollary}
\begin{proof}
Automorphisms take tangent vectors to tangent vectors, and and constant vector fields to constant vector fields.
Hence the constant vector fields tangent to an automorphism orbit at a point are everywhere tangent to that automorphism orbit, and hence are left invariant vector fields on the automorphism group.
By dimension count, they are all such.
Their bracket is the curvature deformed bracket.
\end{proof}
\Danger{} Kobayashi \cite{Kobayashi:1995} Theorem 3.2 states that every automorphism group orbit is a closed submanifold, from which he derives its smooth structure, but he only proves it is a closed subset invariant under flows of the complete infinitesimal automorphisms, hence a closed union of their orbits.
He also has a smooth structure on the automorphism group (\cite{Kobayashi:1995} Theorem 3.1), so that the automorphism group is of dimension equal to the dimension of the space of complete infinitesimal automorphisms.
But it is not clear that this smooth structure gives a topology for which the action is proper, or for which the automorphism group embeds to its orbits.
Sternberg \cite{Sternberg:1983} p. 347 theorem 4.2 quotes Kobayashi to establish the smoothness and properness of the action, as do many authors, even recently.
Gromov \cite{Gromov1988} p. 84 1.5.B proves a result which implies Kobayashi's, but the proof is not clear to me; even Gromov's definition of the topology of the automorphism group is not clear to me.
d'Ambra and Gromov \cite{dAmbra/Gromov:1991} p. 70 5.12 assert that the Lie group structure on the automorphism group of a rigid geometric structure is due to Lie, with no citation, and that a proof of properness is in \cite{Gromov1988}.
I think that Dmitri Zaitsev brought this problem to my attention, probably around 2006, and I have heard it mentioned by Andreas \v{C}ap, Sorin Dumitrescu, Charles Frances, Niky Kamran and Chris Wendl \cite{chriswendl}.
\subsection{Point set topology detour}
\Danger{} There are two definitions of manifold.
Everyone agrees that a topological \(n\)-manifold is a paracompact Hausdorff space, locally homeomorphic to \(\R^n\).
(Almost everyone: some eschew paracompactness \cite{Gauld:2014,Helgason:1978}.
Palais alone, in the unusual but important \cite{Palais:1957}, eschews Hausdorffness.)
Most authors require \emph{second countability},\SubIndex{second countability} i.e. a countable basis of open sets, a stronger condition, equivalent to paracompactness together with any one of:
\begin{itemize}
\item
a countable set of components, or
\item
separability, i.e. containing a countable dense subset, or
\item
existence of an embedding into Euclidean space \cite{Munkres:2000} p. 315.
\end{itemize}
Geometers prefer theorems to remain agnostic about whether they require second countability or merely paracompactness, i.e. which definition of manifold we use.
For any homogeneous space \((X,G)\) with \(X\) and \(G\) second countable, and manifold \(M\) with finitely many components, one might conjecture that the automorphism group of any \((X,G)\)-geometry on \(M\) is second countable.
\v{C}ap and Slov\`ak \cite{Cap/Slovak:2009} p. 97 note that this conjecture has not been proved, even if \(X,G\) and \(M\) are connected.
We are not specifying which definition of manifold we use, so all of our results here admit the same statement and proof for either definition of manifold, in particular theorem~\vref{thm:strong} and corollary~\vref{corollary:finite.components}, proving the conjecture.
I am not aware of a proof of this conjecture even for Riemannian metrics on connected manifolds; the usual proofs do not consider second countability \cite{Matveev.Troyanov:2017}, \cite{Montgomery.Zippin:1974} chapter 5, \cite{Myers.Steenrod:1939}.

It is not hard to prove that a faithful and transitive Lie group action on a manifold with finitely many components is proper just when it is the isometry group of an invariant Riemannian metric \cite{Duistermaat.Kolk:2000} p. 105 Proposition 2.5.2, \cite{Meinrenken:2003} p. 12 Corollary 1.27, \cite{Palais:1961} p. 316 Theorem 4.3.1.
Using the Riemannian metric as a Cartan geometry, if a manifold has finitely many components then any Lie group acting faithfully and properly on it is second countable.
Again, I don't know of a proof of this fact in the literature.

A connected effective homogeneous space \((X,G)\) is \emph{frantic}\define{frantic!homogeneous space}\define{homogeneous!space!frantic} if, for every \((X,G)\)-geometry on any connected manifold \(M\), the automorphisms of that geometry, thought of as maps of \(M\), form a closed subgroup of the homeomorphism group, in the compact-open topology on \(M\).
Frances and Melnick (a.k.a. FRANces eT melnICk) \cite{Frances.Melnick:2019} proved that various homogeneous spaces are frantic; the classification of frantic spaces is unknown.
\subsection{Automorphisms of effective geometries}
\begin{proposition}[Sharpe \cite{Sharpe:1997}]
Every isomorphism \(\G\to\G'\) of Cartan geometries with effective model is determined by its map \(M\to M'\) on the underlying base manifolds. 
\end{proposition}
\begin{proof}
Any two isomorphisms agree up to automorphism, and any two automorphisms which agree on the base agree up to bundle automorphism leaving every fiber invariant.
In other words, it suffices to determine the kernel of the obvious group morphism taking each automorphism of a Cartan geometry to the associated diffeomorphism of the base manifold.
Hence we derive our result from proposition~\vref{proposition:characterize.verticals}.
\end{proof}
Let \(K\) be the kernel of the model \((X,G)\) of Cartan geometry \(\G\to M\) as~\vpageref{subsection:effective.homog}.
Take the linear subspace \(\LieG^K\subseteq\LieG\) consisting of those \(A\in\LieG\) on which all elements \(k\in K\) satisfy \(\Ad_k A=A\).
\begin{proposition}\label{proposition:characterize.verticals}
Take a Cartan geometry \(\G\to M\).
Pick one point \(p_i\in\G\) over each component \(M_i\subseteq M\), and an element \(k_i\in K\) of the kernel of the model.
There is at most one automorphism which is trivial on \(M\) with \(p_i\mapsto p_ik_i\).
Hence every automorphism which is trivial on \(M\) is represented by an element
\[
\set{k_i}\in K^{\pi_0(M)}.
\]
Two isomorphisms 
\[
\begin{tikzcd}
\G'\arrow[d]&\G\arrow[l]\arrow[r]\arrow[d]&\G'\arrow[d]\\
M'&M\arrow[l]\arrow[r]&M'
\end{tikzcd}
\]
of Cartan geometries which agree on \(M\) agree up to an automorphism which is trivial on \(M\).

If every element of
\[
\set{k_i}\in K^{\pi_0(M)}
\]
arises from an automorphism trivial on \(M\) then the curvature of the Cartan geometry is valued in \(\Lm*{2}{\LieG/\LieH}\otimes(\LieG^K)\).

Conversely, if the curvature of the Cartan geometry is valued in \(\Lm*{2}{\LieG/\LieH}\otimes(\LieG^K)\), there is a morphism
\[
\prod_i\fundamentalGroup{M_i}\to K
\]
so that, on the covering space \(\hat{M}\to M\) given by the kernel of this morphism, \(K^{\pi_0(M)}\) is the group of automorphisms acting trivially on \(\hat{M}\).
\end{proposition}
\begin{proof}
We can assume the model connected.
Take an automorphism \(\G\xrightarrow{\Phi}\G\) leaving every fiber invariant.
Our automorphism thus agrees, on each fiber, with right action of some element of \(H\), maybe a different element at each point: \(\Phi(p)=ph(p)\) for a unique map \(\G\xrightarrow{h}H\).
Check that \(H\)-invariance is \(h(pg)=g^{-1}h(p)g\) for every \(g\in H\).
Denote right translation by \(h\in H\) as \(\RT{h}\).
The Cartan connection pulls back, computed on a constant vector field, to
\[
(\Phi^*\omega)(A(p))
=
\left(\Phi'(p)A(p)\right)\hook\omega.
\]
\begin{problem}{kernel.automorphisms}
Prove that
\[
\Phi^*\omega=\Ad_h^{-1}\omega+h^*\omega_H.
\]
\end{problem}
\begin{answer}{kernel.automorphisms}
Compute
\begin{align*}
\Phi'(p)A(p)
&=
\left.\frac{d}{dt}\right|_{t=0}
(\fl{tA}p)h(\fl{tA}p),
\\
&=
\left.\frac{d}{dt}\right|_{t=0}
(\fl{tA}p)h(p)
+
\left.\frac{d}{dt}\right|_{t=0}
ph(\fl{tA}p),
\\
&=
\left.\frac{d}{dt}\right|_{t=0}
\RT{h(p)}\fl{tA}p
+
\left.\frac{d}{dt}\right|_{t=0}
ph(\fl{tA}p),
\\
&=
\RT{h(p)*}
A
+
\left.\frac{d}{dt}\right|_{t=0}
ph(\fl{tA}p).
\end{align*}
On the first term,
\begin{align*}
\left(\RT{h(p)*}A\right)\hook\omega
&=
A\hook \RT{h}^*\omega,
\\
&=A\hook\Ad_h^{-1}\omega,
\\
&=
\Ad_h^{-1}A.
\end{align*}
On the second, we move by \(h(\fl{tA}p)\), inside the fiber, so \(\omega=\omega_H\):
\begin{align*}
\left(\left.\frac{d}{dt}\right|_{t=0}
ph(\fl{tA}p)\right)
\hook\omega
&=
\left(\left.\frac{d}{dt}\right|_{t=0}
h(\fl{tA}p)\right)
\hook\omega_H,
\\
&=
A\hook h^*\omega_H.
\end{align*}
So finally,
\[
\Phi^*\omega=\Ad_h^{-1}\omega+h^*\omega_H.
\]
\end{answer}
So a smooth map \(\G\xrightarrow{\Phi}\G\) preserving every fiber is an automorphism just when 
\[
\omega=\Ad_h^{-1}\omega+h^*\omega_H,
\]
where \(\Phi(p)=ph(p)\).
This is equivalent to asking that we have a map \(\G\xrightarrow{h}H\) so that
\[
(\Ad_h^{-1}-I)\omega\in\LieH.
\]
Since \(\omega\) is onto \(\LieG\), the map \(\G\xrightarrow{h} H\) lies in the subgroup \(H_1\subseteq H\) of \(h\in H\) for which
\[
(\Ad_h-I)\LieG\subseteq\LieH,
\]
as~\vpageref{subsection:effective.homog}.
Inductively, \(h\) lies in the kernel, \(\G\xrightarrow{h} K\), so we write \(h\) as \(k\).

Pick any \(k_0\in K\), and define a map \(\G\xrightarrow{k}K\) by \(k(p_0)=k_0\) at some point \(p_0\), \(k(p_0g)=g^{-1}k(p_0)g\) for \(g\in H\), and the equation
\[
k^*\omega_K=\omega-\Ad_k^{-1}\omega.
\]
On the manifold \(\G\times K\), with points written as \((p,k)\), this equation is satisfied by integral manifolds of \(\vartheta:=\omega_K-\omega+\Ad_k^{-1}\omega\).
Note that \(\vartheta\) is a linear submersion from each tangent space of \(\G\times K\) to \(\LieK\) (just looking at the \(\omega_K\) term) so its kernel consists of linear subspaces of the tangent spaces of \(\G\times K\) of dimension equal to the dimension of \(\G\).
In particular, by the Frobenius theorem, if \(d\vartheta\ne 0\) on such a subspace, there is no local integral manifold, i.e. this value of \(k\in K\) does not arise from an automorphism trivial on \(M\).
Even if \(d\vartheta=0\) on some such subspace, there is at most one integral manifold through each point.

Suppose we find a linear subspace of a tangent space of \(\G\times K\) on which \(\vartheta=0\).
Compute that on that subspace,
\[
d\vartheta=\frac{1}{2}(\Ad_k^{-1}-I)K\sigma\wedge\sigma.
\]
Since \(k\) lies in the kernel \(K\) of \((X,G)\), \(\Ad_k-I\) is valued in \(\LieK\).
So \(d\vartheta\) is valued in \(\Lm*{2}{\LieG/\LieH}^*\otimes \LieK\).
Clearly \(d\vartheta=0\) on every subspace where \(\vartheta=0\) everywhere just exactly when the curvature is a \(2\)-form valued in \(\LieG^K\).
Apply the Frobenius theorem: the curvature is a \(2\)-form valued in \(\LieG^K\) just when there is a unique foliation of \(\G\times K\) by graphs of local solutions \(\text{open }\subseteq\G\xrightarrow{k}K\).

The solutions extend globally to graphs of components of \(\G\) because, just as we have global solvability of Lie equations, the equation \(0=\vartheta\) is linear in any locally faithful representation of \(K\), so is globally solvable along curves in \(\G\) with any initial condition.
The solutions are defined on a covering space of \(\G\).
Each local solution is a local section of the bundle \(\amal{\G}{H}{K}\to M\), under the \(H\)-conjugation on \(K\).
So monodromy occurs on \(M\), not on \(\G\): a morphism \(\fundamentalGroup{M}\to K\).
\end{proof}
How many derivatives of an automorphism do we need to compute before the automorphism is determined?
In a geometry with effective model on a connected manifold, if the structure group \(H\) has a maximal compact subgroup, we only need to compute derivatives at enough points of \(M\) to be able to rigidly move a generic simplex in \(M_c\), fewer derivatives than one might guess.
\subsection{Automorphisms and covering spaces}
Suppose that \(\G\to M\) is a Cartan geometry with model \((X,G)\) on a connected manifold \(M\).
Take the universal covering space \(\tilde{M}\to M\) and let \(\tilde\G\to\G\) be the pullback Cartan geometry.
Let \(\pi:=\pi_1(M)\).
A \emph{extended automorphism}\define{extended automorphism}\define{automorphism!extended} of the geometry on \(M\) is a \(\pi\)-equivariant automorphism of the geometry on \(\tilde{M}\).
If $\tilde\Phi$ is an extended automorphism, define maps by
\[
\begin{tikzpicture}[
back line/.style={},
cross line/.style={preaction={draw=white, -,
line width=6pt}}]
\matrix (m) [matrix of math nodes,
row sep=3em, column sep=3em,
text height=1.5ex,
text depth=0.25ex]{
& \tilde\G & & \G \\
\tilde\G  & & \G \\
& \tilde{M} & & M \\
\tilde{M} & & M \\
};
\path[->]
(m-1-2) edge (m-1-4)
edge  node[above] {$\tilde\Phi$} (m-2-1)
edge [back line] (m-3-2)
(m-1-4) edge (m-3-4)
edge node[above] {$\Phi$} (m-2-3)
(m-2-1) edge [cross line] (m-2-3)
edge (m-4-1)
(m-3-2) edge [back line] (m-3-4)
edge [back line] node[above] {$\tilde\varphi$} (m-4-1)
(m-4-1) edge (m-4-3)
(m-3-4) edge node[above] {$\varphi$} (m-4-3)
(m-2-3) edge [cross line] (m-4-3);
\end{tikzpicture}
\]
\begin{example}
Clearly if $\varphi$ and $\Phi$ are the identity, then $\tilde\varphi$ can be any element of the fundamental group of $M$, and  $\tilde\Phi$ is determined.
\end{example}
A diffeomorphism 
$\map[\psi]{\tilde{M}}{\tilde{M}}$
is the lift $\psi=\tilde\varphi$
of a diffeomorphism 
$\map[\varphi]{M}{M}$
just when there is a group automorphism 
\(
\mapto{\gamma \in \fundamentalGroup{M}}{\gamma' \in \fundamentalGroup{M}}
\)
so that $\psi \circ \gamma = \gamma' \circ \psi$.
\begin{example}
Every automorphism $\left(\varphi,\Phi\right)$
of a Cartan geometry lifts to some
extended automorphism, by taking 
$\tilde\varphi$ to be any lift to a map
on the universal covering space.
\end{example}
Let $\Aut[M]$ be the automorphism group,
and
\[
1 \to \fundamentalGroup{M} \to \ExtAut{M} \to \Aut[M] \to 1
\]
be the obvious exact sequence of the extended automorphism group.

Suppose that $\G \to M$ is a flat Cartan geometry with model \((X,G)\).
Pick a point $p_0 \in\G$ and let $m_0 \in M$ be its image. 
Pick $\tilde\Phi \in \ExtAut{M}$, and as above define $\varphi,\tilde\varphi,\Phi$.
Pullback $\omega$ to $\tilde\G$.
Suppose that
\(
\map[\Delta]{\tilde\G}{G}
\)
is the developing map of $\G$, so that $\Delta^*\omega_G = \omega_{\tilde\G}$.
Define
\[
h\left(\tilde\Phi\right)
=
\Delta \circ \tilde\Phi \left(\tilde{p}_0\right)
\]
and call the map  $\map[h]{\ExtAut{M}}{G}$ the \emph{extended holonomy morphism}.\define{holonomy!morphism!extended}\define{extended holonomy morphism}
\begin{lemma}
The extended holonomy morphism of any developing map on any connected manifold is a Lie group morphism, injective on Lie algebras, which agrees on $\fundamentalGroup{M}$ with the holonomy morphism.
\end{lemma}
\begin{proof}
Pick two extended automorphisms  $\Phi_1, \Phi_2$ and let $\Phi=\Phi_2 \circ \Phi_1$.
Let 
\(
g_j = \Delta\left(\Phi_j\left(p_0\right)\right),
\)
and
\(
g = 
\Delta\left(\Phi\left(p_0\right)\right).
\)
We need to prove that
\(
g=g_2 g_1.
\)
Start by assuming that
$\Phi_1\left(p_0\right)$
and
$\Phi_2\left(p_0\right)$
lie in the same path component as $p_0$.
Take smooth paths $p_1(t), p_2(t) \in \G$ so that
$p_1(0)=p_2(0)=p_0$ and 
$p_1(1)=\Phi_1\left(p_0\right)$
and
$p_2(1)=\Phi_2\left(p_0\right)$
and let $p(t)=p_2(t)$ for $0 \le t \le 1$
and then $p(t)=\Phi_2 \circ p_1(t-1)$
for $1 \le t \le 2$.
Let $g_1(t)=\Delta\left(p_1\left(t\right)\right)$
and $g_2(t)=\Delta\left(p_2\left(t\right)\right)$
and $g(t)=\Delta\left(p(t)\right)$.
Then for $t \ge 1$, 
\[
\dot{g}(t)\hook\omega_G = \dot{g}_1(t-1)\hook\omega_G 
\]
so
\(
g(t) = g_2(1) g_1(t-1).
\)
Let $t=2$.

If these $\Phi_1\left(p_0\right)$ and $\Phi_2\left(p_0\right)$ do not lie in the same path component as $p_0$, then we make use of $H$-equivariance to
get them to.

The kernel on Lie algebras consists of vector fields on \(\tilde\G\) mapping to the zero vector field on \(G\) by the developing map, a local diffeomorphism.
\end{proof}
\section{Automorphisms of submanifolds}
Take a manifold \(M\) with a Cartan geometry \(\G\to M\) with model \((X,G)\).
Take an immersed submanifold \(S\) of \(M\). 
An \emph{automorphism}\define{automorphism!of a submanifold}\define{submanifold automorphism} of this submanifold is a diffeomorphism \(S\to S\) and an automorphism of the Cartan geometry making
\[
\begin{tikzcd}
S\arrow[r]\arrow[d]&M\arrow[d]\\
S\arrow[r]&M
\end{tikzcd}
\]
commute.
\begin{example}
Consider a surface \(S\) in \(3\)-dimensional Euclidean space \(M\).
Each automorphism is a diffeomorphism of the surface and a rigid motion of Euclidean space which matches the diffeomorphism along the surface.
\end{example}
\begin{example}
If \(S\) is a closed and embedded submanifold of \(M\) the automorphisms of \(S\) in any Cartan geometry on \(M\) form the obvious closed subgroup of \(\Aut[M]\), and hence a Lie subgroup \cite{Mimura/Toda:1991} p. 44. 
\end{example}
\begin{example}
Again consider a surface \(S\) in \(3\)-dimensional Euclidean space \(M\), for example take \(S\) to be a pair of disjoint spheres, but immersed into \(M\) as the same unit sphere, so a \(2\)-to-\(1\) immersion to a single embedded sphere in \(M\).
The automorphism group is not a subgroup of the rigid motions.
Pick a dense subset of each component sphere, so that the two subsets are carried by the diffeomorphism to disjoint sets in \(M\); our immersion is injective on the union of those subsets, a dense subset of \(S\).
If we take a countable collection of spheres instead of just two, we find an automorphism group with uncountable many components.
\end{example}
\begin{example}
Again consider a surface \(S\) in \(3\)-dimensional Euclidean space \(M\): suppose that \(S\) consists of all planes parallel to and at a rational distance from a given plane. 
The automorphism group of \(S\) acts properly on \(S\) as a Lie group, but does not act properly on \(M\), even though \(S\) is injectively immersed.
\end{example}
\begin{example}
Consider a densely winding real line \(S\) in a torus \(M\).
Equip \(M\) with its usual Cartan geometry, modelled on \((X,G)=(M,M)\).
The automorphisms are the translations of \(S\), acting as a subgroup of \(M\), hence a Lie group \(S\) mapping \(S\to M\) as an immersed subgroup.
\end{example}
\begin{example}
Take a homogeneous space \(X=G/H\).
Let \(X_{\delta}\) be \(X\) with the discrete topology.
Let \(H_{\delta}:=H\) with its usual topology.
Let \(G_{\delta}\) be \(G\) as a group with the topology whose open sets are unions of translates \(gU\) of open sets \(U\subseteq H\).
The identity map \(X_{\delta}\to X\) is equivariant for the identity map \(G_{\delta}\to G\), an equivariant immersion which has automorphism group \(G_{\delta}\).
If \(X\) has positive dimension then this automorphism group is an immersed but not embedded subgroup of \(G\).
\end{example}
\begin{example}
Again let \(M\) be \(3\)-dimensional Euclidean space, and suppose that \(S\) consists of all planes parallel to and at a rational distance from a given plane in \(M\).
Equip \(M\) with its standard flat conformal geometry. 
Suppose that the origin belongs to one of our planes in \(S\).
Rescalings of \(M\) by rational numbers act on \(S\) as automorphisms.
A sequence of such automorphisms, by rational numbers converging to an irrational, acts as automorphisms on the plane through the origin, but has no limit in \(S\) on any of the other planes.
We want a topology on the automorphism group of \(S\) so that this sequence won't converge; in this example, we want the discrete topology on those rational rescalings.
\end{example}
\begin{example}
Let \(S\) be the real number line and \(M\) be the oriented Euclidean plane, and map \(\theta\in S\mapsto e^{i\theta}\in M\).
For each real constant \(c\), the map \(\theta\mapsto\theta+c\) on \(S\) and the rotation of \(M\) around the origin by angle \(c\) commute, and these are the automorphisms of this immersion.
So the automorphism group of \(S\) is \(\R\) equipped with an immersion
\[
\R\to \SO{2}\ltimes\R,
\]
which is not injective.
\end{example}
Denote the immersion as \(S\xrightarrow{\iota} M\).
Let \(\G_S:=\iota^*\G\), i.e. the set of pairs \((s,p)\) for \(s\in S\) and \(p\in\G\) both mapping to the same point of \(M\).
So \(\G_S\to\G\) is an immersion as well, with \(\dim\G_S=\dim S+\dim H\), and we have a commutative diagram
\[
\begin{tikzcd}
\G_S\arrow[d]\arrow[r]&\G\arrow[d]\\
S\arrow[r]&M.
\end{tikzcd}
\]
\begin{theorem}
Take a Cartan geometry \(\G\to M\) on a connected manifold \(M\).
Take a manifold \(S\), and an immersion \(S\to M\) which is injective on a dense open subset of \(S\).
Then the automorphisms of \(S\) form a finite dimensional Lie group for a unique Lie group structure for which
\[
\begin{tikzcd}
\Aut[S]\arrow[r]&\G_S\arrow[d]\\
&\Aut[S]\backslash\G_S
\end{tikzcd}
\]
is a principal bundle.
This Lie group acts smoothly on \(\G_S\), \(S\), \(\G\) and \(M\).
The Lie group morphism \(\Aut[S]\to\Aut[M]\) is an immersion.
\end{theorem}
\begin{proof}
Let \(\Gamma\) be the group of diffeomorphisms of \(\G_S\) preserving \(\omega\) and commuting with \(H\)-action, with the topology of uniform convergence on compact sets with all derivatives.
So \(\Aut[S]\) is a subgroup of \(\Gamma\times\Aut[M]\).
Take a linear surjection \(\LieG\xrightarrow{\pi}\LieG'\) to a vector space \(\LieG'\) of the same dimension as \(\G_S\).
Let \(\omega':=\pi\circ\omega\in\Omega^1_{\G_S}\).
Let \(\G'\subseteq\G_S\) be the set of points at which \(\omega'\) is a linear isomorphism of tangent spaces, an open set.
Taking various choices of linear surjection, we cover \(\G_S\) in these open sets \(\G'\).
On each \(\G'\), \(\omega'\) is a coframing, i.e. a Cartan connection for a Cartan geometry with infinitesimal model \((\LieG,\LieG)\), so with base space \(\G'\) (not \(S\)).
Since \(\omega'\) is a linear isomorphism on each tangent space of \(\G'\), \(\omega=a\omega'\) for a unique smooth map
\[
\G_S\xrightarrow{a}\LieG^{\prime\vee}\otimes\LieG.
\]
So \(\omega\) is determined by the decorated Cartan geometry of \(\G'\) with decoration \(a\).
Add the composition
\[
\G'\subseteq\G_S\to\G\to\G/\Aut[M]
\]
as another decoration to \(\G'\).
Inside the automorphism group of this decorated Cartan geometry, \(\Aut[S]\) is the closed subgroup commuting with the \(H\)-action.

If some automorphism in \(\Aut[S]\) fixes a point of \(\G_S\) then it fixes the image of that point in \(\G\), so acts on \(\G\) trivially, so acts on \(M\) trivially.
It then fixes every point of the dense open subset of \(S\) on which \(S\to M\) is injective.
So it fixes every point of \(S\) and of \(\G\) and hence of \(\G_S\subseteq S\times\G\).
In particular, \(\Aut[S]\) acts freely on \(\G_S\) and so on \(\G'\).
Apply theorem~\vref{thm:strongly.effective} to find that
\[
\begin{tikzcd}
\Aut[S]\arrow[r]&\G'\arrow[d]\\
&\Aut[S]\backslash\G'
\end{tikzcd}
\]
is a principal bundle, for each of the \(\Aut[S]\)-invariant open sets \(\G'\subseteq\G_S\).
If \(\Aut[S]\to\Aut[M]\) maps some vector to zero, it maps some infinitesimal automorphism vector field on \(\G_S\) to zero in \(\G\), but \(\G_S\to\G\) is an immersion.
\end{proof}
\begin{corollary}
Take a Cartan geometry \(\G\to M\) on a manifold \(M\) with finitely many components.
Take a manifold \(S\), and an immersion \(S\to M\) which is injective on a dense open subset of \(S\).
Then the automorphisms of \(S\) form a finite dimensional Lie group acting smoothly on \(\G_S\), \(S\), \(\G\) and \(M\), and mapping to \(\Aut[M]\) by a Lie group morphism which is an immersion.
\end{corollary}

\section{Example: homogeneous affine surfaces}\label{section:homog.aff.surf}
Consider how we might classify homogeneous surfaces modelled on the affine plane, i.e. with an affine connection invariant under a group of automorphisms acting transitively on the surface.
The automorphisms preserve the affine connection, so preserve its torsion tensor, and preserve the associated torsion free affine connection obtained by subtracting off the torsion tensor.
Hence we may assume that the affine connection is torsion free.

If we let \(V:=\R^2\) then the model of affine geometry is \((X,G)=(V,\GL{V}\ltimes V)\).
Write affine transformations \(x\mapsto hx+v\) as matrices
\[
g=
\begin{pmatrix}
h&v\\
0&0
\end{pmatrix}
\]
where the affine group is \(G=H\ltimes V\) with \(H=\GL{V}\).
The Cartan connection \(\omega\) then becomes
\[
\omega=
\begin{pmatrix}
\sigma&\gamma\\
0&0
\end{pmatrix}
\]
Torsion freedom is precisely
\[
d\sigma^i+\gamma^i_j\wedge\sigma^j=0.
\]
So the curvature is
\[
d\gamma^i_j+\gamma^i_k\wedge\gamma^j_k=R^i_j\sigma^1\wedge\sigma^2.
\]
\subsection{Flat surfaces and subordinacy}
In case \(R=0\), the affine connection is flat and torsion-free, so the affine Cartan geometry is flat, so locally isomorphic to the model.
It could still be homogeneous without being globally isomorphic to the model.
But then it develops to the model by a developing map which is a local isomorphism.
Assuming that \(M\) is connected, each automorphism then extends uniquely to an automorphism of the model, i.e. the affine plane.
So then \(M\) is a covering of an open set in the plane, homogeneous under a subgroup of the affine group.
A homogeneous Cartan geometry \(\G\to M\) is \emph{subordinate}\define{subordinate!homogeneous Cartan geometry}\define{Cartan!geometry!homogeneous subordinate} to another \(\G'\to M'\), with the same model, if there is a local diffeomorphism \(\G\to\G'\), equivariant for the structure group and for a morphism of the automorphism groups.
\begin{example}
The hyperbolic plane has projective connection, induced from its Riemannian geometry, equivariantly embedded into the projective plane, with the model projective connection, hence subordinate.
\end{example}
We want to find the maximal homogeneous affine surfaces, ordered by subordinacy.
We leave the reader to classify the possible homogeneous affine surfaces subordinate to these.
We can thus suppose from here on that \(R\ne 0\).
\subsection{Normalizing curvature}
Any matrix \(R=(R^i_j)\) could perhaps show up in some homogeneous affine surface.
But on any automorphism group orbit \(R\) is constant.
The automorphism group commutes with the \(H\)-action, so the \(H\)-action permutes orbits while transforming \(R\) in the obvious \(H\)-representation.
Take the structure equations:
\[
0=d\gamma+\gamma\wedge\gamma-R\sigma^1\wedge\sigma^2.
\]
Apply right action by an element \(h\in H\), using
\[
r_h^*\omega=\Ad_h^{-1}\omega
\]
which we expand out to give
\begin{align*}
r_h^*\gamma&=h^{-1}\gamma h,\\
r_h^*\sigma&=h^{-1}\sigma,
\end{align*}
to find that
\begin{align*}
0
&=
r_h^*(d\gamma+\gamma\wedge\gamma-R\sigma^1\wedge\sigma^2),
\\
&=
h^{-1}(d\gamma+\gamma\wedge\gamma)h-(r_h^*R)(1/\det h)\sigma^1\wedge\sigma^2,
\\
&=
(h^{-1}Rh-(r_h^*R)(1/\det h))\sigma^1\wedge\sigma^2,
\end{align*}
so that
\[
r_h^*R=(\det h)h^{-1}Rh.
\]
So we can reduce \(R\) to real Jordan normal form up to a scale factor, constant on the orbit.
By picking \(h=\lambda I\) for any \(\lambda\ne 0\), we can move from an orbit where \(R\) has some value to one where \(R\) has the same value scaled by \(\lambda^2\).
So we can assume that \(\det R=-1,0\) or \(1\).
If \(R\) is real diagonalizable, we can swap the two eigenvalues using
\[
h=
\begin{pmatrix}
0&1\\
-1&0
\end{pmatrix}
\]
or use
\[
h=
\begin{pmatrix}
1&0\\
0&-1
\end{pmatrix}
\]
to swap the signs of both eigenvalues, so we can assume that they are in order, say
\[
R=
\begin{pmatrix}
r&0\\
0&\pm 1/r
\end{pmatrix}
\]
with \(r\ge1\) or
\[
R=
\begin{pmatrix}
1&0\\
0&0
\end{pmatrix}
\]
or \(R=0\).
If \(R\) is complex diagonalizable but not real diagonalizable then \(\det R>0\) so we can assume that \(\det R=1\), so 
\[
R=e^{i\theta}
\]
for some \(-\pi<\theta<\pi\) with \(\theta\ne 0\).
If \(R\) is not complex diagonalizable then it has a generalized eigenspace, so a double eigenvalue, positive determinant, so \(\det R=0\) or \(1\).
For explicit calculation, if we write
\[
h=
\begin{pmatrix}
a&b\\
c&d
\end{pmatrix}
\]
then
\[
r_h^*R=
\begin{pmatrix}
d&-b\\
-c&a
\end{pmatrix}
R
\begin{pmatrix}
a&b\\
c&d
\end{pmatrix}.
\]
After a bit of linear algebra,
\[
R=
\begin{pmatrix}
1&\pm 1\\
0&1
\end{pmatrix},
\begin{pmatrix}
0&1\\
0&0
\end{pmatrix}
\text{ or }
\begin{pmatrix}
0&-1\\
0&0
\end{pmatrix}.
\]
Let \(H_0\subseteq H\) be the subgroup preserving our normal form for \(R\), i.e. the \(h\in H\) for which
\[
hR=(\det h)Rh.
\]
The set of points \(\G_0\subseteq\G\) at which \(R\) takes on this normal form is a principal right \(H_0\)-bundle.

If \(\det R\ne 0\), so \(\det R=\pm 1\), taking determinant tells us that the elements of \(H_0\) have determinant \(\pm 1\), and those with \(\det h=1\) are precisely those which commute with \(R\), while those with \(\det h=-1\) are those which anticommute with \(R\).
The possible groups \(H_0\), of dimensions \(1,2,3\) or \(4\), are:
\[
%\bgroup
%\def\arraystretch{2.5}%
\begin{array}{rl}
\toprule
R&H_0\\
\cmidrule(r){1-1}\cmidrule(l){2-2}
0&\GL{2}\\
I&\SL{2}\\
i&\C^{\times}\sqcup\bar{\mathbb{C}}^{\times}\\
-i&\C^{\times}\sqcup\bar{\mathbb{C}}^{\times}\\
e^{i\alpha}
&\C^{\times}\\
\begin{pmatrix}
1&0\\
0&-1
\end{pmatrix}
&
\Set{\begin{pmatrix}
a&0\\
0&1/a
\end{pmatrix}
}
\sqcup
\Set{\begin{pmatrix}
0&b\\
1/b&0
\end{pmatrix}
}
\\[10pt]
\begin{pmatrix}
1&0\\
0&0
\end{pmatrix}
&
\Set{\begin{pmatrix}
a&0\\
0&1/a
\end{pmatrix}
}
\\[10pt]
\begin{pmatrix}
1&\pm 1\\
0&1
\end{pmatrix}
&
\Set{\pm\begin{pmatrix}1&b\\0&1\end{pmatrix}}
\\[10pt]
\begin{pmatrix}
0&1\\
0&0
\end{pmatrix}&\Set{\begin{pmatrix}a&b\\0&\pm 1\end{pmatrix}}\\[10pt]
\begin{pmatrix}
0&-1\\
0&0
\end{pmatrix}&\Set{\begin{pmatrix}a&b\\0&\pm 1\end{pmatrix}}\\[10pt]
\begin{pmatrix}
r&0\\
0&1/r
\end{pmatrix}
&
\Set{\begin{pmatrix}
a&0\\
0&1/a
\end{pmatrix}
}
\\[10pt]
\begin{pmatrix}
r&0\\
0&-1/r
\end{pmatrix}
&
\Set{\begin{pmatrix}
a&0\\
0&1/a
\end{pmatrix}
}
\\
\bottomrule
\end{array}
%\egroup
\]
for any \(r>1\) and \(-\pi/2<\alpha<\pi/2\) and \(\alpha\ne 0\).
\subsection{Differentiate}
On the orbit, differentiate the structure equations:
\[
d\gamma+\gamma\wedge\gamma=R\sigma^1\wedge\sigma^2
\]
taking \(R\) to be constant to find
\[
0=(R\gamma-\gamma R+R(\gamma^1_1+\gamma^2_2))\wedge\sigma^1\wedge\sigma^2.
\]
\subsection{A special case}
Consider the first nonzero case: \(R=I\).
Our equation
\[
0=(R\gamma-\gamma R+R(\gamma^1_1+\gamma^2_2))\wedge\sigma^1\wedge\sigma^2
\]
simplifies to
\[
0=(\gamma^1_1+\gamma^2_2)\wedge\sigma^1\wedge\sigma^2
\]
which, by Cartan's lemma, tells us that
\[
\gamma^1_1+\gamma^2_2=a_i\sigma^i,
\]
for some \(a_1,a_2\).

By invariance of the Cartan connection under the automorphism group, these \(a_i\) are constant on every automorphism group orbit.
Differentiate the equation
\[
\gamma^1_1+\gamma^2_2=a_i\sigma^i,
\]
to find
\[
0=2\sigma^1\wedge\sigma^2+a_i\gamma^i_j\wedge\sigma^j.
\]
But \(\sigma^1,\sigma^2\) are linearly independent on the orbit, because it projects by submersion to the homogeneous affine surface \(M\).
So at most one of these \(a_i\) vanish.

Under the \(H_0\)-action, we want to see how these \(a_i\) transform.
We start with
\[
r_h^*\gamma=\Ad_h^{-1}\gamma=h^{-1}\gamma h,
\]
and
\[
r_h^*\sigma=h^{-1}\sigma.
\]
By invariance of the trace of a matrix,
\[
r_h^*\gamma^i_i=\gamma^i_i.
\]
But then 
\begin{align*}
0&=
r_h^*(\gamma^i_i-a_i\sigma^i),
\\
&=
\gamma^i_i-(r_h^*a_i)(h^{-1})^i_j\sigma^j,
\\
&=
(a_j-(r_h^*a_i)(h^{-1})^i_j)\sigma^j.
\end{align*}
By homogeneity of the automorphism group on the manifold \(M\), these \(\sigma^j\) remain linearly independent on the automorphism group orbit in \(\G\).
Hence
\[
r_h^*a_i=a_jh^j_i.
\]
Since \(H_0=\SL{2}\), we can transform these \(a_i\) by any linear transformation with unit determinant.
Since not both \(a_i\) vanish, we can arrange by \(H_0\)-action that \(1=a_1\) and \(0=a_2\) on some automorphism group orbit.
The subgroup \(H_1\subset H_0\) on which the equations \(1=a_1\), \(0=a_2\) are preserved is precisely the group of \(2\times 2\) matrices of the form
\[
\begin{pmatrix}
1&0\\
c&1
\end{pmatrix}.
\]
The set of points \(\G_1\subset\G_0\subset\G\) at which \((a_1,a_2)=(1,0)\) is therefore a principal right \(H_1\)-bundle.

Recall that the automorphism group embeds into \(\G\) as its orbit.
So at this stage, we can say that if the curvature \(R\) is a multiple of the identity, on a homogeneous affine surface which is not flat, then the automorphism group has dimension at most \(3\), since we have reduced the bundle \(\G\) to a bundle \(\G_1\) with \(1\)-dimensional structure group \(H_1\).

On the orbit, 
\begin{align*}
\gamma^1_1+\gamma^2_2=\sigma^1.
\end{align*}
Plug this in above to find
\[
0=2\sigma^1\wedge\sigma^2+\gamma^1_1\wedge\sigma^1+\gamma^1_2\wedge\sigma^2.
\]
Distributing the first term more neatly as
\[
0=(\gamma^1_1-\sigma^2)\wedge\sigma^1+(\gamma^1_2+\sigma^1)\wedge\sigma^2,
\]
we find, again by Cartan's lemma, that
\begin{align*}
\gamma^1_1&=a_{11}\sigma^1+(a_{12}+1)\sigma^2,\\
\gamma^1_2&=(a_{21}-1)\sigma^1+a_{22}\sigma^2,\\
\gamma^2_2&=(1-a_{11})\sigma^1-(a_{12}+1)\sigma^2
\end{align*}
for unique functions \(a_{ij}=a_{ji}\).
Again these are invariants, so constant on orbits.

\begin{problem}{homogeneous.affine.surface.structure.eqns}
Take exterior derivative of each of these equations in turn to see what further conditions they impose on these constants \(a_{ij}\).
\end{problem}
\begin{answer}{homogeneous.affine.surface.structure.eqns}
Already the first one yields:
\begin{align*}
0
&=
d(\gamma^1_1-a_{11}\sigma^1-(a_{12}+1)\sigma^2),
\\
&=
d\gamma^1_1-a_{11}d\sigma^1-(a_{12}+1)d\sigma^2,
\\
&=
-\gamma^1_2\wedge\gamma^2_1+\sigma^1\wedge\sigma^2\\
&\qquad-a_{11}(-\gamma^1_1\wedge\sigma^1-\gamma^1_2\wedge\sigma^2)\\
&\qquad-(a_{12}+1)(-\gamma^2_1\wedge\sigma^1-\gamma^2_2\wedge\sigma^2),
\\
&=
-((a_{21}-1)\sigma^1+a_{22}\sigma^2)\wedge\gamma^2_1+\sigma^1\wedge\sigma^2\\
&\qquad-a_{11}(-[a_{11}\sigma^1+(a_{12}+1)\sigma^2]\wedge\sigma^1-[(a_{21}-1)\sigma^1+a_{22}\sigma^2]\wedge\sigma^2)\\
&\qquad-(a_{12}+1)(-\gamma^2_1\wedge\sigma^1-[(1-a_{11})\sigma^1-(a_{12}+1)\sigma^2]\wedge\sigma^2),
\\
&=
\gamma^2_1\wedge((a_{21}-1)\sigma^1+a_{22}\sigma^2)+\sigma^1\wedge\sigma^2\\
&\qquad-a_{11}(-(a_{12}+1)\sigma^2\wedge\sigma^1-(a_{21}-1)\sigma^1\wedge\sigma^2)\\
&\qquad-(a_{12}+1)(-\gamma^2_1\wedge\sigma^1-(1-a_{11})\sigma^1\wedge\sigma^2),
\\
&=
(a_{21}-1)\gamma^2_1\wedge\sigma^1
+a_{22}\gamma^2_1\wedge\sigma^2+\sigma^1\wedge\sigma^2\\
&\qquad-a_{11}((a_{12}+1)\sigma^1\wedge\sigma^2-(a_{21}-1)\sigma^1\wedge\sigma^2)\\
&\qquad+(a_{12}+1)\gamma^2_1\wedge\sigma^1
+(a_{12}+1)(1-a_{11})\sigma^1\wedge\sigma^2,
\\
&=
(a_{21}-1)\gamma^2_1\wedge\sigma^1
+a_{22}\gamma^2_1\wedge\sigma^2+\sigma^1\wedge\sigma^2\\
&\qquad-2a_{11}\sigma^1\wedge\sigma^2\\
&\qquad+(a_{12}+1)\gamma^2_1\wedge\sigma^1
+(a_{12}+1)(1-a_{11})\sigma^1\wedge\sigma^2),
\\
&=
2a_{12}\gamma^2_1\wedge\sigma^1
+
a_{22}\gamma^2_1\wedge\sigma^2
+
(1-2a_{11}+(a_{12}+1)(1-a_{11}))\sigma^1\wedge\sigma^2.
\end{align*}
Hence \(0=a_{12}=a_{22}\) and \(a_{11}=2/3\).

Plugging these in
\begin{align*}
\gamma^1_1&=\frac{2}{3}\sigma^1+\sigma^2,\\
\gamma^1_2&=-\sigma^1,\\
\gamma^2_2&=\frac{1}{3}\sigma^1-\sigma^2.
\end{align*}
from which
\begin{align*}
d\sigma^1&=2\sigma^1\wedge\sigma^2,\\
d\sigma^2&=-\gamma^2_1\wedge\sigma^1-\frac{1}{3}\sigma^1\wedge\sigma^2,\\
\end{align*}
and then take exterior derivatives of all of our equations to find that they yield no new identities.
\end{answer}
After the reader finishes this exercise, we can write out the structure equations:
\begin{align*}
d\sigma^1&=2\sigma^1\wedge\sigma^2,\\
d\sigma^2&=-\gamma^2_1\wedge\sigma^1-\frac{1}{3}\sigma^1\wedge\sigma^2,\\
d\gamma^2_1&=-\gamma^2_1\wedge(\frac{1}{3}\sigma^1+2\sigma^2).
\end{align*}
At this point, we want to claim that we have found an example.
But how can we conclude that there is a \(3\)-dimensional Lie group with these structure equations?

There is an abstract approach, which we try first, because of its great generality.
Then we will try a more concrete approach.
\subsection{Differential closure ensures existence}\label{page:diff.closure}
Imagine that the variables \(\sigma^1,\sigma^2,\gamma^2_1\) are just abstract symbols, and consider the associative algebra they generate.
Write the multiplication of that algebra as wedge product, but it is still just a formal multiplication of abstract symbols.
We impose the condition that these symbols anticommute (i.e. to be very precise, we quotient out our algebra by the ideal generated by anticommutators \(x\wedge y+y\wedge x\) for each pair \(x,y\) of our symbols).
But we want to write down equations involving exterior derivatives of our symbols, so we need to add some operator to our algebra to represent an ``abstract exterior derivative''.
Add to our algebra some abstract symbols \(d\sigma^1,d\sigma^2,d\gamma^2_1\), and take again the abstract algebra generated by our symbols.
Impose not only that \(\sigma^1,\sigma^2,\gamma^2_1\) anticommute with one another, but also that they commute with \(d\sigma^1,d\sigma^2,d\gamma^2_1\) and that these commute with one another.
Clearly this new algebra is spanned, as a vector space, by all formal wedge products of these symbols
\[
\sigma^1,\sigma^2,\gamma^2_1,d\sigma^1,d\sigma^2,d\gamma^2_1
\]
which can only have at most one of the \(\sigma^1,\sigma^2,\gamma^2_1\) symbols without being zero, but can have as many of the others as we like.
Hence our algebra is infinite dimensional.

Construct an abstract linear map \(d\) on this algebra, by asking that \(d(\sigma^1)=d\sigma^1\), and so on, and that \(d\) satisfy the Leibnitz identity.
Looking at the basis we have written down, it is clear that the Leibnitz identity uniquely determines \(d\) on each basis element in terms of other basis elements.

When we write ``consider the equations
\begin{align*}
d\sigma^1&=2\sigma^1\wedge\sigma^2,\\
d\sigma^2&=-\gamma^2_1\wedge\sigma^1-\frac{1}{3}\sigma^1\wedge\sigma^2,\\
d\gamma^2_1&=-\gamma^2_1\wedge(\frac{1}{3}\sigma^1+2\sigma^2)\text{''}
\end{align*}
in this algebra, we mean simply that we consider the ideal generated by
\begin{align*}
d\sigma^1&-2\sigma^1\wedge\sigma^2,\\
d\sigma^2&+\gamma^2_1\wedge\sigma^1+\frac{1}{3}\sigma^1\wedge\sigma^2,\\
d\gamma^2_1&+\gamma^2_1\wedge(\frac{1}{3}\sigma^1+2\sigma^2).
\end{align*}
Taking exterior derivatives of both sides of all of these equations, and then substituting in the right hand sides for the left hand sides, is precisely working in the quotient by this ideal.
Finding that we arrive at \(0=0\) in all such equations is precisely saying that this ideal is a \(d\)-closed ideal, also known as a \emph{differential ideal}.\define{differential ideal}

Let us consider more abstractly what this means.
Take any abstract symbols \(\omega^i\) and associate to them abstract symbols \(d\omega^i\), and take the associative algebra generated by their formal wedge products, with the \(\omega^i\) anticommuting with one another, and commuting with the \(d\omega^i\), which commute with one another.
Define \(d\) by \(d(\omega^i)=d\omega^i\) and the Leibnitz identity.

Suppose we have some equations determining every \(d\omega^i\) as some constant coefficient linear combination of wedge products of the \(\omega^j\).
(A physicist might say that the \(d\) operation is \emph{closed}, i.e. we know how to compute it on generators.)
We can write those equations out as the ideal generated by some expressions
\[
d\omega^i+c^i_{jk}\omega^j\wedge\omega^k.
\]
Note that these \(c^i_{jk}\) can be assumed anticommutative, i.e.
\[
0=c^i_{jk}+c^i_{kj}.
\]
Consider also a different algebra, perhaps not associative: let \(\LieG\) be \(\R^n\) with basis \(e_1,\dots,e_n\) and with multiplication
\[
\lb{x}{y}:=c^i_{jk}x^jy^ke_i.
\]

The reader can then check that our ideal is a differential ideal precisely when the multiplication in \(\LieG\) satisfies the Jacobi identity.
But the Jacobi identity is precisely the condition that \(\LieG\) is a Lie algebra.
By Lie's Third Theorem (theorem~\vref{theorem:Lie.3.a}), this is precisely that \(\LieG\) is the Lie algebra of some Lie group \(G\) (which can be assumed to be connected and simply connected).
But then our formal equations
\[
d\omega^i=-c^i_{jk}\omega^j\wedge\omega^k
\]
are precisely the Maurer--Cartan structure equations of \(G\).

Hence we can say that, in the formal algebra defined above, a set of equations determining every \(d\omega^i\) as some combination of the \(\omega^j\) are the Maurer--Cartan structure equations of some Lie group precisely when, formally taking \(d\) on both sides of the equations, and plugging the original equations back in, yields nothing but \(0=0\).

Another way to look at this: it is purely a computational convenience.
Writing out everything in terms of abstract symbols \(\omega^i\), which we always think of as if they were \(1\)-forms on some manifold, is really just a convenient way to organize these \(c^i_{jk}\) constants.
Rather than having \(n^2(n-1)/2\) constants \(c^i_{jk}\), we have \(n\) formal symbols \(\omega^i\).

Applied to our example, there must be a \(3\)-dimensional Lie group with the structure equations
\begin{align*}
d\sigma^1&=2\sigma^1\wedge\sigma^2,\\
d\sigma^2&=-\gamma^2_1\wedge\sigma^1-\frac{1}{3}\sigma^1\wedge\sigma^2,\\
d\gamma^2_1&=-\gamma^2_1\wedge(\frac{1}{3}\sigma^1+2\sigma^2).
\end{align*}

For more information on algebras generated by formal symbols with some notion of \(d\) (also known as \emph{differential graded algebras})\SubIndex{differential graded algebra} see \cite{Yekutieli_2019} chapter 3.
But we warn the reader that the formal theory has nothing further to offer us; the route to prove existence of Cartan geometries with local conditions on their curvature is via exterior differential systems.
We will not discuss exterior differential systems here; see \cite{BCGGG:1991,Ivey/Landsberg:2003,McKay}.

\subsection{Recognizing the Lie group}
Now that we are convinced that our structure equations are those of a \(3\)-dimensional Lie group, we have to ask what Lie group it is, and whether it really acts as the automorphism group of a homogeneous affine surface.
We start to see something when we note that 
\[
d\gamma^2_1=-\gamma^2_1\wedge(\frac{1}{3}\sigma^1+2\sigma^2)
\]
makes it look like perhaps there is something important about 
\[
\frac{1}{3}\sigma^1+2\sigma^2,
\]
or some suitable multiple thereof.
Let us replace the role of \(\sigma^2\) with that of
\[
\zeta:=\sigma^2+\frac{1}{6}\sigma^1.
\]
This simplifies the structure equations and makes them not so hard to recognize.

We will cheat: start with the prior knowledge, from studying projective structures on curves, of the structure equations of \(\SL{2}\): the Lie algebra \(\LieSL{2}\) consists of the traceless \(2\times 2\) matrices, so the Maurer--Cartan \(1\)-form has the form
\[
\omega=
\begin{pmatrix}
\zeta&\xi\\
\eta&-\zeta
\end{pmatrix}
\]
so that the Maurer--Cartan equations \(0=d\omega+\lb{\omega}{\omega}\) expand to
\begin{align*}
d\xi&=-2\zeta\wedge\xi,\\
d\zeta&=-\xi\wedge\eta,\\
d\eta&=2\zeta\wedge\eta.
\end{align*}

Let
\begin{align*}
\xi:=\sigma^1, 
\eta:=\gamma^2_1, 
\zeta:=\sigma^2+\frac{1}{6}\sigma^1.
\end{align*}
Check that these satisfy the structure equations of \(\SL{2}\).
So our homogeneous affine surface has automorphism group with identity component \(\SL{2}\), at least up to covering.
Note that our subgroup \(H_2\) is exactly embedded as the closed subgroup of matrices
\[
\begin{pmatrix}
1&0\\
c&1
\end{pmatrix}.
\]
But this is precisely the stabilizer of the point 
\[
x_0:=
\begin{pmatrix}
0\\
1
\end{pmatrix}
\]
So one such homogeneous affine surface \((X,G)\) is \(X=\R^2-0\) with \(G=\SL{2}\).
Note that the affine geometry is not flat, even though \(X\subset\R^2\) has a \(G\)-invariant flat affine geometry.
Note that \(G=\SL{2}\to X\) is a principal right \(H_1\)-bundle, and that \(G\) embeds into \(\G_1\) as each of its orbits.
But \(\G_1\to X\) is also a principal right \(H_1\)-bundle, hence \(G=\G_1\): the automorphism group of the affine geometry on \(X\) is precisely \(G=\SL{2}\).

We can make the \(1\)-forms more explicit.
Write each matrix \(g\in\SL{2}\) as
\[
g=
\begin{pmatrix}
a&x\\
c&y
\end{pmatrix}
\]
and compute
\[
g^{-1}dg
=
\begin{pmatrix}
y\,da-x\,dc&y\,dx-x\,dy\\
a\,dc-c\,da&a\,dy-c\,dx
\end{pmatrix}
=
\begin{pmatrix}
\zeta&\xi\\
\eta&-\zeta
\end{pmatrix}
=
\begin{pmatrix}
\sigma^2+\frac{1}{6}\sigma^1&\sigma^1\\
\gamma^2_1&-\sigma^2-\frac{1}{6}\sigma^1
\end{pmatrix}.
\]
Hence
\begin{align*}
\sigma^1&=y\,dx-x\,dy,\\
\sigma^2&=y\,da-x\,dc-\frac{1}{6}(y\,dx-x\,dy),\\
\gamma^2_1&=a\,dc-c\,da.
\end{align*}

Suppose we find another connected homogeneous affine surface \((X',G')\) with curvature a nonzero multiple of the identity, so the same structure equations, so that in particular \(G'\) has the same Lie algebra as \(G=\SL{2}\) and some point of \(X'\) has the same stabilizer \(H_1\) as the point \((0,1)\in X\).
Up to covering, \(X'\) is connected and simply connected.
The group \(H_1\) is contractible, so \(G'\) is connected and simply connected.
So \(G'\) is the universal covering group of \(G=\SL{2}\) and \(X'\) is the universal covering surface of \(X\).
The connected homogeneous affine surface with curvature a nonzero multiple of the identity is thus unique up to covering.

\section{Computing homogeneous examples}
A Cartan geometry \(\G\to M\) is \emph{homogeneous}\define{homogeneous!Cartan geometry}\define{Cartan!geometry!homogeneous} when its automorphism group acts transitively on \(M\).
How can we find homogeneous Cartan geometries of a given model, perhaps even classify all of them?
Let us organize the lessons we learned from the example of homogeneous affine surfaces.
By theorem~\vref{thm:strong}, the automorphism group embeds into the total space.
Since the Cartan geometry is invariant under automorphisms, the curvature \(k\) is constant along each automorphism group orbit.
We can alter that constant by picking a different orbit, i.e. by action of the structure group, say \(H\); specifically, for any \(h\in H\),
\[
r_h^*k=\rho(h)k,
\]
where
\[
H\xrightarrow{\rho}\LieG\otimes\Lm{2}{(\LieG/\LieH)}^*
\]
is the obvious \(H\)-representation, i.e. the curvature module.
So we use the structure group to normalize the curvature into some normal form, some element \(k_0\) of the curvature module.
Let \(H_1\subseteq H\) be the subgroup preserving that normal form, the set \(\G_1\subseteq\G\) on which \(k=k_0\).
Then \(\G_1\to M\) is a principal right \(H_1\)-bundle.
There will be some automorphism group orbit inside \(\G_1\).
On that automorphism group, the structure equation
\[
0=d\omega+\frac{1}{2}\lb{\omega}{\omega}-\frac{1}{2}k_0\sigma\wedge\sigma
\]
differentiates to expand out into expressions involving only \(d\omega,d\sigma\).
But \(\sigma=\omega+\LieH\), so we expand out everything into \(d\omega\), which the structure equation turns back into an expression in \(\omega\) and \(k_0\).
This could give additional algebraic conditions on \(k_0\).

There is another source of equations.
On \(\G_1\), the Cartan connection \(\omega\) is perhaps no longer a trivialization of the tangent bundle.
Rather, if \(H_1\subseteq H\) has positive codimension, then \(\omega+\LieH_1\) is semibasic, say
\[
\omega+\LieH_1=a\,\sigma.
\]
This expression
\[
\G_1\xrightarrow{a}(\LieG/\LieH_1)\otimes(\LieG/\LieH)^*
\]
is \(H_1\)-equivariant.
Clearly it is also invariant under the automorphism group, hence constant on our orbit.
Repeating the process, pick some constant value \(a_0\) for \(a\) and let \(H_2\subseteq H_1\) be the subgroup fixing
\[
a_0\in(\LieG/\LieH_1)\otimes(\LieG/\LieH)^*.
\]
We then let \(\G_2\subseteq\G_1\) be the subset of points on which \(a=a_0\), a principal right \(H_2\)-subbundle.
Again, when applying \(d\) to the structure equations, they expand out to equations in the components of \(\omega,a_0,k_0\).
Again, these equations might force algebraic relations among \(a_0,k_0\).

On \(\G_2\), 
\[
\omega+\LieH_2=b\,\sigma,
\]
for a unique
\[
\G_2\xrightarrow{b}(\LieG/\LieH_2)\otimes(\LieG/\LieH)^*
\]
which projects to \(a_0\) under the obvious map
\[
\LieG/\LieH_2\to\LieG/\LieH_1
\]
and is \(H_2\)-equivariant.

Continue in this way until we find that the subgroups \(H_{k+1}\subset H_k\) are of equal dimension, so that we don't get any new invariants showing up.
(In particular, if some \(H_j\) is a linear algebraic group, then all subsequent
\[
H_{j+1},H_{j+2},\dots
\]
are linear algebraic, and each contained in the previous.
By the ascending chain condition, the sequence eventually stabilizes, say all equal to some \(H_0\subseteq H\).)
Eventually, \(\omega+\LieH_0=c_0\sigma\) is a constant, and plugging in to the structure equations
\[
0=d\omega+\frac{1}{2}\lb{\omega}{\omega}-\frac{1}{2}k_0\sigma\wedge\sigma
\]
yields algebraic equations on \(c_0,k_0\) which are satisfied by our choice of \(c_0,k_0\).
These equations then determine structure constants for the remaining linearly independent components \(\omega^i\) of \(\omega\).
By our discussion above, these are the structure equations of a Lie group.
Of course, we only see from the structure equations the Lie algebra of the Lie group and the embedding of the various subalgebras.
We have to determine the Lie group of automorphisms of this homogeneous geometry (if the geometry really exists!).
There is still the possibility that the resulting structure equations are those of a Lie group for which the Lie algebra of \(H_0\) does not exponentiate to a closed subgroup.
Moreover it may not be clear how to find a more explicit geometric description of this Lie group.

\chapter{Mutation}
\section{Changing model: charts and transition maps}
We begin by considering locally homogeneous structures related by an elementary ``change of models''.
Take a morphism \((X,G)\xrightarrow{(\varphi,\Phi)}(X',G')\) of homogeneous spaces.
By equivariance, \(X\to X'\) is a local diffeomorphism just when it is a covering map to an invariant open subset of \(X'\).
If this happens then every \((X,G)\)-structure induces an \((X',G')\)-structure by composing its charts with \(X\to X'\).
Let us make this explicit.

Cover our manifold \(M\) in open sets \(M_a\subseteq M\), on which we have our \(X\)-charts  \(M_a\xrightarrow{\varphi_a}X_a\) defined.
Take a pair of \(X\)-charts \(\varphi_a,\varphi_b\) from the \((X,G)\)-structure.
Their transition map is a constant element \(g_{ab}\in G\).
The preimage of \(M_a\) in the bundle \(\G\) is an open set \(\G_a\) identified with the open set \(G_a\subseteq G\) which is the preimage of \(X_a\).
Our transition map glues the open sets \(G_a\subseteq G\) together by
\[
g\in G_a\mapsto g_{ab}g\in G_b.
\]
invariant under right \(H\)-action.

We can replace the open sets \(M_a\) by smaller open sets so that \(X\to X'\) restricts to each \(X_a\) to a diffeomorphism \(X_a\to X_a':=\varphi(X_a)\subseteq X'\).
Applying \(G\to G'\) gives transition maps \(g'_{ab}:=\Phi g_{ab}\) which glue together the \((X',G')\)-structure with \(X'\)-charts \(\varphi'_a:=\varphi\circ\varphi_a\), giving a bundle \(\G'\to M\) out of the open sets \(G'_a:=\Phi(G_a)H'\) glued by those transition maps.

It will turn out to be important and not easy to find a useful expression for the Cartan connection on \(\G'\); we return to this problem~\vpageref{page:changing.model:form}.

Take an \((X',G')\)-structure \(\G'\to M\).
We want to ask whether this structure is induced from an \((X,G)\)-structure.
Factoring \((X,G)\to(X',G')\) into 
\[
(X,G)\to(X',h(G))\to(X',G'),
\]
our question splits into one in which \(G\xrightarrow{h}G'\) is surjective, and one in which it is injective.

If the bundle \(\G'\to M\) arises from some smooth principal \(H\)-bundle \(\G\to M\), and if \(\LieG\to\LieG'\) is an isomorphism, we just let \(\omega:=\Phi^{-1}\circ\pi^*\omega'\).
\begin{example}
If \(M\) is a circle then (by picking any connection) all bundles on \(M\) are quotients \(\G=(\R\times H)/\sim\) where \((x,h)\sim(x+1,h_0h)\), with \(h_0\) unique up to conjugacy.
So every element of \(H'\) is conjugate to an element in the image \(H\to H'\) just when all \(H'\)-bundles arise from \(H\)-bundles.
\end{example}
\begin{example}
If \(M\) is contractible then every \((X',G')\)-structure lifts to an \((X,G)\)-structure, since all bundles on \(M\) are trivial.
\end{example}
\begin{example}
If \(M=S^n\) is a sphere, it has a covering by two contractible open sets (the sphere punctured at the north or south pole).
The overlap contracts to the equator.
So if \(\pi_{n-1}(H)\to\pi_{n-1}(H')\) is onto, we can lift every \(H'\)-bundle on \(S^n\) to an \(H\)-bundle.
\end{example}
\section{Induced bundles}
To understand more clearly how bundles change when we change models of locally homogeneous structures, we take a brief detour to the theory of principal bundles.
If \(\Bun=\Bun_G\to M\) is a principal right \(G\)-bundle, and \(G\xrightarrow{\rho}G'\) is a Lie group morphism, let \(\Bun'=\Bun_{G'}=\Bun_G \times^G G'\) be the quotient of \(\Bun_G\times G'\) by the \(G\)-action
\[
(p,g')g=(pg,\rho(g)^{-1}g');
\]
this \(\Bun'\) is the \emph{associated principal right \(G'\)-bundle}\define{associated!principal bundle}\define{principal bundle!associated} of \(\Bun\), with right \(G'\)-action descending from
\[
(p,g')g''=(p,g'g'')
\]
for \(g',g''\in G'\), and with smooth \(G\)-equivariant map 
\[
\begin{tikzcd}[cramped]
\Bun\arrow[dr]\arrow[rr]&&\Bun'\arrow[dl]\\
&M&
\end{tikzcd}
\]
descending from
\[
\begin{tikzcd}[cramped,column sep=0em]
p\arrow[d,maps to]&\in&\Bun\arrow[d]\\
(p,1)\arrow[d,maps to]&\in&\Bun\times G'\arrow[d]\\
(p,1)G&\in&\Bun'.
\end{tikzcd}
\]
\begin{theorem}
Suppose that \(G\to G'\) is a morphism of Lie groups and that \(G\to\Bun\to M\) is a principal \(G\)-bundle.
Every principal \(G'\)-bundle \(G'\to\Bun'\to M\) which admits a principal bundle morphism
\[
\begin{tikzcd}
\Bun\arrow[rr]\arrow[dr]&&\Bun'\arrow[dl]\\
&M&
\end{tikzcd}
\]
is canonically isomorphic to \(\amal{\Bun}{G}{G'}\).
\end{theorem}
\begin{proof}
Denote the group morphism as \(G\xrightarrow{\rho}G'\).
The bundle map, denoted \(\Bun\xrightarrow{\Phi}\Bun'\), extends to a map
\[
(p,g')\in\Bun\times G'\xmapsto{\Psi}\Phi(p)g'\in\Bun'.
\]
It is invariant under the right \(G\)-action
\[
(p,g')g:=(pg,\rho(g)^{-1}g'),
\]
which is right \(G'\)-equivariant, hence descends to a morphism of bundles, which is an isomorphism, as it is a \(G'\)-bundle morphism which is the identity on the base.
\end{proof}
\section{Changing model: bundle and form}\label{page:changing.model:form}
We return to the problem of how the bundle and form change when we change the model of a locally homogeneous structure.
We want to describe the new bundle, and its associated form, in a language more suited to our goal of studying Cartan geometries:
\begin{lemma}\label{lemma:change.model}
Take a morphism \((X,G)\xrightarrow{(\varphi,\Phi)}(X',G')\) of homogeneous spaces.
Denote the associated Lie algebra morphism also by \(\Phi\).
Pick corresponding points \(x_0\in X\), \(x_0':=\varphi(x_0)\in X'\).
Let \(H:=G^{x_0},H':={G'}^{x_0'}\), \(H'':=\Phi^{-1}H'\), \(X'':=\varphi(X)\).

Then \(H''\to G\to X''\) is a principal bundle and induces the pullback bundle \(G_{H'}=\left.G'\right|_{X''}\to X''\).
In more detail: define
\[
(g,h')\in G\times H'\xrightarrow{\Psi}g':=gh'\in G'.
\]
Then at any point \((g,h')\in G\times H'\),
\[
\Psi^*\omega_{G'}=\omega_{H'}+\Ad_{h'}^{-1}\Phi\omega_G.
\]
The map \(\Psi\) is a surjective \(H''\)-bundle map to \(\left.G'\right|_{X''}\) under the right \(H''\)-action
\[
(g,h')h:=(gh,\Phi(h)^{-1}h'),
\]
\(H'\)-equivariant under the right \(H'\)-action
\[
(g,h')k':=(g,h'k')
\]
and \(G\)-equivariant under the left \(G\)-action
\[
a(g,h')=(ag,h').
\]
Quotienting by that right \(H''\)-action gives an \(H'\)-bundle isomorphism
\[
\begin{tikzcd}
G_{H'}=\amal{G}{H''}{H'}\arrow[rr]\arrow[dr]&&\left.G'\right|_{X''}\arrow[dl]\\
&X''.&
\end{tikzcd}
\]
\end{lemma}
\begin{proof}
This map \(\Psi\) is clearly invariant and equivariant as stated:
\begin{align*}
\LT{\Phi(g)}\Psi&=\Psi\LT{g},\\
\RT{h'}\Psi&=\Psi\RT{h'},\\
\RT{h}\Psi&=\Psi\RT{h}.
\end{align*}
Therefore \(\Psi\) descends as stated to a map \(\amal{G}{H''}{H'}\to G'\).
It image is clearly the pullback \(\left.G'\right|_{X''}\) of \(G'\to X'\) to \(X''\subset X'\).
By \(G\)-equivariance, \(X''\) is a \(G\)-orbit in \(X'\), so a smooth submanifold.
Two elements \((g,h),(\bar{g},\bar{h})\in G\times H'\) map to the same point \(g'\in G'\)  precisely when
\[
\Phi(g)h=\Phi(\bar{g})\bar{h}
\]
i.e.
\[
\Phi(g^{-1}\bar{g})=h\bar{h}^{-1},
\]
so \(\bar{g}=gk\in\Phi^{-1}H'=H''\) and \(\Phi(k)=h\bar{h}^{-1}\) so \(\bar{h}=\Phi(k)^{-1}h\).
So the map \(\Psi\) descends to an injection
\[
\amal{G}{H''}{H'}\to G'.
\]
The subgroup \(H''=\Psi^{-1}H'\subseteq G\) is a closed Lie subgroup since \(H'\subseteq G'\) is closed. 
Hence it acts on \(G\) freely and properly, hence acts on \(G\times H'\) freely and properly, so the quotient is a smooth manifold and
\[
G\times H\to\amal{G}{H''}{H'}
\]
is a smooth principal \(H''\)-bundle, and the map to \(G'\) is smooth.
The left \(G\)-action and right \(H'\)-action on \(G\times H'\) are smooth, and invariant under the right \(H''\)-action, so quotient to smooth \(\Psi\)-equivariant actions.

Thus 
\[
\amal{G}{H''}{H'}\to \left.G'\right|_{X''}
\]
is an equivariant bijection of homogeneous spaces of \(G\times H'\), hence a diffeomorphism, and a bundle isomorphism over \(X'\) by equivariance.
\begin{problem}{check.homogeneous.Maurer.Cartan}
Check the equation on Maurer--Cartan forms.
\end{problem}
\begin{answer}{check.homogeneous.Maurer.Cartan}
Take some \(u\in T_g G\) and some \(v\in T_{h'}H'\); let \(A:=u\hook\omega_G\) and \(B:=v\hook\omega_{H'}\).
Compute
\begin{align*}
\Psi^*\omega_{G'}(u,v)
&=
\omega_{G'}\Psi'(g,h')(u,v),
\\
&=
\LT{g'*}^{-1}\Psi'(g,h')(u,v),
\\
&=
\LT{\Psi(g)h'*}^{-1}\Psi'(g,h')(u,v),
\\
&=
\LT{h'*}^{-1}
\LT{\Psi(g)*}^{-1}\Psi'(g,h')(u,v),
\\
&=
\LT{h'*}^{-1}
\LT{\Psi(g^{-1})*}\Psi'(g,h')(u,v),
\\
&=
\LT{h'*}^{-1}
(\LT{\Psi(g^{-1})}\Psi)'(g,h')(u,v),
\\
&=
\LT{h'*}^{-1}
(\Psi\LT{g^{-1}})'(g,h')(u,v),
\\
&=
\LT{h'*}^{-1}
\Psi'(1,h')(\LT{g^{-1}})'(g,h')(u,v),
\\
&=
\LT{h'*}^{-1}
\Psi'(1,h')(\LT{g*}^{-1}u,v),
\\
&=
\LT{h'*}^{-1}
\Psi'(1,h')(A,v),
\\
&=
\LT{h'*}^{-1}
\Psi'(1,h')(A,0)+
\LT{h'*}^{-1}
\Psi'(1,h')(0,v),
\\
&=
\LT{h'*}^{-1}
\RT{h'*}A+
\LT{h'*}^{-1}v,
\\
&=
\Ad_{h'}^{-1}A+B.
\end{align*}
\end{answer}
\end{proof}
\begin{corollary}
Take a morphism \((X,G)\xrightarrow{(\varphi,\Phi)}(X',G')\) of homogeneous spaces with \(X\to X'\) a local diffeomorphism.
Denote the associated Lie algebra morphism also by \(\Phi\).
Pick corresponding points \(x_0\in X\), \(x_0':=\varphi(x_0)\in X'\).
Let \(H:=G^{x_0},H':={G'}^{x_0'}, H'':=\Phi^{-1}H'\).

Take an \((X,G)\)-structure with bundle \(H\to\G\to M\) and with associated form \(\omega\).
Then the associated \((X',G')\)-structure has bundle
\[
\G'=\amal{\G}{H''}{H'},
\]
and its Cartan connection pulls back to \(\G\times H'\) to become
\[
\omega_{\G'}=\omega_{H'}+\Ad_{h'}^{-1}\Phi\omega_{\G}.
\]
\end{corollary}
\begin{proof}
From lemma~\vref{lemma:change.model}, this is true on each of the \(H\)-invariant open sets in the bundle \(\G\) which are identified by \(X\)-charts with \(H\)-invariant open subsets of \(G\).
By \(G\)-equivariance, this description glues together under the transition maps.
\end{proof}

\section{Changing model}
Take a morphism \((X,G)\xrightarrow{(\varphi,\Phi)}(X',G')\) of homogeneous spaces, say with \(X=G/H\) and \(X'=G'/H'\).
By equivariance, \(X\to X'\) has differential, at some point, a linear isomorphism, just when it is a local diffeomorphism and a covering map to its image.
If \(X \to X'\) is a local diffeomorphism and a covering map to its image then every \((X,G)\)-geometry \(\G\to M\) induces an \((X',G')\)-geometry \(\G'\to M\) by taking bundle \(\G':=\amal{\G}{H}{H'}\) and form \(\omega'\) which pulls back to \(\G\times H'\) to be, at each point \((p,h')\in\G\times H'\),
\[
\omega'=\omega_{H'}+\Ad_{h'}^{-1}\Phi\omega
\]
where \(\LieG\xrightarrow{\Phi}\LieG'\) is the induced Lie algebra morphism.
Note that this is much like what we did with locally homogeneous structures, but not quite, because we quotiented by \(\Phi^{-1}H'\) there.
We can't do that here, because \(\Phi^{-1}H'\subset G\) doesn't have an action on the \(H\)-bundle \(\G\to M\).
Since \(X\to X'\) is a covering map \(\Phi^{-1}H'/H\) is the fiber, so discrete.
So quotienting by \(H\) in place of \(\Phi^{-1}H'\) changes very little in the theory of locally homogeneous structures, and is the best we can do in the theory of Cartan geometries.

Take an \((X',G')\)-geometry \(\G'\to M\).
We want to ask whether this geometry is induced from an \((X,G)\)-geometry.
Factoring \((X,G)\to(X',G')\) into 
\[
(X,G)\to(X',h(G))\to(X',G'),
\]
our question splits into one in which \(G\xrightarrow{h}G'\) is surjective, and one in which it is injective.
If the bundle \(\G'\to M\) arises from some smooth principal \(H\)-bundle \(\G\to M\), and if \(\LieG\to\LieG'\) is an isomorphism, we just let \(\omega:=\Phi^{-1}\circ\pi^*\omega'\).
We will return to this question~\vpageref{page:Cartan.geometry.holonomy}.

\section{Infinitesimal models}
It is convenient at times to drop reference to models; we really only use ``infinitesimal data'' from the homogeneous model.
An \emph{infinitesimal model}\define{infinitesimal model} is a pair \((\g,H)\) of a Lie group \(H\), with Lie algebra \(\LieH\), and a finite dimensional \(H\)-module \(\g\) containing \(\LieH\) as an \(H\)-submodule.
It is \emph{Lie}\define{infinitesimal model!Lie}\define{Lie infinitesimal model} if \(\g\) is a finite dimensional Lie algebra \(\g=\LieG\) and \(H\) acts on \(\LieG\) as Lie algebra automorphisms extending its adjoint action on \(\LieH\).
Every homogeneous space \((X,G)\) has Lie infinitesimal model \((\LieG,H)\) where \(H:=G^{x_0}\) is the stabilizer of a point \(x_0\in X\).
\begin{problem}{cf}
Find a Lie infinitesimal model \((\LieG,H)\) which is \emph{not} the infinitesimal model of \emph{any} homogeneous space \((X,G)\); geometries with such models are occasionally encountered in the literature.
\end{problem}
\begin{problem}{cg}
Define the kernel of an infinitesimal model, generalizing the definition for a homogeneous space.
\end{problem}
As far as I know, the effective infinitesimal models with \(\g/\LieH\) of low dimension, even of dimension \(1\), are not classified.
If \(\g/\LieH\) has dimension zero, \(\g=\LieH\) and \((\g,H)=(\LieH,H)\) is just a Lie group \(H\) with its Lie algebra; it is effective just when \(H=\set{1}\), almost effective just when \(H\) has dimension zero.
\begin{example}
If \((X,G)\) is homogeneous and \(X\) is connected, the homogeneous spaces \((X,G)\) and \((\tilde{X},\tilde{G})\) have the same infinitesimal model.
Every \((X,G)\)-geometry is thus a \((\tilde{X},\tilde{G})\)-geometry.
\end{example}
\section{Affine forms}
An infinitesimal model \((\g,H)\) is \emph{reductive}\define{reductive!infinitesimal model}\define{infinitesimal model!reductive} if \(\g\) splits \(\g=\LieH\oplus\redComplement\) into a sum of \(H\)-modules.
To each reductive infinitesimal model, we associate its \emph{affine form} \((X,G)\) where \(X:=\redComplement,G:=H\ltimes\redComplement\): a homogeneous model.
\begin{example}
For example, if \(H\) is compact or semisimple, then every infinitesimal model \((\g,H)\) is reductive, so has an affine form, hence a reductive homogeneous model.
\end{example}
\begin{example}
The affine form of hyperbolic space is Euclidean space.
\end{example}
\begin{example}
The affine form of a smooth affine quadric \(X=\set{a \in \LieSL{2,\C{}}|\det a = 1}\) under its group of regular algebraic automorphisms \(G = \PSL{2}\) is \(X'=\C^3\), \(G'=\C^{\times} \times \C^2\), affine space acted on by rescaling and translations.
\end{example}
\section{Cartan geometries: a broader definition}
A \((\g,H)\)-geometry, or \emph{Cartan geometry infinitesimally modelled on \((\g,H)\)}\define{Cartan!geometry!with infinitesimal model} is a right principal \(H\)-bundle \(\G\to M\) with an \(H\)-equivariant \(\g\)-valued \(1\)-form \(\omega\), the \emph{Cartan connection}, giving a linear isomorphism of all tangent spaces of \(\G\) with \(\g\), and agreeing with the Maurer--Cartan form on the fibers.
Clearly an \((X,G)\)-geometry is a \((\LieG,H)\)-geometry.
All of our discussion above admits obvious generalizations to \((\g,H)\)-geometries, and henceforth we employ this new definition of Cartan geometry.
\section{Mutations of models}
A \emph{mutation}\define{mutation!of infinitesimal models}\define{infinitesimal model!mutation} of infinitesimal models 
\[
\begin{tikzcd}
(\g,H)\arrow[r,"{(\varphi,\Phi)}"]&(\g',H')
\end{tikzcd}
\]
is 
\begin{itemize}
\item
a Lie group morphism \(H\xrightarrow{\Phi}H'\) which yields an isomorphism of Lie algebras \(\LieH\to\LieH'\) and 
\item
a \(\Phi\)-equivariant linear map \(\g\xrightarrow{\varphi}\g'\) agreeing with \(\Phi\) on \(\LieH\).
\end{itemize}
A \emph{mutation}\define{mutation!of homogeneous spaces}\define{homogeneous space!mutation} of homogeneous spaces is a mutation of their infinitesimal models.
\Danger{} Sharpe's definition of mutation is stricter, but he runs into some errors \cite{Lotta2004}.
\section{Mutations of Cartan geometries}
Given a \((\g,H)\)-geometry \(\G\to M\) with Cartan connection \(\omega\), let \(\G':= \amal{\G}{H}{H'}\), so that \(\omega':=\varphi\circ\omega\) descends to be defined on \(\G'\), a \((\g',H')\)-geometry: the \((\varphi,\Phi)\)-\emph{mutation}\define{mutation!of a Cartan geometry}\define{Cartan!geometry!mutation} of \(\G\to M\).
\begin{example}
\[
\begin{array}{r@{\,}ll}
(X,G)&=(\R^n,\Orth{n}\ltimes\R^n)&\text{Euclidean space}\\
(X',G')&=(\mathbb{H}^n,\Orth{n,1})&\text{hyperbolic space}\\
(X'',G'')&=(S^n,\Orth{n+1})&\text{the sphere}
\end{array}
\]
Euclidean space is the usual model of Riemannian geometry.
Hyperbolic space and the sphere are two \emph{other} models of Riemannian geometry, because we have the isomorphic mutations \(H\xrightarrow{\text{id}}H\) on \(H=\Orth{n}\), extending to \(\LieG\) as \(H\)-module isomorphisms:
\[
\begin{pmatrix}
A&v\\
0&0
\end{pmatrix}
\in\LieG
\xrightarrow{\varphi}
\begin{pmatrix}
A&v\\
v^t&0
\end{pmatrix}
\in\LieG'
\text{ or }
\begin{pmatrix}
A&v\\
-v^t&0
\end{pmatrix}
\in\LieG''
\]
\end{example}
\begin{example}
The conformal geometry of the sphere admits no mutation besides isomorphisms of the homogeneous model.
\end{example}
\begin{example}
Spin geometry mutates to Riemannian geometry by \(H=\Spin{n}\mapsto H'=\SO{n}\), \(\LieG\) the same for both.
\end{example}
Mutation preserves and reflects vastness and completeness, because the constant vector fields are the same vector fields.
\begin{example}
Any reductive infinitesimal model has an affine form, a canonical reductive homogeneous model.
Hence all Cartan geometries with reductive infinitesimal model are naturally viewed by mutation as Cartan geometries modelled on their affine forms.
For example, any Cartan geometry modelled on complex projective space, under its group of biholomorphic isometries, is naturally also modelled on its affine form: complex Euclidean space, under the group generated by translations and unitary linear transformations.
\end{example}
\begin{problem}{ch}
If two homogeneous spaces \((X,G)\), \((X',G')\) have a mutation \((\varphi,\Phi)\) of their infinitesimal models, an \((X,G)\)-geometry, say with curvature \(k\), mutates to an \((X',G')\)-geometry, say with curvature \(k'\); use lemma~\vref{lemma:bracket} to relate \(k'\) to \(k\).
\end{problem}
Since \(\omega\) identifies all tangent spaces with \(\g\), \(d\omega\) is uniquely expressed as \(d\omega=k\omega\wedge\omega\) for a function \(\G\xrightarrow{k}\g\otimes\Lm{2}{\g}^*\), the \emph{curvature}\define{curvature!of infinitesimally modelled Cartan geometry} of the Cartan geometry.
Then \(dk=k'\omega\) for a unique \(\G\xrightarrow{k'}\g\otimes\Lm{2}{\g}^*\otimes \g^*\), and so on to all orders, giving \emph{curvature functions},\define{curvature functions} with which we can easily generalize Cartan's rigidity theorems.
If \(H\) is a linear algebraic group and \(\g\) is a regular \(H\)-module then we can similarly generalize Gromov's rigidity theorems.
On the other hand, if the infinitesimal model is Lie, we define the curvature \(k\) by
\[
d\omega + \frac{1}{2} \lb{\omega}{\omega}=\frac{1}{2}k\sigma\wedge\sigma.
\]
\begin{theorem}\label{thm:mutate}
A Cartan geometry with infinitesimal model has constant curvature just when it is a mutation of a flat Cartan geometry with Lie infinitesimal model, perhaps without a homogeneous model.
\end{theorem}
\begin{proof}
Clearly any mutation of a constant curvature geometry has constant curvature.
Suppose that \(\G\to M\) has constant curvature, an \(H\)-invariant element of \(\g\otimes\Lm{2}{\g}^*\).
Differentiate the equation 
\[
d\omega=\frac{1}{2}k\omega\wedge\omega
\]
to get the \emph{Bianchi identity}.\SubIndex{Bianchi identity}
Put in that the curvature is constant, to get the Jacobi identity for the Lie bracket of the mutation, so it is a Lie algebra.
(The same proof works for Cartan geometries whose infinitesimal models are Lie.)
\end{proof}
\begin{corollary}
Suppose that \((X,G)\) is a homogeneous space with invariant Riemannian metric.
Any constant curvature \((X,G)\)-geometry on any compact manifold is a complete Klein manifold, perhaps with some other homogeneous model.
\end{corollary}
\begin{proof}
Mutate to a flat infinitesimally modelled geometry, by theorem~\vref{thm:mutate}.
As above, by compactness of \(H\), if \(X=G/H\), there is a homogeneous model \((X',G')\) for the mutation.
Apply corollary~\vref{corollary:compact.flat}.
\end{proof}
\begin{lemma}\label{lemma:local.dropping}
Take a Cartan geometry \(\G\to M\) with model \((X,G)\) with \(X=G/H\), and with Cartan connection \(\omega\).
Let \(H_0\subset G\) be the local expected structure group of the Cartan geometry.
Consider the infinitesimal model \((\LieG,H_0)\).
Then \(\G\) is covered by open sets \(U_a\), each of which embeds as an open set \(U_a\subset\G_a\) into the total space of a \((\LieG,H_0)\)-Cartan geometry \(\G_a\to M_a\), with Cartan connection \(\omega_a\) pulling back to \(U_a\) to equal \(\omega\).
\end{lemma}
\begin{proof}
Pick a point \(p_0\in\G\).
Pick a submanifold \(M_0\subseteq\G\) containing \(p_0\), with tangent space complementary to the constant vector fields of \(\LieH_0\) at \(p_0\).
Since the constant vector fields form a finite dimensional vector space, varying with choice of point \(p\in\G\) smoothly, the tangent spaces of \(M_0\) at all points \(p\in M_0\) near \(p_0\) are also complementary to the constant vector fields.
We can replace \(M_0\) by a smaller submanifold, so that the constant vector fields are complementary to all tangent spaces of \(M_0\), and so that \(M_0\subset\G\) is an embedded submanifold.

The map
\[
h\in H_0,p\in\G\mapsto\fl{\log h}p\in\G
\]
is defined on a suitable neighborhood of \(1\times\G\) by existence and uniqueness for ordinary differential equations.
For any compact set \(K\subset\G\) around \(p_0\), we can find a relatively compact set open set \(U\subseteq H_0\) so that this map is defined on \(K\times U\).
We replace \(M_0\) by a smaller embedded submanifold so that
\[
h\in U,m\in M_0\xmapsto{\varphi}\fl{\log h}m\in\G
\]
is defined.
Since this map is clearly of full rank at each point \(m\in M_0\), we can shrink \(M_0\) and \(U\) to arrange that the map is a local diffeomorphism.

Denote \(\varphi^*\omega\) also as \(\omega\).
Clearly \(\RT{h}^*\omega=\Ad_h^{-1}\omega\) for all \(h\in U\).
We can therefore extend \(\omega\) unambiguously by the same rule for all \(h\in H_0\).
So \(\omega\) extends uniquely to a \((\LieG,H_0)\)-Cartan connection on \(\G_0:=M_0\times H_0\).
\end{proof}

\section{Infinitesimal models and homogeneous spaces}
When does a Lie infinitesimal model \((\LieG,H)\) have a model \((X,G)\)?
If a model exists, then a model exists with \(X\) connected and simply connected.
Given only the Lie infinitesimal model, let \(\tilde{G}\) be the connected and simply connnected Lie group with Lie algebra \(\LieG\). 
Since \(H\) acts as Lie algebra automorphisms of \(\LieG\), it acts as Lie group automorphisms of \(\tilde{G}\).
The Lie algebra inclusion \(\LieH\subseteq\LieG\) gives an injective and immersive Lie group morphism \(H'\to \tilde{G}\) from a unique connected Lie group \(H'\) with Lie algebra \(\LieH\).
Let \(\tilde{H}\to H^0\) be the universal covering Lie group, covering all connected Lie groups with Lie algebra \(\LieH\).
The covering maps sit in exact sequences
\begin{align*}
1\to \pi_1(H^0)\to&\tilde{H}\to H^0\to 1,\\
1\to \pi_1(H')\to&\tilde{H}\to H'\to 1.
\end{align*}
\begin{theorem}
A Lie infinitesimal model \((\LieG,H)\) has a model \((X,G)\) just when
\begin{itemize}
\item
\(H'\subseteq \tilde{G}\) is closed and
\item
\(\pi_1(H')\subseteq\pi_1(H^0)\) and
\item
the quotient
\[
\pi_1(H^0)/\pi_1(H')\subseteq H'\subseteq\tilde{G}
\]
is an \(H\)-invariant central subgroup of \(\tilde{G}\).
\end{itemize}
If there is a model \((X,G)\) then there is a unique model with \(X\) connected and simply connected, up to isomorphism of homogeneous spaces.
\end{theorem}
\begin{proof}
Suppose that there is a model \((X,G)\).
Every connected Lie group with Lie algebra \(\LieG\) is a quotient \(G^0=\tilde{G}/\Gamma\) for some discrete central subgroup \(\Gamma\subseteq\tilde{G}\).
So \(G^0=\tilde{G}/\Gamma\) for a discrete central subgroup \(\Gamma\) and \(H^0\subseteq G^0\).
So \(\Gamma\subseteq\tilde{G}\cap H'\) is discrete and central in \(\tilde{G}\) with \(H'/\Gamma=H^0\).
So \(\pi_1(H')\subseteq\pi_1(H^0)\).
Then \(\Gamma=\pi_1(H^0)/\pi_1(H')\), a subgroup of \(H'\subseteq\tilde{G}\), a central subgroup of \(\tilde{G}\).
The action of \(H\) on \(\LieG\) by Lie algebra automorphisms determines an action of \(H\) on \(\tilde{G}\) by Lie group automorphisms, since \(\tilde{G}\) is connected.
Since \(H'/\Gamma=H^0\) is closed in \(G^0\), \(H'\) is closed in \(\tilde{G}\).
Since \(\Gamma\) is central it commutes with the \(H^0\)-action.
Since we can quotient, \(\Gamma\) is \(H\)-invariant: our three conditions above are satisfied.
Write the action of \(H\) on \(G^0\) as \(G^0\xrightarrow{\alpha_h}G^0\) for \(h\in H\).

We can assume \(X\) is connected and simply connected, so \(\pi_0(H)=\pi_0(G)\) by homotopy exact sequence (\cite{Steenrod:1999} Theorem 17.4, \cite{tomDieck:2008} p. 123 Theorem 6.1.2), hence \(G=G^0\times\pi_0(H)\) as a manifold.
Choose a set theoretic section \(a\in\pi_0(H)\mapsto \bar{a}\in H\) of the obvious \(H\to\pi_0(H)\).
Naturally we can arrange, if \(b=a^{-1}\) then \(\bar{b}=\bar{a}^{-1}\) and \(\bar{1}=1\), but we can't expect multiplication to match.
Every element of \(G\) is a unique product \(g\bar{a}\) for some \(g\in G^0\), \(a\in\pi_0(H)\).
So as a set, \(G=G^0\times\pi_0(H)\), a manifold with the product smooth structure, since \(\pi_0(H)\) is a discrete group.
Let
\[
h_{ab}:=\bar{a}\bar{b}(\overline{ab})^{-1}\in H^0.
\]
In particular, \(h_{1b}=h_{a1}=1\) for any \(a,b\in\pi_0(H)\).

If we take \(g_0,g_1\in G^0\) then
\begin{align*}
(g_0\bar{a})(g_1\bar{b})
&=g_0(\bar{a}g_1\bar{a}^{-1})\bar{a}\bar{b},
\\
&=(g_0\alpha_{\bar{a}}(g_1)h_{ab})\overline{ab}.
\end{align*}
Expand out to find, for any \(g\in G^0\),
\[
(g\bar{a})^{-1}=\alpha_{\bar{a}^{-1}}(g^{-1})\bar{a}^{-1}.
\]

The \(H\)-action on \(\tilde{G}\) is by diffeomorphisms, and so therefore is the \(H\)-action on \(G^0\). 
The product smooth structure on the set \(G:=G^0\times\pi_0(H)\) is invariant under any group structure on \(G\) extending that of \(H\) and \(G^0\), a smooth manifold structure.
Define operations
\[
(g_0,a)(g_1,b)=
(g_0\alpha_{\bar{a}}(g_1)h_{ab},ab),
\]
and inverse by
\[
(g,a)^{-1}
=
(\alpha_{\bar{a}^{-1}}(g^{-1}),a^{-1}).
\]
This is the Lie group structure of \(G\).

Conversely, suppose that the three conditions of the theorem are satisfied.
Define \(\Gamma:=\pi_1(H^0)/\pi_1(H')\) and \(G^0:=\tilde{G}/\Gamma\).
Repeat our steps above: pick any set theoretic section of \(H\to\pi_0(H)\) and define the binary operation above.
Note that the inverse operation satisfies the inverse axiom of a group.
We only need to test for associativity.
Expand out to check that associativity is precisely
\[
h_{ab}h_{ab,c}=\alpha_{\bar{a}}(h_{bc})h_{a,bc}.
\]
But this occurs entirely inside \(H^0\), not inside \(G^0\), so the associativity is guaranteed, and the operations form a group: there is a model.

(To be fancy, we can define expressions
\[
z_{abc}=\alpha_{\bar{a}}(h_{bc})h_{a,bc}h^{-1}_{ab,c}h^{-1}_{ab}.
\]
If we let \(Z=Z_{G^0}\) be the center of \(G^0\), then \(z_{abc}\in Z\) and \(z=\set{z_{abc}}\) is a \(3\)-cocycle, determining a cohomology class \([z]\in H^3(\pi_0(H),Z)\) and vanishes in group cohomology just when the multiplication is that of a group \cite{Maclane:1995} p. 124, IV.8.
But we don't need to do this: our \(z\) is valued in the center of \(H^0\), and we already have a group structure on \(H\).)
\end{proof}

\chapter{Automorphism group orbits}
\section{Cutting out the orbits in the bundle}
Take a Cartan geometry \(\G\to M\) with infinitesimal model \((\g,H)\).
On \(\G\), write the curvature \(k\), and all of its derivatives
\[
dk=k'\omega, dk'=k''\omega, \dots
\]
in terms of their component functions in an arbitrary basis of \(\LieH\) and \(\g\), which we call \emph{curvatures}.\define{curvatures}
On a dense open subset \(\Greg\subseteq\G\) the curvatures have locally constant rank,\ an \(H\)-invariant condition, so projecting to a dense open set in \(\Mreg\subseteq M\), the \emph{regular set}.\define{regular!set}
\begin{example}
Take this ellipsoid of revolution
\[
\includegraphics[width=2cm]{ellipsoid}
\]
with its Riemannian metric as a Cartan geometry.
Its Gauss curvature \(G\) has nonzero differential except at the north and south pole (maxima) and the equator (minima).
Any invariant functions are invariant under the isometry of reflection in the horizontal plane, and are functions of \(G\).
The regular set \(\Mreg\) is the complement of the poles and equator.
\end{example}
\begin{theorem}%
[Regular partition \cite{Cartan:161} p. 122, \cite{Olver:1995} chapter 14]%
\label{theorem:smooth.strat.I}%
\define{theorem!Cartan regular partition}
\define{Cartan!regular partition theorem}
For any Cartan geometry \(\G\to M\), above each component of the regular set, the local automorphism orbits in \(\G\) are the level sets of the component functions of the curvature and finitely many of its derivatives, and these level sets are the fibers of a smooth submersion.
The local infinitesimal automorphisms have flows acting locally transitively on those level sets. 
\end{theorem}
\begin{example}
The automorphism group orbits lie inside these level sets, but might have lower dimension than any of the level sets.
The only global isometry of a sphere punctured at \(3\) points in general position is the identity, but there are three dimensions of local automorphisms: many local automorphisms do not extend globally.
\end{example}
\begin{proof}
Replace \(M\) by \(\Mreg\) without loss of generality.
Consider the open set \(\G'\subseteq\G\) on which a particular collection \(\bar{k}=(k^1,\dots,k^{\ell})\) of chosen curvatures have maximal rank among curvatures.
Clearly \(\G\) is the union of such open sets.
This open set \(\G'\) is invariant under local automorphisms.
By the implicit function theorem, every curvature is locally a function of these curvatures.
If we have two points of \(\G'\) at which all curvatures of first order agree locally as functions of \(\bar{k}\), then their derivatives agree by just differentiating them, so all curvatures agree near those two points of \(\G'\).

Let \(R\subseteq\G\times\G\) be the set of pairs \((p,q)\in\G\times\G\) which have the same values of all curvatures of all orders.
The rate of change \(A_{\G}k\) of any curvature along any constant vector field is given in terms of the  derivative of \(k\), so is encoded in curvatures and components of \(A\).
Hence the rank of each curvature is encoded in the curvatures.
In other words, every point of \(R\) lies in some open set \(\G'\times\G'\).
If some local automorphism of the Cartan geometry takes a point \(p\in\G\) to a point \(q\in\G\) then \((p,q)\in R\).
Note that \(R\subseteq\G\times\G\) is locally a submanifold of codimension equal to the maximal rank of the curvatures and projects to either factor of \(\G\) by a submersion, since \(R\) is cut out by equations specifying invariants of \(q\) as functions of those of \(p\), or vice versa, so we can specify \(p\) or \(q\) freely.

The curvatures of at most a given order transform among one another, under \(H\)-action.
Those which form a maximal rank collection on our open set \(\G'\) transform into linear combinations giving another choice of maximal rank collection on each open set \(\G'h\), for \(h\in H\).
So \(R\) is invariant under diagonal \(H\)-action.

At any points \(p,q\in\G\) with \((p,q)\in R\), curvature at \(p\) agrees with curvature at \(q\), so the brackets of constant vector fields agree.
The expressions \(A_{\G}\bar{k}\) are expressed in terms of the curvatures, by taking derivatives, so agree at \(p\) and at \(q\).
Hence the vector fields
\[
(A_{\G},A_{\G}), \quad A\in\LieG
\]
are tangent to \(R\).

From the Orbit Theorem~\vref{thm:orbit}, the orbits of these vector fields are submanifolds of \(R\).
Since the brackets match up at the corresponding points, each orbit projects to an orbit in each factor \(\G\) of \(\G\times\G\).
But the vector fields \(A_{\G}\) span the tangent spaces of \(\G\).
So each orbit is the leaf of a foliation, with tangent space the span of these \((A_{\G},A_{\G})\), locally the graph of a map between open sets of \(\G\) matching the Cartan connections.
The projection of \(\G\times\G\) to either factor of \(\G\) projects each leaf by a local diffeomorphism, invariant under the flows of the constant vector fields.
So each leaf lies in a unique \(H\)-folio, projecting to \(\G\) by \(H\)-equivariant local diffeomorphism, matching constant vector fields, so the graph of a local automorphism.

Take a curve smooth \(p(t)\in\G\) consisting of points all of which lie in the same local automorphism group orbit.
There are local automorphisms \(\varphi_t\) with \(\varphi_t(p_0)=p(t)\), where \(p_0:=p(0)\), hence a path \((p_0,p(t))\in R\).
Conversely, any path \((p_0,p(t))\in R\) with \(p(0)=p_0\) has all \(p(t)\) in the local automorphism group orbit.
Each local automorphism \(\varphi_t\) is locally uniquely determined: its graph is an open subset of a unique leaf, so it varies smoothly with \(t\).
Define a vector field \(Z\) by 
\[
Z(p):=\left.\frac{d}{dt}\right|_{t=0}\varphi_t(p),
\]
defined on some open set around \(p_0\).
We can assume that this open set is \(H\)-invariant, as \(R\) is.
Clearly \(Z\) is a local infinitesimal automorphism, and every local infinitesimal automorphism arises this way.
Since the local automorphisms have orbits in \(\G\) the level sets of some curvatures, this is the only constraint on \(p'(0)=Z(p_0)\), i.e. the infinitesimal automorphisms at each point \(p_0\in\G\) span the tangent spaces of the level sets of the curvatures.
\end{proof}
\Danger{} Cartan \cite{Cartan:161} and Olver \cite{Olver:1995} give proofs of this theorem, but they state and prove it more intuitively and less precisely.
\begin{example}
A disjoint union of infinitely many unit circles with their usual Riemannian metrics has local automorphisms with orbits of dimension \(1\), but an infinite dimensional group of global isometries, with uncountably many components.
\end{example}
\begin{example}
The conformal structure of the real projective space of dimension \(3\) or more has the same local automorphisms as the sphere, all of which extend globally on the sphere, but extend to the real projective space as two-valued maps typically.
Again there are ``too many'' local automorphisms.
\end{example}
\begin{example}
Take a flat torus of dimension \(3\) or more, and perturb it in a little bump to be not conformally flat.
The local automorphisms away from the bump are those of the sphere, but extend only away from the bump, typically with infinite multivaluedness: ``too many'' local automorphisms.
\end{example}
\begin{theorem}%
[Analytic partition]%
\label{thm:Cartan.analytic.partition}%
\define{theorem!Cartan analytic partition}
\define{Cartan!analytic partition theorem}
If \(\G\to M\) is a real analytic Cartan geometry on a compact manifold \(M\), some finite set of curvatures on \(\G\) distinguish which points can be carried to which other points by local automorphisms of the Cartan geometry.
So the local automorphisms act with closed real analytic subsets of \(\G\) as orbits.
Away from a closed nowhere dense analytic subset of \(M\) these local automorphism orbits form the fibers of a real analytic submersion of \(\G\).
\end{theorem}
\begin{proof}
The same proof as for the previous theorem works on a dense open set.
Where it fails, curvature degenerates, so it fails on an analytic subset of \(\G\), but \(H\)-invariant, so an \(H\)-bundle over an analytic subset of \(M\), hence a compact subset.
Try again using the same tactic, succeeding to find suitable invariants on an open set, the complement of a smaller dimensional analytic set.
Repeat by induction on dimension.
\end{proof}
By the same reasoning, given two real analytic Cartan geometries, on compact manifolds, some finite set of curvatures determines which points of one can be carried by local isomorphisms to which points of the other: the local equivalence problem can be solved.

Recall that \(\Aut\to\G\to\G/\Aut\) is a principal bundle, so there are smooth invariants (precisely the smooth functions on \(\G/\Aut\)) that distinguish automorphism group orbits, but it is not clear how to find them, since local invariants can only distinguish local automorphism group orbits.
By corollary~\vref{corollary:analytic.vast},
\begin{corollary}\label{corollary:analytic.auts}
On any connected manifold \(M\), the automorphism group orbits in \(\G\) of a vast, real analytic Cartan geometry \(\G\to M\) with trivial holonomy of infinitesimal automorphisms are the local automorphism group orbits.
\end{corollary}
\subsection{Decorations}\SubIndex{decoration}\SubIndex{Cartan!geometry!decorated}
Cartan's rigidity theorems extend trivially to decorated Cartan geometries: just replace the phrase \emph{Cartan geometry} with the phrase \emph{decorated Cartan geometry}, and interpret the terms \emph{constant vector field} and \emph{complete} as referring to the underlying Cartan geometry.
The component functions of the decorations (in any embedding of the target manifold in Euclidean space) together with the curvatures have locally constant rank on a dense open set.
\begin{example}
The Cartan connection identifies any vector field \(Z\) on \(M\) with an \(H\)-equivariant map \(A:=Z\hook\omega\colon\G\to\LieG/\LieH\), a decoration.
Similarly, any tensor field is a decoration, and any finite dimensional Lie algebra action.
\end{example}
\begin{example}
If \(\Gamma\subseteq\Aut\) is a closed subgroup of the automorphism group of a decorated Cartan geometry on a connected manifold, the quotient map \(\G\to\Gamma\backslash\G\) is another decoration, and the automorphism group, with this decoration included, is \(\Gamma\).
\end{example}
\section{Cutting out the orbits in the base manifold}
\subsection{Scalar invariants}
A homogeneous space \((X,G)\) is \emph{infinitesimally algebraic}\define{infinitesimally algebraic} if \(X=G/H\) and \(H\) admits the structure of a linear algebraic group, refining its structure as a Lie group, so that \(\LieG\) is a regular \(H\)-module.
The various associated \(H\)-modules in which we find curvature and its derivatives form an infinite sequence of regular modules.
To any finite sum of these \(H\)-modules, we associate rational invariants, as in theorem~\vref{theorem:Rosenlicht}.
Recall that their values distinguish \(H\)-orbits except along a ``bad set'', an invariant algebraic hypersurface on which we can find further rational invariants, and so on.
On any \((X,G)\)-geometry \(\G\to M\), plugging in curvature and its derivatives into these rational invariants defines scalar invariants on \(\G\), perhaps not defined everywhere (and perhaps some not defined anywhere!), rational in the curvature and its derivatives.
These scalar invariants are \(H\)-invariant, so descend to functions on \(M\), each defined on some open (perhaps empty!) set.
As we raise the number of derivatives we allow ourselves to take, we raise the number of scalar invariants we must take to separate \(H\)-orbits.

On the scaffold, we can separate out orbits to as high an order as we like, using invariant \emph{polynomials} in the curvature and these derivatives, so smooth functions everywhere defined on the scaffold.
\subsection{Gromov's rigidity theorems}
The scalar invariants of a given order tell us which points of \(M\) can be matched up to one another over which some bundle isomorphism matches up Cartan connections to some given order.
Let \(\Mreg\subseteq M\) be the set of points above which the curvature and derivatives of some order stay locally inside some nested bad set and outside of the next one and some particular set of invariants are thereby defined near there, and collectively have maximal rank among our scalar invariants defined near there.
Clearly \(\Mreg\) is invariant under local automorphisms.
By Sard's theorem \cite{Sard:1965}, \(\Mreg\subseteq M\) is a dense open set.
\begin{theorem}%
[Regular partition 
\cite{Benoist:1997}, 
\cite{dAmbra/Gromov:1991} p. 72, 5.14.B, 
\cite{Melnick:2011}, 
\cite{Pecastaing:2016}, 
\cite{Sternberg:1983} p. 347 corollary 4.1, 
\cite{Zeghib:2000}]%
\label{theorem:smooth.strat}%
\define{theorem!Gromov regular partition}
\define{Gromov's regular partition theorem}
Take a Cartan geometry \(\G\to M\) with infinitesimally algebraic model.
Over each component of \(\Mreg\), the local automorphism orbits are the level sets of some scalar invariants with linearly independent differentials.
The infinitesimal automorphisms defined above open sets of \(\Mreg\) have flows acting locally transitively on those level sets. 
\end{theorem}
\begin{proof}
Return to the proof of Cartan's theorem~\vref{theorem:smooth.strat.I}.
As before, we can replace \(M\) by \(\Mreg\), so assume \(M=\Mreg\).
We can even replace \(M\) by an open set on which the relevant nested bad set is constant.
As in the proof of Cartan's theorem, locally above \(M\) we can choose some curvatures of locally maximal rank.
We can replace \(M\) by an open set in which this rank is constant.
Consider the order of our chosen curvatures; call it \emph{the} order.\define{order!of a Cartan geometry}\define{Cartan!geometry!order}
Consider the curvature and its derivatives up to the order.
Above each point of \(M\), they lie in some nested bad set, where our scalar invariants might not separate orbits and on which some of them might not be defined.
Above \(M\), we can also pick some scalar invariants, with locally maximal rank, which we can assume is maximal rank, in some nested bad set, but perhaps expressed using more derivatives of curvature.
Since these are locally of maximal rank, every one of our scalar invariants, of any order, defined somewhere in our open set is locally a function of the scalar invariants we picked, by the implicit function theorem.
So we can take the orders to be the same, and hence the nested bad sets are the same.
Our finite set of scalar invariants separate points of that nested bad set.
Two points of \(M\) have the same values for those invariants just when there are points above them in \(\G\) where all curvatures up to the order agree, so points of \(R\).
So all curvatures of all orders match at those points just when our maximal rank collection of scalar invariants do.
Apply theorem~\vref{theorem:smooth.strat.I}.
\end{proof}
A fundamental theorem of Kruglikov and Lychagin generalizes this work and describes syzygies and invariant differential operators among the differential invariants \cite{Kruglikov.Lychagin:2016}.
\begin{theorem}%
[Analytic partition 
\cite{Benoist:1997}, 
\cite{dAmbra/Gromov:1991} p. 35, 1.11.A, 
\cite{Melnick:2011}, 
\cite{Pecastaing:2016}, 
\cite{Zeghib:2000}]%
\label{thm:Gromov.analytic.partition}
\define{theorem!Gromov analytic partition}
\define{Gromov's analytic partition theorem}
Suppose that \(\G\to M\) is a real analytic Cartan geometry, with infinitesimally algebraic model, on a compact manifold \(M\).
Then some finite set of scalar invariants distinguish which points of \(M\) can be carried to which other points by local automorphisms of the Cartan geometry.
So the local automorphisms of open sets act with locally closed real analytic subsets of \(M\) as orbits.
Away from a closed proper analytic subset (in particular, nowhere dense), these local automorphism orbits form the fibers of a real analytic submersion.
At least one orbit in \(M\) of the local automorphisms is a compact real analytic embedded submanifold.
\end{theorem}
\begin{proof}
The same proof as for the previous theorem works on a dense open set.
Where it fails, curvature degenerates into some bad set, so it fails on an analytic subset.
On the smooth points of that set, we try again using the same tactic.
Repeat by induction.
When we get to the last induction step, the orbit is closed.
\end{proof}
\subsection{The open dense theorem}
\begin{corollary}%
[The open-dense theorem 
\cite{Benoist:1997}, 
\cite{dAmbra/Gromov:1991} p. 35, 1.11.A, 
\cite{Melnick:2011}, 
\cite{Pecastaing:2016}, 
\cite{Zeghib:2000}]%
\define{theorem!Gromov open dense}
\define{Gromov's open dense theorem}
Take a Cartan geometry \(\G\to M\) with infinitesimally algebraic model.
If local automorphisms of \(\G\to M\) act on \(M\) with a somewhere dense orbit, then that orbit is open in \(M\).
In particular, if some set of automorphisms acts on \(M\) with a dense orbit, then it generates a subgroup of the automorphism group whose closure in the automorphism group has a dense open orbit.
\end{corollary}
\begin{example}
Surprisingly, there are real analytic torsion-free affine connections on the plane whose automorphism group acts transitively on a dense open set, but not transitively \cite{Dumitrescu:2014}.
Their scalar invariants, expressed in curvature and derivatives of some perhaps very high order, are constant on that dense open set.
They extend continuously (as constants) to the whole plane.
So we might expect the automorphism group to act transitively.
But in fact, the automorphism group orbits in \(\G\) rise higher and higher, in any local trivialization, roughly like hyperbolae:
\[
\begin{tikzpicture}
\tikzset{
    elli/.style args={#1:#2and#3}{
        draw,
        shape=ellipse,
        rotate=#1,
        minimum width=2*#2,
        minimum height=2*#3,
        outer sep=0pt,
    }
}
%
% #1 optional parameters for \draw
% #2 angle of rotation in degrees
% #3 offset of center as (pointx, pointy) or (name-o-coordinate)
% #4 length of plus (semi)axis, that is axis which hyperbola crosses
% #5 length of minus (semi)axis
% #6 how much of hyperbola to draw in degrees, with 90 you’d reach infinity
%
\newcommand\tikzhyperbola[6][thick]{%
    \draw [#1, rotate around={#2: (0, 0)}]
        plot [variable = \t, samples=1000, domain=-#6:#6] ({#4 / cos( \t )}, {#5 * tan( \t )});
    \draw [#1, rotate around={#2: (0, 0)}]
        plot [variable = \t, samples=1000, domain=-#6:#6] ({-#4 / cos( \t )}, {#5 * tan( \t )});
}
\fill[white] (-3,-3) rectangle (3,3);
\fill[gray!20] (-2,-2) rectangle (2,2);
\tikzhyperbola[thick]{45}{(0,0)}{.5cm}{.5cm}{70}
\node[right] (G) at (2,0) {\(\G\)};
\draw[very thick,blue] (0,-2) -- (0,2);
\draw[very thick] (-2,-2.5) -- (2,-2.5);
\node[right] (M) at (2,-2.5) {\(M\)};
\draw[-latex] (G) to (M);
\fill[blue] (0,-2.5) circle (1.5pt);
\end{tikzpicture}
\]
The boundary of our dense open set consists of points where the curvature and derivatives enter a bad set.
\end{example}
\begin{corollary}
Take a vast real analytic Cartan geometry \(\G\to M\) with infinitesimally algebraic model on a connected manifold \(M\).
Suppose that \(M\) has trivial infinitesimal holonomy and that some set of local automorphisms act on \(M\) with a dense orbit.
Then the automorphism group acts on \(M\) with a dense open orbit.
\end{corollary}
\begin{proof}
Apply corollary~\vref{corollary:analytic.auts} and theorem~\vref{thm:Cartan.analytic.partition}.
\end{proof}
\begin{corollary}
Suppose that \(M\) is a compact connected manifold.
Take a real analytic Cartan geometry \(\G\to M\) with infinitesimally algebraic model and trivial holonomy of infinitesimal automorphisms.
Suppose that some set of local automorphisms act on \(M\) with a somewhere dense orbit.
Then the automorphism group acts on \(M\) with a dense union of finitely many open orbits.
If the geometry is holomorphic, the automorphism group acts on \(M\) with a dense open orbit.
\end{corollary}
\begin{proof}
The local automorphisms are driven by local infinitesimal automorphisms.
By theorem~\vref{theorem:smooth.strat}, these local infinitesimal automorphisms act locally transitively on a submanifold of \(M\) which contains a dense open subset of \(M\).
By theorem~\vref{theorem:analytic.extension}, the local infinitesimal automorphisms extend to global infinitesimal automorphisms.
These are complete vector fields, because \(M\) is compact, so they generate a Lie subgroup of the automorphism group, acting locally transitively on \(\Mreg\).
The irregular points \(M-\Mreg\) form a real analytic set, so separate \(M-\Mreg\) locally into finitely many open sets, hence globally by compactness.
If the geometry is holomorphic, the complement \(M-\Mreg\) is a complex analytic set, so has real codimension \(2\) or more, so \(\Mreg\) is connected.
\end{proof}
\begin{theorem}
Suppose that \(\G\to M\) is a Cartan geometry.
Suppose that the model admits an invariant metric.
Then the Cartan geometry is homogeneous if and only if its automorphisms have a dense orbit.
If the geometry is real analytic then the Cartan geometry is homogeneous if and only if its automorphisms have an orbit somewhere dense in every component of \(M\).
\end{theorem}
\begin{proof}
The stabilizer of the model is compact, so linear algebraic.
The scalar invariants can be taken as polynomials in the curvature and derivatives, so smooth \(H\)-invariant functions on \(\G\), so the scalar invariants are constant on a dense open set in \(M\) just when constant on \(\G\).
\end{proof}
\begin{example}
The automorphisms of a Riemannian manifold have a dense orbit just when it is homogeneous
\end{example}
\begin{example}
If the automorphisms of a conformal structure on a manifold \(M\) have a dense orbit in the scaffold \(M_c\) then the conformal structure is homogeneous.
The scaffold is the space of pairs \((m,g)\) of point \(m\in M\) and inner product \(g\) on \(T_m M\) in the conformal class.
So if we can approximate any point and inner product from any other by automorphisms, our conformal structure is homogeneous.
\end{example}
\begin{example}
If \((X,G)\) is a strong and effective homogeneous space with finitely many components and the stabilizer of a point has finitely many components then \(G\) is its automorphism group so \(G\) has finitely many components.
\end{example}
\begin{theorem}[Gromov \cite{Gromov1988} p. 100 theorem 3.5.C]
Suppose that \(M\) is a compact manifold with a real analytic Cartan geometry \(\G\to M\) and its model is infinitesimally algebraic with finite kernel.
Suppose that \(\G\to M\) has trivial holonomy of infinitesimal automorphisms on each component of \(M\).
The automorphism group of \(\G\to M\) has finitely many components, and hence so do all of its orbits in \(\G\) and in \(M\).
At least one orbit in \(M\) is a compact real analytic embedded homogeneous submanifold.
\end{theorem}
\begin{example}
Take any real analytic Cartan geometry on the \(3\)-sphere.
Suppose that its model is infinitesimally algebraic with finite kernel.
Its automorphism group has finitely many components and 
\begin{itemize}
\item
a finite orbit or 
\item
an orbit which is real analytic compact curve or
\item
an orbit which is a real analytic compact surface or
\item
acts transitively.
\end{itemize}
\end{example}
\begin{example}
The \(\Z\)-action on Euclidean space generated by \(x\mapsto x/2\) is conformal, so extends to a conformal action on the sphere; it is a closed subgroup of the automorphism group of a real analytic Cartan geometry, and it is thus the automorphism group of a decorated Cartan geometry, but is not the automorphism group of any real analytic Cartan geometry with infinitesimally algebraic model with finite kernel.
\end{example}
\begin{proof}
Without loss of generality, \(M\) is connected and the model \((X,G)\), \(X=G/H\), is effective.
Inside a fiber of \(\G\to M\), say over a point \(m_0\), pick some point \(p_0\) and consider the subset of points of that fiber on which the curvatures (in the sense of the proof of theorem~\vref{theorem:smooth.strat}) agree with those at \(p_0\).
This subset is cut out by finitely many curvatures near each point of that fiber.
Each fiber of \(\G\to M\) is a copy of \(H\), i.e. is a linear algebraic group, as we saw in lemma~\vref{lemma:stabilzer.closed}.
That subset intersects that fiber in a linear algebraic subgroup: the elements of \(H\) preserving those various curvatures.
So our subset is algebraic (and also a real analytic submanifold) and has therefore finitely many components in each fiber of \(\G\to M\).

Local infinitesimal automorphisms of \(\G\to M\) act freely and transitively on the components of that subset, as in our analytic partition theorem.
These local infinitesimal automorphisms extend, by simple connectivity, along \(M\).
They are complete vector fields because \(M\) is compact.
So they generate a subgroup in \(\Aut^{m_0}\).
On the other hand, the Lie algebra of \(\Aut^{m_0}\) consist of the infinitesimal automorphisms vanishing at \(m_0\) on \(M\), hence has flow preserving our subset, so belonging to that subgroup.
Therefore the identity component \(\Aut^{m_0}\) is precisely the identity component of that linear algebraic group, and \(\Aut^{m_0}\) has finitely many components.
\end{proof}
Perhaps the same result is true for compact manifolds with finite holonomy of infinitesimal automorphisms.
\begin{corollary}
On a compact manifold, take any real analytic Cartan geometry with infinitesimally algebraic model and finite kernel and trivial holonomy of infinitesimal automorphisms.
If there is no nonzero global infinitesimal automorphism then the automorphism group is finite.
\end{corollary}
\begin{example}
As we have seen, there are complete Cartan geometries, even complete Riemannian geometries on surfaces, whose automorphism group is one dimensional with infinitely many components, or discrete with infinitely many components.
\end{example}
{\centering\tiny
\begin{tabular}{p{5cm}p{5cm}}
\includegraphics[width=5cm]{delaunay.png} & \includegraphics[width=5cm]{Schwarz_P_surface} \\
GeometrieWerkstatt CC BY-NC-SA-3.0 &
By Anders Sandberg \\ & Own work \\ & CC BY-SA 3.0 \\
& {\verb!commons.wikimedia.org/w/index.php?curid=20737176!}
\end{tabular}
\par
}
\begin{example}
In dimension \(3\) or more, we can similarly construct conformal structures with infinite discrete automorphism group; they admit no real analytic conformal compactification; in dimension \(2\), every conformal structure conformally compactifies by the Poincar\'e uniformization theorem \cite{Saint-Gervais:2010}.
\end{example}
\begin{example}
For any lattice \(\Lambda\subseteq V\) in a finite dimensional vector space \(V\), let \(\GL{\Lambda}\subseteq\GL{V}\) be its linear automorphisms.
Affine space \((X,G)=(V,\GL{V}\ltimes V)\) has quotient space \(M:=V/\Lambda\) the torus, with automorphism group \(\GL{\Lambda}\ltimes(V/\Lambda)\) with infinitely many components.
To be concrete, for \(n=2\) these automorphisms include the shearing transformations like
\[
  \begin{tikzpicture}[scale=.5]
    \coordinate (Origin)   at (0,0);
    \coordinate (XAxisMin) at (-3,0);
    \coordinate (XAxisMax) at (5,0);
    \coordinate (YAxisMin) at (0,-2);
    \coordinate (YAxisMax) at (0,5);
    \draw [thin, gray,-latex] (XAxisMin) -- (XAxisMax);% Draw x axis
    \draw [thin, gray,-latex] (YAxisMin) -- (YAxisMax);% Draw y axis

    \clip (-3,-2) rectangle (5cm,5cm); % Clips the picture...
    \pgftransformcm{1}{0.6}{0.7}{1}{\pgfpoint{0cm}{0cm}}
          % This is actually the transformation matrix entries that
          % gives the slanted unit vectors. You might check it on
           % MATLAB etc. . I got it by guessing.
    \coordinate (Bone) at (0,2);
    \coordinate (Btwo) at (2,-2);
    \draw[gray!50,style=help lines] (-14,-14) grid[step=2cm] (14,14);
          % Draws a grid in the new coordinates.
          %\filldraw[fill=gray, fill opacity=0.3, draw=black] (0,0) rectangle (2,2);
              % Puts the shaded rectangle
%    \foreach \x in {-7,-6,...,7}{% Two indices running over each
%      \foreach \y in {-7,-6,...,7}{% node on the grid we have drawn 
%        \node[draw,gray,circle,inner sep=.5pt,fill] at (2*\x,2*\y) {};
%            % Places a dot at those points
%      }
%    }
%    \draw [ultra thick,-latex,red] (Origin)
%        -- (Bone) node [above left] {};
%    \draw [ultra thick,-latex,red] (Origin)
%        -- (Btwo) node [below right] {};
%    \draw [ultra thick,-latex,red] (Origin)
%        -- ($(Bone)+(Btwo)$) node [below right] {};
%    \draw [ultra thick,-latex,red] (Origin)
%        -- ($2*(Bone)+(Btwo)$) node [above left] {};
    \filldraw[fill=brown, fill opacity=0.2, draw=black] (Origin)
        rectangle ($2*(Bone)+(Btwo)$);
    \filldraw[fill=blue, fill opacity=0.2, draw=black] (Origin)
        -- ($(Bone)$) -- ($(Bone)+(Btwo)$) -- ($(Btwo)$) -- cycle; 
    %\draw [thin,-latex,red, fill=gray, fill opacity=0.3] (0,0)
        % -- ($2*(0,2)+(2,-2)$)
        % -- ($3*(0,2)+2*(2,-2)$) -- ($(0,2)+(2,-2)$) -- cycle;
  \end{tikzpicture}
\]
given by the matrix
\[
\begin{pmatrix}
1&2\\
0&1
\end{pmatrix}.
\]
\end{example}
Gromov proves more results like the theorems above, for a class of geometric structures on manifolds much larger than the class of Cartan geometries.
His papers are difficult; proofs easier to read of these and many related results are now available \cite{An:2012,Benoist:1997,Feres:2002,Feres.Lampe:2000,Melnick:2011,Pecastaing:2016,Zeghib:2000}.
It does not seem to be known whether there is a compact manifold with finite fundamental group admitting a real analytic Cartan geometry whose automorphism group has infinitely many components.

\chapter{\texorpdfstring{Soldering forms and $G$-structures}{Soldering forms and G-structures}}\label{section:soldering}
\section{The problem}
Take a Lie infinitesimal model \((\LieG,H)\).
Take a manifold \(M\) of dimension \(\dim{\LieG}-\dim{H}\).
We wonder when a principal right \(H\)-bundle \(\Bun\to M\) is the principal bundle of a \((\LieG,H)\)-geometry.
\section{Soldering forms on bundles}
Take a Lie group \(H\) and an \(H\)-module \(V\).
Take a manifold \(M\) of the same dimension as \(V\).
A \emph{soldering form}\define{soldering form} \(\sigma\) on a principal \(H\)-bundle \(\Bun\to M\) is an \(H\)-equivariant \(1\)-form \(\sigma\in\Omega^1_{\Bun}\otimes V\) vanishing precisely on the vertical vectors, i.e. the tangent spaces of the fibers of \(\Bun\to M\).
\begin{problem}{soldering.equiv}
Explain why a soldering form is equivalent to a choice of vector bundle isomorphism \(TM\cong\amal{\Bun}{H}{V}\).
\end{problem}
\section{The frame bundle}
Take a finite dimensional vector space \(V\) and a manifold \(M\) with \(\dim M=\dim V\).
A \emph{\(V\)-valued frame}\define{frame} on \(M\) at a point \(m\in M\) is a linear isomorphism \(T_m M\xrightarrow{u}V\).
The \(V\)-valued \emph{standard frame bundle}\define{frame bundle!standard}\define{standard frame bundle} \(\framebundle{M}\) of \(M\) is the set of pairs \((m,u)\) so that \(m\in M\) and \(u\) is a \(V\)-valued frame at \(m\).
Let \((m,u)\in\framebundle{M}\xmapsto{\pi}m\in M\).
The frame bundle is thus a principal right \(\GL{V}\)-bundle under the right action
\[
(m,u)h:=(m,h^{-1}u),
\]
also denoted \(\RT{h}(m,u)\).
The frame bundle has a soldering form: if we let \(H:=\GL{V}\) and \(\Bun:=\framebundle{M}\), we have a soldering form defined by
\[
v\hook\sigma:=u(\pi'(m,u)v),
\]
for any \(v\in T_{(m,u)}\framebundle{M}\).
\section{First order structures}
Suppose that \(H\) is a Lie group and \(V\) is a finite dimensional \(H\)-module.
Take a manifold \(M\) of the same dimension as \(V\).
An \emph{\(H\)-structure} on \(M\) is a principal right \(H\)-bundle \(\Bun\) and an \(H\)-equivariant bundle map
\(
\Bun\to\framebundle{M}.
\)
\Danger{} an \(H\)-structure is \emph{not} the same use of the word \emph{structure} as in \emph{locally homogeneous structure}.
Any \(H\)-structure determines a soldering form \(\sigma\) by pullback of the soldering form from the frame bundle.
If \(H\subseteq\GL{V}\) is an immersed subgroup then each \(H\)-structure is an immersed subbundle of the frame bundle.
If \(H\subseteq\GL{V}\) is a closed subgroup then each \(H\)-structure is an embedded subbundle of the frame bundle, in which case \(H\)-structures are precisely sections of the fiber bundle \(\framebundle{M}/H\to M\).
\section{Quotienting out trivia}
Let \(H_1\subseteq H\) be the subgroup acting trivially on \(V\), which we recall from lemma~\vref{lemma:eff.hom}.
Every \(H\)-structure descends to an \(H/H_1\)-structure \(\Bun/H_1\to\framebundle{M}\); most often it is this \(H/H_1\)-structure one encounters in the literature.
Note that \(H/H_1\subseteq\GL{V}\) is an immersed subgroup.
\section{\texorpdfstring{Soldering form $\cong$ first order structure}{Soldering form equivalent to first order structure}}
The soldering form on the frame bundle is universal in the following sense.
Suppose that \(\Bun\xrightarrow{\pi}M\) is a principal right \(H\)-bundle with a soldering form \(\sigma\).
Pick a point \(p\in\Bun\) and let \(m:=\pi(p)\).
The covector \(\sigma_p\in T^*_p\Bun\otimes V\) vanishes on the fiber \(T_p(\Bun_p)\). 
There is a unique linear isomorphism
\[
T_m M\cong T_p\Bun/T_p(\Bun_p)\xrightarrow{u}V,
\]
defined by \(u\circ\pi'(p)=\sigma_p\); denote this \(u\) as \(u=\Phi(p)\), defining an \(H\)-equivariant map
\(
\Bun\xrightarrow{\Phi}\framebundle{M}
\).
The reader can check that \(\Phi^*\sigma=\sigma\).
Summing up: a soldering form determines an \(H\)-structure, and vice versa.
Hence a principal \(H\)-bundle \(\Bun\to M\) has a soldering form exactly when there is an isomorphism
\[
\amal{\Bun}{H}{\GL{V}}\cong\framebundle{M}
\]
of principal \(\GL{V}\)-bundles.
If \(M\) is homotopy equivalent to a finite \(CW\)-complex (for example, if \(M\) is a compact manifold, or the interior of a compact manifold with boundary and corners), the existence of such an isomorphism is expressible in the language of characteristic classes \cite{Steenrod:1999} p. 178.
In particular, every \((X,G)\)-geometry, or \((H,\LieG)\)-geometry, gives rise to an \(H/H_1\)-structure, where as above \(H_1\subseteq H\) is the subgroup acting trivially on \(V=\LieG/\LieH\).
\begin{example}
If \((X,G)=(\RP{n},\PGL{n+1})\), an \((X,G)\)-geometry is a projective connection, and \(H=G^{x_0}\) is the subgroup of matrices
\[
\begin{pmatrix}
a&b\\
0&d
\end{pmatrix}
\]
split into blocks, \(a\in\R^{\times}\), \(d\in\GL{n}\), \(b\in\R^{n*}\), defined up to rescaling the matrix by any nonzero constant.
The subgroup \(H_1\) consists of the matrices
\[
\begin{pmatrix}
1&b\\
0&I
\end{pmatrix}
\]
So \(H/H_1=\GL{V}\), where \(V=\R^n\).
Hence, on any manifold \(M\), any projective connection has associated \(H/H_1\)-structure just the frame bundle itself.
\end{example}
\begin{example}
A conformal structure in dimension \(3\) or more is an \((X,G)\)-geometry where \((X,G)=(S^n,\PO{n+1,1})\).
Check that \(H/H_1\) is the group of conformal linear transformations of the tangent space \(T_{x_0} X\) at the north pole, i.e. linear transformations preserve the standard metric on the sphere, up to a constant factor.
Hence the associated \(H/H_1\)-structure is a Riemannian metric defined up to scaling by a positive smooth function.
\end{example}
\begin{theorem}%
[Ehresmann \cite{Ehresmann:1951} p. 43]%
\define{theorem!Ehresmann soldering}
\define{Ehresmann's soldering theorem}
Suppose that \((\LieG,H)\) is a Lie infinitesimal model, for example that of a homogeneous space.
Let \(V:=\LieG/\LieH\).
Take a manifold \(M\) of the same dimension as \(V\).
Take a principal right \(H\)-bundle \(\Bun\to M\).
The following are equivalent:
\begin{itemize}
\item 
This bundle is the bundle of a \((\LieG,H)\)-geometry.
\item
This bundle has a soldering form.
\item
There is an isomorphism \(\amal{\Bun}{H}{\GL{V}}\cong\framebundle{M}\) of principal \(\GL{V}\)-bundles.
\item
There is an isomorphism \(\amal{\Bun}{H}{V}\cong TM\) of vector bundles.
\end{itemize}
\end{theorem}
\begin{proof} 
We first ask whether, at some point of \(\Bun\), we can construct a single covector \(\omega\in T_p^*\Bun\otimes\LieG\) so that \(A\hook\omega=A\) for \(A\in\LieH\) and \(\omega+\LieH=\sigma\).
Locally, we can write \(\Bun\) as a product \(M\times H\), and any point can be made to become \((m,1)\).
The soldering form \(\sigma\) at that point is a linear isomorphism \(T_m M\to V:=\LieG/\LieH\).
We want to construct a covector \(\omega\) so that \(\omega(0,A)=A\) and \(\omega(\dot{m},0)+\LieH=\sigma(\dot{m})\).
So we need \(\omega(\dot{m},A)=\omega(0,A)+\omega(\dot{m},0)=A+\xi(\dot{m})\) where \(\xi\in T^*_m M\otimes\LieG\) can be any linear map lifting \(\sigma\).
But clearly there is at least one such linear map, just by taking any linear splitting \(\LieG=\LieH\oplus V\).
Moreover, any two such agree up to a covector in \(T_m^*M \otimes \LieH\).

So the set of pairs \((p,\omega)\) with \(p\in\Bun\) and \(\omega\in T_p^*\Bun\otimes\LieG\) so that \(A\hook\omega=A\) and \(\omega+\LieH=\sigma\) is an affine bundle \(\vb{A}\to\Bun\) modelled on the adjoint vector bundle \(\vbh\to\Bun\).
Let \(H\) act on this set \(\vb{A}\) of covectors by
\[
\omega h:=\Ad_h \RT{h}^*\omega,
\]
acting freely and properly on \(\vb{A}\) since it lifts the free and proper action on \(\Bun\), so the quotient is a manifold.
It preserves the affine structure on each fiber, so the quotient \(\vb{A}/H\to M\) is a smooth bundle of affine spaces.
Every smooth bundle of affine spaces has a global smooth section, by partition of unity.
\end{proof}
Kobayashi \cite{Kobayashi:1956} Theorem 2, Barakat \cite{Barakat:2004} proved some cases of this theorem.
\begin{example}
There are many homogeneous models \((X,G)\) for which any two \((X,G)\)-geometries on the same manifold have isomorphic principal bundles \cite{McKay:2016}.
They are not classified.
\end{example}
\section{Higher order structures}
We can iterate the construction of an \(H/H_1\)-structure from a Cartan geometry.
We have constructed from any \((X,G)\)-geometry \(\G\to M\) a quotient map 
\[
pH_1\in\G/H_1\mapsto \omega_p+\LieH\in\framebundle{M}.
\]
This map has image some \(H/H_1\)-bundle in \(\framebundle{M}\), which we might prefer to denote \(M^{(1)}\).
We can then view \(\G\to M^{(1)}\) as a Cartan geometry modelled on \((X^{(1)},G)\) where \(X^{(1)}:=G/H_1\).
We can repeat this process.

A \emph{higher order structure}\define{higher order structure} is a sequence of structures, each a bundle over the last: for a sequence of Lie groups \(H\), \(H^{(1)}\), \(H^{(2)}\), \dots.
It has \emph{finite type}\define{finite type}\define{higher order structure!finite type} if the sequence terminates in some trivial group \(H^{(N)}=\set{1}\).
A \(\set{1}\)-structure is precisely a trivialization, i.e. a coframing or a framing, it is also a \((G,G)\)-Cartan geometry for any Lie group \(G\) of dimension equal to the manifold on which we have the \(\set{1}\)-structure.
\begin{theorem}
Every finite type higher order structure, on any manifold with finitely many components, has automorphism group a finite dimensional Lie group acting smoothly.
\end{theorem}
\begin{proof}
By finite type, the higher order structure induces a \(\set{1}\)-structure, hence a Cartan geometry, so the automorphism group is a closed subgroup of the automorphism group of that Cartan geometry.
\end{proof}
\section{Effecting automorphisms}
We continue the theory of the soldering form.
Take a Cartan geometry \(\G\to M\) with model \((X,G)\), \(X=G/H\).
Let \(M_1:=\G/H_1\), a lift \(\G\to M_1\) to a geometry modelled on \(X_1:=G/H_1\).
Repeat the process, building bundles 
\[
\begin{tikzcd}
\dots\arrow[r]\arrow[d]&M_3\arrow[r]\arrow[d]     &M_2\arrow[r]\arrow[d]    &M_1\arrow[r]\arrow[d]&M\\
\dots\arrow[r]         &\framebundle{M_2}\arrow[r]&\framebundle{M_1}\arrow[r]&\framebundle{M}\arrow[ur]
\end{tikzcd}
\]
\begin{theorem}%
[Pecastaing \cite{Pecastaing:2016}]%
\define{theorem!Pecastaing}
\define{Pecastaing's theorem}
If \((X,G)\) is a connected, effective homogeneous space, any local automorphism of \((X,G)\)-geometries is determined by its value and first \(k\) derivatives as a map of the base manifolds, where \(k\le\dim H\) and \(k\) is at most the number of times that dimensions drop in
\[
\begin{tikzcd}
0\arrow[r]&\dots\arrow[r]&\LieH_3\arrow[r]&\LieH_2\arrow[r]&\LieH_1\arrow[r]&\LieH.
\end{tikzcd}
\]
\end{theorem}
\begin{proof}
By lemma~\vref{lemma:kernel}, if \(X\) is connected then the sequence of groups
\[
\begin{tikzcd}
\dots\arrow[r]&H_3\arrow[r]&H_2\arrow[r]&H_1\arrow[r]&H
\end{tikzcd}
\]
converges to the kernel \(K\).
At each step, if \(H_{i+1}\subseteq H_i\) is an open subgroup, i.e. a union of components, i.e. \(\dim H_{i+1}=\dim H_i\), i.e. \(\LieH_{i+1}=\LieH_i\), then \(\LieG/\LieH_{i+1}=\LieG/\LieH_i\), so we have the same frame bundle \(\framebundle{M_{i+1}}=\framebundle{M_i}\).
Moreover, \(H_{i+2}\) is the subset of \(H_{i+1}\) acting trivially on \(\LieG/\LieH_{i+1}=\LieG/\LieH_i\), i.e. \(H_{i+2}=H_{i+1}\), so we hit the kernel at most one step after the dimensions of the \(\LieH_i\) cease to drop.
Every local diffeomorphism
\[
\text{open }\subseteq M\xrightarrow{\varphi}\text{open }\subseteq M' 
\]
lifts to a map of frame bundles by
\[
(m,u)\mapsto (\varphi(m),u\circ\varphi'(m)^{-1}),
\]
depending on one derivative.
So the first derivative of an automorphism determines the zero order derivative on the frame bundle, so on the immersed submanifold \(M_1\subseteq\framebundle{M}\).
Apply induction.
If \(\LieH_i=\LieH_{i+1}\), these are the same frame bundles, so require the same number of derivatives down on \(M\) to specify the automorphism on \(M_{i+1}\).
\end{proof}
\begin{example}
Pseudo-riemannian geometry has \((X,G)=(\R^{p,q},\SO{p,q}\ltimes\R^{p,q})\), \(H=\SO{p,q}\), \(H_1=1\), so any two local automorphisms which agree in value and first derivative at a point agree on any connected set containing that point.
\end{example}
\begin{example}
Projective connections have \((X,G)=(\RP{n},\PGL{n+1})\), with \(H\) consisting of matrices
\[
\begin{pmatrix}
a&b\\
0&d
\end{pmatrix}
\]
up to scaling, so elements of \(\LieG/\LieH\) look like
\[
\begin{pmatrix}
*&*\\
X&*
\end{pmatrix}
\]
on which elements of \(H\) act by \(X\mapsto a^{-1}dX\), so \(H_1\) consists of matrices
\[
\begin{pmatrix}
1&b\\
0&I
\end{pmatrix}
\]
so \(\LieG/\LieH_1\) consists of matrices 
\[
\begin{pmatrix}
x&*\\
Y&Z
\end{pmatrix}
\]
on which elements of \(H_1\) act by
\[
\begin{pmatrix}
x\\
Y\\
Z
\end{pmatrix}
\mapsto
\begin{pmatrix}
x+bY\\
Y\\
Z-Yb
\end{pmatrix}.
\]
So \(H_2=1\).
Hence any two local automorphisms of a projective connection which agree in value and two derivatives at a point agree on any connected set containing that point.
\end{example}

\chapter{Holonomy}
\section{Holonomy of a loop}
Take a principal bundle \(G\to\Bun\to M\) with connection \(\omega\).
Every absolutely continuous curve in \(M\) lifts to a horizontal curve, uniquely determined by the starting point where we lift the initial point: in any local trivialization, the horizontal curves lifting the given curve are the solutions to a Lie equation, and we can cover the curve in \(M\) by finitely many open sets on which a trivialization exists.
The lift of a loop in \(M\) might not be a loop in \(\Bun\), but returns to the same fiber, so to a point which differs from the original point by a unique element of \(G\): the \emph{holonomy}\define{holonomy} of the loop.
If we replace the starting point in \(\Bun\), the new starting point differs from the old one by a unique element of \(G\).
The new lift and the old agree all along by this element of \(G\), since the connection has \(G\)-invariant horizontal spaces.
So the holonomy of the lift conjugates when we change the initial point.
Reparameterizing the loop, preserving orientation, we reparameterize its horizontal lift.
Hence the holonomy is unchanged.
\begin{example}
Take a cone and puncture the tip:
\[
\begin{tikzpicture}
\newcommand{\radiusx}{.75}
\newcommand{\radiusy}{.25}
\newcommand{\height}{1.5}
\newcommand{\Rx}{(\radiusx)*sqrt(1-(\radiusy/\height)*(\radiusy/\height))}
\newcommand{\Ry}{\radiusy*(\radiusy/\height)}
\coordinate (a) at (-{\Rx},{\Ry});
\coordinate (b) at ({\Rx},{\Ry});
%\coordinate (a) at (-{\radiusx*sqrt(1-(\radiusy/\height)*(\radiusy/\height))},{\radiusy*(\radiusy/\height)});
%\coordinate (b) at ({\radiusx*sqrt(1-(\radiusy/\height)*(\radiusy/\height))},{\radiusy*(\radiusy/\height)});
\draw[gray!50,thick,fill=gray!30] (a)--(0,\height)--(b)--cycle;
\fill[gray!50] circle ({\radiusx} and \radiusy);
\begin{scope}
\clip ([xshift=-2mm]a) rectangle ($(b)+(1mm,-2*\radiusy)$);
\draw[thick,gray!50] circle ({\radiusx} and \radiusy);
\end{scope}
\begin{scope}
\clip ([xshift=-2mm]a) rectangle ($(b)+(1mm,2*\radiusy)$);
\end{scope}
\end{tikzpicture}
\]
Cut it straight from the tip to a point of the base,
\[
\begin{tikzpicture}
\newcommand{\radiusx}{.75}
\newcommand{\radiusy}{.25}
\newcommand{\height}{1.5}
\newcommand{\Rx}{\radiusx*sqrt(1-(\radiusy/\height)*(\radiusy/\height))}
\newcommand{\Ry}{\radiusy*(\radiusy/\height)}
\coordinate (a) at (-{\Rx},{\Ry});
\coordinate (b) at ({\Rx},{\Ry});
\draw[gray!50,thick,fill=gray!30] (a)--(0,\height)--(b)--cycle;
\fill[gray!50] circle ({\radiusx} and \radiusy);
\begin{scope}
\clip ([xshift=-2mm]a) rectangle ($(b)+(1mm,-2*\radiusy)$);
\draw[thick,gray!50] circle ({\radiusx} and \radiusy);
\end{scope}
\begin{scope}
\clip ([xshift=-2mm]a) rectangle ($(b)+(1mm,2*\radiusy)$);
\end{scope}
\draw[thick,black!50] (0,\height) -- ({\Rx*cos(40)},{-\radiusy*sin(40)});
\end{tikzpicture}
\]
Unravel it in the plane, to see 
\[
\begin{tikzpicture}
\draw[thick,gray!50,fill=gray!30] (0,0) --  (220:1.5) arc(220:120:1.5) -- cycle;
\end{tikzpicture}
\]
Imagine rolling the cone in the plane, leaving this trace behind.
The punctured cone is the quotient of the universal covering space of the punctured plane by the group of rotations by the cone angle.
The flat connection on the oriented orthonormal frame bundle of the plane pulls back to the universal covering space of the punctured plane, and quotients to the cone.
The holonomy of the loop around the tip is the rotation by the cone angle.
\end{example}
\begin{example}
Take the trivial bundle \(\Bun:=M\times G\) with the standard flat connection \(\omega=\omega_G\).
The horizontal curves are the curves \(p(t)=(m(t),g_0)\) for constant \(g_0\in G\), hence the holonomy of each loop is \(1\).
\end{example}
\begin{example}
Take a connected manifold \(M\), let \(\pi:=\fundamentalGroup{M}\), and take some group morphism
\[
\pi\xrightarrow{\rho}G
\]
to a Lie group \(G\).
Take the universal covering space \((\tilde{M},\tilde{m}_0)\xrightarrow{p}(M,m_0)\).
The standard flat connection \(\omega=\omega_G\) on \(\tilde{M}\times G\) is \(\pi\)-invariant, so descends to a flat connection on \(\Bun:=\amal{\tilde{M}}{\pi}{G}\).
Every continuous curve in \(M\) from a point \(m_0\) lifts to a curve in \(\tilde{M}\).
On \(\tilde{M}\times G\), the horizontal curves are again \(p(t)=(m(t),g_0)\).
But then on the quotient, \(p(1)=(m(1),g(1))\) is identified with
\[
(\gamma m(1),\rho(\gamma)g(1)),
\]
for every \(\gamma\in\pi\).
So \(m(t)\) drops to a loop in \(M\) just when \(m(1)=\gamma m(0)\) for some \(\gamma\in\pi\), necessarily unique, and the holonomy of this loop is \(\rho(\gamma)\in G\).
\end{example}
\section{Curvature as infinitesimal holonomy}
Take a principal bundle \(G\to\Bun\to M\) and a connection \(\omega\) on that bundle.
Take a contractible loop and contract it (smoothly, to avoid technical complications).
Let us pick some notation for the contraction.
Take a rectangle in the \(s,t\)-plane.
Map each horizontal line segment in our rectangle, of constant \(t\), into a loop in \(M\).
All the loops start and end at the same point:
\[
\begin{tikzpicture}
\begin{scope}[xshift=-3cm,yshift=-1cm]
\fill[gray!20] (0,0) rectangle (2,2);
%\draw[very thick,black] (0,0) -- (1,0);
\coordinate[label=right:\(s\)] (s) at (2,0);
\coordinate[label=above:\(t\)] (t) at (0,2);
\foreach \i in {0,...,20}
{
	\pgfmathsetmacro\clr{5*(\i)}
	\draw[thick,blue!\clr!brown] (0,{.1*\i}) -- (2,{.1*\i});
}
\draw[ultra thick,brown] (2,0) -- (0,0) -- (0,2);
\draw[ultra thick,brown] (2,0) -- (2,2);
\end{scope}
\coordinate (ar) at (-.75,0);
\coordinate (en) at (.75,0);
\path[-stealth,ultra thick] (ar) edge[bend left] (en);
\begin{scope}[xshift=2cm]
\fill[gray!20] (0,0) circle (1cm);
\foreach \i in {0,...,20}
{
	\pgfmathsetmacro\clr{5*(\i)}
	\draw[thick,blue!\clr!brown] ({1-.05*\i},0) circle ({.05*\i});
}
\end{scope}
\end{tikzpicture}
\]
The top edge maps to the original loop.
The bottom edge is the finished contraction: it goes to a single point, say \(m_0\in M\).
All of our loops start and stop at that same point \(m_0\).
So the left and right sides also go to that same point, since they are where each loop starts and stops.
So the left and right hand sides, and the bottom, of the rectangle all go to the same point \(m_0\in M\).

Pick a point \(p_0\in\Bun_{m_0}\).
Lift each horizontal curve in the rectangle to a horizontal curve in the bundle \(\Bun\) (but now the word ``horizontal'' means that \(\omega=0\) on it), to start at \(p_0\).
In other words, the left side of the rectangle maps to \(p_0\).
The bottom edge is still mapped to a point, since it is the lift of a point, so it goes to \(p_0\).
The bottom right corner maps to \(p_0\), by continuity.

Careful with the right side.
When we mapped to \(M\), the right side mapped to \(m_0\).
So when we lift, the right side maps to the fiber \(\Bun_{m_0}\).
(We will suppose a right principal bundle, but it is easy to change notation for a left principal bundle.)
Mapping
\[
g\in G\mapsto p_0g\in\Bun_{m_0}
\]
identifies \(G\) with \(\Bun_{m_0}\), pulling back the connection to the Maurer--Cartan form.
So the holonomy along the horizontal loops, thought of as a map taking the right side to \(G\), solves the Lie equation
\[
\omega_G=(\partial_t)\hook\omega_{\Bun}
\]
for a curve \(g\in G\) with initial condition \(g=1\).

We can solve that same Lie equation (without the same interpretation as holonomy) on all of the vertical line segments (of constant \(s\)), starting with \(g=1\) at the bottom:
\[
\begin{tikzpicture}
\begin{scope}[xshift=-3cm,yshift=-1cm]
\fill[gray!20] (0,0) rectangle (2,2);
%\draw[very thick,black] (0,0) -- (1,0);
\coordinate[label=right:\(s\)] (s) at (2,0);
\coordinate[label=above:\(t\)] (t) at (0,2);
\foreach \i in {0,...,20}
{
	\pgfmathsetmacro\clr{5*(\i)}
	\draw[thick,blue!\clr!brown] ({.1*\i},0) -- ({.1*\i},2);
}
\end{scope}
\coordinate (ar) at (-.75,0);
\coordinate (en) at (.75,0);
\path[-stealth,ultra thick] (ar) edge[bend left] (en);
\begin{scope}[xshift=2cm]
\fill[gray!20] (0,0) circle (1cm);
\foreach \i in {0,...,20}
{
	\pgfmathsetmacro\clr{5*(\i)}
	\draw[thick,blue!\clr!brown] (1,0) -- ({cos(360*\i/20)},{sin(360*\i/20)});
}
\end{scope}
\end{tikzpicture}
\]
We produce a map \(g(s,t)\in G\) with \(g=1\) along the left and bottom sides, equal to the holonomy along the right side.
Denoting partial derivatives with subscripts,
\[
g_t\hook\omega_G=\partial_t\hook\omega_{\Bun}.
\]
\begin{align*}
\partial_s(g_t\hook\omega_G)
&=
\partial_s(\partial_t\hook\omega),
\\
&=
\partial_t\hook\LieDer_{\partial_s}\omega,
\\
&=
\partial_t\hook(\partial_s\hook d\omega+d(\partial_s\hook\omega)),
\\
&=
\partial_t\hook(\partial_s\hook \Omega+0),
\\
&=
\Omega(\partial_s,\partial_t).
\end{align*}
Integrate in \(s\): along the right side, the holonomy \(g\) satisfies
\[
g_t\hook\omega_G=\int\Omega(\partial_s,\partial_t)ds.
\]
This equation gives a meaning to ``curvature is infinitesimal holonomy''.
\Danger{} here we are not calculating the holonomy \(g\) itself, but only its left logarithmic derivative \(g_t\hook\omega_G\).
Calculating \(g\) requires solving this equation: the Lie equation associated to this left logarithmic derivative.
\section{The holonomy group}
If we follow one loop, and then another, clearly we multiply their holonomies.
If we follow a loop backwards, we invert the holonomy.
Hence the holonomies which arise at a given point of the bundle form a group.
For notation: take a point \(p_0\in\Bun\) and its image point \(m_0\in M\).
The \emph{holonomy group}\define{holonomy!group} \(\holonomyGroup{G}[p_0]\subset G\) at \(p_0\) is the set of all holonomies of lifts of loops in \(M\) starting and ending at \(m_0\).
(The traditional notation for the holonomy group is \(\operatorname{Hol}\), but we will also need a holonomy reduction and a holonomy algebra, so it could become unclear which of these the notation \(\operatorname{Hol}\) refers to.)
As we have seen above, under \(G\)-action, the holonomy group conjugates.
\begin{lemma}\label{lemma:holonomy.horizontally.invariant}
Points connected by a horizontal curve have the same holonomy group.
To be precise, take a principal bundle \(\Bun\to M\) with a connection.
Suppose that \(p_0,p_1\in\Bun\) are connected by an absolutely continuous horizontal curve.
Then \(\holonomyGroup{G}[p_0]=\holonomyGroup{G}[p_1]\).
\end{lemma}
\begin{proof}
Let \(m_0,m_1\in M\) be the images of \(p_0,p_1\).
Choose a path from \(m_0\) to \(m_1\) that has a horizontal lift from \(p_0\) to \(p_1\).
Picture taking a loop at \(m_1\) and then drawing a ``lasso''\define{lasso} from \(m_0\):
\[
\begin{tikzpicture}
\coordinate (m0) at (0,0);
\coordinate (m1) at (4,2);
\begin{scope}[very thick,decoration={
    markings,
    mark=at position 0.5 with {\arrow{stealth}}}
    ] 
\draw[postaction={decorate},thick] (m0) to[curve through={(1,.4) .. (2,-.1) .. (3,1.3)}] (m1);
\draw[postaction={decorate},thick] (m1) to[curve through={(5,3) .. (6,2) .. (4,1.3)}] (m1);
\end{scope}
\fill[draw=white] (m0) circle (1.6pt);
\fill[draw=white] (m1) circle (1.6pt);
\node[below] at (m0) {\(m_0\)}; 
\node[above left] at (m1) {\(m_1\)}; 
\end{tikzpicture}
\]
following our chosen path from \(m_0\) to \(m_1\), then the loop at \(m_1\), then the chosen path backwards to \(m_0\).
This associates to every loop at \(m_1\) a loop at \(m_0\).

For the sake of notation, suppose that our bundle is a left \(G\)-bundle.
Take our loop at \(m_1\) and lift it up to a path from \(p_1\), say with holonomy \(g\); so ending at \(gp_1\).

Lift the lasso: walk along our horizontal path from \(p_0\) to \(p_1\), then that lifted path from \(p_1\) to \(gp_1\), and then the \(g\)-translate of our horizontal path, backwards, takes us back to \(gp_0\).
So the holonomy of the lasso is the holonomy of the loop: \(g\).
Hence the holonomy group of \(p_0\) contains \(g\), i.e. contains that of \(p_1\).
Reverse the roles of \(p_0,p_1\) to see that the holonomy groups are equal.
\end{proof}
\begin{corollary}
Take a principal bundle with connection, on a connected base manifold. 
Above any two points of the base manifold, there are points of the bundle connected by a horizontal absolutely continuous path, and hence having the same holonomy group.
\end{corollary}
\begin{corollary}
Take a principal bundle with connection, on a connected base manifold. 
The holonomy groups at any two points of the bundle are conjugate subgroups.
\end{corollary}
\begin{theorem}\label{thm:conjugate.hol}
Take a principal bundle with connection \(G\to\Bun\to M\) on a connected manifold \(M\).
If some group belongs to the holonomy group at a chosen point of \(\Bun\) then, at every point of \(\Bun\), some conjugate belongs to the holonomy group.
If that group is normal in the structure group then it belongs to the holonomy group at every point.
If the holonomy group is normal in \(G\) at some point, it is constant throughout the bundle.
\end{theorem}
We write \(\holonomyGroup{G}\) to mean \(\holonomyGroup{G}[p_0]\) for an unspecified point \(p_0\in\Bun\).
\begin{lemma}
A principal bundle with connection, on a connected manifold, is isomorphic to the trivial bundle with the standard flat connection just when its holonomy is trivial at some, hence any, point.
\end{lemma}
\begin{proof}
The trivial principal bundle with the standard flat connection has trivial holonomy: \(\holonomyGroup{G}[p_0]=\set{1}\), as the horizontal paths are just
\[
p(t)=(m(t),g_0).
\]
On the other hand, if some connection has \(\holonomyGroup{G}[p_0]=\set{1}\), we can parallel transport along any paths in \(M\), and the resulting parallel transport identifies all fibers.
The identification is independent of the path, by triviality of the holonomy.
Hence the principal bundle is trivial.
The trivialization makes the horizontal spaces precisely the spaces where we don't move in the fiber direction.
So our trivialization gives trivial parallel transport and the standard flat connection.
\end{proof}
\begin{theorem}\label{theorem:induced.holonomy}
Take a Lie group morphism \(G\xrightarrow{\rho}G'\) and a principal \(G\)-bundle \(G\to\Bun\to M\) with associated \(G'\)-bundle \(\Bun'\), with bundle morphism \(\Bun\to\Bun'\).
Take a point \(p\in\Bun\) and the associated point \(p'\in\Bun'\).
For each connection on \(\Bun\), the induced connection on \(\Bun'\) has holonomy group
\[
\dot{G}'_{p'}=\rho \dot{G}_p.
\]
\end{theorem}
\begin{proof}
Recall from theorem~\vref{theorem:induced.connection} that the induced connection has horizontal space the image of the horizontal space.
\end{proof}
\section{Restricted holonomy}
The \emph{restricted holonomy group}\define{restricted holonomy group}\define{holonomy!group!restricted} is the subgroup \(\restrictedHolonomyGroup{G}[p_0]\subseteq\holonomyGroup{G}[p_0]\) generated by the holonomies of contractible absolutely continuous loops in \(M\)  starting and ending at the image point \(m_0\in M\) of \(p_0\in\Bun\).
Clearly the restricted holonomy group is a normal subgroup of the holonomy group, by conjugating any contraction of a loop.
\begin{lemma}
Any two absolutely continuous paths in a manifold are homotopic with fixed endpoints if and only if they are absolutely continuously homotopic with fixed endpoints.
\end{lemma}
\begin{proof}
We give an outline but leave the details to the reader; for a complete proof, see \cite{McKay:diff.geom} p. 352 lemma 38.11.
Also see \cite{Kobayashi/Nomizu:1996} p. 284 Appendix 7: the Factorization Lemma and \cite{Lee:2013} p. 142 theorem 6.29 for essentially the same proof of essentially the same result.
Picture a continuous map of a box into a manifold, and suppose that we have already arranged the map to be absolutely continuous, or even smooth, in certain variables.
For each face, we have also arranged the map to be absolutely continuous, or even smooth, in some other variables, on that face.
Suppose that our manifold is a convex domain in Euclidean space.
We leave the reader to argue, by local convolutions, that we can smooth out the interior, by a small smoothing homotopy, fixing every point of whichever of the faces we like, but smoothing the faces we choose not to fix, preserving absolute continuity or smoothness in the variables we wish.
In fact, we can do much better with convolutions \cite{Fefferman:2009,Fefferman.Israel:2020,McKay:diff.geom,Whitney:1934a,Whitney:1934b,Whitney:1934c}.

Now consider a map of a box into a more general manifold \(M\); we can assume our manifold \(M\) is connected.
Cover \(M\) in open sets diffeomorphic to convex domains.
By compactness of the box, we can divide it into a grid of smaller boxes, meeting along their faces, on each of which our map takes that box into one of those open sets.
We can therefore smooth on each box, by small homotopy, one box at a time, fixing the map on faces as needed.
The same proof easily extends to maps between manifolds with boundary and corners.
\end{proof}
\begin{corollary}
The restricted holonomy group of any connection is a connected subgroup of the structure group.
\end{corollary}
\begin{corollary}
The restricted holonomy group of any connection on any connected manifold is the holonomy group of the pullback to some covering space, in particular to the universal covering space.
\end{corollary}
\section{Lie group structure}
Recall that an \emph{immersed submanifold}\define{immersion} \(S\subseteq M\) is a subset which is somehow equipped with a topology and smooth structure for which the inclusion mapping is an immersion; there may be more than one such topology, for example if \(S\) is a figure \(8\)  in the plane.
An immersed submanifold \(S\subseteq M\) of a manifold \(M\) is \emph{initial}\define{initial!submanifold} if every smooth map of manifolds \(N\to M\) with image in \(S\) is a smooth map \(N\to S\).
The figure \(8\) is not initial for any of its topologies as an immersed submanifold.
\begin{lemma}\label{lemma:initial.unique}
A subset which is an immersed submanifold for some topology is initial for at most one topology and for at most one smooth structure.
\end{lemma}
\begin{proof}
The immersed submanifold \(S\), in all of its topologies, maps to \(M\) smoothly, hence to itself in any initial topology smoothly.
So in any two initial topologies, it maps to itself smoothly from each to the other, a diffeomorphism.
\end{proof}
An \emph{initial subgroup}\define{initial!subgroup} \(G\subseteq H\) of a Lie group \(H\) is an initial submanifold which is also a subgroup.
\begin{theorem}%
[Initial subgroup theorem]%
\label{theorem:initial.subgroup}%
\define{theorem!initial subgroup}
\define{initial!subgroup theorem}
Every subgroup \(G\subseteq H\) of a Lie group with a countable set of connected components is an initial subgroup for a unique smooth structure.
Its connected components in that topology are its path connected components as a subset of \(H\).
The Lie algebra of \(G\) is the subset of the Lie algebra of \(H\) whose flows lie in \(G\).
\end{theorem}
\begin{proof}
Hilgert and Neeb \cite{Hilgert.Neeb:2012} p. 356 prove that there is such a smooth structure and that its connected components and Lie algebra are as stated.
By lemma~\vref{lemma:initial.unique}, it is unique.
\end{proof}
\begin{theorem}%
[Yamabe 
\cite{Hilgert.Neeb:2012} p. 347 Theorem 9.6.1]%
\label{theorem:Yamabe}%
\define{theorem!Yamabe}
\define{Yamabe's theorem}
Any connected subgroup of a Lie group is a connected initial Lie subgroup for a unique smooth structure.
\end{theorem}
Hence the restricted holonomy group is an injectively immersed Lie subgroup of the structure group.
\begin{lemma}\label{lemma:G.zero}
Take a group \(G\) and a subgroup \(G_0\subseteq G\).
Suppose that \(G_0\) has a smooth structure in which it is a connected Lie group.
Suppose that \(G\) is the union of countably many left translates of \(G_0\).
Then there is a Lie group structure on \(G\) making \(G_0\) its identity component if and only if 
\(G_0\subseteq G\) is normal and  each element of \(G\) acts by adjoint action on \(G_0\) continuously.
That Lie group structure is then unique.
\end{lemma}
\begin{proof}
If \(G\) is a Lie group and \(G_0\) is the identity component of \(G\), clearly \(G_0\) is normal and each element of \(G\) acts by adjoint action on \(G_0\) continuously.
Suppose the converse.
Every continuous automorphism of a Lie group is smooth, because its graph is a closed subgroup of the product.
So adjoint action of any element of \(G\) on \(G_0\) is a smooth Lie group automorphism.

Write \(G\) as the disjoint union of left translates of \(G_0\).
Make each of these left translates of \(G_0\) diffeomorphic to \(G_0\) by that left translation.
Suppose we change the choice how we left translate: say you left translate by \(g_1\) and I left translate by \(g_2\).
The left translates \(g_1G_0\) and \(g_2G_0\) are disjoint or equal.
If equal, the smooth structures are related by
\[
g\in G_0\mapsto g_2^{-1}g_1g \in G_0, 
\]
which is a diffeomorphism of \(G_0\), since \(g_2^{-1}g_1\in G_0\).
Thus \(G\) becomes a manifold independent of how we choose to left translate.
But also this is the unique smooth structure on \(G\) in which \(G_0\) is the identity component.

We need the multiplication to be smooth.
Pick points \(x_0,y_0\in G\).
We need to prove that multiplication \(x,y\mapsto xy\) is smooth for \(x\) near \(x_0\) and \(y\) near \(y_0\).
Write each \(x\) near \(x_0\) as \(x=x_1x_0\), \(y=y_1y_0\) for \(x_0,y_0\in G_0\) near \(1\).
Then
\[
xy=x_1y_1x_0'y_0
\]
where \(x_0'=y_1^{-1}x_0y_1\).
So this is smooth as a function of \(x_0',y_0\), and hence as a function of \(x_0,y_0\).
\end{proof}
\begin{corollary}
Take a subgroup \(G\subseteq H\) of a Lie group \(H\), and a path connected normal subgroup \(G_0\subseteq G\).
Suppose that \(G\) is the union of countably many left translates of \(G_0\).

Then \(G_0\) is the identity path component of \(G\) as a subset of \(H\).
In the unique smooth structure on \(G\) so that \(G\) is an initial subgroup of \(H\), \(G_0\) is its identity component, so this smooth structure is the same one as in lemma~\vref{lemma:G.zero}.
\end{corollary}
\begin{proof}
Since \(G_0\) is path connected, it lies in a single path connected component \(G'_0\subseteq G\).
Since \(G\) lies in countably many translates of \(G_0\), \(G\) also lies in countably many translates of \(G'_0\).
By lemma~\vref{lemma:G.zero}, there is a unique smooth structure on \(G\) for which \(G'_0\) is the identity component.
By theorem~\vref{theorem:initial.subgroup}, this is the unique smooth structure on \(G\) for which \(G\subseteq H\) is initial.

By Yamabe's theorem, \(G_0\subseteq H\) is a connected initial Lie subgroup.
Countably many translates of \(G_0\) cover \(G\).
By Sards' theorem \cite{Sard:1965}, \(\dim G_0\ge\dim G\).
But \(G\) contains \(G_0\), so again by Sard's theorem, \(\dim G_0\le \dim G\), hence \(\dim G_0=\dim G=\dim G_0'\).
Since \(G_0\) and \(G_0'\) are initial, the inclusion map \(G_0\to H\) is a smooth injective Lie group morphism \(G_0\to G_0'\) of equal dimensional Lie groups, so \(G_0\) is an open subgroup of \(G_0'\), so a connected component, so \(G_0=G_0'\).
\end{proof}
\begin{theorem}
Take a connection on a principal bundle over a connected manifold.
The holonomy group has the restricted holonomy as its identity path component, in the topology induced from the structure group.
There is a unique smooth structure on the holonomy group in which the restricted holonomy is its identity component.
This is the unique smooth structure in which the holonomy  group is a initial Lie subgroup of the structure group.
\end{theorem}
\begin{proof}
The restricted holonomy group is a normal subgroup.
Pick generators \(\set{\gamma_{\alpha}}\) of the fundamental group \(\pi_1(M)\).
Note that the fundamental group \(\pi_1(M)\) of any manifold is countably generated.
(To see this, cover the manifold in a countable collection of connected and simply connected open sets, each labelled by some symbol.
Any path is determined up to homotopy by a string of finitely many symbols representing open sets it enters and which cover its image; note that the string need not be unique.)
Therefore countably many translates of the restricted holonomy group cover the holonomy group.
By lemma~\vref{lemma:G.zero}, there is a unique Lie group structure on the holonomy group making it a Lie group with the restricted holonomy group as identity component.
Moreover, the restricted holonomy is the path component of the identity.
\end{proof}
\begin{example}
This theorem has a surprising consequence. 
Take some loop \(\gamma\) on \(M\).
Let \(g\) be its holonomy.
Then \(g\) is connected by a continuous path inside \(\holonomyGroup{G}\) to the identity if and only if there is some other loop \(\gamma'\), perhaps very far from \(\gamma\), which is contractible, and has the same holonomy \(g\).
\end{example}
\section{The orbit theorem and holonomy}
The \emph{holonomy reduction}\define{holonomy!reduction} at a point \(p_0\) is the set \(\holonomyReduction{\Bun}[p_0]\subset\Bun\) of points of \(\Bun\) connected to a point \(p_0\) by an absolutely continuous horizontal curve.

Given any vector field \(v\) on \(M\), its \emph{horizontal lift}\define{horizontal lift!of vector field} \(\hat{v}\) is the horizontal vector field on \(\Bun\) projecting to it.
Let \(\VF\) be the set of all lifts of horizontal vector fields.
By the Orbit Theorem, every orbit of \(\VF\) is a connected initial submanifold containing all lifted vector fields among its tangent vector fields.
Every horizontal vector arises by lifting some vector field.
So the tangent spaces to the orbit contain the horizontal spaces of the connection.
Hence every absolutely continuous horizontal curve lies inside a single orbit.
So each holonomy reduction lies in a single orbit.

On the other hand, the orbit is the union of compositions of flows of lifted vector fields.
Hence every point of the orbit is on an absolutely continuous horizontal curve.
So the holonomy reduction through a point is precisely the orbit of that point.
Therefore each holonomy reduction \(\holonomyReduction{\Bun}[p_0]\), say for a point \(p_0\in\Bun\) projecting to a point \(m_0\in M\), is the union of end points of the lifts of all absolutely continuous curves in \(M\) starting at \(m_0\).
But we can also replace those paths by piecewise smooth flows of lifted vector fields, which are also horizontal paths, and still (by definition of the orbit) reach the same points.
Hence, when we define holonomy groups and holonomy reductions, there is no loss of generality in employing piecewise smooth paths in place of absolutely continuous paths.
The holonomy reduction \(\holonomyReduction{\Bun}[p_0]\) intersects the fiber \(\Bun_{m_0}\) precisely in the set of translates of \(p_0\) by elements of the holonomy group.

\section{Holonomy of analytic connections}
\begin{theorem}
Take a real analytic principal bundle \(G\to\Bun\to M\) with a real analytic connection.
Take a locally slack set of real analytic vector fields on \(M\), with a single orbit in each path component of \(M\). 
Then the holonomy reduction through any point of \(\Bun\) is precisely the orbit in \(\Bun\) of the horizontal lifts of all of the vectors fields in that set.
\end{theorem}
\begin{proof}
We may assume that \(M\) is path connected.
By theorem~\vref{theorem:orbit.maps}, each orbit in \(\Bun\) projects by a surjective submersion to an orbit in \(M\).
The Lie saturate pushes down to the Lie saturate.
By the Hermann--Nagano theorem (theorem~\vref{thm:Hermann.Nagano}), each holonomy reduction \(\holonomyReduction{\Bun}[p_0]\) has tangent space at \(p_0\) precisely the values \(v(p_0)\) of the vector fields \(v\) in the Lie saturate of the real analytic horizontal lifts.
Again by the Hermann--Nagano theorem, the Lie saturate in \(M\) spans the orbit, so spans every tangent space of \(M\).
Hence its horizontal lifts span every horizontal space.
The brackets of the horizontal lifts agree up to vertical vector fields.
So the Lie saturate of the horizontal lifts upstairs also spans the horizontal spaces.
By the orbit theorem, the orbits upstairs are the holonomy reductions.
\end{proof}
\begin{corollary}
Take a real analytic principal bundle with a real analytic connection.
By the Orbit Theorem, the holonomy reduction through any point is precisely the orbit of any of the following sets of vector fields on \(\Bun\):
\begin{itemize}
\item
the set of all horizontal lifts of real analytic vector fields defined on open subsets of \(M\),
\item
the set of all real analytic horizontal vector fields, defined on open sets in \(\Bun\),
\item
the set of all smooth horizontal vector fields,
\item
the set of all horizontal lifts of smooth vector fields defined on open subsets of \(M\)
\item
the set of all horizontal lifts of some locally slack set of real analytic vector fields, each defined on an open subset of \(M\), with a single orbit in each path component of \(M\).
\end{itemize}
\end{corollary}
The set of real analytic vector fields on \(M\) is locally slack (see the definition~\vpageref{page:define.slack}), clearly.
Hence the set of real analytic horizontal lifts is locally slack.
\section{The holonomy reduction}
Every horizontal vector is the value of the lift of a complete vector field \(v\) on \(M\).
Complete vector fields lift to complete vector fields, again by global solvability of Lie equations.
So every orbit is the orbit of these complete vector fields.
By theorem~\vref{theorem:orbit.maps}, each orbit \(\holonomyReduction{\Bun}[p_0]\) maps to \(M\) by a smooth fiber bundle map.
The fibers are precisely the orbits of the holonomy group.
By horizontal invariance of holonomy groups (lemma~\vref{lemma:holonomy.horizontally.invariant}), the holonomy group \(\holonomyGroup{G}[p]\) is the same for every point \(p\) in the holonomy reduction.
The holonomy group is thereby immersed in the structure group, with the restricted holonomy as identity path component.
\begin{theorem}[The reduction theorem I]%
\define{theorem!reduction I}
\define{reduction theorem I}
Take a connection on a principal bundle \(G\to\Bun\to M\) over a connected manifold \(M\).
The holonomy reduction \(\holonomyReduction{\Bun}[p_0]\subset\Bun\) is a principal right \(\holonomyGroup{G}\)-bundle over \(M\).
\end{theorem}
This is a slight improvement on the usual reduction theorem \cite{Joyce:2007} p. 30.
\begin{proof}
Being a subgroup of the structure group, \(\holonomyGroup{G}\) acts freely on \(\Bun\), hence on \(\holonomyReduction{\Bun}:=\holonomyReduction{\Bun}[p_0]\subseteq\Bun\).
The \(\holonomyGroup{G}\)-action on \(\holonomyReduction{\Bun}\), as a map
\[
\holonomyGroup{G}\times\holonomyReduction{\Bun}\to\holonomyReduction{\Bun}
\]
factors through the \(G\)-action on \(\Bun\), so is smooth.
We will assume a right principal bundle, but the obvious change of notation works for a left principal bundle.
Since \(\holonomyReduction{\Bun}\to M\) is a smooth fiber bundle, we can assume it is trivial, say \(\holonomyReduction{\Bun}=M\times F\) for some fiber \(F\).
The group \(\holonomyGroup{G}\) acts smoothly and simply transitively on \(F\).
Since \(\holonomyGroup{G}\) has countably many components (see \vpageref{page:G.mod.H.trouble}), \(F\cong\holonomyGroup{G}/\set{1}\), so \(\holonomyReduction{\Bun}\) is a principal \(\holonomyGroup{G}\)-bundle.
\end{proof}
Note that holonomy reduction depends on choice of a point \(p_0\in\Bun\), but any two are related by the \(G\)-action on \(\Bun\), conjugating their structure groups, i.e. the holonomy groups.
\begin{theorem}%
[The reduction theorem II]%
\define{theorem!reduction II}
\define{reduction theorem II}
Take a connection on a principal bundle over a connected manifold \(M\).
The connection pulls back to any of its holonomy reductions to be a connection.
\end{theorem}
\begin{proof}
The pullback \(\omega\) of the connection is a \(1\)-form valued in \(\LieG\), perhaps not in \(\LieH\), transforming in the adjoint representation, and agreeing with the Maurer--Cartan form on the vertical vectors.
It vanishes on the horizontal vector fields that we lift from \(M\).
By definition of \(\holonomyReduction{\Bun}\) as an orbit of such vector fields, these vector fields are tangent to \(\holonomyReduction{\Bun}\).
All put together, these horizontal vector fields span, at each point of \(\Bun\), a linear subspace of dimension equal to that of \(M\).
But they are tangent to \(\holonomyReduction{\Bun}\), so that linear subspace is tangent to \(\holonomyReduction{\Bun}\) and \(\omega=0\) on it.
Together with the vertical directions, these horizontal directions split the tangent space to \(\holonomyReduction{\Bun}\), as it is a principal \(\holonomyGroup{G}\)-bundle.
But \(\omega=0\) on the horizontal, equals the \(\holonomyGroup{G}\)-Maurer--Cartan form on the vertical, and together these split the tangent space, so \(\omega\) is valued in the Lie algebra of \(\holonomyGroup{\LieG}\) on all tangent vectors to \(\holonomyReduction{\Bun}\).
\end{proof}
\section{The holonomy algebra}
The \emph{holonomy algebra}\define{holonomy!algebra} \(\holonomyAlgebra{G}\) is the Lie algebra of the holonomy group \(\holonomyGroup{G}\).
\begin{corollary}\label{corollary:normal.sub.hol}
Take a principal bundle \(G\to\Bun\to M\) with connection on a connected manifold \(M\).
Suppose that at some point of \(\Bun\), some \(G\)-invariant subspace in \(\LieG\) lies in the holonomy algebra. 
Then the subgroup of \(G\) generated by that subspace lies in the restricted holonomy group of every point of \(\Bun\).
\end{corollary}
\begin{proof}
A \(G\)-invariant subspace of \(\LieG\) is precisely a \(G\)-invariant ideal of \(\LieG\).
The holonomy group at the given point \(p_0\) contains the connected subgroup generated by that ideal.
This subgroup is normal: each point is a product of exponentials of elements of the ideal, and these are carried to one another by adjoint \(G\)-action.
By horizontal invariance, the subgroup also lies in the holonomy algebra of any point \(p\) connected to \(p_0\) by a horizontal path.
The points \(p\) connected to \(p_0\) by a horizontal path project onto \(M\), because the horizontal space at each point \(p\) projects by linear isomorphism to a tangent space of \(M\) and \(M\) is connected.
By horizontal invariance, that subgroup belongs to the holonomy group of some point in each fiber of \(\Bun\to M\).
Conjugate by the action of the structure group: since the subgroup is normal, it belongs to the restricted holonomy group of every point.
\end{proof}
\begin{corollary}
At each point of any principal bundle with connection, the curvature of the connection is a \(2\)-form valued in the holonomy algebra of that point.
\end{corollary}
\begin{proof}
The connection pulls back to a connection on the holonomy reduction, so the curvature pulls back to the curvature.
Every tangent space to the holonomy reduction contains the horizontal spaces of the connection, and the curvature depends only on horizontal components.

We can also see this directly from the explicit formula
\[
\Omega(\partial_s,\partial_t)=\partial_s(g_t\hook\omega_G)
\]
which we computed above, by constructing a suitable contracting loop.
\end{proof}
\begin{corollary}
For any Lie group \(G\) and initial Lie subgroup \(H\subseteq G\), a principal \(G\)-bundle \(\Bun_G\to M\) with connection is isomorphic to the \(G\)-bundle
\[
\Bun_G\cong\amal{\Bun_H}{H}{G}
\]
associated to some principal \(H\)-bundle \(\Bun_H\to M\), and has the induced connection, if and only if it has holonomy in \(H\).
\end{corollary}
\begin{proof}
If \(\Bun_G\cong\amal{\Bun_H}{H}{G}\), just take any connection \(\omega\) on \(\Bun_H\), and then the differential form
\[
\omega_G+\Ad_g^{-1}\omega
\]
descends (as we have seen) to a connection on \(\Bun_G\), with \(\Bun_H\subseteq\Bun_G\) containing the horizontal directions at each of its points.
Hence the holonomy orbit lies in \(\Bun_H\) and the holonomy lies in \(H\).
On the other hand, if \(\Bun_H\) is the holonomy reduction, it is a principal \(H\)-subbundle and the connection restricts to it to be a connection.
One easily identifies \(\Bun_G\cong\amal{\Bun_H}{H}{G}\).
\end{proof}
\begin{corollary}
For any Lie group \(G\) and initial Lie subgroup \(H\subseteq G\), a principal \(G\)-bundle \(\Bun_G\to M\) is isomorphic to the \(G\)-bundle
\[
\Bun_G\cong\amal{\Bun_H}{H}{G}
\]
associated to some principal \(H\)-bundle \(\Bun_H\to M\) if and only if there is a connection on \(\Bun_G\) with holonomy in \(H\).
\end{corollary}
\section{The Ambrose--Singer theorem}
\begin{theorem}%
[Ambrose--Singer \cite{Ambrose.Singer:1953}]%
\label{thm:Ambrose.Singer}%
\define{theorem!Ambrose--Singer}%
\define{Ambrose--Singer theorem}%
Take a principal bundle with a connection.
The holonomy algebra at a point is precisely the span of all values \(\Omega(u,v)\) of the curvature form on tangent vectors \(u,v\) tangent to the holonomy reduction at any point of the holonomy reduction.
It is also precisely the span of all values \(\Omega'(u,v)\) of all parallel translates of the curvature form on tangent vectors \(u,v\) tangent to the holonomy reduction at a single point of the holonomy reduction.
\end{theorem}
Cartan was aware of this \cite{Cartan:1926} top of page 5, no doubt with the same proof.
\begin{proof}
We can replace the bundle by its holonomy reduction, so assume that the holonomy group is \(G\) and the bundle is \(G\to\Bun\to M\), a right principal \(G\)-bundle, and that \(\Bun\) is its own holonomy reduction through each of its points.
Let \(V\subseteq\LieG\) be the span of all values \(\Omega(u,v)\) of the curvature on tangent vectors \(u,v\) tangent to \(\Bun\).
Denote the connection by \(\omega\).
Consider the \(1\)-form \(\omega+V\in\nForms{1}{\Bun}\otimes(\LieG/V)\).
From the equation
\[
\Omega=d\omega+\frac{1}{2}\lb{\omega}{\omega},
\]
we see that, modulo \(\omega+V\),
\[
d\omega+V=\Omega+V=0
\]
since \(\Omega\) is valued in \(V\).
By the Frobenius theorem, \(\Bun\) is foliated by submanifolds whose tangent spaces consist in the tangent vectors on which \(\omega\in V\).
In particular, the horizontal directions \(\omega=0\) lie tangent to the leaves.
Hence each leaf contains the orbit of the horizontal vector fields, i.e. the holonomy reduction.
So there is a unique leaf: all of \(\Bun\).
But on \(\Bun\), \(\omega\) is onto \(\LieG\), since \(\omega\) equals the Maurer--Cartan form on fibers.
Hence \(V=\LieG\).

Pick any complete vector field \(v\) on \(M\), and lift it to a horizontal vector field \(\hat{v}\) on \(\Bun\).
Since the horizontal spaces of any connection are invariant under the structure group action, \(\hat{v}\) is invariant under the structure group action, so the flow of \(\hat{v}\) is by bundle automorphisms of \(\Bun\). 
Move the connection \(\omega\) by the flow of that lifted vector field to produce another connection.
The holonomy group is unchanged, as the flow preserves all holonomy orbits, i.e. all holonomy reductions.
We find therefore that when we move the curvature by the flow, it is now the curvature of some other connection with the same holonomy.
So its values still have the same span.
We can repeat this, composing such flows.
So the span of curvature values at all points of \(\Bun\) is the same as the span of curvature values of all of these moved curvature tensors at any one point of \(\Bun\).

We can approximate any piecewise smooth curve in \(M\) by flows of several vector fields, flowing along one vector field for some time, then another, and so on.
The resulting parallel transported curvature is therefore approximated by the curvature from the flows, and so lies in the closure, hence in the span, of those curvature values. 
\end{proof}
\section{Reductions of structure group}
\begin{lemma}
Take a principal bundle \(G\to\Bun\to M\), on a connected manifold \(M\) of dimension \(2\) or more.
Take a connection on that bundle, flat in some nonempty open set, with restricted holonomy group \(H\) at some point \(p_0\).
Suppose that \(H\) is not the identity component of \(G\).
Then there is another connection, flat in some nonempty open set contained in the previous open set, agreeing with the previous connection outside a compact set, with restricted holonomy group at \(p_0\) (1) containing \(H\) and (2) of larger dimension than \(H\).
\end{lemma}
\begin{proof}
Take some element \(A\in\LieG-\LieH\).
Let \(m_0\in M\) be the image of \(p_0\in\Bun\).
Pick a point \(m_1\ne m_0\) and a point \(p_1\) in the fiber of \(m_0\) in the holonomy reduction of \(p_0\).
Near \(m_1\), use the flat connection to trivialize the bundle, so that the connection becomes the standard flat connection.
Denote points \(p=(m,h)\) in the holonomy reduction, and \(p=(m,g)\) in the \(G\)-bundle, over some neighborhood of \(m_1\in M\).
We can suppose that \(p_1=(m_1,1)\). 

Take a horizontal curve from \(p_0\) to \(p_1\), say \(p(t)\), \(0\le t\le 1\).
For \(t\) near \(1\), this curve is \(p(t)=(m(t),1)\).
We can suppose that this curve \(t\mapsto m(t)\) is immersed near \(t=1\), and hence in some coordinates is, for \(t\) near \(1\):
\[
t\mapsto x=(x^1,x^2,\dots,x^n)=(1-t,0,\dots,0).
\]
Let
\[
\Gamma(x)=A\,x^1 f(x) dx^2,
\]
where \(f(x)\) is some bump function equal to \(1\) near the origin and equal to zero not far away.
Our curve \(p(t)\) is still horizontal for the old connection \(\omega\), but also remains horizontal for the new connection
\[
\omega'=\omega+\Ad_g^{-1}\Gamma
\]
on the \(G\)-bundle.
The curvature of that connection takes the value \(A\) on points of the \(\omega\)-holonomy reduction of \(p_0\), i.e. at \(p_1\).
By the Ambrose--Singer theorem, the \(\omega'\)-restricted holonomy algebra contains \(A\).

We claim the \(\omega'\)-restricted holonomy group also contains that of \(\omega\).
Picture a lasso\SubIndex{lasso} leaving \(m_1\), heading out backwards along \(m(t)\), away from the open set where \(\omega\) is flat, to pick up holonomy, and then going back to \(m_1\) along the same path by which it came.
\[
\begin{tikzpicture}[xscale=-1]
\coordinate (m0) at (0,0);
\coordinate (m1) at (4,2);
\begin{scope}[very thick,decoration={
    markings,
    mark=at position 0.5 with {\arrow{stealth}}}
    ] 
\draw[postaction={decorate},thick] (m0) to[curve through={(1,.4) .. (2,-.1) .. (3,1.3)}] (m1);
\draw[postaction={decorate},thick] (m1) to[curve through={(5,3) .. (6,2) .. (4,1.3)}] (m1);
\end{scope}
\fill[draw=white] (m0) circle (1.6pt);
\fill[draw=white] (m1) circle (1.6pt);
\node[below] at (m0) {\(m_1\)}; 
\node[left] at (m1) {\(m_0\)}; 
\end{tikzpicture}
\]

All restricted \(\omega\)-holonomy at \(p_1\) arises this way.
Since \(p(t)\) is horizontal for \(\omega\) and for \(\omega'\), the holonomy along the lasso agrees for \(\omega\) and for \(\omega'\).
The new holonomy algebra contains the contributions of these lassos, which are the same dimension for either connection, so contains all restricted holonomy of \(\omega\).
\end{proof}
\begin{lemma}\label{lemma:make.holonomy.big}
Take a principal bundle \(G\to\Bun\to M\), on a connected manifold \(M\) of dimension \(2\) or more.
That bundle has a connection with restricted holonomy precisely the identity component of \(G\).
For any connection \(\omega\) and any nonempty open subset \(U\subseteq M\), there is a connection \(\omega'\) with restricted holonomy precisely the identity component of \(G\) which agrees with \(\omega\) outside of a compact subset of \(U\).
\end{lemma}
\begin{proof}
It suffices to prove the result for \(M\) connected.

Suppose that we can prove the result for trivial bundles.
Then we can pick \(U\) in which the bundle is trivial, and the result follows.
So suppose that the bundle is trivial.
Take the standard flat connection \(\omega\).

If \(\LieG=0\) this suffices, so suppose not.

Take some point \(m_1\in M\) and a curvature bundle element \(\Omega_1\) taking on some nonzero value \(\Omega_1(u,v)\) at a point \(p_1\) in the fiber over \(m_1\).
Take a relatively compact open set \(U_1\) around \(m_1\), whose closure is not \(M\).
Pick a connection over \(U_1\) whose curvature at \(p_1\) is \(\Omega_1\).
Use a partition of unity to glue this connection to \(\omega\), so that it agrees with \(\omega\) outside of the preimage in \(\Bun\) of \(U_1\).
Call the resulting connection \(\omega\).
The holonomy algebra now has dimension at least \(1\).

After \(k\) steps, suppose that the curvature algebra is of dimension at least \(k\).
Suppose that there is some point at which the holonomy algebra is still not \(\LieG\).
Then this is true at every point, because the holonomy algebra changes only by adjoint action from one point to another.
Take some point \(m_{k+1}\in M\) near which \(\omega\) is flat, a point \(p_{k+1}\) in the fiber of \(m_{k+1}\), and a curvature bundle element \(\Omega_{k+1}\) taking on some value \(\Omega_{k+1}(u,v)\) at \(p_{k+1}\) outside of the holonomy algebra.
Take some relatively compact open set \(U_{k+1}\) around \(m_{k+1}\), with closure not intersecting any of the closures of the previously chosen open sets \(U_i\).
Over \(U_{k+1}\), pick a connection whose curvature at \(p_{k+1}\) is \(\Omega_{k+1}\).
Use a partition of unity to glue this connection to \(\omega\), so that it agrees with \(\omega\) outside of the preimage in \(\Bun\) of \(U_{k+1}\).

\end{proof}
\begin{theorem}
Take a connected manifold \(M\) of dimension \(2\) or more.
For any principal bundle \(G\to\Bun_G\to M\), point \(p\in\Bun\) and Lie subgroup \(H\subseteq G\), the following are equivalent:
\begin{itemize}
\item
there is a connection on \(\Bun_G\) with holonomy group at \(p\) contained in \(H\)
\item
\(
\Bun_G\cong\amal{\Bun_{H}}{H}{G}
\)
for some \(H\)-bundle \(\Bun_{H}\).
\item
there is a connection on \(\Bun_G\) with holonomy group at \(p\) contained in \(H\) and with restricted holonomy group at \(p\) precisely the identity component of \(H\).
\end{itemize}
\end{theorem}
It is rare, even in explicit examples, that we know how to find the holonomy reduction, as it requires solving Lie equations.
Hence we might not know its tangent spaces, and we cannot directly apply the Ambrose--Singer theorem.
By purely local calculation, we find the smallest Lie algebra \(\LieN\subseteq\LieG\) containing all values of the curvature at all points, at least constraining the holonomy algebra.
Since the curvature transforms in the adjoint representation, this Lie algebra contains all \(G\)-adjoint conjugates of the holonomy algebra.
\begin{theorem}%
[Cartan holonomy envelope \cite{Cartan:1926} pp. 6,8]%
\label{thm:Cartan.holonomy.normal}
\define{theorem!Cartan holonomy envelope}
\define{Cartan!holonomy envelope theorem}
Suppose that \(G\to\Bun\to M\) is a principal bundle with connection on a connected manifold \(M\).
Let \(\LieN\subseteq\LieG\) be the span of the values of the curvature at all points of \(\Bun\).
Then \(\LieN\subseteq\LieG\) is a \(G\)-invariant ideal containing the holonomy algebra.
Suppose that \(N\subseteq G\) is a normal subgroup with Lie algebra \(\LieN\subseteq\LieG\).
The restricted holonomy group at every point of \(\Bun\) lies in \(N\).
If the holonomy group at some point of \(\Bun\) lies in \(N\) then the holonomy group at every point of \(\Bun\) lies in \(N\).
\end{theorem}
\begin{proof}
For each point \(p\in\Bun\), let \(V_p\subset\LieG\) be the span of the values
\[
\Omega(u,v)\in\LieG
\]
of curvature for \(u,v\in T_p\Bun\).
So \(\LieN\) is the span of
\[
\bigcup_{p\in\Bun} V_p\subset\LieG.
\]
The curvature transforms in the adjoint representation, so \(V_p\) does as well, and so \(\LieN\subseteq\LieG\) is an adjoint invariant linear subspace.
Differentiating the adjoint action, \(\LieN\subseteq\LieG\) is a Lie subalgebra and an ideal.

We assume a right bundle.
Holonomy is conjugated by \(G\)-action:
\[
\holonomyGroup{G}[pg]=\Ad_g^{-1}\holonomyGroup{G}[p].
\]
By normalcy of \(N\subseteq G\), the holonomy lies in \(N\) at some point just when it lies in \(N\) at every other point of the same fiber.
But then we can follow horizontal curves to get from fiber to fiber, because \(M\) is connected.
So the holonomy group lies everywhere in \(N\).

Consider the \(1\)-form \(\omega+\LieN\in\nForms{1}{\Bun}\otimes(\LieG/\LieN)\).
From the equation
\[
\Omega=d\omega+\frac{1}{2}\lb{\omega}{\omega},
\]
we see that, modulo \(\omega+\LieN\),
\[
d\omega+\LieN=\Omega+\LieN=0
\]
since \(\Omega\) is valued in \(\LieN\).
By the Frobenius theorem, \(\Bun\) is foliated by submanifolds whose tangent spaces consist in the tangent vectors on which \(\omega\in\LieN\).
In particular, the horizontal directions \(\omega=0\) lie tangent to the leaves.
Hence each leaf contains the orbit of the horizontal vector fields, i.e. the holonomy reduction.
\end{proof}
Suppose that we wish to determine whether the holonomy group of a principal bundle \(G\to\Bun\to M\) with connection, on a connected and simply connected manifold \(M\), is up to conjugacy some group \(L\subseteq G\).
First, check that at every point of \(\Bun\), the curvature values must lie in some conjugate of the Lie algebra \(\LieL\) of \(L\).
Second, look at the subset \(\Bun'\subseteq\Bun\) of points at which the curvature is valued in the Lie algebra \(\LieL\) of \(L\).
Then the holonomy reduction must lie in \(\Bun'\).
We can compute \(\Bun'\) by a local calculation, without knowing how to compute the holonomy group, i.e. without solving any Lie equations, just by calculating the curvature.
Replace \(\Bun'\) by the set of points of \(\Bun'\) at which all flow lines of all vector fields \(v\) with \(v\hook\omega\in\LieL\) stay, at least for a short time, inside \(\Bun'\).
(It might be trickier to find these points, but at least it is still determined by a local condition.)
Repeat this process of replacing \(\Bun'\) by that subset of \(\Bun'\).
After we repeat infinitely often, if \(\Bun'\) is not empty then, as in the proof of theorem~\vref{thm:Cartan.holonomy.normal}, the equations \(\omega+\LieL=0\) satisfy the requirements of the Frobenius theorem so that the horizontal directions are then tangent to the leaves, so contain a holonomy reduction and the holonomy group lies in \(L\), at least locally.
\section{Covariant derivatives}
Suppose that \(G\to\Bun\xrightarrow{\pi}M\) is a principal right (or left) \(G\)-bundle.
A \emph{vertical vector}\define{vertical!vector} is a tangent vector to \(\Bun\) in the kernel of \(\pi'(p)\).
Each vertical vector is uniquely expressed as \(A_{\Bun}(p)\) where
\[
A\in\LieG\mapsto A_{\Bun}(p):=\left.\frac{d}{dt}pe^{tA}\right|_{t=0}\in T_p\Bun.
\]

Pick a connection \(\omega\) on \(\Bun\to M\).
At each point \(p\in\Bun\), the connection determines a \emph{horizontal space}\define{horizontal space} \((0=\omega)\subseteq T_p\Bun\).
We can split every tangent vector uniquely into a sum of horizontal and vertical.
Denote the horizontal part of a tangent vector \(v\in T_p\Bun\) by \(\horizontalPart{v}\).

We are not interested in all possible forms on \(\Bun\), but we are interested in vector-valued forms.
Take a \(G\)-module \(G\xrightarrow{\rho}\GL{V}\); for simplicity we might write \(gv\) to mean \(\rho(g)v\), for \(g\in G\), \(v\in V\).
A form \(\xi\) valued in \(V\) is \emph{equivariant}\define{equivariant form} if
\[
r_g^*\xi=g^{-1}\xi.
\]
To each \(V\)-valued \(k\)-form \(\xi\) on \(\Bun\), we associate its \emph{horizontal part}:\define{horizontal part}
\[
\horizontalPart{\xi}(v_1,\dots,v_k):=
\xi(\horizontalPart{v}_1,\dots,\horizontalPart{v}_k).
\]
The \emph{covariant derivative}\define{covariant derivative} is the horizontal part of the exterior derivative
\[
\covder\xi:=\horizontalPart{(d\xi)}.
\]
\begin{example}
One motivation for this definition: on the connection itself
\[
d\omega=-\frac{1}{2}\lb{\omega}{\omega}+\Omega
\]
so the covariant derivative is
\[
\covder\omega=\Omega,
\]
since \(\lb{\omega}{\omega}\) is purely vertical, as \(\omega=0\) precisely on the horizontal vectors.
\end{example}
Take a principal bundle \(G\to\Bun\to M\) and \(G\)-module \(V\).
Denote by \(\vb{V}^p\) the set of equivariant \(V\)-valued semibasic \(p\)-forms on \(\Bun\).
\begin{lemma}
On \(\vb{V}^{\complex}\),
\[
\nabla \xi:=d\xi+\omega\wedge\xi.
\]
\end{lemma}
\begin{proof}
Note that
\[
\horizontalPart{(d\xi+\omega\wedge\xi)}=
\horizontalPart{(d\xi)}=\nabla\xi.
\]
We need only prove that \(d\xi+\omega\wedge\xi\) is horizontal, i.e. 
\[
0=A_{\Bun}\hook(d\xi+\omega\wedge\xi)
\]
for every \(A\in\LieG\).
Differentiating in \(g\) the equivariance condition, i.e. writing \(g=e^{tA}\) for some \(A\in\LieG\) and differentiating in \(t\) at \(t=0\),
\[
\LieDer_{A_{\Bun}}\xi=-A\xi.
\]
The Cartan formula gives another computation of the Lie derivative:
\[
\LieDer_{A_{\Bun}}\xi=A_{\Bun}\hook d\xi+d(A_{\Bun}\hook\xi).
\]
\end{proof}
Roughly speaking, the covariant derivative takes advantage of the extra data of being equivariant, cancelling out terms in the exterior derivative that are only due to equivariance.
\begin{example}
Take the special case where \(\xi\) is a \(0\)-form, i.e. an equivariant function
\[
\Bun\xrightarrow{f}V.
\]
Every equivariant function determines a section of the associated vector bundle \(\vb{V}:=\amal{\Bun}{G}{V}\).
The covariant derivative is \(\nabla f=df+\omega f\).
\end{example}
\begin{lemma}
Take a principal bundle \(\Bun\to M\) with connection on a connected manifold.
An equivariant function 
\[
\Bun\xrightarrow{f}V
\]
has vanishing covariant derivative just when it is constant on every holonomy reduction.
When this occurs, that constant is invariant under the holonomy group.
\end{lemma}
\begin{proof}
Suppose that \(\nabla f=0\), i.e. \(df=-\rho(\omega)f\).
Then \(f\) is constant along flow lines of horizontal vector fields since \(\omega=0\) on these.
The orbits of the horizontal vector fields are the holonomy reductions.
Being constant, but also equivariant under the holonomy group, \(f\) has value on each holonomy reduction being an element of \(V\) under the holonomy group.
\end{proof}

\section{Vector bundle valued forms}
Take a principal bundle \(G\to\Bun\xrightarrow{\pi}M\) and an associated vector bundle \(\vb{V}\to M\), say \(\vb{V}=\amal{\Bun}{G}{V}\) for some \(G\)-module \(V\).
Denote by \(\vb{V}^p\) the set of equivariant \(V\)-valued semibasic \(p\)-forms on \(\Bun\).
For example, the curvature of a connection is a section of \(\vbg^2\).
\begin{problem}{ck}
Define a vector bundle isomorphism \(\vb{V}^p=\vb{V}\otimes\Lm{p}{T^*M}\).
Prove that, for any \(G\)-equivariant map of principal \(G\)-bundles 
\[
\begin{tikzcd}
\Bun\arrow[r]\arrow[d]&\Bun'\arrow[d]\\
M\arrow[r]&M'
\end{tikzcd}
\]
over a smooth map \(M\to M'\), your vector bundle isomorphisms commute with the obvious pullbacks of differential forms and associated vector bundles.
\end{problem}
\begin{answer}{ck}
First, we define a map \(\vb{V}^p\to\vb{V}\otimes\nForms{p}{M}\), and prove that it is a vector bundle isomorphism.
We will assume that \(G\to\Bun\xrightarrow{\pi}M\) is a left principal bundle; the case of a right principal bundle is a slight change of notation, by defining a left \(G\)-action on \(\Bun\) by
\[
gp:=pg^{-1},
\]
for \(p\in\Bun'\).
Every differential form on \(M\) pulls back to a \(G\)-equivariant semibasic form on \(\Bun\), so the result we want to prove holds when \(V=\R\) is a trivial \(G\)-module.

Take an associated vector bundle \(\vb{V}:=\amal{\Bun}{G}{V}\).
We have a quotient map
\[
\Bun\times V\xrightarrow{q}\vb{V},
\]
which is a linear isomorphism on each fiber over \(\Bun\).
Write each point \(\vb{v}\in\vb{V}\) as a quotient \(\vb{v}=q(p,v)=G(p,v)\), and the points of \(\pi^*\vb{V}\) as pairs \((p,\vb{v})\).
Let \(\underline{V}:=\Bun\times V\) be the trivial vector bundle over \(\Bun\).
So we have a map
\[
\begin{tikzcd}
\underline{V}\arrow[rr,"\varphi"]\arrow[dr]&&\pi^*\vb{V}\arrow[dl]\\
&\Bun
\end{tikzcd}
\]
defined by \(\varphi(p,v):=(p,\vb{v})\) where \(\vb{v}:=q(p,v)\).
This map is a linear isomorphism on each fiber.
By problem~\vref{problem:v.b.fiber.pres}, \(\varphi\) is a vector bundle isomorphism of vector bundles on \(\Bun\).
In particular, for each point \(p\in\Bun\), mapping to a point \(m:=\pi(p)\in M\), 
\[
V\xrightarrow{\varphi_p}\vb{V}_m
\]
is a linear isomorphism.

Take a \(G\)-equivariant \(V\)-valued semibasic form \(\xi\) on \(\Bun\), i.e. \(\xi\) is a section of \(\vb{V}^p\).
To say that \(\xi\) is \(G\)-equivariant here means precisely that, for any point \(p\in\Bun\) and any tangent vectors
\[
u_1,\dots,u_p\in T_p\Bun
\]
and any \(g\in G\), 
\[
\xi_{gp}(g_*u_1,\dots,g_*u_p)=
\xi_p(u_1,\dots,u_p).
\]
Being semibasic means precisely that
\[
\xi_p(u_1,\dots,u_p)
\]
doesn't change if we replace any \(u_i\) by \(u_i+z_i\), for any tangent vector \(z_i\in T_p\Bun\) with \(\pi_*z_i=0\).

Clearly \(\varphi\xi\) is \(\pi^*\vb{V}\)-valued, and still a \(G\)-equivariant semibasic form on \(\Bun\).
Pick a point \(p\in\Bun\) and let \(m:=\pi(p)\in M\).
Take any tangent vectors \(w_1,\dots,w_p\in T_m M\).
Since \(\Bun\xrightarrow{\pi}M\) is a fiber bundle map, it is a submersion, so there are tangent vectors
\[
u_1,\dots u_p\in T_p\Bun
\]
so that
\[
w_1=\pi_* u_1, \dots, w_k=\pi_* u_k.
\]
These are perhaps not unique, but each \(u_i\) is unique up to adding vertical vectors, i.e. replacing by \(u_i+z_i\) for any tangent vectors \(z_i\in T_p\Bun\) with \(\pi_*z_i=0\).
But the expression
\[
\varphi\xi(u_1,\dots,u_p)\in\vb{V}_m
\]
is independent of the choice of these \(u_1,\dots,u_k\), since they are unique up to adding vertical vectors, i.e. vectors tangent to the fibers of \(\Bun\xrightarrow{\pi}M\), and \(\xi\) is semibasic.
Therefore we can set
\[
\bar\xi(w_1,\dots,w_p):=\varphi\xi(u_1,\dots,u_p)\in\vb{V}_m.
\]

Claim: This \(\bar\xi\) is also independent of the choice of point \(p\), depending only on 
\[
m,w_1,\dots,w_p.
\]
Proof: change \(p\) to some other point \(q\in\Bun\) in the same fiber \(\Bun_m\), then \(q=gp\) for a unique element \(g\in G\).
Since \(\pi g=\pi\), for any tangent vector \(u\in T_p\Bun\), 
\[
\pi_*g_*z=(\pi g)_* z=\pi_*u.
\]
Hence if we set \(u'_i:=g_*u_i\in T_q\Bun\), then
\[
\pi_*u'_i=\pi_* u_i=w_i.
\]
So then if we pick any tangent vectors \(u''_i\in T_q\Bun\) with \(\pi_*u''_i=w_i\), we must have
\[
\pi_*(u''_i-u'_i)=0,
\]
so \(u''_i-u'_i\) is vertical, say \(u''_i-u'_i=z_i\), and therefore
\begin{align*}
\xi_q(u''_1,\dots,u''_p)
&=
\xi_q(u'_1+z_1,\dots,u'_p+z_p),
\\
&=
\xi_q(u'_1,\dots,u'_p),
\\
&=
\xi_{gp}(g_*u_1,\dots,g_*u_p),\\
&=
\xi_p(u_1,\dots,u_p).
\end{align*}
By taking a smooth local section of \(\Bun\to M\), we see that \(\bar\xi\) is a smooth differential form.

Taking a smooth map of manifolds \(M\to M'\) and \(G\)-equivariant maps \(\Bun\to\Bun'\) covering that map, we reverse our steps above to see that the differential form pullback and vector bundle pullback map the vector bundle isomorphisms commute.
\end{answer}

\section{Inducing vector bundle connections}
\begin{problem}{vb.onto}
Take a vector bundle morphism
\[
\begin{tikzcd}
V\arrow[rr,"\varphi"]\arrow[dr]&&W\arrow[dl]\\
&M.&
\end{tikzcd}
\]
To each section \(v\) of \(V\) we associate the section \(w:=\varphi\circ v\).
Prove that this linear map of sections to sections is onto just when, for every point \(m\in M\), the linear map
\[
V_m\xrightarrow{\varphi_m}W_m
\]
of fibers is onto.
\end{problem}
\begin{answer}{vb.onto}
Suppose that the linear map of sections to sections is onto.
Take a vector \(w_m\in W_m\).
Using a local trivialization and a bump function, we easily construct a section \(w\) taking on this value at \(m\).
By hypothesis, there is some section \(v\) so that \(w=\varphi\circ v\).
But then \(w_m=\varphi_m v_m\).
So \(\varphi_m\) is onto for each \(m\in M\).

Next we suppose that \(\varphi_m\) is onto for each \(m\in M\).
Suppose first that \(V\) and \(W\) are trivial vector bundles \(V=M\times V_0\), \(W=M\times W_0\) for vector spaces \(V_0,W_0\).
Then \(\varphi\) is just a linear map at each point
\[
V_0\xrightarrow{\varphi_m}W_0,
\]
smoothly varying in \(m\).
Picking bases, \(\varphi\) is represented by a matrix with full rank and at least as many columns as row, say
\[
\varphi_m
=
\begin{pmatrix}
A_m & B_m 
\end{pmatrix}
\]
for matrices \(A_m,B_m\) depending smoothly on \(m\in M\), with \(A_m\) square.
Pick a point \(m_0\in M\).
Permuting basis vectors, we can arrange that \(A_{m_0}\) is invertible.
Therefore \(A_m\) is invertible near \(m=m_0\).
Writing elements of \(V\) as
\[
\begin{pmatrix}
x\\
y
\end{pmatrix}
\]
we can solve
\[
w(m)=A_mx+B_my,
\]
by taking \(y=0\) and inverting \(A_m\), for \(m\) near \(m_0\).
Therefore there is a local solution to our problem near each point of \(M\).

Returning to the general problem, for any vector bundles, locally we can trivialize, so there are local solutions \(v_{\alpha}\) in open sets \(M_{\alpha}\) covering \(M\).
Use a partition of unity subordinate to this open cover.
\end{answer}

Take a Lie group \(G\) and a \(G\)-module \(V\).
So \(G\) acts on \(V\) by a morphism of Lie groups
\[
G\xrightarrow{\rho}\GL{V}.
\]
Differentiate to get a morphism \(\LieG\xrightarrow{\rho}\LieGL{V}\) of Lie algebras.
Consider the exact sequence of \(G\)-modules
\[
0\to\LieK\to\LieG\to\LieGL{V}\to Q\to 0
\]
where \(\LieK\subseteq\LieG\) is the subalgebra of elements \(A\in\LieG\) so that \(\rho(A)=0\) and 
\[
Q:=\LieGL{V}/\rho(\LieG)
\]
is a \(G\)-module.
Take a principal bundle \(G\to\Bun\to M\) and the associated vector bundle  \(\vb{V}=\amal{\Bun}{G}{V}\).
\begin{problem}{Koszul.to.Ehresmann.III}
Suppose that a vector bundle connection on \(\vb{V}\) is associated to a connection on \(\Bun\). 
Prove that that connection on \(\Bun\) is unique precisely up to adding a section of \(\vbk^1\).
Hence the connection on \(\Bun\) inducing a given vector bundle connection is unique if and only if \(\LieK=0\), i.e. \(\LieG\subseteq\LieGL{V}\) is a subalgebra, i.e. \(G\) acts on \(V\) with only a discrete subgroup of \(G\) acting trivially.
\end{problem}
\begin{answer}{Koszul.to.Ehresmann.II}
Take two connections \(\omega,\omega'\) on \(\Bun\).
Take a section \(s\) of \(\vb{V}\), which we identify with a \(G\)-equivariant map \(\Bun\xrightarrow{s}V\).
The associated vector bundle connections \(\nabla,\nabla'\) on \(\vb{V}\) are expressed as operations
\[
\nabla s=ds+\rho(\omega)s, \nabla's=ds+\rho(\omega')s.
\]
So these agree, for all \(s\), just when 
\[
\rho(\omega-\omega')=0
\]
i.e. just when \(\delta:=\omega-\omega'\) is a \(G\)-equivariant semibasic \(1\)-form valued in \(\LieK\subseteq\LieG\).
Conversely, given any \(G\)-equivariant semibasic \(1\)-form \(\delta\) valued in \(\LieK\), we can let \(\omega':=\omega+\delta\).
\end{answer}
\begin{problem}{Koszul.to.Ehresmann.II}
Define a section of \(\vb{Q}^1\) associated to each vector bundle connection \(\nabla\), the ``error'' of \(\nabla\), so that the error vanishes just when \(\nabla\) is associated to a connection on \(\Bun\).
In particular, prove that every vector bundle connection on any vector bundle is associated to a unique connection on the frame bundle of that vector bundle.
\end{problem}
\begin{answer}{Koszul.to.Ehresmann.II}
Take a connection \(\nabla\) on \(\vb{V}\).
Pick a basis of local sections \(s_1,\dots,s_N\) of \(\vb{V}\) defined on some open set \(U\subseteq M\).
Each becomes some \(V\)-valued function \(f_1,\dots,f_N\) on \(\Bun\).
Each has covariant derivative \(\nabla s_i\) a section of \(T^*M\otimes\vb{V}\), which becomes some \(V\)-valued semibasic \(1\)-form \(\nabla f_i\) on \(\Bun\).
We have to ask if there is a connection \(\omega'\) inducing \(\nabla\), i.e. if there is some \(\LieG\)-valued \(G\)-equivariant \(1\)-form \(\omega'\) so that
\[
\rho(\omega')f_i = \nabla f_i-df_i, \quad i=1,2,\dots, N.
\]
Pick a connection \(\omega\) on \(\vb{V}\) over \(U\).
To find a connection \(\omega'\) inducing \(\nabla\), we need precisely that there is a section \(\delta\) of \(\vbg^1\), i.e. a \(\LieG\)-valued \(G\)-equivariant semibasic \(1\)-form, so that
\[
\rho(\delta)f_i=\nabla f_i-df_i-\rho(\omega)f_i, \quad i=1,2,\dots, N,
\]
and then we can set \(\omega'=\omega+\delta\).
We can change basis, from \(s_i\) to some \(s'_i\) say, and this will change these \(f_i\) to some new \(f'_i\), say.
But then \(s'_i=g_i^js_j\) for a unique invertible matrix \(g\in\GL{N}\) at each point of \(M\), and hence \(f'_i=g_i^jf_j\) and
\[
\nabla s'_i=g_i^j\nabla s_j+s_j dg^j_i,
\]
giving
\[
\nabla f'_i=g_i^j\nabla f_j+f_j dg^j_i,
\]
So then
\begin{align*}
\nabla f'_i-df'_i-\rho(\omega)f'_i
=
g^j_i(\nabla f_i-df_i-\rho(\omega)f_i),
\end{align*}
allowing precisely the same possibilities for \(\delta\).

Define \(\varepsilon\in\LieGL{V}\otimes\nForms{1}{\Bun}\) by, for \(v\in T_p\Bun\),
\[
\varepsilon(v)f_i=\nabla f_i(v)-df_i(v)-\rho(\omega(v))f_i, \quad i=1,2,\dots, N.
\]
Clearly \(\varepsilon\) is independent of the choice of basis of local sections \(f_i\).
Moreover, \(\varepsilon(v)=0\) for \(v\) vertical, and \(\varepsilon\) is \(G\)-equivariant.
Consider the projection
\[
a\in\LieGL{V}\otimes T^*M\mapsto\bar{a}\in(\LieGL{V}/\rho(\LieG))\otimes T^*M 
\]
If we let \(W:=\LieGL{V}/\rho(\LieG)\), with associated vector bundle \(\vb{W}:=\amal{\Bun}{G}{W}\), we can define the ``error'' \(\bar\varepsilon\). 
Note that while \(\varepsilon\) depends on \(\omega\), \(\bar\varepsilon\) does not, and therefore is defined on \(M\), not just on \(U\).
Clearly if \(\nabla\) is induced by some connection  on \(\Bun\) then \(\bar\varepsilon=0\).

Claim: \(\bar\varepsilon=0\) just when \(\nabla\) is induced by a connection on \(\Bun\).
Proof: suppose that \(\bar\varepsilon=0\).
So then \(\varepsilon\) is valued in \(\rho(\LieG)\).
By problem~\vref{problem:vb.onto}, there is some section \(\delta\) of \(\vbg^1\) mapping to \(\varepsilon\).
\end{answer}
\begin{lemma}
On \(\vb{V}^{\complex}\),
\[
\nabla^2=\rho(\Omega)\wedge
\]
i.e., on any principal bundle \(\Bun\to M\) with connection, for any equivariant \(p\)-form valued in a \(G\)-module \(G\xrightarrow{\rho}\GL{V}\),
\[
\nabla(\nabla\xi)=\rho(\Omega)\wedge\xi,
\]
where \(\Omega\) is the curvature of the connection.
\end{lemma}
\begin{problem}{nabla.2}
Prove it.
\end{problem}
\begin{answer}{nabla.2}
We leave the reader to fill in some details of the algebra:
\begin{align*}
\nabla^2\xi
&=
(d+\omega\wedge)(d\xi+\omega\wedge\xi),
\\
&=
d^2\xi+d(\omega\wedge\xi)+\omega\wedge d\xi+\omega\wedge\omega\wedge\xi,
\\
&=
d\omega\wedge\xi+\omega\wedge\omega\wedge\xi,
\\
&=
\Omega\wedge\xi.
\end{align*}
\end{answer}
\begin{problem}{difference.of.connections}
Take a principal bundle \(G\to\Bun\xrightarrow{\pi}M\).
Recall that \(\vbg:=\amal{\Bun}{G}{\LieG}\) and that \(\vbg^1=\vbg\otimes T^*M\).
Explain why the difference of any two connections on the bundle is identified with a section of \(\vbg^1\), and conversely, if we fix one connection, any section of \(\vbg^1\) is the difference between that connection and some other connection.
In other words, the space of smooth connections is an affine space modelled on the vector space of sections of \(\vbg^1\).
\end{problem}
\begin{answer}{difference.of.connections}
Let \(\delta:=\omega'-\omega\) be the difference of two connections \(\omega,\omega'\) on \(\Bun\).
Then \(\delta\) is a \(G\)-equivariant \(\LieG\)-valued \(1\)-form on \(\Bun\) vanishing on the fibers.
By problem~\vref{problem:ck}, \(\delta\) is identified with a section of \(\vbg^1\).
\end{answer}
\section{Holonomy of Cartan geometries}\label{page:Cartan.geometry.holonomy}
Suppose now that \(\Bun=\G_G\) is the \(G\)-bundle of a Cartan geometry \(\G\to M\).
Recall that development is parallel transport, i.e. horizontal lifting, i.e. the horizontal lift of a curve in \(M\) to a curve in \(\G_G\) is the projection of a curve in \(\G\times G\) on which the Cartan connection agrees with the Maurer--Cartan form.
The holonomy group of \(\G_G\) at \(p_0\) is the group \(\holonomyGroup{G}\) of values at the end point \(g(b)\) of absolutely continuous curves \(g(t)\in G\), \(a\le t\le b\), given by requiring \(g(a)=1\in G\) and that \(g(t)\) is the development of some absolutely continuous curve \(p(t)\in\G\), i.e. , projecting to a loop \(m(t)\in M\) with \(p(a)=p_0\).
\begin{example}
The holonomy group of the model is \(\set{1}\).
\end{example}
\begin{theorem}
Take a morphism \((X,G)\to(X',G')\) of homogeneous spaces.
Suppose that \(G\to G'\) is injective.
An \((X',G')\)-geometry is induced from an \((X,G)\)-geometry just when there is some point of the total space at which the holonomy lies in the image of \(G\to G'\).
\end{theorem}
For each point \(p_0\in\G\), let \(\mathscr{O}\subset\G\times G\) be set of points \((p(t),g(t))\) which occur on absolutely continuous paths with \(p(0)=p_0\) and \(g(0)=1\).
Note that \(\mathscr{O}\) is thus the orbit through \((p_0,1)\) of the obvious family of vector fields, i.e. those on which \(\omega_{\G}=\omega_G\).
The holonomy reduction at \(p_0\) of \(\G_G\), which we denote by \(\holonomyReduction{\G_G}\), is the \(H\)-quotient
\[
\holonomyReduction{\G_G}=\mathscr{O}/H\subseteq\amal{\G}{H}{G}=\G_G.
\]
\begin{example}
Suppose that \(\G\to M\) is a flat Cartan geometry modelled on a homogeneous space \((X,G)\).
Pick a point \(m_0\in M\).
Let \((\tilde{M},\tilde{m}_0)\xrightarrow{p}(M,m_0)\) be the universal covering space and \(\pi:=\fundamentalGroup{M}\) be the fundamental group.
Denote the pullback Cartan geometry on \(\tilde{M}\) by \(\tilde\G:=p^*\G\), with the pullback Cartan connection.
Since the Cartan connection is flat, it is locally isomorphic to the standard flat connection, and then globally on \(\tilde{M}\) by simply connectivity: \(\tilde\G_G=\tilde{M}\times G\).
Every Cartan geometry on \(\tilde{M}\) whose associated \(G\)-bundle is \(\tilde{G}_G\) and whose associated Cartan connection is the standard flat connection is given precisely by a reduction of structure group to \(H\subseteq G\) which is nowhere tangent to the horizontal, i.e. to \(T_{\tilde{m}}\tilde{M}\oplus 0\).
A reduction of structure group to \(H\) is an \(H\)-invariant submanifold of \(\tilde{G}_G\) whose every fiber is precisely a single \(H\)-orbit, hence of the form \(gH\) for some \(g\in G\).
So a reduction of structure group to \(H\) is precisely a smooth map \(\tilde{M}\to G/H=X\).
This reduction is a Cartan geometry just when it is nowhere tangent to the horizontal, i.e. just when its tangent spaces map by linear isometry to tangent spaces of \(X\), i.e. just when it is a local diffeomorphism \(\tilde{M}\to X\), a developing map.
It then quotients to \(M\) under \(\pi\) just when \(\pi\) acts by automorphisms, i.e. when the developing map is equivariant for a group morphism \(\pi\xrightarrow{h}G\).
We have proved theorem~\vref{thm:flat.first.appearance}.
But we have also found that the holonomy group is the image of \(\pi\xrightarrow{h}G\).
\end{example}
For any Cartan geometry \(\G\to M\), any principal \(H\)-subbundle \(\G'\subset\G_G\) which is nowhere tangent to the horizontal also defines a Cartan geometry, with the same connection on \(\G_G\), so the same holonomy reduction through each point of \(\G_G\).
So we can't guarantee that \(\G\) intersects the holonomy reduction of \(\G_G\) at any more than one point in general, or guarantee that the intersection is smooth.
These \(\G,\holonomyReduction{\G_G}[p_0]\subseteq\G_G\) are, roughly speaking, ``unrelated'' given by the independent data of the reduction of structure group and the Cartan connection, which are related only by the required that the reduction of structure group be nowhere horizontal for the Cartan connection, an open condition.
\begin{theorem}
Take a homogeneous space \((X,G)\), say with stabilizer \(H:=G^{x_0}\) for some point \(x_0\in X\).
Take a flat Cartan geometry \(H\to\G\to M\) on a connected manifold \(M\).
Every connected Lie subgroup of \(G\) is the restricted holonomy group of a Cartan geometry on \(M\) with total space \(\G\), hence is the holonomy group of a Cartan geometry on a covering space of \(M\).
\end{theorem}
\begin{proof}
Pick a connected subgroup \(L\subseteq G\).
Recall our construction of a connection with restricted holonomy \(L\) from lemma~\vref{lemma:make.holonomy.big}: we started with any flat connection \(\omega\) on \(\G\).
We picked an open set \(U\subseteq M\) over which the bundle is trivialized by the flat connnection; any simply connected open set will do, for example.
Trivialize over \(U\), so write points of \(\G_{G,U}\) as \(p=(m,g)\) for some \(g\in G\).
We constructed a new connection \(\omega'\) of the form
\[
\omega'=\omega+\Ad_g^{-1}\Gamma,
\]
where \(\Gamma\) is a \(1\)-form on \(M\), valued in \(\LieL\), supported on a compact subset of \(M\).

The generic choice of \(\Gamma\) will have restricted holonomy \(L\), once we have seen that there is some \(\Gamma\) that has.
Looking back at the construction, we could just as well have used \(\varepsilon\Gamma\), for any \(\varepsilon\ne 0\).
In other words, the set of choices of \(\Gamma\) to make \(\omega'\) have restricted holonomy \(L\) is an open cone, dense in the topology of uniform convergence on compact sets with any number of derivatives.
In particular, we can arrange that \(\Gamma\) is as close as we like to zero, uniformly in any number of derivatives.
Hence we can suppose that the horizontal spaces of \(\omega'\) are uniformly as close as we like to those of \(\omega\).
Thus the new connection \(\omega'\) is still a Cartan connection for the old \(H\)-bundle \(\G_H\).
\end{proof}
\begin{corollary}
Take a homogeneous space \((X,G)\).
Every connected Lie subgroup of \(G\) is the holonomy group of a Cartan geometry.
\end{corollary}
This completes the classification of restricted holonomy groups of Cartan geometries.
But if we ask for some special class of Cartan geometries, for example satisfying some invariant condition on curvature, then it is not clear how to find which holonomy groups can arise from Cartan geometries satisfying that condition.
\begin{theorem}\label{theorem:open.holonomy}
Take a Cartan geometry \(\G\to M\), say with model \((X,G)\), on a connected manifold.
Suppose that there is a connected open set \(U\subseteq\G\) at which the holonomy group contains some subset \(S\subseteq G\).
Let \(G'\) be the smallest subgroup generated by \(S\) and normalized by the identity component of \(G\).
Then the holonomy group of every point in the path component of \(U\) contains \(G'\).
\end{theorem}
\begin{proof}
Suppose that our model is \((X,G)\) and our Cartan geometry is \(\G\to M\).
Let \(G_0\) be the group generated by that set.
Pick a point \(p_0\).
Pick any \(g_1\in G\) close enough to the identity, and there will be a path \(p(t)\) starting at \(p_0\), developing to \(G\) to a path \(g(t)\) with \(g(0)=1\) and \(g(1)=g_1\).
(To prove this: developing along radial lines in some coordinates on \(\G\) centered at \(p_0\), we get curves in \(G\) with arbitrary velocities, so cover an open set in \(G\).)
The holonomy groups are conjugate:
\[
\holonomyGroup{G}[p_0]=g_1^{-1}\holonomyGroup{G}[p_1]g_1
\]
so contains \(g_1^{-1}G_0g_1\).
This holds for all \(g_1\) near \(1\), so holds for all \(g_1\) in the identity component of \(G\).
\end{proof}
From theorem~\vref{theorem:induced.holonomy}:
\begin{theorem}
The Cartan connection of any reductive geometry splits as usual into an affine connection and the soldering form.
The holonomy of the reductive geometry maps to the holonomy of the affine connection.
\end{theorem}
\begin{proof}
We can mutate the model to the affine form, so assume the model is some \((X,G)=(V,H\ltimes V)\).
Then the map \(\G_G\to\G_H=\G_G/V\) takes connection to connection.
\end{proof}
\begin{example}
The holonomy group, as Cartan geometry, of any pseudo-Riemannian metric, maps to the holonomy group of its Levi--Civita connection.
\end{example}
\section{The curvature algebra}
Denote the Cartan connection as \(\omega:=\omega_{\G}\) and its curvature as \(\Omega:=\Omega_{\G}\).
As we saw in lemma~\vref{lemma:forms.on.product}, on \(\mathscr{O}\), the connection form on \(\G_G\) pulls back to
\[
\tilde\omega:=\Ad_g(\omega-\omega_G)
\]
under the map we called \(\pi'\circ\iota\).
Taking exterior derivative on \(\mathscr{O}\), the curvatures are thus related by
\[
\tilde\Omega:=\Ad_g\Omega.
\]
Writing
\[
\Omega=\frac{1}{2}k\sigma\wedge\sigma
\]
we have \(\Ad_g k\in\holonomyAlgebra{G}\otimes\Lm*{2}{\LieG/\LieH}^*\) at every point 
\[
(p,g)\in\mathscr{O}\subseteq\G\times G.
\]
\begin{lemma}
Suppose that \(\G\to M\) is a Cartan geometry modelled on a homogeneous space \((X,G)\) with \(X=G/H\).
Take a \(G\)-module \(V\) and a section \(f\) of the associated vector bundle on \(M\).
Pick a point \(p_0\in\G\) and let \(\holonomyGroup{G}:=\holonomyGroup{G}[p_0]\) 
Suppose that \(\dot{V}\subseteq V\) is a \(\holonomyGroup{G}\)-submodule and that \(gf\) is valued in \(\dot{V}\) on the preimage \(\mathscr{O}\subset\G\times G\) of the holonomy reduction of \(p_0\).
Then, on \(\G\), the covariant derivative 
\[
\nabla^{\G} f=df+\omega f=f'\sigma
\]
has value \(f'(p_0)\) lying in \(\dot{V}\otimes(\LieG/\LieH)^*\).
\end{lemma}
\begin{problem}{covariant.derivative.holonomy}
Give a proof.
\end{problem}
\begin{answer}{covariant.derivative.holonomy}
Since \(gf\) is valued in \(\dot{V}\) on \(\mathscr{O}\), trivially \(d(gf)\) is valued there too, on any tangent vector to \(\mathscr{O}\).
Compute on \(\mathscr{O}\):
\begin{align*}
d(gf)
&=
dg\,f+g\,df,
\\
&=
(\Ad_g\omega_G)gf+g\,(\nabla^{\G} f-\omega f),
\\
&=
(\Ad_g\omega-\tilde\omega)gf+gf'\sigma-(\Ad_g\omega)f,
\\
&=
-\tilde\omega gf+gf'\sigma.
\end{align*}
The first term, \(\tilde\omega gf\), we already know is valued in \(\dot{V}\) on any tangent vector to \(\mathscr{O}\).
So \(gk'\sigma\) is valued in \(\dot{V}\) on any tangent vector to \(\mathscr{O}\).
Note that \(\sigma\) is a \(1\)-form on \(\mathscr{O}\) pulled back from \(\G\) via the submersion
\[
\mathscr{O}\subseteq\G\times G\to\G.
\]
So at the point \((p,g)=(p_0,1)\), \(f'\) is valued in \(V'\otimes(\LieG/\LieH)^*\).
\end{answer}
The \emph{curvature algebra}\define{curvature!algebra} of a Cartan geometry \(\G\to M\) at a point \(p_0\in\G\) is the set of values
\[
k(u_1,u_2),k'(u_1,u_2,u_3),\dots
\]
of all covariant derivatives at \(p_0\) of curvature of all orders, for any
\[
u_1,u_2,\dots\in\LieG/\LieH.
\]
\begin{corollary}
The curvature algebra of any Cartan geometry lies in the holonomy algebra.
\end{corollary}
\section{Example: surfaces with Riemannian metric}
Suppose that \((X,G)\) is the Euclidean plane \(X=\R^2\) under its orientation preserving isometry group 
Recall that development is precisely rolling, for surfaces with Riemannian metric.
The holonomy of the Cartan geometry of a surface with Riemannian metric is thus \(G\) precisely when we can rolling the surface on the plane \(X\), and get from any position and rotation on the plane to any other.
We will see we can do this precisely when the surface is not flat.

Write \(G=\SO{2}\ltimes\R^2\).
We can identify \(\R^2=\C\) and then \(\SO{2}\) is the unit circle in the complex plane.
As usual we write each isometry \(g\in G\) as a \(2\times 2\) matrix
\[
g=
\begin{pmatrix}
h&v\\
0&1
\end{pmatrix}
\]
with \(h=e^{i\theta}\in \SO{2}\) and \(v\in V=\C\).
Let \(\G\to M\) be the Riemannian geometry of an oriented surface with Riemannian metric.
We will prove that the restricted holonomy is \(\{1\}\) or \(G\).

Denoting the Cartan connection as
\[
\omega=
\begin{pmatrix}
-i\gamma&\sigma\\
0&0
\end{pmatrix}
\]
as usual, we have the structure equations
\begin{align*}
d\sigma-i\gamma\wedge\sigma&=0,\\
d\gamma&=K\,dA
\end{align*}
of surface geometry where
\[
dA:=\sigma^1\wedge\sigma^2
\]
is the area form on the sphere \(M\) and \(K\) is the Gauss curvature.
The restricted holonomy is trivial \(\holonomyGroup{G}=\set{1}\) just when the surface develops to the Euclidean plane isometrically without holonomy, so \(M\) has some immersion to the plane, by which its orientation and Euclidean metric are pulled back.

The curvature of the Cartan connection is
\[
\Omega
=
d\omega+\omega\wedge\omega=
\begin{pmatrix}
-i&0\\
0&0
\end{pmatrix}K\,dA.
\]
So the curvature function is
\[
k=
\begin{pmatrix}
-i&0\\
0&0
\end{pmatrix}K.
\]
Hence the curvature algebra, at every point where the Gauss curvature \(K\ne 0\), contains
\[
\begin{pmatrix}
i&0\\
0&0
\end{pmatrix}.
\]
The holonomy algebra contains the curvature algebra.
The holonomy group contains the exponentials of the holonomy algebra.
So at points where the Gauss curvature is not zero, the holonomy group contains
\[
\begin{pmatrix}
e^{i\varphi}&0\\
0&1
\end{pmatrix}
\]
for any \(\varphi\).
By theorem~\vref{theorem:open.holonomy}, the holonomy group at every point contains all conjugates of these:
\begin{align*}
\Ad_g
\begin{pmatrix}
e^{i\varphi}&0\\
0&1
\end{pmatrix}
&=
\begin{pmatrix}
e^{i\theta}&v\\
0&1
\end{pmatrix}
\begin{pmatrix}
e^{i\varphi}&0\\
0&1
\end{pmatrix}
\begin{pmatrix}
e^{-i\theta}&-e^{-i\theta}v\\
0&1
\end{pmatrix},
\\
&=
\begin{pmatrix}
e^{i\varphi}&(1-e^{i\varphi})v\\
0&1
\end{pmatrix}.
\end{align*}
Being \(3\)-dimensional, this is already an open subset of \(G\), which is connected, so the holonomy group is \(G\).

Let us look at this problem from another direction.
The curvature on \(\G_G\) pulls back to \(\G\times G\) to become
\begin{align*}
\tilde\Omega
&=
\Ad_g\Omega,
\\
&=
\begin{pmatrix}
e^{i\theta}&v\\
0&1
\end{pmatrix}
\begin{pmatrix}
-i&0\\
0&0
\end{pmatrix}K\,dA
\begin{pmatrix}
e^{-i\theta}&-e^{-i\theta}v\\
0&1
\end{pmatrix},
\\
&=
\begin{pmatrix}
-ie^{i\theta}&0\\
0&0
\end{pmatrix}
\begin{pmatrix}
e^{-i\theta}&-e^{-i\theta}v\\
0&1
\end{pmatrix}
K\,dA,
\\
&=
\begin{pmatrix}
-i&-iv\\
0&0
\end{pmatrix}
K\,dA.
\end{align*}
Notice that nonzero Gauss curvature causes both rotation and translation components to the curvature on \(\G_G\), which we did not see on \(\G=\G_H\).

Note that \((p,g)=(p_0,I)\in\mathscr{O}\) by hypothesis.
Recall that \(\mathscr{O}\) is the set of points we can arrive at on developments of paths in \(\G\).
In particular, there is no constraint on the possible velocities of \(g\) at any point of \(\mathscr{O}\) in the tangent space \(T_g G\).
Hence we can move
\[
g=
\begin{pmatrix}
e^{i\theta}&v\\
0&1
\end{pmatrix}
\]
to any point of \(G\), i.e. make \(\theta\) and \(v\) move independently and reach any values, while staying inside \(\mathscr{O}\).

At each point of \(\mathscr{O}\), \(\tilde\Omega\) is valued in \(\holonomyAlgebra{G}\).
The holonomy algebra contains the span of
\[
\begin{pmatrix}
-i&-iv\\
0&0
\end{pmatrix}
\]
for each point \((p,g)\in\mathscr{O}\) wherever \(G(p)\ne 0\).
Since any point \(g\in G\) arises from some point \((p,g)\in\mathscr{O}\), the holonomy algebra contains an open set of elements of the form
\[
\begin{pmatrix}
-i&-iv\\
0&0
\end{pmatrix}
\]
and so contains their span: \(\holonomyAlgebra{G}=\LieG\).
\section{Example: the complex affine line as model}
Let \((X,G)\) be the complex affine line \(X=\C\) with \(G\) the group of complex affine transformations
\[
x\in X\mapsto ax+b\in X,
\]
for \(a,b\in\C\) with \(a\ne 0\).
We can write each element of \(G\) as a matrix
\[
g=
\begin{pmatrix}
a&b\\
0&1
\end{pmatrix}.
\]
Every oriented surface with a Riemannian metric has an induced \((X,G)\)-geometry.
So we have a large collection of examples, with holonomy contained in the isometry group of the oriented Euclidean plane.

Take an \((X,G)\)-geometry on a connected surface \(M\).
The Cartan connection is valued in \(\LieG\) so has the form
\[
\omega
=
\begin{pmatrix}
\gamma&\sigma\\
0&0
\end{pmatrix}
\]
where \(\gamma,\sigma\) are complex valued \(1\)-forms.
The semibasic forms are generated by \(\sigma,\bar\sigma\).
The curvature is
\begin{align*}
\Omega&=d\omega+\omega\wedge\omega,\\
&=
d
\begin{pmatrix}
\gamma&\sigma\\
0&0
\end{pmatrix}
+
\begin{pmatrix}
\gamma&\sigma\\
0&0
\end{pmatrix}
\wedge
\begin{pmatrix}
\gamma&\sigma\\
0&0
\end{pmatrix},
\\
&=
\begin{pmatrix}
A&\bar{B}\\
0&0
\end{pmatrix}
\frac{i}{2}\sigma\wedge\bar\sigma
\end{align*}
for some complex valued functions \(A,B\) on \(\G\).
We can write this as
\begin{align*}
d\gamma&=A\frac{i}{2}\sigma\wedge\bar\sigma,\\
d\sigma+\gamma\wedge\sigma&=\bar{B}\frac{i}{2}\sigma\wedge\bar\sigma.
\end{align*}
Note that \(B\) is the torsion.
The Cartan geometry is complex analytic (hence no \(\bar\sigma\) terms) just when flat.
\begin{problem}{cl}
Prove that any homogeneous \((X,G)\)-geometry with automorphism group of real dimension \(4\) is flat; classify them globally.
\end{problem}
Expand out 
\begin{align*}
r_h^*\omega&=\Ad_h^{-1}\omega,\\
r_h^*\Omega&=\Ad_h^{-1}\Omega,
\end{align*}
to find, for 
\[
h=
\begin{pmatrix}
a&0\\
0&1
\end{pmatrix},
\]
\begin{align*}
r_h^*\gamma&=\gamma,\\
r_h^*\sigma&=a^{-1}\sigma,\\
r_h^*A&=|a|^2A,\\
r_h^*B&=aB.
\end{align*}
In particular, the \(1\)-form \(\tau:=B\sigma\) is basic, a \((1,0)\)-form on the Riemann surface \(M\), which we can think of as the torsion \((1,0)\)-form.

Clearly \(0=A=B\) just for flat \((X,G)\)-geometries, which we know are given by developing maps.
Equivalently, they are holomorphic affine connections on Riemann surfaces.
For example, there are no flat \((X,G)\)-geometries on the sphere, since there is no local diffeomorphism to the complex affine line \(X=\C\).
Moreover, there are none on compact surfaces other than the torus, since a flat holomorphic affine connection forces vanishing of all characteristic classes \cite{Atiyah:1957}.
The flat complex affine structures on tori are well known; we saw them~\vpageref{example:elliptic.curve.affine}.

At every point of the orbit \(\mathscr{O}\), 
\[
\tilde\Omega
=
g
\begin{pmatrix}
A&\bar{B}\\
0&0
\end{pmatrix}
g^{-1}
=
\begin{pmatrix}
A&a\bar{B}-bA\\
0&0
\end{pmatrix}.
\]
\subsection{Example: the complex affine as model: torsion without curvature}
Continuing with \((X,G)\) the complex affine line, suppose now that \(A=0\) everywhere but \(B\ne 0\) at some point of \(\G\).
Throughout the orbit, we find the values of curvature 
\[
\begin{pmatrix}
A&a\bar{B}-bA\\
0&0
\end{pmatrix}=
\begin{pmatrix}
0&a\bar{B}\\
0&0
\end{pmatrix}
\]
span all complex matrices of the form
\[
\begin{pmatrix}
0&B'\\
0&0
\end{pmatrix}.
\]
Hence the restricted holonomy group is the group of matrices of the form
\[
\begin{pmatrix}
1&b\\
0&1
\end{pmatrix}
\]
for \(b\in\C\).

We can locally explicitly solve for all such.
Note that \(d\gamma=0\).
By the Poincar\'e lemma, we can locally construct a complex valued function \(Z\) so that \(\gamma=Z^{-1}dZ\).
Compute
\begin{align*}
d(Z\sigma)
&=
dZ\wedge\sigma+Z\,d\sigma,
\\
&=
Z^{-1}dZ\wedge Z\sigma+Z(-\gamma\wedge\sigma+\bar{B}\frac{i}{2}\sigma\wedge\bar\sigma),
\\
&=
\frac{i\bar{B}}{2\bar{Z}}(Z\sigma)\wedge(\bar{Z}\bar\sigma).
\end{align*}
So \(Z\sigma\) is basic, say \(Z\sigma=w(z)dz\) for some local complex coordinate \(z\) on \(M\), so
\[
\sigma=\frac{w}{Z}dz,
\]
and we compute that
\[
\tau=B\sigma=-\frac{2i}{\bar{w}}\overline{\frac{\partial w}{\partial\bar{z}}}dz.
\]

Conversely, take any \((1,0)\)-form \(\tau\) on a Riemann surface \(M\), and in any local coordinates, identifying \(M\) with an open set \(U\subseteq\C\), write \(\tau\) as
\[
\tau=F\,dz
\]
for some smooth complex valued function \(F=F(z)\).
By the Cauchy--Pompeiu formula \cite{Griffiths.Harris:1978} p. 2, if we take a relatively compact domain \(U\) in the domain of the complex coordinate \(z\), with piecewise \(C^1\) boundary, and let
\[
W(z):=h(z)-\frac{1}{4\pi}\int_U \frac{d\zeta}{\zeta-z}\wedge\bar\tau,
\]
where \(h(z)\) can be any holomorphic function in \(U\), then \(W\) solves
\[
\frac{\partial W}{\partial\bar{z}}=-\frac{i\bar{F}}{2},
\]
in \(U\).
Let \(w:=e^W\), and we have solved the linear elliptic partial differential equation
\[
\frac{\partial w}{\partial\bar{z}}=-\frac{iw\bar{F}}{2}.
\]
for an unknown complex valued function \(w(z)\ne 0\).
Let \(\G:=U\times\C^{\times}\), with coordinates \((z,Z)\) and forms
\begin{align*}
\gamma&:=\frac{dZ}{Z},\\
\sigma&:=\frac{w}{Z}dz.
\end{align*}
We have found all local forms of the \((X,G)\)-geometries with \(A=0\), hence with holonomy in the group of matrices of the form
\[
\begin{pmatrix}
1&b\\
0&1
\end{pmatrix}
\]
modulo the problem of local solvability of the differential equation.
The restricted holonomy group is exactly this group unless \(\tau=0\) everywhere.
Note that the zero locus of \(\tau\) on the surface \(M\) is an invariant of the geometry, so there are clearly infinitely many examples which are not locally isomorphic.
\begin{problem}{cm}
If
\[
\gamma':=\gamma+\frac{i}{2}\bar\tau-\frac{i}{2}\tau
\]
check that \(\gamma'\) is a torsion-free connection with real valued curvature.
\end{problem}
\subsection{Example: the complex affine as model: curvature without torsion}
Again we consider geometries modelled on the complex affine line.
Suppose that \(B=0\) everywhere, a torsion-free \((X,G)\)-geometry.
We can suppose that \(A\ne 0\) on some open set.
As above, the holonomy algebra of any point in that open set contains
\[
\begin{pmatrix}
A&b\\
0&0
\end{pmatrix}
\]
for that \(A\) and any \(b\in\C\).
So the restricted holonomy algebra is \(G\) or is the group of matrices of the form
\[
\begin{pmatrix}
e^{tA_0}&B\\
0&1
\end{pmatrix},
\]
where \(A_0\in\C\) is our curvature at that point \(A_0:=A(p_0)\), and \(B\) can be any complex number.
The subgroup of matrices of the form
\[
\begin{pmatrix}
e^{tA_0}&0\\
0&1
\end{pmatrix}
\]
is normal in \(G\).
By corollary~\vref{corollary:normal.sub.hol}, this subgroup lies in the restricted holonomy of every point.
But then if some point \(p_1\) has curvature \(A_1:=A(p_1)\) not a real multiple of \(A_0\), the holonomy algebra at \(p_1\) contains
\[
\begin{pmatrix}
A_0&0\\
0&0
\end{pmatrix},
\begin{pmatrix}
A_1&0\\
0&0
\end{pmatrix},
\]
and so contains their span.
By the same argument, so do all points, so the holonomy is \(G\).
Hence we can suppose that \(B=0\) and that has the same form, i.e. \(A(p_1)\) is a real multiple of \(A(p_0)\).

Call an \((X,G)\)-geometry with \(B=0\) and holonomy group a strict subgroup of \(G\) \emph{nondegenerate} if this multiple \(r\) is positive.
Every  such geometry is nondegenerate at least on a nonempty open set.
So \(A=re^{i\alpha_0}\) for a real constant \(\alpha_0\) and a smooth function \(r>0\) on \(\G\).
Under the structure group action, check that
\[
r_h^*r=|h|^2r.
\]
Therefore \(r=1\) on a circle of points in each fiber of \(\G\to M\): the orbit of the unit circle in \(H=\C^{\times}\).
Differentiate
\[
d\gamma=r e^{i\alpha_0}\frac{i}{2}\sigma\wedge\bar\sigma,
\]
to find that 
\[
0=(dr-r(\gamma+\bar\gamma))\wedge\sigma\wedge\bar\sigma,
\]
so by Cartan's lemma,
\[
dr=r(\gamma+\bar\gamma)-2is\sigma+2i\bar{s}\bar\sigma
\]
for some complex valued function \(s\).

Therefore the set of points at which \(r=1\) is a smooth principal \(S^1\)-bundle \(\G_0\subseteq\G\).
On that bundle, write \(\gamma\) in real and imaginary parts as \(\gamma=\alpha-i\beta\).
Compute that
\[
\gamma+\bar\gamma=2is\sigma-2i\bar{s}\bar\sigma,
\]
so that
\[
\alpha=is\sigma-i\bar{s}\bar\sigma.
\]
Let
\[
\beta':=\beta+s\sigma+\bar{s}\bar\sigma.
\]
Check that
\[
d\sigma=i\beta'\wedge\sigma.
\]
So \(\beta'\) is the Levi--Civita connection of an oriented Riemannian metric on \(M\), in the conformal class of the \((X,G)\)-geometry.
It is convenient to write \(dA\) for the area form of the metric
\[
dA:=\frac{i}{2}\sigma\wedge\bar\sigma.
\]

Note that
\begin{align*}
d\gamma&=e^{i\alpha_0}\frac{i}{2}\sigma\wedge\bar\sigma,\\
&=e^{i\alpha_0}dA,\\
&=\cos(\alpha_0)dA+i\sin(\alpha_0)dA,\\
&=d\alpha-i\,d\beta,
\end{align*}
so that
\begin{align*}
d\alpha&=\cos(\alpha_0)dA,\\
d\beta&=-\sin(\alpha_0)dA.
\end{align*}

Differentiate \(\alpha=is\sigma-i\bar{s}\bar\sigma\) to find that
\begin{align*}
d\alpha&=i(ds+is\beta')\wedge\sigma-i(d\bar{s}-i\bar{s}\beta')\wedge\bar\sigma,\\
&=\cos(\alpha_0)dA.
\end{align*}
So \(ds+is\beta'\) is semibasic; it is convenient to write it, for some complex \(p,q\), as
\[
ds+is\beta'=p\sigma-\frac{q}{4}\bar\sigma.
\]
Plug in to find that
\[
q_0=\cos(\alpha_0).
\]
Note that \(\xi:=s\sigma\) is a semibasic \((1,0)\)-form, i.e. a complex linear form, with
\begin{align*}
d\xi&=d(s\sigma),\\
&=ds\wedge\sigma+s\,d\sigma,\\
&=(-is\beta'+p\sigma-\frac{q}{4}\bar\sigma)\wedge\sigma+s\,i\beta'\wedge\sigma,
\\
&=
\frac{q}{4}\sigma\wedge\bar\sigma,
\\
&=
-i\frac{q}{2}\,dA.
\end{align*}
In particular, \(\xi\) is basic, and so \(q\) is also basic.
Note that \(\beta'=\beta+\xi+\bar\xi\).

The curvature of the metric associated to \(\beta'\) is
\begin{align*}
d\beta'&=d\beta+d\xi+d\bar\xi,
\\
&=-\sin(\alpha_0)dA-i\frac{q}{2}\,dA+i\frac{\bar{q}}{2}\,dA,
\\
&=
(-\sin(\alpha_0)-\frac{i}{2}(q-\bar{q}))dA,
\\
&=
(-\sin(\alpha_0)+q_1)dA.
\end{align*}
So the Gauss curvature is \(q_1-\sin(\alpha_0)\).

Suppose we start with any oriented surface \(M\) with Riemannian metric.
Denote its oriented orthonormal frame bundle as \(\G_0\to M\).
Write its structure equations as
\begin{align*}
d\sigma'&=i\beta'\wedge\sigma',\\
d\beta'&=K'\,dA
\end{align*}
where again \(dA:=\frac{i}{2}\sigma'\wedge\bar\sigma'\).
Pick any unit complex number \(e^{i\alpha_0}\).
Suppose we can find a \((1,0)\)-form \(\xi\) on \(M\) so that
\[
d\xi=\frac{1}{2}(K'-ie^{i\alpha_0})dA.
\]
Consider the bundle \(G:=\G_0\times\R^+\).
On that bundle, we let
\begin{align*}
\alpha:=\frac{dr}{r}+i(\xi-\bar\xi),\\
\beta:=\beta'-\xi-\bar\xi,\\
\gamma:=\alpha-i\beta, \\
\sigma:=r\sigma'.
\end{align*}
Tracing back our steps, we see that we have reconstructed the torsion-free \((X,G)\)-geometries with restricted holonomy group precisely the set of matrices
\[
\begin{pmatrix}
e^{te^{i\alpha_0}}&B\\
0&1
\end{pmatrix}
\]
for \(t\in\R\), \(B\in\C\).
So the construction of all torsion-free \((X,G)\)-geometries reduces to solving
\[
2d\xi=(K'-ie^{i\alpha_0})dA.
\]
Locally, this has a solution \(\xi\) by the \(\bar\partial\)-Poincar\'e lemma.

Consider the global problem, for \(M\) compact.
We can assume that \(M\) is connected.
Integration gives
\[
0=2\pi\chi_M-ie^{i\alpha_0}\area{M}.
\]
In particular, the Euler characteristic can't vanish for any compact nondegenerate \((X,G)\)-geometry.
Note that \(e^{i\alpha_0}=\pm i\), say \(e^{i\alpha_0}=(-1)^{\varepsilon}i\) where \(\varepsilon=0\) or \(1\). 
Our Euler characteristic is
\[
\chi_M=\frac{(-1)^{1+\varepsilon}\area{M}}{2\pi}
\]
which we can write in terms of the genus \(g_M\) as
\[
g_M=1+(-1)^{\varepsilon}\frac{\area{M}}{4\pi}.
\]
So if \(M\) is the sphere, \(\varepsilon=1\), \(\area{M}=4\pi\), while if \(M\) is of genus \(2\) or more,
\[
\area{M}=4\pi(g_M-1).
\]
The holonomy group, for \(M\) the sphere or of higher genus, is then
\[
\begin{pmatrix}
e^{it}&b\\
0&1
\end{pmatrix}
\]
for \(t\in\R\).

Conversely, take any oriented compact surface with Riemannian metric, and with genus \(0\) and area \(4\pi\), or genus at least \(2\) and area \(4\pi(g_M-1)\).
Denote its Gauss curvature by \(K'\).
Recall that on \((1,0)\)-forms on a Riemann surface, \(d=\bar\partial\) \cite{Griffiths.Harris:1978}.
By the Hodge isomorphism, and the computation of Dolbeault cohomology, there is a \((1,0)\)-form \(\xi\) with
\[
2\bar\partial\xi=(K'+(-1)^{\varepsilon})dA
\]
where \(\varepsilon=0\) for \(M\) a sphere and \(\varepsilon=1\) for \(M\) of genus \(2\) or more.
Moreover, this \(\xi\) is unique on the sphere, unique up to adding a holomorphic \(1\)-form on any surface of genus \(2\) or more.
We can make explicit homogeneous examples by taking constant curvature oriented surfaces, scaled to have suitable area, with \(\xi=0\).

\begin{appendices}
%\appendix
\chapter{Principal bundles and vector bundles}
\section{Example: frames of a surface}
Take a \(2\)-dimensional vector space \(V\) with a positive definite inner product and fix an orientation of \(V\).
A \emph{frame} on \(V\) is a positively oriented ordered orthonormal basis.
(We will, at times, use the term \emph{frame} to mean many similar constructions.)
There is no canonical choice of frame on \(V\).
But if we pick one frame, then all others are given uniquely by rotating it by some angle.
So the set of frames is a homogeneous space of the group \(G=\SO{2}\) with trivial stabilizer.
Picture the frames as a ``copy'' of \(G\): once we fix one frame, all others are given uniquely by rotating that one.
Since there is no canonical choice of frame on an abstract vector space, there is no canonical choice of how to identify each frame with a rotation matrix, i.e. a point of \(G\).

Take  an oriented smooth surface \(M\) smoothly immersed in \(3\)-dimensional Euclidean space.
A \emph{frame} at a point \(m\in M\)  is a frame of the vector space \(T_m M\).
\[
\includegraphics[width=8cm]{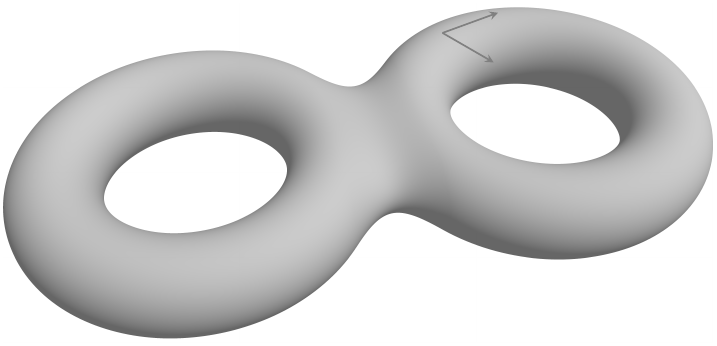}
\]
The \emph{frame bundle}\define{frame bundle} \(\Bun\) of the surface is the set  of all its frames.
(We will, at times, use the term \emph{frame bundle} for many similar constructions.)
The group \(G:=\SO{2}\) of rotations of the plane acts on the frame bundle, by rotating each frame by some angle.
The frames at a given point \(m\) are rotations of one another, so the set of frames at a given point \(m\) is a homogeneous space of \(G\) with trivial stabilizer.
Picture the frames at a point \(m\) as a ``copy'' of \(G\): once we fix one frame, all others at the same point \(m\) are given uniquely by rotating that one.
Since there is no canonical choice of frame on a surface, there is no canonical choice of how to identify each frame with a rotation matrix, i.e. a point of \(G\).

Locally, on some sufficiently small open set of our surface, we can choose a pair of smoothly varying linearly independent vector fields.
Applying Gram--Schmidt orthogonalization, we can arrange that they form a positively oriented orthonormal basis at all points of that open set: a \emph{framing}.\define{framing}
At each point of that open set, every frame becomes a rotation of the framing at that point, hence has an associated rotation matrix.
Globally, it might be impossible to pick a framing, for example on the sphere \(M=S^2\).
Even locally, you and I might make different choices of framing, as there is no canonical choice.

We want to capture this idea of a family of ``copies'' of a Lie group \(G\), parameterized by the points \(m\in M\) of a  manifold \(M\), but without a canonical choice of identity element.
Such a family is, roughly speaking, what we call a \emph{principal bundle}.
The precise definition will take some effort to reach.
We can very roughly picture any principal bundle as if its base manifold is a surface, and each point of its total space is a frame at a point of the surface, as Cartan did \cite{Cartan:1926}.
\section{Definition}
\begin{example}
Consider the trivial example.
Take a Lie group \(G\) and a manifold \(M\); the \emph{standard trivial principal \(G\)-bundle}\define{standard trivial principal bundle}\define{principal bundle!trivial}\define{trivial principal bundle!standard} is \(\Bun:=M\times G\), with the smooth map \(\Bun\xrightarrow{\pi}M\) given by \(\pi(m,g)=m\).
If we equip it with the left \(G\)-action \(g_0(m,g_1)=(m,g_0g_1)\), it is the \emph{left standard trivial principal \(G\)-bundle}.
If, instead, we equip it with the right \(G\)-action \((m,g_0)g_1=(m,g_0g_1)\), it is the \emph{right standard trivial principal \(G\)-bundle}.
\end{example}
Take a Lie group \(G\) and a manifold \(M\).
A \(G\)-\emph{prebundle}\define{prebundle} on \(M\) is a smooth map \(\Bun\xrightarrow{\pi}M\) invariant under an action of \(G\) on \(\Bun\).
The \emph{total space} is \(\Bun\); the \emph{base space} is \(M\); the \emph{structure group} is \(G\).
An \emph{morphism}\define{morphism!principal}\define{principal bundle morphism} of \(G\)-prebundles \(\Bun\to M\) and \(\Bun'\to M'\) is a commutative diagram of smooth maps
\[
\begin{tikzcd}
\Bun\arrow[r, "\Phi"]\arrow[d]&\Bun'\arrow[d] \\
M\arrow[r,"\varphi"]&M'
\end{tikzcd}
\]
so that \(\Phi\) is \(G\)-equivariant.
It is an \emph{isomorphism} if \(\Phi\) and \(\varphi\) are diffeomorphisms.
A \emph{trivialization}\define{trivialization} of a prebundle is an isomorphism with the standard trivial bundle over some manifold.
A prebundle which admits a trivialization is \emph{trivial}.\define{trivial!principal bundle}\define{principal bundle!trivial}
The \emph{restriction}\define{restriction!principal bundle}\define{principal bundle!restriction} of a prebundle \(\Bun\xrightarrow{\pi}M\) over an open set \(U\subseteq M\) is the prebundle with total space \(\Bun_U:=\pi^{-1}U\), base space \(U\), and map the restriction
\[
\Bun_U\xrightarrow{\left.\pi\right|_U}U.
\]
A prebundle is \emph{locally trivial}\define{locally trivial} if every point of the base space lies in an open set over which the restriction is trivial.
A \emph{principal bundle}\define{principal bundle}\define{bundle!principal} is a locally trivial prebundle.
(We often write \emph{bundle} to mean \emph{principal bundle}.)
\section{An equivalent definition}
A \emph{section}\define{section} of a smooth map \(X\to Y\) of manifolds is a smooth map \(Y\to X\) so that the composition \(Y\to X\to Y\) of the two maps is the identity map.
A \emph{local section} of a smooth map \(X\to Y\) is a smooth map \(U\to X\) of an open subset \(U\subseteq Y\) so that the composition \(U\to X\to Y\) of the two maps is the inclusion map \(U\to Y\).
A smooth map \(X\to Y\) \emph{admits local sections} if the domains of local sections cover \(Y\).
\begin{example}
Every surjective submersion admits local sections, by the implicit function theorem.
\end{example}
\begin{problem}{submer}
Prove that a smooth map \(X\to Y\) has a local section near a point \(y\in Y\) if and only if there is a point of \(X\) mapped to \(y\)  near which the map is a submersion.
Give the simplest example you can which is not a submersion.
\end{problem}
The \emph{fiber}\define{fiber} \(X_y\) of a smooth map \(X\xrightarrow{\varphi}Y\) over a point \(y\in Y\) is the set \(X_y:=\varphi^{-1}\set{y}\).
\begin{problem}{connections.loc.triv}
Prove that a prebundle is a principal bundle just when all of following hold:
\begin{itemize}
\item
the group acts freely
\item 
the orbits of the group are precisely the fibers of the map
\item
the orbits are homogeneous spaces of the group and
\item
the map admits local sections.
\end{itemize}
\end{problem}
\begin{answer}{connections.loc.triv}
Clearly a local trivial prebundle admits local sections, by trivializing first and then composing with \(m\mapsto (m,1)\).

Suppose that \(\Bun\xrightarrow{\pi}M\) is a prebundle satisfying those conditions.
Each point \(m\in M\) lies in an open set \(U\subseteq M\) on which there is a local section.
The composition \(U\to\Bun\to M\) is the inclusion map \(U\to M\), so is injective with full rank, hence \(U\to\Bun\) is injective with full rank, i.e. is an injective immersion.
It is injective because \(U\to\Bun\to M\)
We can cover \(M\) in open sets \(M_a\subseteq M\) so that each is the domain of a local section \(M_a\xrightarrow{s_a}\Bun\).
Let \(\Bun_a:=\pi^{-1}M_a\), so our local sections are maps \(M_a\xrightarrow{s_a}\Bun_a\).
For simplicity of notation, assume that \(\Bun\to M\) is a right \(G\)-prebundle.
Define smooth maps
\[
(m,g)\in M_a \times G\xmapsto{t_a} s_a(m)g\in \Bun_a.
\]
We need only prove that this \(t_a\) is a diffeomorphism.
The composition 
\[
M_a \times G\xrightarrow{t_a} \Bun_a\xrightarrow{\pi} M_a
\]
is the standard trivial right principal bundle on \(M_a\).
In particular, this composition is a surjective submersion, so \(M_a\times G\xrightarrow{t_a} \Bun_a\) has full rank on each slice \(M_a\times\set{g}\).

If \(t_a(m_0,g_0)=t_a(m_1,g_1)\) then, composing with projection to \(M_a\), \(m_0=m_1\), and so 
\[
s_a(m_0)g_0=s_a(m_0)g_1.
\]
But \(G\) acts freely on the fibers of \(\Bun\to M\), so \(g_0=g_1\).
So \(t_a\) is injective.

Pick any point \(p\in \Bun_a\).
Let \(m:=\pi(p)\).
Let \(q:=s_a(m,1)\).
So then \(\pi(q)=m\).
Hence \(p,q\in \Bun_m\) are in the same fiber.
The fibers are the \(G\)-orbits.
Therefore \(p=qg\) for some \(g\in G\).
Therefore \(p=s_a(m,1)g=s_a(m,g)\).
So \(t_a\) is onto.

On the fibers, \(G\xrightarrow{t_a} \Bun_m\) is \(G\)-equivariant, i.e. just the \(G\)-action, which is free.
The fiber \(\Bun_m\) is the \(G\)-orbit, so this map is onto.
The orbits are homogeneous spaces (i.e. a countable set of components of \(G\) act transitively on some component of \(\Bun_m\)).
So for each fiber \(\Bun_m\), \((\Bun_m,G)\) is a homogeneous space with trivial stabilizer.
Let \(p_0:=t_a(m,1)\).
So \(g\in G\xmapsto{t_a}p_0g\in \Bun_m\) is a diffeomorphism.

So \(t_a\) is a diffeomorphism on each fiber, hence has full rank on fibers.
Its composition with \(\Bun_a\to M_a\) is a submersion.
So \(M_a\times G\to \Bun\) has full rank, i.e. is a local diffeomorphism.
Being a bijection, it is a diffeomorphism.
\end{answer}
A principal bundle is a \emph{left}\define{principal bundle!left}\define{bundle!principal!left} principal bundle if the \(G\)-action is left, a \emph{right}\define{principal bundle!right}\define{bundle!principal!right} principal bundle if the \(G\)-action is right.
A \emph{morphism} of principal bundles is just a prebundle morphism between principal bundles.
\begin{problem}{top.row}
Prove that any morphism of principal bundles is an isomorphism just when it is a diffeomorphism of their base manifolds.
\end{problem}
\section{Bundles as quotients}
Take a Lie group \(G\) acting on a manifold \(M\).
The action is \emph{proper}\define{action!proper}\define{proper!group action} if the map
\[
(g,m)\in G\times M\mapsto
\left.
\begin{cases}
(gm,m)&\text{ for a left action},\\
(mg,m)&\text{ for a right action}
\end{cases}
\right\}
\in M\times M
\]
is proper.
\begin{problem}{proper.group.action}
An action of a Lie group \(G\) on a manifold \(M\) is proper just when, for any convergent sequence \(x_1,x_2,\dots \to x\) in \(M\) and any sequence \(g_1,g_2,\dots\) in \(G\), if \(g_1 x_1, g_2 x_2, \dots\) converges in \(M\), then after perhaps replacing by a subsequence, \(g_1,g_2,\dots\) converges in \(G\).
\end{problem}
\begin{problem}{ae}
Every compact Lie group acting smoothly acts properly.
\end{problem}
\begin{example}
Irrational winding of a torus is a free improper action.
\end{example}
\begin{problem}{ad}
Prove that a Lie group action of a Lie group \(G\) on a manifold \(\Bun\) is free and proper just when there is a smooth structure on the quotient \(M:=\Bun/G\), so that \(\Bun\to M\) is a principal \(G\)-bundle, and this smooth structure is unique \cite{Duistermaat.Kolk:2000} p.53 theorem 1.11.4.
\end{problem}
\begin{problem}{make.bundle}
Conversely, suppose that \(G\to\Bun\to M\) is a principal bundle.
Since \(\Bun\to M\) is \(G\)-invariant, it quotients to a map \(\Bun/G\to M\).
Prove that this map is a diffeomorphism, and that the quotient map \(\Bun\to\Bun/G\) is the pullback of \(\Bun\to M\), so every principal bundle is isomorphic to a free and proper quotient.
\end{problem}
\begin{problem}{free.and.proper}
For any principal bundle \(G\to\Bun\to M\), prove that, for each point \(p_0\in\Bun\), the map \(g\in G\mapsto gp_0\in\Bun\) is a diffeomorphism to a fiber.
\end{problem}
\begin{example}
The circle acts on the unit length vectors in any complex inner product space by multiplication by unit complex numbers, hence the sphere \(S^{2n+1}\) is a principal \(S^1\)-bundle over complex projective space, the \emph{Hopf fibration}.\define{Hopf fibration}
\end{example}
{
\NewDocumentCommand\pG{}{{\vphantom{G}}'\!G}
\NewDocumentCommand\pg{}{{\vphantom{g}}'\!g}
\NewDocumentCommand\pM{}{{\vphantom{M}}'\!M}
\NewDocumentCommand\pmm{}{{\vphantom{m}}'\!m}
\NewDocumentCommand\pn{}{{\vphantom{n}}'\!n}
\begin{problem}{biquot}
Suppose that \(\pG,G'\) are Lie groups acting on the left and the right on a manifold \(M\), and commuting.
Let \(\pG\times G'\) act on \(M\) by 
\[
(g,h)m:=gmh^{-1}.
\]
Let
\begin{align*}
\pM&:=\pG\backslash M,\\
M'&:=M/G',\\
\pM'&:=(\pG\times G')\backslash M
\end{align*}
as topological quotient spaces.
Write down homeomorphisms between
\[
\pM', (\pG\backslash M)/G',\text{ and } \pG\backslash(M/G').
\]
Prove that the following are equivalent:
\begin{itemize}
\item
\(\pG\times G'\) acts freely and properly on \(M\)
\item
\(\pG\) acts freely and properly on \(M\) and  \(G'\) acts freely and properly on \(\pM\)
\item
\(G'\) acts freely and properly on \(M\) and  \(\pG\) acts freely and properly on \(M'\).
\end{itemize}
Prove that if any, hence all, of these happen then then all of the quotient maps
\[
\begin{tikzcd}
& M\arrow[dl]\arrow[dr]\arrow[dd]& \\
\pM\arrow[dr]& {} & M'\arrow[dl] \\
& \pM' &
\end{tikzcd}
\]
are smooth principal bundle maps, and \(\pG\) acts on \(M\to M'\) by smooth principal bundle automorphisms, and \(G'\) acts on \(M\to\pM\) by smooth principal bundle automorphisms.
\end{problem}
}
\begin{answer}{biquot}
\NewDocumentCommand\pG{}{{\vphantom{G}}'\!G}
\NewDocumentCommand\pg{}{{\vphantom{g}}'\!g}
\NewDocumentCommand\pM{}{{\vphantom{M}}'\!M}
\NewDocumentCommand\pmm{}{{\vphantom{m}}'\!m}
\NewDocumentCommand\pn{}{{\vphantom{n}}'\!n}
Suppose that \(\pG\) acts freely and properly on \(M\) and that \(G'\) acts on \(\pM\) freely and properly.
So \(\pG\to M\to\pM\) is a right principal bundle and \(G'\) acts on \(M\) by smooth principal bundle automorphisms, equivariant for \(M\to\pM\).
So \(G'\) acts on \(\pM\) smoothly.
If \(gmh=m\) for some \(m\in M\), \(g\in\pG\), \(h\in G'\) then apply \(M\to\pM\) to see that \(h=1\), so \(G'\) acts freely on \(M\).
If \(m_i\to m\) and \(g_im_ih_i\to n\) then apply \(M\to\pM\) to see that some subsequence of of \(h_i\) converges, say to \(h\).
Replace \(m_i\) by \(m_ih^{-1}\), \(m\) \(mh^{-1}\) and \(n\) by \(nh^{-1}\) to arrange \(m_i\to m\), \(g_im_i\to n\).
So \(g_i\) has a convergent subsequence.
So \(\pG\times G'\) acts freely and properly on \(M\).

By obvious change of notation, the same steps apply with the left and right actions swapped.

Suppose that \(\pG\times G'\) acts freely and properly on \(M\).
So \(M\to\pM'\) is a principal \(\pG\times G'\)-bundle.
Being a subgroup, \(\pG\) acts freely and properly on \(M\) too.
Pick an open set \(U\subset\pM'\) containing that point, over which \(M\to\pM'\) is trivial, say identified with \(U\times\pG\times G'\).
Then the preimage in \(\pM\) is identified with the quotient \(U\times G'\), and the preimage in \(M'\) is identified with the quotient \(U\times\pG\).
So the action of \(G'\) on \(\pM\) is free and proper with quotient \(\pM'\).
\end{answer}
\section{Constructions}
\subsection{Homogeneous spaces}
By problem~\vref{problem:n}, for every homogeneous space \((X,G)\), if we pick a point \(x_0\in X\) and let \(H:=G^{x_0}\) be the stabilizer of that point, then the map 
\[
g\in G\mapsto gx_0\in X
\]
is a principal right \(H\)-bundle:
\[
\begin{tikzcd}
G\ar[d]&H\arrow{l}\\
X.
\end{tikzcd}
\]
For example, working over the real or complex numbers or the quaternions, if we take \(X=\Proj{n}\) and \(G=\PGL{n}\) then \(H\) is the group of projection transformations fixing a point, say
\[
x_0=[1:0:\dots:0],
\]
so \(H\) is the group of matrices, up to scaling by nonzero constants, of the form:
\[
h=
\begin{bmatrix}
a&b\\
0&d
\end{bmatrix}
\]
with block sizes
\[
\begin{bNiceMatrix}[margin,first-row,first-col]
  & 1 & n \\
1 & \ & \ \\
n & \ & \ \\
\end{bNiceMatrix}
\]
\subsection{Pullback}
For any smooth map of manifolds \(M'\xrightarrow{\varphi}M\), and prebundle \(\Bun\xrightarrow{\pi}M\), the \emph{pullback} prebundle has total space
\[
\Bun'=\varphi^*\Bun:=\set{(m,p)|m\in M', p\in\Bun, \varphi(m)=\pi(p)},
\]
and map \((m,p)\mapsto m\).
The pullback has \(G\)-action \((m,p)g:=(m,pg)\) if \(\Bun\to M\) is a right \(G\)-prebundle, and \(g(m,p):=(m,gp)\) if \(\Bun\to M\) is a left \(G\)-prebundle.
The pullback has morphism
\[
\begin{tikzcd}
\Bun\arrow[r, "\Phi"]\arrow[d,"\pi"']&\Bun'\arrow[d, "\pi'"] \\M\arrow[r, "\varphi"']& M'
\end{tikzcd}
\]
where \(\Phi(m,p)=p\).
\begin{problem}{pullb}
Prove that the pullback of any principal bundle is a principal bundle.
\end{problem}
\subsection{Whitney sum}
If \(G\to\Bun\xrightarrow{\pi}M\) and \(G'\to\Bun'\xrightarrow{\pi}M\) are principal bundles, their \emph{Whitney sum}\define{Whitney sum} \(G\times G'\to\Bun\times_M\Bun'\to M\) has total space \(\Bun\times_M\Bun'\) the set of pairs \((x,x')\in\Bun\times\Bun'\) so that \(x\) and \(x'\) project to the same point \(m\in M\), i.e. \(m=\pi(x)=\pi'(x')\), with obvious group action
\[
(g,g')(x,x')=(gx,g'x'),
\]
for left principal bundles, and
\[
(x,x')(g,g')=(xg,x'g'),
\]
for right principal bundle, and so on.
\begin{problem}{Whitney.sum}
Prove that the Whitney sum of any two principal bundles is a principal bundle.
\end{problem}
\subsection{Induced bundles}
Suppose that \(G\to\Bun\to M\) and that \(G\xrightarrow{\varphi}G'\) is a Lie group morphism.
Consider the group action of \(G\) on \(\Bun\times G'\):
\begin{align*}
g(x,g')&:=(gx,g'g^{-1}), \text{if \(\Bun\) is left principal bundle,}\\
(x,g')g&:=(xg,g^{-1}g'), \text{if \(\Bun\) is a right principal bundle.}
\end{align*}
Consider the group action of \(G'\) on \(\Bun\times G'\):
\begin{align*}
k'(x,g')&:=(gx,k'g'), \text{if \(\Bun\) is left principal bundle,}\\
(x,g')k'&:=(xg,g'k'), \text{if \(\Bun\) is a right principal bundle.}
\end{align*}
\begin{problem}{principal.bundles:commuting.actions}
Prove that these commute.
\end{problem}
\begin{problem}{principal.bundles:commuting.actions.2}
Prove that the \(G'\)-action on \(\Bun\times G'\) commutes with the \(G\)-action.
\end{problem}
Let \(\Bun':=\amal{\Bun}{G}{G'}\) be the quotient of \(\Bun\times G'\) by the \(G\)-action.
\begin{problem}{principal.bundles:commuting.actions.3}
Prove that the \(G'\)-action descends to the quotient, giving a principal \(G'\)-bundle. 
\end{problem}
Roughly speaking, principal bundle morphisms are induced:
\begin{problem}{principal.bundles:commuting.actions.4}
Suppose that \(G\to\Bun\to M\) and \(G'\to\Bun'\to M\) are principal bundles over the same manifold \(M\).
Fix a Lie group morphism \(G\to G'\).
Suppose that \(\Bun\to\Bun'\) is a smooth map, equivariant for \(G\to G'\), taking the fiber over any point to the fiber over the same point.
Prove that the map extends to a unique principal bundle isomorphism \(\amal{\Bun}{G}{G'}\to\Bun'\).
Conversely, prove that if we define \(\Bun':=\amal{\Bun}{G}{G'}\), then the composition of the obvious maps
\[
\Bun\to\Bun\times G'\to\amal{\Bun}{G}{G'}=\Bun'
\]
is smooth, equivariant for \(G\to G'\), and takes the fiber over any point to the fiber over the same point.
\end{problem}
\section{Associated vector bundles}
Suppose that \(G\to\Bun\xrightarrow{\pi}M\) is a principal right \(G\)-bundle and that \(V\) is a \(G\)-module, i.e. \(V\) is a finite dimensional real (or complex or quaternionic) vector space equipped with a Lie group morphism \(G\xrightarrow{\rho}\GL{V}\).
The \emph{associated vector bundle}\define{associated!vector bundle}\define{vector!bundle!associated} is \(\vb{V}:=\amal{\Bun}{G}{V}\).
\begin{problem}{a.vb.fiber}
Why is the associated vector bundle a vector bundle?
\end{problem}
\begin{answer}{a.vb.fiber}
The map \((p,v)\in\Bun\times V\mapsto m:=\pi(p)\) is \(G\)-invariant, hence descends to a map \(\vb{V}\to M\); denote the fiber of this map over a point \(m\in M\) by \(\vb{V}_m\).
Each point \(w\in\vb{V}_m\) is the equivalence class \((p,v)G\) of a pair \((p,v)\in\Bun\times V\).
Each fiber is the quotient \(\amal{\Bun_m}{G}{V}\).
Picking any point \(p\in\Bun_m\), we map 
\[
v\in V\mapsto (p,v)G\in\amal{\Bun_m}{G}{V}.
\]
If we change the choice of \(p\), say to \(pg\), we change this to
\[
(pg,v)G=(p,g^{-1}v)G,
\]
a linear isomorphism on \(V\), so the fibers are vector spaces invariantly, although the identification with \(V\) depends on the choice of point \(p\).
The act of \(G\) on \(\Bun\) is free and proper, so the action on \(\Bun\times V\) is also free and proper, so \(\Bun\times V\to\amal{\Bun}{G}{V}\) is a principal right \(G\)-bundle and so \(\vb{V}\) is a smooth manifold.
The map \(\Bun\to V\to\Bun\) is smooth, so the quotient by \(G\) is the smooth map \(\amal{\Bun}{G}{V}\to M\).
So \(\amal{\Bun}{G}{V}\) is a smooth manifold, smoothly mapped to \(M\), with vector space fibers.
We need to check local triviality.
Replace \(M\) by an open subset on which \(G\to\Bun\to M\) is trivial.
It suffices to prove the triviality of the associated vector bundle of the trivial bundle \(\Bun=M\times G\), for which \(\amal{M\times G}{G}{V}=M\times V\) by taking inside each orbit \((m,g,v)G\) the point \((m,1,v)\).
\end{answer}
\begin{problem}{moebius}
How is the M\"obius strip the associated line bundle of a principal bundle?
\end{problem}
\[
\includegraphics[width=5cm]{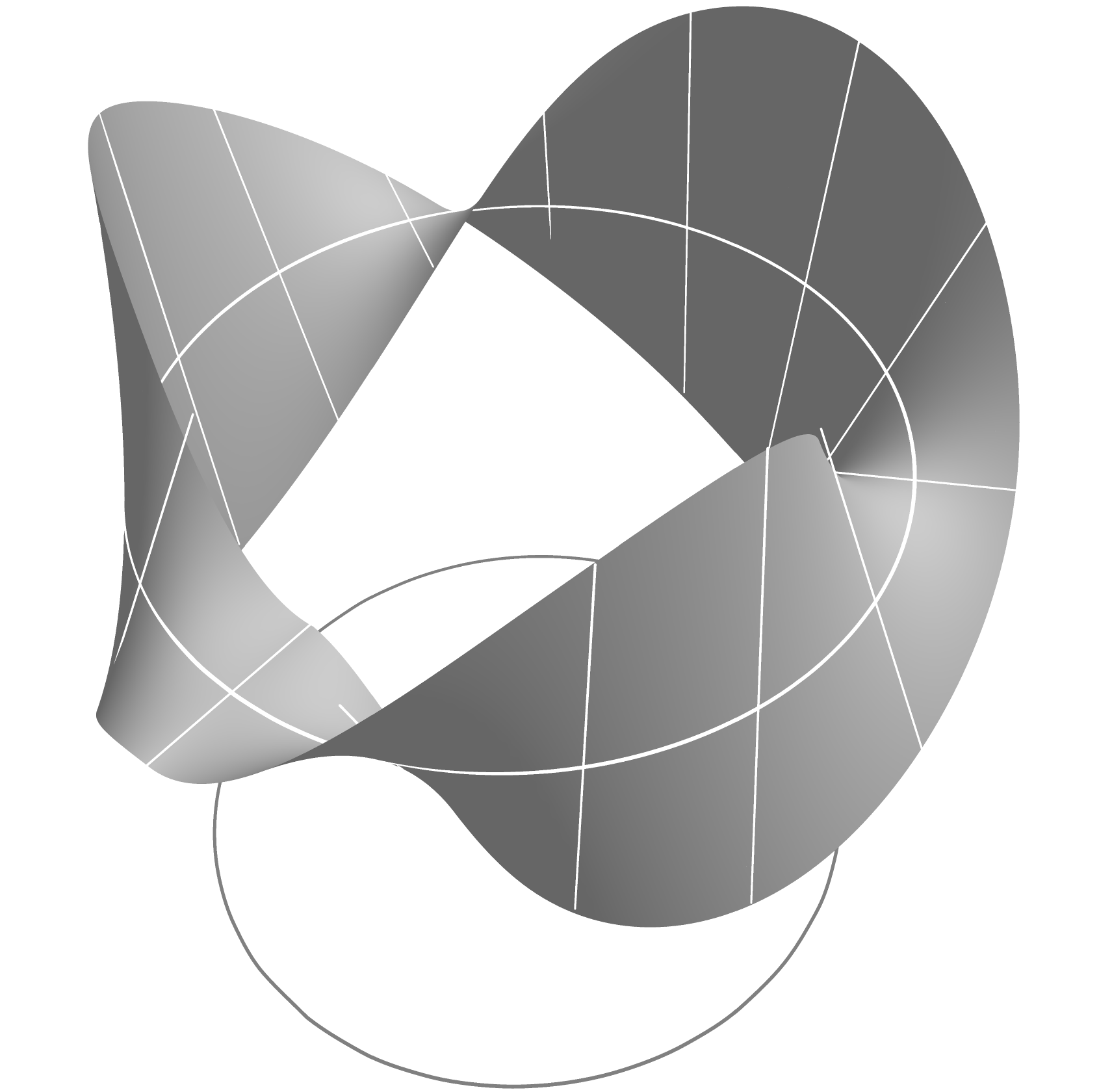}
\]
Given any smooth \(G\)-equivariant function \(\Bun\xrightarrow{f}V\), its graph is \(G\)-invariant, and drops to a section of \(\vb{V}\to M\); conversely all smooth sections arise this way.
\begin{example}
Take \((X,G)=(\Proj{n},\SL{n})\) over the field \(k=\R\) or \(\C\) of the real or complex numbers.
So \(H\) is the set of matrices of the form
\[
h=
\begin{pmatrix}
a&b\\
0&d
\end{pmatrix}
\]
with block sizes
\[
\begin{bNiceMatrix}[margin,first-row,first-col]
  & 1 & n \\
1 & \ & \ \\
n & \ & \ \\
\end{bNiceMatrix}
\]
and with determinant \(1\), i.e. \(a\det d=1\).
Since \(H\subseteq G\) is precisely the stabilizer of the line spanned by
\[
e_0=(1,0,\dots,0),
\]
\(X=G/H\) is the set of lines through the origin in \(k^{n+1}\).
Let \(\OO(d)\) be the vector bundle associated to the \(H\)-module \(V=k\) with the \(H\)-representation
\[
\rho(h)t=t/a^d
\]
for \(t\in V\).
Given any polynomial \(p(x)\) homogeneous of degree \(d\) in variables, construct a function \(G\xrightarrow{f}k\) by \(f(g)=p(ge_0)\).
Then 
\[
f(gh)=p(ghe_0)=p(gae_0)=a^dp(ge_0)=a^df(g)=\rho(h)^{-1}f(g).
\]
So \(f\) is a section of \(\OO(d)\).
\end{example}
\begin{example}
Consider a more common description of \(\OO(d)\).
Write each point of \(X\) as \(L\in X\), so \(L\subseteq k^{n+1}\) is a line through the origin.
Let \(\mathscr{L}\) be the set of pairs \((L,q)\) where \(L\in X\) and 
\begin{itemize}
\item
if \(d>0\), \(q\in (L^*)^{\otimes d}\) is a degree \(d\) polynomial on \(L\),
\item
if \(d=0\), \(q\in k\),
\item
if \(d<0\), \(q\in L^{\otimes |n|}\).
\end{itemize}
The map
\[
(L,q)\in \mathscr{L}\mapsto L\in X
\]
has fibers 
\[
\mathscr{L}_L:=
\begin{cases}
(L^*)^{\otimes d},&d>0,\\
k,&d=0,\\
L^{\otimes |d|},&d<0.
\end{cases}
\]
Clearly by taking the usual affine chart on projective space \(X=\Proj{n}\), we can then make the obvious section of \(\mathscr{L}\to X\), so it is a line bundle.
The group \(G\) acts on \(\mathscr{L}\) by
\[
g(L,q)=(gL,g_* q)
\]
where
\[
g_*q=
\begin{cases}
q\circ g^{-1},&d>0,\\
q,&d=0,\\
gq,&d<0.
\end{cases}
\]
Let 
\[
q_0:=
\begin{cases}
x\mapsto (x^0)^d,&d>0,\\
1,&d=0,\\
e_0^{\otimes d},&d<0.
\end{cases}.
\]
Map 
\[
(g,c)\in G\times k\xrightarrow{\varphi} (gL_0,cg_*q_0)\in \mathscr{L}.
\]
The fibers consists of the \(H\)-orbits, and \(H\) acts freely, so the map descends to a \(G\)-equivariant diffeomorphism
\[
\OO(d)=\amal{G}{H}{k}\cong \mathscr{L},
\]
which is linear on the fibers.
\end{example}
\begin{problem}{v.b.fiber.pres}
Suppose that \(V,W\to M\) are vector bundles and that \(V\xrightarrow{\varphi}W\) is a smooth map taking the fiber \(V_m\) over each point \(m\in M\) to the fiber \(W_m\), a linear isomorphism on each fiber.
Prove that \(\varphi\) is a vector bundle isomorphism.
\end{problem}
\begin{answer}{v.b.fiber.pres}
For trivial bundles, it is clear.
All vector bundles are locally trivial.
\end{answer}
Given any vector bundle \(\vb{V}\to M\), and a vector space \(V\) of dimension equal to the rank of \(\vb{V}\), the \(V\)-valued \emph{frame bundle}\define{frame bundle} \(\framebundle{\vb{V}}\to M\) of \(\vb{V}\) is the set of all pairs \((m,u)\) with \(m\in M\) and \(\vb{V}_m\xrightarrow{u}V\) a linear isomorphism.
\begin{problem}{frame.bundle.GLV}
Prove that the frame bundle is a principal right \(\GL{V}\)-bundle under the action
\[
(m,u)g:=(m,g^{-1}u),
\]
for \((m,u)\in\framebundle{\vb{V}}\) and \(g\in\GL{V}\).
\end{problem}
\begin{answer}{frame.bundle.GLV}
For a trivial bundle, it is clear.
If you glue together just two trivial bundles, you already see how the story unfolds: the transition maps of the vector bundle is that of the frame bundle as a principal bundle.
\end{answer}
Every vector bundle isomorphism yields an isomorphism of frame bundles as right principal bundles.
Every right principal bundle isomorphism yields an isomorphism of all associated vector bundles.
Consider the map
\[
((m,u),v)\in\framebundle{\vb{V}}\times V\mapsto u^{-1}v\in\vb{V}.
\]
This map is invariant under the \(\GL{V}\)-action
\[
((m,u),v)g:=((m,g^{-1}u),g^{-1}v).
\]
so descends to the quotient
\[
\amal{\framebundle{\vb{V}}}{\GL{V}}{V}\to\vb{V}.
\]
Since the frame bundle is a principal bundle, so locally trivial, we can look in a local trivialization to see that this map is a fiber preserving diffeomorphism.
Giving each fiber the obvious vector space structure, it is a vector bundle isomorphism.
The inverse map takes each section \(s\) of \(\vb{V}\) to graph of the \(G\)-equivariant function
\[
(m,u)\in\framebundle{\vb{V}}\xmapsto{f}u(s(m))\in V.
\]
\begin{problem}{equiv.fn}
Every \(G\)-equivariant function
\(
\framebundle{\vb{V}}\xrightarrow{f}V
\)
arises uniquely from a section of \(\vb{V}\) in this way.
\end{problem}
\begin{problem}{quotient.v.b}
Suppose that \(V\to M\) is a vector bundle, and \(G\) is a Lie group acting freely and properly on the right on \(M\), and acting on \(V\) as vector bundle isomorphisms over \(M\), commuting with the bundle map \(V\to M\).
Prove that the quotients \(\bar{M}:=V/G\) and \(\bar{V}:=V/G\) are smooth manifolds and the quotient maps \(M\to\bar{M}\), \(V\to\bar{V}\) are smooth principal right \(G\)-bundles and that \(\bar{V}\to\bar{M}\) has a unique structure of smooth vector bundle so that, on each fiber \(V_m\to\bar{V}_{\bar{m}}\) is a linear isomorphism.
\end{problem}
\begin{answer}{quotient.v.b}
Apply problem~\vref{problem:biquot} to the frame bundle of \(\vb{V}\) and the \(G\)-action on it to see that \(G\) acts freely any properly on the frame bundle, and that a quotient bundle \(\GL{V}\to\bar\Bun\to\bar{M}\) exists over the quotient manifold \(\bar{M}=G\backslash M\).
Let \(\bar{V}\) be the associated vector bundle \(\bar{V}:=\amal{\bar\Bun}{\GL{V}}{V}\).
By construction, \(\bar{V}=\bar\Bun/(G\times\GL{V})\), so again by problem~\vref{problem:biquot}, \(\bar{V}=G\backslash V\).
\end{answer}

\chapter{Basic theory of ordinary differential equations}
\section{Time varying vector fields}
A \emph{time varying vector field}\define{time varying vector field} \(X\) on a manifold \(M\) associates to each point \(m\in M\) and time \(t\) in some interval \(I\subseteq\R\) a vector \(X_t(m)\in T_m M\), locally integrable in \(t\), and locally Lipschitz in \(m\) with a local Lipschitz constant (in some coordinates) bounded over time by a locally integrable function of time.
\section{Weak solutions}
A locally absolutely continuous curve \(x(t)\), perhaps not differentiable, is a \emph{flow line} of \(X_t(x)\), also called a \emph{solution} of the associated equation \(\dot{x}=X_t\), if, for each time \(t=t_0\), in some local coordinates with  \(x(t_0)\) in the domain of those coordinates, and for any \(t\) close enough to \(t_0\):
\[
x(t)=x(t_0)+\int_{t_0}^t X_t(x(t))\,dt.
\]
\begin{problem}{ode.coords}
Justify replacing the expression \emph{some local coordinates} above by \emph{any local coordinates}.
\end{problem}
\begin{answer}{ode.coords}
If \(X_t(x)\) is absolutely continuous in \(t\) we can differentiate to find
\[
\dot{x}=X_t(x),
\]
and the usual change of variables for vector fields applies to prove coordinate independence.
If \(X_t(x)\) is merely locally integrable in \(t\), we can take a sequence of time varying vector fields \(X_{1t},X_{2t},\dots\) which approach \(X_t\) locally in some coordinates, \(L^1\) locally in \(t\), and with local Lipschitz bound approach that of \(X_t\).
For example, we can mollify with a smooth bump function, in some local coordinates.
Apply existence and uniqueness to these, before and after any change of coordinates.
\end{answer}
\section{Existence and uniqueness}
\begin{theorem}[Picard's existence and uniqueness theorem \cite{Sontag:1998} C.3.8 p. 482]\label{thm:Picard}
There is a local solution which passes through a given point of the manifold at a given time, unique if we ask that it be defined on a maximal time interval.
If the time varying vector field is \(C^k\), some \(k=1,2,\dots\) or \(C^{\infty}\) or \(C^{\omega}\) then its flows lines are \(C^{k+1}\), \(C^{\infty}\) or \(C^{\omega}\).
\end{theorem}
\section{Continuity}
Suppose that \(M\) is a manifold with a Riemannian metric \(g\).
Denote the \(g\)-distance between points as \(d\).
Suppose that \(p,q\in M\) and that \(q\) lies at a distance from \(p\) less than the injectivity radius at \(p\).
There is a unique unit speed minimal geodesic starting at \(p\) at time zero and travelling to \(q\).
Denote parallel transport along that path by \((p\to q)\).
\begin{theorem}\label{thm:continuity}
Take a time varying vector field \(X\) on a Riemannian manifold \(M\).
Suppose that \(I\xrightarrow{x}M\) solves \(\dot x=X_t(x)\) on a compact interval \(I\subset\R\).
Pick some \(t_0\in I\).
Let \(x_0:=x(t_0)\).
Take some positive number \(R>0\) so that at every point of \(x(t)\), \(R\) is less than the injectivity radius of \(M\) at \(x(t)\).

Take a time varying vector field \(Y\) on \(M\).
There is an integrable function \(\beta(t)>0\) for \(t\in I\)so that, for every time \(t\),
\[
\beta(t)\ge\frac{|(y\to x(t))(X+Y)_t(y)-(X+Y)_t(x(t))|}{d(x(t),y)}
\]
for all \(y\in M\) with \(d(x(t),y)\le R\).
Pick any such function.

Pick a point \(y_0\in M\).
Suppose that
\[
d(x_0,y_0)+\int_I|Y_s(x(s))|\,ds
\le
Re^{-\int_I \beta(s)\,ds}.
\]

Let \(y(t)\) be the solution of
\[
\dot y=X_t(y)+Y_t(y)
\]
with \(y(t_0)=y_0\), on the maximal time interval on which this is defined.
Then \(y(t)\) is defined for all \(t\in I\) and 
\[
d(x(t),y(t))
\le
\left(d(x_0,y_0)+\int_{t_0}^t|Y_s(x(s))|\,ds\right)
\exp\int_{t_0}^t \beta(s)\,ds.
\]
\end{theorem}
\begin{proof}
The function \(\beta\) exists, because our definition of time varying vector field forces \(\beta\) to exist locally in charts, in the Euclidean metric, and we only need finitely many charts to cover the path \(x(t)\), \(t\in I\), and then some constants to relate the Euclidean metric to \(g\).

Suppose that \(y(t)\) is defined and stays at distance less than \(R\) from \(x(t)\) on some interval \(J\subseteq I\) containing \(t_0\).
Let \(U\subseteq J\) be the set of times \(t\) at which \(x(t)\ne y(t)\).
At every time \(t\in J-U\), our result holds trivially.
Each point \(t\in U\) lies inside an interval \(J_t\subseteq J\) with each end either at an end of \(J\) or at a point outside \(U\).
It suffices to prove the result on each interval \(J_t\).
So we can assume that \(x(t)\ne y(t)\) except perhaps at the end points of \(J\).

Denote distance from a point \(p\) by \(d_p(q):=d(p,q)\).
For \(q\) inside the injectivity radius of \(p\), but with \(q\ne p\), \(\nabla d_p(q)\) is a unit vector, and is the velocity of the unique minimal geodesic from \(p\) to \(q\).

We have a unique unit speed minimal geodesic \(x(t)\to y(t)\), which has an initial unit velocity \(u(t)=-d_q(x(t))\) at \(x(t)\), and a final unit velocity \(v(t)=d_p(y(t))\) at \(y(t)\).
Note that \((x(t)\to y(t))u(t)=v(t)\).
Since \(x(t),y(t)\) are absolutely continuous, and the exponential map is smooth, \(u(t),v(t)\) are absolutely continuous.
\begin{align*}
d(x(t),y(t))-d(x_0,y_0)
&=
\int_{t_0}^t
\frac{d}{ds}d(x(s),y(s))\,ds,
\\
&=
\int_{t_0}^t
\frac{d}{ds}d(x(s),y(s))\,ds,
\\
&=
\int_{t_0}^t
\left<\dot x(s),(-u)\right>\,ds
+
\int_{t_0}^t
\left<\dot y(s),v\right>
\,ds,
\\
&=
\int_{t_0}^t
-\left<X_s(x),u\right>\,ds
+
\int_{t_0}^t
\left<(X+Y)_s(y),v\right>
\,ds,
\\
&=
\int_{t_0}^t
-\left<X_s(x),u\right>\,ds
+
\int_{t_0}^t
\left<(X+Y)_s(y),(x\to y)u\right>
ds,
\\
&=
\int_{t_0}^t
-\left<X_s(x)+Y_s(x),u\right>\,ds
+
\int_{t_0}^t
\left<Y_s(x),u\right>
\,ds,
+
\int_{t_0}^t
\left<(y\to x)(X+Y)_s(y),u\right>
ds
\\
&=
\int_{t_0}^t
\left<(y\to x)(X+Y)_s(y)-(X+Y)_s(x),u\right>\,ds
+
\int_{t_0}^t
\left<Y_s(x(s)),u\right>
\,ds,
\\
&\le
\int_{t_0}^t
\left|(y\to x)(X+Y)_s(y)-(X+Y)_s(x)\right|\,ds
+
\int_{t_0}^t
\left|Y_s(x(s))\right|
\,ds,
\\
&\le
\int_{t_0}^t
\beta(s)d(x(s),y(s))\,ds
+
\int_{t_0}^t
\left|Y_s(y(s))\right|
\,ds.
\end{align*}
By Gronwall's inequality \cite{Sontag:1998} p.~475,
\[
d(x(t),y(t))
\le
\left(d(x_0,y_0)+\int_{t_0}^t|Y_s(x(s))|\,ds\right)
\exp\int_{t_0}^t \beta(s)\,ds.
\]
By our hypotheses, the path \(y(t)\) remains in the injectivity radius of \(x(t)\) on the closure of \(J\).
Therefore the path \(y(t)\) exists and satisfies the same estimates on a longer time interval than \(J\) unless \(J=I\).
\end{proof}
\begin{theorem}\label{thm:lift.continuity}
Suppose that \(M\) is a manifold, \(V\subseteq TM\) a vector subbundle.
Suppose that \(M\xrightarrow{\varphi}\bar{M}\) is a smooth submersion and that \(TM=V\oplus\ker\varphi'\).
Take a sequence of absolutely continuous paths \(x_1,x_2,\dots\) in \(\bar{M}\), all defined on the same closed interval \(I\subseteq\R\), converging uniformly to an absolutely continuous path \(x\) in \(\bar{M}\), and with velocities \(\dot x_1,\dot x_2,\dots\) bounded in \(L^1\) norm.
Take points \(m_1,m_2,\dots\to m\in M\) so that, if we let \(\bar m_j:=\varphi(m_j)\), then
\[
\bar m_j=x_j(t_0), 
\]
and \(\bar m=x(t_0)\).
Suppose that there is a \(V\)-tangent absolutely continuous path \(p\) on \(M\), defined on \(I\), lifting \(x\), with \(p(t_0)=m\).

Then, for all but finitely many \(x_j\), there is a unique \(V\)-tangent absolutely continuous path \(p_j\) on \(M\), defined on \(I\), lifting \(x_j\), with \(p_j(t_0)=m_j\).
Moreover, \(p_j\) converges uniformly to \(p\).
The \(\dot p_j\) converge in \(L^1\) to \(\dot p\) if and only if the \(\dot x_j\)  converge in \(L^1\) to \(\dot x\).
\end{theorem}
\begin{proof}
To be precise, the bound in \(L^1\) norm can be computed in any Riemannian metric, defined perhaps only in some small neighborhood of \(x\).

By the implicit  function theorem, we can cover \(M\) and \(\bar{M}\) in domains of coordinate charts in which \(M\xrightarrow{\varphi}\bar{M}\) is a linear projection map.
We can cover the curve \(p\) in finitely many of these. 
We will build \(V\)-tangent lifts \(p_j\) for all of the \(x_j\) that stay inside these chart domains in \(\bar{M}\), which is all \(x_j\) but finitely many.
So we can assume that \(M\) is the union of finitely many such chart domains.
Write \(I\) as the union of finitely many closed subintervals, on each of which \(p\) stays in such a chart domain.
We only need to prove the result for each subinterval, and glue, by existence and uniqueness.
So we can work entirely in a single chart domain. 

So we can assume that \(M\) is an open set in Euclidean space with coordinates \((x,y)\), for \(x,y\) points in Euclidean spaces of some dimensions, and that \(\varphi(x,y)=x\).
Hence \(V\) is uniquely expressed as \(dy=f(x,y)dx\), for some smooth function \(f\).

We can think of \(V\) as a differential equation
\[
\dot y=f(x(t),y)\dot x
\]
for an unknown \(y\), but a known \(x\).
Local solutions exist by Picard's theorem.
By hypothesis, there is a global solution \(y(t)\) corresponding to the given \(x(t)\), \(t\in I\), with \(y(t_0)\) given.
Our aim is to prove the existence of global solutions \(y_i(t)\) of all but finitely many equations
\[
\dot y_i=f(x_i(t),y_i(t))\dot x_i
\]
with given \(y_i(t_0)\) near \(y(t_0)\).

To apply Picard's theorem, we need to control the local Lipschitz constant for \(f(x,y)\).
Since \(f(x,y)\) is smooth, we can acheive this control by perhaps replacing our charts by finitely many charts on precompact domains in the domains of our previously given charts.

We also need to control
\[
\bar{y}_i(t):=
\int_{t_0}^t
\left\{
f(x_i(s),y(s))\dot{x}_i-
f(x(s),y(s))\dot{x}
\right\}ds,
\]
formed from integrals computed along the given solution \(y(t)\).
To be precise, we need to make sure that
\[
\sup_{t\in I}|\bar{y}_i(t)|
\]
gets sufficiently small for large enough \(i\) \cite{Sontag:1998} p.~486, Theorem~55.

Compute
\begin{align*}
\bar{y}_i(t)
&=
\int_{t_0}^t
\left\{
f(x_i(s),y(s))\dot{x}_i-
f(x(s),y(s))\dot{x}
\right\}ds,
\\
&=
\int_{t_0}^t
f(x_i(s),y(s))\left(\dot{x}_i-\dot{x}\right)
ds\\
&\qquad+
\int_{t_0}^t
\left\{
f(x_i(s),y(s))-
f(x(s),y(s))
\right\}\dot{x}ds,
\\
&=
\left.
f(x_i(s),y(s))(x_i-x)
\right|_{s=t_0}^{s=t}
-
\int_{t_0}^t
\left\{f_x\dot x_i+f_y \dot y\right\}(x_i-x)
ds\\
&\qquad+
\int_{t_0}^t
\left\{
f(x_i(s),y(s)-
f(x(s),y(s))\right\}\dot{x}
ds.
\end{align*}
The first term goes to zero by uniform convergence \(x_i\to x\), the second by \(L^1\) boundedness of \(x_i\), and uniform convergence \(x_i\to x\), the third by uniform converge \(x_i\to x\).
Therefore the solutions \(y_i(t)\) exist for \(t\in I\) and approach \(y(t)\) uniformly \cite{Sontag:1998} p.~486 Theorem~55.
\end{proof}
\section{Reparameterization}
Take a manifold \(M\).
A curve in \(M\) is \emph{absolutely continuous} or \emph{locally Lipschitz} if, in any chart, on every interval on which the curve is defined and stays in the chart domain, it is absolutely continuous or locally Lipschitz in the chart.
We topologize the continuous paths by the compact open topology.
Cover \(M\) by a countable, locally finite, collection of precompact open sets, each of which has closure lying in the domain of a chart.
Take an interval \(I\subseteq\R\), not necessarily open, closed or bounded.
Take a continuous curve in \(I\xrightarrow{x}M\).
Cover \(I\) in a countable set of closed intervals, on each of which \(x\) remains in one of our precompact open sets.
We topologize the absolutely continuous maps by refining the compact open topology using the Frech\'et topology on the locally integrable functions applied to the velocity of the path, as measured in each of our charts, on these intervals \cite{Meise.Vogt:1997} p.~40.
We topologize the locally Lipschitz maps by refining the compact open topology using the Frech\'et topology on the bounded measureable functions applied to the velocity of the path, as measured in each of our charts, on these intervals.
\begin{theorem}\label{thm:reparam}
Take a manifold \(M\) with a Riemannian metric \(g\).
The length of any absolutely continuous path \(I\to M\) is defined and invariant under composition \(J\to I\to M\) with any absolutely continuous surjective increasing or decreasing (perhaps not strictly) map \(J\to I\).
For each point of \(I\), there is a canonical continuous choice of reparameterization making the path unit speed.
\end{theorem}
\begin{proof}
The change of variables formula
\[
\int \left|\dot p(t)\right|_g \, dt=
\int \left|\dot p(\rho(t))\right|_g \, |\dot \rho(t)|dt
\]
is valid for any absolutely continuous increasing (perhaps not strictly increasing) surjective map of intervals \cite{Rudin:1987} p.~156.

By definition,
\[
\dot p\in L^1_{\text{loc}}(I,p^*TM).
\]
Let 
\[
\rho(t):=\int_{t_0}^t|\dot p(s)|\,ds,
\]
an absolutely continuous increasing (perhaps not strictly) function
\[
I\xrightarrow{\rho}\R.
\]
Let \(J:=\rho(I)\subseteq\R\) be its image.
By continuity of \(\rho\), \(J\) is connected, so an interval.
If \(\rho(t_1)=\rho(t_2)\) for some \(t_1\le t_2\) then
\[
0=\int_{t_1}^{t_2}|\dot p(s)|\,ds
\]
so \(p(t)\) is constant on \(t_1\le t\le t_2\).
Therefore the map
\[
P(s):=p(t)
\]
for \(s=\rho(t)\) is defined, for \(s\in J\),
\[
J\xrightarrow{P}M.
\]
Since \(P\) is a reparameterization of \(p\), every piece of \(P\) has the same length as any piece of \(p\) which we can reparameterize to it, as above, so finite.
So \(P\) is also absolutely continuous.
Take an interval of \(P\),say \(P(s)\) on \(s_1\le s\le s_2\).
Take the largest \(t_1\) and smallest \(t_2\) mapped by \(\rho\) to \(s_1,s_2\) respectively.
The distance between the end points is
\begin{align*}
d(P(s_1),P(s_2))
&=
d(p(t_1),p(t_2)),
\\
&=
\int_{t_1}^{t_2}|\dot p(t)|_g\,dt,
\\
&=
\int_{t_1}^{t_2}\rho(t)\,dt,
\\
&=
s_2-s_1.
\end{align*}
So \(P\) has unit speed, and hence is Lipschitz.

Take an open set, given by uniform local open estimates on velocity, and open sets controlling position.
If \(p\) maps to some \(P\) satisfying these estimates, then every curve \(q\) near \(p\), in the topology on absolutely continuous curves, stays near \(p\), so has image \(Q\) staying near \(P\), and has almost the same reparameterization since 
\[
\int ||\dot p|-|\dot q||\le\int |\dot p-\dot q|
\]
is small on any of our intervals.
So \(p\mapsto P\) is continuous.
\end{proof}

\chapter{Lie equations}\label{appendix:Lie.equations}
This appendix doesn't discuss Cartan geometries, but builds some essential tools.
Each curve \(g(t)\in G\) in a Lie group \(G\) has \emph{Darboux derivative}\define{Darboux derivative} or \emph{left logarithmic derivative}\define{left logarithmic derivative} \(A(t):=\dot{g}(t)\hook\omega_G\).
Conversely, take a curve \(A(t)\in\LieG\) in the Lie algebra of a Lie group \(G\).
The \emph{Lie equation}\define{Lie!equation} of that curve is the equation \(A(t)=\dot{g}(t)\hook\omega_G\).
The \emph{solutions} of a Lie equation are the curves \(g(t)\in G\) with Darboux derivative \(A(t)\) \cite{Bryant:1991}.
Each Lie equation \(A(t)=\dot{g}\hook\omega_G\) is equivalent to the equation of flow \(\dot{g}=X_t\) of the time varying vector field \(X_t(g)=\LT{g*}A(t)\).

Picard's theorem (theorem~\vref{thm:Picard}) constructs, for any \(C^k\) Darboux derivative \(A(t)\), a local \(C^{k+1}\) solution \(g(t)\) with any initial condition \(g(t_0)=g_0\). 
Since \(\omega_G\) is left invariant, left translations take solutions to solutions.
Glue local solutions together by left translations to make a global solution.
\begin{theorem}\label{theorem:Lie.equations}
Every \(L^1_{\text{loc}}\) Lie equation has a unique global solution, with any given initial condition.
In other words, for any \(L^1_{\text{loc}}\) Darboux derivative, defined on an interval of the real number line, the associated Lie equation admits a locally absolutely continuous solution, defined on the same interval, with any given initial condition at any time in that interval.
Any other solution defined on an interval is the same solution perhaps restricted to a subinterval and perhaps left translated.
If the Darboux derivative is \(C^k\) then the solution is \(C^{k+1}\), for any integer \(k\ge 1\).
If the Darboux derivative is real analytic then the solution is real analytic.
\end{theorem}
\begin{proof}
The Darboux equation is smooth in position, locally uniformly in time, and locally integrable in time.
By Picard's existence and uniqueness theorem in any coordinates \cite{Sontag:1998} C.3.8 p. 482, there is a local solution passing through a given point of the Lie group at a given time, unique if we ask that it be defined on a maximal time interval.

In other words, we can cover the domain of the Darboux derivative \(A(t)\) with time intervals, each open in the domain of \(A(t)\), on each of which there is a solution \(g(t)\) to our Lie equation, with some initial condition at some time in that time interval.
We can assume that each interval is maximal, subject to existence and uniqueness of a solution with that initial condition on that interval.

Fix one such maximal time interval.
By left translation of any one of these solutions, i.e. replacing \(g(t)\) by \(g_0g(t)\) for some constant \(g_0\in G\), we get another solution.
Varying \(g_0\), we can choose to have any initial condition we like, at any time in the time interval.
By local existence and uniqueness, again the interval is maximal for existence and uniqueness of a solution, for all possible initial conditions at at time in that time interval.

Take two such maximal time intervals, overlapping.
Take a solution on each time interval.
They agree on the overlap, up to left translation.
After that left translation, they glue together to a solution on the union of the intervals.
By maximality, the two intervals are equal.

Take one of these maximal intervals \(I\).
Take some \(t\in\bar{I}\) lying in the domain of the Darboux derivative \(A(t)\).
Then \(t\) also lies in one of our intervals, say \(J\).
Both being open, \(I\) intersects \(J\), so they are equal.
So \(I\) is closed and open in the domain of the Darboux derivative, hence a component of the domain, which is an interval, so connected, hence \(I\) is this domain.

By smoothness of local solutions, \(C^k\) Darboux derivative ensures a \(C^{k+1}\) solution \cite{Sontag:1998} C.3.11 p. 484.
\end{proof}
\begin{example}
If \(G=\GL{V}\) for some vector space \(V\),  write \(\omega_G=g^{-1}\,dg\), so our differential equation is
\[
\frac{dg}{dt}=gA(t),
\]
a linear ordinary differential equation.
Instead of our theorem above, we can just directly employ Picard's theorem \cite{Sontag:1998} C.4 p. 487 to prove existence, uniqueness and smoothness of solution of any linear ordinary differential equation.
\end{example}
\begin{example}
If we can solve a Lie equation in the identity component of \(G\), we can solve the same Lie equation in \(G\), by left translation.
So we can effectively assume that  \(G\) is connected when we try to solve Lie equations.
\end{example}
\begin{example}
Take a discrete normal subroup \(\Gamma\subseteq G\), and the quotient
\[
1\to\Gamma\to G\to\bar{G}\to 1.
\]
Then \(G\to\bar{G}\) is a covering map, so solve in \(\bar{G}\) and lift to \(G\), or solve in \(G\) and project to \(\bar{G}\).
So to solve Lie equations, it is enough to assume that \(G\) is connected and simply connected.
On the other hand, we could just as well quotient by any discrete normal subgroup if we prefer.
\end{example}
\begin{example}
If \(G\) is a Lie subgroup of a Lie group \(G'\), we can solve any Lie equation for \(G\) inside \(G'\) instead.
By the existence and uniqueness from theorem~\vref{theorem:Lie.equations}, the solution lies in \(G\).
By Ado's theorem \cite{Hilgert.Neeb:2012} p. 189 theorem 7.4.1, there is a faithful representation of the Lie algebra \(\LieG\) of \(G\).
By Lie's Third Theorem \cite{Hilgert.Neeb:2012} p. 334 9.4.11 (which follows easily from Ado's theorem), if \(G\) is connected and simply connected, there is a representation \(G\to\GL{V}\) with discrete kernel.
Hence the solution of the Lie equation can be found inside \(\GL{V}\), i.e. as the linear system of ordinary differential equations \(\dot{g}=gA(t)\), for a given matrix \(A(t)\) and an unknown matrix \(g(t)\).
\end{example}
\begin{example}
If \(G\) is abelian and connected, its universal covering Lie group is Euclidean space.
The Lie equation in Euclidean space is \(\dot{g}=A(t)\), with solution \(g(t)=\int A(t)\,dt\).
\end{example}
\begin{example}
If \(G\) is a semidirect product \(G=N\rtimes Q\), say with multiplication
\[
(n,q)(n',q')=(nqn'q^{-1},qq')
\]
then left translation in \(Q\) is independent of that in \(N\).
Decomposing each \(A\) in the Lie algebra of \(G\) into \(A=(A^N,A^Q)\), our Lie equation \(\dot{g}=\LT{g*}A\) becomes a Lie equation \(\dot{q}=\LT{q*}A^Q\) in \(Q\), and then a ``coupled'' equation 
\[
\dot{n}=\LT{nq*}\RT{q*}^{-1}A^N
\]
in \(N\).
If we solve the equation in \(Q\), we can rewrite the coupled equation as a Lie equation
\[
\dot{n}=\LT{n}B,
\]
in \(N\), with
\[
B:=\Ad_q A^N.
\]
By repeated application of this trick, we can solve any Lie equation in any solvable Lie group in quadratures, reducing to abelian Lie group factors.
More generally, we can use this trick to solve any Lie equation whose Darboux derivative lies in a solvable subalgebra.
By Levi decomposition, Lie equations reduce to those of semisimple Lie groups, and hence of simple Lie groups, hence of simple Lie groups in adjoint form.
\end{example}
\begin{problem}{cq}
The \(G\) be the group of linear transformations of a \(2\)-dimensional real or complex vector space which preserve a particular line through the origin. 
Solve the general Lie equation for \(G\).
\end{problem}
\begin{example}
On the other hand, we don't know of any trick to integrate Lie equations explicitly for any simple Lie group.
The lowest dimensional example is the group \(\SO{3}\) of rotations of Euclidean space: given the angular momentum \(A(t)\) of a rigid body, we don't see how to compute explicitly its rotation \(g(t)\), the solution of the Lie equation \(g^{-1}\dot{g}=A\), unless all \(A(t)\) fix the same point.
Similarly, given a one parameter family of infinitesimal projective transformations of the projective line, we don't know how to find the resulting projective transformation to which they integrate, unless they all fix the same point.
\end{example}
\section{Continuity}
We won't ever need the continuity theory for Lie equations, but we include it as a trivial example of applying continuity for time varying vector fields.
\begin{corollary}
Pick a Lie group \(G\) with Lie algebra \(\LieG\) and left invariant Maurer--Cartan form \(\omega\).
Pick a left invariant Riemannian metric \(g\) on \(G\).
Denote \(g\)-distance as \(d\).
Let \(R>0\) be a constant smaller than the injectivity radius of the metric \(g\) at the identity element, hence at every element of \(G\).
Denote parallel transport along a minimal geodesic from points \(x\) to \(y\) by \((x\to y)\).
Let
\[
C:=\sup_{d(x,1)\le R}\frac{\|(x\to 1)L_{x*}-I\|}{d(x,1)},\]
where the numerator is an operator norm.

Pick a locally integrable function
\[
t\in I\xmapsto{A}A_t\in\LieG
\]
on an interval \(I\subseteq\R\), perhaps unbounded, perhaps not open, perhaps not closed.
Denote its \(L^1\)-norm by \(\left\|A\right\|_1\).
Pick a time \(t_0\in I\).
Pick a point \(x_0\in G\) and let \(x(t)\) be the solution of the Lie equation \(\dot x\hook\omega_G=A\) with \(x(t_0)=x_0\).

Pick a locally integrable function \(I\xrightarrow{B}\LieG\) on the same interval \(I\subseteq\R\).
Pick a point \(y_0\in G\).
If
\[
d(x_0,y_0)+\left\|B\right\|_1
\le
Re^{-C\left\|A+B\right\|_1},
\]
then the solution \(y(t)\) to the Lie equation
\[
\dot y\hook\omega=A+B
\]
with \(y(t_0)=y_0\) is defined on \(I\) and satisfies
\[
d(x(t),y(t))
\le
\left(
d(x_0,y_0)
+
\int_{t_0}^t |B_s|ds
\right)
e^{C\int_{t_0}^t |A_s+B_s|ds}.
\]
\end{corollary}
\begin{proof}
By left invariance of the metric,
\[
L_{a*}(b\to c)=(a^{-1}b\to a^{-1}c).
\]
We have \(X_t\hook\omega=A_t\), \(Y_t\hook\omega=B_t\), so
\begin{align*}
\frac{\left|(y\to x(t))(X+Y)_t(y)-(X+Y)_t(x(t))\right|}{d(x(t),y)}
&=
\frac{\left|(y\to x)L_{y*}(A+B)_t-L_{x*}(A+B)_t\right|}{d(x(t),y)},
\\
&\le
\frac{\left|(x^{-1}y\to 1)L_{x^{-1}y*}-I\right|}{d(x(t),y)}|(A+B)_t|,
\\
&=C|(A+B)_t|.
\end{align*}
Apply theorem~\vref{thm:continuity}.
\end{proof}
We can do a bit better: let \(G^0\subseteq G\) the identity component and \(\tilde{G}^0\to G^0\) the universal covering group.
By left invariance, the estimates above are identical in \(G\) and in \(G^0\).
The flow lines lift up and project down by \(\tilde{G}^0\to G^0\).
So it suffices to compute \(C\) and \(R\) in \(\tilde{G}^0\) rather than in \(G\).
But the estimates above, computed in \(\tilde{G}^0\) could be better, and are no worse, than in \(G\).

\chapter{Vector field orbits}
This section doesn't discuss Cartan geometries, but builds more essential tools.
\section{The orbit theorem}
Take a set \(\VF\) of smooth vector fields on a manifold \(M\), each vector field only defined on some open set, different vector fields perhaps defined on different open sets.
The \(\VF\)-\emph{orbit}\define{orbit} of a point \(m_0\in M\) is the smallest subset of \(M\) containing \(m_0\) and invariant under the flows (perhaps only locally defined!) of all vector fields in \(\VF\).
\[
\includegraphics[width=4cm]{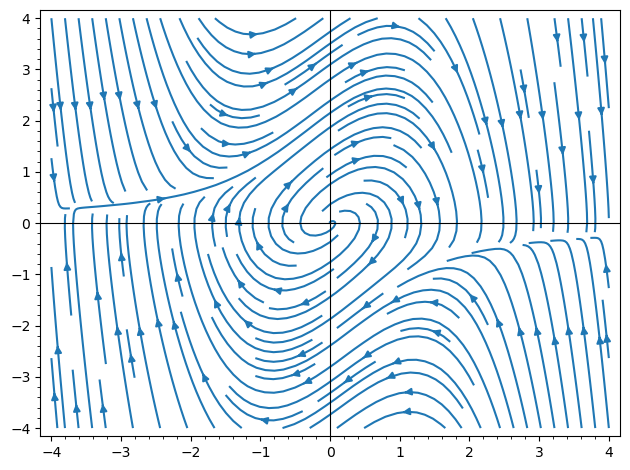}
\]
There might be points of \(M\) near which none of the vector fields in \(\VF\) are defined; these points are themselves orbits, by definition.
If two orbits intersect then they are equal.
The \emph{radical}\define{radical} \(\VF*\) of \(\VF\) is the set of all smooth vector fields, defined on various open subsets of \(M\), whose flows, where defined, preserve the orbits of \(\VF\). 

If \(X,Y\) are vector fields and if, for some real number \(s\), the diffeomorphism \(e^{sX}\) is defined on the domain of \(Y\) then \(e^{sX}_*Y\) is defined on the image of that domain.
We refer to this process as \emph{pushing around}\define{pushing around} vector fields.
\begin{lemma}\label{lemma:radical}
Any set of smooth vector fields has the same orbits as its radical does.
The radical is closed under pushing around.
If a vector field belongs to the radical, so does its rescaling by any smooth function defined on its domain.
The radical is a sheaf of sets: 
\begin{itemize}
\item
If a vector field belongs to the radical, so does its restriction to any open subset.
\item
if some vector fields, defined on various open sets, belong to the radical, and agree where the open sets overlap, then they are the restrictions to those open sets of an element of the radical.
\end{itemize}
\end{lemma}
\begin{proof}
Suppose that our set of vector fields is \(\VF\), and they live on a manifold \(M\).
We can assume that \(M\) is connected, and hence second countable.
If we multiply an element of \(\VF\) by a smooth function, we reparameterize the orbits, but don't enlarge them.
So the radical is closed under scaling by smooth functions.
If \(X,Y\in\VF*\) and the vector field \(Z:=e^X_* Y\) is defined then the flow of \(e^X_*Y\) is 
\[
e^{tZ} = e^X e^{tY},
\]
which preserves orbits.
So \(Z\in\VF*\): \(\VF*\) is closed under pushing around.
\end{proof}
By rescaling by smooth functions:
\begin{lemma}
For any point of the manifold, and any element of the radical defined near that point, there is a globally defined and compactly supported element of the radical which agrees with that element near that point.
Any set of vector fields has the same orbits as the globally defined and compactly supported elements of its radical.
\end{lemma}
Let \(\check\VF\) be the smallest set of locally defined vector fields containing \(\VF\) and closed under constant rescaling and pushing around.
So \(\VF\) and \(\check\VF\) have the same orbits.

For each point \(m_0 \in M\), take as many vector fields as possible \(X_1,\dots,X_k\), out of \(\check\VF\), which are linearly independent at \(m_0\). 
Refer to the number \(k\) of vector fields as the \emph{orbit dimension}.\define{orbit!dimension}  
Push around: the orbit dimension is a constant throughout each orbit.
The map
\[
\left(t_1, \dots, t_k\right) \in \text{open } \subseteq \R^k \mapsto
e^{t_1 X_1} \dots e^{t_k X_k} m_0 \in M
\]
is a \emph{flow parameterization}.\define{flow parameterization}
\[
\includegraphics{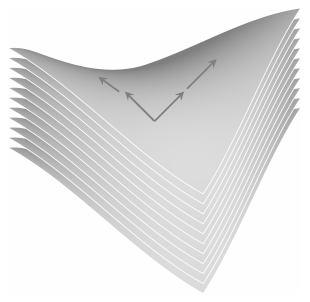}
\]
We assume it to be defined in some open set on which it is a smooth embedding.
Its inverse is a \emph{flow chart}\define{flow chart} and its image is a
\emph{flow set}.\define{flow set}
The \emph{orbit topology} on an orbit is the coarsest topology in which all flow charts are continuous.
\begin{theorem}[Sussmann's Orbit Theorem \cite{Sussmann:1973}]%
\label{thm:orbit}%
\define{theorem!orbit}%
\define{orbit theorem}
Take a set of smooth vector fields on a manifold.
The orbit of any point under any set of smooth vector fields is a smooth immersed submanifold, for a unique topology and smooth structure so that every flow line of any of the vector fields is continuous.
In this topology, the orbit is connected and second countable.
In the smooth structure, it is initial.
A dense open set of points lie on an orbit of locally maximal dimension.
The vector fields of the set, and those of the radical, restrict to be smooth vector fields on every orbit.
The vector fields of the radical span every tangent space to every orbit.
The radical is a sheaf of Lie algebras of vector fields, and a sheaf of modules over the sheaf of smooth functions, closed under pushing around, and closed under pointwise convergence.
A vector field, defined on an open set, belongs to the radical just when it is tangent to all of the orbits.
\end{theorem}
\begin{proof}
The tangent space to a flow set at each of its points is spanned by the linearly independent vector fields
\[
X_1,e^{t_1 X_1}_* X_2, \dots,e^{t_1 X_1}_* \dots e^{t_{k-1} X_{k-1}}_*
X_k,
\] 
which belong to \(\check\VF\), since they are just pushed around copies of the \(X_j\).
Let \(U\) be a flow set.
Suppose that \(Y\in\check\VF\) is a vector field, which is not tangent to \(U\). 
Then at some point of \(U\), \(Y\) is not a multiple of those pushed around vector fields, so the orbit dimension exceeds \(k\), a contradiction.

Therefore all vector fields in \(\check\VF\) are tangent to all flow sets. 
Each flow set is an embedded submanifold of \(M\), with tangent spaces having these linearly independent vector fields as a basis.
Each flow chart is a smooth chart on an embedded submanifold.
So every smooth map to \(M\) with image inside a single flow set is smooth in that flow chart.

Since the vector fields of \(\check\VF\) are all tangent to the flow sets, which are embedded submanifolds of \(M\), any point inside any flow set stays inside that flow set under the flow of any vector field in \(\check\VF\), at least for a short time.
Any point stays inside the flow set under compositions of flows of the vector fields, at least for a short time. 
Therefore a point belonging to two flow sets remains in both of them under the flows that draw out either of them, at least for short times. 
Therefore that point belongs to a smaller flow set lying inside both of them, given by restricting one or the other flow parameterization of one or the other flow set. 
Therefore the intersection of flow sets is a union of flow sets.
An open set of an orbit is thus precisely a union of flow sets.

The flow charts make the orbit locally homeomorphic to Euclidean space.  
Every open subset of \(M\) intersects every flow set in a flow set, so intersects every open set of the orbit in an open set of the orbit. 
So the inclusion mapping of the orbit into \(M\) is a continuous injection.
So the orbit topology is Hausdorff.
In this topology, the orbit is path connected, since the flow lines of the vector fields are continuous.

Next we will prove that the orbit is second countable, hence a topological manifold.
To prove this, it suffices to cover the orbit in a countable set of second countable open sets.

A \emph{splitting} is a diffeomorphism
\[
W\to U\times V
\]
of an open set \(W\subseteq M\) to a product of manifolds \(U,V\).
\[
\includegraphics{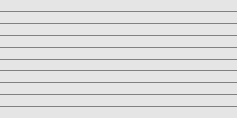}
\]
A \emph{slice}\define{slice!of a family of vector fields} \(S_{v_0}\) is the preimage in \(W\) of a set
\[
U\times\set{v_0}\subseteq U\times V.
\]
\[
\includegraphics{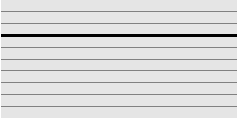}
\]
For any positive integer \(k\), a \emph{\(k\)-splitting}\define{splitting!of a family of vector fields} for the vector fields \(\VF\) is a splitting so that every \(k\)-dimensional orbit of \(\VF\) intersects \(W\) in a union of \(k\)-dimensional slices, each an open subset of that orbit.

Pick a point \(m_0\in M\), say with orbit dimension \(k\).
Pick vector fields \(X_1,\dots,X_k\in\VF_1\), linearly independent at \(m_0\).
Take an embedded submanifold \(V\subseteq M\) through \(m_0\) whose tangent space at \(m_0\) is complementary to the span of the \(X_j\).
After perhaps shrinking \(V\), we can arrange that all \(X_j\) are linearly independent at every point of \(V\).
So every point \(v\in V\) has orbit dimension at least \(k\).
The map
\[
(t,v)\mapsto e^{t_1X_1}\dots e^{t_kX_k}v\in M
\]
is a diffeomorphism, at least near \((0,m_0)\).
The inverse of this map is a splitting, perhaps after restricting to some smaller domain.
The orbit dimension at \(v=m_0\) is \(k\), and is at least \(k\) at every point \(v\in V\).
For points \(v\in V\) whose orbit dimension exceeds \(k\), the slice \(S_v\) is only an embedded smooth submanifold lying inside a single orbit.
But if the orbit dimension at a point \(v\in V\) is \(k\), then by definition of the flow charts, this splitting has slice \(S_v\) a flow set, so open.
So the splitting is a \(k\)-splitting.
So every point lies in the domain of a \(k\)-splitting, where \(k\) is the dimension of the orbit of that point.

Once we pick a particular dimension \(k\), we can cover \(M\) in domains of \(k\)-splittings.
Since \(M\) is second countable, we can pick a countable open covering by domains of \(k\)-splittings
\[
M_j\to U_j \times V_j.
\]

Consider the overlaps of two splittings:
\[
\includegraphics{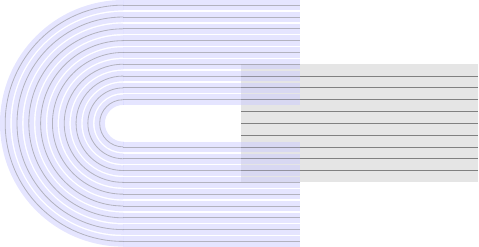}
\]
Pick a slice \(S=S_v\subseteq M_i\).
Consider its intersection with some other \(M_j\).
\[
\includegraphics{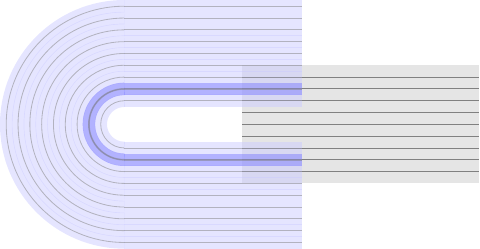}
\]
The intersection \(S\cap M_j\) is a union of open sets of slices in \(M_j\).
But the intersection is also an open subset of \(S\).
Since \(S\) is a second countable manifold, every open set of \(S\) has only countably many components.
So \(S\cap M_j\) intersects countably many slices of \(M_j\).

Start at a point \(m_0\).
As we move along any path in its orbit, continuous in the orbit topology, we pass through open sets \(M_{j_1},M_{j_2},\dots,M_{j_N}\) along a sequence of slices
\[
S_{v_{j_1}}, 
S_{v_{j_2}}, 
\dots
S_{v_{j_N}}.
\]
The union of these slices is a open set of the orbit, covering that path.

For any given initial point \(m_0\), the first slice is determined by \(m_0\), the second by a choice from among countably many, and so on.
So the set of such sequences, starting at a given point \(m_0\), is countable.
So the orbit is covered by countably many slices.
Each slice is second countable and open in the orbit.
So the orbit is second countable.

Therefore each orbit is a second countable topological manifold.
The flow charts give the orbit a smooth structure, so it is a smooth manifold.
The orbit is path connected.

The inclusion mapping of the orbit into \(M\) is a smooth immersion since it is locally the inclusion of slices.

Take a smooth (or continous) curve in \(M\) which lies in the domain of a \(k\)-splitting.
Only countably many slices of the \(k\)-splitting lie in the orbit.
If the curve travels from one slice to another, it passes through uncountably many slices, so leaves the orbit.
Conversely, if the curve stays in the orbit, it stays in a single slice.
Each slice is an open set in the orbit, so the curve maps smoothly (or continuously) to the orbit.

Any curve in \(M\) locally lies in the domain of a \(k\)-splitting.
So any smooth (or continuous) curve in \(M\) which stays in a single orbit is smooth (or continuous) as a map to that orbit.

Take any smooth (or continuous) map \(Z\to M\) from a manifold \(Z\) (or a locally path connected topological space).
Following paths inside it, the map is smooth (or continuous) to the orbit.

Apply this to the flow of some element of \(\VF*\): every element of \(\VF*\) has all its flow lines from a given point staying inside a slice, at least for short time.
Every element of \(\VF*\) restricts to every orbit to remain a smooth vector field.
The orbits are immersed submanifolds tangent to all elements of \(\VF*\).

Conversely, take a vector field on an open subset of \(M\), which is tangent to all orbits.
By the same argument, its flow preserves all orbits, so the vector field belongs to \(\VF*\).
So \(\VF*\) is the set of vector fields tangent to all orbits.
If two vector fields are tangent to an immersed submanifold, then their rescalings by smooth functions, their sum, and their bracket, are also tangent to that immersed submanifold.
So \(\VF*\) is a sheaf of modules over smooth functions, and a sheaf of Lie algebras.
Since any immersion is locally an embedding, the restriction of any element of \(\VF*\) to any orbit remains smooth.

Take a sequence \(X_1,X_2,\dots\in\VF*\) converging pointwise to a vector field \(X\).
At each point \(m\in M\), the values \(X_j(m)\) lie in tangent space to the orbit, which is a finite dimensional vector space.
So \(X(m)\), at that point, lies there too.
Pick an open subset \(U\) of the orbit, embedded into \(M\).
Then \(X\) is tangent to \(U\) at every point of \(U\).
So the flow of \(X\) preserves \(U\) for short times.
So \(X\) preserves the orbit, so \(X\in\VF*\).

Suppose we find a new smooth structure so that the orbit is smoothly immersed in the ambient manifold \(M\).
Let \(\mathscr{O}\) be the orbit in the old smooth structure, \(\mathscr{O}'\) the orbit in the new smooth structure; as subsets of \(M\) they are equal, but perhaps have different topologies or smooth structures.
Since \(\mathscr{O}\) is initial, \(\mathscr{O}'\to\mathscr{O}\) is a smooth bijection.

Suppose that every flow line of any element of \(\VF\) is continuous in \(\mathscr{O}'\).
Since the flow lines are continuous, each one stays for a little while in an open set on which the inclusion map \(\mathscr{O}'\to M\) is a smooth embedding.
The flow lines are smooth in \(M\), each either an immersed curve or a fixed point.
Therefore they are smooth in those smoothly embedded submanifolds, so in \(\mathscr{O}'\).
So the vector fields of \(\VF\) restrict to smooth vector fields.
Pushing around, every element of \(\check\VF\) restricts to be smooth in \(\mathscr{O}'\).

The flow charts are thus smooth bijections in this new smooth structure.
The flow charts are charts for \(\mathscr{O}\).
So the identity map \(\mathscr{O}\to\mathscr{O}'\) is a smooth bijection.

The vector fields in \(\VF*\) span the orbit tangent space at each point.
Take any maximal set of the vector fields which are linearly independent at some point.
Then they are still defined and linearly independent nearby, hence the orbits of all nearby points are of at least as large a dimension.
Hence if a point has locally maximal dimensional orbit, so do all nearby points.
\end{proof}
\Danger{}~Many proofs of the Orbit Theorem do not prove that the orbits are second countable.
\begin{problem}{Sussman:compact}
Suppose that \(\VF\) is a set of smooth vector fields, globally defined, on a manifold \(M\).
Let \(\VF_c\) be the set of compactly supported smooth vector fields whose flows preserve the orbits of \(\VF\).
Prove that \(\VF_c\) is a Lie algebra of vector fields, and a module over the algebra of smooth functions, and has the same orbits as \(\VF\).
Prove that the group of diffeomorphisms of \(M\) preserving the orbits of \(\VF\) acts with the same orbits as \(\VF\), even if no element of \(\VF\) is a complete vector field.
\end{problem}
\section{Analytic vector fields}
Replace every instance of \emph{smooth} by \emph{real analytic}: lemma~\vref{lemma:radical} and the Orbit Theorem (theorem~\vref{thm:orbit}) remain true, with the same replacements in the proofs. 
The orbits of any set of real analytic vector fields thus form real analytic immersed second countable initial submanifolds.

The Orbit Theorem does not tell us how to find the tangent spaces of the orbits.
It only says that the tangent space of the orbit through each point is the span of the values of vector fields obtained from all possible pushings around.
There are very few examples where we can computing a pushing around.
We will give a practical recipe applicable to many real analytic examples: the Hermann--Nagano theorem.
\section{The Lie saturate}
Suppose that \(\VF\) is a set of smooth vector fields defined on open subsets of a manifold \(M\).
The \emph{Lie saturate}\define{Lie saturate} \(\VF'\) of \(\VF\) is the smallest set of smooth vector fields, defined on various open subsets of \(M\), 
\begin{itemize}
\item
containing \(\VF\) and
\item
containing the zero vector field and
\item
closed under constant rescaling and
\item
closed under addition and brackets: for any two vector fields of \(\VF'\), their sum and Lie bracket, defined on the intersections of their domains, also belong to \(\VF'\).
\end{itemize}
\begin{example}
If all vector fields in \(\VF\) are globally defined, \(\VF'\) is the Lie algebra generated by \(\VF\).
\end{example}
For each point \(m\in M\), the \emph{values at the point} \(m\) of \(\VF\) are
\[
\VF_m:=\set{v(m)|v\in\VF}.
\]
We will find a criterion for a set of real analytic vector fields so that the values of the Lie saturate at each point form the tangent space to the orbit at that point.
In particular, the values are of constant dimension along orbits, which can help to find the orbits.
\begin{example}
Let \(\VF\) be the set consisting of two vector fields in the plane: \(\partial_x\) defined on the plane, and \(\partial_y\) defined on the right half plane.
Then \(\VF'\) is the set of all vector fields of the form \(a\partial_x\), defined on the plane, or of the form \(a\partial_x+b\partial_y\) defined on the right half plane.
The plane is a single orbit of \(\VF\), since we can use \(\partial_x\) to get to the right half plane, then use \(\partial_y\) to reach the desired \(y\) coordinate, then use \(\partial_x\) again to reach the desired \(x\) coordinate.
In this example, the values of the Lie saturate are smaller than the tangent spaces to the orbit at all points in the closed left half plane.
\end{example}
Calculation of \(\VF'\) from \(\VF\) only requires differentiating, not solving differential equations, so is often easier than computing \(\VF*\).
\section{Generic regularity}
The \emph{rank}\define{rank!family of vector fields} of \(\VF\) at \(m\) is the dimension of the span of its values \(\VF_m\subseteq T_m M\).
On a connected manifold, a set of vector fields is \emph{generically regular}\define{generically regular!family of vector fields} if its rank is maximal on a dense open set.
On a disconnected manifold, a set of vector fields is \emph{generically regular}\define{generically regular!family of vector fields} if the restriction of those vector fields to any component of the manifold is generically regular.
\begin{example}
Any set of globally defined real analytic vector fields on a manifold is generically regular: take a finite set which have some rank at a point, and they have at least that rank at every point.
\end{example}
\begin{example}
By the Orbit Theorem, for any set \(\VF\) of vector fields on a manifold, the set \(\check\VF\) given by pushing around the vector fields of \(\VF\) is generically regular.
\end{example}
\begin{example}
Take a real analytic vector subbundle \(V\subseteq TM\) of the tangent bundle of a manifold \(M\).
The set of real analytic local sections of \(V\) is generically regular.
\end{example}
\begin{example}
Take a real analytic vector bundle \(V\to M\) and a real analytic vector bundle morphism \(V\to TM\).
Claim: the set of vector fields on \(M\) which are locally in the image of this morphism is generically regular.
Proof: \(V\) is locally trivial.
On each open set on which \(V\) is trivial, map a basis of local sections of \(V\). 
The morphism has maximal rank on a dense open set.
\end{example}
\begin{lemma}\label{lemma:tgt.orb}
Take a set \(\VF\) of real analytic vector fields defined on open subsets of a manifold.
Suppose that its Lie saturate is generically regular.
Then the values of the Lie saturate at each point form the tangent space to the orbit of \(\VF\) through that point.
In particular, the values of the Lie saturate are of constant dimension along orbits.
\end{lemma}
\begin{proof}
We can assume that \(M\) is connected without loss of generality.
The orbits are constructed via flow charts directly from the vector fields in \(\VF\), which are real analytic, so the orbits are real analytic in that real analytic smooth structure.
As above, let \(\check\VF\) be the smallest set of vector fields containing \(\VF\) and closed under pushing around.
We have already seen that the orbit \(\mathscr{O}\) through any point \(m\in M\) has tangent space
\[
T_m\mathscr{O}=\operatorname{span}\check\VF_m.
\]
So we need to show that \(\check\VF\) and \(\VF'\) span the same values in each tangent space of \(M\).
\begin{problem}{analytic.pushforward}
Take vector fields \(v,w\) defined on open sets of a manifold \(M\).
On the intersection of the open sets on which \(v\) and \(w\) are defined, prove that the Taylor series for \(e^{tv}_*w\) in \(t\) consists of repeated Lie brackets:
\[
e^{tv}_*w=w+t\LieDer_v w+\frac{t^2}{2}\LieDer^2_v w+\dots
\]
\end{problem}
\begin{answer}{analytic.pushforward}
One point that could be confusing: note that this is not a Taylor series on the manifold \(M\).
Once we fix a point \(m\in M\), the map
\[
t\mapsto (e^{tv}_*w)(m):=(e^{tv})'(e^{-tv}m)w(e^{-tv}m)\in T_m M
\]
is a curve in a single tangent space, a single vector space, so has a well defined Taylor series.
Consider for which \(m\) and \(t\) this is defined.
We have to have the flow of \(v\) defined for time \(-t\) through \(m\), and then have \(w\) defined at the image of that point.
So if \(v\) is defined on \(V\) and \(w\) on \(W\), and \(V_t\subseteq V\) is the set of points where \(v\) has flow defined for time \(t\), then we need
\[
m\in V_{-t}\cap e^{tv}(W\cap V_t).
\]
Since \(V_t\) gets smaller as \(|t|\) gets larger, and every point of \(V\) lies in \(V_t\) for some interval of \(t\) near \(t=0\), we can ensure that this is defined anywhere in \(V\) if we make \(t\) small enough.

On the interior of the set where \(v\) vanishes, the result is obvious.
By continuity, it suffices to prove the result on open sets covering the set where \(v\ne 0\).
By the flow box theorem, it suffices to prove the result for \(v\) a translation vector field on Euclidean space, where the result is obvious.
\end{answer}
Each of the terms in the infinite series on the right hand side lies in \(\VF'\).
Because \(v,w\) are analytic, the series converges to \(e^{tv}_*w\), at least on some open set, for all small enough \(t\).
In any finite dimensional vector space, the span of a sequence of vectors can only increase in dimension finitely many times.
So at each point of any open set where \(v,w,e^{tv}_*w\) are defined, and where the series converges, the span of the limit is already contained in the span of some finite number of terms.
In other words, if the series converges at a point \(m\in M\), then \(v,w,e^{tv}_*w\) all lie in \(\VF'_m\).
(Careful: pushing around may allow us to move the domain of \(e^{tv}_*w\) outside of the domain of \(w\), where the series is not defined.)
We can replace \(M\) by the dense open subset of \(M\) on which \(\VF_m\) has constant dimension, and still find an open set of points \(m\in M\) where \(e^{tv}_*w(m)\in\VF'_m\) for all small enough \(t\).

By the constancy of dimension, any local sections \(v_1,\dots,v_k\) of \(\VF'\) which form a basis at one point form a basis at every nearby point, a basis of local sections.
So 
\[
v_1\wedge v_2\wedge\dots\wedge v_k\wedge e^{tv}_*w
\]
on some open set.
The same holds wherever this wedge product is defined.
The same holds when we replace these local sections by any other basis of local sections.
So \(e^{tv}_*w(m)\in\VF'_m\) for every \(m\) where \(e^{tv}_*w\) is defined.
(Note that, by the Frobenius theorem, and bracket closure of \(\VF'\), the spaces \(\VF_m\) are the tangent spaces of the leaves of a foliation, the orbits of \(\VF\), throughout the dense open set of maximal dimensional \(\VF'_m\).)
\end{proof}
\section{Locally slack}
We still don't have a useful criterion to decide when the Lie saturate of a set of vector fields has values equal to the tangent spaces of the orbits.

We invent some terminology.
Take a set \(\VF\) of real analytic vector fields, defined on open sets in a manifold \(M\).
The set \(\VF\) is \emph{slack}\define{slack} on an open set \(U\subseteq M\) if, for any point \(m\in U\) and any vector field \(v\in\VF\) defined near \(m\), there are real analytic functions \(f_1,\dots,f_n\) defined near \(m\) and vector fields \(v_1,\dots,v_n\in\VF\), each defined on an open set containing \(U\), so that, on some open set around \(m\),
\[
v=f_1v_1+\dots+f_nv_n.
\]
(Roughly speaking, the vector fields globally defined on \(U\) give rise to all of the ones defined on any open subset of  \(U\).
Slack is a weaker condition than being a flabby sheaf, since the functions \(f_i\) are only defined near \(m\).
It differs from local freedom in that the coefficient functions \(f_i\) might not be arbitrary.)
The set \(\VF\) is \emph{locally slack}\define{locally slack!family of vector fields}\label{page:define.slack} if \(M\) is covered in open sets on which \(\VF\) is slack.
\begin{example}
The set of all real analytic sections of any vector subbundle of the tangent bundle is locally slack.
\end{example}
\begin{example}
Any finite dimensional set of global real analytic vector fields is locally slack.
\end{example}
\begin{lemma}\label{lemma:v.b.slack}
If \(V\to TM\) is a real analytic vector bundle morphism then the set of all local sections of \(TM\) which are locally in its image is locally slack.
\end{lemma}
\begin{proof}
Locally we can trivialize \(V\), and then take a basis of local sections of \(V\).
Every section in the image is locally the image of a linear combination of these.
\end{proof}
\begin{lemma}
If some real analytic globally defined vector fields on a connected manifold are linearly independent at a point, they are linearly independent on a dense open set.
\end{lemma}
\begin{proof}
Their wedge product vanishes everywhere if and only if it vanishes on some open set, by analyticity.
\end{proof}
\begin{corollary}
Any finite set of real analytic vector fields on a connected manifold have maximal rank on a dense open set.
\end{corollary}
\begin{lemma}\label{lemma:gen.reg}
Any locally slack set of vector fields is generically regular.
\end{lemma}
\begin{proof}
Take a locally slack set \(\VF\) of vector fields.
Take an open set \(U\) on which \(\VF\) is slack.
The rank of \(\VF\) is at most the dimension of the manifold, so bounded.
Pick a point \(m\in U\) of maximal rank of \(\VF\).
Pick out a finite number of vector fields from \(\VF\) whose values at \(m\) span all of the values \(\VF_m\).
Pick finitely sections defined on \(U\), so that the given sections are expressed as linear combinations of them, with real analytic coefficients, defined near \(m\).
These have maximal rank on a dense open subset of \(U\), hence so does \(\VF\).
\end{proof}
\begin{theorem}[Hermann--Nagano]\label{thm:Hermann.Nagano}\define{theorem!Hermann--Nagano}\define{Hermann--Nagano!theorem}
Take a locally slack set \(\VF\) of real analytic vector fields, defined on open subsets of a manifold.
Its Lie saturate \(\VF'\) is locally slack, being slack on the same open sets as \(\VF\), so is generically regular.
The vector spaces \(\VF'_m\) are the tangent spaces to the orbits of \(\VF\), and in particular are of constant dimension along orbits.
\end{theorem}
\begin{proof}
Take an open set \(U\) on which \(\VF\) is slack.
Consider the set of all iterated brackets
\[
\ad_{v_{i_1}}
\ad_{v_{i_2}}
\dots
\ad_{v_{i_{k-1}}}
v_{i_k}
\]
of elements of \(\VF\) which are defined on \(U\).
Every element of \(\VF'\), defined near some point \(m\in U\), is expressed by a finite sum of real analytic functions multiplied by iterated brackets of elements of \(\VF\) defined near \(m\), so by iterated brackets of elements of \(\VF\) defined on \(U\).
So \(\VF'\) is slack on \(U\).
Being locally slack, \(\VF'\) is generically regular by lemma~\vref{lemma:gen.reg}.
Apply lemma~\vref{lemma:tgt.orb}.
\end{proof}
\section{Maps matching vector fields}
Take a smooth map \(M\xrightarrow{\varphi}M'\), and smooth vector fields \(X,X'\) on \(M,M'\).
Write \(X\xrightarrow{\varphi}X'\) to mean that \(\varphi'\left(m\right) X\left(m\right) =X'(m')\) for any \(m\in M\) where \(m':=\varphi(m)\), i.e. \(X\) pushes forward to \(X'\).
For sets \(\VF,\VF'\) of vector fields on \(M,M'\), a \emph{correspondence}\define{correspondence!of set of vector fields} \(\VF\xrightarrow{\varphi}\VF'\) is a smooth map \(M\xrightarrow{\varphi}M'\) so that
\begin{enumerate}
\item for any \(X \in \VF\) there is an \(X'\in\VF*'\) with \(X\xrightarrow{\varphi}X'\) and
\item for any \(X' \in \VF'\) there is an \(X\in\VF*\) so that \(X\xrightarrow{\varphi}X'\).
\end{enumerate}
Say that \(\varphi\) \emph{pushes down}\define{push down!set of vector fields} \(\VF\) to \(\VF'\).
A vector field on \(M\) is \emph{equicomplete}\define{equicomplete} over \(M'\) if, for every fiber of \(\varphi\), the flow of the vector field is defined on some time interval \(-\varepsilon<t<\varepsilon\) through every point of that fiber.
A correspondence \(\VF\xrightarrow{\varphi}\VF'\) is \emph{equicomplete} if, for every \(X'\in\VF'\), there is a vector field \(X\in\VF*\) equicomplete over \(M'\) so that \(X\xrightarrow{\varphi}X'\).
\begin{theorem}%
[Ehresmann \cite{Ehresmann1961}, 
Sharpe \cite{Sharpe:1997} p. 93 proposition 8.8, 
McKay \cite{McKay:2004a}]%
\label{theorem:orbit.maps}%
\define{theorem!orbit maps}%
\define{orbit maps theorem}
If a set of vector fields pushes down via some smooth map then the map takes orbits to orbits by surjective submersions.
If the set of vector fields is equicomplete then orbits are fiber bundles over orbits.
\end{theorem}
\begin{proof}
Clearly \(\VF\xrightarrow{\varphi}\VF'\) takes each orbit to an orbit.
Given a path
\[
e^{t_1 X_1} \dots e^{t_k X_k}m_0
\]
down in \(M'\), we can lift it to one in \(M\), so \(\varphi\) maps each orbit onto an orbit, and maps onto every orbit.
We can ``push around'' even by incomplete vector fields, but only locally, preserving orbits.
Pushing upstairs in \(M\) corresponds to pushing downstairs in \(M'\), so without loss of generality, both \(\VF\) and \(\VF'\) are closed under ``pushing around''.
As in the orbit theorem~\vref{thm:orbit}, the flow charts are formed via pushing around, so \(\varphi\) is a submersion.

Suppose the family is equicomplete.
Multiply vector fields on \(M'\) by bump functions, and those on \(M\) by their pullbacks: we can assume all vectors fields in \(\VF'\) are globally defined and complete.
By restricting to an orbit in \(M\), we may assume that \(M\) and \(M'\) are orbits.
As in the proof of theorem~\vref{thm:orbit}, for each point \(m' \in M'\), construct
a flow chart
\[
\left(t_1, \dots, t_k\right) \mapsto e^{t_1 X'_1} \dots e^{t_k X'_k}
m'.
\]
Pick equicomplete vector fields \(X_k\) with \(X_k\xrightarrow{\varphi}X'_k\), hence complete.
Let \(U'\subseteq M'\) be the associated flow set; on \(U'\) these \(t_j\) are now coordinates.  
Let \(U:=\varphi^{-1} U'\subseteq M\). 
Let \(Z\) be the fiber of \(M\xrightarrow{\varphi}M'\) above the origin of the flow chart. 
Map
\[
u \in U \mapsto \left(u',z\right) \in U' \times Z
\]
by \(u'=\varphi(u)\) and
\[
z = e^{-t_k X_k} \dots e^{-t_1 X_1} u.
\]
This gives \(M\) the local structure of a product.  
The transition maps have a similar form, composing various flows, so \(M\to M'\) is a fiber bundle.
\end{proof}
\begin{corollary}\label{corollary:Ehresmann.thm}
A proper submersion mapping to a connected manifold is a fiber bundle mapping.
\end{corollary}
\begin{proof}
Take the family consisting of any vector field upstairs which pushes down to some vector field downstairs.
\end{proof}
\begin{problem}{cr}
Prove that a smooth map of manifolds is a fiber bundle map just when it pushes down an equicomplete family of vector fields.
\end{problem}
\begin{problem}{cs}
Prove that the composition of fiber bundle maps is a fiber bundle map.
\end{problem}
\section{An equivariant orbit theorem}
If \(\VF\) is a set of vector fields on a manifold \(M\) and \(H\) is a Lie group acting on \(M\), the \emph{\(H\VF\)-orbits}\define{orbit} are the minimal nonempty \(H\)-invariant unions of \(\VF\)-orbits.
\begin{theorem}%
[Blumenthal \cite{Blumenthal1986}]%
\label{cor:orbitMapEquiv}%
\define{theorem!Blumenthal morphism}
\define{Blumenthal morphism theorem}
Suppose that 
\begin{itemize}
\item
\(\G\xrightarrow{\varphi}\G'\) is a smooth map of manifolds
\item  
\(\VF\xrightarrow{\varphi}\VF'\) for some sets of vector fields
\item
\(\VF\) is equicomplete over \(\G'\)
\item
Lie groups \(H,H'\) act on \(\G,\G'\)
\item
\(\G\to\G'\) is equivariant for a Lie group morphism \(H\to H'\).
\end{itemize}
Then the \(H\VF\)-orbits fiber equivariantly over their images in the \(H'\VF'\)-orbits.
If \(H\to H'\) is onto, the image of each \(H\VF\)-orbit is an \(H'\VF'\)-orbit.
If \(H,H'\) act freely and properly on \(\G,\G'\) and \(H\to H'\) has closed image then the \(H\)-quotients of the \(H\VF\)-orbits fiber smoothly over their images in the \(H'\)-quotients of the \(H'\VF'\)-orbits.
\end{theorem}
\begin{proof}
Each \(H\)-translate of a single equicomplete vector field on \(\G\) is equicomplete with the same survival time of all flow lines, and maps to the corresponding \(H'\)-translated vector field.
Taking unions of these translates, we can suppose that \(\VF,\VF'\) are \(H,H'\)-invariant.

Let \(H''\) be the image of \(H\) in \(H'\), so isomorphic to the quotient of \(H\) by the kernel, hence a Lie group.
The \(H\)-action on \(\G\) gives a Lie algebra \(\LieH\)-action on \(\G\).
Add all vector fields of that \(\LieH\)-action to \(\VF\), and simultaneously those of the \(\LieH''\)-action to \(\VF'\).
We preserve equicompleteness, since these vector fields are complete.

By theorem~\vref{theorem:orbit.maps}, \(\VF\)-orbits fiber over \(\VF'\)-orbits.
Any \(H\)-translate of any \(\VF\)-orbit is an \(\VF\)-orbit; any two are equal or disjoint.
Translating by elements of the identity component, every \(\VF\)-orbit is invariant because \(\LieH\subseteq\VF\).
Hence the \(H\VF\)-orbit of a point is a disjoint union of \(\VF\)-orbits, all of the same dimension, at most one for each component of \(H\), so an immersed submanifold.
Each \(H\VF\)-orbit maps to an \(H''\VF'\)-orbit, which is a submanifold of an \(H'\VF'\)-orbit.
This map is also a fiber bundle map, with each fiber a union of components of \(H\) times a fiber of the \(\VF\)-orbit over the \(\VF'\)-orbit.

We can assume that \(H\) and \(H'\) act freely and properly and that \(H\to H'\) has closed image \(H''\subseteq H'\).
We can replace \(\G\) by a single \(H\VF\)-orbit and \(\G'\) by a single \(H'\VF'\)-orbit.
Our problem reduces to:
\end{proof}
\begin{proposition}
Suppose that \(\G\xrightarrow{\Phi}\G'\) is a fiber bundle mapping, equivariant for free and proper actions of Lie groups \(H,H'\) on manifolds \(\G\) and \(\G'\) for some Lie group morphism \(H\to H'\) with closed image.
Then \(M:=\G/H\to M':=\G'/H'\) is a smooth fiber bundle.
\end{proposition}
\begin{proof}
The actions are free and proper, so the quotients \(M:=\G/H, M':=\G'/H'\) are smooth manifolds with smooth maps
\[
\begin{tikzcd}
\G\arrow[r,"\Phi"]\arrow[d,"\pi"]&\G'\arrow[d,"\pi'"]\\
M\arrow[r,"\varphi"]&M'.
\end{tikzcd}
\]
Denote by \(K\) the kernel of \(H\to H'\).
Let \(M'':=M/K\), \(\G'':=\G/K\), \(H'':=H/K\).
Then \(M\to M''\) is a fiber bundle mapping, since \(K\) acts freely and properly.
So it suffices to prove that \(M''\to M'\) is a fiber bundle mapping.
So we can replace \(M\) by \(M''\), \(H\) by \(H''\), \(\G\) by \(\G''\), i.e. we can assume that \(H\subseteq H'\) is a closed subgroup.

So we suppose that \(\G\xrightarrow{\Phi}\G'\) is a fiber bundle mapping, equivariant for free and proper actions of Lie groups \(H,H'\) on manifolds \(\G\) and \(\G'\) for some closed subgroup \(H\subseteq H'\) of a Lie group \(H'\).
Let \(M'':=\G'/H\).
\[
\begin{tikzcd}
\G\arrow[r,"\Phi"]\arrow[d,"\pi"]&\G'\arrow[r,"\operatorname{id}"]\arrow[d,"\pi"]&\G'\arrow[d,"\pi'"]\\
M\arrow[bend right,rr,"\varphi"]\arrow[r]&M''\arrow[r]&M'.
\end{tikzcd}
\]
Clearly the right side,\(M''\to M'=\G'/H'\), is a fiber bundle with fibers \(H'/H\), since it is just the quotient map \(\G'/H\to\G'/H'\).
The composition of fiber bundles is a fiber bundle.
So it suffices to prove that \(M\to M''\) is a fiber bundle.
We will see that it has the same fibers as \(\G\to\G'\).

So we can suppose that \(\G\xrightarrow{\Phi}\G'\) is a fiber bundle mapping, equivariant for free and proper actions of a Lie group \(H\), and let \(M:=\G/H\) and \(M':=\G'/H\).
Each point of \(M'\) lies in the domain of a local section of the principal \(H\)-bundle \(\G'\to M'\).
The preimage in \(\G\) of the image of this local section is a fiber bundle over the domain in \(M'\), and strikes each \(H\)-orbit in at most one point.
By dimension count, it is a local section of \(\G\to M\). 
\end{proof}
\section{The Frobenius theorem}
A \emph{distribution}\define{distribution} or \emph{plane field}\define{plane field} is a vector subbundle \(V\subseteq TM\) of the tangent bundle of a manifold \(M\).
A \emph{\(V\)-tangent vector field} is a section of \(V\) as a vector bundle.
A distribution \(V\) is \emph{bracket closed}\define{bracket closed} when the set of \(V\)-tangent vector fields are closed under Lie bracket, hence form a Lie algebra.
A \emph{Pfaffian system}\define{Pfaffian system} \(I\) is a collection of \(1\)-forms whose span in each cotangent space has constant dimension; it is \emph{Frobenius}\define{Frobenius} if \(dI=0\) modulo the ideal generated by \(I\) in the algebra of differential forms; \(I\) \emph{cuts out} the distribution \(V=I^{\perp}\subset TM\).
\begin{theorem}[Frobenius]\label{theorem:Frobenius}\define{theorem!Frobenius}\define{Frobenius!theorem}
A distribution \(V\) is bracket closed just when some, hence any, Pfaffian system cutting out \(V\) is Frobenius, and this occurs just when \(V\) is the bundle of tangent vectors to the leaves of a foliation \(F\); we write \(V=TF\).
These leaves are the maximal connected integral submanifolds of \(V\), and are also are the orbits of the \(V\)-tangent vector fields.
\end{theorem}
\[
\includegraphics[width=5cm]{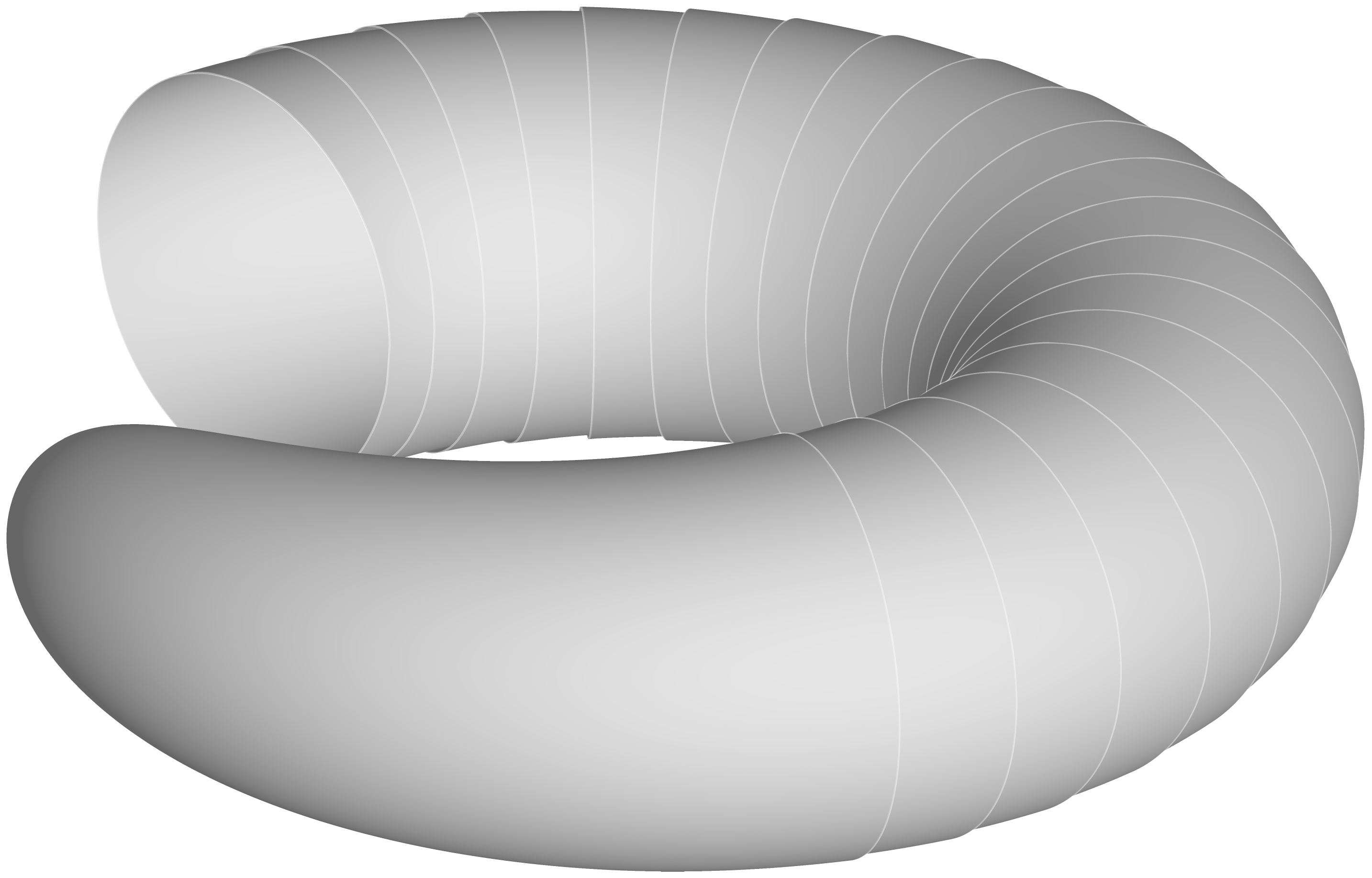}
\includegraphics[width=5cm]{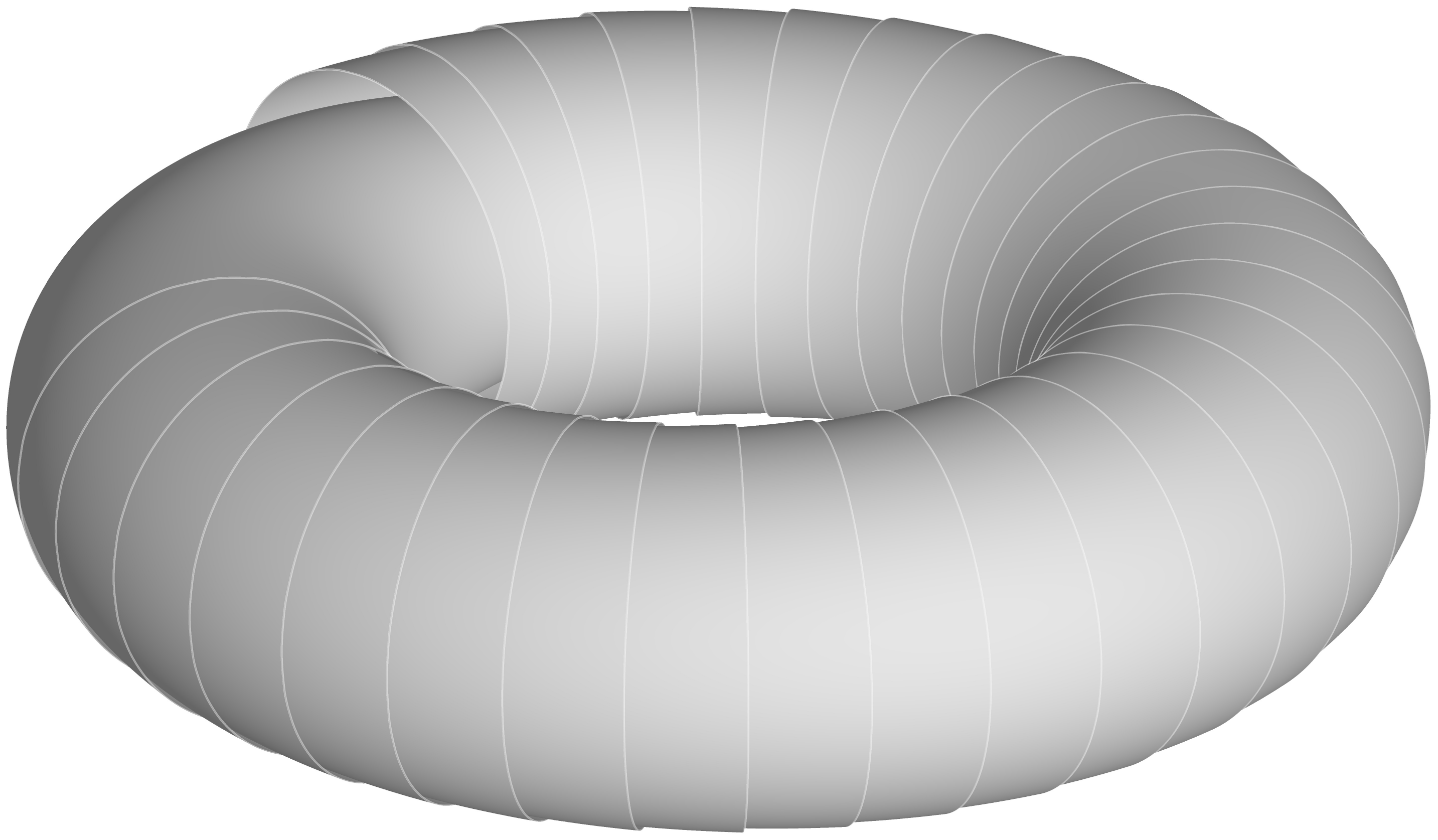}
\]
\begin{problem}{ct}
Use Sussmann's Orbit Theorem to prove the Frobenius theorem.
\end{problem}
\section{An equivariant Frobenius theorem}
Suppose that a Lie group \(H\) acts smoothly on a manifold \(M\) preserving a vector subbundle \(V\subseteq TM\), and that the vector fields of the Lie algebra of \(H\) are \(V\)-tangent; say that \(V\) is an \emph{\(H\)-distribution}.\define{distribution}
If \(V\) is bracket closed, say that the associated foliation \(F\) is an \emph{\(H\)-foliation}.
A \emph{folio}\define{folio} of \(F\) is the union of the \(H\)-translates of a leaf of \(F\); we may also call it an \emph{\(H\)-folio} or \emph{\(F\)-folio} or \emph{\(V\)-folio}.
\begin{corollary}[Equivariant Frobenius theorem]\label{corollary:equivariant.Frobenius}\define{Frobenius theorem!equivariant}\define{equivariant Frobenius theorem}\define{theorem!equivariant Frobenius}
Suppose that a Lie group \(H\) acts smoothly on a manifold \(M\) and that \(V\) is a bracket closed \(H\)-distribution.
Then \(V\) is the tangent bundle of a unique \(H\)-foliation \(F\).
Each point of \(M\) lies in a unique folio.
The folios are immersed submanifolds of dimension equal to the dimension of the leaves of \(F\).
If \(M\to \bar{M}:=M/H\) is a submersion of manifolds then \(F\) descends to a foliation on \(\bar{M}\), whose leaves are precisely the \(H\)-quotients of the folios.
\end{corollary}

\chapter{Semibasic differential forms}\label{appendix:semibasic}
Take a smooth map of manifolds \(P\xrightarrow{\varphi}Q\).
A smooth differential form \(\xi\) on \(P\) is \emph{basic}\define{basic differential form} if \(\xi=\varphi^*\eta\) for some smooth differential form \(\eta\) on \(Q\).
A smooth differential form \(\xi\) on \(P\) is \emph{semibasic}\define{semibasic differential form} if, for each point \(p_0\in P\), if we let \(q_0:=\varphi(p_0)\), then there is an alternating form \(\eta_{q_0}\) on \(T_{q_0} Q\) for which \(\xi_{p_0}=\varphi^*\eta_{q_0}\).
We might say that to be semibasic is to be \emph{pointwise basic}.
Basic forms are semibasic.
Our aim in this appendix is to find local and global conditions on a semibasic form under which it is basic.
For each point \(p_0\in P\), a vector \(v\in T_{p_0} P\) is \emph{vertical}\define{vertical!vector} if it is tangent to the fiber of \(P\to Q\).
\begin{example}
If \(P\to Q\) is
\[
\begin{tikzcd}
(x,y,z)\in\R^3\arrow[r]&x\in\R
\end{tikzcd}
\]
then 
\[
\begin{array}{@{}llll@{}}
\text{degree}&\text{forms}&\text{semibasic}&\text{basic}\\
\cmidrule(r){1-1}\cmidrule(lr){2-2}\cmidrule(lr){3-3}\cmidrule(l){4-4}
0&f(x,y,z)&\text{all}&f(x)\\[7pt]
1&f(x,y,z)dx&f(x,y,z)dx&f(x)dx\\
&\quad+g(x,y,z)dy\\
&\quad+h(x,y,z)dz\\[7pt]
2&f(x,y,z)dy\wedge dz&0&0\\
&\quad +g(x,y,z)dz\wedge dx\\
&\quad +h(x,y,z)dx\wedge dy\\[7pt]
3&f(x,y,z)dx\wedge dy\wedge dz&0&0
\end{array}
\]
\end{example}
We allow notation like \(dx^{ij}:=dx^i\wedge dx^j\), and so on.

Take a smooth submersion \(P\xrightarrow{\pi}Q\).
Take a point \(p_0\in P\).
Recall that the \emph{vertical vectors} in \(T_{p_0} P\) are those on which \(\pi'(p_0)=0\).
Functions \(y^1,\dots,y^p\) defined near \(p_0\) are \emph{nondegenerate} for the submersion if their differentials are linearly independent on vertical vectors and vanish on a vertical vector just when it is the zero vertical vector.
Take some nondegenerate functions \(y^1,\dots,y^p\) on \(P\) near a point \(p_0\).
Take any coordinates \(x^1,\dots,x^q\) on \(Q\) defined near \(\pi(p_0)\).
Denote their pullbacks to \(P\) also by \(x^1,\dots,x^q\).
Then 
\[
x^1,\dots,x^q,y^1,\dots,y^p
\] 
are coordinates near \(p_0\), by the implicit function theorem.
Call these \emph{adapted coordinates} to the submersion \(P\to Q\).
\begin{theorem}
Take a smooth submersion \(P\xrightarrow{\pi}Q\) and a differential form \(\xi\) on \(P\).
The following are equivalent:
\begin{itemize}
\item \(\xi\) is semibasic,
\item \(v\hook\xi=0\) for all vertical vectors \(v\),
\item \(\xi\) is a section of the vector bundle \(\pi^*\Lambda^*Q\),
\item In any adapted coordinates, \(\xi=a_I(x,y)dx^I\) for unique smooth functions \(a_I(x,y)\).
\end{itemize}
\end{theorem}
\begin{proof}
Suppose that \(\xi\) is semibasic \(k\)-form.
At an arbitrary point \(p_0\in P\), write 
\[
\xi_{p_0}=\pi(p_0)^*\eta_{q_0}.
\]
In other words, for any tangent vectors \(v_1,\dots,v_k\in T_{p_0} P\), if we let
\[
w_j:=\pi'(p_0)v_j,
\]
then
\[
\xi(v_1,\dots,v_k)=\eta(w_1,\dots,w_k).
\]
In particular, if any one of these \(v_j\) is vertical, the left hand side vanishes.
So \(v\hook\xi=0\) for any vertical vector \(v\), i.e. with \(\pi'(p_0)v=0\).

Any differential form can be expanded in any coordinates:
\[
\xi=a_{IJ}dx^I\wedge dy^J,
\]
Our vertical vectors are the span of the
\[
\frac{\partial}{\partial y^j}.
\]
Hence \(v\hook\xi=0\) for any vertical vector \(v\) just when \(\ell=k\), i.e. no \(dy^j\) factors for any \(j\).
Hence the coordinate expression as stated is equivalent to vanishing wedge with any vertical vector.

The pullback bundle \(\pi^*\Lambda^*Q\) has sections precisely the choices of an element of
\(
\Lambda_q^*Q
\)
at each point \(q=\pi(p)\), for each point \(p\in P\), smoothly varying with \(p\).
In our coordinates, this is a choice of 
\[
b_I dx^I,
\]
at each point \(p=(x,y)\), smoothly varying in \(x,y\), i.e. precisely such an expression
\(
\xi=a_Idx^I.
\)

Finally, if \(\xi\) admits an expression as \(\xi=a_Idx^I\), then at each point \(p=(x,y)\), clearly it is the pullback
\(
\pi^*\eta
\)
at a point \((x,y)\in P\) of 
\[
\eta = b_I dx^I,
\]
precisely if we take 
\[
b_I=a_I(x,y).
\]
\end{proof}
\begin{example}
Let \(P\) be the interior of an annulus in the plane, and \(Q\) the horizontal axis, \(P\xrightarrow{\varphi}Q\) an orthogonal linear projection:
\[
\begin{tikzpicture}[background rectangle/.style={fill=white,rounded corners}, show background rectangle]
\fill[gray!20,draw=gray!10] (0,1.1) circle (1);
\node at (1.3,1.1) {\(P\)};
\fill[white,draw=gray!10] (0,1.1) circle (.5);
\draw[gray!40] (-1,0) -- (1,0) node[right,black] {\(Q\)};
\end{tikzpicture}
\]
Let \(U_Q\subseteq Q\) be the open interval of points \(q\in Q\) over which the fibers \(P_q:=\varphi^{-1}\set{q}\) consist of two disjoint open intervals.
Let \(U_P:=\varphi^{-1}U_Q\).
\[
\begin{tikzpicture}[background rectangle/.style={fill=white,rounded corners}, show background rectangle]
\fill[gray!20,draw=gray!10] (0,1.1) circle (1);
\node at (1.3,1.1) {\(P\)};
\fill[white,draw=gray!10] (0,1.1) circle (.5);
\draw[gray!40] (-1,0) -- (1,0) node[right,black] {\(Q\)};
\draw[gray!50,very thick] (-.5,0) -- (.5,0);
\draw[gray,fill=gray!10] (-.5,0) circle (1pt);
\draw[gray,fill=gray!10] (.5,0) circle (1pt);
\begin{scope}
\clip(-.5,.1) rectangle (.5,2.1);
\fill[gray!50,draw=gray!10] (0,1.1) circle (1);
\fill[white,draw=gray!10] (0,1.1) circle (.5);
\end{scope}
\end{tikzpicture}
\]
In the usual Cartesian coordinates on the plane, all \(0\)-forms \(f(x,y)\) on \(P\) are semibasic; the localdly basic are those with \(\partial f/\partial y=0\).
The semibasic \(1\)-forms on \(P\) dare precisely the \(1\)-forms of the form \(f(x,y)dx\); again the locally basic are those with \(\partial f/\partial y=0\).
Zero is the only semibasic \(2\)-form on \(P\).
Take a bump function \(Q\xrightarrow{h_Q}\R\) vanishing outside \(U_Q\).
let
\[
h_P(x,y):=
\begin{cases}
h_Q(x),&\text{ in one component of \(U_P\)},\\
0,&\text{ in the other}.
\end{cases}
\]
\[
\begin{tikzpicture}[background rectangle/.style={fill=white,rounded corners}, show background rectangle]
\fill[gray!20,draw=gray!10] (0,1.1) circle (1);
\node at (1.3,1.1) {\(P\)};
\fill[white,draw=gray!10] (0,1.1) circle (.5);
\draw[gray!40] (-1,0) -- (1,0) node[right,black] {\(Q\)};
\draw[gray!50,very thick] (-.5,0) -- (.5,0);
\draw[gray,fill=gray!10] (-.5,0) circle (1pt);
\draw[gray,fill=gray!10] (.5,0) circle (1pt);
\shade[left color=gray!50,
                       right color=gray!50,
                       middle color=white] % <---
                       (-.1,-.01) rectangle (.2,.01);
\begin{scope}
\clip(-.5,.1) rectangle (.5,2.1);
\fill[gray!50,draw=gray!10] (0,1.1) circle (1);
\clip (0,1.1) circle (1);
\shade[left color=gray!50,
                       right color=gray!50,
                       middle color=white] % <---
                       (-.1,1.1) rectangle (.2,2.1);
\fill[white,draw=gray!10] (0,1.1) circle (.5);
\end{scope}
\end{tikzpicture}
\]
Then \(h(x,y)\) is a semibasic \(0\)-form on \(P\) and \(h(x,y)dx\) is a semibasic \(1\)-form on \(P\).
Both are locally basic, but neither is basic.
\end{example}
\begin{example}
For any covering map \(P\to Q\), for example the usual \(S^n\to\RP{n}\), every form is semibasic, but the basic forms are those invariant under the local deck transformations.
For a normal covering map, such as \(S^n\to\RP{n}\), we need only check the global covering transformations.
\end{example}
A \emph{vertical infinitesimal automorphism}\define{vertical!infinitesimal automorphism} of a differential form \(\xi\) on \(P\) is a vertical vector field on \(P\) with \(0=\LieDer_X \xi\).
A \emph{vertical automorphism}\define{vertical!automorphism} of a differential form \(\xi\) on \(P\) is a diffeomorphism of \(P\) preserving \(\xi\) and preserving the fibers of \(P\to Q\); let \(\Aut[\xi/Q]\) be the group of vertical automorphisms.
Consider the subgroup \(\Aut[\xi/Q]^0\) generated by the flows of complete vertical infinitesimal automorphisms.
Intuitively, we picture \(\Aut[\xi/Q]^0\) as like an identity component of \(\Aut[\xi/Q]\), but we make no attempt to justify that picture.
\begin{proposition}\label{thm:local.semibasic}
Take a surjective submersion \(P\xrightarrow{\pi}Q\).
Pick a semibasic form \(\xi\) on \(P\).
Then the following are equivalent:
\begin{itemize}
\item
The form \(d\xi\) is semibasic.
\item
The vertical infinitesimal automorphisms span the vertical vectors at every point of \(P\).
\item
The group \(\Aut[\xi/Q]^0\) acts transitively on every path component of every fiber of \(P\to Q\).
\item
The differential form \(\xi\) is locally basic.
To be precise: we can cover \(Q\) in open sets \(Q_a\) and \(P\) in open sets \(P_a\), so that \(P\to Q\) takes \(P_a\) onto \(Q_a\), and \(\xi\) is basic for \(P_a\to Q_a\), i.e. we can pick differential forms \(\eta_a\) on \(Q_a\) so that \(\xi=\pi^*\eta_a\) on \(P_a\).
\item
The same as the previous statement, but we can further choose these \(P_a\) so that each is \(\Aut[\xi/Q]\)-invariant, and hence intersects each fiber of \(P\to Q\) in a union of components of that fiber.
\end{itemize}
\end{proposition}
\begin{proof}
Take any vector field \(X\).
By Cartan's formula,
\[
\LieDer_X\xi=X\hook d\xi+d(X\hook \xi).
\]
Suppose that \(X\) is vertical.
If \(\xi\) and \(d\xi\) are semibasic then \(X\hook\xi=0\) and also \(X\hook d\xi=0\), so
\[
\LieDer_X\xi=0,
\]
i.e. \(\xi\) is invariant under vertical vector fields.
In particular, \(\xi\) is invariant under the flows of all complete vertical vector fields, which we can construct trivially in local coordinates on \(P\), by the implicit function theorem.
Hence the one parameter subgroups of \(\Aut[\xi/Q]\) act locally transitively on the fibers.
Hence their orbits are the components of the fibers.

Suppose that the vertical infinitesimal automorphisms span the vertical vectors.
Take any vertical infinitesimal automorphism \(X\) so
\[
0=\LieDer_X\xi=X\hook d\xi+d(X\hook \xi)=X\hook d\xi,
\]
since \(\xi\) is semibasic.
Hence \(d\xi\) is semibasic.

Suppose that the one parameter subgroups of \(\Aut[\xi/Q]\) act on every fiber of \(P\to Q\), with the subgroup they generate acting transitively on every path component of that fiber.
Their velocity vector fields are precisely the complete vertical infinitesimal automorphisms.
The largest family \(\bar{\mathfrak{F}}\) of complete vector fields whose flows preserve these orbits is precisely that same family.
Hence they form a Lie algebra and a module over the smooth functions and are invariant under one another's flows (pushing around) by theorem~\vref{thm:orbit}.
So their span is the same throughout their orbits, giving the tangent spaces to the orbits, i.e. they span the vertical vectors.

Take adapted coordinates \(x^i,y^j\).
Then
\[
\xi=a_I(x,y)dx^I,
\]
for unique smooth functions \(a_I(x,y)\), since \(\xi\) is semibasic.

Hence
\[
d\xi=
\partial_{x^i} a_I dx^{iI}
+
\partial_{y^j} a_I dy^j \wedge dx^I.
\]
Cancellations can occur in the first term, but the \(dx,dy\) are linearly independent, so \(d\xi\) is semibasic just when all of these vanish, i.e. just when
\[
\frac{\partial a_I}{\partial y^j}=0
\]
for all \(j\), i.e. just when \(\xi\) depends only on \(x\) in our coordinates.
Alternatively, we can say that \(\xi\) is invariant under the flows of vertical vector fields, hence under ``locally arbitrary'' motions in the \(y\) variables, so depends only on the \(x\) variables.
In other words, if \(\xi\) and \(d\xi\) are semibasic, then \(\Aut[\xi/Q]\) contains diffeomorphisms acting locally transitively on the fibers (as we already knew) and equivalently we can write \(\xi\) in \(x\) variables entirely, i.e. as a form on an open set in \(Q\).
Hence we can cover \(P\) in open sets \(P_a\subseteq P\) on each of which \(\xi=\pi^*\eta_a\) for some differential form \(\eta\) on the open set \(Q_a:=\pi P_a\).

Take one such \(P_a\) and some diffeomorphism \(\varphi\in\Aut[\xi/Q]\).
Then \(\varphi\) moves each point of \(P_a\) up its fiber of \(P\to Q\), acting trivially on \(\xi\) and on \(\eta_a\).
So on \(P_a\cup\varphi P_a\), \(\xi=\pi^*\eta_a\); we can replace \(P_a\) by \(P_a\cup\varphi P_a\).
By the same argument, we can replace \(P_a\) by
\[
\bigcup_{\varphi\in\Aut[\xi/Q]} \varphi P_a.
\]
\end{proof}
\begin{lemma}\label{lemma:G.inv.basic}
Take a surjective submersion \(P\xrightarrow{\pi}Q\).
Pick a semibasic form \(\xi\) on \(P\).
If some Lie group \(G\) acts on \(P\) preserving \(\xi\) and preserving every fiber of \(P\to Q\), and acts on each fiber with open orbits, then the group \(\Aut[\xi/Q]^0\) acts transitively on every path component of every fiber, and hence proposition~\vref{thm:local.semibasic} applies.
\end{lemma}
\begin{proof}
We can replace \(G\) by its identity component, so assume \(G\) is connected, so the one parameter subgroups of \(G\) generate the action and belong to \(\Aut[\xi/Q]^0\).
The orbits on any fiber are open, but then so are their complements, being unions of orbits.
\end{proof}
\begin{theorem}\label{theorem:basic}
Take a surjective submersion \(P\xrightarrow{\pi}Q\) and a semibasic form \(\xi\) on \(P\).
Suppose that, for any two components of any fiber of \(\pi\), some element of \(\Aut[\xi/Q]\) takes one of these components to the other.
Then the following are equivalent:
\begin{itemize}
\item
The form \(\xi\) is basic.
\item
The form \(\xi\) is \emph{uniquely} basic: there is a unique smooth differential form \(\eta\) on \(Q\) with \(\xi=\pi^*\eta\).
\item
The form \(d\xi\) is semibasic.
\item
The group \(\Aut[\xi/Q]^0\) acts on every fiber of \(P\to Q\), transitively on every path component of that fiber.
\item
The vertical infinitesimal automorphisms span the vertical vectors at every point of \(P\).
\end{itemize}
\end{theorem}
\begin{proof}
Suppose that \(\xi\) is basic, so \(\xi=\pi^*\eta\).
Then clearly \(\xi\) is semibasic and
\[
d\xi=d\pi^*\eta=\pi^*d\eta
\]
is semibasic too.
But \(\xi\) is invariant under any diffeomorphism preserving \(\pi\), i.e. \(\Aut[\xi/Q]\) is the set of all diffeomorphisms of \(P\) preserving \(\pi\), hence acts transitively on the components of the fibers.

By theorem~\vref{thm:local.semibasic}, \(d\xi\) is semibasic just when \(\Aut[\xi/Q]\) acts locally transitively on the fibers of \(P\to Q\), and this just when
the vertical infinitesimal automorphisms span the vertical vectors at every point of \(P\).

By hypothesis, for any two components of any fiber of \(\pi\), some element of \(\Aut[\xi/Q]\) takes one of these components to the other.
The inverse reverses them, so \(\Aut[\xi/Q]\) interchanges any two components of any fiber, so acts transitively on the components of each fiber.
Hence \(\Aut[\xi/Q]\) acts locally transitively on every fiber just when its acts transitively on every fiber.
We can assume that each \(P_a\) is \(\Aut[\xi/Q]\)-invariant, so a union of fibers of \(P\to Q\), i.e. \(P_a=\pi^{-1}Q_a\).
We have now covered \(Q\) in open sets \(Q_a\) on each of which we have a smooth differential form \(\eta_a\) so that \(\pi^*\eta_a=\xi\) on \(P_a=\pi^{-1}Q_a\).
Since \(\pi\) is a surjective submersion, \(\pi^*\) is injective.
But \(\pi^*\eta_a=\xi=\pi^*\eta_b\) over \(Q_a\cap Q_b\), hence \(\eta_a=\eta_b\) on \(Q_a\cap Q_b\), so there is a unique smooth form \(\eta\) on \(Q\) with \(\xi=\pi^*\eta\).
\end{proof}

\chapter{Invariant theory}\label{appendix:invariant.theory}
We aim to understand how to construct invariants of Lie group actions.
In many important examples, we can write down the Lie group action explicitly in polynomial expressions.
We will see that this ensures, roughly speaking, that there are enough invariant rational functions to separate most of the group orbits.
\begin{example}
The group \(G\) of matrices of the form
\[
\begin{pmatrix}
t&0\\
0&\frac{1}{t}
\end{pmatrix}
\]
acts on the plane \(\R^2\), taking each point \((x,y)\) and rescaling \(x\) by \(t\), \(y\) by \(1/t\).
The function \(xy\) is invariant.
The orbits:
\[
\includegraphics[width=4cm]{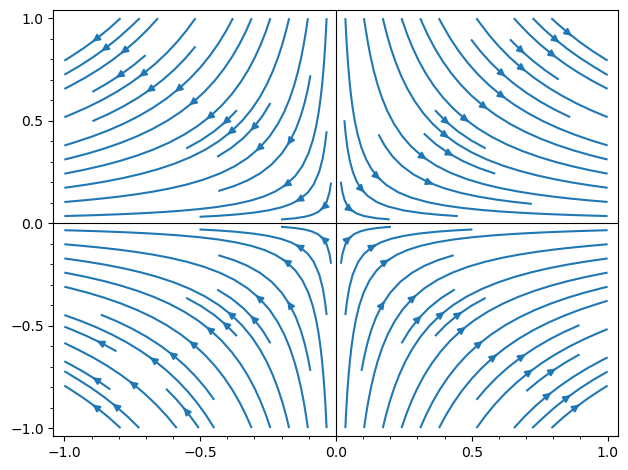}
\]
are 
\begin{itemize}
\item
the origin, 
\item
the horizontal axis punctured at the origin, 
\item 
the vertical axis punctured at the origin, 
\item
for any constant \(c\ne 0\), the hyperbola \(xy=c\).
\end{itemize}
So the invariant function \(xy\) ``separates orbits'', i.e. takes different values at points just when they lie on different orbits, except along the horizontal and vertical axes.

Taking small nonzero values of \(t\), we can take any point of the horizontal axis and bring it as close as we like to the origin.
So any continuous function invariant under the group action takes on the same value at every point of the horizontal axis as it does at the origin.
Take large values of \(t\), hence small values of \(1/t\), and the same argument applies: any continuous function invariant under the group action takes on the same value at every point of the horizontal axis and vertical axes as it does at the origin.
So there is no continuous invariant ``better'' than \(xy\), i.e. no continuous invariant can tell us which of these ``bad points'' of the two axes lie on which orbit.
We will want to cut out the bad points.
\end{example} 
An \emph{algebraic set}\define{algebraic set} is a set cut out by polynomial equations in finitely many variables, which we will assume are real or complex variables.
By the Hilbert basis theorem \cite{Kunz:2013} p. 11 proposition 2.3, any algebraic set cut out by infinitely many polynomial equations is already cut out by taking only finitely many of those polynomial equations.
Any algebraic set is a union of finitely many submanifolds, and has finitely many path components; see \cite{Whitney:1957} for an elementary proof.
Its \emph{dimension}\define{dimension} is the largest dimension of those submanifolds.
A \emph{regular map} of an algebraic set is a polynomial map restricted to that set \cite{Shafarevich:2013} p. 27.
A \emph{linear algebraic group}\define{linear algebraic group} is an algebraic set of matrices forming a group under matrix multipication \cite{Harris:1995} p. 114 lecture 10, \cite{Procesi:2007} chapter 7, \cite{OnishchikVinberg:1990}.
Every linear algebraic group is, by definition, a closed subgroup of the general linear group, and hence a closed embedded Lie subgroup \cite{Mimura/Toda:1991} p. 44.
\begin{example}
The orthogonal group, the symplectic group, the circle of unit complex numbers, the general and special linear groups, the upper triangular invertible matrices, the strictly upper triangular matrices, the unitary and special unitary groups, the automorphism group of any finite dimensional algebra or Lie algebra are linear algebraic groups.
Products of linear algebraic groups are linear algebraic groups.
\end{example}
\begin{example}
An elliptic curve is a complex algebraic group, but not isomorphic to any complex linear algebraic group; it is isomorphic as a real Lie group to the torus.
Densely winding subgroups of a torus are not linear algebraic.
The groups \(\GL{n}\R,\SO{n,1}\) are linear algebraic groups, but their identity components are not.
\end{example}
\begin{example}
By a theorem of Chevalley, every compact Lie group is a linear algebraic group \cite{Brocker.tom.Dieck:1995} III.8. 
For the curious: this holds for a unique real linear algebraic group each of whose complex components contains a real point and which contains no linear algebraic Zariski closed \(\R^{\times}\)-subgroup; moreover the linear algebraic group is semisimple just when the compact Lie group has finite center \cite{Brocker.tom.Dieck:1995} III.8.
From a different point of view, the classification of connected compact Lie groups makes fairly clear that they are linear algebraic \cite{Procesi:2007} chapter 10, section 7.2, theorem 4, page 380.
\end{example}
A \emph{regular morphism}\define{regular!morphism!of linear algebraic groups}\define{morphism!regular!of linear algebraic groups} of linear algebraic groups is a morphism of groups so that output matrix entries are polynomial in input matrix entries.
Take a linear algebraic group \(G\).
A \emph{regular \(G\)-module}\define{regular!module} is a finite dimensional vector space \(V\) with a regular morphism \(G\to\GL{V}\).
\begin{example}
Every finite dimensional representation of a compact Lie group is regular (for connected compact Lie groups, see \cite{Procesi:2007} p. 228 theorem 3; more generally see \cite{Brocker.tom.Dieck:1995} p. 149 theorem 7.5, \cite{casselman_2015}, \cite{Chevalley:1955.3} chapter VI). 
\end{example}
Given a set of functions invariant under a group action, each function perhaps only defined on an invariant open set, we say that these functions \emph{separate orbits}\define{separate orbits}\define{orbit!separate} of the action if any two points in different orbits have different value of at least one of the functions which has values defined at both points.
\begin{example}
Any densely winding subgroup of the torus is not algebraic.
Any continuous function on the torus, invariant under the subgroup, is constant.
\end{example}
\begin{problem}{ci}
Find invariant rational functions for the action of the Lorentz group on Minkowski space.
Do they separate orbits?
\end{problem}
A \emph{quasialgebraic subset}\define{quasialgebraic!subset}\define{subset!quasialgebraic} of an algebraic set \(X\) is a difference \(X-Y\) where \(Y\subset X\) is an algebraic set of lower dimension.
By definition, every quasialgebraic subset is an open subset.
A \emph{fat subset} of an algebraic set \(X\) is a union of topological components of a quasialgebraic subset.

Take a linear algebraic group \(G\), a regular \(G\)-module \(V\), and a \(G\)-invariant algebraic set \(Z\subseteq V\).
A \emph{good set} \(U\subseteq Z\) is a subset so that
\begin{itemize}
\item
\(U\) is \(G\)-invariant and
\item
\(U\subseteq Z\) is a quasialgebraic subset, hence an open subset, and
\item
\(U\subseteq V\) is a real analytic submanifold and
\item
there is a finite set of rational functions on \(V\), \(G\)-invariant on \(Z\), which separate the \(G\)-orbits on \(U\) and
\item
the quotient map \(U\to\bar{U}:=G\backslash U\) is a real analytic smooth fibration of real analytic manifolds and
\item
\(\bar{U}\) is a fat subset of some algebraic set \(\bar{Z}\), hence an open subset, so that
\item
the quotient map \(U\to\bar{U}\) composes with the inclusion \(\bar{U}\to\bar{Z}\) to give an \(G\)-invariant regular map \(U\to\bar{Z}\).
\end{itemize}
The set \(Z-U\) is a \emph{bad set}.\define{bad set}
\begin{theorem}%
[Rosenlicht \cite{Rosenlicht:1963}, 
\cite{Grosshans:1997} p. 107 theorem 19.5, 
\cite{Shafarevich:1989} p. 155, theorem 2.3]%
\label{theorem:Rosenlicht}%
\define{theorem!Rosenlicht}\define{Rosenlicht theorem}
Take a linear algebraic group \(G\) and a regular \(G\)-module \(V\).
Every \(G\)-invariant algebraic set \(Z\subseteq V\) contains a good set.
\end{theorem}
\begin{example}
The bad set \(Z':=Z-U\) is then also a \(G\)-invariant algebraic set.
Applying this theorem to \(Z'\) in place of \(Z\), we produce a good set for \(Z'\), say \(U'\subseteq Z'\).
Invariant rational functions on \(Z'\) separate orbits on \(U'\), and so on by induction.
\end{example}
\begin{example}
The group \(G\) of matrices of the form
\[
\begin{pmatrix}
a&0\\
0&b
\end{pmatrix}
\]
with \(ab=1\) acts on \(V=\R^2\); consider the algebraic set \(Z=V\).
The orbits:
\[
\includegraphics[width=4cm]{hyperbolic-orbits}
\]
are ``mostly'' cut out by the invariant rational function \(f(x,y)=xy\), except on the bad set 
\[
Z'=\set{(x,y)|xy=0}=(x=0)\cup(y=0),
\]
which has invariant rational functions \(f_1(x,y)=x/y\) and \(f_2(x,y)=y/x\), and bad set 
\[
Z''=\set{(0,0)},
\]
which has invariant rational function \(f_3(x,y)=0\) and empty bad set.
\end{example}
\begin{example}
The group \(G:=\set{\pm 1}\) acts by sign change on the real number line \(Z=V=\R\).
We can take the bad set to be \(Z':=\set{0}\), so \(U:=\R^{\times}\) has quotient space \(\bar{U}:=\R^+\).
We can take \(\bar{Z}:=\R\) and \(\bar{Z}':=\set{0}\), so \(\bar{U}\subset\bar{Z}-\bar{Z}'\) is \(\R^+\subset\R^{\times}\), a single topological component.
The quotient map is
\[
x\in U\mapsto x^2\in\R^+=\bar{U}\subset\bar{Z}=\R,
\]
a regular map.
\end{example}
\begin{example}
The rescalings \(G=\R^{\times}I\) act on \(V=\R^n\), with all orbits approaching the origin.
So \(Z=V\) has invariant rational functions \(x_i/x_j\).
We can take \(Z':=\set{0}\) with rational function \(0\), \(Z''\) empty.
Note that at some points of the good set \(U=\R^n-\set{0}\), some of these rational functions might not be defined.
We can still separate orbits: any distinct orbits have enough invariant rational functions defined near both of them so that we can separate them.
\end{example}
\begin{example}
For any compact Lie group \(G\) and finite dimensional \(G\)-module \(V\), all invariant rational functions can be expressed as functions of invariant polynomials \cite{Procesi:2007}, so we don't have to worry about denominators vanishing, but we still encounter bad sets.
\end{example}
\section{Proof of Rosenlicht's theorem}
Rosenlicht's original proof \cite{Rosenlicht:1963} uses big machinery from algebraic geometry, and Gromov \cite{Gromov1988} p. 89 only gives an outline of a proof, so it seems worthwhile to expand Gromov's outline into a complete proof using elementary algebraic geometry.
Let us briefly review some elementary algebraic geometry, only over the real or complex numbers, for which we recommend \cite{Harris:1995,Shafarevich:2013}; these are only references we will depend on.
An algebraic set is \emph{irreducible}\define{algebraic set!irreducible}\define{irreducible algebraic set} if not the union of two proper algebraic subsets.
An irreducible algebraic set is also called an \emph{affine variety}\define{affine variety}\define{variety!affine} \cite{Reid:1988}.
By the Noetherian property of polynomials, every algebraic set is a unique union of affine varieties.
Every affine variety contains a smooth point, i.e. a point near which the variety is cut out by polynomials with linearly independent differentials \cite{Shafarevich:2013} p. 93.
An \emph{algebraic set} in projective space is a set which intersects each affine subspace (i.e. complement of a projective hypersurface) in an algebraic set.
Repeating the definition of irreducible, an irreducible algebraic set in projective space is a \emph{projective variety}.\define{variety!projective}\define{projective variety}
Every affine variety has closure in projective space a projective variety \cite{Shafarevich:2013} p. 45.
Every projective variety is compact, because it is a closed subset of projective space.
Hence the image of any projective subvariety under any morphism of varieties is compact; in fact, it is easy to see that it is a projective subvariety \cite{Shafarevich:2013} p. 57.
Any affine variety is its projective closure minus the set of points at infinity, which is a projective subvariety of lower dimension.
So under any morphism of varieties, the image of any affine variety is the image of its projective closure minus the image of the set of points at infinity.
Hence the image of any morphism of affine algebraic sets is a difference of affine algebraic sets and is Zariski dense in its Zariski closure; see \cite{Borel:1991} p. 21 for proof over any field.
In particular, every orbit of a regular action of a linear algebraic group is a difference of a variety minus an affine algebraic subset of lower dimension.
The closure of an orbit is thus a union of the orbit with various lower dimensional orbits;  again \cite{Borel:1991} p. 53 proves this over any field, but we stick with the real or complex numbers.
For any Lie group action, every orbit is a smooth submanifold.
Every orbit of a linear algebraic group acting on an algebraic set is a smooth embedded submanifold whose boundary is thus a union of lower dimensional orbits.
In particular, minimal dimensional orbits are smooth embedded submanifolds; see \cite{Humphreys:1975} p. 60 section 8.3 for proof, over any field, that they are smooth varieties.

More generally, given a linear algebraic group \(G\) and a regular \(G\)-module \(V\), and a subvariety \(S\subseteq V\), the union of the \(G\)-orbits through \(S\) is the image of
\[
(g,s)\in G\times S\mapsto gs\in V,
\]
so is a difference of \(G\)-invariant subvarieties.

For any affine algebraic set (or affine algebraic variety) \(X\), the symmetric product \(S^n X\) is an affine algebraic set (affine algebraic variety), whose points are the unordered \(n\)-tuples \([x_1,\dots,x_n]\) of points of \(X\) (with multiplicities), with coordinate ring the polynomials invariant under permutation of points \cite{Artin:2022} p.~80, section 2.7, \cite{Harris:1995} section 10.23.

We prove Rosenlicht's theorem (theorem~\vref{theorem:Rosenlicht}), following the proof outlined by Gromov \cite{Gromov1988} p. 89:
\begin{proof}
Any affine algebraic set has unique decomposition into a union of affine varieties.
So the decomposition of \(Z\) into affine varieties is \(G\)-invariant.
If Rosenlicht's theorem applies to each of these affine varieties, it applies to \(Z\).
So it suffices to prove Rosenlicht's theorem for \(Z\) irreducible, i.e. an affine variety.
The generic point of \(Z\) is smooth, with maximal local dimension.
The singular points of \(Z\) are \(G\)-invariant, and form an algebraic subvariety \(Z'\subset Z\) of lower dimension.
The Lie algebra of \(G\) is realized as vector fields on the smooth points of \(Z\), of locally maximal rank at the generic point.
Since \(Z\) is irreducible, this locally maximal rank is constant away from some \(G\)-invariant subvariety, which we can put into \(Z'\) too.
So on \(U:=Z-Z'\), every point has smallest dimensional \(G\)-stabilizer, so maximal dimensional \(G\)-orbit.
As discussed above, each \(G\)-orbit is a difference of \(G\)-invariant affine algebraic sets \(A-B\), so open in its closure \(A\); its boundary \(B\) is therefore a union of lower dimensional \(G\)-orbits.
The subset of \(Z\) on which the \(G\)-orbits are of lower dimension lies in \(Z'\), so all \(G\)-orbits in \(U\) are Zariski closed in \(U\) and of the same dimension.

Take a \(G\)-orbit \(\mathscr{O}:=Gz_0\subseteq U\) through a point \(z_0\in U\).
Take an affine subspace \(A\subseteq V\) of complementary dimension to \(\mathscr{O}\).
By B\'ezout's theorem \cite{Shafarevich:2013} p. 236, there are finitely many intersection points of \(\mathscr{O}\) and \(A\).
If we replace real variables by complex variables, complexifying the varieties involved, then that number of intersection points is the same for any two generics choices of \(A\).
Therefore, for the moment, let us prove the Rosenlicht theorem for complex algebraic varieties.
We will reconsider the real case later on.

For generic choice of \(A\), \(A\) and \(\mathscr{O}\) don't intersect at infinity, or on the boundary of \(\mathscr{O}\), so the expected number of intersections occur on complex points of \(\mathscr{O}\).
By Sard's theorem \cite{Sard:1965}, for generic choice of \(A\), the points at which \(T:=Z\cap A\) strikes \(\mathscr{O}\) are transverse intersections, and remain so through all nearby orbits.
Loss of transversality is an algebraic condition, so transversality of \(T\) to the \(G\)-orbits fails on an algebraic subset \(T'\subset T\) of lower dimension.
As above, \(GT'\) is a difference of \(G\)-invariant affine algebraic sets.
Moreover, \(GT'\) does not intersect any \(G\)-orbit near \(\mathscr{O}\).
So \(GT'\) lies in \(G\)-invariant algebraic subset of \(Z\) of lower dimension than \(Z\).
Hence we can put \(GT'\) into \(Z'\).
So \(T\) intersects every \(G\)-orbit in \(U\) transversely.

Similarly, the set of points of \(Z\) lying on \(G\)-orbits which intersect \(T\) outside of \(U\) is closed and nowhere dense, so lies in an algebraic subset, which we can place into \(Z'\).
So \(T\) intersects every \(G\)-orbit in \(U\) transversely and always at the same number of complex points.

Write down the complex intersection points of \(T\) and \(\mathscr{O}\) in some order.
For nearby orbits, the intersection points move only a little, by the implicit function theorem.
So over the complex numbers, we have a holomorphic map define near each point of \(T\), from an open subset of \(T\) to \(Z\), mapping each point of \(T\) to that same point but thought of as a point of \(Z\), and this map is tranverse to orbits, where it is defined.
Suppose that there are \(n\) intersection points.

Let \(\hat{Z}_1\) be the set of tuples \((z,g)\) where \(z\in Z\) and \(g\in S^n G\) 
\[
g=[g_1,\dots,g_n]
\]
so that \(g_iz\in T\).
Clearly \(\hat{Z}_1\) is an algebraic subset of \(Z\times S^n G\).

Let \(\hat{Z}\subset\hat{Z}_1\) be the algebraic closure of the Zariski open set of points at which all \(g_i\) are distinct.
For \(z\in U\), the \(G\)-equivariant map
\[
(z,g)\in\hat{Z}\xmapsto{z}z\in Z
\]
is bijective, so an local isomorphism of complex manifolds.
Note that \(\hat{Z}_1\subseteq Z\times S^nG\) has Zariski closure in a suitable projective space sitting inside the product of the projective closures of \(Z\) and \(S^n G\).
In particular, \(\hat{Z}_1\xmapsto{z} Z\) extends to the projective closures.
By finiteness of the generic fiber, \(\hat{Z}\to Z\) is a finite map in the sense of algebraic geometry \cite{Harris:1995} p. 178 lemma 14.8.
By the implicit function theorem for finite maps of algebraic varieties \cite{Harris:1995} p. 179 corollary 14.10, the map \(z\) restricts to \(U\) to become an isomorphism of algebraic varieties.

Consider the map
\[
(z,g)\in\hat{Z}_1\xmapsto{t}t=[g_1z,\dots,g_nz]\in S^nT,
\]
which is invariant under the \(G\)-action 
\[
g_0(z,g):=
(g_0z,[g_1g_0^{-1},\dots,g_ng_0^{-1}]).
\]
Over \(U\), this map drops to a holomorphic \(G\)-invariant map
\[
U\xrightarrow{t}S^nT.
\]
Different orbits are disjoint so intersect \(T\) at disjoint sets of points, hence distinct points of \(S^n T\).
So this map separates points of \(U\).

The rational functions on \(S^n T\) pull back to \(\hat{Z}\), hence defined on \(Z\), so arise from rational functions on \(V\).
They are \(G\)-invariant on \(\hat{Z}\) so also on \(Z\).

We can put into \(Z'\) any points near which \(U\xrightarrow{t}S^n T\) is not of locally maximal rank differential.
So we can assume that \(t\) has constant rank.
The image \(\bar{U}\) of \(t\) is, once again, a difference \(\bar{Z}-\bar{Z}'\) of affine varieties \(\bar{Z},\bar{Z}'\), with \(\bar{Z}'\subseteq\bar{Z}\) of lower dimension.
We can assume that \(\bar{Z}\) is the Zariski closure of \(\bar{U}\).
We can put the preimage of the singular locus of \(\bar{Z}\) into \(Z'\), so assume that \(\bar{Z}\) is a smooth affine algebraic variety away from \(\bar{Z}'\).
Our map \(U\xrightarrow{t}\bar{U}\) is of constant rank and surjective, so a submersion.
Since \(T\) intersects every orbit of \(U\) transversely, some Euclidean open subset of \(T\) intersects transversely and at a single point, a local section, and these local sections strike every orbit in \(U\).
So \(U\xrightarrow{t}\bar{U}\) is a holomorphic submersion with local holomorphic sections, hence a holomorphic fibration.

Finally, let us return to the original problem, over the real numbers.
We first complexify and follow the argument above.
At each step, when we construct \(Z'\) above, we use only real equations, so \(Z'\) is also a real algebraic variety.
Choose our affine subspace \(A\) to be real. 
The number of complex intersection points in B\'ezout's theorem is constant. 
The number of real intersection points might not be constant, but only locally constant, since pairs of complex conjugate intersection points can collide to produce a pair of real intersection points as we move, but the collisions are singular, so lie in \(Z'\).
In particular, if there are no real intersection points of \(T\) above a real point of \(\bar{Z}-\bar{Z}'\), then \(\bar{U}\)  does not intersect that component of \(\bar{Z}-\bar{Z}'\).
\end{proof}

\end{appendices}

\providecommand{\bysame}{\leavevmode\hbox to3em{\hrulefill}\thinspace}
\providecommand{\MR}{\relax\ifhmode\unskip\space\fi MR }
% \MRhref is called by the amsart/book/proc definition of \MR.
\providecommand{\MRhref}[2]{%
  \href{http://www.ams.org/mathscinet-getitem?mr=#1}{#2}
}
\providecommand{\href}[2]{#2}
\begin{thebibliography}{100}

\bibitem{Abraham/Marsden:1978}
Ralph Abraham and Jerrold~E. Marsden, \emph{Foundations of mechanics},
  Benjamin/Cummings Publishing Co. Inc. Advanced Book Program, Reading, Mass.,
  1978, Second edition, revised and enlarged, With the assistance of Tudor Ra\c
  tiu and Richard Cushman. \MR{81e:58025}

\bibitem{Allcock:2018}
Daniel Allcock, \emph{Spherical space forms revisited}, Trans. Amer. Math. Soc.
  \textbf{370} (2018), no.~8, 5561--5582. \MR{3803143}

\bibitem{Allcock.Carlson.Toledo:2011}
Daniel Allcock, James~A. Carlson, and Domingo Toledo, \emph{The moduli space of
  cubic threefolds as a ball quotient}, Mem. Amer. Math. Soc. \textbf{209}
  (2011), no.~985, xii+70. \MR{2789835}

\bibitem{Ambrose.Singer:1953}
W.~Ambrose and I.~M. Singer, \emph{A theorem on holonomy}, Trans. Amer. Math.
  Soc. \textbf{75} (1953), 428--443. \MR{63739}

\bibitem{Amores:1979}
A.~M. Amores, \emph{Vector fields of a finite type {$G$}-structure}, J.
  Differential Geom. \textbf{14} (1979), no.~1, 1--6 (1980),
  \url{http://projecteuclid.org/euclid.jdg/1214434847}. \MR{577874}

\bibitem{An:2012}
Jinpeng An, \emph{Rigid geometric structures, isometric actions, and algebraic
  quotients}, Geom. Dedicata \textbf{157} (2012), 153--185,
  \url{https://doi-org.ucc.idm.oclc.org/10.1007/s10711-011-9603-2}.
  \MR{2893482}

\bibitem{Arbarello.Cornalba.Griffiths.Harris:1985}
E.~Arbarello, M.~Cornalba, P.~A. Griffiths, and J.~Harris, \emph{Geometry of
  algebraic curves. {V}ol. {I}}, Grundlehren der mathematischen Wissenschaften
  [Fundamental Principles of Mathematical Sciences], vol. 267, Springer-Verlag,
  New York, 1985. \MR{770932}

\bibitem{Arnold:1989}
V.~I. Arnolʹd, \emph{Mathematical methods of classical mechanics}, second ed.,
  Graduate Texts in Mathematics, vol.~60, Springer-Verlag, New York, 1989,
  Translated from the Russian by K. Vogtmann and A. Weinstein. \MR{997295}

\bibitem{Artin:2022}
Michael Artin, \emph{Algebraic geometry---notes on a course}, Graduate Studies
  in Mathematics, vol. 222, American Mathematical Society, Providence, RI,
  [2022] \copyright 2022. \MR{4465876}

\bibitem{Atiyah:1957}
M.~F. Atiyah, \emph{Complex analytic connections in fibre bundles}, Trans.
  Amer. Math. Soc. \textbf{85} (1957), 181--207. \MR{MR0086359 (19,172c)}

\bibitem{Barakat:2004}
Mohamed Barakat, \emph{The existence of {C}artan connections and geometrizable
  principle bundles}, Arch. Math. (Basel) \textbf{83} (2004), no.~2, 159--163,
  \url{https://doi-org.ucc.idm.oclc.org/10.1007/s00013-004-1048-5}.
  \MR{2104944}

\bibitem{Baston/Eastwood:1989}
Robert~J. Baston and Michael~G. Eastwood, \emph{The {P}enrose transform},
  Oxford Mathematical Monographs, The Clarendon Press Oxford University Press,
  New York, 1989, Its interaction with representation theory, Oxford Science
  Publications. \MR{MR1038279 (92j:32112)}

\bibitem{Bavard.Mounoud:2013}
Christophe Bavard and Pierre Mounoud, \emph{Sur les surfaces lorentziennes
  compactes sans points conjugu\'{e}s}, Geom. Topol. \textbf{17} (2013), no.~1,
  469--492. \MR{3039767}

\bibitem{Beem.Ehrlich.Easley:1996}
John~K. Beem, Paul~E. Ehrlich, and Kevin~L. Easley, \emph{Global {L}orentzian
  geometry}, second ed., Monographs and Textbooks in Pure and Applied
  Mathematics, vol. 202, Marcel Dekker, Inc., New York, 1996. \MR{1384756}

\bibitem{Benoist:1997}
Yves Benoist, \emph{Orbites des structures rigides (d'apr\`es {M}. {G}romov)},
  Integrable systems and foliations/{F}euilletages et syst\`emes
  int\'{e}grables ({M}ontpellier, 1995), Progr. Math., vol. 145, Birkh\"{a}user
  Boston, Boston, MA, 1997,
  \url{https://doi-org.ucc.idm.oclc.org/10.1007/978-1-4612-4134-8_1},
  pp.~1--17. \MR{1432904}

\bibitem{Blumenthal1986}
R.~A. Blumenthal, \emph{Mappings between manifolds with {C}artan connections},
  Differential geometry, {{P}e\~{n}\'{\i}scola} 1985 \cite{Penis.Cola:1986},
  Lecture Notes in Math., vol. 1209, Springer, Berlin, 1986,
  \url{https://doi-org.ucc.idm.oclc.org/10.1007/BFb0076622}, pp.~94--99.
  \MR{863747}

\bibitem{Borel:1963}
Armand Borel, \emph{Compact {C}lifford-{K}lein forms of symmetric spaces},
  Topology \textbf{2} (1963), 111--122,
  \url{https://doi-org.ucc.idm.oclc.org/10.1016/0040-9383(63)90026-0}.
  \MR{146301}

\bibitem{Borel:1991}
\bysame, \emph{Linear algebraic groups}, second ed., Graduate Texts in
  Mathematics, vol. 126, Springer-Verlag, New York, 1991. \MR{92d:20001}

\bibitem{Brocker.tom.Dieck:1995}
Theodor Br\"{o}cker and Tammo tom Dieck, \emph{Representations of compact {L}ie
  groups}, Graduate Texts in Mathematics, vol.~98, Springer-Verlag, New York,
  1995, Translated from the German manuscript, Corrected reprint of the 1985
  translation. \MR{1410059}

\bibitem{Bryant:1991}
Robert~L. Bryant, \emph{An introduction to {L}ie groups and symplectic
  geometry}, Geometry and quantum field theory ({P}ark {C}ity, {UT}, 1991),
  IAS/Park City Math. Ser., vol.~1, Amer. Math. Soc., Providence, RI, 1995,
  \url{https://doi-org.ucc.idm.oclc.org/10.1090/pcms/001/02}, pp.~5--181.
  \MR{1338391}

\bibitem{BCGGG:1991}
Robert~L. Bryant, Shiing-Shen Chern, R.~B. Gardner, H.~L. Goldschmidt, and
  P.~A. Griffiths, \emph{Exterior differential systems}, Springer-Verlag, New
  York, 1991. \MR{92h:58007}

\bibitem{Bryant.Griffiths:1995}
Robert~L. Bryant and Phillip~A. Griffiths, \emph{Characteristic cohomology of
  differential systems. {II}. {C}onservation laws for a class of parabolic
  equations}, Duke Math. J. \textbf{78} (1995), no.~3, 531--676. \MR{1334205}

\bibitem{Calabi.Markus:1962}
E.~Calabi and L.~Markus, \emph{Relativistic space forms}, Ann. of Math. (2)
  \textbf{75} (1962), 63--76. \MR{133789}

\bibitem{Canary.Epstein.Green:2006}
R.~D. Canary, D.~B.~A. Epstein, and P.~L. Green, \emph{Notes on notes of
  {T}hurston}, Fundamentals of Hyperbolic Manifolds: Selected Expositions,
  London Mathematical Society Lecture Note Series, vol. 238, London Math. Soc.,
  London, UK, 2006, pp.~1--115.

\bibitem{Cap/Slovak:2009}
Andreas {\v{C}}ap and Jan Slov{\'a}k, \emph{Parabolic geometries. {I}},
  Mathematical Surveys and Monographs, vol. 154, American Mathematical Society,
  Providence, RI, 2009, Background and general theory. \MR{2532439
  (2010j:53037)}

\bibitem{Cartan:1924}
E.~Cartan, \emph{Sur les vari\'{e}t\'{e}s \`a connexion projective}, Bull. Soc.
  Math. France \textbf{52} (1924), 205--241. \MR{1504846}

\bibitem{Cartan:1926}
{\'E}.~Cartan, \emph{Les groupes d'holonomie des espaces
  g{\'e}n{\'e}ralis{\'e}s.}, Acta Math. \textbf{48} (1926), 1--42 (French).

\bibitem{Cartan:1910}
Elie Cartan, \emph{Les syst\`emes de {P}faff, \`a cinq variables et les
  \'{e}quations aux d\'{e}riv\'{e}es partielles du second ordre}, Ann. Sci.
  \'{E}cole Norm. Sup. (3) \textbf{27} (1910), 109--192, Also in
  \cite{Cartan:II}, pp. 927--1010. \MR{1509120}

\bibitem{Cartan:68}
\'Elie Cartan, \emph{Les espaces {\`a} connexion conforme}, Ann. Soc. Polon.
  Mat. \textbf{2} (1923), 171--221, Also in \cite{Cartan:III2}, pp. 747--797.

\bibitem{Cartan1930}
{\'E}lie Cartan, \emph{La th{\'e}orie des groupes finis et continus et
  l'analysis situs.}, M{\'e}morial des sciences math{\'e}matiques. {Fasc}. 42.
  {Paris}: {Gauthier}-{Villars}. 62 p. (1930)., 1930.

\bibitem{Cartan:136bis}
\'Elie Cartan, \emph{Sur la g\'eometrie pseudo--conforme des hypersurfaces de
  l'espace de deux variables complexes, {II}}, Annali Sc. Norm. Sup. Pisa
  \textbf{1} (1932), 333--354, Also in \cite{Cartan:III2}, pp. 1217--1238.

\bibitem{Cartan:136}
Elie Cartan, \emph{Sur la g\'{e}om\'{e}trie pseudo-conforme des hypersurfaces
  de l'espace de deux variables complexes}, Ann. Mat. Pura Appl. \textbf{11}
  (1933), no.~1, 17--90, Also in \cite{Cartan:II}, pp. 1232--1305. \MR{1553196}

\bibitem{Cartan:161}
\'Elie Cartan, \emph{Les probl{\`e}mes de equivalence}, S{\'e}minaire de Math.
  \textbf{expos{\'e} D} (1937), 113--136, Also in \cite{Cartan:II}, pp.
  1311--1334.

\bibitem{Cartan:1938}
\bysame, \emph{Les espaces g\'en\'eralis\'es et l'integration de certaines
  classes d'\'equations differentielles}, C. R. Acad. Sci. Paris S\'er. I Math.
  (1938), no.~206, 1689--1693, also in {\OE}uvres {C}ompletes, vol. III, partie
  2, pp. 1621-1636.

\bibitem{Cartan:174}
\bysame, \emph{La geometria de las ecuaciones diferencials de tercer orden},
  Rev. Mat. Hispano-Amer. \textbf{4} (1941), 1--31, also in {\OE}uvres
  {C}ompl{\`e}tes, Partie III, Vol. 2, 174, p. 1535--1566.

\bibitem{Cartan:II}
Elie Cartan, \emph{\oe uvres compl\`etes. {P}artie {II}. {V}ol. 1. {A}lg\`ebre,
  formes diff\'{e}rentielles, syst\`emes diff\'{e}rentiels. {V}ol. 2. {G}roupes
  infinis, syst\`emes diff\'{e}rentiels, th\'{e}ories d'\'{e}quivalence},
  Gauthier-Villars, Paris, 1953. \MR{58523}

\bibitem{Cartan:III2}
\'{E}lie Cartan, \emph{\oe uvres compl\`etes. {P}artie {III}. {V}ol. 2}, second
  ed., \'{E}ditions du Centre National de la Recherche Scientifique (CNRS),
  Paris, 1984, G\'{e}om\'{e}trie diff\'{e}rentielle. Divers. [Differential
  geometry. Miscellanea], With biographical material by Shiing Shen Chern,
  Claude Chevalley and J. H. C. Whitehead. \MR{753093}

\bibitem{Cartan:1992}
{\'E}lie Cartan, \emph{Le\c cons sur la g\'eom\'etrie projective complexe. {L}a
  th\'eorie des groupes finis et continus et la g\'eom\'etrie diff\'erentielle
  trait\'ees par la m\'ethode du rep\`ere mobile. {L}e\c cons sur la th\'eorie
  des espaces \`a connexion projective}, Les Grands Classiques
  Gauthier-Villars. [Gauthier-Villars Great Classics], \'Editions Jacques
  Gabay, Sceaux, 1992, Reprint of the editions of 1931, 1937 and 1937.
  \MR{1190006 (93i:01030)}

\bibitem{casselman_2015}
William Casselman, \emph{Compact groups as algebraic groups}, Aug 2015,
  \url{https://personal.math.ubc.ca/~cass/research/pdf/Compact-algebraic.pdf}.

\bibitem{Chavel:2006}
Isaac Chavel, \emph{Riemannian geometry}, second ed., Cambridge Studies in
  Advanced Mathematics, vol.~98, Cambridge University Press, Cambridge, 2006, A
  modern introduction. \MR{2229062}

\bibitem{Chern:1989}
Shiing-shen Chern (ed.), \emph{Global differential geometry}, MAA Studies in
  Mathematics, vol.~27, Mathematical Association of America, Washington, DC,
  1989. \MR{90d:53003}

\bibitem{Chevalley:1955.3}
Claude Chevalley, \emph{Th\'eorie des groupes de {L}ie. {T}ome {III}.
  {T}h\'eor\`emes g\'en\'eraux sur les alg\`ebres de {L}ie}, Actualit\'es Sci.
  Ind. no. 1226, Hermann \& Cie, Paris, 1955. \MR{0068552}

\bibitem{Clarke:1979}
C.~J.~S. Clarke, \emph{Boundary definitions}, Gen. Relativity Gravitation
  \textbf{10} (1979), no.~12, 977--980,
  \url{https://doi-org.ucc.idm.oclc.org/10.1007/BF00776514}. \MR{550360}

\bibitem{Clifton:1966}
Yeaton~H. Clifton, \emph{On the completeness of {C}artan connections}, J. Math.
  Mech. \textbf{16} (1966), 569--576. \MR{MR0205183 (34 \#5017)}

\bibitem{Connelly:1971}
Robert Connelly, \emph{A new proof of {B}rown's collaring theorem}, Proceedings
  of the American Mathematical Society \textbf{27} (1971), no.~1, 180--182.

\bibitem{Conrad:2024}
Keith Conrad, \emph{${SL}_2\mathbb{Z}$}, unpublished, April 2024.

\bibitem{dAmbra/Gromov:1991}
G.~D'Ambra and M.~Gromov, \emph{Lectures on transformation groups: geometry and
  dynamics}, Surveys in differential geometry, Lehigh Univ., Bethlehem, PA,
  1991, From lectures in Cambridge MA, 1990, pp.~19--111. \MR{1144526
  (93d:58117)}

\bibitem{Saint-Gervais:2010}
Henri~Paul de~Saint-Gervais, \emph{Uniformisation des surfaces de {R}iemann},
  ENS \'{E}ditions, Lyon, 2010.

\bibitem{Demailly2012}
J.-P. Demailly, \emph{Complex analytic and differential geometry}, unpublished,
  June 2012.

\bibitem{Dieudonne:1949}
Jean Dieudonn\'{e}, \emph{Sur une g\'{e}n\'{e}ralisation du groupe orthogonal
  \`a quatre variables}, Arch. Math. \textbf{1} (1949), 282--287,
  \url{https://doi-org.ucc.idm.oclc.org/10.1007/BF00776514}. \MR{29360}

\bibitem{Duistermaat.Kolk:2000}
J.~J. Duistermaat and J.~A.~C. Kolk, \emph{Lie groups}, Universitext,
  Springer-Verlag, Berlin, 2000,
  \url{https://doi-org.ucc.idm.oclc.org/10.1007/978-3-642-56936-4}.
  \MR{1738431}

\bibitem{Dumitrescu:2014}
Sorin Dumitrescu, \emph{An invitation to quasihomogeneous rigid geometric
  structures}, Bridging algebra, geometry, and topology, Springer Proc. Math.
  Stat., vol.~96, Springer, Cham, 2014, pp.~107--123. \MR{3297111}

\bibitem{Earle:1981}
Clifford~J. Earle, \emph{On variation of projective structures}, Riemann
  surfaces and related topics: {P}roceedings of the 1978 {S}tony {B}rook
  {C}onference ({S}tate {U}niv. {N}ew {Y}ork, {S}tony {B}rook, {N}.{Y}., 1978),
  Ann. of Math. Stud., vol.~97, Princeton Univ. Press, Princeton, N.J., 1981,
  pp.~87--99. \MR{624807}

\bibitem{Ehresmann:1936}
Charles Ehresmann, \emph{Sur la notion d'espace complet en g{\'e}om{\'e}trie
  diff{\'e}rentielle}, C. R. Acad. Sci. Paris \textbf{202} (1936), 2033.

\bibitem{Ehresmann:1938}
\bysame, \emph{Sur les arcs analytique d'un espace de {C}artan}, C. R. Acad.
  Sci. Paris (1938), 1433.

\bibitem{Ehresmann:1951}
\bysame, \emph{Les connexions infinit\'esimales dans un espace fibr\'e
  diff\'erentiable}, Colloque de topologie (espaces fibr\'es), Bruxelles, 1950,
  Georges Thone, Li\`ege, 1951, pp.~29--55. \MR{MR0042768 (13,159e)}

\bibitem{Ehresmann1961}
\bysame, \emph{Structures feuillet{\'e}es}, Proc. Vth Can. Congress (Toronto)
  (E.~M. Rosenthall, ed.), Univ. of Toronto Press, 1961.

\bibitem{Ehrlich:2006}
Paul~E. Ehrlich, \emph{A personal perspective on global {L}orentzian geometry},
  Analytical and numerical approaches to mathematical relativity, Lecture Notes
  in Phys., vol. 692, Springer, Berlin, 2006, pp.~3--34. \MR{2222545}

\bibitem{Erickson}
Jacob~W. Erickson, \emph{A visual introduction to {C}artan geometries},
  \url{www.math.umd.edu/~jwericks}.

\bibitem{Fefferman:2009}
Charles Fefferman, \emph{Whitney's extension problems and interpolation of
  data}, Bull. Amer. Math. Soc. (N.S.) \textbf{46} (2009), no.~2, 207--220.
  \MR{2476412}

\bibitem{Fefferman.Israel:2020}
Charles Fefferman and Arie Israel, \emph{Fitting smooth functions to data},
  CBMS Regional Conference Series in Mathematics, vol. 135, American
  Mathematical Society, Providence, RI, [2020] \copyright 2020. \MR{4235100}

\bibitem{Feres.Lampe:2000}
R.~Feres and P.~Lampe, \emph{Cartan geometries and dynamics}, Geom. Dedicata
  \textbf{80} (2000), no.~1-3, 29--41,
  \url{https://doi-org.ucc.idm.oclc.org/10.1023/A:1005219805069}. \MR{1762497}

\bibitem{Feres:2002}
Renato Feres, \emph{Rigid geometric structures and actions of semisimple {L}ie
  groups}, Rigidit\'{e}, groupe fondamental et dynamique, Panor. Synth\`eses,
  vol.~13, Soc. Math. France, Paris, 2002, pp.~121--167. \MR{1993149}

\bibitem{Frances:2012}
Charles Frances, \emph{About geometrically maximal manifolds}, J. Topol.
  \textbf{5} (2012), no.~2, 293--322,
  \url{https://doi-org.ucc.idm.oclc.org/10.1112/jtopol/jts003}. \MR{2928078}

\bibitem{Frances:2014}
\bysame, \emph{Removable and essential singular sets for higher dimensional
  conformal maps}, Comment. Math. Helv. \textbf{89} (2014), no.~2, 405--441,
  \url{https://doi-org.ucc.idm.oclc.org/10.4171/CMH/323}. \MR{3225453}

\bibitem{Frances.Melnick:2019}
Charles Frances and Karin Melnick, \emph{Topology of automorphism groups of
  parabolic geometries}, Geom. Topol. \textbf{23} (2019), no.~1, 135--169,
  \url{https://doi-org.ucc.idm.oclc.org/10.2140/gt.2019.23.135}. \MR{3921318}

\bibitem{Fricke.Klein:2017}
Robert Fricke and Felix Klein, \emph{Lectures on the theory of automorphic
  functions. {V}ol. 2}, CTM. Classical Topics in Mathematics, vol.~4, Higher
  Education Press, Beijing, 2017, Translated from the German original [
  MR0183872] by Arthur M. DuPre. \MR{3838412}

\bibitem{Gardner:1989}
Robert~B. Gardner, \emph{The method of equivalence and its applications},
  CBMS-NSF Regional Conference Series in Applied Mathematics, vol.~58, Society
  for Industrial and Applied Mathematics (SIAM), Philadelphia, PA, 1989.
  \MR{MR1062197 (91j:58007)}

\bibitem{Garling:2011}
D.~J.~H. Garling, \emph{Clifford algebras: an introduction}, London
  Mathematical Society Student Texts, vol.~78, Cambridge University Press,
  Cambridge, 2011. \MR{2816665}

\bibitem{Gauld:2014}
David Gauld, \emph{Non-metrisable manifolds}, Springer, Singapore, 2014.
  \MR{3244277}

\bibitem{Goldman:1999}
William~M. Goldman, \emph{Complex hyperbolic geometry}, Oxford Mathematical
  Monographs, The Clarendon Press, Oxford University Press, New York, 1999,
  Oxford Science Publications. \MR{1695450}

\bibitem{Goldman:2009}
\bysame, \emph{Trace coordinates on {F}ricke spaces of some simple hyperbolic
  surfaces}, Handbook of {T}eichm\"{u}ller theory. {V}ol. {II}, IRMA Lect.
  Math. Theor. Phys., vol.~13, Eur. Math. Soc., Z\"{u}rich, 2009, pp.~611--684.
  \MR{2497777}

\bibitem{Goldman:2010}
\bysame, \emph{Locally homogeneous geometric manifolds}, Proceedings of the
  {I}nternational {C}ongress of {M}athematicians. {V}olume {II} (New Delhi),
  Hindustan Book Agency, 2010, pp.~717--744. \MR{2827816 (2012g:57034)}

\bibitem{Goldman:2022}
\bysame, \emph{Geometric structures on manifolds}, Graduate Studies in
  Mathematics, vol. 227, American Mathematical Society, Providence, RI, [2022]
  \copyright 2022. \MR{4500072}

\bibitem{Griffiths.Harris:1978}
Phillip Griffiths and Joseph Harris, \emph{Principles of algebraic geometry},
  Pure and Applied Mathematics, Wiley-Interscience [John Wiley \& Sons], New
  York, 1978. \MR{507725}

\bibitem{Gromov1988}
Michael Gromov, \emph{Rigid transformations groups}, G\'{e}om\'{e}trie
  diff\'{e}rentielle ({P}aris, 1986), Travaux en Cours, vol.~33, Hermann,
  Paris, 1988, pp.~65--139. \MR{955852}

\bibitem{Gromov:2007}
Misha Gromov, \emph{Metric structures for {R}iemannian and non-{R}iemannian
  spaces}, Modern Birkh\"{a}user Classics, Birkh\"{a}user Boston, Inc., Boston,
  MA, 2007, Based on the 1981 French original, With appendices by M. Katz, P.
  Pansu and S. Semmes, Translated from the French by Sean Michael Bates.
  \MR{2307192}

\bibitem{Grosshans:1997}
Frank~D. Grosshans, \emph{Algebraic homogeneous spaces and invariant theory},
  Lecture Notes in Mathematics, vol. 1673, Springer-Verlag, Berlin, 1997,
  \url{https://doi-org.ucc.idm.oclc.org/10.1007/BFb0093525}. \MR{1489234}

\bibitem{Guediri1995}
M.~Guediri and J.~Lafontaine, \emph{Sur la complétude des variétés
  pseudo-riemanniennes}, Journal of Geometry and Physics \textbf{15} (1995),
  no.~2, 150--158.

\bibitem{Guillemin.Pollack.2010}
Victor Guillemin and Alan Pollack, \emph{Differential topology}, AMS Chelsea
  Publishing, Providence, RI, 2010, Reprint of the 1974 original. \MR{2680546}

\bibitem{Hammerl:2007}
Matthias Hammerl, \emph{Homogeneous {C}artan geometries}, Arch. Math. (Brno)
  \textbf{43} (2007), no.~5, 431--442. \MR{2381786}

\bibitem{Harris:1995}
Joe Harris, \emph{Algebraic geometry}, Graduate Texts in Mathematics, vol. 133,
  Springer-Verlag, New York, 1995, A first course, Corrected reprint of the
  1992 original. \MR{1416564}

\bibitem{Helgason:1978}
Sigurdur Helgason, \emph{Differential geometry, {L}ie groups, and symmetric
  spaces}, Graduate Studies in Mathematics, vol.~34, American Mathematical
  Society, Providence, RI, 2001, Corrected reprint of the 1978 original.
  \MR{MR1834454 (2002b:53081)}

\bibitem{Hilbert1971}
David Hilbert, \emph{Foundations of geometry}, Second edition. Translated from
  the tenth German edition by Leo Unger, Open Court, LaSalle, Ill., 1971.
  \MR{0275262 (43 \#1019)}

\bibitem{Hilgert.Neeb:2012}
Joachim Hilgert and Karl-Hermann Neeb, \emph{Structure and geometry of {L}ie
  groups}, Springer Monographs in Mathematics, Springer, New York, 2012,
  \url{http://0-dx.doi.org.library.ucc.ie/10.1007/978-0-387-84794-8}.
  \MR{3025417}

\bibitem{Hubbard:1981}
John~H. Hubbard, \emph{The monodromy of projective structures}, Riemann
  surfaces and related topics: {P}roceedings of the 1978 {S}tony {B}rook
  {C}onference ({S}tate {U}niv. {N}ew {Y}ork, {S}tony {B}rook, {N}.{Y}., 1978),
  Ann. of Math. Stud., vol.~97, Princeton Univ. Press, Princeton, N.J., 1981,
  pp.~257--275. \MR{624819}

\bibitem{Humphreys:1975}
James~E. Humphreys, \emph{Linear algebraic groups}, Springer-Verlag, New
  York-Heidelberg, 1975, Graduate Texts in Mathematics, No. 21. \MR{0396773 (53
  \#633)}

\bibitem{Ivey/Landsberg:2003}
Thomas~A. Ivey and J.~M. Landsberg, \emph{Cartan for beginners: differential
  geometry via moving frames and exterior differential systems}, American
  Mathematical Society, Providence, RI, 2003.

\bibitem{Jacobowitz1990}
Howard Jacobowitz, \emph{An introduction to {CR} structures}, Mathematical
  Surveys and Monographs, vol.~32, American Mathematical Society, Providence,
  RI, 1990. \MR{1067341}

\bibitem{Joyce:2007}
Dominic~D. Joyce, \emph{Riemannian holonomy groups and calibrated geometry},
  Oxford Graduate Texts in Mathematics, vol.~12, Oxford University Press,
  Oxford, 2007. \MR{2292510}

\bibitem{Kapovich:1989}
M.~Kapovich, \emph{Deformation spaces of flat conformal structures},
  Proceedings of the {S}econd {S}oviet-{J}apan {J}oint {S}ymposium of
  {T}opology ({K}habarovsk, 1989), vol.~8, {S}oviet-{J}apan {J}oint {S}ymposium
  of {T}opology, no.~1, 1990, pp.~253--264. \MR{1043223}

\bibitem{Kazdan}
Jerry~L. Kazdan, \emph{Applications of partial differential equations to
  problems in geometry}, to appear.

\bibitem{Klein:2008}
Felix~C. Klein, \emph{A comparative review of recent researches in geometry},
  2008, \url{https://arxiv.org/abs/0807.3161}.

\bibitem{Klingen:1990}
Helmut Klingen, \emph{Introductory lectures on {S}iegel modular forms},
  Cambridge Studies in Advanced Mathematics, vol.~20, Cambridge University
  Press, Cambridge, 1990. \MR{1046630}

\bibitem{Klingenberg:1995}
Wilhelm P.~A. Klingenberg, \emph{Riemannian geometry}, second ed., De Gruyter
  Studies in Mathematics, vol.~1, Walter de Gruyter \& Co., Berlin, 1995.
  \MR{1330918}

\bibitem{Klingler:1996}
Bruno Klingler, \emph{Compl\'{e}tude des vari\'{e}t\'{e}s lorentziennes \`a
  courbure constante}, Math. Ann. \textbf{306} (1996), no.~2, 353--370.
  \MR{1411352}

\bibitem{Kobayashi:1954}
Sh{\^o}shichi Kobayashi, \emph{Espaces \`a connexion de {C}artan complets},
  Proc. Japan Acad. \textbf{30} (1954), 709--710. \MR{MR0069570 (16,1053d)}

\bibitem{Kobayashi:1956}
Sh\^{o}shichi Kobayashi, \emph{On connections of {C}artan}, Canadian J. Math.
  \textbf{8} (1956), 145--156. \MR{77978}

\bibitem{Kobayashi:1957}
\bysame, \emph{Theory of connections}, Ann. Mat. Pura Appl. (4) \textbf{43}
  (1957), 119--194. \MR{96276}

\bibitem{Kobayashi:1995}
Sh{\^o}shichi Kobayashi, \emph{Transformation groups in differential geometry},
  Springer-Verlag, Berlin, 1995, Reprint of the 1972 edition. \MR{96c:53040}

\bibitem{KobayashiNagano:1964}
Sh{\^o}shichi Kobayashi and Tadashi Nagano, \emph{On projective connections},
  J. Math. Mech. \textbf{13} (1964), 215--235. \MR{28 \#2501}

\bibitem{Kobayashi/Nomizu:1996}
Sh{\^o}shichi Kobayashi and Katsumi Nomizu, \emph{Foundations of differential
  geometry. {V}ol. {I}}, Wiley Classics Library, John Wiley \& Sons, Inc., New
  York, 1996, Reprint of the 1963 original, A Wiley-Interscience Publication.
  \MR{1393940}

\bibitem{Komrakov.Churyumov.Doubrov.1993}
B.~Komrakov, A.~Churyumov, and B.~Doubrov, \emph{Two-dimensional homogeneous
  spaces}, Pure Mathematics preprint of Matematisk Institutt, Universitetet i
  Oslo, \url{https://www.duo.uio.no/handle/10852/43320}, 1993.

\bibitem{Kruglikov.Lychagin:2016}
Boris Kruglikov and Valentin Lychagin, \emph{Global {L}ie-{T}resse theorem},
  Selecta Math. (N.S.) \textbf{22} (2016), no.~3, 1357--1411. \MR{3518554}

\bibitem{Kruglikov.The:2017}
Boris Kruglikov and Dennis The, \emph{The gap phenomenon in parabolic
  geometries}, J. Reine Angew. Math. \textbf{723} (2017), 153--215,
  \url{https://doi-org.ucc.idm.oclc.org/10.1515/crelle-2014-0072}. \MR{3604980}

\bibitem{Kuga:1993}
Michio Kuga, \emph{Galois' dream: group theory and differential equations},
  Birkh\"{a}user Boston, Inc., Boston, MA, 1993, Translated from the 1968
  Japanese original by Susan Addington and Motohico Mulase. \MR{1199112}

\bibitem{Kunz:2013}
Ernst Kunz, \emph{Introduction to commutative algebra and algebraic geometry},
  Modern Birkh\"auser Classics, Birkh\"auser/Springer, New York, 2013,
  Translated from the 1980 German original [MR0562105] by Michael Ackerman,
  With a preface by David Mumford, Reprint of the 1985 edition [MR0789602].
  \MR{2977456}

\bibitem{Lee:2013}
John~M. Lee, \emph{Introduction to smooth manifolds}, second ed., Graduate
  Texts in Mathematics, vol. 218, Springer, New York, 2013. \MR{2954043}

\bibitem{Loray/MarinPerez:2009}
Frank Loray and David Mar{\'{\i}}n~P{\'e}rez, \emph{Projective structures and
  projective bundles over compact {R}iemann surfaces}, Ast\'erisque (2009),
  no.~323, 223--252. \MR{2647972}

\bibitem{Lotta2004}
Antonio Lotta, \emph{On model mutation for reductive {C}artan geometries and
  non-existence of {C}artan space forms}, Kodai Math. J. \textbf{27} (2004),
  no.~2, 174--188,
  \url{https://doi-org.ucc.idm.oclc.org/10.2996/kmj/1093351324}. \MR{2069768}

\bibitem{Maclane:1995}
Saunders Mac~Lane, \emph{Homology}, Classics in Mathematics, Springer-Verlag,
  Berlin, 1995, Reprint of the 1975 edition. \MR{1344215}

\bibitem{Maskit:1988}
Bernard Maskit, \emph{Kleinian groups}, Grundlehren der mathematischen
  Wissenschaften [Fundamental Principles of Mathematical Sciences], vol. 287,
  Springer-Verlag, Berlin, 1988. \MR{959135}

\bibitem{Massey:1967}
William~S. Massey, \emph{Algebraic topology: {A}n introduction}, Harcourt,
  Brace \& World, Inc., New York, 1967. \MR{35 \#2271}

\bibitem{Matveev.Troyanov:2015}
Vladimir~S. Matveev and Marc Troyanov, \emph{Completeness and incompleteness of
  the {B}inet-{L}egendre metric}, Eur. J. Math. \textbf{1} (2015), no.~3,
  483--502, \url{https://doi-org.ucc.idm.oclc.org/10.1007/s40879-015-0046-4}.
  \MR{3401902}

\bibitem{Matveev.Troyanov:2017}
\bysame, \emph{The {M}yers-{S}teenrod theorem for {F}insler manifolds of low
  regularity}, Proc. Amer. Math. Soc. \textbf{145} (2017), no.~6, 2699--2712,
  \url{https://doi-org.ucc.idm.oclc.org/10.1090/proc/13407}. \MR{3626522}

\bibitem{McKay:diff.geom}
Benjamin McKay, \emph{Lectures on differential geometry},
  \url{https://ben-mckay.github.io/benmckay.github.io/}.

\bibitem{McKay:2004a}
\bysame, \emph{Sussmann's orbit theorem and maps}, Differential Geom. Appl.
  \textbf{25} (2007), no.~3, 277--280,
  \url{https://doi.org/10.1016/j.difgeo.2006.11.005}. \MR{2330456}

\bibitem{McKay2011b}
\bysame, \emph{Holomorphic {C}artan geometries on uniruled surfaces}, C. R.
  Math. Acad. Sci. Paris \textbf{349} (2011), no.~15-16, 893--896. \MR{2835898}

\bibitem{McKay:2016}
\bysame, \emph{Holomorphic geometric structures on {K}\"ahler-{E}instein
  manifolds}, Manuscripta Math. \textbf{153} (2017), no.~1-2, 1--34.
  \MR{3635971}

\bibitem{McKay}
\bysame, \emph{Introduction to exterior differential systems}, unpublished,
  2022.

\bibitem{Meinrenken:2003}
Eckhard Meinrenken, \emph{Group actions on manifolds},
  \url{www.math.toronto.edu/mein/teaching/LectureNotes/action.pdf}, 2003.

\bibitem{Meise.Vogt:1997}
Reinhold Meise and Dietmar Vogt, \emph{Introduction to functional analysis},
  Oxford Graduate Texts in Mathematics, vol.~2, The Clarendon Press, Oxford
  University Press, New York, 1997, Translated from the German by M. S.
  Ramanujan and revised by the authors. \MR{1483073}

\bibitem{Melnick:2011}
Karin Melnick, \emph{A {F}robenius theorem for {C}artan geometries, with
  applications}, Enseign. Math. (2) \textbf{57} (2011), no.~1-2, 57--89,
  \url{http://dx.doi.org/10.4171/LEM/57-1-3}. \MR{2850584}

\bibitem{Melnick:2021}
\bysame, \emph{Rigidity of transformation groups in differential geometry},
  Notices Amer. Math. Soc. \textbf{68} (2021), no.~5, 721--732,
  \url{https://doi-org.ucc.idm.oclc.org/10.1090/noti2279}. \MR{4249428}

\bibitem{Mimura/Toda:1991}
Mamoru Mimura and Hirosi Toda, \emph{Topology of {L}ie groups, i and ii},
  Translations of Mathematical Monographs, vol.~91, American Mathematical
  Society, Providence, Rhode Island, 1991.

\bibitem{Misner/Thorne/Wheeler:1973}
Charles~W. Misner, Kip~S. Thorne, and John~Archibald Wheeler,
  \emph{Gravitation}, W. H. Freeman and Co., San Francisco, CA, 1973.
  \MR{418833}

\bibitem{Molzon/Mortensen:1996}
Robert Molzon and Karen~Pinney Mortensen, \emph{The {S}chwarzian derivative for
  maps between manifolds with complex projective connections}, Trans. Amer.
  Math. Soc. \textbf{348} (1996), no.~8, 3015--3036. \MR{MR1348154 (96j:32028)}

\bibitem{Montgomery.Zippin:1974}
Deane Montgomery and Leo Zippin, \emph{Topological transformation groups},
  Robert E. Krieger Publishing Co., Huntington, N.Y., 1974, Reprint of the 1955
  original. \MR{0379739}

\bibitem{Mostow:1950}
George~Daniel Mostow, \emph{The extensibility of local {L}ie groups of
  transformations and groups on surfaces}, Ann. of Math. (2) \textbf{52}
  (1950), 606--636. \MR{MR0048464 (14,18d)}

\bibitem{Munkres:2000}
James~R. Munkres, \emph{Topology}, Prentice Hall, Inc., Upper Saddle River, NJ,
  2000, Second edition of [ MR0464128]. \MR{3728284}

\bibitem{Myers.Steenrod:1939}
S.~B. Myers and N.~E. Steenrod, \emph{The group of isometries of a {R}iemannian
  manifold}, Ann. of Math. (2) \textbf{40} (1939), no.~2, 400--416.
  \MR{1503467}

\bibitem{Penis.Cola:1986}
A.~M. Naveira, A.~Ferr\'{a}ndez, and F.~Mascar\'{o} (eds.), \emph{Differential
  geometry, {P}e\~{n}\'{\i}scola 1985}, Lecture Notes in Mathematics, vol.
  1209, Springer-Verlag, Berlin, 1986. \MR{863741}

\bibitem{Needham:2023}
Tristan Needham, \emph{Visual complex analysis}, anniversary ed., Oxford
  University Press, Oxford, 2023, With a foreword by Roger Penrose.
  \MR{4577812}

\bibitem{Nomizu:1954}
Katsumi Nomizu, \emph{Invariant affine connections on homogeneous spaces},
  Amer. J. Math. \textbf{76} (1954), 33--65,
  \url{https://doi-org.ucc.idm.oclc.org/10.2307/2372398}. \MR{59050}

\bibitem{Nomizu:1978}
\bysame, \emph{Kinematics and differential geometry of submanifolds}, Tohoku
  Math. J. (2) \textbf{30} (1978), no.~4, 623--637, Rolling a ball with a
  prescribed locus of contact. \MR{516894}

\bibitem{Olver:1995}
Peter~J. Olver, \emph{Equivalence, invariants, and symmetry}, Cambridge
  University Press, Cambridge, 1995. \MR{96i:58005}

\bibitem{ONeill:1983}
Barrett O'Neill, \emph{Semi-{R}iemannian geometry}, Pure and Applied
  Mathematics, vol. 103, Academic Press, Inc. [Harcourt Brace Jovanovich,
  Publishers], New York, 1983, With applications to relativity. \MR{719023}

\bibitem{OnishchikVinberg:1990}
A.~L. Onishchik and {\`E}.~B. Vinberg, \emph{Lie groups and algebraic groups},
  Springer Series in Soviet Mathematics, Springer-Verlag, Berlin, 1990,
  Translated from the Russian and with a preface by D. A. Leites.
  \MR{91g:22001}

\bibitem{Palais:1957}
Richard~S. Palais, \emph{A global formulation of the {L}ie theory of
  transformation groups}, Mem. Amer. Math. Soc. No. \textbf{22} (1957),
  iii+123. \MR{22 \#12162}

\bibitem{Palais:1961}
\bysame, \emph{On the existence of slices for actions of non-compact {L}ie
  groups}, Ann. of Math. (2) \textbf{73} (1961), 295--323. \MR{126506}

\bibitem{Papadopoulos.Troyanov.2014}
Athanase Papadopoulos and Marc Troyanov (eds.), \emph{Handbook of {H}ilbert
  geometry}, IRMA Lectures in Mathematics and Theoretical Physics, vol.~22,
  European Mathematical Society (EMS), Z\"{u}rich, 2014. \MR{3309067}

\bibitem{Pecastaing:2016}
Vincent Pecastaing, \emph{On two theorems about local automorphisms of
  geometric structures}, Ann. Inst. Fourier (Grenoble) \textbf{66} (2016),
  no.~1, 175--208, \url{http://aif.cedram.org/item?id=AIF_2016__66_1_175_0}.
  \MR{3477874}

\bibitem{Penrose:2005}
Roger Penrose, \emph{The road to reality}, Alfred A. Knopf, Inc., New York,
  2005, A complete guide to the laws of the universe. \MR{2116746}

\bibitem{Petersen:2016}
Peter Petersen, \emph{Riemannian geometry}, third ed., Graduate Texts in
  Mathematics, vol. 171, Springer, Cham, 2016. \MR{3469435}

\bibitem{Procesi:2007}
Claudio Procesi, \emph{Lie groups}, Universitext, Springer, New York, 2007, An
  approach through invariants and representations. \MR{MR2265844 (2007j:22016)}

\bibitem{Reid:1988}
Miles Reid, \emph{Undergraduate algebraic geometry}, London Mathematical
  Society Student Texts, vol.~12, Cambridge University Press, Cambridge, 1988.
  \MR{982494}

\bibitem{Rosenlicht:1963}
Maxwell Rosenlicht, \emph{A remark on quotient spaces}, An. Acad. Brasil. Ci.
  \textbf{35} (1963), 487--489. \MR{0171782}

\bibitem{Rudin:1987}
Walter Rudin, \emph{Real and complex analysis}, third ed., McGraw-Hill Book
  Co., New York, 1987. \MR{924157}

\bibitem{Rudin:2008}
\bysame, \emph{Function theory in the unit ball of {$\mathbb{C}^n$}}, Classics
  in Mathematics, Springer-Verlag, Berlin, 2008, Reprint of the 1980 edition.
  \MR{2446682}

\bibitem{Sard:1965}
Arthur Sard, \emph{Hausdorff measure of critical images on {B}anach manifolds},
  Amer. J. Math. \textbf{87} (1965), 158--174. \MR{173748}

\bibitem{Schmidt:1974}
B.~G. Schmidt, \emph{A new definition of conformal and projective infinity of
  space-times}, Comm. Math. Phys. \textbf{36} (1974), 73--90,
  \url{http://projecteuclid.org.ucc.idm.oclc.org/euclid.cmp/1103859662}.
  \MR{339775}

\bibitem{Schmidt:1979}
\bysame, \emph{Remarks about modifications of the {$b$}-boundary definition},
  Gen. Relativity Gravitation \textbf{10} (1979), no.~12, 981--982,
  \url{https://doi-org.ucc.idm.oclc.org/10.1007/BF00776515}. \MR{550361}

\bibitem{Shafarevich:1989}
I.~R. Shafarevich (ed.), \emph{Algebraic geometry. {IV}}, Encyclopaedia of
  Mathematical Sciences, vol.~55, Springer-Verlag, Berlin, 1994, Linear
  algebraic groups. Invariant theory, A translation of {\emph{Algebraic
  geometry. 4}} (Russian), Akad. Nauk SSSR Vsesoyuz. Inst. Nauchn. i Tekhn.
  Inform., Moscow, 1989 [ MR1100483 (91k:14001)], Translation edited by A. N.
  Parshin and I. R. Shafarevich.
  \url{https://doi-org.ucc.idm.oclc.org/10.1007/978-3-662-03073-8}.
  \MR{1309681}

\bibitem{Shafarevich:2013}
Igor~R. Shafarevich, \emph{Basic algebraic geometry. 1}, third ed., Springer,
  Heidelberg, 2013, Varieties in projective space. \MR{3100243}

\bibitem{Sharpe:1997}
Richard~W. Sharpe, \emph{Differential geometry}, Graduate Texts in Mathematics,
  vol. 166, Springer-Verlag, New York, 1997, Cartan's generalization of Klein's
  Erlangen program, With a foreword by S. S. Chern. \MR{98m:53033}

\bibitem{Sharpe:2002}
\bysame, \emph{An introduction to {C}artan geometries}, Proceedings of the 21st
  Winter School ``Geometry and Physics'' (Srn\'\i, 2001), no.~69, 2002,
  pp.~61--75. \MR{2004f:53023}

\bibitem{Siegel:1943}
Carl~Ludwig Siegel, \emph{Discontinuous groups}, Annals of Mathematics
  \textbf{44} (1943), no.~4, 674--689.

\bibitem{Siegel:3}
\bysame, \emph{Topics in complex function theory. {V}ol. {III}}, Wiley Classics
  Library, John Wiley \& Sons, Inc., New York, 1989, Abelian functions and
  modular functions of several variables, Translated from the German by E.
  Gottschling and M. Tretkoff, With a preface by Wilhelm Magnus, Reprint of the
  1973 original, A Wiley-Interscience Publication. \MR{1013364}

\bibitem{Sikora:2012}
Adam~S. Sikora, \emph{Character varieties}, Trans. Amer. Math. Soc.
  \textbf{364} (2012), no.~10, 5173--5208,
  \url{https://doi-org.ucc.idm.oclc.org/10.1090/S0002-9947-2012-05448-1}.
  \MR{2931326}

\bibitem{Singer.Sternberg:1965}
I.~M. Singer and Shlomo Sternberg, \emph{The infinite groups of {L}ie and
  {C}artan. {I}. {T}he transitive groups}, J. Analyse Math. \textbf{15} (1965),
  1--114. \MR{217822}

\bibitem{Sontag:1998}
Eduardo~D. Sontag, \emph{Mathematical control theory}, second ed., Texts in
  Applied Mathematics, vol.~6, Springer-Verlag, New York, 1998, Deterministic
  finite-dimensional systems.
  \url{https://doi-org.ucc.idm.oclc.org/10.1007/978-1-4612-0577-7}.
  \MR{1640001}

\bibitem{Spivak:1979b}
Michael Spivak, \emph{A comprehensive introduction to differential geometry.
  {V}ol. {II}}, second ed., Publish or Perish Inc., Wilmington, Del., 1979.
  \MR{82g:53003b}

\bibitem{Steenrod:1999}
Norman Steenrod, \emph{The topology of fibre bundles}, Princeton University
  Press, Princeton, NJ, 1999, Reprint of the 1957 edition, Princeton
  Paperbacks. \MR{2000a:55001}

\bibitem{Sternberg:1983}
Shlomo Sternberg, \emph{Lectures on differential geometry}, second ed., Chelsea
  Publishing Co., New York, 1983, With an appendix by Sternberg and Victor W.
  Guillemin. \MR{MR891190 (88f:58001)}

\bibitem{Sullivan/Thurston:1983}
Dennis Sullivan and William Thurston, \emph{Manifolds with canonical coordinate
  charts: some examples}, Enseign. Math. (2) \textbf{29} (1983), no.~1-2,
  15--25. \MR{702731}

\bibitem{Sussmann:1973}
H{\'e}ctor~J. Sussmann, \emph{Orbits of families of vector fields and
  integrability of distributions}, Trans. Amer. Math. Soc. \textbf{180} (1973),
  171--188. \MR{47 \#9666}

\bibitem{tomDieck:2008}
Tammo tom Dieck, \emph{Algebraic topology}, EMS Textbooks in Mathematics,
  European Mathematical Society (EMS), Z\"{u}rich, 2008,
  \url{https://doi-org.ucc.idm.oclc.org/10.4171/048}. \MR{2456045}

\bibitem{Wang1958}
Hsien-Chung Wang, \emph{On invariant connections over a principal fibre
  bundle}, Nagoya Math. J. \textbf{13} (1958), 1--19.

\bibitem{chriswendl}
Chris Wendl, \emph{Palais's and {K}obayashi's theorems on automorphism groups
  of geometric structures}, MathOverflow, \url{mathoverflow.net/q/417762}
  (version: 2022-03-09).

\bibitem{Weyl1938}
Hermann Weyl, \emph{Book {R}eview: {L}a {T}h\'{e}orie des {G}roupes {F}inis et
  {C}ontinus et la {G}\'{e}om\'{e}trie {D}iff\'{e}rentielle trait\'{e}es par la
  {M}\'{e}thode du {R}ep\`ere {M}obile}, Bull. Amer. Math. Soc. \textbf{44}
  (1938), no.~9, 598--601,
  \url{https://doi-org.ucc.idm.oclc.org/10.1090/S0002-9904-1938-06789-4}.
  \MR{1563803}

\bibitem{Whitney:1934a}
Hassler Whitney, \emph{Analytic extensions of differentiable functions defined
  in closed sets}, Trans. Amer. Math. Soc. \textbf{36} (1934), no.~1, 63--89.
  \MR{1501735}

\bibitem{Whitney:1934b}
\bysame, \emph{Differentiable functions defined in closed sets. {I}}, Trans.
  Amer. Math. Soc. \textbf{36} (1934), no.~2, 369--387. \MR{1501749}

\bibitem{Whitney:1934c}
\bysame, \emph{Functions differentiable on the boundaries of regions}, Ann. of
  Math. (2) \textbf{35} (1934), no.~3, 482--485. \MR{1503174}

\bibitem{Whitney:1957}
\bysame, \emph{Elementary structure of real algebraic varieties}, Ann. of Math.
  (2) \textbf{66} (1957), 545--556,
  \url{https://doi-org.ucc.idm.oclc.org/10.2307/1969908}. \MR{95844}

\bibitem{Wolf:1967}
Joseph~A. Wolf, \emph{Spaces of constant curvature}, sixth ed., AMS Chelsea
  Publishing, Providence, RI, 2011. \MR{2742530}

\bibitem{Yekutieli_2019}
Amnon Yekutieli, \emph{Derived categories}, Cambridge University Press,
  November 2019.

\bibitem{Zeghib:2000}
Abdelghani Zeghib, \emph{On {G}romov's theory of rigid transformation groups: a
  dual approach}, Ergodic Theory Dynam. Systems \textbf{20} (2000), no.~3,
  935--946, \url{https://doi-org.ucc.idm.oclc.org/10.1017/S0143385700000511}.
  \MR{1764937}

\end{thebibliography}
\end{document}